\documentclass[amssymb,amsfonts,reqno]{amsart}

 \usepackage{xcolor}
 \usepackage[citebordercolor={green}]{hyperref}

\usepackage{color}
\usepackage{latexsym}
\usepackage{amsmath}
\usepackage{amssymb}
\usepackage{amsfonts}
\usepackage{esint}
\usepackage{graphicx}
\usepackage{accents}
\newcommand{\Lie}{{\mathcal{L}}}

\newcommand{\der}{\nabla}

\newcommand{\les}{\lesssim}
\newcommand{\bea}{\begin{eqnarray}}
\newcommand{\eea}{\end{eqnarray}}

\newcommand{\derm}{ { \der^{(\bf{m})}} }
\newcommand{\rderm}{   {\mbox{$\nabla \mkern-13mu /$\,}^{(\bf{m})} }   }

\newcommand{\rpa}{\mbox{${\pa} \mkern-09mu /$}}

\newcommand{\eps}{{\varepsilon}}\newcommand{\R}{{\mathbb R}}

\newcommand{\E}{{\cal E}}

\newcommand{\la}{\langle}\newcommand{\si}{\sigma}\renewcommand{\b}{\beta}

\newcommand{\cal}{\mathcal}

\def\a{\alpha}\def\ga{\gamma}\def\de{\delta}
\def\bm{\left( \begin{array}{cc}}
\def\endm{\end{array}\right)}\newcommand{\eq}{\end{equation}}
\def\a{\alpha}\def\b{\beta}
\def\ga{\gamma}\def\de{\delta}\def\Box{\square}\def\pa{\partial}

\def \rectangle#1#2{\hbox{\vrule\vbox to #2 {\hrule\hbox to #1{\hfil}\vfil\hrule}\vrule}}
\def\Lb{\underline{L}}
\def\a{\alpha}\def\b{\beta}\def\ga{\gamma}
\def\de{\delta}\def\Box{\square}\def\pa{\partial}
\def\Lb{\underline{L}}

\def\pa{\partial}

\def\beaa{\begin{eqnarray*}}
\def\eeaa{\end{eqnarray*}}
\def\pa{\partial}
\def\a{{\alpha}}
\def\b{{\beta}}
\def\ga{\gamma}

\def\de{\delta}

\def\eps{\epsilon}
\def\la{\lambda}

\def\si{\sigma}

\def\Om{\Omega}

\def\th{\theta}

\def\Lb{{\underline{L}}}

\def\g{{\bf g}}

\def\SSS{{\Bbb S}}

\def\nn{{\mathbb N}}

\def\R{{\mathbb R}}

\def\12{\frac{1}{2}}

\def\N{{\mathcal N}}

\def\bep{\begin{proposition}}
\def\eep{\end{proposition}}

\def\4{\frac{1}{4}}

\def\12{\frac{1}{2}}

\def\N{\nn}

\def\bep{\begin{proposition}}
\def\eep{\end{proposition}}

\def\bm#1{\boldsymbol{#1}} 

\def\nin{\notin}

\def\nin{\notin}

\def\build#1_#2^#3{\mathrel{\mathop{\kern 0pt#1}\limits_{#2}^{#3}}}
\def\4{\frac{1}{4}}

\def\<{\langle}
\def\>{\rangle}

\theoremstyle{plain}
\newtheorem{theorem}{Theorem}
\newtheorem{proposition}{Proposition}
\newtheorem{lemma}{Lemma}
\newtheorem{corollary}{Corollary}

\theoremstyle{remark}
\newtheorem{remark}{Remark}

\theoremstyle{definition}
\newtheorem{definition}{Definition}
\numberwithin{equation}{section}
\numberwithin{proposition}{section}
\numberwithin{definition}{section}
\numberwithin{lemma}{section}
\numberwithin{corollary}{section}
\numberwithin{remark}{section}
\begin{document}
\include{psfig}
\title[Exterior Stability of Minkowski for Einstein-Yang-Mills]{Exterior stability of the $\bf{(1+3)}$-dimensional Minkowski space-time solution to the Einstein-Yang-Mills equations}
\author{Sari Ghanem}
\address{University of Lübeck}
\email{sari.ghanem@uni-luebeck.de}
\maketitle

\begin{abstract}
We prove the exterior stability of the Minkowski space-time, $\R^{1+3}$\,, solution to the Einstein-Yang-Mills system in both the Lorenz and harmonic gauges, where the Yang-Mills fields are valued in any arbitrary Lie algebra $\cal G$\,, associated to any compact Lie group $G$\,. We start with an arbitrary sufficiently small initial data, defined in a suitable energy norm for the perturbations of the Yang-Mills potential and of the Minkowski space-time, and we show the well-posedness of the Cauchy development in the exterior of the fully coupled Einstein-Yang-Mills equations in the Lorenz gauge and in wave coordinates, and we prove that this leads to solutions converging to the zero Yang-Mills curvature and to the Minkowski space-time. Furthermore, we obtain dispersive estimates in wave coordinates on the Yang-Mills potential in the Lorenz gauge and on the metric, as well as on the gauge invariant norm of the Yang-Mills curvature. This provides a new proof to the exterior stability result by P. Mondal and S-T. Yau, based on an alternative approach, by using a null frame decomposition that was first used by H. Lindblad and I. Rodnianski for the case of the Einstein vacuum equations. In this third paper of a series, we detail all the new material concerning our proof so as to provide lecture notes for Ph.D. students wanting to learn non-linear hyperbolic differential equations and stability problems in mathematical General Relativity.
\end{abstract}

\setcounter{page}{1}
\pagenumbering{arabic}


\section{Brief summary of the main novelties of this paper}\label{BrielfSummary}

We start here with a brief summary of the main novelties meant for specialised experts: readers who are unfamiliar with the subject can start by reading the introduction in Section \ref{Introduction}, and can skip for the moment this short summary, and should get first somewhat familiar with the introductions of our previous papers \cite{G4} and \cite{G5}.

We shall mention, in this short section, difficulties that do not exist at all neither for the Einstein vacuum equations nor for the Einstein-Maxwell equations. Let us first mention quickly, nevertheless, that the argument that we are going to run to prove the non-linear stability of the fully coupled Einstein-Yang-Mills system, is a continuity argument on the energy, defined in \eqref{definitionoftheenergynorm}: we start with an a prior estimate (see \eqref{aprioriestimate}) and we want to upgrade this estimate to ultimately prove that it is a true estimate.

\subsection{The system of non-linear hyperbolic partial differential equations}\

Using the Lorenz gauge condition, that is a condition that breaks the gauge invariance of the Einstein-Yang-Mills equations (see proofs of Lemmas \ref{estimateoncovariantgradientofAL} and \ref{LorenzgaugeestimateforgradientofLiederivativesofAL}), where one can see that these are not only estimates, but one can express literally the terms on the left hand side of the inequality as a combination of terms in the right hand side of the inequality), and using the wave coordinates, that is a system of coordinates for the Einstein-Yang-Mills equations satisfying some condition, namely the harmonic gauge condition, that breaks the diffeomorphism invariance of the Einstein-Yang-Mills equations (see proofs of Lemmas \ref{waveconditionestimateonzeroLiederivativeofmetric} and \ref{wavecoordinatesestimateonLiederivativesZonmetric}), one can then write the Einstein-Yang-Mills equations as a coupled system of non-linear hyperbolic equations on the Einstein-Yang-Mills potential $A$\;, that is a one-tensor valued in the Lie algebra ${\cal G}$\;, and on the perturbation metric $h^1 = g - m - h^0$\;, that is a scalar valued two-tensor, where $g$ is the unknown metric solution to the evolution problem for the Einstein-Yang-Mills system, and where $h^0$ is the spherically symmetric Schwarzschildian part of the perturbation.

In addition, using a null-frame tetrad $\{\, \underline{L}\,, L\,, e_a\,, a \in \{1, 2 \}\, \}$\;, that is a frame constructed using wave coordinates (see Definition \ref{definitionofthenullframusingwavecoordinates}), where $\cal T:= \{\, L\,, e_a\,, a \in \{1, 2 \}\, \}$ is a set of three vectors of the frame tangent to the outgoing null-cone for the Minkowski metric $m$\;, where $m$ is defined to be Minkowski metric in wave coordinates, and where $\cal U := \{\, \underline{L}\,, {\cal T}\,  \}$ are the full components of the frame, we get a coupled system of non-linear wave equations on the Einstein-Yang-Mills potential $A$ and on the perturbation metric $h^1$, that is schematically speaking the following (see Subsection \ref{StructureofthewaveequationforthebadcomponentsofAandhone}):
\bea\label{waveequationonbadtermfortheEinsteinYangMillspoential}
   \notag
 &&   g^{\la\mu} \derm_{\la}   \derm_{\mu}   A_{{\underline{L}}}     \\
 \notag
    &=& \derm h  \cdot  \rderm A        + \rderm  h \cdot  \derm A  +    A   \cdot     \rderm A   +  \derm  h  \cdot  A^2   +  A^3 \\
    \notag
  && + O( h \cdot  \derm h \cdot  \derm A) + O( h \cdot  A \cdot \derm A)  + O( h \cdot  \derm h \cdot  A^2) + O( h \cdot  A^3) \\
         \notag
 && +  A_L   \cdot   \derm A      +  A_{e_a}  \cdot     \derm A_{e_a}     \; ,\\
  \eea
 and
    \bea
\notag
&&    g^{\alpha\beta} \derm_\alpha \derm_\beta h^1_{{\underline{L}} {\underline{L}} }     \\
\notag
            &=&   \rderm h \cdot \derm h +  h \cdot ( \derm h )^2      +   \rderm A    \cdot  \derm  A   +  A^2  \cdot  \derm A     +  A^4 \\
    \notag
     && + O \big(h \cdot  (\derm A)^2 \big)   + O \big(  h  \cdot  A^2 \cdot \derm A \big)     + O \big(  h   \cdot  A^4 \big)   \\
         \notag
 && + (  \derm h_{\cal TU} ) ^2 +   ( \derm A_{e_{a}} )^2    \\
 \notag
 && +     g^{\alpha\beta} \derm_\alpha \derm_\beta h^0 \; . \\
\eea
and for the “good” components of $A$ and $h^1$ (see Subsection \ref{StructureofthewaveequationforthegoodcomponentsofAandhone}), we have the “better” structure that is
  \bea
   \notag
 &&   g^{\la\mu} \derm_{\la}   \derm_{\mu}   A_{{\cal T}}     \\
 \notag
    &=&   \derm h  \cdot  \rderm A      +  \rderm  h  \cdot  \derm A   +    A   \cdot   \rderm A  +  \derm  h  \cdot  A^2   +  A^3 \\
    \notag
  && + O( h \cdot  \derm h \cdot  \derm A) + O( h \cdot  A \cdot \derm A)  + O( h \cdot  \derm h \cdot  A^2) + O( h \cdot  A^3)  \; ,\\
  \eea
  and
 \bea
\notag
 &&   g^{\alpha\beta} \derm_\alpha \derm_\beta h^1_{ {\cal T} {\cal U}}    \\
 \notag
    &=&  \rderm h \cdot \derm h + h \cdot (\derm h)^2      +  \rderm A   \cdot  \derm  A   +    A^2 \cdot \derm A     + A^4  \\
\notag
     && + O \big(h \cdot  (\derm A)^2 \big)   + O \big(  h  \cdot  A^2 \cdot \derm A \big)     + O \big(  h   \cdot  A^4 \big)   \\
  \notag
  && +    g^{\alpha\beta} \derm_\alpha \derm_\beta h^0  \; .\\
\eea

We see in the troublesome wave equation in \eqref{waveequationonbadtermfortheEinsteinYangMillspoential} on $A_{{\underline{L}}}$\;, that there exist “bad” terms which are $ A_{e_a}  \cdot     \derm A_{e_a} $\;, where $a \in \{1, 2 \}$\,, and $ A_L   \cdot   \derm A    $\;.

\subsection{Dealing with the “bad” term $A_{e_a}  \cdot     \derm A_{e_a} $\;}\label{dealingiwthbadtermAaderAa}\
 
 \textbf{The problems:}\
 
 Unlike the case of the Einstein-vacuum equations or even the Einstein-Maxwell system, we need to deal here with not a square of derivatives, but with a product of the derivative and the zeroth derivative of the field. This induces the difficulty that while we have a good decay estimate on $|\derm A_{e_a}|$ (see Lemma \ref{estimatefortangentialcomponentspotential}), we need to also translate it into a good decay estimate on $|A_{e_a}|$\,, and we do not know how to do that without having an $\eps$ loss in the time decay rate for the leading term (that is the zeroth Lie derivative) for which we cannot afford such a loss.
 
 The reason that we would have such a loss in the time decay rate for $| A_{e_a}|$\;, as opposed to $|\derm A_{e_a}|$\;, is because an estimate on the gradient could be translated into an estimate on the field itself, only if we integrate (as far as our knowledge and approach are concerned) along specific lines of integration and use the asymptotic behaviour of the initial data (see the construction of such integration in Definition \ref{definitionequationforintegralalongthenullcoordinateplusboundaryterm} and see Remark \ref{estimategoodcomponentspotentialandmetric}). However, for this to be successful in maintaining for $|A_{e_a}|$ the time decay rate in $t$ that we established originally for $|\derm A_{e_a}|$ in Lemma \ref{estimatefortangentialcomponentspotential}, we need to have more decay in $q:= r - t$ (that measures the distance to the outgoing light cone in Minkowski space-time). But if we want to have a higher decay rate in $|q|$ in the estimate of decay for $|\derm A_{e_a}|$\;, the only estimate that we can think of in order to establish such an upgraded estimate (with more decay in $|q|$) is the weighted estimate established by Lindblad and Rodnianski in Corollary 7.2 of \cite{LR10} (see \eqref{LinfinitynormestimateongradientderivedbyLondbladRodnianski}). The problem is: any weighted estimate, with a non-trivial weight (that is needed here), would necessarily imply no better (as far as our knowledge is concerned) than a Grönwall inequality on the gradient in order to conclude a result that would give the desired upgraded estimate. But such a Grönwall inequality would necessarily mean an $\eps$ loss for the time decay rate in $t$ for  $|\derm A_{e_a}|$ (see Lemma \ref{weightedestimatefortangentialcomponentspotential}), and thus, we would have already ruined our argument. 
 
That was one of the main difficulties for this term (we will discuss later, in this brief summary, another main difficulty for this exact term): it comes from the fact that in order to upgrade the a priori estimate on the energy, using a Grönwall type inequality, we need to control \textit{correctly} the product $| A_{e_a} | \cdot  |   \derm A_{e_a}|$ using the dispersive a prior estimates. By controlling “correctly”, we mean in a way where the bootstrap argument could close, and in order for us to achieve this closure: we need to control each term in the product separately, namely the term $|A_{e_a}|$ and then the term $|\derm A_{e_a}|$\;. Let us explain:

Regarding \textbf{the term $|A_{e_a}|$} in the above mentioned product: we want to use a \textit{suitable} a priori estimate on $|A_{e_a}|$\;, that would give the right factor for the other term $| \derm A_{e_a}|$\;. For this, we \textit{need} to establish a decay in $\frac{1}{t}$ for $|A_{e_a}|$\;. The only way we could think of, given our approach, in order to establish such an estimate for $|A_{e_a}|$\;, is to integrate an a priori estimate on the gradient $|\derm A_{e_a}|$ -- in fact, since the energy for a system of non-linear wave equations is at the level of the gradient, the a priori estimates that we start with are for the gradient. However, if we want to keep the decay in  $\frac{1}{t}$ for $|A_{e_a}|$\;, this means that we need to have a decay not only in $\frac{1}{t}$ for $| \derm A_{e_a}|$\;, but also an additional decay in $|q|$\;, so that when we integrate properly, we could still conserve the decay in $\frac{1}{t}$ for $| A_{e_a}|$ as well.
 
 The problem here is: if we want to have more decay in $|q|$ in the estimate of decay for $|\derm A_{e_a}|$\;, this comes at the expense (at the best of our knowledge) of having an $\eps$ loss in the time decay rate for $|\derm A_{e_a}|$ (see Lemma \ref{weightedestimatefortangentialcomponentspotential}), and this destroys our argument in getting the right factor for $| \derm A_{e_a}|$\;. In fact, we can afford an $\eps$ loss in the time decay rate for the Lie derivatives in the direction of Minkowski vector fields, namely $| \derm  ( \Lie_{Z^I} A_{e_a} ) |$\;, but we \textit{cannot} afford this for the \textit{zeroth} Lie derivative. 

In its turn, \textbf{the term $| \derm A_{e_a}|$} needs to be controlled in a way suitable enough so that we would have the right factor for $|A_{e_a}|$\;. While one can obtain for the term $|\derm A_{e_a}|$\;, a \textit{relatively} good estimate (see Lemma \ref{estimatefortangentialcomponentspotential}), yet it is \textit{not enough} since the other term in the product is $| A_{e_a} |$\,, and \textit{not the gradient} (the gradient is actually the term that we would need to control, given the expression of the energy), and therefore, we need a more powerful estimate on $|\derm A_{e_a}|$ (than the one obtained in Lemma \ref{estimatefortangentialcomponentspotential}) so that we could apply on $| A_{e_a} |$ the Hardy type inequality in the exterior (that we established in Corollary \ref{HardytypeinequalityforintegralstartingatROm}), while keeping the decay rate as $\frac{1}{t}$ (which is the most crucial condition for a successful estimate on the zeroth Lie derivative).

The problem again is: in order to have a factor for $|A_{e_a}|$ with more decay in $|q|$\;, we need to have a weighted estimate in $|q|$ for $|\derm A_{e_a}|$\;, and then integrate properly to translate this into an estimate on $|A_{e_a}|$\;. Again, this would lead only to a Grönwall inequality that would necessarily imply an $\eps$ loss in the time decay rate for $|\derm A_{e_a}|$ (see Lemma \ref{weightedestimatefortangentialcomponentspotential}) that we cannot recover for $|A_{e_a}|$\;.

 \textbf{The solutions:}\
 
 To control correctly the term $ |A_{e_a} | \cdot |   \derm A_{e_a} |$\;, we use first an energy estimate that involves only the component $A_{e_a}$\;, an energy estimate that we established in \cite{G5} (see Lemma \ref{Theveryfinal }). This is unlike the case of the Einstein vacuum equations, where an energy estimate for all components of the metric, without establishing an energy estimate for each component, would be sufficient. In fact, this is a key feature in this paper, that we have to deal separately at several occasions with some tangential components of the Einstein-Yang-Mills potential, namely $A_{e_a}$\;, before dealing with all components of the potential $A$\;, unlike the components of the metric which could be all dealt with together simultaneously at these occasions, which is a key difference not only with the Einstein vacuum equations, but also with the Einstein-Maxwell equations, where such “troublesome” structure does not exist -- we shall explain this later.

Such an energy estimate (see Lemma  \ref{Theveryfinal }), would allow us to control a weighted $L^2$ norm of $|  \derm A_{e_a} |$ (see Lemma \ref{ControlontheweightedLtwonormofdermAea}), by space-time integrals that do \textit{not} involve the “bad” term $ A_{e_a}  \cdot    \derm A_{e_a}$\;, because this “bad” term does not appear in the source terms for the wave equation for the “good” terms  $A_{e_a}$\;, $a \in \{1, 2 \}$\;. This in turn, using a weighted Klainerman-Soboloev estimate in the exterior region (exterior to an outgoing null cone for the Minkowski metric, defined to be the Minkowski metric in wave coordinates), would translate into pointwise decay on $|  \derm A_{e_a} |$ (see Lemma \ref{upgradedestimateontheA_acomponentsthatconservesthedecayasoneoverttimesintegralsoftheenergysoastousegeneralizedgronwall}).

Now, by using the very special fact that it is an $A_{e_a}$ component, and not just any component, it has the special feature that a partial derivative in the direction of $r$, is also a covariant derivative (covariant with respect to the Minkowski metric $m$) in the direction of $r$ of that component, i.e. $$\pa_r  \Lie_{Z^I} A_{e_{a}}  = \derm_r  ( \Lie_{Z^I} A_{e_{a}} )$$ (see \eqref{specialfactaboutAeacomponentsthatexpolitsthefactthatitisnotonlyAtangentialbutspecialone} and \eqref{partialderivativeindirectionofrofAeacomponenentiscovariantderivativeofAeacomponent}). This in turn would allow us, by special integration, as we detailed in \cite{G4}, to translate the pointwise estimate on $|  \derm A_{e_a} |$ into pointwise estimate on $| A_{e_a} |$ (see Lemma \ref{upgradedestimateontheA_acomponentsthatconservesthedecayasoneoverttimesintegralsoftheenergysoastousegeneralizedgronwall}).

This would enable us to control the term $$\frac{ (1+t)}{\eps} \cdot \Big( \sum_{|K| +|J| \leq |I|}   | \Lie_{Z^K}  A_{e_{a}}   | \cdot |\derm ( \Lie_{Z^J}  A_{e_{a}} )  | \Big)^{2}$$ with the right decay factors in time and space (see Lemma \ref{estimateonthebadtermproductA_atimesdermA_ausingbootstrapsoastoclosethegronwallinequalityonenergy}), so that the \textit{weighted} integral in space of such quantity, namely $$\int_{\Sigma^{ext}_{\tau} }  \frac{ (1+t)}{\eps} \cdot \Big( \sum_{|K| +|J| \leq |I|}   | \Lie_{Z^K}  A_{e_{a}}   | \cdot |\derm ( \Lie_{Z^J}  A_{e_{a}} )  | \Big)^{2}   \cdot w(q) \; ,$$ could be controlled (see Lemma \ref{HardyinequalityontheestimateonthebadtermproductA_atimesdermA_ausingbootstrapsoastoclosethegronwallinequalityonenergy}) by a quantity that it is a time integral multiplied with the right factor so that the time integral $$\int_{t_1}^{t_2} \int_{\Sigma^{ext}_{\tau} }  \frac{ (1+t)}{\eps} \cdot \Big( \sum_{|K| +|J| \leq |I|}   | \Lie_{Z^K}  A_{e_{a}}   | \cdot |\derm ( \Lie_{Z^J}  A_{e_{a}} )  | \Big)^{2}   \cdot w(q) \cdot dt $$ would be controlled by a time integral of a time integral, with all terms entering with the right factors so that they could be controlled only by one single time integral, that gives a suitable control to be able to apply later a Grönwall inequality on all components.

\subsection{Dealing with the “bad” term $ A_L   \cdot   \derm A $\;}\label{dealingwithALderA}\

 \textbf{The problems:}\

It might look troublesome to estimate this term, because, to start with, it concerns a product of not two derivatives, but a product of a derivative $|\derm A |$ with a zeroth derivative $|A_L |$ (and not $|\derm A_L |$). Again, the reason we need derivatives, is because the expression of the energy that we want to upgrade involves derivatives (given our natural choice of energy for non-linear hyperbolic wave equations). However, an estimate on $|\derm A_L |$  (and not on $|A_L |$) could be translated into an estimate on $|A_L |$\;, if we integrate the gradient as in the specific integration constructed in Definition \ref{definitionequationforintegralalongthenullcoordinateplusboundaryterm}. However, integration along these lines means significant loss in the decay rate in $|q|$ for $|A_L|$\;, a loss which in the interior region (interior to the outgoing light cone for the Minkowski metric $m$) would destroy entirely the desired factor that we are seeking for the other term $|\derm A|$\;. In fact, in the interior region, we cannot put additional weights in $|q|$ (as opposed to what we can do in the exterior region) because this would change the sign of the derivative of the weight (derivative with respect to $q$) in a way that ruins completely the energy estimate that we established in \cite{G5}: an increasing weight in $|q|$ in the interior region where $q \in ( -\infty, q_0)$\;, would mean that the derivative of the weight $\frac{\pa \widetilde{w} (q)}{\pa q} < 0$ (see Definition \ref{defwidetildew} for the weight) and this enters with the wrong sign in the energy estimate (see Lemma \ref{energyestimatewithoutestimatingthetermsthatinvolveBIGHbutbydecomposingthemcorrectlysothatonecouldgettherightestimatewithtildew}
and Corollary \ref{energyestimatewithoutestimatingthetermsthatinvolveBIGHbutbydecomposingthemcorrectlysothatonecouldgettherightestimate}). This issue regarding the weight is the reason why we consider the problem in the exterior region (and not everywhere) -- we shall elaborate on this more in the introduction.

Another well-known difficulty, is that even if this term $ A_L   \cdot   \derm A $ could be estimated by a product of two derivatives, say by $|\derm A|^2 $\;, there is no hope that such an estimate could work without taking into consideration the very specific structure of the equations, and that’s due to the famous counter-example of John, see \cite{John}.

  \textbf{The solutions:}\
  
It turns out that the Lorenz gauge condition, implies a very good estimate in the exterior on the zeroth Lie derivative of $A_L$ (see Lemma \ref{simpleestimateonlyonthe ALcomponenetwithOUTLiederivative}, or see also Lemma \ref{estimategoodcomponentspotentialandmetric}, where the sum over $|I|-1$ is absent if $|I| = 0$) and a very good estimate on the partial derivative of that component $\pa A_L$ (see Lemma \ref{estimateonpartialderivativeofALcomponent}), which are estimates both needed to control successfully the term $ | A_L  | \cdot    | \derm A   | $\;.

In fact, the leading term of Lie derivatives in the direction of Minkowski vector fields of this product, namely $$ \sum_{|K| = |I| }  \Big(  |  A_L |  \cdot     |  \derm ( \Lie_{Z^K} A )  |  +  |  \Lie_{Z^K} A_L |  \cdot     |  \derm A   |   \Big) \;,$$ contains on one hand, a term  $ \sum_{|K| = |I| }  |  \Lie_{Z^K} A_L | $ that enters with a factor that seems to be the “wrong” factor (see Corollary \ref{TheleeadingtermfortroublesomecomponentsintehsourcesforA}), and on the other hand, a term $ \sum_{|K| = |I| } |  \derm ( \Lie_{Z^K} A )  | $ that enters with the right factor in order to apply a Grönwall inequality on the energy. However, applying a Hardy type inequality in the exterior (see Corollary \ref{HardytypeinequalityforintegralstartingatROm}), we can estimate the term with the “wrong” factor, namely $$ \sum_{|K| = |I| }  |  \Lie_{Z^K} A_L | \; , $$ by its partial derivative instead, namely $| \pa  \Lie_{Z^K} A_L  | $\;, with again a seemingly “wrong” factor, but thanks to Lemma \ref{estimateonpartialderivativeofALcomponent}, this partial derivative can then be estimated by tangential derivatives of the potential $|\rderm A|$\;, with some good error terms, for which we have a very good control on (in Lemma \ref{estimateonthetermthatcontainspartialderivatoveofALinthettimessquareofbadcomponenentsinsourcetermsforAthatcontainstheALcomponent}), despite the “wrong” decaying factor for $|\rderm A|$\;, and this is thanks to our energy estimate (see Lemma  \ref{Theveryfinal }) that controls a space-time integral for such tangential derivatives.

The lower order terms, lower in number of Lie derivatives, namely $$\sum_{|J| + |K| \leq |I| -1}    | \Lie_{Z^J}  A_L |  \cdot     |  \derm ( \Lie_{Z^K} A )  | \;, $$ have a factor for $  |  \derm ( \Lie_{Z^K} A )  |$ which this time involves the Lie derivatives of that special component, namely $|\Lie_{Z^K}  A_L |$\;, and not the “nice” zeroth Lie derivative as opposed to what we explained above. To control the term $  |  \derm ( \Lie_{Z^K} A )  |$ that enters with the “wrong” factor, we use the fact that we have only an $\eps$ loss in the time decay rate (see Lemma \ref{estimateonlowerordertermsforthepotentialAinthesourcetermsforenergy}) thanks to our upgrading of the dispersive estimates for the Lie derivatives of the fields (see Section \ref{UpgradedestimatesfortheLiederivativesoftheEinstein-Yang-Mills fields}, that also presents new challenges that we shall explain later). To control the lower order terms $| \Lie_{Z^J}  A_L |$\;, we plan carefully to use again, in Lemma \ref{HardytpyeestimatesontehtermscontainingALinthebadstructureofsourcestermsforAusefultoobtainenergyestimates}, the Hardy type inequality in the exterior (given in Corollary \ref{HardytypeinequalityforintegralstartingatROm}) to translate them into partial derivatives, but this time, we do not care any more if these are partial derivative of special components. In fact, since all of these are lower order terms, we can close the argument by choosing $\eps$ small enough, small depending on the number of derivatives that we want to control (and depending on other parameters that we shall all show so that one can follow carefully our argument) -- and obviously we must control at least some derivatives so that our applications of the weighted Klainerman-Sobolev inequalities would make sense.

Thereby, we estimate this term $ | \Lie_{Z^I}    \big(  A_L   \cdot     \derm A    \big)    | $ in a way that gives a suitable control in Corollary \ref{thestructureofthebadtermALdermAusingboostrapassumptionanddecompistionofthesumandlowerordertermsexhibitedtodealwithAL}, in such a manner that we would have the correct control in Lemma \ref{estimateonthebadtermproductALtimesdermAusingbootstrapsoastoclosethegronwallinequalityonenergy} on $$\frac{ (1+t)}{\eps}  \cdot \Big(  \sum_{|K| \leq |I|}  | \Lie_{Z^K}    \big(  A_L   \cdot     \derm A    \big)    |   \Big)^{2} $$ with the \textit{right} factors so that the \textit{weighted} space integral $$\int_{\Sigma^{ext}_{\tau} }  \frac{ (1+t)}{\eps}  \cdot \Big(  \sum_{|K| \leq |I|}  | \Lie_{Z^K}    \big(  A_L   \cdot     \derm A    \big)    |   \Big)^{2}  \cdot w(q)  $$ would be suitably controlled in Lemma \ref{HardytpyeestimatesontehtermscontainingALinthebadstructureofsourcestermsforAusefultoobtainenergyestimates} after applying the Hardy type inequality (of Corollary \ref{HardytypeinequalityforintegralstartingatROm}).

In conclusion, both the leading term and the lower order terms have the right factors in Lemma \ref{HardyinequalityontheestimateonthebadtermproductALtimesdermAusingbootstrapsoastoclosethegronwallinequalityonenergy}, that would allow us to inject this successfully in the energy estimate of Lemma \ref{Theveryfinal }.

\subsection{Upgrading the dispersive estimates for the Lie derivatives}\label{dealingwithLiederivativesupgrade}\

 \textbf{The problems:}\

We want to upgrade the dispersive estimates for the Lie derivatives in the direction of Minkowski vector fields of the Einstein-Yang-Mills potential $A$\;, and also of the metric $h$\;, with a tolerable $\eps$ loss in the decay rate in time $t$ (a loss that we could not tolerate for the zeroth Lie derivative), where the $\eps$ loss would ultimately enter in the energy estimate for lower order terms, that are lower in order of number of Lie derivative in comparison to the leading term -- in fact, the leading term that has the highest number of Lie derivatives will only see the factor generated from zeroth Lie derivative which does not have the $\eps$ loss. 

The way in which we are going to upgrade the dispersive estimates for the Lie derivatives, is by using the estimate derived by Lindblad and Rodnianski in Corollary 7.2 of \cite{LR10} (see \eqref{LinfinitynormestimateongradientderivedbyLondbladRodnianski}), but this time, we need to apply it for Lie derivatives of the fields, and for this, we need to estimate correctly the commutator term (see Lemma \ref{theexactcommutatortermwithdependanceoncomponents}). The trouble is that the estimate that was used in the literature to estimate the commutator term for the Einstein vacuum and for Einstein-Maxwell does \textit{not} work.

The difficulty is that when we upgrade for the Lie derivatives using the estimate derived by Lindblad and Rodnianski (see \eqref{LinfinitynormestimateongradientderivedbyLondbladRodnianski}), we actually aim to get a Grönwall inequality, in the form of an integral in time, where the integrand is $\frac{1}{t}$\;, from which we could conclude the time decay rate of $\frac{1}{t^{1-\eps}}$\;, with an $\eps$ loss. However, for the full components of $\derm A$\;, we do not know how to get such a Grönwall type inequality precisely because in the sources for the non-linear wave equation on $A_{{\underline{L}}}$\;, there is the term  $A_{e_a}  \cdot     \derm A_{e_a} $\;, that is not a product for derivatives, and therefore, on one hand, we cannot get the right factor of $\frac{1}{t}$ in the Grönwall inequality because we cannot get the desired decay rate in time for $|A_{e_a}|$\;, and on the other hand, the presence of $ | A_{e_a}| $ does not allow us to establish a Grönwall type inequality on the gradient even if we use the good estimate on $|\derm A_{e_a}|$ (of Lemma \ref{estimatefortangentialcomponentspotential}). We are being therefore \textit{forced} to estimate both terms in the product simultaneously, in a way that does not generate a Grönwall type integral. But, if we want to estimate the term $$  \sum_{|K| + |I| \leq |J| }   | \Lie_{Z^K}  A_{e_a}   |  \cdot | \derm   ( \Lie_{Z^I} A_{e_a} )  | \; ,$$
this means that we need to upgrade the estimate for $| \derm   ( \Lie_{Z^I} A_{e_a} )  |$ first, then aim to translate this, if we can, into an upgraded estimate on $| \Lie_{Z^I} A_{e_a}   |$\;, before upgrading for the component $|\derm A_{{\underline{L}}} |$\;. 

The problem is: the gradient that appears in the Grönwall inequality generated from the previously used estimates for the commutator term injected in the estimate of Lindblad-Rodnianski (of \eqref{LinfinitynormestimateongradientderivedbyLondbladRodnianski}), is a gradient for the \textit{full} components, and therefore would not allow to close an argument, through Grönwall lemma, on $| \derm   ( \Lie_{Z^I} A_{e_a} )|$\;, by appealing only to the structure of the sources of the non-linear terms satisfied by the wave equation on $A_{e_a}$ only.

Another problem: even if one upgrades an estimate on the \textit{covariant} gradient of \textit{one} specific component, this does not mean that one can through integration translate this into an estimate on the zeroth derivative of that \textit{specific} component, since this is a tensor. In fact, when we integrate to transform an estimate on the gradient into an estimate on the component itself, we need to control the partial derivative $\pa_r$ of that component, and this is not necessarily dominated by the covariant gradient of that very specific component, since the covariant derivatives of the vectors of our frame are not necessarily vanishing.

 \textbf{The solutions:}\

As explained, unlike the case of the Einstein vacuum equations, and also unlike the case of the Einstein-Maxwell system, in the case of the Einstein-Yang-Mills equations, we need to have a \textit{more suitable} estimate for the following commutator term $$| \Lie_{Z^I}  ( g^{\la\mu} \derm_{\la}   \derm_{\mu}     A_{e_a} ) - g^{\la\mu}    \derm_{\la}   \derm_{\mu}  (  \Lie_{Z^I} A_{e_a}  ) | \; .$$

We notice that the term in the estimate of commutator that is behind the imposed Grönwall inequality, with an integral involving the gradients of the full components as previously used in the literature, is the term that appears with the weak factor $\frac{1}{|q|}$\;, and not the one that appears with the strong factor $\frac{1}{t}$\;. We realize that this troublesome term could be estimated in a more refined fashion as follows (see \eqref{estimateforgradientoftensorestaimetedbyLiederivativesoftensorwithrefinedcomponentsforthebadfactor}),
 \beaa
 \notag
 |\derm A_{\cal T} | &\les& \sum_{|I| \leq 1}  \frac{1}{(1+t+|q|)} \cdot | \Lie_{Z^I} A |  +  \sum_{ V^\prime \in \cal T } \sum_{|I| \leq 1}  \frac{1}{(1+|q|)} \cdot  | \Lie_{Z^I} A_{ V^\prime} | \; , \\
 \eeaa
 where now, we get an estimate where term with the weak factor is insensitive to the bad component $A_{{\underline{L}}}$ that bothered us (as explained above), which when injected to estimate the commutator term, leads to an estimate that is more refined (see Lemma \ref{commutationformaulamoreprecisetoconservegpodcomponentsstructure}), where the terms with the weak factors do not see the bad component $A_{{\underline{L}}}$ (see Lemma \ref{Thecommutationformulawithpossibilityofsperationoftangentialcomponentsaswell}).
 
 We now inject this estimate on the commutator term in the estimate of Lindblad-Rodnianski that we need in order to upgrade the estimate on the Lie derivatives (see Lemma \ref{estimatethatallowsupgradeincorporatingtermsfromthecommutationformula}). As mentioned, the terms that have the good decay factor in $t$ do not generate a Grönwall type integral (see Lemma \ref{estimateonthegooddecayingpartofthecommutatorterm}). However, the terms that have the weak decaying factor in $|q|$ generate a Grönwall integral (see Lemma \ref{estimateonthebadpartofthecommutatortermtobtainaGronwallforcomponents}), but this time, it does not involve the full components, but only the good ones $A_{\cal T}$\;, for which the non-linear wave equations do not have in their sources the troublesome term $A_{e_a}  \cdot     \derm A_{e_a} $\;. This allows us to establish a Grönwall type inequality that enables us to upgrade first for the good components separately, and consequently for $\derm A_{e_a} $ (see Lemma \ref{GronwallinequalitongradientofAandh1withliederivativesofsourceterms}).
 
However, in order for us to apply the upgrade on $\derm A_{e_a} $ in an upgrade for $\derm A_{{\underline{L}}}$\;, we estimate the product $$  \sum_{|K| + |I| \leq |J| }   | \Lie_{Z^K}  A_{e_a}   |  \cdot | \derm   ( \Lie_{Z^I} A_{e_a} )  | \; ,$$ which involves estimating $ | \Lie_{Z^K}  A_{e_a}   | $ (which is not the covariant gradient). So how to translate an estimate on the covariant gradient of specific components $|\derm ( \Lie_{Z^I}  A_{\cal T} )|$\;, into an estimate on the partial derivative $| \pa_r ( \Lie_{Z^I}  A_{e_a} ) |$\;, which is the term that is actually needed to be estimated in order to integrate to estimate $\Lie_{Z^I}  A_{e_a} $? We use the specific fact that it is not just any component, but an  $A_{e_a}$ component, and we notice that $\derm_{r} e_a = 0 $\;. This allows us to get the right estimate on this product (see Lemma \ref{upgradeddecayestimateonthetermAe_agradAe_a}), which we shall use to upgrade for the full components (see Lemma \ref{SteptwoforinductonforA}).

We point out that here, the term $ \sum_{|K| + |I| \leq |J| }  | \Lie_{Z^K} A_L  | \cdot    | \derm  ( \Lie_{Z^I}  A  ) | $ is not a problem this time, since we can estimate the product successfully using the a priori estimates (see Lemma \ref{upgradeddecayestimateonthetermALgradientA}) and this is thanks to the Lorenz gauge estimate on $A_L$ (see Lemma \ref{estimategoodcomponentspotentialandmetric}).

Finally, we will prove the theorem.

\subsection{The theorem}\label{Thetheoremofexteriorstabilityfornequalthree}\

\begin{theorem}
Assume that we are given an initial data set $(\Sigma, \overline{A}, \overline{E}, \overline{g}, \overline{k})$ for \eqref{EYMsystemforintro}. We assume that $\Sigma$ is diffeomorphic to $\R^3$\;. Then, there exists a global system of coordinates $(x^1, x^2, x^3) \in \R^3$ for $\Sigma$\;. We define
\bea
r := \sqrt{ (x^1)^2 + (x^2)^2 +(x^3)^2  }\;.
\eea
Furthermore, we assume that the data $(\overline{A}, \overline{E}, \overline{g}, \overline{k}) $ is smooth and asymptotically flat. 

Let $\chi$ be a smooth function such as
 \bea\label{defXicutofffunction}
\chi (r)  := \begin{cases} 1  \quad\text{for }\quad r \geq \frac{3}{4} \;  ,\\
0 \quad\text{for }\quad r \leq \frac{1}{2} \;. \end{cases} 
\eea
Let $M$ be the mass such as $ 0 < M \leq \eps^2 \leq 1 $\;, let $\de_{ij}$ be the Kronecker symbol, and let $\overline{h}^1_{ij} $ be defined in this system of coordinates $x^i$\;, by
 \bea
\overline{h}^1_{ij} := \overline{g}_{ij} - (1 + \chi (r)\cdot  \frac{ M}{r}  ) \de_{ij} \; .
\eea

We then define the weighted $L^2$ norm on $\Sigma$\;, namely $\overline{\E}_N$\;, for $\ga > 0$\;, by
 \bea\label{definitionoftheenergynormforinitialdata}
 \notag
\overline{\E}_N &:=&  \sum_{|I|\leq N} \big(   \| (1+r)^{1/2 + \ga + |I|}   \overline{D} (  \overline{D}^I  \overline{A}    )  \|_{L^2 (\Sigma)} +  \|(1+r)^{1/2 + \ga + |I|}    \overline{D}  ( \overline{D}^I \overline{h}^1   )  \|_{L^2 (\Sigma)} \big) \\
\notag
&:=&  \sum_{|I|\leq N}      \big(   \sum_{i=1}^{n}  \|(1+r)^{1/2 + \ga + |I|}     \overline{D} (  \overline{D}^I  \overline{A_i}    )  \|_{L^2 (\Sigma)} +  \sum_{i, j =1}^{n}  \|(1+r)^{1/2 + \ga + |I|}     \overline{D}  ( \overline{D}^I \overline{h}^1_{ij}   )  \|_{L^2 (\Sigma)} \big) \; , \\
\eea
where the integration is taken on $\Sigma$ with respect to the Lebesgue measure $dx_1 \ldots dx_n$\;, and where $\overline{D} $ is the Levi-Civita covariant derivative associated to the given Riemannian metric $\overline{g}$\;.

We also assume that the initial data set $(\Sigma, \overline{A}, \overline{E}, \overline{g}, \overline{k})$ satisfies the Einstein-Yang-Mills constraint equations, namely
\bea
\notag
  \mathcal{R}+ \overline{k}^i_{\,\, \, i} \overline{k}_{j}^{\,\,\,j}  -  \overline{k}^{ij} \overline{k}_{ij}   &=&    \frac{4}{(n-1)}   < \overline{E}_{i}, \overline{E}^{ i}>   \\
 \notag
 && +  < \overline{D}_{i}  \overline{A}_{j} - \overline{D}_{j} \overline{A}_{i} + [ \overline{A}_{i},  \overline{A}_{j}] ,\overline{D}^{i}  \overline{A}^{j} - \overline{D}^{j} \overline{A}^{i} + [ \overline{A}^{i},  \overline{A}^{j}] >  \;  ,\\
 \notag
\overline{D}_{i} \overline{k}^i_{\,\,\, j}    - \overline{D}_{j} \overline{k}^i_{\, \,\,i}  &=&  2 < \overline{E}_{i}, \overline{D}_{j}  \overline{A}^{i} - \overline{D}^{i} \overline{A}_{j} + [ \overline{A}_{j},  \overline{A}^{i}]  >  \;  ,\\
\notag
\overline{D}^i \overline{E}_{ i} + [\overline{A}^i, \overline{E}_{ i} ]  &=& 0  \;  .
\eea

For any $N \geq 11$\;, there exists a constant $ \overline{c} (\cal{K}, N, \ga) $\;, that depends on $\cal{K}$\;, on $\ga$\;, and on $N$\;, such that if

\bea\label{Assumptiononinitialdataforglobalexistenceanddecay}
\overline{\E}_{N+2} &\leq&   \overline{c} (\cal{K}, N, \ga)  \; ,\\
M  &\leq&   \overline{c} (\cal{K}, N, \ga)  \; ,
\eea
then there exists a solution $(\cal{M}, A, g)$ to the Cauchy problem for the fully coupled Einstein-Yang-Mills system \ref{EYMsystemforintro} in the future of the whole causal complement of any compact $\cal{K} \subset \Sigma$ \;, converging to the null Yang-Mills potential and to the Minkowski space-time in the following sense: if we define the metric $ m_{\mu\nu}$ to be the Minkowski metric in wave coordinates $(x^0, x^1, x^3)$\, and define $t = x^0$\;, and if we define in this system of wave coordinates 
\bea
h^{1}_{\mu\nu} := g_{\mu\nu} - m_{\mu\nu}- h^0_{\mu\nu}  \;  ,
\eea
where for $t > 0$, 
\bea\label{guessonpropagationofthesphericallsymmetricpart}
h^0_{\mu\nu} := \chi(r/t) \cdot \chi(r)\cdot \frac{M}{r}\de_{\mu\nu}  \;  ,
\eea
and where for $t=0$\;,  \bea\label{definitionofthesphericallysymmtericpartofinitialdata}
h^0_{\mu\nu} ( t= 0) := \chi(r)\cdot \frac{M}{r}\cdot \de_{\mu\nu}  \; ,
\eea
then, for $\overline{h}^1_{ij} $ and $\overline{A}_i$ decaying sufficiently fast, we have the following estimates on $h^1$\;, and on $A$ in the Lorenz gauge, for the norm constructed using wave coordinates by taking the sum over all indices in wave coordinates, given in \eqref{boundonNLiederivativeofgradient}, \eqref{boundonNLiederivativeofzerothderivative}, \eqref{boundonNLiederivativeotheYangMillscurvature}. That there exists a constant $C(N)$ to bound  $N-2$ Lie derivative of the fields in direction of Minkowski vector fields, and to bound the growth of $\E _{N} (\cal K) (t)$ in \eqref{theboundinthetheoremonEnbyconstantEN}, and that there exists a constant $\eps$ that depends on $\overline{c} (\cal{K}, N, \ga)$\;, on  $\cal{K} \subset \Sigma$\;, on $N$ and on $\ga$\;, such that we have the following estimates in the whole complement of the future causal of the compact $\cal{K}  \subset \Sigma$\;, for all $|I| \leq N -2$\;, 

\bea\label{boundonNLiederivativeofgradient}
 \notag
  && \sum_{\mu= 0}^{n} |\derm  ( \Lie_{Z^I}  A_{\mu} ) (t,x)  |     +   \sum_{\mu, \nu = 0}^{n}  |\derm  ( \Lie_{Z^I}  h^1_{\mu\nu} ) (t,x)  |   \\
   \notag
    &\leq&  C(\cal{K} ) \cdot   C(N)  \cdot \frac{\eps }{(1+t+|r-t|)^{1-\eps} (1+|r-t|)^{1+\gamma}} \; ,\\
      \eea
and
 \bea\label{boundonNLiederivativeofzerothderivative}
 \notag
  \sum_{\mu= 0}^{n}  |\Lie_{Z^I} A_{\mu} (t,x)  | +   \sum_{\mu, \nu = 0}^{n}  |\Lie_{Z^I}  h^1_{\mu\nu} (t,x)  |  &\leq&C(\cal{K} ) \cdot C(N) \cdot c(\ga) \cdot \frac{\eps }{(1+t+|r-t|)^{1-\eps} (1+|r-t|)^{\gamma}} \; , \\
      \eea
      where $Z^I$ are the Minkowski vector fields (see Subsection \ref{TheMinkowksivectorfieldsdefinition}).

      In particular, the gauge invariant norm on the Yang-Mills curvature decays as follows, for all $|I| \leq N - 2$\;,
 \bea\label{boundonNLiederivativeotheYangMillscurvature}
 \notag
 \sum_{\mu, \nu = 0}^{n}  |\Lie_{Z^I} F_{\mu\nu}  (t,x) |  &\leq&C(\cal{K} ) \cdot   C(N)    \cdot \frac{\eps }{(1+t+|r-t|)^{1-\eps} (1+|r-t|)^{1+\gamma}}  \\
  \notag
&& +     C(\cal{K} ) \cdot C(N) \cdot c(\ga)   \cdot \frac{\eps }{(1+t+|r-t|)^{2-2\eps} (1+|r-t|)^{2\gamma}} \; . \\
      \eea
      
      Furthermore, if one defines $w$ as follows, 
\bea
w(r-t):=\begin{cases} (1+|r-t|)^{1+2\gamma} \quad\text{when }\quad r-t>0 \;, \\
         1 \,\quad\text{when }\quad r-t<0 \; , \end{cases}
\eea
and if we define  $\Sigma_t^{ext} (\cal K)$ as being the time evolution in wave coordinates of $\Sigma$ in the future of the causal complement of $\cal{K}$\;, then for all time $t$\;, we have
\bea\label{theboundinthetheoremonEnbyconstantEN}
\notag
\E _{N} (\cal K) (t) &:=&  \sum_{|J|\leq N} \big( \|w^{1/2}   \derm ( \Lie_{Z^J} h^1   (t,\cdot) )  \|_{L^2 (\Sigma_t^{ext} (\cal K)) } +  \|w^{1/2}   \derm ( \Lie_{Z^J}  A   (t,\cdot) )  \|_{L^2 (\Sigma_t^{ext} (\cal K))  }  \big) \\
\notag
&\leq& C(N) \cdot \eps \cdot (1+t)^{\eps} \; . \\
\eea

More precisely, for any constant $E(N)$ and for any $\ga > 0$\;, there exists a constant $\de (\ga) > 0$ that depends on $\ga$\;, and there exists a constant $\overline{c}_1 (\cal{K}, N, \ga, \de)$\;, that depends on  $\cal{K} \subset \Sigma$\;, on $E(N)$\;, on $N$\;, on $\ga$ and on $\de$\;, such that if

\bea\label{Assumptiononinitialdataforglobalexistenceanddecay}
 \overline{\E}_{N+2} (0) &\leq&   \overline{c}_1 (\cal{K}, E (N), N, \ga, \de)   \; ,\\
 M &\leq&  \overline{c}_1 (\cal{K}, E (N), N, \ga, \de)  \; ,
 \eea
then, we have in the whole complement of the future causal of the compact $\cal{K}  \subset \Sigma$\;, for all time $t$\,, 
\bea\label{Theboundontheglobaleenergybyaconstantthatwechoose}
 \E_{N}(\cal K)  (t) \leq E(N) \cdot (1+t)^{\de} \; .
\eea

\end{theorem}

\section{Introduction}\label{Introduction}

We prove the exterior stability of the $(1+3)$-Minkowski space-time governed by the evolution problem in General Relativity with matter, which in this case is the Einstein-Yang-Mills system. In fact, in $(1+3)$ dimensions, the Einstein-Yang-Mills equations read (see \cite{G4}) the following system on the unknown $(\cal M, A, \g)$\;,
\bea
\begin{cases} \label{EYMsystemforintro}
R_{ \mu \nu}  =& 2 < F_{\mu\b}, F_{\nu}^{\;\;\b} >- \frac{1}{2} \cdot g_{\mu\nu } < F_{\a\b },F^{\a\b } >\; ,   \\
0 \;\;\;\, \,=&  \der_{\alpha} F^{\a\b}  + [A_{\alpha}, F^{\a\b} ]    \;, \\
F_{\a\b} =& \der_{\a}A_{\b} - \der_{\b}A_{\a} + [A_{\a},A_{\b}] \; ,  \end{cases} 
      \eea
where $\cal M$ is the unknown manifold, where $A$ is the unknown Yang-Mills potential valued in the Lie algebra $\cal G$ associated to any compact Lie group $G$\;, where $\g$ is the unknown Lorentzian metric, and where $\der_{\alpha}$ is the unknown space-time covariant derivative of Levi-Civita prescribed by $\g$ and $R_{ \mu \nu} $ is the Ricci tensor. We start with a general initial data for the Yang-Mills potential $A$ chosen to be small, and with general initial data for the metric $\g$ that is asymptotical flat and chosen to be close to the Minkowski initial data. Then, we show that in the Lorenz gauge and in wave coordinates, the perturbations disperse in time in the complement of the future causal domain of a compact set, and lead to a solution that is converging to the Minkowski space-time in the exterior.

To achieve this, we recast the Einstein-Yang-Mills system, \ref{EYMsystemforintro}, in the Lorenz gauge and in wave coordinates, as a coupled system of non-linear covariant wave equations on both the unknown Yang-Mills potential and the unknown metric (see \cite{G4}). The goal is to use a bootstrap argument on a higher order energy norm, so as to make a statement on both the existence and the dispersion of the solutions (see  Subsection \ref{Thebootstrapargumentandnotationonboundingtheenergy}). 

However, in the case of $n=3$\,, where $n$ is the space dimension, there is a serious problem due to the lack of integrability of the decay generated from the Klainerman-Sobolev inequality combined with the bootstrap assumption on the exterior energy. First, let us point out that we need to work with an energy defined only as an integral on the exterior of the outgoing light cone, due to the lack of integrability in the interior region to the light cone of some of the factors in the source terms which already fail in the case $n=4$\,. In fact, as far as this remark is concerned, one could conjecture that one does not have stability at all in the interior, where there could be concentration of energy, as pointed out by Mondal and Yau, in \cite{MY1}, for the Einstein-Yang-Mills equations. Yet, our approach is totally different and independent of theirs, and relies on using a null frame decomposition used by Lindblad and Rodnianski, in \cite{LR2} and \cite{LR10}, for the Einstein vacuum equations in wave coordinates in their studying of the stability problem of the Minkowski space-time in vacuum, which we carry out, here, to the coupled Einstein-Yang-Mills system in the Lorenz gauge. In other words, we study the non-linear structure using a null frame decomposition (see Lemma \ref{TheEinstein-Yang-Millssysteminanullframe}), as a starting point to address the difficulties that arise in three space-dimensions.

First, let us start by pointing out, most importantly, the fact that if one does not allow a polynomial growth of the type $(1+t)^\de$ of the higher order energy, one would then need to have integrable factors in the Grönwall inequality, in order to close the continuity argument. However, the Klainerman-Sobolev inequality in $(1+3)$-dimensions, would then give a decay rate in time as $\frac{1}{t}$ which when integrated, would give a logarithmic growth and therefore would lead to polynomial growth of the energy with a rate of the type $(1+t)^\de$\;. Thus, proving boundedness of the energy is not possible to achieve using the Klainerman-Sobolev inequality alone.

Instead, one would allow a polynomial growth of the energy norm, at the expense of losing even more from the decay on the fields after applying the Klainerman-Sobolev inequality. A growth at the rate $(1+t)^\de$ of the energy (see \eqref{bootstrap}), implies then through the Klainerman-Sobolev inequality a decay rate in time on the gradient of the fields at the rate of  $\frac{1}{t^{1-\de}}$ (see Lemma \ref{apriordecayestimatesfrombootstrapassumption}), which is not at all integrable, and therefore a Grönwall inequality on the energy using such a decay rate would even worse, give exponential growth!, which is even far worse than what we have started with, that is a polynomial growth. This already gives one an insight about why the problem is difficult.

Yet, one can wish to improve the decay rate on the gradient of the fields, by any mean whatsoever that works, in order to improve the decay rate from  $\frac{1}{t^{1-\de}}$ to  $\frac{1}{t}$\;, and therefore, a Grönwall inequality on the energy would then allow one to conclude that the polynomial growth that we started with will be conserved, and if one starts with an initial data small enough, then one could also improve it by a multiplication factor and hence, close the continuity argument.

However, it turns out that exploiting the very special structure of the Einstein-Yang-Mills equations (see Lemma \ref{TheEinstein-Yang-Millssysteminanullframe}) along special estimates generated from the Lorenz gauge (see Lemma \ref{LorenzgaugeestimateforgradientofLiederivativesofAL}) and from the wave coordinates condition (see Lemma \ref{wavecoordinatesestimateonLiederivativesZonmetric}), would allow us to improve the decay rate of the gradient, for the zeroth Lie derivative, of some special components, but not for all components. Although, this does not extend to the Lie derivatives of even the special components, it turns out that this is sufficient to derive a suitable Grönwall inequality, exploiting of course the structure of the Einstein-Yang-Mills equations. Let us explain:

In wave coordinates, one can define a null frame, that is a null frame for the Minkowski metric $m$\;, since $m$ is constructed using these same wave coordinates (see Definition \ref{definitionofthenullframusingwavecoordinates}). In this null frame, there are three vectors tangential to the outgoing light cone $ \{L, e_1, e_2\}$\;, where $L$ is null, tangent and orthogonal to the outgoing light cone (keeping in mind that it is a Lorentzian metric). A fourth vector $\underline{L}$ is conjugate to $L$\;, null and is in fact tangential to an ingoing null cone. Now, an important observation is that derivatives in the direction of the outgoing null cone are expected to decay faster than those in the direction of the ingoing null cone. In fact, that is because for wave equations, one expects the fields to disperse more in time when one goes away from their support. We are indeed working with initial data that is in a certain weighted Sobolev space, not necessarily with compact support though, however this imposes already that the gradient of the fields decay sufficiently enough at spatial infinity and therefore, there is a direction, that of the outgoing light cone, that goes away from the possible concentration of energy of the initial data, and another direction, that of the ingoing light cone, that goes into the region where the energy could be concentrated. Hence, the vectors $ \{L, e_1, e_2\}$ are already in a good outgoing direction and therefore derivatives of fields in that direction are expected to decay faster that in the direction of the ingoing light cone $\underline{L}$.

One can see this from the following: if one defines $\rpa$ to be the projection of the derivatives at some point in space-tine $(t, x) \in \R^{1+3} $ on the outgoing light cone passing through that point, projecting using an euclidian metric, that is the usual metric on $\R^{1+3}$ (see Definition \ref{defrestrictedderivativesintermsofZ}), then we have the following inequality (see \eqref{goodderivative} and Lemmas \ref{betterdecayfortangentialderivatives} and \ref{decayrateforfullderivativeintermofZ}),
\bea\label{boundinggradientandtangentialderivativesbyMinkwoskiLiederivatives}
\notag
(1 + t + | (|x| -t) | ) \cdot |\rpa \phi | (t, x) + (1 +  |(|x| -t)  | ) \cdot  |\pa \phi | (t, x) &\les& \sum_{|I| = 1} | \Lie_{Z^I} \phi | (t, x) \; . \\
\eea
If $\phi$ satisfies, say a free wave equation on the Minkowski background, namely $m^{\a\b} \derm_\a \derm_\b \phi = 0 $\;, then since the vectors $Z^I$ are the Killing vector fields of Minkowski, then applying $\Lie_{Z^I}$ to the wave operator in hand, either commutes or generates a term proportional to the wave operator in hand, we therefore get also $ m^{\a\b} \derm_\a \derm_\b (  \Lie_{Z^I} \phi ) = 0 $\;. However, by energy estimates, if one bounds $\phi$\;, then the same would be true for  $\Lie_{Z^I} \phi$\;, since the equation is preserved for the Lie derivatives. Consequently, the right hand side of the inequality \eqref{boundinggradientandtangentialderivativesbyMinkwoskiLiederivatives} would then be bounded and as a result, we see already due to the above mentioned inequality that we get a decay rate for both the gradient of the fields and the tangential derivatives of the fields.

However, we immediately notice from the above inequality \eqref{boundinggradientandtangentialderivativesbyMinkwoskiLiederivatives}, that derivatives tangential to the outgoing light cone decay better in time $t$ than the whole derivatives. This was provided that  $\Lie_{Z^I} \phi $ are bounded, for which we used the fact that Minkowski vector fields commute, or almost commute, with the wave operator governing $\phi$\;. Hence, in a bootstrap argument, if one assumes a certain control on either the boundedness or the growth of the energy, one would have then also assumed a better decay rate for tangential derivatives. Yet, one still needs to close the argument by improving such a bound on the energy. If the non-linear structure of the wave equation in hand, involves “good” derivatives, i.e. the tangential derivatives, which already have a better decay rate than the one that we already started with by using the Klainerman-Sobolev inequality combined with the assumption on the growth of the energy, then one can try to exploit this better decay rate in the non-linear structure in order to improve the bound on the growth of the energy, and thereby to conclude that the a priori estimate that we started with was indeed a true estimate on the energy.

One can in fact do so, i.e. run successfully such a continuity argument for non-linear wave equations that have a null structure: see, for instance, the work of Christodoulou in \cite{Chr1} and Klainerman in \cite{Kl2}. For null structures, one can bound the source terms by a “good” tangential derivative multiplied by a “bad” whole derivative. Since tangential derivatives decay fast enough, this would already offer the sufficient enough decay rate of the type  $\frac{1}{(1+t )^{1+\la}}$\;, $\la > 0$ (that is good) multiplied by the whole derivative, which then one could transform into $\frac{1}{(1+t )^{1+\la}}$ multiplied by the energy norm. However, we are left with the bad decay rate  $\frac{1}{(1+t )^{1-\de}}$ multiplied by the tangential derivative; this does not have a sufficient decay factor in order to be transformed into $\frac{1}{(1+t )^{1+\la}}$ multiplied by the energy norm. However, for non-linear wave equations, one can control a space-time integral of the square of the tangential derivatives, or differently speaking, one can show that the $L^2$ norm of tangential derivatives is integrable in time (see for instance, \cite{LR10} and \cite{G5}, or see Lemma \ref{Theveryfinal }). Thus, one does not have to worry about the fact that tangential derivatives do not have a good decaying factor in front of them in the non-linear structure of the wave equation. 

One could then hope to generalize such an argument to non-linear wave operators, as the one in study here, $ g^{\a\b} \derm_\a \derm_\b \phi $\;, as long as the commutation of the Minkowski Lie derivatives with such an operator has a good structure and as long as the Minkowski Lie derivatives of the source terms has a good structure as well (for example, see Lemma \ref{theexactcommutatortermwithdependanceoncomponents}).

This being said, one must keep in mind that if one had a square of a full derivative in the non-linear structure, then one could expect a blow-up in finite time. This was, for instance, proven by John, see \cite{John}, for $\Box \phi = (\pa_t \phi)^2$\;. Thus, a square of derivative could be disastrous for global existence of solutions. One needs to study more. Yet, the Einstein-Yang-Mills equations do not satisfy the null condition, which is a sufficient criteria, but not necessary, to prove global existence of solutions. Hence, one has to exploit more the structure. In the case of the Einstein vacuum equations, without matter, this was dealt with by Lindblad and Rodnianski in \cite{LR10} using wave coordinates and a null frame decomposition.

The hope there, however, was to generalize this also to components of a tensor, and not only derivatives of a tensor. In other words, one hopes that there are components of the field strength solution of the wave operator $ g^{\a\b} \derm_\a \derm_\b \phi $\;, which decay better than other components. It turns out again that considering on one hand, components tangential to the outgoing light cone $\phi_{\cal T}$\;, where $\cal T := \{L, e_1, e_2\}$\;, and on the other hand, the component which is non-tangential to the outgoing light cone, namely $\phi_{\underline{L}}$\;, following the philosophy of \cite{LR10}, could do the job also in the case of the Einstein-Yang-Mills equations, as shown here, provided, however, that the source terms for these “good” components $\phi_{\cal T}$ and “bad” component $\phi_{\underline{L}}$\;, have a good structure that we exploit here.

First, let us point out that in $(1+3)$-space-time dimensions, the initial data for the metric has a spherically symmetric part, that carries the mass $M$ (see \ref{definitionofthesphericallysymmtericpartofinitialdata}), that does not decay fast enough to have finite weighted energy, such as for the energy we consider (see Remark \ref{thesphericallysymmmetricparthasinfiniteenergy}). Hence, the right quantity to look at, for the wave equation, is the metric from which one subtracts this spherically symmetric part $h^0$ of the initial data that behaves as $\frac{M}{r}$\;, keeping in mind that the rest $h^1 := h - h^0$ could be imposed on to decay faster. This would then mean that one should know how $h^0$ propagates. This was dealt with by Lindblad and Rodnianski by “a guess” on how $g^{\la\mu} \derm_{\la}   \derm_{\mu}  h^0$ should look like, and one can then prove that the remaining quantity $g^{\la\mu} \derm_{\la}   \derm_{\mu}  h^1$ is well-posed, meaning any remaining propagation in time of $h^0$ from the guess that they made (see \eqref{guessonpropagationofthesphericallsymmetricpart}), could be injected in $h^1$\;, while being consistent with the initial data for $h^1$\;, and could be proven as a global solution defined for all time, in the Einstein vacuum case (see \cite{LR10}).

In the case of the Einstein vacuum equations studied in \cite{LR10}, as shown by Lindblad and Rodnianski, one could obtain the following estimates for the “good” components, where $\cal U := \{  \underline{L}, L, e_1, e_2 \}$\;, $\cal T := \{ L, e_1, e_2 \}$\;,
 \bea
\notag
 && | g^{\alpha\beta} \derm_\alpha \derm_\beta h^1_{ {\cal T} {\cal U}} |   \\
 \notag
    &\les&  |\rderm h| \cdot |\derm h|  +  |h| \cdot |\derm h|^2 +   | g^{\alpha\beta} \derm_\alpha \derm_\beta h^0 | \; ,
\eea
and for the “bad” component,
 \bea
\notag
 && | g^{\alpha\beta} \derm_\alpha \derm_\beta h^1_{{\underline{L}} {\underline{L}} } |   \\
 \notag
    &\les&  |\rderm h| \cdot |\derm h|  +  |h| \cdot |\derm h|^2 +   | g^{\alpha\beta} \derm_\alpha \derm_\beta h^0 | \\
    \notag
&&    +  |\derm h_{\cal TU} |^2  \; .
\eea
The term $|\rderm h| \cdot |\derm h|$ comes out from the null structure and could be dealt with easily by already known heuristics. The term  $|h| \cdot |\derm h|^2 $ decays fast enough not to cause any trouble whatsoever. However, the term  $|\derm h_{\cal TU} |^2$ might look, at first instance, as the counter example for global existence by John in \cite{John}, however, the idea of Lindblad and Rodnianski is that in fact, the terms $h^1_{ {\cal T} {\cal U}}$ satisfy wave equations with a null structure -- it is exactly the reason why such terms are called “good” components --, and therefore, they could be shown to decay better than the decay rate that one obtains from the Klainerman-Sobolev inequality and the bootstrap assumption. This is achieved by obtaining an estimate on the derivatives of solutions to non-linear wave equations, provided that one controls the non-linear structure good enough. This is again, as mentioned earlier, a generalization of having improved decay rate for tangential derivatives (or for “good” derivatives), provided that one bounds the Lie derivatives of the field (say by using energy estimates for wave equations). Yet, here one wishes to have improved decay rate for \textit{all} derivatives, considering the fact that it is derivative of a solution to a non-linear wave equation, provided that one controls good enough the non-linearity structure -- which could be done for “good” components.

The estimate derived by Lindblad and Rodnianski (see Corollary 7.2 of \cite{LR10}), that is at the heart of their argument for the case of the Einstein vacuum equations, is essentially the following (see \eqref{LinfinitynormestimateongradientderivedbyLondbladRodnianski}),
\beaa
 (1 + t + | (|x| -t) | ) \cdot |\pa \phi |  &\les& \sum_{|I| \leq 1} | \Lie_{Z^I} \phi | \\
 &&+ \quad \text{a time integral of terms involving $g^{\alpha\beta} \derm_\alpha \derm_\beta\phi$}\\
  &&+  \quad  \text{a time integral of terms involving the field $\phi$ or $\pa\phi$}\; .
\eeaa
These time integrals could be bounded using the standard Klainerman-Sobolev inequality with the bootstrap assumption on the energy, or the weighted energy. Since in the case of Einstein vacuum equations, the good components have an excellent structure, namely a null structure, then that time integral could be bounded and one gets a decay rate of the type $\frac{1}{t}$ for the good components. Since these good components enter as product, namely as a square $|\derm h_{\cal TU} |^2$\;, in the non-linear structure for the “bad” components, keeping in mind the discussion above, we would have a factor of $\frac{1}{t}$ multiplied by the gradient, which is good enough to obtain a Grönwall inequality on the energy with the right factor to close the bootstrap argument on the energy, assuming a polynomial growth. This is essentially the argument of Lindblad and Rodnianski in \cite{LR10}.

Of course, such a discussion is simplified. On one hand, one needs to study the structure of the commutation between the wave operator and the Minkowski Lie derivatives. As shown by Lindblad and Rodnianski in \cite{LR10}, using the wave coordinates condition, one could get good estimates for such commutation that is useful for the case of the Einstein vacuum equations. On the other hand, the Lie derivatives of the “good” components of the metric do not decay in an integrable fashion: this only holds true for the zeroth-Lie derivative. However, in the Grönwall inequality on the energy, these terms with “bad” decaying factors, only appear as lower order terms (lower in the number of Lie derivatives) and thus, one could still close the argument.

In the case of the Einstein-Yang-Mills system in the Lorenz gauge and in wave coordinates, \eqref{EYMsystemforintro}, the equations could be estimated in the null frame that was constructed using wave coordinates (see Definition \ref{definitionofthenullframusingwavecoordinates}), as in the following for the “good” components (see Lemma \ref{estimateonthesourcetermsforgoodcomponentofPoentialAandgoodcompometrich}),
  \bea\label{estimateonnonlinearstructureofgoodcomponentforAinintrrod}
   \notag
 &&  | g^{\la\mu} \derm_{\la}   \derm_{\mu}   A_{{\cal T}}   | \\
 \notag
    &\les& | \derm h | \cdot  |\rderm A |       + | \rderm  h | \cdot  |\derm A |  +   | A  | \cdot    | \rderm A  | +  | \derm  h | \cdot  | A  |^2   +  | A  |^3 \\
    \notag
  && + O( h \cdot  \derm h \cdot  \derm A) + O( h \cdot  A \cdot \derm A)  + O( h \cdot  \derm h \cdot  A^2) + O( h \cdot  A^3) \; , \\
  \eea
  and
 \bea\label{estimateonnonlinearstructureofgoodcomponentforhoneinintrrod}
\notag
 && | g^{\alpha\beta} \derm_\alpha \derm_\beta h^1_{ {\cal T} {\cal U}} |   \\
 \notag
    &\les&   |\rderm h| \cdot |\derm h|+|h| \cdot |\derm h|^2   +  | \rderm A   | \cdot  |\derm  A  |  +    | A |^2 \cdot | \derm A   |  + | A |^4  \\
\notag
     && + O \big(h \cdot  (\derm A)^2 \big)   + O \big(  h  \cdot  A^2 \cdot \derm A \big)     + O \big(  h   \cdot  A^4 \big)  \\
  && +   | g^{\alpha\beta} \derm_\alpha \derm_\beta h^0 | \; , 
\eea

and as in the following for the “bad” component (see Lemma \ref{structureofthesourcetermsofthewaveoperatoronAandh}),
\bea\label{estimateonnonlinearstructureofbadcomponentforAinintrrod}
   \notag
 && | g^{\la\mu} \derm_{\la}   \derm_{\mu}   A_{{\underline{L}}}  |   \\
 \notag
    &\les& | \derm h | \cdot  |\rderm A |       + | \rderm  h | \cdot  |\derm A |  +   | A  | \cdot    | \rderm A  | +  | \derm  h | \cdot  | A  |^2   +  | A  |^3 \\
    \notag
  && + O( h \cdot  \derm h \cdot  \derm A) + O( h \cdot  A \cdot \derm A)  + O( h \cdot  \derm h \cdot  A^2) + O( h \cdot  A^3) \\
         \notag
 && + | A_L  | \cdot    | \derm A   |   + | A_{e_a}  | \cdot    | \derm A_{e_a}  |     \; , \\
  \eea
 and
    \bea\label{estimateonnonlinearstructureofbadcomponentforhoneinintrrod}
\notag
&&  | g^{\alpha\beta} \derm_\alpha \derm_\beta h^1_{{\underline{L}} {\underline{L}} } |    \\
\notag
            &\les& |\rderm h| \cdot |\derm h| + |h| \cdot |\derm h|^2   +  | \rderm A   | \cdot  |\derm  A  | +  | A |^2  \cdot | \derm A   |   +  |A |^4 \\
    \notag
     && + O \big(h \cdot  (\derm A)^2 \big)   + O \big(  h  \cdot  A^2 \cdot \derm A \big)     + O \big(  h   \cdot  A^4 \big)   \\
      \notag
 &&+ |\derm h_{\cal TU} |^2  +   |\derm A_{e_{a}}  |^2 \\
 && +   | g^{\alpha\beta} \derm_\alpha \derm_\beta h^0 | \; . 
\eea

If one looks at the non-linear structure of this Einstein-Yang-Mills system in the Lorenz gauge and in this null frame constructed using wave coordinates, obviously, the first terms which appear, of the type $| \derm h | \cdot  |\rderm A |$ and $ | \rderm  h | \cdot  |\derm A |$ come out from a null structure, which everyone knows how to deal with by already established heuristics, as discussed.

Yet, even for the so called “good” components $A_{\cal T} $\;, we see already terms that fail to decay fast enough in the interior region (that is interior to the outgoing light cone whose tip is the origin), which are $| A  | \cdot    | \rderm A  |$ and  $| A |^3$\;, which failed already to decay sufficiently enough in the interior even for $n= 4$\;. If we look at the product $| A  | \cdot    | \rderm A  |$\;, one has a term which consists of a factor multiplied by $| \rderm A  |$, and another term that consists of a factor multiplied by $|A|$\;. Hence, one of the terms is $|A|$ multiplied by a factor that does not decay good enough in the interior: it is a term of the type $ |A| \cdot \begin{cases}  \frac{1 }{(1+t+|q|)^{2-\delta} (1+|q|)^{\ga}} \;,\quad\text{when }\quad q>0 \;,\\
      \frac{ (1+|q|)^{\frac{1}{2} }}{ (1+t+|q|)^{2-\delta} } \;,\quad\text{when }\quad q<0 \; , \end{cases} \\ $ where here $q:= r-t$. This means that in the interior, we have a term of the type $|A| \cdot\frac{1}{(1+t+|q|)^{\frac{3}{2}-\delta}}$\;, whereas we need a factor of the type $\frac{1}{(1+t)\cdot (1+|q| )}$ to have sufficient control. (Yet, one then can see that in $n=4$\;, this interior factor would become $|A| \cdot\frac{1}{(1+t+|q|)^{\frac{5}{2}-\delta}}$ which is good enough, provided that we use the fact that the other term is in fact a tangential derivative.)

However, the term $|A|^3$ is far much more problematic in the interior region: it could be estimated as $|A| \cdot \begin{cases}   \frac{1}{ (1+ t + | q | )^{2-2\delta }  (1+| q |   )^{2\ga}} \;  ,\quad\text{when }\quad q>0 \; ,\\
\notag
    \frac{ (1+| q |   )}{ (1+ t + | q | )^{2-2\delta }  }  \; , \,\quad\text{when }\quad q<0 \; , \end{cases} \\$ which clearly does not have the desired $\frac{1}{(1+t)\cdot (1+|q| )}$ factor in the interior. That’s even in the case of $n=4$, where that factor would be $  \frac{ (1+| q |   )}{ (1+ t + | q | )^{3-2\delta }  } $: if $\de = 0$ we would then need better than  $\frac{1}{(1+t)\cdot (1+|q| )}$  as discussed earlier, and if $\de > 0$, we would do not have what is needed neither, in order to close the argument. And, this was only dealing with the non-linear structure of the good components $A_{\cal T}$.
    
 The wave operator on the bad component $A_{\underline{L}}$ have source terms, additional than ones in the source terms for the good component, which are  $ | A_L  | \cdot    | \derm A   | $ and $ | A_{e_a}  | \cdot    | \derm A_{e_a}  | $, which both have their problems even when one looks at the problem from now on in the exterior region only. In fact, if one ignores the components and looks at them as  $ | A  | \cdot    | \derm A  | $, then one would have to deal with terms which are of the type $$ | \derm A  |  \cdot  \begin{cases} \frac{1}{ (1+ t + | q | )^{1-\delta }  (1+| q |   )^{\ga}} \; ,\quad\text{when }\quad q>0 \;,\\
  \frac{(1+| q |   )^{\frac{1}{2} }}{ (1+ t + | q | )^{1-\delta }  } \;  , \,\quad\text{when }\quad q<0\; , \end{cases} $$
    and 
$$ |A| \cdot \begin{cases}  \frac{1 }{(1+t+|q|)^{1-\delta} (1+|q|)^{1+\ga}} \; ,\quad\text{when }\quad q>0 \; ,\\
        \frac{1  }{(1+t+|q|)^{1-\delta}(1+|q|)^{\frac{1}{2} }} \; , \quad\text{when }\quad q<0 \; . \end{cases} $$
    
But what we are looking for in order to close the continuity argument, are terms of the type $\frac{1}{(1+t)} \cdot | \derm A  |$ and $\frac{1}{ (1+t) \cdot (1+|q|) } \ | A  |$, which fail by large to be satisfied in the interior region. Even if we look at the exterior region only, we would then have $\frac{1}{(1+t)^{1-\de}} \cdot | \derm A  |$ and $\frac{1}{ (1+t)^{1-\de} \cdot (1+|q|) } \cdot |  A  |$ which falls quite short to close any argument of this kind using Grönwall inequality, due to the $\de > 0$ loss: no matter how small $\de$ is, this would give exponential growth of the energy. Again, if $\de = 0$\;, this would then give polynomial growth of the energy, but in that case, we have already assumed a uniform bound on the energy (by taking $\de = 0$), and a polynomial growth would not be able to improve such an a priori bound. On has also to always keep in mind the counter example of John, \cite{John}: if such an argument could close, then it would also work for the counter example $\Box \phi = (\pa_t \phi )^2$ for which we already know that solutions blow up in finite time. Hence, the argument should take into consideration the special structure of the non-linearity of these wave equations.
 
Looking, from now on, only in the exterior region (by exterior here we now mean the region $q \geq q_0$, where $q_0$ could be negative), by then we do not have any loss of decay in time $t$ in the interior of the outgoing light cone (whose tip is the origin, prescribed by $q < 0$), because by then, $q$ is bounded in that interior, i.e. $|q| \leq |q_0|$, for $q < 0$.

It turns out that in the Lorenz gauge condition and in the null frame constructed using wave coordinates, one can deal with all the terms in the non-linear structure of the covariant wave equations on $A$ and on $h^1$ (in \eqref{estimateonnonlinearstructureofgoodcomponentforAinintrrod}, \eqref{estimateonnonlinearstructureofgoodcomponentforhoneinintrrod}, \eqref{estimateonnonlinearstructureofbadcomponentforAinintrrod} and in \eqref{estimateonnonlinearstructureofbadcomponentforhoneinintrrod}) and always by separating “good” components and “bad” components. In particular, one can deal with the most troublesome terms, namely $A_{e_a}  \cdot     \derm A_{e_a} $ (see Subsection \ref{dealingiwthbadtermAaderAa}) and $ A_L   \cdot   \derm A $ (see Subsection \ref{dealingwithALderA}) as well as from the new problems generated in upgrading the a priori dispersive estimates for the Lie derivatives (see \ref{dealingwithLiederivativesupgrade}): problems which do not arise for the Einstein vacuum equation nor for the Einstein-Maxwell system. This is the specific structure of the Einstein-Yang-Mills equations.

Finally, we will prove Theorem \ref{Thetheoremofexteriorstabilityfornequalthree}. We refer the reader to the Section \ref{BrielfSummary}.

\section{Ingredients for the proof}

\subsection{Definitions for the perturbations of the Minkowski metric}\

  \begin{definition}\label{definitionoftheMinkwoskimetricminwavecoordinates}
In wave coordinates, we define $m$ to be Minkowski metric. This means that in the system of wave coordinates $\{x^0, ..., x^n\}$\;, we prescribe $m$ by
\beaa
m_{00}&=&-1 \;,\qquad  m_{ii}=1\;,\quad \text{if}\quad i=1, ...,n,\\
\quad\text{and}\quad m_{\mu\nu}&=&0\;,\quad \text{if}\quad \mu\neq \nu \quad \text{for} \quad  \mu, \nu \in \{0, 1, ..., n\} .
\eeaa
\end{definition}

  \begin{definition}
We define $h$ as the 2-tensor given by:
\bea
h_{\mu\nu} = g_{\mu\nu} - m_{\mu\nu} \, .
\eea
\end{definition}

  \begin{definition}
Let $m^{\mu\nu}$ be the inverse of $m_{\mu\nu}$. We define
\bea\label{definitionofsmallh}
h^{\mu\nu} &=& m^{\mu\mu^\prime}m^{\nu\nu^\prime}h_{\mu^\prime\nu^\prime} \\
\label{definitionsofbigH}
H^{\mu\nu} &=& g^{\mu\nu}-m^{\mu\nu}.
\eea

\end{definition}
  \begin{definition}
 We define $\der^{(\bf{m})}$ to be the covariant derivative associated to the flat metric $m$\;. Given the definition of $m$ in Definition \ref{definitionoftheMinkwoskimetricminwavecoordinates}), the Christoffel symbols are vanishing in wave coordinates, and therefore, for all $\mu, \nu \in \{0, 1, \ldots, n \}$\;,
 \bea
{ \der^{(\bf{m})}}_{e_{\mu}} e_\nu := 0 \; ,
 \eea
 where
 \bea
 e_\mu = \frac{\pa}{\pa x_\mu} \;,
 \eea
 and where $\{x^0, x^1, \ldots, x^n \}$ are the wave coordinates.
\end{definition}

\subsection{The Minkowski vector fields}\label{TheMinkowksivectorfieldsdefinition} \

Let $x^{\mu}$ be the system of wave coordinates, and let $ x^{0} = t$\;. We define 
\beaa
x_{\b} = m_{\mu\b} x^{\mu} \, ,
\eeaa
where we raised and lowered indices with respect to the Minkowski metric $m$, defined in wave coordinates to be the Minkowski metric. We define
\beaa
Z_{\a\b} = x_{\b} \pa_{\a} - x_{\a} \pa_{\b}  \; ,\\
S = t \pa_t + \sum_{i=1}^{3} x^i \pa_{i} \; .
\eeaa
The Minkowski vector fields will be denoted by $Z$\;, and they are defined such that
\bea
Z \in {\cal Z}  := \big\{ Z_{\a\b}\,,\, S\,,\, \pa_{\a} \, \,  | \, \,   \a\,,\, \b \in \{ 0, \ldots, 3 \} \big\} \; .
\eea

\begin{definition}\label{definitionofZI}
The family ${\cal Z}$ has  $11$ vector fields: $6$ vectors for the Lorentz boosts and rotations, $4$ space-time translations and one scaling vector field. One can order them and assign to each vector an $11$-dimensional integer index $(0, \ldots, 1, \ldots,0)$. Hence, a collection of $k$ vector fields from the family ${\cal Z}$, can be described by the set $I=(\iota_1, \ldots,\iota_k)$, where each $\iota_i$ is an $11$-dimensional integer, where $|I|=k = \sum_{i=1}^{k} | \iota_i |$, with $|\iota_i|=1$\;. We define
\bea
Z^I :=Z^{\iota_1}\ldots Z^{\iota_k} \quad \text{for} \quad I=(\iota_1, \ldots,\iota_k),  
\eea
where $\iota_i$ is an 11-dimensional integer index, with $|\iota_i|=1$, and $Z^{\iota_i}$ representing each a vector field from the family ${\cal Z}$.

For a tensor $T$, of arbitrary order, either a scalar or valued in the Lie algebra, we define the Lie derivative as
\bea
\Lie_{Z^I} T :=\Lie_{Z^{\iota_1}} \ldots \Lie_{Z^{\iota_k}} T \quad \text{for} \quad I=(\iota_1, \ldots,\iota_k) .
\eea

Furthermore, we make the following definition: when we write $I = I_1+I_2$, it means that we divided the set $I$ into two sets $I_1$ and in $I_2$, while preserving the order of $I$ in $I_1$ and in $I_2$, i.e., if $I=(\iota_1, \ldots,\iota_k)$, then $I_1=(\iota_{i_1}, \ldots,\iota_{i_n})$ and $I_2=(\iota_{i_{n+1}}, \ldots,\iota_{i_k})$, where $i_1< \ldots<i_n$ and $i_{n+1}< \ldots<i_k$. By a sum $\sum_{I_1+I_2=I}  $\;, we mean that we make the sum over all such partitions for a given $I$. With this convention, the Leibnitz rule holds and reads for sufficiently smooth functions $f$ and $g$, 
\bea
Z^I(f \cdot g)=\sum_{I_1+I_2=I} (Z^{I_1} f) \cdot (Z^{I_2} g) \; .
\eea

\end{definition}
We also refer the reader to our previous paper \cite{G4} for more details.

\subsection{Notations and estimates on the perturbations}\

\begin{definition}\label{definitionofbigOforLiederivatives}

For a family of tensors Let $\Lie_{Z^{I_1}}K^{(1)}, \ldots, \Lie_{Z^{I_m}} K^{(m)}$, where each tensor $ K^{(l)}$ is again either $A$ or $h$ or $H$\;, or $\derm A$\;, $\derm h$ or $\derm H$\;, we define
\bea
\notag
&& O_{\mu_{1} \ldots \mu_{k} } (\Lie_{Z^{I_1}}K^{(1)} \cdot \ldots \cdot \Lie_{Z^{I_m}} K^{(m)}  ) \\
\notag
&:=& \prod_{l=1}^{m} \Big[  \prod_{|J_l| \leq |I_l|} Q_{1}^{J_l} ( \Lie_{Z^{J_l}} K^{(l)} ) \cdot \Big( \sum_{n=0}^{\infty} P_{n}^{J_l}  ( \Lie_{Z^{J_l}} K^{(l)} ) \Big) \Big] \; . \\
\eea
where again $P_n^{J_l} (K^{l} )$ and $Q_1^{J_l} (K)$, are tensors that are Polynomials of degree $n$ and $1$, respectively, with $Q^{J_l}_1 (0) = 0$ and $Q^{J_l}_1 \neq 0$\;, of which the coefficients are components in wave coordinates of the metric $\textbf m$ and of the inverse metric $\textbf m^{-1}$, and of which the variables are components in wave coordinates of the covariant tensor $\Lie_{Z^{J_l}} K^{l}$, leaving some indices free, so that at the end the whole product $$  \prod_{l=1}^{m} \Big[  \prod_{|J_l| \leq |I_l|}  Q_{1}^{J_l} ( \Lie_{Z^{J_l}} K^{(l)} ) \cdot \Big( \sum_{n=0}^{\infty} P_{n}^{J_l}  ( \Lie_{Z^{J_l}} K^{(l)} ) \Big) \Big]$$ gives a tensor with free indices $\mu_{1} \ldots \mu_{k}$\;. To lighten the notation, we shall drop the indices and just write $O (\Lie_{Z^{I_1}}K^{(1)} \cdot \ldots \cdot \Lie_{Z^{I_m}} K^{(m)}  )$.

\end{definition}

We have the following lemma from a previous paper, see \cite{G4}.
\begin{lemma}\label{linkbetweenbigHandsamllh}
We have
\bea
H^{\mu\nu}= -h^{\mu\nu}+O^{\mu\nu}(h^2) \, .
\eea
\end{lemma}

\subsection{The Einstein-Yang-Mills equations in the harmonic and Lorenz gauges}\

We have proved the following lemma in \cite{G4}.
\begin{lemma}\label{Einstein-Yang-MillssysteminLorenzwavegauges}

Let,  
\bea
 P(\pa_\mu h,\pa_\nu h) :=\frac{1}{4} m^{\alpha\alpha^\prime}\pa_\mu h_{\alpha\alpha^\prime} \, m^{\beta\beta^\prime}\pa_\nu h_{\beta\beta^\prime}  -\frac{1}{2} m^{\alpha\alpha^\prime}m^{\beta\beta^\prime} \pa_\mu h_{\alpha\beta}\, \pa_\nu h_{\alpha^\prime\beta^\prime} \, ,
\eea and
 
\bea
\notag
&& Q_{\mu\nu}(\pa h,\pa h) \\
\notag
&:=& \pa_{\alpha} h_{\beta\mu}\, \, m^{\alpha\alpha^\prime}m^{\beta\beta^\prime} \pa_{\alpha^\prime} h_{\beta^\prime\nu} -m^{\alpha\alpha^\prime}m^{\beta\beta^\prime} \big(\pa_{\alpha} h_{\beta\mu}\,\,\pa_{\beta^\prime} h_{\alpha^\prime \nu} -\pa_{\beta^\prime} h_{\beta\mu}\,\,\pa_{\alpha} h_{\alpha^\prime\nu}\big)\\
\notag
&& +m^{\a\a'}m^{\b\b'}\big (\pa_\mu h_{\a'\b'}\, \pa_\a h_{\b\nu}- \pa_\a h_{\a'\b'} \,\pa_\mu h_{\b\nu}\big )  + m^{\a\a'}m^{\b\b'}\big (\pa_\nu h_{\a'\b'} \,\pa_\a h_{\b\mu} - \pa_\a h_{\a'\b'}\, \pa_\nu h_{\b\mu}\big )\\
\notag
&& +\frac 12 m^{\a\a'}m^{\b\b'}\big (\pa_{\b'} h_{\a\a'}\, \pa_\mu h_{\b\nu} - \pa_{\mu} h_{\a\a'}\, \pa_{\b'} h_{\b\nu} \big ) +\frac 12 m^{\a\a'}m^{\b\b'} \big (\pa_{\b'} h_{\a\a'}\, \pa_\nu h_{\b\mu} - \pa_{\nu} h_{\a\a'} \,\pa_{\b'} h_{\b\mu} \big ) \, ,\\
\eea
and
\bea
G_{\mu\nu}(h)(\pa h,\pa h) :=  O  (h \cdot (\pa h)^2) \, ,
\eea
i.e. $G_{\mu\nu}(h)(\pa h,\pa h) $ is a quadratic form in $\pa h$ with
coefficients smoothly dependent on $h$ and vanishing when $h$
vanishes: $G_{\mu\nu}(0)(\pa h,\pa h)=0$.

The Einstein-Yang-Mills equations in both the harmonic and Lorenz gauges are the following equations given in \eqref{nonlinearwaveequationonAusingpartialderivativesinwavecoordinatessystemwitharemarkonwritingthesources}
and \eqref{nonlinearwaveequationonhusingpartialderivativesinwavecoordinatessystemwithPQG}. 

That is on one hand, a non-linear wave equation on the Yang-Mills potential $A$,
      \bea\label{nonlinearwaveequationonAusingpartialderivativesinwavecoordinatessystemwitharemarkonwritingthesources}
   \notag
g^{\la\mu} \pa_{\la}   \pa_{\mu}   A_{\si}      &=&  m^{\a\ga} m ^{\mu\la}  (  \pa_{\si}  h_{\ga\la} )   \pa_{\a}A_{\mu}       +   \frac{1}{2}  m^{\a\mu}m^{\b\nu}   \big(   \pa_\a h_{\b\si} + \pa_\si h_{\b\a}- \pa_\b h_{\a\si}  \big)   \cdot  \big( \pa_{\mu}A_{\nu} - \pa_{\nu}A_{\mu}  \big) \\
 \notag
&& +      \frac{1}{2}  m^{\a\mu}m^{\b\nu}   \big(   \pa_\a h_{\b\si} + \pa_\si h_{\b\a}- \pa_\b h_{\a\si}  \big)   \cdot   [A_{\mu},A_{\nu}] \\
 \notag
 && -  m^{\a\mu} \big(  [ A_{\mu}, \pa_{\a} A_{\si} ]  +    [A_{\alpha},  \pa_{\mu}  A_{\si} - \pa_{\si} A_{\mu} ]    +    [A_{\alpha}, [A_{\mu},A_{\si}] ]  \big)  \\
 \notag
  && + O( h \cdot  \pa h \cdot  \pa A) + O( h \cdot  \pa h \cdot  A^2) + O( h \cdot  A \cdot \pa A) + O( h \cdot  A^3) \, ,\\
  \eea
            \begin{remark}
      We could have written the term $m^{\a\mu} \big(  [ A_{\mu}, \pa_{\a} A_{\si} ]  +    [A_{\alpha},  \pa_{\mu}  A_{\si} - \pa_{\si} A_{\mu} ] \big) $ as $m^{\a\mu} \big(  2 [ A_{\mu}, \pa_{\a} A_{\si} ]  -   [A_{\alpha},   \pa_{\si} A_{\mu} ] \big) $\;, but we prefer to keep it in this form.
      \end{remark}

On the other hand, a non-linear wave equation on the metric $h$\;,
  \bea\label{nonlinearwaveequationonhusingpartialderivativesinwavecoordinatessystemwithPQG}
\notag
 g^{\alpha\beta}\pa_\alpha\pa_\beta
h_{\mu\nu} &=& P(\pa_\mu h,\pa_\nu h) + Q_{\mu\nu}(\pa h,\pa h)  +G_{\mu\nu}(h)(\pa h,\pa h)  \\
\notag
 &&   -4   m^{\si\b} \cdot  <   \pa_{\mu}A_{\b} - \pa_{\b}A_{\mu}  ,  \pa_{\nu}A_{\si} - \pa_{\si}A_{\nu}  >    \\
 \notag
 &&   + \frac{2}{(n-1)} m_{\mu\nu }  m^{\si\b}  m^{\a\la}    \cdot  <  \pa_{\a}A_{\b} - \pa_{\b}A_{\a} , \pa_{\la}A_{\si} - \pa_{\si}A_{\la} >   \\
 \notag
&&           -4 m^{\si\b}  \cdot  \big( <   \pa_{\mu}A_{\b} - \pa_{\b}A_{\mu}  ,  [A_{\nu},A_{\si}] >   + <   [A_{\mu},A_{\b}] ,  \pa_{\nu}A_{\si} - \pa_{\si}A_{\nu}  > \big)  \\
\notag
&& + \frac{2}{(n-1)} m_{\mu\nu }  m^{\si\b}  m^{\a\la}    \cdot \big(  <  \pa_{\a}A_{\b} - \pa_{\b}A_{\a} , [A_{\la},A_{\si}] >    +  <  [A_{\a},A_{\b}] , \pa_{\la}A_{\si} - \pa_{\si}A_{\la}  > \big) \\
\notag
 &&  -4 m^{\si\b}  \cdot   <   [A_{\mu},A_{\b}] ,  [A_{\nu},A_{\si}] >      +  \frac{2}{(n-1)} m_{\mu\nu }  m^{\si\b}  m^{\a\la}   \cdot   <  [A_{\a},A_{\b}] , [A_{\la},A_{\si}] >  \\
     && + O \big(h \cdot  (\pa A)^2 \big)   + O \big(  h  \cdot  A^2 \cdot \pa A \big)     + O \big(  h   \cdot  A^4 \big)    \, .
\eea
\end{lemma}
            \begin{remark}
      In three space dimensions, we are actually interested in the wave equation on the perturbation $h^1 := h - h^0$\;, where $h^0$ is defined in \eqref{definitionofthesphericallysymmtericpartofinitialdata} (see Subsection \ref{subsectionontheenergynormwithdiscussiononwhyweconsiderhoneandnotfullh}). 
      \end{remark}

\subsection{The weights}\

Let $q := r - t $\,, where $t$ and $r$ are defined using the wave coordinates, as we explained in \cite{G4}. We recapitulate the following definitions from \cite{G5}.

\begin{definition}\label{defoftheweightw}
We define $w$ by
\beaa
w(q):=\begin{cases} (1+|q|)^{1+2\gamma} \quad\text{when }\quad q>0 , \\
         1 \,\quad\text{when }\quad q<0 . \end{cases}
\eeaa
for some $\gamma > 0$\,.

\end{definition}

\begin{definition}\label{defwidehatw}
We define $\widehat{w}$ by 
\beaa
\widehat{w}(q)&:=&\begin{cases} (1+|q|)^{1+2\gamma} \quad\text{when }\quad q>0 , \\
        (1+|q|)^{2\mu}  \,\quad\text{when }\quad q<0 , \end{cases} \\
        &=&\begin{cases} (1+ q)^{1+2\gamma} \quad\text{when }\quad q>0 , \\
      (1 - q)^{2\mu}  \,\quad\text{when }\quad q<0 ,\end{cases} 
\eeaa

for $\ga > 0$ and $\mu < 0$\,. Note that the definition of $\widehat{w}$\,, gives for $\ga \neq - \frac{1}{2} $ and $\mu \neq 0$ (which is assumed here), that 
\beaa
\widehat{w}^{\prime}(q) \sim \frac{\widehat{w}(q)}{(1+|q|)} \; ,
\eeaa
(see Lemma \ref{derivativeoftildwandrelationtotildew}).

\end{definition}
\begin{remark}\label{remarkonwhyweneedmustrictlynegative}
We take $\mu$ non-zero so that we would have for $q<0$, the derivative $\frac{\pa \widehat{w}}{\pa q}$ to be non-vanishing, so that to have a control on certain space-time integral of tangential derivatives, which is needed here (see Corollary \ref{theenergyestimatewithcontrolontangentialderivativesusinghatwoveroneplusq}).
Also, we take $\mu$ negative because we want the derivative $\frac{\pa \widehat{w}}{\pa q} > 0$\,, as we will see that this is what we need in order to obtain an energy estimate on the fields (see Corollary \ref{energyestimatewithoutestimatingthetermsthatinvolveBIGHbutbydecomposingthemcorrectlysothatonecouldgettherightestimate}). In other words, $\mu < 0$ is a necessary condition to ensure that $\widehat{w}^{\prime} (q)$ enters with the right sign in the energy estimate (see Corollary \ref{WeightedenergyestimateusingHandnosmallnessyet}).
\end{remark}

\begin{definition}\label{defwidetildew}
We define $\widetilde{w}$ by 
\beaa
\widetilde{w} ( q)&:=&  \widehat{w}(q) + w(q) \\
&:=&\begin{cases} 2 (1+|q|)^{1+2\gamma} \quad\text{when }\quad q>0 , \\
       1+  (1+|q|)^{2\mu}  \,\quad\text{when }\quad q<0 . \end{cases} \\    
\eeaa
Note that the definition of $\widetilde{w}$ is constructed so that Lemma \ref{equaivalenceoftildewandtildeandofderivativeoftildwandderivativeofhatw} holds, which we need in order to obtain Corollary \ref{energyestimatewithoutestimatingthetermsthatinvolveBIGHbutbydecomposingthemcorrectlysothatonecouldgettherightestimate}.
\end{definition}

We recapitulate the following lemmas from \cite{G5}.
\begin{lemma}\label{equaivalenceoftildewandtildeandofderivativeoftildwandderivativeofhatw}
We have
\bea
\widetilde{w}^{\prime}  &\sim & \widehat{w}^{\prime } \; .
\eea
Furthermore, for $\mu < 0$\,(which is the case here), we have
\bea
\widetilde{w} ( q)& \sim & w(q) \; .
\eea

\end{lemma}

\begin{lemma}\label{derivativeoftildwandrelationtotildew}
Let $\widetilde{w}$ as in Definition \ref{defwidetildew},
We have, for $\ga \neq - \frac{1}{2} $ and $\mu \neq 0$\,,  
\bea
\widehat{w}^{\prime}(q) \sim \frac{\widehat{w}(q)}{(1+|q|)} \; .
\eea

\end{lemma}

\subsection{The exterior region}\

\begin{definition}\label{definitionofSigmatexterasintersectedwithcomplementofinteriorofcausaldomainofK}
 Let $\cal{K}$ be a compact region of the initial slice $\Sigma$ and let $C(\cal K) $ be the future causal domain of dependance of $\cal{K}$ for the metric $g$\;. Now, let $\overline{C}(\cal K)$ be the complement of $C$ in the future domain of dependance $\Sigma_{t_0}$ for the metric $g$\;. 
 We define for a given fixed $t$\;,
   \bea
\Sigma^{ext}_{t}(\cal K)   &:=&  \Sigma_t  \cap \overline{C} (\cal K) \, .
\eea
We shall write $\Sigma^{ext}_{t}$ to denote $\Sigma^{ext}_{t}(\cal K)$\;.  

We define $N_{t_1}^{t_2}$ as the null boundary (null for the metric $g$) of $\overline{C}$ restricted to times between $t_1$ and $t_2$\;.
  \begin{remark}
Yet, in Definition \ref{definitionofSigmaexteriorattequalconstantusingthestressenergymomentumtensorTgforwaveeq}, we give a different definition for $\Sigma^{ext}_{t}$ and in Definition \ref{definitionofNandofNtruncatedbtweentwots}, we give a different definition for $N_{t_1}^{t_2}$\;. However, we note that these definitions are equivalent for our purposes, i.e. in our context, the estimates that we get on $\Sigma^{ext}_{t}$ from one definition (of the two Definitions \ref{definitionofSigmatexterasintersectedwithcomplementofinteriorofcausaldomainofK}
 and \ref{definitionofSigmaexteriorattequalconstantusingthestressenergymomentumtensorTgforwaveeq}) could be translated into the same estimates with the other definition (see Remark \ref{RemarkaboutthedifferenceofdefinitionofSigmaexterior}). Thus, for moment, we can think of $\Sigma^{ext}_{t}$ as defined in Definition \ref{definitionofSigmatexterasintersectedwithcomplementofinteriorofcausaldomainofK}. 
  \end{remark}

\end{definition}

\subsection{Weighted Klainerman-Sobolev inequality in the exterior}\

We will now state a Klainerman-Sobolev inequality which holds true in the exterior region $\overline{C}$\;, the complementary of $C$ (the future causal domain of dependance for the metric $g$\;, of the compact $\cal K \subset \Sigma$). 

The integration for the $L^2$ will be supported only on the exterior regions $\Sigma^{ext}_{t} = \Sigma_{t} \cap \overline{C}$\;. We have globally the following pointwise estimate in the exterior region $\overline{C}$ for any smooth scalar function $\phi$ vanishing at spatial infinity, i.e. $\lim_{r \to \infty} \phi (t, x^1, x^2, x^3) = 0$\;,
\bea
\notag
|\phi(t,x)| \cdot (1+t+|q|) \cdot \big[ (1+|q|) \cdot w(q)\big]^{1/2} \leq
C \cdot \sum_{|I|\leq 2 } \|\big(w(q)\big)^{1/2} Z^I \phi(t,\cdot)\|_{L^2 (\Sigma^{ext}_{t} ) } \; , \\
\eea
where here the $L^2(\Sigma^{ext}_{t} )$ norm is taken on $\Sigma^{ext}_{t} $ slice with respect to the Lebesgue measure $dx_1 \ldots dx_n$ obtained from wave coordinates.

\subsection{Definition of the norms}\

We define the tensor $E$ as the euclidian metric in wave coordinates, that is prescribed by
\bea
\notag
E_{\mu\nu} = m(\frac{\pa}{\pa x^\mu}, \frac{\pa}{\pa x^\nu})+2 m(\frac{\pa}{\pa x^\mu}, \frac{\pa}{\pa t}) \cdot m(\frac{\pa}{\pa x^\nu}, \frac{\pa}{\pa t}) \, .
\eea
For a tensor of arbitrary order, say $K_\a$\;, we define
\bea
 | K  |^2= E^{\mu\nu}  < K_{\mu} ,   K_{\nu} >  \, .
\eea
Therefore, we have
\bea
 | K  |^2 &=&  \sum_{\a \in  \{t, x^1, \ldots, x^n \}} | K_\a |^2   \, ,
\eea

We then define for a  tensor of arbitrary order, for example $K_{\a}$\;, 
\bea
 | \pa K |^2 := | \derm K |^2 :=  E^{\a\b} E^{\mu\nu}    { \der^{(\bf{m})}}_{\mu} K_\a \ \cdot    { \der^{(\bf{m})}}_{\nu} K_\b \; .
 \eea
We then have
 \bea
 \notag
| \pa K |^2  &=&   |  {\der^{(\bf{m})}}_{t} K |^2  +  |  {\der^{(\bf{m})}}_{x^1} K |^2 +\ldots  + |  {\der^{(\bf{m})}}_{x^n} K |^2 \\
&=&  \sum_{\a,\; \b \in  \{t, x^1, \ldots, x^n \}} |  \pa_{\a} K_\b |^2   \, ,
\eea
since in wave coordinates, the Minkowski covariant derivative of coordinate vector fields are vanishing.

We recapitulate the following Lemma from \cite{G4}.
\begin{lemma}
At a point $p$ of the space-time, let $x^{\mu}$  be the wave coordinate system. For a sufficiently smooth function $f$ and for a norm $\mid .\mid $\,, we define for all $I$ and $Z^I$ as previously defined, the following norm in the wave coordinates system $\{t, x^1, \ldots, x^n \}$\,,
\bea
| Z^I \pa f | := \sqrt{ | Z^I \pa_{t} f  |^2 + \sum_{i=1}^{n} | Z^I  \pa_{i} f  |^2 } \;  .
\eea
Then, we have,
\bea
|  Z^I \pa  f |  \leq C (|I| ) \cdot    \sum_{|J| \leq |I| }  | \pa  (  Z^J  f ) | \; ,
\eea
where $C (|I| )$ is a constant that depends only on $|I|$\,.
\end{lemma}

\subsection{The energy norm}\label{subsectionontheenergynormwithdiscussiononwhyweconsiderhoneandnotfullh}

We recall that we are given an initial data set which we write as $(\Sigma, \overline{A}, \overline{E}, \overline{g}, \overline{k})$\;, and that $\Sigma$ is diffeomorphic to $\R^n$\;, and therefore there exists a global system of coordinates $(x^1, ..., x^n) \in \R^n$ for $\Sigma$\;. We define
\bea
r := \sqrt{ (x^1)^2 + ...+(x^n)^2  }\;.
\eea
We assume that the initial data set is smooth and asymptotically flat. Now, we define a smooth function $\chi$\;, given by
 \bea\label{defXicutofffunction}
\chi (r)  := \begin{cases} 1  \quad\text{for }\quad r \geq \frac{3}{4} \;  ,\\
0 \quad\text{for }\quad r \leq \frac{1}{2} \;. \end{cases} 
\eea
We let $\de_{ij}$ be the Kronecker symbol, and we define $\overline{h}^1_{ij} $ in this system of coordinates $x^i$\;, for $i \in \{1, \ldots, n \}$\;, by
 \bea
\overline{h}^1_{ij} := \overline{g}_{ij} - (1 + \chi (r)\cdot  \frac{2\cdot M}{r}  ) \de_{ij} \; .
\eea
We define
 \bea
\overline{h}^0_{ij} :=  \chi (r)\cdot  \frac{2\cdot M}{r}   \de_{ij} \; .
\eea
Thus, the initial data can be written as
\bea
\overline{g}= \overline{h}^1 +\overline{h}^0 +  \de_{ij} \; .
\eea

We foresee that the propagation of $\de_{ij}$  will be the Minkowksi metric in wave coordinates. Now, we do not know how the propagation of $\overline{h}^0$ will look like. We will follow a guess made by Lindblad and Rodnianksi in \cite{LR10}, by defining $h^0$ for $t >0$ as 
\bea\label{guessonpropagationofthesphericallsymmetricpart}
h^0_{\mu\nu} := \chi(r/t) \cdot \chi(r)\cdot \frac{M}{r}\de_{\mu\nu}  \;  ,
\eea
and for $t=0$\;,  \bea\label{definitionofthesphericallysymmtericpartofinitialdata}
h^0_{\mu\nu} ( t= 0) := \chi(r)\cdot \frac{M}{r}\cdot \de_{\mu\nu}  \; ,
\eea

We then define $h^1$ as
\bea
h^{1}_{\mu\nu} := g_{\mu\nu} - m_{\mu\nu}- h^0_{\mu\nu}  \;  .
\eea

We define the higher order energy norm as the following $L^2$ norms on $A$ and $h^1$ in the exterior, 
 \bea\label{definitionoftheenergynorm}
 \notag
\E_{N} (t) :=  \sum_{|I|\leq N} \big(   \|w^{1/2}   \derm ( \Lie_{Z^I}  A   (t,\cdot) )  \|_{L^2 (\Sigma_t^{ext}) } +  \|w^{1/2}   \derm ( \Lie_{Z^I} h^1   (t,\cdot) )  \|_{L^2 (\Sigma_t^{ext})} \big) \, , \\
\eea
where the integration is taken with respect to the Lebesgue measure $dx_1 \ldots dx_n$\:.

\subsection{The bootstrap argument} \label{Thebootstrapargumentandnotationonboundingtheenergy}\

It is also called a continuity argument. We start with a local solution defined on a maximum time interval $[0, T_{\text{loc}})$ and that is well-posed in the energy norm $\E_{N} (t)$ for some $N \in \N$\;. This means that the time dependance of the energy $\E_{N} (t)$ is continuous. Furthermore, by maximality of $T_{\text{loc}}$ and the well-posedness of the solution, the time interval for the local solution must be excluding $T_{\text{loc}}$\,, otherwise the energy will be finite at $t=T_{\text{loc}}$ and this means that we could extend the local solution again beyond the time $t=T_{\text{loc}} $ by repeating the argument for establishing a local solution starting at time $t = T_{\text{loc}} $. 

Differently speaking, the maximal $T_{\text{loc}}$ is characterised by
\beaa
\lim_{t \to T_{\text{loc}} } \E_{ N } (t)  =  \infty \;.
\eeaa
We look at any time $T \in [0, T_{\text{loc}})$, such that for all $t$ in the interval of time $[0, T]$, we have
\bea\label{aprioriestimate}
\E_{ N } (t)  \leq E (N )  \cdot \eps \cdot (1 +t)^\delta \;,
\eea
where $E ( N )$\, is a constant that depends on $ N $\,, where $\eps > 0$ is a constant to be chosen later small enough, and where $\delta \geq 0$ is to be chosen later. In addition, we start with an initial data such that this estimate holds true for $t=0$\,, i.e.
\bea \label{theboostrapimposestheconditiononinitialdatasothatnonemptyset}
\E_{ N } (0)  \leq E (N ) \cdot  \eps \;,
\eea
and thus we know that such a $T$ exists, since at least $T=0$ satisfies the estimate.

We will then show that for $t \in [0, T]$\,, the same estimate holds true but with $\eps$ replaced with $\frac{\eps}{2}$\,, i.e. we then prove that for all $t$ in the time interval $[0, T]$\,,
\bea\label{improvedapriori}
\E_{ N } (t)  \leq E ( N )  \cdot \frac{\eps}{2}   \cdot (1 +t)^\delta \; .
\eea

Consequently, we would have shown that the set
$$\{ T  \; \;  | \; \text{for all} \quad  t \in [0, T]\,,  \quad  \E_{ N } (t)  \leq E (N )  \cdot \eps  \cdot (1 +t)^\delta  \}$$ is relatively open in $[0, T_{\text{loc}})$\,, non-empty since $0$ belongs to the set, and we know that it is relatively closed in $[0, T_{\text{loc}})$ since the map $t \to \E_{N} (t)$  is continuous, and thus, the set is the whole interval $[0, T_{\text{loc}})$\,.

As a result, we would have shown that for all $t \in [0, T_{\text{loc}})$\,, we have
\beaa
\E_{ N } (t)  \leq E ( N )  \cdot \frac{\eps}{2}  \cdot (1 +t)^\delta \;.
\eeaa
As a result, we have
\beaa
\lim_{t \to T_{\text{loc}} } \E_{ N } (t)  \leq E ( N )  \cdot \frac{\eps}{2}  \cdot (1 +T_{\text{loc}})^\delta < \infty \;.
\eeaa
By continuity of the energy, this means that  $\E_{ N } (T_{\text{loc}})$ is finite and we can then repeat the argument for establishing a local solution starting at time $t = T_{\text{loc}} $ which would lead to a local solution defined beyond the time $t=T_{\text{loc}} $\,, which contradicts the maximality of $T_{\text{loc}} $\,.

\subsection{The bootstrap assumption}\label{Theassumptionforthebootstrapandthenotation}\

To run our continuity argument, we start by assuming that for all $k \leq N \in \N$\,, where $N$ is to be determined later, we have
  \bea\label{bootstrap}
\E_{ k } (t)  \leq E (k )  \cdot \eps \cdot (1 +t)^\delta \;.
\eea
We choose 
\bea\label{epsissmallerthan1}
0 < \eps \leq 1\;,
\eea
so that any powers of $\eps$ are in fact bounded by $\eps$\,. To lighten the notation, we also choose here
\bea\label{assumptionontheconstantboundsforenergy}
E ( k ) \leq 1\;,
\eea
so that any sum of powers of $E ( k )$ is in fact bounded by a constant multiplied by $E ( k )$\,. Furthemore, we choose
\bea\label{assumptionontheorderingoftheconstantboundsforenergy}
E ( k_1 )  \leq E(k_2)\;,
\eea
for all $k_1 \leq k_2$\,, with $k_1\, , k_2 \in \N$\,, given the fact that $\E ( k_1 )  \leq \E(k_2)$. The reason we choose to put the constants $E ( k )$\,, rather than want an $\eps$ to be fixed, is to show in the estimates the dependance on the energy and mainly, on the number of Lie derivatives involved. In other words, these constants $E ( k)$ are not needed but are there to make clearer in the argument the number of Lie derivatives of fields for which we use the bootstrap assumption.

\subsection{A priori decay estimates}\

The a priori estimates are decay estimates that are generated from the weighted Sobolev inequality combined with the bootstrap assumption \eqref{bootstrap}. The a priori estimates have nothing to do with the Einstein-Yang-Mills equations, but they come from the fact that we chose the energy to be in the form of what dominates the right hand side of the Klainerman-Sobolev inequality when applied to  $\derm \Lie_{Z^I} A$ and $ \derm \Lie_{Z^I} h^1$\;. The fact that these estimates are generated from the bootstrap assumption, and are not proven yet to be true estimates, is the reason why we call them “a priori decay estimates”.

\begin{lemma}\label{apriordecayestimatesfrombootstrapassumption}
Under the bootstrap assumption \eqref{bootstrap}, taken for $N =  |I| +  2$, if for all $\mu, \nu \in (t, x^1, \ldots, x^n)$, and for any functions $ \pa_\mu \Lie_{Z^I} h^1_\nu \, ,  \pa_\mu \Lie_{Z^I} A_\nu  \in C^\infty_0(\R^n)$, then we have in the exterior region $\overline{C} \subset \{ q \geq q_0 \}$,
 \bea
 \notag
|\derm  ( \Lie_{Z^I}  A ) (t,x)  |    &\leq & \begin{cases} C ( |I| ) \cdot E ( |I| + 2)  \cdot \frac{\eps }{(1+t+|q|)^{1-\delta} (1+|q|)^{1+\ga}},\quad\text{when }\quad q>0,\\
 \notag
       C ( |I| ) \cdot E ( |I| + 2)  \cdot \eps \cdot \frac{ (1+|q|)^{\frac{1}{2} }}{ (1+t+|q|)^{1-\delta} } \,\quad\text{when }\quad q<0 , \end{cases} \\
        \notag
       &\leq& C(q_0) \cdot C ( |I| ) \cdot E ( |I| + 2 )  \cdot \frac{\eps }{(1+t+|q|)^{1-\delta} \cdot  (1+|q|)^{1+\gamma}} \; , \\
      \eea
and 
 \bea
 \notag
|\derm ( \Lie_{Z^I}  h^1 ) (t,x)  |   &\leq & \begin{cases} C ( |I| ) \cdot E ( |I| + 2)  \cdot \frac{\eps }{(1+t+|q|)^{1-\delta} (1+|q|)^{1+\ga}},\quad\text{when }\quad q>0,\\
 \notag
       C ( |I| ) \cdot E ( |I| + 2)  \cdot \eps \cdot \frac{ (1+|q|)^{\frac{1}{2} }}{ (1+t+|q|)^{1-\delta} } \,\quad\text{when }\quad q<0 , \end{cases} \\
        \notag
 &\leq& C(q_0) \cdot C ( |I| ) \cdot E ( |I| + 2)  \cdot \frac{\eps }{(1+t+|q|)^{1-\delta} \cdot  (1+|q|)^{1+\gamma}}  \; . \\
      \eea

\end{lemma}

\subsection{The spatial asymptotic behaviour of $ \Lie_{Z^I} A (t, x) $ at $t=0$}\

\begin{lemma}\label{estimateonZderivativeofafuncionbypartialderivativeoff}
We have for all vector $Z \in  {\cal Z}$, and for all sufficiently smooth function $f$, the following estimate for $t \geq 0$,
\beaa
 | Z f  |&\les& ( 1 + t +  |x  | ) \cdot | \pa f  | .
\eeaa
\end{lemma}

\begin{lemma}\label{spatialdecayoninitialdeataforanapriorestimateontheEinsteinYangMillsfieldswithoutgradiant}
If the factor $\gamma$ in the weight is such that $\gamma > \max\{0, \delta -1\} $, then under the bootstrap assumption \eqref{bootstrap}, taken for $k =  |I| +  \lfloor  \frac{n}{2} \rfloor  +1$, we have for all $ | I  | \geq 1$, 
\bea\label{boundonthehyperplanetequalzeroforzeroderivativeoftehfields}
\notag
 | \Lie_{Z^I} A (0,x)  | +   | \Lie_{Z^I}  h^1 (0,x)  |   &\les&  C ( |I| ) \cdot E ( |I| + \lfloor  \frac{n}{2} \rfloor  +1)  \cdot \frac{\eps  }{(1+r)^{1+\gamma-\delta}}  \; ,\\
\eea
and 
\bea\label{blimitatrgoestoinfinityathyperplanetequalzerozeroforzeroderivativeoftehfields}
 \lim_{ r \to \infty }  \big(  | \Lie_{Z^I}  A (0,x)  | +     | \Lie_{Z^I}  h^1 (0,x)  |  \big) &=& 0 \;.
\eea
Also, we choose to take the initial data such that \eqref{boundonthehyperplanetequalzeroforzeroderivativeoftehfields} is also true for $ | I  | = 0$, which implies \eqref{blimitatrgoestoinfinityathyperplanetequalzerozeroforzeroderivativeoftehfields}.
\end{lemma}

\begin{lemma} \label{asymptoticbehaviouratteqzeroforallfields}
 For $\gamma > \max\{0, \delta -1\} $, we have
 \bea\label{boundonthehyperplanetequalzeroforthefullderivativeoftheLieZfields}
\notag
 | \derm   ( \Lie_{Z^I} A ) A (0,x)  | +   |\derm   ( \Lie_{Z^I} A )  h^1 (0,x)  |   &\les&  C ( |I| ) \cdot E ( |I| + \lfloor  \frac{n}{2} \rfloor  +1)  \cdot \frac{\eps  }{(1+r)^{1+\gamma-\delta}}  \; ,\\
\eea
and 
 \bea
\notag
 \lim_{ r = \to \infty } \big(   |\derm ( \Lie_{Z^I}  A ) (0,x)  | +  |\derm  ( \Lie_{Z^I}   h^1 ) (0,x)  |  \big)  &=&  0  \; .\\
\eea
\end{lemma}

 \subsection{Estimates on $ \Lie_{Z^I} A $ and $ \Lie_{Z^I} h^1 $ for $t>0$.}\

Now, we will use \eqref{boundonthehyperplanetequalzeroforzeroderivativeoftehfields} in Lemma \ref{asymptoticbehaviouratteqzeroforallfields} to estimate the Lie derivatives in the direction of Minkowski vector fields of the Einstein-Yang-Mills fields $ \Lie_{Z^I} A $ and $ \Lie_{Z^I} h^1 $, for $t > 0$. This will be done by specific integration till we reach the hyperplane prescribed by $t=$ and then use  \eqref{boundonthehyperplanetequalzeroforzeroderivativeoftehfields}.

\begin{lemma}\label{apriorestimateontheEinsteinYangMillsfieldswithoutgradiant}
Under the bootstrap assumption \eqref{bootstrap}, taken for $N =  |I| + 2$, and with initial data such that
\beaa
|  A (0,x) | +  |  h^1 (0,x)|   &\les&   \frac{\eps }{ (1+|q|)^{1+\gamma-\delta}} \, ,
\eeaa
we then have for all $ | I  | $, in the exterior region $\overline{C} \subset \{ q \geq q_0 \}$,
 \bea
 \notag
| \Lie_{Z^I}  A  (t,x)  |        &\leq& \begin{cases} c (\gamma) \cdot  C ( |I| ) \cdot E ( |I| +  2)  \cdot \frac{\eps }{(1+t+|q|)^{1-\delta} (1+|q|)^{\gamma}},\quad\text{when }\quad q>0,\\
       C ( |I| ) \cdot E ( |I| + 2 )  \cdot \frac{\eps \cdot (1+| q |   )^\frac{1}{2} }{(1+t+|q|)^{1-\delta} }  \,\quad\text{when }\quad q<0 . \end{cases} \\
      \notag
      &\leq& C(q_0) \cdot c (\gamma) \cdot  C ( |I| ) \cdot E ( |I| + 2)  \cdot \frac{\eps }{(1+t+|q|)^{1-\delta} \cdot  (1+|q|)^{\gamma}} \; , \\
      \eea
and
\bea
\notag
| \Lie_{Z^I}  h^1  (t,x)  |      &\leq& \begin{cases} c (\gamma) \cdot  C ( |I| ) \cdot E ( |I| +  2)  \cdot \frac{\eps }{(1+t+|q|)^{1-\delta} (1+|q|)^{\gamma}},\quad\text{when }\quad q>0,\\
       C ( |I| ) \cdot E ( |I| + 2 )  \cdot \frac{\eps \cdot (1+| q |   )^\frac{1}{2} }{(1+t+|q|)^{1-\delta} }  \,\quad\text{when }\quad q<0 . \end{cases} \\
      \notag
      &\leq& C(q_0) \cdot c (\gamma) \cdot  C ( |I| ) \cdot E ( |I| + 2)  \cdot \frac{\eps }{(1+t+|q|)^{1-\delta} \cdot (1+|q|)^{\gamma}} \; . \\
      \eea
      
\end{lemma}

\begin{proof}

We note that for $h$ small, the exterior region $\overline{C}$ has the property that lines prescribed by $\{ (\tau, r), \quad \tau + r = constant \}$ intersect the boundary of $\overline{C}$ only once. Differently speaking, each point $(t, |x|)$ in $\overline{C}$, for $t \geq 0$, can be joined to the hyperplane $\tau = 0$ through a line of the form $\tau + r = t + |x| $ which is totally contained in $\overline{C}$. Since the integration of the gradient is over lines $\tau + r = constant$ which is in the exterior region and therefore the estimates on the gradient hold true, we could then integrate as in \cite{G4}. Then, we note that under the bootstrap assumption, and therefore under the a priori estimates in Lemmas \ref{apriordecayestimatesfrombootstrapassumption} and \ref{asymptoticbehaviouratteqzeroforallfields}, we have that for any $q_0 \in \R$\;, there exists a point $(t, r=0)$ such that $N$ (defined in Definition \ref{definitionofNandofNtruncatedbtweentwots}) whose tip is $(t, r=0)$ is contained in the region $\{(t, x) \;|\; q:= r-t \leq q_0 \}$\;, and also for any $N$ there exists a $q_0$ such that the exterior region contains $\{ q \geq q_0 \}$\;. Thus, working in the exterior region $\overline{C}$ is equivalent, for our purposes, to looking at the region $\{ q \geq q_0 \}$\;.

\end{proof}

\section{Improved decay estimates for good derivates}

From the a priori estimates generated from the bootstrap assumption combined with a weighted Sobolev inequality, one can derive improved decay estimates for “good” derivates, which are derivatives tangential to the outgoing light cone at a point in the space-time. To explain, let us first start by establishing a relation between the space-time derivatives of a function on one hand, and the Lorentz boosts, rotations and scaling vector field applied to a function on the other hand. Since we already established an estimate on the latter applied to the Yang-Mills potential $A$ in wave coordinates and to $h^1$, we could then see how this could translate into estimate on space-time derivatives of these, and extract better decay estimates for certain good derivatives.

\subsection{Estimating derivatives in terms of the Minkowski vector fields}\

\begin{lemma}
We have for $i \in \{1, \ldots, n \}$\;, 
\bea
  \pa_{t}  &=&  \frac{ t S - x^i  Z_{0i} }{   t^2- r^2  } \;,\\
   \pa_{i}  &=&\frac{  - x_i  S + t Z_{0i}  - x^j Z_{ij}   }{   t^2- r^2  } \;, \\
\pa_r &=& \frac{ (- r S + t \frac{x^{i}}{r}  Z_{0i}   )  }{   t^2- r^2  } \; ,
\eea
and we also have
\bea
 \pa_{i}  &=&   \frac{  - x_i  S    }{   t^2- r^2  }   + \frac{x_i x^j Z_{0j}   }{ t(  t^2- r^2 )  }   +  \frac{ Z_{0i}   }{ t  }\; .
\eea

\end{lemma}

\begin{proof}

Computing
\beaa
\notag
 t S - x^i  Z_{0i}  &=&t ( t \pa_{t} + x^{i} \pa_{i} )  - x^i   (  x_{i} \pa_{t} + t \pa_{i} ) \; .\\
 \eeaa

We have, for a summation on $i$ over spatial indices,
\beaa
x^i x_{i} = m_{\a i} x^i x^{\a} = m_{i i} x^i x^i = \sum_{i=1}^n (x^i)^2 = r^2 \;,
\eeaa
thus
 \beaa
  \notag
 t S - x^i  Z_{0i}  &=& ( t^2 - r^2) \pa_{t} \;.
\eeaa
Hence,
\bea
  \pa_{t}  &=&  \frac{ t S - x^i  Z_{0i} }{   t^2- r^2  } \;.
\eea

We compute
\beaa
\notag
 - x_i  S + t Z_{0i}  - x^j Z_{ij}  &=&  - x_i ( t \pa_{t} + x^{j} \pa_{j} )  + t (  x_{i} \pa_{t} + t \pa_{i} ) - x^j (  x_{j} \pa_{i} -  x_{i}  \pa_{j} ) \\
 \notag
 &=&     t^2 \pa_{i}  - x^j   x_{j} \pa_{i} =   (  t^2 - r^2 ) \pa_{i} \; .
\eeaa
Thus,
\bea
 \pa_{i}  &=&\frac{  - x_i  Z_0 + t Z_{0i}  - x^j Z_{ij}   }{   t^2- r^2  } \; .
\eea

We have
\bea
\pa_r = \frac{x^{i}}{r} \pa_{i}\; .
\eea
We compute
\beaa
\notag
\pa_r &=& \frac{x^{i}}{r} \pa_{i} = \frac{x^{i}}{r} \cdot \frac{ (- x_i  S + t Z_{0i}  - x^j Z_{ij} )  }{( t- r) (t +r) } \\
&=& \frac{ (- r S + t \frac{x^{i}}{r}  Z_{0i}  - \frac{x^{i}}{r}  x^j Z_{ij} )  }{   t^2- r^2  }\; .
\eeaa
We have
\beaa
 x^i  x^j Z_{ij} = - x^j x^i  Z_{ji} \;,
\eeaa
hence
\beaa
 x^i  x^j Z_{ij} =0 \; .
\eeaa
As a result
\bea
\pa_r &=& \frac{ (- r S + t \frac{x^{i}}{r}  Z_{0i}   )  }{   t^2- r^2  } \;.
\eea

We showed that
\beaa
 \pa_{i}  &=&\frac{  - x_i  Z_0    }{   t^2- r^2  } + \frac{   t Z_{0i}   }{   t^2- r^2  }  + \frac{- x^j Z_{ij}   }{   t^2- r^2  } \; .
\eeaa

We compute
\beaa
t Z_{0 i} = t ( x_{i} \pa_{0} - x_{0} \pa_{i} ) = t ( x_{i} \pa_{t} + t \pa_{i} )  = t  x_{i} \pa_{t} + t^2 \pa_{i} \;,
\eeaa
and
\beaa
- x^j Z_{ij} = - x^j  ( x_{j} \pa_{i} - x_{i} \pa_{j} ) = - x^j  x_{j} \pa_{i}  + x^j x_{i} \pa_{j}  = - r^2 \pa_{i}  + x_{i} x^j  \pa_{j} \;,
\eeaa
we obtain
\beaa
t Z_{0 i}  - x^j Z_{ij} &=& t  x_{i} \pa_{t} + t^2 \pa_{i}  - r^2 \pa_{i}  + x_{i} x^j  \pa_{j} \\
&=& t  x_{i} \pa_{t} + ( t^2  - r^2 ) \pa_{i}  + x_{i} x^j  \pa_{j} \;.
\eeaa
On the other hand,
\beaa
\frac{x_i x^j Z_{0j}  + (t^2 - r^2) Z_{0i} }{t} &=& \frac{ x_i x^j ( x_{j} \pa_{t} - x_{0} \pa_{j} )  + (t^2 - r^2)  ( x_{i} \pa_{t} - x_{0} \pa_{i} ) }{t} \\
&=& \frac{ x_i r^2 \pa_{t} - x_{0}x_i x^j  \pa_{j}  + t^2  x_{i} \pa_{t} -  t^2  x_{0} \pa_{i}   - r^2  x_{i} \pa_{t} + r^2 x_{0} \pa_{i}  }{t} \\
&=& \frac{ - x_{0}x_i x^j  \pa_{j}  + t^2  x_{i} \pa_{t} -  ( t^2 -r^2 ) x_{0} \pa_{i}  }{t} \\
&=& \frac{  t x_i x^j  \pa_{j}  + t^2  x_{i} \pa_{t} +  ( t^2 -r^2 ) t \pa_{i}  }{t} \\
&=&   x_i x^j  \pa_{j}  + t  x_{i} \pa_{t} +  ( t^2 -r^2 )  \pa_{i} \;  .\\
\eeaa
Consequently,
\bea
t Z_{0 i}  - x^j Z_{ij} &=& \frac{x_i x^j Z_{0j}  + (t^2 - r^2) Z_{0i} }{t}\;,
\eea
and injecting this in the expression of $\pa_i$, we get
\bea
 \pa_{i}   &=&  \frac{  - x_i  S    }{   t^2- r^2  }   + \frac{x_i x^j Z_{0j}   }{ t(  t^2- r^2 )  }   +  \frac{ Z_{0i}   }{ t  } \;.
\eea

\end{proof}

\begin{definition}\label{defrestrictedderivativesintermsofZ}

We define for $i \in \{1, \ldots, n\} $\;,
\bea
\notag
\rpa_{i} &:=& \pa_i - E ( \pa_{i} ,  \pa_{r}  ) \cdot \pa_{r} \\
\notag
&=& \pa_i - E ( \pa_{i} ,  \frac{x^j}{r} \pa_{j}  ) \cdot  \pa_{r} \\
&=& \pa_i - \frac{x_i}{r} \pa_{r} \; ,
\eea
and
\bea
\rpa_{0} &:=& \frac{ \pa_t +  \pa_r }{2} \;.
\eea
\end{definition}

\begin{lemma}\label{restrictedderivativesintermsofZ}

For restricted partial derivatives defines as in Definition \ref{defrestrictedderivativesintermsofZ}, we have for $i \in \{1, \ldots, n\} $\;,
\bea
\rpa_{i} &=& \frac{  - \frac{x_i}{r} \frac{x^j}{r}   Z_{0j}  +  Z_{0i}    }{  t  } \; , \\
\rpa_{0}  &=&  \frac{S + \frac{x^i}{r}  Z_{0i}   }{  2( t + r ) } \;, 
\eea
and we also have
 \bea
\rpa_i  &=&   \frac{x^j}{r^2}   Z_{ij} \; .
 \eea

\end{lemma}

\begin{proof}

We compute
\beaa
\rpa_{i} &=& \pa_i - \frac{x_i}{r} \pa_{r}\\
&=&  \frac{  - x_i  S   }{   t^2- r^2  }   + \frac{x_i x^j Z_{0j}   }{ t(  t^2- r^2 )  }   +  \frac{ Z_{0i}   }{ t  }  - \frac{x_i}{r} \frac{ (- r S + t \frac{x^{j}}{r}  Z_{0j}   )  }{   t^2- r^2  } \\
&=& \frac{  - t x_i   S + x_i x^j Z_{0j}  + (  t^2- r^2  )  Z_{0i}  + t x_i  S  -  t^2  \frac{x_i}{r} \frac{x^{j}}{r}  Z_{0j}    }{  t(  t^2- r^2  ) }  \\
&=& \frac{   x_i x^j Z_{0j}  + (  t^2- r^2  )  Z_{0i}    -  t^2 \frac{x_i}{r} \frac{x^{j}}{r}  Z_{0j}   }{  t(  t^2- r^2  ) }  \\
&=& \frac{   x_i x^j Z_{0j}  -  \frac{t^2}{r^2} x_i x^{j}  Z_{0j}  + (  t^2- r^2  )  Z_{0i}      }{  t(  t^2- r^2  ) } \\
&=& \frac{     - (t^2 - r^2) \frac{x_i x^j}{r^2} Z_{0j}  +  (  t^2- r^2  )  Z_{0i}    }{   t(  t^2- r^2  ) } \\
&=& \frac{     -  \frac{x_i x^j}{r^2} Z_{0j}  +  Z_{0i}    }{   t } .
\eeaa
Consequently,
\bea
\rpa_{i} &=& \frac{  - \frac{x_i}{r} \frac{x^j}{r}   Z_{0j}  +  Z_{0i}    }{  t  }  
\eea

Also, computing
\beaa
   \frac{x^j}{r^2}   Z_{ij} &=&     \frac{x^j}{r^2}   ( x_j \pa_i  -  x_i \pa_j )  =   \frac{x^j}{r^2}    x_j \pa_i  -  \frac{x^j}{r^2} x_i \pa_j \\
   &=&   \pa_i  -  \frac{x_i}{r} \frac{x^j}{r} \pa_j =  \pa_i  -  \frac{x_i}{r}  \pa_r \, ,
 \eeaa
we get,
 \bea
\rpa_i  &=&   \frac{x^j}{r^2}   Z_{ij} \, .
 \eea
  
We compute
 \beaa
\rpa_{0} &=& \frac{ \pa_t +  \pa_r }{2} \\
&=& \frac{ t S - x^i  Z_{0i} }{  2( t^2- r^2 ) } + \frac{ (- r S + t \frac{x^{i}}{r}  Z_{0i}   )  }{ 2(  t^2- r^2 ) } \\
&=& \frac{( t -r ) S + ( t \frac{x^{i}}{r}  - x^i ) Z_{0i}   }{  2( t^2- r^2 ) } = \frac{( t -r ) S + \frac{x^i}{r} ( t- r  ) Z_{0i}   }{  2( t^2- r^2 ) }\\
&=&  \frac{ S + \frac{x^i}{r}  Z_{0i}   }{  2( t + r ) } ,
\eeaa
and we obtain
\bea
\rpa_{0}  &=&  \frac{Z_0 + \frac{x^i}{r}  Z_{0i}   }{  2( t + r ) } .
\eea

\end{proof}

We will now point out the following lemma, which will be used through out the manuscript, which shows how decay rates in different variables, $t$, $r$ and $q$, relate together depending on the different regions $q \leq 0$ and $q \geq 0$.
\begin{lemma}\label{decayfactorsforregions}
We consider in this lemma the region $t \geq 0$.

For $q \leq 0$, we have
\bea
 \frac{1}{1+2t} &\leq&  \frac{1}{1 + t + |q|  } \; ,
  \eea
  and
  \bea
  \frac{1}{1 + r + 2 |q|  } &= & \frac{1}{1 + t + |q|  }\; ,\\
 \frac{1}{1 + t } &= & \frac{1}{1 + r + |q|  }\;  .
  \eea
For $q \geq 0$, we have
\bea
\frac{1}{1 + r   }  &= & \frac{1}{1 + t + |q|  }\;  .
\eea

Also, for all $q$, we have
\bea
\frac{1}{1 + t + r   }  &\les & \frac{1}{1 + t + |q|  }\;  .
\eea

\end{lemma}

\begin{proof}

For $q \leq 0$, then $ |q| = - q = t-r$. We then compute
\beaa
1+t+ |q|= 1+t +  t-r = 1+2t - r \leq 1+2t ,
\eeaa
and hence,
\bea
 \frac{1}{1+2t} &\leq&  \frac{1}{1 + t + |q|  } \; .
  \eea
 We also compute
 \bea
1+t+ |q|= 1+t +  t-r = 1+ 2t -r = 1+ r+ 2t -2r = 1 + r + 2 |q| ,
\eea
 and 
  \bea
1+r+ |q|= 1+r +  t-r = 1+t .
\eea
The last equality also gives
\beaa
1+t + r =  1+2r+ |q| ,
\eeaa
which leads to 
\beaa
1+t + r \geq   1+r+ |q|  &\geq&  \frac{1}{2} ( 2+2r+ 2|q|  ) \geq  \frac{1}{2} ( 1+r+ 2|q|  ) ,
\eeaa
and hence
\beaa
  \frac{1}{ 2( 1+t + r )  }  &\leq&  \frac{1}{ 1+r+ 2|q|  }  ,
\eeaa
and consequently
\beaa
\frac{1}{ 1+t + r  }  &\les&    \frac{1}{ 1+r+ 2|q|  } .
\eeaa
But we proved that in the region $q<0$, we have $1+r+ 2|q|  = 1 + t + |q|$. Thus, we get for $q<0$,
\bea
\frac{1}{ 1+t + r  }  &\les&    \frac{1}{ 1+t + |q|  } .
\eea

For $q \geq 0$, we have $ |q| =  q = r-t$. We compute
 \bea
1+t+ |q|= 1+t + r- t-=  1 + r ,
\eea
which also gives
 \bea
 1 + t +r = 1+2t+ |q| \geq 1+t+ |q|  ,
\eea
and therefore, for $q \geq 0$, we get
\bea
\frac{1}{ 1+t + r  }  &\leq&    \frac{1}{ 1+t+ |q|  } .
\eea
\end{proof}

\begin{lemma}\label{betterdecayfortangentialderivatives}

We define
\bea
 |\rpa f |^2 :=     |\rpa_0 f |^2 + \sum_{i=1}^n   |\rpa_i f |^2   .
  \eea
Then, we have the following inequality for all $t \geq 0$ and for all $q \in \R$,
\bea
(1 + t + |q| ) \cdot |\rpa f | \les \sum_{|I| = 1} |Z^I f |  .
\eea
\end{lemma}

\begin{proof}

We consider first the region where either $t \geq 1 $ or $r \geq 1$. In that region, we have $2t  \geq 1 + t$ or $2r \geq 1 + r$ and therefore, either 
\beaa
\frac{1}{t   } \les  \frac{1}{ 1+t  } \; ,
\eeaa
or
\beaa
\frac{1}{r   } \les  \frac{1}{ 1+r   } \;.
\eeaa
We also have in that region  $t+r \geq 1$ and therefore $2t+2r \geq 1 + t +r $ and hence,
\beaa
\frac{1}{ t + r   } \les  \frac{1}{ 1+t+ r   } \; .
\eeaa
We now evaluate $  |\rpa_0 f |^2 + \sum_{i=1}^n   |\rpa_i f |^2   $ . We have
\beaa
|  \rpa_{0} f|   &=&|  \frac{S  f+ \frac{x^i}{r}  Z_{0i} f   }{  2( t + r ) } | \les  \frac{ | S f | + | \frac{x^i}{r}  Z_{0i} f  | }{  1 + t + r  } .
\eeaa
Given that
\beaa
| \frac{x^i}{r} |  =  | \frac{x^j  }{  \sqrt{(x^1)^2 + (x^2)^2 +(x^3)^2 } }  |  \leq 1 ,
\eeaa
we obtain
\bea
|  \rpa_{0} f|   &\les & \sum_{|I| = 1}  \frac{ |Z^I f | }{  1 + t + r  } \les  \sum_{|I| = 1}  \frac{ |Z^I f | }{  1 + t + | q|   }  .
\eea

To evaluate $  \sum_{i=1}^n  |\rpa_i f |^2 $, we look at the cases for the sign of $q$:

\textbf{Case $q = r-t \leq 0$:}\\

In that case, being in the region $t \geq 1 $ or $r \geq 1$, would impose for $q \leq 0$ that in fact $t \geq 1 $. We have for all $i \in \{1, \ldots, n \}$, 
\beaa
| \rpa_{i} f | &=& |  \frac{  - \frac{x_i}{r} \frac{x^j}{r}   Z_{0j} f +  Z_{0i}  f  }{  t  }  |   \\
& \les&    \frac{ | - \frac{x_i}{r} \frac{x^j}{r}   Z_{0j} f +  Z_{0i}  f | }{  1 +t  }  \\
&& \text{(considering that $t \geq 1 $ when we look in the region $t \geq 1 $ or $r \geq 1$, for $q \leq 0$) } \\
& \les&    \frac{ | - \frac{x_i}{r} \frac{x^j}{r}   Z_{0j} f |  + |  Z_{0i}  f | }{  1 +t +  | q| }  \\
&& \text{(considering that we are in the region $q \leq 0 $) } \\
& \les&  \sum_{|I| = 1}  \frac{ |Z^I f | }{  1 + t + | q|   }  .
\eeaa

\textbf{Case $q = r-t > 0$:}\\

 Being in the region $t \geq 1 $ or $r \geq 1$, would impose for $q  > 0$ that in fact $r \geq 1 $. For all $i \in \{1, \ldots, n \}$
\beaa
| \rpa_{i} f | &=& | \frac{x^j}{r^2}   Z_{ij} f | = \frac{ | \frac{x^j}{r}   Z_{ij} f | }{r} \leq  \frac{ | \frac{x^j}{r}   Z_{ij} f | }{1 +r}  \\
&& \text{(where we used the fact that $r \geq 1$ in the region $t \geq 1 $ or $r \geq 1$, for $q > 0$) } \\
&\les&  \frac{ | \frac{x^j}{r}   Z_{ij} f | }{1 +t +  | q|  }  \\
&& \text{(where we used the fact that $q \geq 0$) } \\
& \les&  \sum_{|I| = 1}  \frac{ |Z^I f | }{  1 + t + | q|   } .
\eeaa

Thus, in the region $t \geq 1 $ or $r \geq 1$, we have
\beaa
 |\rpa f | &=& \sqrt{ |\rpa_0 f |^2 + \sum_{i=1}^3   |\rpa_i f |^2   } \leq |\rpa_0 f  | + \sum_{i=1}^3   |\rpa_i f | \\
&& \text{(using $\sqrt{ x +y } \leq \sqrt{x} + \sqrt{y}$) } \\
& \les&  \sum_{|I| = 1}  \frac{ |Z^I f | }{  1 + t + | q|   }  .
\eeaa

Finally, we consider the region $0 \leq t \leq 1 $ and $r \leq 1$. In this region, we have
\beaa
1 + t +  | q| = 1+ t +| r-t| \leq 1 + t +  r +| t |  \leq 1 + 2t + r \leq 1 +2 + 1 \leq 4
\eeaa 
 and therefore
  \beaa
1  \les  \frac{1}{ 1+t+ | q|   } .
\eeaa
Now, we evaluate  $ |\rpa_r f  |$. We recall that by definition, we have
 \beaa
| \rpa_{0} f | &=& | \frac{ \pa_t +  \pa_r }{2}  f | \les |  \pa_t f |  + |   \pa_r f | \\
&\les&  \sum_{|I| = 1}   |Z^I f |\\
&& \text{(since $\pa_\a \in \cal Z $ for $\a \in \{0, 1,\ldots, n \}$ )}  .
\eeaa

For $i \in \{1, \ldots, n\} $, by definition of the derivatives tangential to the light cone, we have
\beaa
\notag
| \rpa_{i} f |&=& | \pa_i f - \frac{x^i}{r} \pa_{r}  f | \leq  | \pa_i f |  + |  \frac{x^i}{r} \pa_{r}  f | \\
&\leq&  | \pa_i f |  + |   \pa_{r}  f | .
\eeaa
Since 
\bea
\pa_r f= \frac{x^{i}}{r} \pa_{i} f ,
\eea
we get
\beaa
| \rpa_{i} f  |  &\leq&  | \pa_i f |  + |   \frac{x^{i}}{r} \pa_{i}   f | \\
&\leq&  | \pa_i f |  + |  \pa_{i}   f | \\
&\les&  \sum_{|I| = 1}   |Z^I f |\\
&& \text{(because $\pa_\a \in \cal Z $ for $\a \in \{0, 1, \ldots, n \}$ }  .
\eeaa
This proves that 
\beaa
| \rpa f  | &\les&  \sum_{|I| = 1}   |Z^I f | \les   \sum_{|I| = 1}  \frac{ |Z^I f | }{  1 + t + | q|   } \\
&& \text{(considering that we are in the region $0 \leq t \leq 1 $ and $r \leq 1$) }.
\eeaa
\end{proof}

\begin{lemma}\label{decayrateforfullderivativeintermofZ}

We have the following inequality for all $t \geq 0$ and for all $q$, 
\bea
(1 +  |q| ) \cdot |\pa f | &\les& \sum_{|I| = 1} |Z^I f | \; .
\eea
\end{lemma}

\begin{proof}
First, we consider points in the space-time such that $q\neq 0$ and $t \geq 0$. We have
\beaa
 |  \pa_{t} f   | &=&  | \frac{ t S f - x^i  Z_{0i} f }{   t^2- r^2  }  | \leq  \frac{  | t S f  | +  |- x^i  Z_{0i} f  | }{   | ( t + r ) (t-r)  |  } \\
& \leq&  \frac{ ( t + r ) (   \sum_{|I| = 1} |Z^I f | )  }{   ( t + r ) | t-r  |  } \\
  &\leq&  \frac{    \sum_{|I| = 1} |Z^I f |   }{    | q  |  } .
  \eeaa
We also have for  $q\neq 0 $ and $t \geq 0$, for all $i \in \{1, \ldots, n \}$,
  \beaa
    | \pa_{i} f |  &=& \frac{   | - x_i  S f  + t Z_{0i} f  - x^j Z_{ij}   f  |}{   | t^2- r^2  | } \\
   & \leq&  \frac{ ( t + r ) (   \sum_{|I| = 1} |Z^I f | )  }{   ( t + r ) | t-r  |  } \\
  &\leq&  \frac{    \sum_{|I| = 1} |Z^I f |   }{    | q  |  } .
   \eeaa
Hence, for $q\neq 0$ and $t \geq 0$,
\beaa
 |  \pa f   | &\leq&   |  \pa_{t} f   |  +  \sum_{i = 1}^3| \pa_{i} f | \\ 
  &\leq&  \frac{    \sum_{|I| = 1} |Z^I f |   }{    | q  |  } .
  \eeaa
Thus, for $t \geq 0$,
\beaa
  |q| \cdot  |  \pa f  |    &\leq&    \sum_{|I| = 1} |Z^I f |    \\
 && \text{(which holds true also for $q =0$)}.
  \eeaa
  
 Now, we would like to bound $  |  \pa f  | $ on the outgoing light cone prescribed by $q=0$. We note that
\beaa
 |  \pa f  |  &\leq&  |  \pa_{t} f   |  +  \sum_{i = 1}^3| \pa_{i} f |  \\
&\les&  \sum_{|I| = 1}   |Z^I f |\\
&& \text{(because $\pa_\a \in \cal Z $ for $\a \in \{0, 1, \ldots, n \}$)}  .
\eeaa
As a result, for all $t\geq 0$, we have
\beaa
(1 +  |q| ) \cdot |\pa f | \les \sum_{|I| = 1} |Z^I f |  .
\eeaa

 \end{proof}

\subsection{Estimating good derivatives of the fields}\

Thanks to the estimate
\bea \label{goodderivative}
(1 + t + |q| ) \cdot |\rpa f | + (1 +  |q| ) \cdot  |\pa f | &\les& \sum_{|I| = 1} |Z^I f | \, ,
\eea
that we have just shown, in Lemmas \ref{betterdecayfortangentialderivatives} and \ref{decayrateforfullderivativeintermofZ}, we can establish for $ |\rpa f |$ a decay estimate in $t$  even along the outgoing null-come prescribed by $q=0$, however, we do not have such a decay estimate for $|\pa f |$. Thus, we call $\rpa f$ a “good” derivative. Yet, this is provided that one can control $\sum_{|I| = 1} |Z^I f |$. In the case of the Yang-Mills potential $A$ and $h^1$ in wave coordinates, we have already established a control on their respective norms, which are $ |\Lie_{Z^I} A | $ and $|\Lie_{Z^I} h^1 |$, by integrating along $\tau+r = constant$ the decay estimates generated from a direct application of the weighted Sobolev inequality combined with the bootstrap assumption on the growth of the defined energy.

\begin{lemma} We have for $I$ defined as in \ref{definitionofZI},
\bea
\notag
(1 + t + |q| ) \cdot |\rpa ( \Lie_{Z^I} A ) | +(1 + t + |q| ) \cdot  |\rpa ( \Lie_{Z^I} h^1 ) | \les \sum_{|J| \leq |I| + 1} \big( | \Lie_{Z^J} A | + |  \Lie_{Z^J} h^1 |  \big)   . \\
\eea
\end{lemma}

\begin{proof}
By applying the above estimate, \eqref{goodderivative}, for $f$ being taken to be each component of $\Lie_{Z^I} A$ in wave coordinates using the norm on the Lie algebra, and on each component of $\Lie_{Z^I} h^1$ using the scalar norm, we obtain
that for all $\si, \mu, \nu \in (t, x^1, \ldots, x^n)$,

\beaa
(1 + t + |q| )  \cdot |\rpa \Lie_{Z^I} A_\si | +(1 + t + |q| ) \cdot |\rpa \Lie_{Z^I} h^1_{\mu\nu} | \les \sum_{|J| = 1} \big( | Z^J \Lie_{Z^I} A_\si | + | Z^J \Lie_{Z^I} h^1_{\mu\nu} |  \big)   .
\eeaa

Using again the fact that a commutation of two vector fields in $\cal Z$ is a combination of vector fields in $\cal Z$, and that a commutation of a vector field in $\cal Z$ and of a $\pa_\mu$ gives a linear combination of vectors of the form $\pa_\mu$, we get that for all $\si, \mu, \nu \in (t, x^1, x^2, x^3)$,

\beaa
\notag
(1 + t + |q| )  \cdot |\rpa \Lie_{Z^I} A_\si | +(1 + t + |q| )  \cdot |\rpa \Lie_{Z^I} h^1_{\mu\nu} | \les \sum_{|J| \leq |I| + 1} \big( | \Lie_{Z^J} A | + |  \Lie_{Z^J} h^1 |  \big)   ,
\eeaa
and hence,
\bea
\notag
(1 + t + |q| )  \cdot |\rpa ( \Lie_{Z^I} A ) | +(1 + t + |q| )  \cdot |\rpa ( \Lie_{Z^I} h^1 ) | \les \sum_{|J| \leq |I| + 1} \big( | \Lie_{Z^J} A | + |  \Lie_{Z^J} h^1 |  \big)   . \\
\eea
\end{proof}

\begin{lemma}\label{estimategoodderivatives}
Under the bootstrap assumption \eqref{bootstrap}, taken for $N =  |I| + 3 $, we have for all $ I $, in the exterior region $\overline{C} \subset \{ q \geq q_0 \} $,
 \bea
 \notag
 |\rpa ( \Lie_{Z^I}  A ) | + |\rpa ( \Lie_{Z^I}  h^1 ) |     &\leq& \begin{cases} c (\gamma) \cdot C ( |I| +1 ) \cdot E ( |I| +  3) \cdot  \frac{\eps}{ (1+ t + | q | )^{2-\delta }  (1+| q |   )^{\gamma}},\quad\text{when }\quad q>0,\\
  \notag
      C (   |I| +1 ) \cdot E (   |I| +  3) \cdot  \frac{\eps}{ (1+ t + | q | )^{2-\delta }  } (1+| q |   )^{\frac{1}{2} }  \,\quad\text{when }\quad q<0 . \end{cases} \\
       \notag
      &\leq& C(q_0) \cdot  c (\gamma) \cdot C ( |I| +1 ) \cdot E ( |I| +  3) \cdot  \frac{\eps}{ (1+ t + | q | )^{2-\delta }  (1+| q |   )^{\gamma}} \; .\\
      \eea
      
\end{lemma}

\begin{proof}

We proved that for $q<0$, for all $| J|$, we have  
\beaa
\notag
\mid \Lie_{Z^J} A (t,  | x | \cdot \Om)  \mid + \mid \Lie_{Z^J} h^1 (t,  | x | \cdot \Om) \mid &\les&  C ( |J| ) \cdot E ( |J| +  2) \cdot  \frac{\eps}{ (1+ t + | q | )^{1-\delta }  } (1+| q |   )^{\frac{1}{2} }   , 
\eeaa
and in particular, this gives that for all $|J| \leq |I| + 1$,
\beaa
\notag
\mid \Lie_{Z^J} A (t,  | x | \cdot \Om)  \mid + \mid \Lie_{Z^J} h^1 (t,  | x | \cdot \Om) \mid &\les&  C (  |I| + 1  ) \cdot E (   |I|  +  3) \cdot  \frac{\eps}{ (1+ t + | q | )^{1-\delta }  } (1+| q |   )^{\frac{1}{2} }   , 
\eeaa
and consequently, for all $|J| \leq |I| + 1$,
\beaa
\notag
\sum_{|J| \leq |I| + 1}  \big( |\Lie_{Z^J}  A| + |\Lie_{Z^J}  h^1 | \big)   &\les&  C (   |I| +1 ) \cdot E (   |I| +  3) \cdot  \frac{\eps}{ (1+ t + | q | )^{1-\delta }  } \cdot (1+| q |   )^{\frac{1}{2} }   .
\eeaa

This gives, using Lemma \ref{betterdecayfortangentialderivatives}, that for $q<0$, for all $| I|$,
\bea
\notag
 |\rpa ( \Lie_{Z^I}  A) | + |\rpa ( \Lie_{Z^I} h^1) | &\les&  C (   |I| +1 ) \cdot E (   |I| +  3) \cdot  \frac{\eps}{ (1+ t + | q | )^{2-\delta }  } \cdot (1+| q |   )^{\frac{1}{2} }   .\\
\eea

For $q \geq 0$, we proved that
\bea 
\notag
\mid \Lie_{Z^J}  A (t,  | x | \cdot \Om)  \mid + \mid \Lie_{Z^J}  h^1 (t,  | x | \cdot \Om) \mid &\les& c (\gamma)\cdot C ( |J| ) \cdot E ( |J| +  2) \cdot  \frac{\eps}{ (1+ t + | q | )^{1-\delta }\cdot   (1+| q |   )^{\gamma}}   , 
\eea
and hence,
\beaa
\notag
\sum_{|J| \leq |I| + 1}  \big( | \Lie_{Z^J}  A| + | \Lie_{Z^J}  h^1 | \big)   &\les&  c (\gamma)\cdot C ( |I| +1 ) \cdot E ( |I| +3) \cdot  \frac{\eps}{ (1+ t + | q | )^{1-\delta } \cdot  (1+| q |   )^{\gamma}}  .
\eeaa

Using Lemma \ref{betterdecayfortangentialderivatives}, this gives that for $q \geq 0$, for all $| I|$, we have  
\bea
\notag
 |\rpa ( \Lie_{Z^I}  A ) | + |\rpa ( \Lie_{Z^I}  h^1 ) |    &\les&  c (\gamma)\cdot C ( |I| +1 ) \cdot E ( |I| + 3) \cdot  \frac{\eps}{ (1+ t + | q | )^{2-\delta } \cdot  (1+| q |   )^{\gamma}}  . \\
\eea
\end{proof}

 \section{Estimating the spherically symmetric part of the metric}
   
  \begin{definition}
Now, following a guess by Lindblad-Rodnianski in \cite{LR10}, we are going to define for $t > 0$, the following quantity $h^0$, and then define $h^1$ that is a part of $h$ from which we are going to subtract $h^0$ and then we are going to show in this paper that the wave equation for $h^1$, derived from coupled the Einstein-Yang-Mills system, is well-posed:
\beaa
h^1 := h - h^0
\eeaa 
with $h^0$ defined in wave coordinates, i.e. for $\mu,\; \nu \in \{0, 1, \ldots, n\}$, for $t > 0 $, as being 
\bea
h^0_{\mu\nu} := \chi(r/t)\chi(r)\frac{M}{r}\de_{\mu\nu} \; ,
\eea
where $\chi(s)$ is smooth function of $s$ such that
\bea
\chi(s) = \begin{cases} 1 \qquad \text{when} \qquad s\geq 3/4 ,\\
0  \qquad \text{when} \qquad s\leq 1/2 ,\end{cases}
\eea
and for $t=0$,
\bea
h^0_{\mu\nu} ( t= 0) := \chi(r)\frac{M}{r}\de_{\mu\nu} \, ,
\eea
where $M > 0$, and where $\de_{\mu\nu}$ is the Kronecker symbol.
\end{definition}

\begin{lemma}\label{derivativesofsphericalsymmetricpart}
We have
 \beaa
 \notag
|\pa   h^0 (t,x)  |   &\les& \frac{M}{(1+t+|q|)^2 }.
      \eeaa
      
\end{lemma}

\begin{proof}

We are now interested in getting estimates on quantities involving $h^0$. 
We fix $ \mu, \nu \in (t, x^1, \ldots, x^n)$. Computing for $t >0$,
\beaa
 \pa_{t} h^0_{\mu\nu} &=& \pa_{t} \big( \chi(r/t)\chi(r)\frac{M}{r}\de_{\mu\nu} \big) \\
 &=& \chi(r)\frac{M}{r}\de_{\mu\nu} \cdot (-\frac{r}{t^2} ) \chi^\prime (r/t)  = - \frac{r }{t^2} \frac{M}{r} \chi(r) \chi^\prime (r/t) \de_{\mu\nu} \\
&=& - \frac{ M }{t^2}  \chi(r) \chi^\prime (r/t) \de_{\mu\nu} \, .
\eeaa
We have
\beaa
\pa_i r =  \frac{x_i}{r } \,  ,
\eeaa
and thus
\beaa
 \pa_{i} h^0_{\mu\nu} &=& \pa_{i} \big( \chi(r/t)\chi(r)\frac{M}{r}\de_{\mu\nu} \big) \\
&=&   \chi^\prime (r/t) \cdot ( \frac{x_i}{r t} ) \cdot \chi(r)\frac{M}{r}\de_{\mu\nu}  +    \chi (r/t)\chi^\prime (r) \cdot ( \frac{x_i}{r } ) \cdot  \frac{M}{r}\de_{\mu\nu} +    \chi (r/t)\chi(r)\frac{M}{r^2} \cdot ( \frac{- x_i}{r} ) \cdot  \de_{\mu\nu}  \\
&=&   \frac{x_i M}{r^2 t }  \chi^\prime (r/t)  \chi(r) \de_{\mu\nu}  +  \frac{x_i M }{r^2}  \chi (r/t)\chi^\prime (r)  \de_{\mu\nu} +    \frac{- x_i M }{r^3}    \chi (r/t)\chi(r)   \de_{\mu\nu}  . \\
\eeaa

We notice that the term that is not compactly supported behaves like $\frac{M}{r^2}$ and furthermore, the terms which are compactly supported, due to the terms $\chi^\prime (r/t) $ and $\chi^\prime (r) $, have a support either in the region $\frac{1}{2} < \frac{r}{t} < \frac{3}{4}$ or in the region $\frac{1}{2} < r < \frac{3}{4}$, respectively.

\textbf{Terms multiplied with $\chi^\prime (r/t)$}\\

The compactly supported term multiplied with $\chi^\prime (r/t)$ is in the following region, for $t > 0$,
\beaa
\frac{t}{2} < r < \frac{3 }{4} t ,
\eeaa
which in particular means that 
\beaa
\frac{4}{3} < \frac{t}{r} < 2 .
\eeaa
Thus, in that region
\bea
\frac{M}{r} \leq \frac{2M}{t} ,
\eea
and 
\bea
\frac{M}{t} \leq \frac{3M}{4r} .
\eea
Using the fact that 
\beaa
\frac{x_i}{r} \leq 1 ,
\eeaa
we also get
\bea
\frac{x_i}{t} \leq \frac{r}{t} \leq \frac{3}{4} .
\eea

\textbf{Terms multiplied with $\chi (r/t)$}\\

For the terms with $\chi (r/t)$, the support is in the following region, for $t > 0$,
\beaa
 \frac{r}{t} > \frac{1}{2} ,
\eeaa
which means
\bea
 \frac{1}{r} <\frac{2}{t} .
\eea

\textbf{Terms multiplied with $\chi^\prime (r)$}\\

For terms with $\chi^\prime (r)$, we have $\frac{1}{2} < r < \frac{3}{4} $, 
and thus
\beaa
\frac{4}{3} < \frac{1}{r} < 2 \; ,
\eeaa
and hence
\bea
\frac{1}{r} < \frac{3}{4} \frac{1}{r^2} \;.
\eea

\textbf{Terms multiplied with either $\chi (r)$ or $\chi^\prime (r)$}\\

Now, due to the presence of either the term $\chi (r)$ or the term $\chi^\prime (r)$, we have
\beaa
\frac{1}{2} < r ,
\eeaa
and thus $1<2r$, which gives $1+r < 3r$ and consequently, we always have
\bea
\frac{1}{r} < \frac{3}{1+r} ,
\eea
and similarly
\bea
\frac{1}{r^2} \les \frac{1}{1+r^2} .
\eea
In particular, estimates of the type $\frac{1}{r} \les \frac{1}{t}$ which translate as $t \les r$, lead to $1+t \les 1+r \les  3r$ and therefore, they lead to following type estimate
$\frac{1}{r} \les \frac{1}{1+t}$. In other words, we have
\bea
\frac{1}{r} \les \frac{1}{t} \Rightarrow \frac{1}{r} \les  \frac{1}{1+t} \;,
\eea
and similarly
\bea
\frac{1}{r^2} \les \frac{1}{t^2} \Rightarrow \frac{1}{r^2} \les  \frac{1}{1+t^2} \;.
\eea

\textbf{Finally}\\

Estimating terms while considering the expression of the terms and the estimates that we have just established for each of the present terms, we get
\beaa
  \mid \pa_{t} h^0_{\mu\nu} \mid    &\leq&  \mid \frac{ M }{t^2}  \chi(r) \chi^\prime (r/t) \de_{\mu\nu}  \mid \\
  &\les&  \frac{ 1 }{r^2}  \mid  \chi(r) \chi^\prime (r/t) \de_{\mu\nu}  \mid \\
  &\les&  \frac{ 1 }{1 + r^2}  \mid  \chi(r) \chi^\prime (r/t) \de_{\mu\nu}  \mid \\
    &\les&  \frac{ 1 }{1 + t^2}  \mid  \chi(r) \chi^\prime (r/t) \de_{\mu\nu}  \mid \; .
\eeaa

Estimating for $i \in \{1, \ldots, n \}$, 
\beaa
 \mid \pa_{i} h^0_{\mu\nu}\mid  &\leq& \mid  \frac{x_i M}{r^2 t }  \chi^\prime (r/t)  \chi(r) \de_{\mu\nu}   \mid  +   \mid  \frac{x_i M }{r^2}  \chi (r/t)\chi^\prime (r)  \de_{\mu\nu}  \mid  +     \mid  \frac{- x_i M }{r^3}    \chi (r/t)\chi(r)   \de_{\mu\nu}  \mid   \\
 &\leq& \mid  \frac{ 3M}{r^2 }  \chi^\prime (r/t)  \chi(r) \de_{\mu\nu}   \mid  +   \mid  \frac{M }{r}  \chi (r/t)\chi^\prime (r)  \de_{\mu\nu}  \mid  +     \mid  \frac{M }{r^2}    \chi (r/t)\chi(r)   \de_{\mu\nu}  \mid   \\
  &\les& \frac{ 1 }{1 + r^2}   \mid    \chi^\prime (r/t)  \chi(r) \de_{\mu\nu}   \mid  +   \frac{1 }{r^2}  \mid   \chi (r/t)\chi^\prime (r)  \de_{\mu\nu}  \mid  +    \frac{1 }{r^2}  \mid     \chi (r/t)\chi(r)   \de_{\mu\nu}  \mid   \\
    &\les& \frac{ 1 }{1 + r^2}   \mid    \chi^\prime (r/t)  \chi(r) \de_{\mu\nu}   \mid  +   \frac{1 }{1+r^2}  \mid   \chi (r/t)\chi^\prime (r)  \de_{\mu\nu}  \mid  +    \frac{1 }{1+ r^2}  \mid     \chi (r/t)\chi(r)   \de_{\mu\nu}  \mid   \\
        &\les& \frac{ 1 }{1 + t^2}   \mid    \chi^\prime (r/t)  \chi(r) \de_{\mu\nu}   \mid  +   \frac{1 }{1+t^2}  \mid   \chi (r/t)\chi^\prime (r)  \de_{\mu\nu}  \mid  +    \frac{1 }{1+ t^2}  \mid     \chi (r/t)\chi(r)   \de_{\mu\nu}  \mid   \; .
\eeaa

We now look at cases depending on whether we are inside the outgoing light cone or if we are in the outside region of the light cone:

\textbf{Case $q< 0$:}\\

In that case, we have already established that
\beaa
\frac{1}{1+t} \les \frac{1}{1+t+|q|} ,
\eeaa
which implies that
\beaa
(1 +t+|q|)^2 &\les& (1+t)^2 \\
&\les& 1 +  2 t + t^2 \\
&\les& 1 + t^2 + t^2 \\
&&\text{(using $ab \les a^2 + b^2$)} \\
&\les& 1 + t^2 .
\eeaa
Hence, in this region $q< 0$, we have
\bea
\frac{1}{1+t^2} &\les& \frac{1}{(1+t+|q|)^2} .
\eea

The estimates could then be continued as follows
\beaa
  \mid \pa_{t} h^0_{\mu\nu} \mid     &\les&  \frac{ 1 }{1 + r^2}  \mid  \chi(r) \chi^\prime (r/t) \de_{\mu\nu}  \mid \\
      &\les&  \frac{ 1 }{1 + t^2}  \mid  \chi(r) \chi^\prime (r/t) \de_{\mu\nu}  \mid \\
       &\les&  \frac{1}{(1+t+|q|)^2}  \mid  \chi(r) \chi^\prime (r/t) \de_{\mu\nu}  \mid \\
       &&\text{(where we used that we are in $q < 0$)} \\
           &\les&  \frac{1}{(1+t+|q|)^2} .
\eeaa

Estimating for $i \in \{1, \ldots, 3 \}$, 
\beaa
 \mid \pa_{i} h^0_{\mu\nu}\mid      &\les& \frac{ 1 }{1 + r^2}   \mid    \chi^\prime (r/t)  \chi(r) \de_{\mu\nu}   \mid  +   \frac{1 }{1+r^2}  \mid   \chi (r/t)\chi^\prime (r)  \de_{\mu\nu}  \mid  +    \frac{1 }{1+ r^2}  \mid     \chi (r/t)\chi(r)   \de_{\mu\nu}  \mid   \\
        &\les& \frac{ 1 }{1 + t^2}   \mid    \chi^\prime (r/t)  \chi(r) \de_{\mu\nu}   \mid  +   \frac{1 }{1+t^2}  \mid   \chi (r/t)\chi^\prime (r)  \de_{\mu\nu}  \mid  +    \frac{1 }{1+ t^2}  \mid     \chi (r/t)\chi(r)   \de_{\mu\nu}  \mid   \\
                   &\les&  \frac{1}{(1+t+|q|)^2} \\
                         &&\text{(using that we are in the region $q < 0$).} \\
\eeaa

\textbf{Case $q>0$:}\\

In the region $q>0$, we have already established that
\beaa
\frac{1}{1+r} = \frac{1}{1+t+|q|},
\eeaa
which translates as $r= t+ |q|$, which implies $r^2 = ( t+ |q| )^2 $ and therefore it gives
\beaa
1 + r^2 = 1+ ( t+ |q| )^2 .
\eeaa
However,  $(1+  t+ |q| )^2  = 1 + 2( t+ |q| )+ ( t+ |q|)^2  \les 1 + ( t+ |q|)^2$ which implies that for $q>0$, 
\bea
\frac{1}{1 + r^2} = \frac{1}{1+ ( t+ |q| )^2 } \les  \frac{1}{(1+  t+ |q| )^2 } .
\eea

Going back now to our estimates on the terms, we have

\beaa
  \mid \pa_{t} h^0_{\mu\nu} \mid     &\les&  \frac{ 1 }{1 + r^2}  \mid  \chi(r) \chi^\prime (r/t) \de_{\mu\nu}  \mid \\
           &\les&  \frac{1}{(1+t+|q|)^2} \\
                 &&\text{(where we used that we are in $q > 0$)} .\\
\eeaa
For $i \in \{1, \ldots, n \}$, we get
\beaa
 \mid \pa_{i} h^0_{\mu\nu}\mid      &\les& \frac{ 1 }{1 + r^2}   \mid    \chi^\prime (r/t)  \chi(r) \de_{\mu\nu}   \mid  +   \frac{1 }{1+r^2}  \mid   \chi (r/t)\chi^\prime (r)  \de_{\mu\nu}  \mid  +    \frac{1 }{1+ r^2}  \mid     \chi (r/t)\chi(r)   \de_{\mu\nu}  \mid   \\
                   &\les&  \frac{1}{(1+t+|q|)^2} \\
                         &&\text{(using that we are in the region $q > 0$).} \\
\eeaa

\end{proof}

\begin{lemma}\label{Liederivativesofsphericalsymmetricpart}
  We have for all $I$,      
 \beaa
 \notag
|   \Lie_{Z^I} h^0 (t,x)  |   &\leq& c (\delta) \cdot C ( |I| )   \cdot \frac{M }{(1+t+|q|) }  \; .
      \eeaa
      and
 \beaa
 \notag
|\pa (  \Lie_{Z^I} h^0 ) (t,x)  |   &\leq&  C ( |I| )   \cdot \frac{M }{(1+t+|q|)^{2} } \; .
      \eeaa
\end{lemma}

\begin{proof}
We have proved the estimate for $|I| = 0$, in Lemma \ref{derivativesofsphericalsymmetricpart}. Now, we look at the case $|I| =1$: we have then that the vector $Z$ is either $\pa_\a$, for $\a \in \{0, 1, 2, 3\} $ and in that case, we have already proved the estimate in the previous Lemma \ref{derivativesofsphericalsymmetricpart}, or we have $Z$ taking the form of either of the following vector field, which is the case that we will examine here:
\beaa
Z_{ij} &=& x_{j} \pa_{i} - x_{i} \pa_{j} \, , \quad i, j \in \{1, 2, 3\} ; \\
Z_{ti} &=& x_{i} \pa_{t} + t \pa_{i} \, , \quad i \in \{1, 2, 3\} ; \\
S &=&   t \pa_{t}  + \sum_{i=1}^3 x_i \pa_{i}    \, .
\eeaa
We now fix $\mu, \nu \in (t, x^1, x^2, x^3)$. We will show that
 \beaa
 \notag
| Z^I h^0_{\mu\nu} (t,x)  |   &\leq& C ( |I| )   \cdot \frac{M }{(1+t+|q|) } \; , \\
 \notag
|\pa  Z^I h^0_{\mu\nu} (t,x)  |   &\leq& C ( |I| )   \cdot \frac{M }{(1+t+|q|)^{2} } \; .
      \eeaa
Using the fact that a commutation of two vector fields in $\cal Z$ is a combination of vector fields in $\cal Z$, and that a commutation of a vector field in $\cal Z$ and of a $\pa_\mu$ gives a linear combination of vectors of the form $\pa_\nu$, we will then get the result.

However, since $h^0_{\mu\nu}$ is spherically symmetric, we have
\bea
Z_{ij} h^0_{\mu\nu} = 0.
\eea

We showed in the proof of \ref{derivativesofsphericalsymmetricpart}, that or $t >0$,
\beaa
 \pa_{t} h^0_{\mu\nu} &=& - \frac{ M }{t^2}  \chi(r) \chi^\prime (r/t) \de_{\mu\nu} \\
 \pa_{i} h^0_{\mu\nu} &=&   \frac{x_i M}{r^2 t }  \chi^\prime (r/t)  \chi(r) \de_{\mu\nu}  +  \frac{x_i M }{r^2}  \chi (r/t)\chi^\prime (r)  \de_{\mu\nu} +    \frac{- x_i M }{r^3}    \chi (r/t)\chi(r)   \de_{\mu\nu}  . 
\eeaa
Computing
\beaa
Z_{ti} h^0_{\mu\nu} &=& x_{i} \pa_{t} h^0_{\mu\nu} + t \pa_{i} h^0_{\mu\nu}  \\
&=& - \frac{x_{i} M }{t^2}  \chi(r) \chi^\prime (r/t) \de_{\mu\nu}   +  \frac{x_i M}{r^2 }  \chi^\prime (r/t)  \chi(r) \de_{\mu\nu}  +  \frac{x_i M t }{r^2}  \chi (r/t)\chi^\prime (r)  \de_{\mu\nu} \\
&&+    \frac{- x_i M t}{r^3}    \chi (r/t)\chi(r)   \de_{\mu\nu}  
\eeaa
and
\beaa
S h^0_{\mu\nu} &=& t \pa_{t} h^0_{\mu\nu} +  \sum_{i=1}^3 x_i \pa_{i} h^0_{\mu\nu}  \\
&=& - \frac{ M }{t}  \chi(r) \chi^\prime (r/t) \de_{\mu\nu}   +  \sum_{i=1}^3 \big(  \frac{x_i^2 M}{r^2 t }  \chi^\prime (r/t)  \chi(r) \de_{\mu\nu}  +  \frac{x_i^2 M  }{r^2}  \chi (r/t)\chi^\prime (r)  \de_{\mu\nu} \\
&&+    \frac{- x_i^2 M }{r^3}    \chi (r/t)\chi(r)   \de_{\mu\nu}  \big) .
\eeaa

Considering the estimates that we established for each of the terms, we obtain
\beaa
\mid Z_{ti} h^0_{\mu\nu} \mid  &\leq& \mid - \frac{x_{i} M }{t^2}  \chi(r) \chi^\prime (r/t) \de_{\mu\nu}  \mid + \mid  \frac{x_i M}{r^2 }  \chi^\prime (r/t)  \chi(r) \de_{\mu\nu} \mid  + \mid  \frac{x_i M t }{r^2}  \chi (r/t)\chi^\prime (r)  \de_{\mu\nu} \mid  \\
&&+ \mid    \frac{- x_i M t}{r^3}    \chi (r/t)\chi(r)   \de_{\mu\nu} \mid   \\
&\leq& \frac{3}{4} \cdot \frac{3 M }{4r}  \mid \chi(r) \chi^\prime (r/t) \de_{\mu\nu}  \mid + \ \frac{ M}{r }  \mid  \chi^\prime (r/t)  \chi(r) \de_{\mu\nu} \mid  +  2 M  \mid \chi (r/t)\chi^\prime (r)  \de_{\mu\nu} \mid  \\
&&+     \frac{2 M }{r}   \mid  \chi (r/t)\chi(r)   \de_{\mu\nu} \mid   \\
&\les& \frac{ 1 }{r}  \mid \chi(r) \chi^\prime (r/t) \de_{\mu\nu}  \mid  +   \frac{ 1 }{r}   \mid \chi (r/t)\chi^\prime (r)  \de_{\mu\nu} \mid  +     \frac{1 }{r}   \mid  \chi (r/t)\chi(r)   \de_{\mu\nu} \mid  ,
\eeaa
and
\beaa
\mid S h^0_{\mu\nu} \mid &=&  \frac{ M }{t} \mid  \chi(r) \chi^\prime (r/t) \de_{\mu\nu}  \mid  +  \sum_{i=1}^3 \big(  \frac{x_i^2 M}{r^2 t } \mid \chi^\prime (r/t)  \chi(r) \de_{\mu\nu} \mid +  \mid \frac{x_i^2 M  }{r^2}  \chi (r/t)\chi^\prime (r)  \de_{\mu\nu} \mid \\
&&+    \frac{ x_i^2 M }{r^3}  \mid   \chi (r/t)\chi(r)   \de_{\mu\nu}  \mid  \big) \\
&\les &  \frac{ M }{t} \mid  \chi(r) \chi^\prime (r/t) \de_{\mu\nu}  \mid  +    \frac{ M}{ t } \mid \chi^\prime (r/t)  \chi(r) \de_{\mu\nu} \mid +  \mid  M    \chi (r/t)\chi^\prime (r)  \de_{\mu\nu} \mid \\
&&+    \frac{  M }{r}  \mid   \chi (r/t)\chi(r)   \de_{\mu\nu}  \mid  \\
&\les &  \frac{ 1 }{r} \mid  \chi(r) \chi^\prime (r/t) \de_{\mu\nu}  \mid  +  \frac{ 1 }{r} \mid      \chi (r/t)\chi^\prime (r)  \de_{\mu\nu} \mid +    \frac{  1 }{r}  \mid   \chi (r/t)\chi(r)   \de_{\mu\nu}  \mid  \; .
\eeaa
Consequently, we now distinguish the regions inside and outside the light-cone $q \geq 0$ and $q <0$, as before and we get the result. 

Concerning estimating the gradient, we differentiate: differentiating with respect to $\pa_t$, we get
\beaa
 \pa_t Z_{ti} h^0_{\mu\nu} &=&  \frac{2 x_{i} M }{t^3}  \chi(r) \chi^\prime (r/t) \de_{\mu\nu}   + \frac{x_{i} r M }{t^4}   \chi(r) \chi^{\prime\prime} (r/t) \de_{\mu\nu}   -  \frac{x_i r M}{r^2 t^2}  \chi^{\prime\prime} (r/t)  \chi(r) \de_{\mu\nu}  \\
&&  +  \frac{x_i M  }{r^2}  \chi (r/t)\chi^\prime (r)  \de_{\mu\nu} -  \frac{x_i M  }{r t}  \chi^{\prime} (r/t)\chi^\prime (r)  \de_{\mu\nu}  \\
&& +    \frac{- x_i M }{r^3}    \chi (r/t)\chi(r)   \de_{\mu\nu}   -    \frac{- x_i M }{r^2 t}    \chi^\prime (r/t)\chi(r)   \de_{\mu\nu}  \; .
\eeaa
Proceeding similarly for the other terms and distingishing the regions inside and outside the outgoing light-cone whose tip is the origin, as earlier, we get the result.

\end{proof} 
  \begin{lemma} \label{tangentialderivativesphericalsymmetricpart}
 We have for all $I$,      
 \beaa
 \notag
|  \rpa  ( \Lie_{Z^I} h^0 ) (t,x)  |   &\leq&   C ( |I| )   \cdot \frac{ M }{(1+t+|q|)^{2} }  \; .
      \eeaa
\end{lemma}

\begin{proof}
Using following estimate, combined with Lemma \ref{Liederivativesofsphericalsymmetricpart},

   \beaa
 \notag
|  \rpa  ( \Lie_{Z^I} h^0 ) (t,x)  |   &\leq&  |  \pa  ( \Lie_{Z^I} h^0 ) (t,x)  |    \leq   C ( |I| )   \cdot \frac{ M }{(1+t+|q|)^{2} }  \; ,
      \eeaa
we obtain the result.

  \end{proof}

We also have the following lemma from Lindblad-Rodnianski in \cite{LR10}.

\begin{lemma}\label{estimateonthesourcetermsforhzerothesphericallsymmtrpart}
Let $M \leq \eps$, then we have
\beaa
| \Lie_{ Z^I}  ( g^{\la\mu} \derm_{\la}   \derm_{\mu}    h^0 ) | &\les& C ( |I| ) \cdot E ( |I| + 2)  \cdot \begin{cases} c (\gamma) \cdot \frac{\eps }{(1+t+|q|)^{4-\delta}(1+|q|)^{\delta}},\quad q>0 \;, \\
\frac{\eps}{ (1+t+|q|)^{3}} ,\quad q<0 \;.
\end{cases}
\eeaa
Thus, for $\de \leq 1$, and $M \leq \eps$, we have
\beaa
| \Lie_{ Z^I}  ( g^{\la\mu} \derm_{\la}   \derm_{\mu}    h^0 ) | &\les& C ( |I| ) \cdot c (\gamma) \cdot  E ( |I| + 2)  \cdot \frac{\eps }{(1+t+|q|)^{3} } \; .
\eeaa

\end{lemma}

\begin{proof}      
See Lemma 9.9 in \cite{LR10}, and for calculations, see proof of Lemma 8.4 in \cite{Loiz2}. Note that the bootstrap assumption is being used to estimate the Lie derivatives of $H$\,, which appear when taking the Lie derivatives of the wave operator $( H^{\la\mu} + m^{\la\mu} ) \derm_{\la}   \derm_{\mu}    h^0$\, (where we used the definition of $H$ in \ref{definitionsofbigH} to decompose $g$), using Lemma \ref{linkbetweenbigHandsamllh}. Note that the assumption $M \leq \eps$ is being used to estimate $H$ (see Lemma \ref{aprioriestimatesonZLiederivativesonbigH}), by estimating first $h$ (see Lemma \ref{aprioridecayestimates}) and the using Lemma \ref{linkbetweenbigHandsamllh}.
\end{proof}

\begin{remark}\label{thesphericallysymmmetricparthasinfiniteenergy}
The tensor $h^0 $ does not have finite energy. In fact, we have
\beaa
  \mid \pa_{t} h^0_{\mu\nu} \mid    &\sim&  \mid \frac{ M }{t^2}  \chi(r) \chi^\prime (r/t) \de_{\mu\nu}  \mid ,
\eeaa
and therefore,
\beaa
  \mid \pa_{t} h^0_{\mu\nu} \mid^2    &\sim& \frac{ M^2 }{t^4}  \mid   \chi(r) \chi^\prime (r/t) \de_{\mu\nu}  \mid^2 .
\eeaa
On the other hand, we have
\beaa
 \mid \pa_{i} h^0_{\mu\nu}\mid  &\sim& \mid  \frac{x_i M}{r^2 t }  \chi^\prime (r/t)  \chi(r) \de_{\mu\nu}   \mid  +   \mid  \frac{x_i M }{r^2}  \chi (r/t)\chi^\prime (r)  \de_{\mu\nu}  \mid  +     \mid  \frac{- x_i M }{r^3}    \chi (r/t)\chi(r)   \de_{\mu\nu}  \mid   \\
 &\sim& \mid  \frac{ 3M}{r^2 }  \chi^\prime (r/t)  \chi(r) \de_{\mu\nu}   \mid  +   \mid  \frac{M }{r}  \chi (r/t)\chi^\prime (r)  \de_{\mu\nu}  \mid  +     \mid  \frac{M }{r^2}    \chi (r/t)\chi(r)   \de_{\mu\nu}  \mid  ,
\eeaa
and therefore
\beaa
 \mid \pa_{i} h^0_{\mu\nu}\mid^2   &\sim&  \frac{ M^2}{r^4 }  \mid \chi^\prime (r/t)  \chi(r) \de_{\mu\nu}   \mid^2  +    \frac{M^2 }{r^2}  \mid \chi (r/t)\chi^\prime (r)  \de_{\mu\nu}  \mid^2  +      \frac{M^2 }{r^4}   \mid  \chi (r/t)\chi(r)   \de_{\mu\nu}  \mid^2  . \\
\eeaa

As a result, considering the integral $  \|w^{1/2}   \pa   h^0   (t,\cdot) )  \|_{L^2 (\Sigma^{ext}_{t} )}  $ at $t= constant$, we have that $w^{1/2}$ behaves as $1 + r^\frac{1}{2}$ and thus,
\beaa
  && \|w^{1/2}   \pa   h^0   (t,\cdot) )  \|_{L^2 (\Sigma^{ext}_{t} )}  \\
  &\sim&  \int_{\Sigma^{ext}_{t} } \big( r^2  \mid \pa_{t} h^0_{\mu\nu} \mid^2  + r^2  \mid \pa_{i} h^0_{\mu\nu}\mid^2 \big) dx^1 dx^2 dx^3 \\
  &\sim&  \int_{\Sigma^{ext}_{t} } \big(  \frac{ M^2 r^2 }{t^4}  \mid   \chi(r) \chi^\prime (r/t) \de_{\mu\nu}  \mid^2  + \frac{ M^2}{r^2 }  \mid \chi^\prime (r/t)  \chi(r) \de_{\mu\nu}   \mid^2  +    M^2   \mid \chi (r/t)\chi^\prime (r)  \de_{\mu\nu}  \mid^2 \\
  && +      \frac{M^2 }{r^2}   \mid  \chi (r/t)\chi(r)   \de_{\mu\nu}  \mid^2   \big) dx^1 dx^2 dx^3 \\
    &\sim&    \int_{r=r^{ext}}^{\infty}   \int_{\th=0}^{\pi}  \int_{\phi=0}^{2\pi} \big(  \frac{ M^2  }{t^2}  \mid   \chi(r) \chi^\prime (r/t) \de_{\mu\nu}  \mid^2  + \frac{ M^2}{r^2 }  \mid \chi^\prime (r/t)  \chi(r) \de_{\mu\nu}   \mid^2  +    M^2   \mid \chi (r/t)\chi^\prime (r)  \de_{\mu\nu}  \mid^2 \\
  && +      \frac{M^2 }{r^2}   \mid  \chi (r/t)\chi(r)   \de_{\mu\nu}  \mid^2   \big) \cdot r^2 \sin( \th )dr d\th d\phi  \\
  && \text{(where we used polar coordinates $r, \th, \phi$ on $\Sigma_t$)} \\
    &\sim&  \int_{r=r^{ext}}^{\infty} \big(  \frac{ M^2 r^2  }{t^2}  \mid   \chi(r) \chi^\prime (r/t) \de_{\mu\nu}  \mid^2  + M^2  \mid \chi^\prime (r/t)  \chi(r) \de_{\mu\nu}   \mid^2  +    M^2  r^2 \mid \chi (r/t)\chi^\prime (r)  \de_{\mu\nu}  \mid^2 \\
  && +      M^2 \mid  \chi (r/t)\chi(r)   \de_{\mu\nu}  \mid^2   \big)  dr   \\
 && \text{(which clearly blows up).} \\
\eeaa
Similarly, one can show that for the Lie derivatives of the fields as well, the integral $  \|w^{1/2}   \derm  \Lie_{Z^{I}}  h^0   (t,\cdot) )  \|_{L^2  (\Sigma^{ext}_{t} )}  $ blows up. Hence, the energy is not defined for $h$ but only for the tensor $h^1$ that is expected to decay.
\end{remark}

\section{Studying the structure of the source term in the null frame}

We have shown in Lemma \ref{Einstein-Yang-MillssysteminLorenzwavegauges}, that the Einstein-Yang-Mills system in the Lorenz gauge and in wave coordinates is given by a system of non-linear  wave operators on $A$ and on $g$. To study the structure of these non-linearities in the source terms of these hyperbolic operators $g^{\la\mu} \pa_{\la}   \pa_{\mu}   A_{\si} $ and $ g^{\alpha\beta}\pa_\alpha\pa_\beta
g_{\mu\nu}$, we introduce the null frame that will help us decompose the terms into “good” terms and “bad” terms: the good terms will be terms that we could control using estimates that exploit the Lorenz gauge and the wave coordinates conditions -- this will be apparent by the control by “good derivatives” on the “good terms” that we will show.

\begin{definition}\label{definitionofthenullframusingwavecoordinates}
At a point $p$ in the space-time, let
\bea
L &=&  \pa_{t} + \pa_{r} \, =   \pa_{t} +  \frac{x^{i}}{r} \pa_{i} \, , \\
\underline{L} &=&  \pa_{t} - \pa_{r} =  \pa_{t} -  \frac{x^{i}}{r} \pa_{i} \, ,
\eea
and let $\{e_1, e_{2} \}$ be an orthonormal frame on $\SSS^{2}$. 
We define the sets 
\bea
\cal T&=& \{ L,e_1, e_{2}\} \, ,\\
 \cal U&=&\{ \underline{L}, L, e_1, e_{2}\} \, .
\eea
We call the set  $ \cal U$ a null-frame.
\end{definition}

\begin{lemma}
We have for all $A, B \in \{ 1, 2\}$, 
\bea
m_{LL} &=& m_{\underline{L} \underline{L}} = m_{Le_A} =  m_{\underline{L} e_A} = 0 \; , \\
m^{LL} &=& m^{\underline{L} \underline{L}} = m^{Le_A} =  m^{\underline{L} e_A} = 0\; , \\
m_{L \underline{L}} &=& -2 \; , \\
m_{e_A e_B}  &=&  \delta_{AB} \; ,
\eea
and 
\bea
m^{L \underline{L}} &=& - \frac{1}{2} \; ,\\
m^{e_A e_B}  &=&  \delta_{AB}\; .
\eea
\end{lemma}

\begin{proof}

We compute
\beaa
m(L,L)&=& m( \pa_{t} + \pa_{r},\pa_{t} + \pa_{r}) =m( \pa_{t} +  \frac{x^{i}}{r} \pa_{i},  \pa_{t} +  \frac{x^{j}}{r} \pa_{j}) = m( \pa_{t} ,  \pa_{t} ) + m (  \frac{x^{i}}{r} \pa_{i},   \frac{x^{j}}{r} \pa_{j}) \\
&=&m( \pa_{t} ,  \pa_{t} ) +  \frac{x^{i}.x^{i}}{r^2} m (  \pa_{i},  \pa_{i}) = -1+1 =0 ,\\
m(\underline{L},\underline{L})&=& m( \pa_{t} - \pa_{r},\pa_{t} - \pa_{r}) = m ( \pa_{t}, \pa_{t}) + m( \pa_{r}, \pa_{r}) = -1+1 =0 .
\eeaa
We have for all $A \in \{1, 2 \}$,
\beaa
m(L,e_A) =  m(\underline{L}, e_A) = 0 ,
\eeaa
and
\beaa
m(L, \underline{L}) &=& m( \pa_{t} + \pa_{r}, \pa_{t} - \pa_{r} ) = -1-1 = -2 \\
m(e_A, e_B) &=& \delta_{AB} .
\eeaa
 \end{proof}
 
\begin{remark}
To lighten the notation, in what follows, we write $e_A$ to designate a sum over $A \in \{1, 2\}$.
\end{remark}

\begin{lemma}\label{TheEinstein-Yang-Millssysteminanullframe}
In the Lorenz gauge the Yang-Mills potential satisfies the following equations decomposed in a null-frame, i.e. for any $\si \in \cal U$, we have
     \bea
   \notag
&& g^{\la\mu} \derm_{\la}   \derm_{\mu}   A_{\si}   \\ 
  \notag
 &=&  \frac{1}{4} (  \derm_{\si} h_{\underline{L}\underline{L}} ) \cdot  \derm_{L}A_{L}  + \frac{1}{4} (   \derm_{\si}  h_{\underline{L}L} )  \cdot \derm_{L}A_{\underline{L}}   -\frac{1}{2}   (  \derm_{\si}  h_{\underline{L}e_{A}} ) \cdot \derm_{L} A_{e_{A}}  \\
   \notag
 && +  \frac{1}{4} (  \derm_{\si}  h_{L\underline{L}} ) \cdot \derm_{\underline{L} }A_{L} +  \frac{1}{4}  (  \derm_{\si}  h_{LL} ) \cdot \derm_{\underline{L} }A_{\underline{L}} -  \frac{1}{2} (  \derm_{\si}  h_{Le_{A}} ) \cdot  \derm_{\underline{L} }A_{e_{A}} \\
    \notag
&&   -  \frac{1}{2} (  \derm_{\si} h_{e_{A}\underline{L}} )  \cdot \derm_{e_{A}}A_{L}   -  \frac{1}{2} (   \derm_{\si} h_{e_{A}L} ) \cdot  \derm_{e_{A}}A_{\underline{L}} +  ( \derm_{\si}  h_{e_{A}e_{A}} )  \cdot \derm_{e_{A}}A_{e_{A}} \\
   \notag
&& + \frac{1}{4}  \big(    \derm_L h_{\underline{L}\si} -  \derm_{\underline{L}} h_{L \si}  \big)   \cdot  \big( \derm_{\underline{L}}A_{L} -  \derm_{L}A_{\underline{L}}  \big) \\
   \notag
&& - \frac{1}{2}    \big(    \derm_L h_{e_{A}\si} -  \derm_{e_{A}} h_{L \si}  \big)   \cdot   \big(  \derm_{\underline{L}}A_{e_{A}} -  \derm_{e_{A}}A_{\underline{L}}  \big) \\
     \notag
   && - \frac{1}{2}   \big(    \derm_{\underline{L}} h_{e_{A}\si} -  \derm_{e_{A}} h_{\underline{L}\si}  \big)   \cdot \big(  \derm_{L}A_{e_{A}} -  \derm_{e_{A}}A_{L}  \big) \\
 \notag
&& + \frac{1}{4}   \big(    \derm_L h_{\underline{L}\si} -  \derm_{\underline{L}} h_{L \si}  \big)  \cdot   [A_{\underline{L}},A_{L}]    \\
   \notag
&& - \frac{1}{2}   \big(    \derm_L h_{e_{A}\si} -  \derm_{e_{A}} h_{L \si}  \big)   \cdot   [A_{\underline{L}},A_{e_{A}}]   \\
    \notag
   && - \frac{1}{2}   \big(   \derm_{\underline{L}} h_{e_{A}\si} - \derm_{e_{A}} h_{\underline{L}\si}  \big)   \cdot  [A_{L},A_{e_{A}}]    \\
   \notag
 && +  \frac{1}{2} \big(  [ A_{\underline{L}},  \derm_{L} A_{\si} ]  +   [A_{L},   \derm_{\underline{L}}  A_{\si} - \derm_{\si} A_{\underline{L}} ]    +    [A_{L}, [A_{\underline{L}}, A_{\si}] ]  \big)  \\
    \notag
    &&  +  \frac{1}{2} \big(  [ A_{L},  \derm_{\underline{L}} A_{\si} ]  +    [A_{\underline{L}},   \derm_{L}  A_{\si} - \derm_{\si} A_{L} ]    +    [A_{\underline{L}}, [A_{L},A_{\si}] ]  \big)  \\
       \notag
 &&   -   \big(  [ A_{e_{a}}, \derm_{e_{a}} A_{\si} ] + [A_{e_{a}},  \derm_{e_{a}} A_{\si}  -  \derm_{\si} A_{e_{a}} ]     +    [A_{e_{a}}, [A_{e_{a}},A_{\si}] ]  \big) \\
  && + O( h \cdot  \pa h \cdot  \pa A) + O( h \cdot  \pa h \cdot  A^2) + O( h \cdot  A \cdot \pa A) + O( h \cdot  A^3) \, . 
  \eea

The perturbations $h$ of the metric $m$, defined to be the Minkowski metric in wave coordinates, solutions to the Einstein-Yang-Mills equations satisfy the following wave equation, i.e. for any $\mu, \nu \in \cal U$, we have
    \bea
\notag
&& g^{\alpha\beta} \derm_\alpha \derm_\beta h_{\mu\nu}   \\
\notag
     &=&  P(\pa_\mu h,\pa_\nu h) +  Q_{\mu\nu}(\pa h,\pa h)   + G_{\mu\nu}(h)(\pa h,\pa h)  \\
\notag
 &&  +2  <    \derm_{\mu}A_{L} -  \derm_{L}A_{\mu}  ,  \derm_{\nu}A_{ \underline{L} } -  \derm_{ \underline{L} }A_{\nu}  >   \\
 \notag
 &&  +2 <   \derm_{\mu}A_{ \underline{L} } - \derm_{ \underline{L} }A_{\mu}  ,  \derm_{\nu}A_{L} - \derm_{L}A_{\nu}  >    \\
 \notag
  && -4   <    \derm_{\mu}A_{e_A} -  \derm_{e_A }A_{\mu}  ,   \derm_{\nu}A_{e_A} - \derm_{e_A}A_{\nu}  > \\
  \notag
  &&   + \frac{1}{2} m_{\mu\nu }    \cdot  <   \derm_{L}A_{\underline{L}} - \derm_{\underline{L}}A_{L} ,  \derm_{\underline{L}}A_{L} - \derm_{L}A_{\underline{L}} > \\
\notag
&& -  2 m_{\mu\nu }     \cdot  <   \derm_{e_A}A_{\underline{L}} -  \derm_{\underline{L}}A_{e_A} ,  \derm_{e_A}A_{L} - \derm_{L}A_{e_A} > \\
\notag
 &&  +2  \big( <   \derm_{\mu}A_{\underline{L}} - \derm_{\underline{L}}A_{\mu}  ,  [A_{\nu},A_{L}] >   + <   [A_{\mu},A_{\underline{L}}] ,   \derm_{\nu}A_{L} -  \derm_{L}A_{\nu}  > \big) \\
 \notag
&&   + 2  \big( <   \derm_{\mu}A_{L} - \derm_{L}A_{\mu}  ,  [A_{\nu},A_{\underline{L}}] >   + <   [A_{\mu},A_{L}] ,  \derm_{\nu}A_{\underline{L}} - \derm_{\underline{L}}A_{\nu}  > \big)  \\
\notag
 &&  -4     \big( <    \derm_{\mu}A_{e_A} - \derm_{e_A}A_{\mu}  ,  [A_{\nu},A_{e_A}] >    -4  <   [A_{\mu},A_{e_A}] ,   \derm_{\nu}A_{e_A} -  \derm_{e_A}A_{\nu}  > \big)  \\
\notag
&&  +  m_{\mu\nu }   \cdot   <  \derm_{L}A_{\underline{L}} -  \derm_{\underline{L}}A_{L} , [A_{\underline{L}},A_{L}] >   \\
\notag
 &&  -  2 m_{\mu\nu }     \cdot \big(  <   \derm_{e_A }A_{\underline{L}} -  \derm_{\underline{L}}A_{e_A} , [A_{e_A},A_{L}] >    +  <  [A_{e_A},A_{\underline{L}}] ,  \derm_{e_A}A_{L} -  \derm_{L}A_{e_A}  > \big) \\
 \notag
 && + 2   <   [A_{\mu},A_{\underline{L}}] ,  [A_{\nu},A_{L}] >  + 2   <   [A_{\mu},A_{L}] ,  [A_{\nu},A_{\underline{L}}] >  -4    <   [A_{\mu},A_{e_A}] ,  [A_{\nu},A_{e_A}] > \\
 \notag
  && +   \frac{1}{2}   m_{\mu\nu }    \cdot   <  [A_{L},A_{\underline{L}}] , [A_{\underline{L}},A_{L}] >   - 2 m_{\mu\nu }   \cdot   <  [A_{e_A},A_{\underline{L}}] , [A_{e_A},A_{L}] > \\
     && + O \big(h \cdot  (\pa A)^2 \big)   + O \big(  h  \cdot  A^2 \cdot \pa A \big)     + O \big(  h   \cdot  A^4 \big)  \,  . 
\eea

    \end{lemma}
    \begin{proof}
We showed in Lemma \ref{Einstein-Yang-MillssysteminLorenzwavegauges}, that in the Lorenz gauge and in wave coordinates -- in other words for indices running only over wave coordinates, i.e. $\la, \mu, \si, \b, \nu, \a \in \{t, x^1, \ldots, x^n \} $ --, that the Yang-Mills potential satisfies 
          \beaa
   \notag
g^{\la\mu} \pa_{\la}   \pa_{\mu}   A_{\si}      &=&  m^{\a\ga} m ^{\mu\la} (   \pa_{\si}  h_{\ga\la} )   \pa_{\a}A_{\mu}       +   \frac{1}{2}  m^{\a\mu}m^{\b\nu}   \big(   \pa_\a h_{\b\si} + \pa_\si h_{\b\a}- \pa_\b h_{\a\si}  \big)   \cdot  \big( \pa_{\mu}A_{\nu} - \pa_{\nu}A_{\mu}  \big) \\
 \notag
&& +      \frac{1}{2}  m^{\a\mu}m^{\b\nu}   \big(   \pa_\a h_{\b\si} + \pa_\si h_{\b\a}- \pa_\b h_{\a\si}  \big)   \cdot   [A_{\mu},A_{\nu}] \\
 \notag
 && -  m^{\a\mu} \big(  [ A_{\mu}, \pa_{\a} A_{\si} ]  +    [A_{\alpha},  \pa_{\mu}  A_{\si} - \pa_{\si} A_{\mu} ]    +    [A_{\alpha}, [A_{\mu},A_{\si}] ]  \big)  \\
  && + O( h \cdot  \pa h \cdot  \pa A) + O( h \cdot  \pa h \cdot  A^2) + O( h \cdot  A \cdot \pa A) + O( h \cdot  A^3) \, .
  \eeaa

 Since the Christoffel symbols for the connection $\derm$ are vanishing in wave coordinates, we could then write,
            \beaa
   \notag
g^{\la\mu} \derm_{\la}   \derm_{\mu}   A_{\si}     &=&  m^{\a\ga} m ^{\mu\la}  ( \derm_{\si}  h_{\ga\la} ) \cdot  \derm_{\a}A_{\mu}     \\
&&  +   \frac{1}{2}  m^{\a\mu}m^{\b\nu}   \big(   \derm_\a h_{\b\si} +  \derm_\si h_{\b\a}-  \derm_\b h_{\a\si}  \big)   \cdot  \big( \derm_{\mu}A_{\nu} -  \derm_{\nu}A_{\mu}  \big) \\
 \notag
&& +      \frac{1}{2}  m^{\a\mu}m^{\b\nu}   \big(   \derm_\a h_{\b\si} +  \derm_\si h_{\b\a}-  \derm_\b h_{\a\si}  \big)   \cdot   [A_{\mu},A_{\nu}] \\
 \notag
 && -  m^{\a\mu} \big(  [ A_{\mu}, \derm_{\a} A_{\si} ]  +    [A_{\alpha},  \derm_{\mu}  A_{\si} - \derm_{\si} A_{\mu} ]    +    [A_{\alpha}, [A_{\mu},A_{\si}] ]  \big)  \\
  && + O( h \cdot  \derm h \cdot  \derm A) + O( h \cdot  \derm h \cdot  A^2) + O( h \cdot  A \cdot \derm A) + O( h \cdot  A^3) \, .
  \eeaa
 
 Also, we have

 \beaa
 g^{\la\mu} \pa_{\la}   \pa_{\mu}   A_{\si} = g^{\la\mu} \derm_{\la}   \derm_{\mu}   A_{\si} ,
 \eeaa
 
 which is a tensor in $\si$. Thus, the right hand side and the left hand side of the following equation is a tensor in $\si$ and corresponds to a full tensorial contraction on all other indices and hence the expression does not depend on the system of coordinates that we choose. Thus, we can compute the trace in the null-frame $\cal U$.

Decomposing the first term in the null-frame  $\{ \underline{L}, L, e_1, \ldots, e_{(n-1)}\} $, and where we write $e_A$ to designate a sum over $A \in \{1, \ldots, (n-1) \}$ so as to lighten the notation, we get,  
\beaa
  && m^{\a\ga} m ^{\mu\la}  (  \derm_{\si}  h_{\ga\la} )   \derm_{\a}A_{\mu}   \\
 &=&   m^{L\ga} m ^{\mu\la} (  \derm_{\si}  h_{\ga\la} )   \derm_{L}A_{\mu}  +  m^{\underline{L} \ga} m ^{\mu\la} (   \derm_{\si}  h_{\ga\la} )   \derm_{\underline{L} }A_{\mu}   +  m^{e_{A}\ga} m ^{\mu\la}  (  \derm_{\si}  h_{\ga\la} )  \derm_{e_{A}}A_{\mu} \\
 &=&   m^{L\underline{L}} m ^{\mu\la} (  \derm_{\si}  h_{\underline{L}\la} )  \derm_{L}A_{\mu}  +  m^{\underline{L} L} m ^{\mu\la} (  \derm_{\si}  h_{L\la} )   \derm_{\underline{L} }A_{\mu}   +  m^{e_{A}e_{A}} m ^{\mu\la}  (  \derm_{\si}  h_{e_{A}\la} ) \derm_{e_{A}}A_{\mu} \\
 &=&   m^{L\underline{L}} m ^{L\la} (   \derm_{\si} h_{\underline{L}\la} )   \derm_{L}A_{L}  + m^{L\underline{L}} m ^{\underline{L}\la}  ( \derm_{\si}  h_{\underline{L}\la} )  \derm_{L}A_{\underline{L}}   +m^{L\underline{L}} m ^{e_{A}\la}  \derm_{\si} ( h_{\underline{L}\la} )   \pa_{L} A_{e_{A}}  \\
 && +  m^{\underline{L} L} m ^{L\la}  (  \derm_{\si}  h_{L\la} )   \derm_{\underline{L} }A_{L} +  m^{\underline{L} L} m ^{\underline{L}\la} (  \derm_{\si}  h_{L\la} )  \derm_{\underline{L} }A_{\underline{L}} +  m^{\underline{L} L} m ^{e_{A}\la}  (  \derm_{\si}  h_{L\la} )   \derm_{\underline{L} }A_{e_{A}} \\
&&   +  m^{e_{A}e_{A}} m ^{L\la} (   \derm_{\si}  h_{e_{A}\la} )   \derm_{e_{A}}A_{L} +  m^{e_{A}e_{A}} m ^{\underline{L}\la} (  \derm_{\si}  h_{e_{A}\la} )   \derm_{e_{A}}A_{\underline{L}} \\
&& +  m^{e_{A}e_{A}} m ^{e_{A}\la}  (  \derm_{\si}  h_{e_{A}\la} )   \derm_{e_{A}}A_{e_{A}} \\
 &=&   m^{L\underline{L}} m ^{L\underline{L}} (  \derm_{\si}  h_{\underline{L}\underline{L}} )   \derm_{L}A_{L}  + m^{L\underline{L}} m ^{\underline{L}L} (  \derm_{\si}  h_{\underline{L}L} )   \pa_{L}A_{\underline{L}}   +m^{L\underline{L}} m ^{e_{A}e_{A}}  ( \derm_{\si}  h_{\underline{L}e_{A}} )   \derm_{L} A_{e_{A}}  \\
 && +  m^{\underline{L} L} m ^{L\underline{L}}   (  \derm_{\si}  h_{L\underline{L}} )  \derm_{\underline{L} }A_{L} +  m^{\underline{L} L} m ^{\underline{L}L} (  \derm_{\si}  h_{LL} )   \derm_{\underline{L} }A_{\underline{L}} \\
 && +  m^{\underline{L} L} m ^{e_{A}e_{A}} (  \derm_{\si}  h_{Le_{A}} )   \derm_{\underline{L} }A_{e_{A}}  +  m^{e_{A}e_{A}} m ^{L\underline{L}}  (  \derm_{\si}  h_{e_{A}\underline{L}} )   \derm_{e_{A}}A_{L} \\
 && +  m^{e_{A}e_{A}} m ^{\underline{L}L} (  \derm_{\si}  h_{e_{A}L} )  \derm_{e_{A}}A_{\underline{L}}  +  m^{e_{A}e_{A}} m ^{e_{A}e_{A}}  (  \derm_{\si}  h_{e_{A}e_{A}} )   \derm_{e_{A}}A_{e_{A}} .
\eeaa 

For the second term, we compute
\beaa
 && m^{\a\mu}m^{\b\nu}   \big(   \derm_\a h_{\b\si} + \derm_\si h_{\b\a}-  \derm_\b h_{\a\si}  \big)   \cdot  \big( \derm_{\mu}A_{\nu} - \derm_{\nu}A_{\mu}  \big) \\
 &=& m^{L\mu}m^{\b\nu}   \big(  \derm_L h_{\b\si} + \derm_\si h_{\b L}- \derm_\b h_{L \si}  \big)   \cdot  \big( \derm_{\mu}A_{\nu} - \derm_{\nu}A_{\mu}  \big)  \\
 && + m^{\underline{L} \mu}m^{\b\nu}   \big(   \derm_{\underline{L}} h_{\b\si} +\derm_\si h_{\b\underline{L}}- \derm_\b h_{\underline{L}\si}  \big)   \cdot  \big( \derm_{\mu}A_{\nu} - \derm_{\nu}A_{\mu}  \big) \\
 && +  m^{e_{A}\mu}m^{\b\nu}   \big(   \derm_{e_{A}} h_{\b\si} + \derm_\si h_{\b e_{A}}- \derm_\b h_{e_{A} \si}  \big)   \cdot  \big( \derm_{\mu}A_{\nu} -\derm_{\nu}A_{\mu}  \big) \\
 &=& m^{L\underline{L}}m^{\b\nu}   \big(   \derm_L h_{\b\si} + \derm_\si h_{\b L}-\derm_\b h_{L \si}  \big)   \cdot  \big( \derm_{\underline{L}}A_{\nu} - \derm_{\nu}A_{\underline{L}}  \big)  \\
 && + m^{\underline{L} L}m^{\b\nu}   \big(   \derm_{\underline{L}} h_{\b\si} + \derm_\si h_{\b\underline{L}}- \derm_\b h_{\underline{L}\si}  \big)   \cdot  \big( \derm_{L}A_{\nu} - \derm_{\nu}A_{L}  \big) \\
 && +  m^{e_{A}e_{A}}m^{\b\nu}   \big(   \derm_{e_{A}} h_{\b\si} + \derm_\si h_{\b e_{A}}- \derm_\b h_{e_{A} \si}  \big)   \cdot  \big( \derm_{e_{A}}A_{\nu} - \derm_{\nu}A_{e_{A}}  \big) \\
&=& m^{L\underline{L}}m^{L\nu}   \big(   \derm_L h_{L\si} + \derm_\si h_{L L}- \derm_L h_{L \si}  \big)   \cdot   \big( \derm_{\underline{L}}A_{\nu} -\derm_{\nu}A_{\underline{L}}  \big) \\
&& + m^{L\underline{L}}m^{\underline{L}\nu}   \big(   \derm_L h_{\underline{L}\si} + \derm_\si h_{\underline{L} L}-\derm_{\underline{L}} h_{L \si}  \big)   \cdot  \big( \derm_{\underline{L}}A_{\nu} - \pa_{\nu}A_{\underline{L}}  \big) \\
&& + m^{L\underline{L}}m^{e_{A}\nu}   \big(   \derm_L h_{e_{A}\si} + \derm_\si h_{e_{A} L}- \derm_{e_{A}} h_{L \si}  \big)   \cdot   \big( \derm_{\underline{L}}A_{\nu} - \derm_{\nu}A_{\underline{L}}  \big) \\
 && + m^{\underline{L} L}m^{L\nu}   \big(   \derm_{\underline{L}} h_{L\si} +\derm_\si h_{L\underline{L}}- \derm_L h_{\underline{L}\si}  \big)   \cdot  \big(\derm_{L}A_{\nu} -\derm_{\nu}A_{L}  \big) \\
  && + m^{\underline{L} L}m^{\underline{L}\nu}   \big(  \derm_{\underline{L}} h_{\underline{L}\si} + \derm_\si h_{\underline{L}\underline{L}}- \derm_{\underline{L}} h_{\underline{L}\si}  \big)   \cdot  \big( \derm_{L}A_{\nu} - \derm_{\nu}A_{L}  \big) \\
   && + m^{\underline{L} L}m^{e_{A}\nu}   \big(   \derm_{\underline{L}} h_{e_{A}\si} + \derm_\si h_{e_{A}\underline{L}}- \derm_{e_{A}} h_{\underline{L}\si}  \big)   \cdot \big( \derm_{L}A_{\nu} - \derm_{\nu}A_{L}  \big) \\
 && +  m^{e_{A}e_{A}}m^{L\nu}   \big(  \derm_{e_{A}} h_{L\si} +\derm_\si h_{L e_{A}}- \derm_L h_{e_{A} \si}  \big)   \cdot  \big( \derm_{e_{A}}A_{\nu} - \derm_{\nu}A_{e_{A}}  \big)  \\
  && +  m^{e_{A}e_{A}}m^{\underline{L}\nu}   \big(   \derm_{e_{A}} h_{\underline{L}\si} + \derm_\si h_{\underline{L} e_{A}}-\derm_{\underline{L}} h_{e_{A} \si}  \big)   \cdot   \big( \derm_{e_{A}}A_{\nu} - \derm_{\nu}A_{e_{A}}  \big)   \\
   && +  m^{e_{A}e_{A}}m^{e_{A}\nu}   \big(  \derm_{e_{A}} h_{e_{A}\si} + \derm_\si h_{e_{A} e_{A}}- \derm_{e_{A}} h_{e_{A} \si}  \big)   \cdot   \big( \derm_{e_{A}}A_{\nu} - \derm_{\nu}A_{e_{A}}  \big) \\
&=& m^{L\underline{L}}m^{L\underline{L}}   \big(  \derm_L h_{L\si} + \derm_\si h_{L L}- \derm_L h_{L \si}  \big)   \cdot   \big( \derm_{\underline{L}}A_{\underline{L}} - \derm_{\underline{L}}A_{\underline{L}}  \big) \\
&& + m^{L\underline{L}}m^{\underline{L}L}   \big(   \derm_L h_{\underline{L}\si} +\derm_\si h_{\underline{L} L}- \derm_{\underline{L}} h_{L \si}  \big)   \cdot  \big( \derm_{\underline{L}}A_{L} - \pa_{L}A_{\underline{L}}  \big) \\
&& + m^{L\underline{L}}m^{e_{A} e_{A}}   \big(  \derm_L h_{e_{A}\si} + \derm_\si h_{e_{A} L}- \derm_{e_{A}} h_{L \si}  \big)   \cdot   \big( \derm_{\underline{L}}A_{e_{A}} - \derm_{e_{A}}A_{\underline{L}}  \big) \\
 && + m^{\underline{L} L}m^{L\underline{L}}   \big(  \derm_{\underline{L}} h_{L\si} + \derm_\si h_{L\underline{L}}- \derm_L h_{\underline{L}\si}  \big)   \cdot  \big(\derm_{L}A_{\underline{L}} - \derm_{\underline{L}}A_{L}  \big) \\
  && + m^{\underline{L} L}m^{\underline{L}L}   \big(   \derm_{\underline{L}} h_{\underline{L}\si} + \derm_\si h_{\underline{L}\underline{L}}- \derm_{\underline{L}} h_{\underline{L}\si}  \big)   \cdot  \big( \derm_{L}A_{L} - \derm_{L}A_{L}  \big) \\
   && + m^{\underline{L} L}m^{e_{A}e_{A}}   \big(  \derm_{\underline{L}} h_{e_{A}\si} + \derm_\si h_{e_{A}\underline{L}}- \derm_{e_{A}} h_{\underline{L}\si}  \big)   \cdot \big( \derm_{L}A_{e_{A}} - \derm_{e_{A}}A_{L}  \big) \\
 && +  m^{e_{A}e_{A}}m^{L\underline{L}}   \big(  \derm_{e_{A}} h_{L\si} + \derm_\si h_{L e_{A}}- \pa_L h_{e_{A} \si}  \big)   \cdot  \big(\derm_{e_{A}}A_{\underline{L}} - \derm_{\underline{L}}A_{e_{A}}  \big)  \\
  && +  m^{e_{A}e_{A}}m^{\underline{L}L}   \big( \derm_{e_{A}} h_{\underline{L}\si} +\derm_\si h_{\underline{L} e_{A}}- \derm_{\underline{L}} h_{e_{A} \si}  \big)   \cdot   \big( \derm_{e_{A}}A_{L} - \derm_{L}A_{e_{A}}  \big)   \\
   && +  m^{e_{A}e_{A}}m^{e_{A}e_{A}}   \big(  \derm_{e_{A}} h_{e_{A}\si} + \derm_\si h_{e_{A} e_{A}}- \derm_{e_{A}} h_{e_{A} \si}  \big)   \cdot   \big( \derm_{e_{A}}A_{e_{A}} - \derm_{e_{A}}A_{e_{A}}  \big)   .
\eeaa

Similarly, for the third term, we get
   \beaa
&&     m^{\a\mu}m^{\b\nu}   \big(   \derm_\a h_{\b\si} + \derm_\si h_{\b\a}- \derm_\b h_{\a\si}  \big)   \cdot   [A_{\mu},A_{\nu}] \\
&=&  m^{L\underline{L}}m^{\underline{L}L}   \big(   \derm_L h_{\underline{L}\si} + \derm_\si h_{\underline{L} L}- \derm_{\underline{L}} h_{L \si}  \big)  \cdot   [A_{\underline{L}},A_{L}]    \\
&& + m^{L\underline{L}}m^{e_{A} e_{A}}   \big(   \derm_L h_{e_{A}\si} + \derm_\si h_{e_{A} L}- \derm_{e_{A}} h_{L \si}  \big)   \cdot   [A_{\underline{L}},A_{e_{A}}]   \\
 && + m^{\underline{L} L}m^{L\underline{L}}   \big(   \derm_{\underline{L}} h_{L\si} + \derm_\si h_{L\underline{L}}- \derm_L h_{\underline{L}\si}  \big)   \cdot   [A_{L},A_{\underline{L}}]    \\
   && + m^{\underline{L} L}m^{e_{A}e_{A}}   \big(   \derm_{\underline{L}} h_{e_{A}\si} + \derm_\si h_{e_{A}\underline{L}}- \derm_{e_{A}} h_{\underline{L}\si}  \big)   \cdot  [A_{L},A_{e_{A}}]    \\
 && +  m^{e_{A}e_{A}}m^{L\underline{L}}   \big(   \derm_{e_{A}} h_{L\si} + \derm_\si h_{L e_{A}}- \derm_L h_{e_{A} \si}  \big)   \cdot  [A_{e_{A}}, A_{\underline{L}}]    \\
  && +  m^{e_{A}e_{A}}m^{\underline{L}L}   \big(   \derm_{e_{A}} h_{\underline{L}\si} + \derm_\si h_{\underline{L} e_{A}}- \derm_{\underline{L}} h_{e_{A} \si}  \big)   \cdot  [A_{e_{A}}, A_{L}]    .
\eeaa

We compute the fourth term,
\beaa
 && -  m^{\a\mu} \big(  [ A_{\mu}, \derm_{\a} A_{\si} ]  +    [A_{\alpha}, \derm_{\mu}  A_{\si} - \derm_{\si} A_{\mu} ]    +    [A_{\alpha}, [A_{\mu},A_{\si}] ]  \big)  \\
   &=&  -  m^{L \mu} \big(  [ A_{\mu}, \derm_{L} A_{\si} ]  +   [A_{L},  \derm_{\mu}  A_{\si} - \derm_{\si} A_{\mu} ]    +    [A_{L}, [A_{\mu},A_{\si}] ]  \big)  \\
    &&  -  m^{\underline{L} \mu} \big(  [ A_{\mu}, \derm_{\underline{L}} A_{\si} ]  +    [A_{\underline{L}}, \derm_{\mu}  A_{\si} - \derm_{\si} A_{\mu} ]    +    [A_{\underline{L}}, [A_{\mu},A_{\si}] ]  \big)  \\
 &&   -  m^{e_{a} \mu} \big(  [ A_{\mu}, \derm_{e_{a}} A_{\si} ]  +    [A_{e_{a}},  \derm_{\mu}  A_{\si} - \derm_{\si} A_{\mu} ]    +    [A_{e_{a}}, [A_{\mu},A_{\si}] ]  \big)  \\
    &=&  -  m^{L \underline{L}} \big(  [ A_{\underline{L}}, \derm_{L} A_{\si} ]  +    [A_{L},  \derm_{\underline{L}}  A_{\si} - \derm_{\si} A_{\underline{L}} ]    +    [A_{L}, [A_{\underline{L}}, A_{\si}] ]  \big)  \\
    &&  -  m^{\underline{L} L} \big(  [ A_{L}, \derm_{\underline{L}} A_{\si} ]  +    [A_{\underline{L}},  \derm_{L}  A_{\si} -\derm_{\si} A_{L} ]    +    [A_{\underline{L}}, [A_{L},A_{\si}] ]  \big)  \\
 &&   -  m^{e_{a}e_{a}} \big(  [ A_{e_{a}}, \derm_{e_{a}} A_{\si} ]  +    [A_{e_{a}},  \derm_{e_{a}}  A_{\si} - \derm_{\si} A_{e_{a}} ]    +    [A_{e_{a}}, [A_{e_{a}},A_{\si}] ]  \big)  \\
    &=&  -  m^{L \underline{L}} \big(  [ A_{\underline{L}}, \derm_{L} A_{\si} ]  +    [A_{L},  \derm_{\underline{L}}  A_{\si} - \derm_{\si} A_{\underline{L}} ]    +    [A_{L}, [A_{\underline{L}}, A_{\si}] ]  \big)  \\
    &&  -  m^{\underline{L} L} \big(  [ A_{L}, \derm_{\underline{L}} A_{\si} ]  +   [A_{\underline{L}},  \derm_{L}  A_{\si} -\derm_{\si} A_{L} ]    +    [A_{\underline{L}}, [A_{L},A_{\si}] ]  \big)  \\
 &&   -  m^{e_{a}e_{a}} \big( [ A_{e_{a}}, \derm_{e_{a}} A_{\si} ]  +  [A_{e_{a}}, \derm_{e_{a}} A_{\si}  -  \derm_{\si} A_{e_{a}} ]    +    [A_{e_{a}}, [A_{e_{a}},A_{\si}] ]  \big) .
 \eeaa

 This leads to
     \beaa
   \notag
&& g^{\la\mu} \derm_{\la}   \derm_{\mu}   A_{\si}     \\
   \notag
    &=&  m^{L\underline{L}} m ^{L\underline{L}}  (  \derm_{\si}  h_{\underline{L}\underline{L}} )   \derm_{L}A_{L}  + m^{L\underline{L}} m ^{\underline{L}L} (   \derm_{\si}  h_{\underline{L}L} )   \derm_{L}A_{\underline{L}}   +m^{L\underline{L}} m ^{e_{A}e_{A}} (   \derm_{\si}  h_{\underline{L}e_{A}} )   \derm_{L} A_{e_{A}}  \\
   \notag
 && +  m^{\underline{L} L} m ^{L\underline{L}} (  \derm_{\si}  h_{L\underline{L}} )   \derm_{\underline{L} }A_{L} +  m^{\underline{L} L} m ^{\underline{L}L} (   \derm_{\si}  h_{LL} )   \derm_{\underline{L} }A_{\underline{L}} +  m^{\underline{L} L} m ^{e_{A}e_{A}} (  \derm_{\si}  h_{Le_{A}} )   \derm_{\underline{L} }A_{e_{A}} \\
    \notag
&&   +  m^{e_{A}e_{A}} m ^{L\underline{L}} (  \derm_{\si}  h_{e_{A}\underline{L}} )   \derm_{e_{A}}A_{L} +  m^{e_{A}e_{A}} m ^{\underline{L}L}   (  \derm_{\si}  h_{e_{A}L} )  \derm_{e_{A}}A_{\underline{L}} \\
\notag
&& +  m^{e_{A}e_{A}} m ^{e_{A}e_{A}}  \derm_{\si} ( h_{e_{A}e_{A}} )   \derm_{e_{A}}A_{e_{A}} \\
   \notag
&& + \frac{1}{2} m^{L\underline{L}}m^{L\underline{L}}   \big(  \derm_L h_{L\si} + \derm_\si h_{L L}- \derm_L h_{L \si}  \big)   \cdot   \big(\derm_{\underline{L}}A_{\underline{L}} - \derm_{\underline{L}}A_{\underline{L}}  \big) \\
   \notag
&& + \frac{1}{2} m^{L\underline{L}}m^{\underline{L}L}   \big(   \derm_L h_{\underline{L}\si} + \derm_\si h_{\underline{L} L}- \derm_{\underline{L}} h_{L \si}  \big)   \cdot  \big( \derm_{\underline{L}}A_{L} - \derm_{L}A_{\underline{L}}  \big) \\
   \notag
&& + \frac{1}{2} m^{L\underline{L}}m^{e_{A} e_{A}}   \big(   \derm_L h_{e_{A}\si} +\derm_\si h_{e_{A} L}-\derm_{e_{A}} h_{L \si}  \big)   \cdot   \big( \derm_{\underline{L}}A_{e_{A}} - \derm_{e_{A}}A_{\underline{L}}  \big) \\
   \notag
 && + \frac{1}{2} m^{\underline{L} L}m^{L\underline{L}}   \big(  \derm_{\underline{L}} h_{L\si} + \derm_\si h_{L\underline{L}}-\derm_L h_{\underline{L}\si}  \big)   \cdot  \big( \derm_{L}A_{\underline{L}} - \derm_{\underline{L}}A_{L}  \big) \\
    \notag
  && + \frac{1}{2} m^{\underline{L} L}m^{\underline{L}L}   \big( \derm_{\underline{L}} h_{\underline{L}\si} +\derm_\si h_{\underline{L}\underline{L}}-\derm_{\underline{L}} h_{\underline{L}\si}  \big)   \cdot  \big( \derm_{L}A_{L} -\derm_{L}A_{L}  \big) \\
     \notag
   && + \frac{1}{2} m^{\underline{L} L}m^{e_{A}e_{A}}   \big(  \derm_{\underline{L}} h_{e_{A}\si} +\derm_\si h_{e_{A}\underline{L}}- \derm_{e_{A}} h_{\underline{L}\si}  \big)   \cdot \big( \derm_{L}A_{e_{A}} - \derm_{e_{A}}A_{L}  \big) \\
      \notag
 && +  \frac{1}{2} m^{e_{A}e_{A}}m^{L\underline{L}}   \big(  \derm_{e_{A}} h_{L\si} + \derm_\si h_{L e_{A}}- \derm_L h_{e_{A} \si}  \big)   \cdot  \big( \derm_{e_{A}}A_{\underline{L}} - \derm_{\underline{L}}A_{e_{A}}  \big)  \\
    \notag
  && +  \frac{1}{2} m^{e_{A}e_{A}}m^{\underline{L}L}   \big(   \derm_{e_{A}} h_{\underline{L}\si} +\derm_\si h_{\underline{L} e_{A}}- \derm_{\underline{L}} h_{e_{A} \si}  \big)   \cdot   \big( \derm_{e_{A}}A_{L} -\derm_{L}A_{e_{A}}  \big)   \\
     \notag
   && +  \frac{1}{2} m^{e_{A}e_{A}}m^{e_{A}e_{A}}   \big(   \derm_{e_{A}} h_{e_{A}\si} + \derm_\si h_{e_{A} e_{A}}- \derm_{e_{A}} h_{e_{A} \si}  \big)   \cdot   \big( \derm_{e_{A}}A_{e_{A}} - \derm_{e_{A}}A_{e_{A}}  \big) \\
 \notag
&& \frac{1}{2} m^{L\underline{L}}m^{\underline{L}L}   \big(   \derm_L h_{\underline{L}\si} + \derm_\si h_{\underline{L} L}- \derm_{\underline{L}} h_{L \si}  \big)  \cdot   [A_{\underline{L}},A_{L}]    \\
   \notag
&& + \frac{1}{2} m^{L\underline{L}}m^{e_{A} e_{A}}   \big(  \derm_L h_{e_{A}\si} + \derm_\si h_{e_{A} L}- \derm_{e_{A}} h_{L \si}  \big)   \cdot   [A_{\underline{L}},A_{e_{A}}]   \\
   \notag
 && + \frac{1}{2} m^{\underline{L} L}m^{L\underline{L}}   \big(  \derm_{\underline{L}} h_{L\si} +\derm_\si h_{L\underline{L}}- \derm_L h_{\underline{L}\si}  \big)   \cdot   [A_{L},A_{\underline{L}}]    \\
    \notag
   && + \frac{1}{2} m^{\underline{L} L}m^{e_{A}e_{A}}   \big(   \derm_{\underline{L}} h_{e_{A}\si} + \derm_\si h_{e_{A}\underline{L}}- \derm_{e_{A}} h_{\underline{L}\si}  \big)   \cdot  [A_{L},A_{e_{A}}]    \\
      \notag
 && +  \frac{1}{2} m^{e_{A}e_{A}}m^{L\underline{L}}   \big(   \derm_{e_{A}} h_{L\si} + \derm_\si h_{L e_{A}}- \derm_L h_{e_{A} \si}  \big)   \cdot  [A_{e_{A}}, A_{\underline{L}}]    \\
    \notag
  && +  \frac{1}{2} m^{e_{A}e_{A}}m^{\underline{L}L}   \big(   \derm_{e_{A}} h_{\underline{L}\si} +\derm_\si h_{\underline{L} e_{A}}-\derm_{\underline{L}} h_{e_{A} \si}  \big)   \cdot  [A_{e_{A}}, A_{L}] \\
 \notag
 &&  -  m^{L \underline{L}} \big(  [ A_{\underline{L}}, \derm_{L} A_{\si} ]  +    [A_{L},  \derm_{\underline{L}}  A_{\si} - \derm_{\si} A_{\underline{L}} ]    +    [A_{L}, [A_{\underline{L}}, A_{\si}] ]  \big)  \\
    \notag
    &&  -  m^{\underline{L} L} \big(  [ A_{L}, \derm_{\underline{L}} A_{\si} ]  +    [A_{\underline{L}}, \derm_{L}  A_{\si} - \derm_{\si} A_{L} ]    +    [A_{\underline{L}}, [A_{L},A_{\si}] ]  \big)  \\
       \notag
 &&   -  m^{e_{a}e_{a}} \big(   [ A_{e_{a}}, \derm_{e_{a}} A_{\si} ] +  [A_{e_{a}},  \derm_{e_{a}} A_{\si}  - \derm_{\si} A_{e_{a}} ]    +    [A_{e_{a}}, [A_{e_{a}},A_{\si}] ]  \big) \\
  && + O( h \cdot  \pa h \cdot  \pa A) + O( h \cdot  \pa h \cdot  A^2) + O( h \cdot  A \cdot \pa A) + O( h \cdot  A^3) \, .
  \eeaa
  Finally, we get
    \bea
   \notag
&& g^{\la\mu} \derm_{\la}   \derm_{\mu}   A_{\si}   \\ 
  \notag
 &=&  \frac{1}{4} (  \derm_{\si}  h_{\underline{L}\underline{L}} )   \derm_{L}A_{L}  + \frac{1}{4} (   \derm_{\si}  h_{\underline{L}L} )   \derm_{L}A_{\underline{L}}   -\frac{1}{2}  (   \derm_{\si}  h_{\underline{L}e_{A}} )  \derm_{L} A_{e_{A}}  \\
   \notag
 && +  \frac{1}{4} (  \derm_{\si}  h_{L\underline{L}} )  \derm_{\underline{L} }A_{L} +  \frac{1}{4} (   \derm_{\si}  h_{LL} )   \derm_{\underline{L} }A_{\underline{L}} -  \frac{1}{2} (  \derm_{\si}  h_{Le_{A}} )  \derm_{\underline{L} }A_{e_{A}} \\
    \notag
&&   -  \frac{1}{2} (  \derm_{\si}  h_{e_{A}\underline{L}} )   \derm_{e_{A}}A_{L}   -  \frac{1}{2}  (  \derm_{\si}  h_{e_{A}L} )   \derm_{e_{A}}A_{\underline{L}} + (  \derm_{\si}  h_{e_{A}e_{A}} )   \derm_{e_{A}}A_{e_{A}} \\
   \notag
&& + \frac{1}{8}  \big(    \derm_L h_{\underline{L}\si} +  \derm_\si h_{\underline{L} L}-  \derm_{\underline{L}} h_{L \si}  \big)   \cdot  \big( \derm_{\underline{L}}A_{L} -  \derm_{L}A_{\underline{L}}  \big) \\
   \notag
&& - \frac{1}{4}    \big(    \derm_L h_{e_{A}\si} + \derm_\si h_{e_{A} L}-  \derm_{e_{A}} h_{L \si}  \big)   \cdot   \big(  \derm_{\underline{L}}A_{e_{A}} -  \derm_{e_{A}}A_{\underline{L}}  \big) \\
   \notag
 && + \frac{1}{8}    \big(    \derm_{\underline{L}} h_{L\si} +  \derm_\si h_{L\underline{L}}-  \derm_L h_{\underline{L}\si}  \big)   \cdot  \big(  \derm_{L}A_{\underline{L}} -  \derm_{\underline{L}}A_{L}  \big) \\
     \notag
   && - \frac{1}{4}   \big(    \derm_{\underline{L}} h_{e_{A}\si} +  \derm_\si h_{e_{A}\underline{L}}-  \derm_{e_{A}} h_{\underline{L}\si}  \big)   \cdot \big(  \derm_{L}A_{e_{A}} -  \derm_{e_{A}}A_{L}  \big) \\
      \notag
 && -  \frac{1}{4}    \big(    \derm_{e_{A}} h_{L\si} +  \derm_\si h_{L e_{A}}-  \derm_L h_{e_{A} \si}  \big)   \cdot  \big(  \derm_{e_{A}}A_{\underline{L}} -  \derm_{\underline{L}}A_{e_{A}}  \big)  \\
    \notag
  && -  \frac{1}{4}    \big(   \derm_{e_{A}} h_{\underline{L}\si} +  \derm_\si h_{\underline{L} e_{A}}-  \derm_{\underline{L}} h_{e_{A} \si}  \big)   \cdot   \big( \derm_{e_{A}}A_{L} -  \derm_{L}A_{e_{A}}  \big)   \\
 \notag
&& + \frac{1}{8}   \big(    \derm_L h_{\underline{L}\si} + \derm_\si h_{\underline{L} L}-  \derm_{\underline{L}} h_{L \si}  \big)  \cdot   [A_{\underline{L}},A_{L}]    \\
   \notag
&& - \frac{1}{4}   \big(    \derm_L h_{e_{A}\si} + \derm_\si h_{e_{A} L}-  \derm_{e_{A}} h_{L \si}  \big)   \cdot   [A_{\underline{L}},A_{e_{A}}]   \\
   \notag
 && + \frac{1}{8}   \big(    \derm_{\underline{L}} h_{L\si} + \derm_\si h_{L\underline{L}}-  \derm_L h_{\underline{L}\si}  \big)   \cdot   [A_{L},A_{\underline{L}}]    \\
    \notag
   && - \frac{1}{4}   \big(   \derm_{\underline{L}} h_{e_{A}\si} +  \derm_\si h_{e_{A}\underline{L}}- \derm_{e_{A}} h_{\underline{L}\si}  \big)   \cdot  [A_{L},A_{e_{A}}]    \\
      \notag
 && -  \frac{1}{4}    \big(    \derm_{e_{A}} h_{L\si} + \derm_\si h_{L e_{A}}-  \derm_L h_{e_{A} \si}  \big)   \cdot  [A_{e_{A}}, A_{\underline{L}}]    \\
    \notag
  && -  \frac{1}{4}   \big(    \derm_{e_{A}} h_{\underline{L}\si} +  \derm_\si h_{\underline{L} e_{A}}-  \derm_{\underline{L}} h_{e_{A} \si}  \big)   \cdot  [A_{e_{A}}, A_{L}] \\
 \notag
 && +  \frac{1}{2} \big(  [ A_{\underline{L}},  \derm_{L} A_{\si} ]  +   [A_{L},   \derm_{\underline{L}}  A_{\si} - \derm_{\si} A_{\underline{L}} ]    +    [A_{L}, [A_{\underline{L}}, A_{\si}] ]  \big)  \\
    \notag
    &&  +  \frac{1}{2} \big(  [ A_{L},  \derm_{\underline{L}} A_{\si} ]  +    [A_{\underline{L}},   \derm_{L}  A_{\si} - \derm_{\si} A_{L} ]    +    [A_{\underline{L}}, [A_{L},A_{\si}] ]  \big)  \\
       \notag
 &&   -   \big(  [ A_{e_{a}}, \derm_{e_{a}} A_{\si} ] + [A_{e_{a}},  \derm_{e_{a}} A_{\si}  -  \derm_{\si} A_{e_{a}} ]     +    [A_{e_{a}}, [A_{e_{a}},A_{\si}] ]  \big) \\
    \notag
  && + O( h \cdot  \pa h \cdot  \pa A) + O( h \cdot  \pa h \cdot  A^2) + O( h \cdot  A \cdot \pa A) + O( h \cdot  A^3) \, . \\
  \eea

Combining further the terms that repeat and cancelling out some terms using the symmetry of the metric and therefore of the tensor $h$, we obtain the stated result for the wave equation on the Yang-Mills potential.

On the other hand, we showed in Lemma \ref{Einstein-Yang-MillssysteminLorenzwavegauges}, that in wave coordinates, the metric solution to the Einstein-Yang-Mills equations in the Lorenz gauge satisfies the following equation,

  \beaa
\notag
 && g^{\alpha\beta}\pa_\alpha\pa_\beta h_{\mu\nu} \\
  &=& P(\pa_\mu h,\pa_\nu h)  +  Q_{\mu\nu}(\pa h,\pa h)   + G_{\mu\nu}(h)(\pa h,\pa h)  \\
\notag
 &&   -4   m^{\si\b} \cdot  <   \derm_{\mu}A_{\b} - \derm_{\b}A_{\mu}  ,  \derm_{\nu}A_{\si} - \derm_{\si}A_{\nu}  >    \\
 \notag
 &&   + m_{\mu\nu }  m^{\si\b}  m^{\a\la}    \cdot  <  \derm_{\a}A_{\b} - \derm_{\b}A_{\a} , \derm_{\la}A_{\si} - \derm_{\si}A_{\la} >   \\
 \notag
&&           -4 m^{\si\b}  \cdot  \big( <   \derm_{\mu}A_{\b} - \derm_{\b}A_{\mu}  ,  [A_{\nu},A_{\si}] >   + <   [A_{\mu},A_{\b}] ,  \derm_{\nu}A_{\si} - \derm_{\si}A_{\nu}  > \big)  \\
\notag
&& + m_{\mu\nu }  m^{\si\b}  m^{\a\la}    \cdot \big(  <  \derm_{\a}A_{\b} - \derm_{\b}A_{\a} , [A_{\la},A_{\si}] >    +  <  [A_{\a},A_{\b}] , \derm_{\la}A_{\si} - \derm_{\si}A_{\la}  > \big) \\
\notag
 &&  -4 m^{\si\b}  \cdot   <   [A_{\mu},A_{\b}] ,  [A_{\nu},A_{\si}] >      + m_{\mu\nu }  m^{\si\b}  m^{\a\la}   \cdot   <  [A_{\a},A_{\b}] , [A_{\la},A_{\si}] >  \\
     && + O \big(h \cdot  (\pa A)^2 \big)   + O \big(  h  \cdot  A^2 \cdot \pa A \big)     + O \big(  h   \cdot  A^4 \big)  \,  .
\eeaa

 We compute the source terms by decomposing the terms in the null-frame $\{ \underline{L}, L, e_1, \ldots, e_{(n-1)} \}$. Computing the first term
  \beaa
\notag
 &&   -4   m^{\si\b} \cdot  <   \derm_{\mu}A_{\b} - \derm_{\b}A_{\mu}  ,  \derm_{\nu}A_{\si} - \derm_{\si}A_{\nu}  >    \\
 &=& -4  m^{\si L} \cdot  <   \derm_{\mu}A_{L} - \derm_{L}A_{\mu}  ,  \derm_{\nu}A_{\si} - \derm_{\si}A_{\nu}  >  \\
 &&  -4   m^{\si  \underline{L} } \cdot  <   \derm_{\mu}A_{ \underline{L} } - \derm_{ \underline{L} }A_{\mu}  ,  \derm_{\nu}A_{\si} - \derm_{\si}A_{\nu}  >    \\
 && -4   m^{\si e_A} \cdot  <   \derm_{\mu}A_{e_A} - \derm_{e_A }A_{\mu}  ,  \derm_{\nu}A_{\si} - \derm_{\si}A_{\nu}  >  \\
 &=&  -4  m^{ \underline{L}  L} \cdot  <   \derm_{\mu}A_{L} - \derm_{L}A_{\mu}  ,  \derm_{\nu}A_{ \underline{L} } - \derm_{ \underline{L} }A_{\nu}  >   \\
 &&  -4   m^{L  \underline{L} } \cdot  <   \derm_{\mu}A_{ \underline{L} } - \derm_{ \underline{L} }A_{\mu}  ,  \derm_{\nu}A_{L} - \derm_{L}A_{\nu}  >    \\
  && -4   m^{e_A e_A} \cdot  <   \derm_{\mu}A_{e_A} - \derm_{e_A }A_{\mu}  ,  \derm_{\nu}A_{e_A} - \pa_{e_A}A_{\nu}  >  .
   \eeaa
For the second term, we compute
 \beaa
 \notag
 &&     m^{\si\b}  m^{\a\la}    \cdot  < \derm_{\a}A_{\b} - \derm_{\b}A_{\a} , \derm_{\la}A_{\si} - \derm_{\si}A_{\la} >   \\
 \notag
&=& m^{L \b}  m^{\a\la}    \cdot  <  \derm_{\a}A_{\b} - \derm_{\b}A_{\a} ,\derm_{\la}A_{L} - \derm_{L}A_{\la} > \\
&&   + m^{\underline{L} \b}  m^{\a\la}    \cdot  <  \derm_{\a}A_{\b} -\derm_{\b}A_{\a} , \derm_{\la}A_{\underline{L}} - \derm_{\underline{L}}A_{\la} >  \\
&& + m^{e_A \b}  m^{\a\la}    \cdot  <  \derm_{\a}A_{\b} - \derm_{\b}A_{\a} , \derm_{\la}A_{e_A} -\derm_{e_A}A_{\la} > \\
&=& m^{L \underline{L} }  m^{\a\la}    \cdot  <  \derm_{\a}A_{\underline{L}} - \pa_{\underline{L}}A_{\a} , \derm_{\la}A_{L} - \derm_{L}A_{\la} >  \\
&& + m^{\underline{L} L}  m^{\a\la}    \cdot  <  \pa_{\a}A_{L} - \pa_{L}A_{\a} , \pa_{\la}A_{\underline{L}} - \pa_{\underline{L}}A_{\la} >  \\
&& + m^{e_A e_A}  m^{\a\la}    \cdot  <  \derm_{\a}A_{e_A} - \derm_{e_A}A_{\a} , \derm_{\la}A_{e_A} - \derm_{e_A}A_{\la} > \\
&=& m^{L \underline{L} }  m^{L\la}    \cdot  <  \derm_{L}A_{\underline{L}} - \derm_{\underline{L}}A_{L} , \derm_{\la}A_{L} - \derm_{L}A_{\la} > \\
&& + m^{L \underline{L} }  m^{\underline{L}\la}    \cdot  <  \derm_{\underline{L}}A_{\underline{L}} - \derm_{\underline{L}}A_{\underline{L}} ,\derm_{\la}A_{L} - \derm_{L}A_{\la} >  \\
&& + m^{L \underline{L} }  m^{e_A\la}    \cdot  <  \pa_{e_A}A_{\underline{L}} - \pa_{\underline{L}}A_{e_A} , \pa_{\la}A_{L} - \pa_{L}A_{\la} > \\
&& + m^{\underline{L} L}  m^{L\la}    \cdot  <  \derm_{L}A_{L} - \derm_{L}A_{L} , \derm_{\la}A_{\underline{L}} - \derm_{\underline{L}}A_{\la} >   \\
&&+ m^{\underline{L} L}  m^{\underline{L}\la}    \cdot  <  \pa_{\underline{L}}A_{L} - \pa_{L}A_{\underline{L}} , \pa_{\la}A_{\underline{L}} - \pa_{\underline{L}}A_{\la} > \\
&& +  m^{\underline{L} L}  m^{e_A\la}    \cdot  <  \derm_{e_A}A_{L} - \derm_{L}A_{e_A} , \derm_{\la}A_{\underline{L}} -\derm_{\underline{L}}A_{\la} > \\
&& + m^{e_A e_A}  m^{L\la}    \cdot  <  \pa_{L}A_{e_A} - \pa_{e_A}A_{L} , \derm_{\la}A_{e_A} - \derm_{e_A}A_{\la} > \\
&& +  m^{e_A e_A}  m^{\underline{L}\la}    \cdot  <  \derm_{\underline{L}}A_{e_A} -\derm_{e_A}A_{\underline{L}} , \derm_{\la}A_{e_A} - \derm_{e_A}A_{\la} > \\
&& + m^{e_A e_A}  m^{e_A\la}    \cdot  <  \derm_{e_A}A_{e_A} - \derm_{e_A}A_{e_A} , \derm_{\la}A_{e_A} - \derm_{e_A}A_{\la} >  \\
&=& m^{L \underline{L} }  m^{L\underline{L}}    \cdot  < \derm_{L}A_{\underline{L}} - \derm_{\underline{L}}A_{L} , \derm_{\underline{L}}A_{L} - \derm_{L}A_{\underline{L}} > \\
&&+ m^{L \underline{L} }  m^{\underline{L}L}    \cdot  <  \derm_{\underline{L}}A_{\underline{L}} -\derm_{\underline{L}}A_{\underline{L}} , \derm_{L}A_{L} - \derm_{L}A_{L} >  \\
&& + m^{L \underline{L} }  m^{e_A e_A}    \cdot  <  \derm_{e_A}A_{\underline{L}} -\derm_{\underline{L}}A_{e_A} , \derm_{e_A}A_{L} - \derm_{L}A_{e_A} > \\
&& + m^{\underline{L} L}  m^{L\underline{L}}    \cdot  <  \derm_{L}A_{L} - \derm_{L}A_{L} , \derm_{\underline{L}}A_{\underline{L}} - \derm_{\underline{L}}A_{\underline{L}} >  \\
&& + m^{\underline{L} L}  m^{\underline{L}L}    \cdot  <  \derm_{\underline{L}}A_{L} - \derm_{L}A_{\underline{L}} , \derm_{L}A_{\underline{L}} - \derm_{\underline{L}}A_{L} > \\
&& +  m^{\underline{L} L}  m^{e_A e_A}    \cdot  <  \derm_{e_A}A_{L} - \derm_{L}A_{e_A} , \derm_{e_A}A_{\underline{L}} - \derm_{\underline{L}}A_{e_A} > \\
&& + m^{e_A e_A}  m^{L\underline{L}}    \cdot  <  \derm_{L}A_{e_A} - \derm_{e_A}A_{L} , \derm_{\underline{L}}A_{e_A} - \derm_{e_A}A_{\underline{L}} > \\
&&+  m^{e_A e_A}  m^{\underline{L}L}    \cdot  <  \derm_{\underline{L}}A_{e_A} - \derm_{e_A}A_{\underline{L}} , \derm_{L}A_{e_A} - \derm_{e_A}A_{L} > \\
&& + m^{e_A e_A}  m^{e_A e_A}    \cdot  <  \derm_{e_A}A_{e_A} - \derm_{e_A}A_{e_A} , \derm_{e_A}A_{e_A} - \derm_{e_A}A_{e_A} >   .
\eeaa
For the third term, we compute

\beaa
&&           m^{\si\b}  \cdot  \big( <   \derm_{\mu}A_{\b} - \derm_{\b}A_{\mu}  ,  [A_{\nu},A_{\si}] >   + <   [A_{\mu},A_{\b}] ,  \derm_{\nu}A_{\si} - \derm_{\si}A_{\nu}  > \big)  \\
&=& m^{L \b}  \cdot  \big( <   \derm_{\mu}A_{\b} - \derm_{\b}A_{\mu}  ,  [A_{\nu},A_{L}] >   + <   [A_{\mu},A_{\b}] ,  \derm_{\nu}A_{L} -  \derm_{L}A_{\nu}  > \big) \\
&&  +  m^{\underline{L}\b}  \cdot  \big( <   \derm_{\mu}A_{\b} - \derm_{\b}A_{\mu}  ,  [A_{\nu},A_{\underline{L}}] >   + <   [A_{\mu},A_{\b}] ,  \derm_{\nu}A_{\underline{L}} - \derm_{\underline{L}}A_{\nu}  > \big)  \\
 && + m^{e_A\b}  \cdot  \big( <  \derm_{\mu}A_{\b} - \derm_{\b}A_{\mu}  ,  [A_{\nu},A_{e_A}] >   + <   [A_{\mu},A_{\b}] , \derm_{\nu}A_{e_A} - \derm_{e_A}A_{\nu}  > \big)  \\
 &=& m^{L \underline{L}}  \cdot  \big( <   \derm_{\mu}A_{\underline{L}} - \derm_{\underline{L}}A_{\mu}  ,  [A_{\nu},A_{L}] >   + <   [A_{\mu},A_{\underline{L}}] ,  \derm_{\nu}A_{L} - \derm_{L}A_{\nu}  > \big) \\
&&  +  m^{\underline{L}L}  \cdot  \big( <  \derm_{\mu}A_{L} -\derm_{L}A_{\mu}  ,  [A_{\nu},A_{\underline{L}}] >   + <   [A_{\mu},A_{L}] , \derm_{\nu}A_{\underline{L}} -\derm_{\underline{L}}A_{\nu}  > \big)  \\
 && + m^{e_A e_A}  \cdot  \big( <   \derm_{\mu}A_{e_A} - \derm_{e_A}A_{\mu}  ,  [A_{\nu},A_{e_A}] >   + <   [A_{\mu},A_{e_A}] ,  \derm_{\nu}A_{e_A} - \derm_{e_A}A_{\nu}  > \big)  .
\eeaa
For the forth term, 

\beaa
&&   m^{\si\b}  m^{\a\la}    \cdot \big(  <  \derm_{\a}A_{\b} - \derm_{\b}A_{\a} , [A_{\la},A_{\si}] >    +  <  [A_{\a},A_{\b}] , \derm_{\la}A_{\si} - \derm_{\si}A_{\la}  > \big) \\
&=&   m^{L\b}  m^{\a\la}    \cdot \big(  <  \derm_{\a}A_{\b} - \derm_{\b}A_{\a} , [A_{\la},A_{L}] >    +  <  [A_{\a},A_{\b}] , \derm_{\la}A_{L} - \derm_{L}A_{\la}  > \big) \\
&& +   m^{\underline{L}\b}  m^{\a\la}    \cdot \big(  <  \derm_{\a}A_{\b} - \derm_{\b}A_{\a} , [A_{\la},A_{\underline{L}}] >    +  <  [A_{\a},A_{\b}] , \derm_{\la}A_{\underline{L}} - \derm_{\underline{L}}A_{\la}  > \big) \\
&& +   m^{e_A \b}  m^{\a\la}    \cdot \big(  <  \derm_{\a}A_{\b} - \derm_{\b}A_{\a} , [A_{\la},A_{e_A}] >    +  <  [A_{\a},A_{\b}] , \derm_{\la}A_{e_A} - \derm_{e_A}A_{\la}  > \big) \\
&=&   m^{L\underline{L}}  m^{\a\la}    \cdot \big(  <  \derm_{\a}A_{\underline{L}} - \derm_{\underline{L}}A_{\a} , [A_{\la},A_{L}] >    +  <  [A_{\a},A_{\underline{L}}] , \derm_{\la}A_{L} - \derm_{L}A_{\la}  > \big) \\
&& +   m^{\underline{L}L}  m^{\a\la}    \cdot \big(  <  \derm_{\a}A_{L} -\derm_{L}A_{\a} , [A_{\la},A_{\underline{L}}] >    +  <  [A_{\a},A_{L}] , \derm_{\la}A_{\underline{L}} - \derm_{\underline{L}}A_{\la}  > \big) \\
&& +   m^{e_A e_A}  m^{\a\la}    \cdot \big(  <  \derm_{\a}A_{e_A} - \derm_{e_A}A_{\a} , [A_{\la},A_{e_A}] >    +  <  [A_{\a},A_{e_A}] , \derm_{\la}A_{e_A} - \derm_{e_A}A_{\la}  > \big) \\
&=&   m^{L\underline{L}}  m^{L \la}    \cdot \big(  <  \derm_{L}A_{\underline{L}} - \derm_{\underline{L}}A_{L} , [A_{\la},A_{L}] >    +  <  [A_{L},A_{\underline{L}}] ,\derm_{\la}A_{L} - \derm_{L}A_{\la}  > \big) \\
&& + m^{L\underline{L}}  m^{\underline{L}\la}    \cdot \big(  <  \derm_{\underline{L}}A_{\underline{L}} - \derm_{\underline{L}}A_{\underline{L}} , [A_{\la},A_{L}] >    +  <  [A_{\underline{L}},A_{\underline{L}}] , \derm_{\la}A_{L} - \derm_{L}A_{\la}  > \big) \\
 && + m^{L\underline{L}}  m^{e_A\la}    \cdot \big(  < \derm_{e_A }A_{\underline{L}} - \derm_{\underline{L}}A_{e_A} , [A_{\la},A_{L}] >    +  <  [A_{e_A},A_{\underline{L}}] , \derm_{\la}A_{L} - \pa_{L}A_{\la}  > \big) \\
 && +   m^{\underline{L}L}  m^{L\la}    \cdot \big(  <  \derm_{L}A_{L} - \derm_{L}A_{L} , [A_{\la},A_{\underline{L}}] >    +  <  [A_{L},A_{L}] , \derm_{\la}A_{\underline{L}} -\derm_{\underline{L}}A_{\la}  > \big) \\
 && +   m^{\underline{L}L}  m^{\underline{L}\la}    \cdot \big(  <  \derm_{\underline{L}}A_{L} -\derm_{L}A_{\underline{L}} , [A_{\la},A_{\underline{L}}] >    +  <  [A_{\underline{L}},A_{L}] , \derm_{\la}A_{\underline{L}} - \derm_{\underline{L}}A_{\la}  > \big) \\
 && +   m^{\underline{L}L}  m^{e_A \la}    \cdot \big(  < \derm_{e_A }A_{L} - \derm_{L}A_{e_A } , [A_{\la},A_{\underline{L}}] >    +  <  [A_{e_A },A_{L}] , \derm_{\la}A_{\underline{L}} - \derm_{\underline{L}}A_{\la}  > \big) \\
  && +   m^{e_A e_A}  m^{L\la}    \cdot \big(  <  \derm_{L}A_{e_A} - \derm_{e_A}A_{L} , [A_{\la},A_{e_A}] >    +  <  [A_{L},A_{e_A}] , \derm_{\la}A_{e_A} - \derm_{e_A}A_{\la}  > \big) \\
 && +   m^{e_A e_A}  m^{\underline{L}\la}    \cdot \big(  < \derm_{\underline{L}}A_{e_A} - \derm_{e_A}A_{\underline{L}} , [A_{\la},A_{e_A}] >    +  <  [A_{\underline{L}},A_{e_A}] , \derm_{\la}A_{e_A} - \derm_{e_A}A_{\la}  > \big) \\
 && +   m^{e_A e_A}  m^{e_A\la}    \cdot \big(  <  \derm_{e_A}A_{e_A} - \derm_{e_A}A_{e_A} , [A_{\la},A_{e_A}] >    +  <  [A_{e_A},A_{e_A}] , \derm_{\la}A_{e_A} - \derm_{e_A}A_{\la}  > \big) ,
\eeaa
and hence the fourth terms gives
\beaa
&&   m^{\si\b}  m^{\a\la}    \cdot \big(  <  \derm_{\a}A_{\b} - \derm_{\b}A_{\a} , [A_{\la},A_{\si}] >    +  <  [A_{\a},A_{\b}] , \derm_{\la}A_{\si} - \derm_{\si}A_{\la}  > \big) \\
&=&   m^{L\underline{L}}  m^{L \underline{L}}    \cdot \big(  <  \derm_{L}A_{\underline{L}} - \derm_{\underline{L}}A_{L} , [A_{\underline{L}},A_{L}] >    +  <  [A_{L},A_{\underline{L}}] , \derm_{\underline{L}}A_{L} - \derm_{L}A_{\underline{L}}  > \big) \\
&& + m^{L\underline{L}}  m^{\underline{L}L}    \cdot \big(  <  \derm_{\underline{L}}A_{\underline{L}} - \derm_{\underline{L}}A_{\underline{L}} , [A_{L},A_{L}] >    +  <  [A_{\underline{L}},A_{\underline{L}}] , \derm_{L}A_{L} - \derm_{L}A_{L}  > \big) \\
 && + m^{L\underline{L}}  m^{e_A e_A}    \cdot \big(  <  \derm_{e_A }A_{\underline{L}} - \derm_{\underline{L}}A_{e_A} , [A_{e_A},A_{L}] >    +  <  [A_{e_A},A_{\underline{L}}] , \derm_{e_A}A_{L} - \derm_{L}A_{e_A}  > \big) \\
 && +   m^{\underline{L}L}  m^{L\underline{L}}    \cdot \big(  <  \derm_{L}A_{L} - \derm_{L}A_{L} , [A_{\underline{L}},A_{\underline{L}}] >    +  <  [A_{L},A_{L}] ,\derm_{\underline{L}}A_{\underline{L}} - \derm_{\underline{L}}A_{\underline{L}}  > \big) \\
 && +   m^{\underline{L}L}  m^{\underline{L}L}    \cdot \big(  < \derm_{\underline{L}}A_{L} -\derm_{L}A_{\underline{L}} , [A_{L},A_{\underline{L}}] >    +  <  [A_{\underline{L}},A_{L}] , \derm_{L}A_{\underline{L}} - \derm_{\underline{L}}A_{L}  > \big) \\
 && +   m^{\underline{L}L}  m^{e_A e_A}    \cdot \big(  < \derm_{e_A }A_{L} - \derm_{L}A_{e_A } , [A_{e_A},A_{\underline{L}}] >    +  <  [A_{e_A },A_{L}] , \derm_{e_A}A_{\underline{L}} - \derm_{\underline{L}}A_{e_A}  > \big) \\
  && +   m^{e_A e_A}  m^{L\underline{L}}    \cdot \big(  <  \derm_{L}A_{e_A} - \derm_{e_A}A_{L} , [A_{\underline{L}},A_{e_A}] >    +  <  [A_{L},A_{e_A}] , \derm_{\underline{L}}A_{e_A} - \derm_{e_A}A_{\underline{L}}  > \big) \\
 && +   m^{e_A e_A}  m^{\underline{L}L}    \cdot \big(  <  \derm_{\underline{L}}A_{e_A} - \derm_{e_A}A_{\underline{L}} , [A_{L},A_{e_A}] >    +  <  [A_{\underline{L}},A_{e_A}] , \derm_{L}A_{e_A} - \derm_{e_A}A_{L}  > \big) \\
 && +   m^{e_A e_A}  m^{e_A e_A}    \cdot \big(  <  \derm_{e_A}A_{e_A} - \derm_{e_A}A_{e_A} , [A_{e_A},A_{e_A}] >  \\
 &&  +  <  [A_{e_A},A_{e_A}] , \derm_{e_A}A_{e_A} - \derm_{e_A}A_{e_A}  > \big) \\
 &=&   m^{L\underline{L}}  m^{L \underline{L}}    \cdot \big(  <  \derm_{L}A_{\underline{L}} - \derm_{\underline{L}}A_{L} , [A_{\underline{L}},A_{L}] >    +  <  [A_{L},A_{\underline{L}}] , \derm_{\underline{L}}A_{L} - \derm_{L}A_{\underline{L}}  > \big) \\
 && + m^{L\underline{L}}  m^{e_A e_A}    \cdot \big(  <  \derm_{e_A }A_{\underline{L}} - \derm_{\underline{L}}A_{e_A} , [A_{e_A},A_{L}] >    +  <  [A_{e_A},A_{\underline{L}}] , \derm_{e_A}A_{L} - \derm_{L}A_{e_A}  > \big) \\
 && +   m^{\underline{L}L}  m^{\underline{L}L}    \cdot \big(  < \derm_{\underline{L}}A_{L} -\derm_{L}A_{\underline{L}} , [A_{L},A_{\underline{L}}] >    +  <  [A_{\underline{L}},A_{L}] , \derm_{L}A_{\underline{L}} - \derm_{\underline{L}}A_{L}  > \big) \\
 && +   m^{\underline{L}L}  m^{e_A e_A}    \cdot \big(  < \derm_{e_A }A_{L} - \derm_{L}A_{e_A } , [A_{e_A},A_{\underline{L}}] >    +  <  [A_{e_A },A_{L}] , \derm_{e_A}A_{\underline{L}} - \derm_{\underline{L}}A_{e_A}  > \big) \\
  && +   m^{e_A e_A}  m^{L\underline{L}}    \cdot \big(  <  \derm_{L}A_{e_A} - \derm_{e_A}A_{L} , [A_{\underline{L}},A_{e_A}] >    +  <  [A_{L},A_{e_A}] , \derm_{\underline{L}}A_{e_A} - \derm_{e_A}A_{\underline{L}}  > \big) \\
 && +   m^{e_A e_A}  m^{\underline{L}L}    \cdot \big(  <  \derm_{\underline{L}}A_{e_A} - \derm_{e_A}A_{\underline{L}} , [A_{L},A_{e_A}] >    +  <  [A_{\underline{L}},A_{e_A}] , \derm_{L}A_{e_A} - \derm_{e_A}A_{L}  > \big) .
\eeaa
For the fifth term, we compute

\beaa
 &&   m^{\si\b}  \cdot   <   [A_{\mu},A_{\b}] ,  [A_{\nu},A_{\si}] >     \\
 &=& m^{L\b}  \cdot   <   [A_{\mu},A_{\b}] ,  [A_{\nu},A_{L}] >  + m^{\underline{L}\b}  \cdot   <   [A_{\mu},A_{\b}] ,  [A_{\nu},A_{\underline{L}}] > + m^{e_A\b}  \cdot   <   [A_{\mu},A_{\b}] ,  [A_{\nu},A_{e_A}] >  \\
 &=& m^{L\underline{L}}  \cdot   <   [A_{\mu},A_{\underline{L}}] ,  [A_{\nu},A_{L}] >  + m^{\underline{L}L}  \cdot   <   [A_{\mu},A_{L}] ,  [A_{\nu},A_{\underline{L}}] > + m^{e_A e_A}  \cdot   <   [A_{\mu},A_{e_A}] ,  [A_{\nu},A_{e_A}] > .
 \eeaa
For the sixth term, we evaluate
 \beaa
 &&  m^{\si\b}  m^{\a\la}   \cdot   <  [A_{\a},A_{\b}] , [A_{\la},A_{\si}] > \\
 &=&  m^{L\b}  m^{\a\la}   \cdot   <  [A_{\a},A_{\b}] , [A_{\la},A_{L}] > + m^{\underline{L}\b}  m^{\a\la}   \cdot   <  [A_{\a},A_{\b}] , [A_{\la},A_{\underline{L}}] > \\
 && + m^{e_A \b}  m^{\a\la}   \cdot   <  [A_{\a},A_{\b}] , [A_{\la},A_{e_A}] > \\
 &=&  m^{L\underline{L}}  m^{\a\la}   \cdot   <  [A_{\a},A_{\underline{L}}] , [A_{\la},A_{L}] > + m^{\underline{L}L}  m^{\a\la}   \cdot   <  [A_{\a},A_{L}] , [A_{\la},A_{\underline{L}}] > \\
 && + m^{e_A e_A}  m^{\a\la}   \cdot   <  [A_{\a},A_{e_A}] , [A_{\la},A_{e_A}] > \\
  &=&  m^{L\underline{L}}  m^{L\la}   \cdot   <  [A_{L},A_{\underline{L}}] , [A_{\la},A_{L}] > +  m^{L\underline{L}}  m^{\underline{L}\la}   \cdot   <  [A_{\underline{L}},A_{\underline{L}}] , [A_{\la},A_{L}] > \\
  && + m^{L\underline{L}}  m^{e_A\la}   \cdot   <  [A_{e_A},A_{\underline{L}}] , [A_{\la},A_{L}] > \\
&&   + m^{\underline{L}L}  m^{L\la}   \cdot   <  [A_{L},A_{L}] , [A_{\la},A_{\underline{L}}] > + m^{\underline{L}L}  m^{\underline{L}\la}   \cdot   <  [A_{\underline{L}},A_{L}] , [A_{\la},A_{\underline{L}}] > \\
&& + m^{\underline{L}L}  m^{e_A\la}   \cdot   <  [A_{e_A},A_{L}] , [A_{\la},A_{\underline{L}}] > \\
 && + m^{e_A e_A}  m^{L\la}   \cdot   <  [A_{L},A_{e_A}] , [A_{\la},A_{e_A}] > +  m^{e_A e_A}  m^{\underline{L}\la}   \cdot   <  [A_{\underline{L}},A_{e_A}] , [A_{\la},A_{e_A}] > \\
 && + m^{e_A e_A}  m^{e_A\la}   \cdot   <  [A_{e_A},A_{e_A}] , [A_{\la},A_{e_A}] > \\
  &=&  m^{L\underline{L}}  m^{L\underline{L}}   \cdot   <  [A_{L},A_{\underline{L}}] , [A_{\underline{L}},A_{L}] >  + m^{L\underline{L}}  m^{e_A e_A}   \cdot   <  [A_{e_A},A_{\underline{L}}] , [A_{e_A},A_{L}] > \\
&&  + m^{\underline{L}L}  m^{\underline{L}L}   \cdot   <  [A_{\underline{L}},A_{L}] , [A_{L},A_{\underline{L}}] >  + m^{\underline{L}L}  m^{e_A e_A}   \cdot   <  [A_{e_A},A_{L}] , [A_{e_A},A_{\underline{L}}] > \\
 && + m^{e_A e_A}  m^{L\underline{L}}   \cdot   <  [A_{L},A_{e_A}] , [A_{\underline{L}},A_{e_A}] > +  m^{e_A e_A}  m^{\underline{L} L}   \cdot   <  [A_{\underline{L}},A_{e_A}] , [A_{L},A_{e_A}] > .
  \eeaa

Thus, we obtain
    \beaa
\notag
&& g^{\alpha\beta}\derm_\alpha \derm_\beta h_{\mu\nu}  \\
\notag
     &=&  P(\pa_\mu h,\pa_\nu h) +  Q_{\mu\nu}(\pa h,\pa h)   + G_{\mu\nu}(h)(\pa h,\pa h)  \\
     \notag
 &&  -4  m^{ \underline{L}  L} \cdot  <  \derm_{\mu}A_{L} - \derm_{L}A_{\mu}  , \derm_{\nu}A_{ \underline{L} } - \derm_{ \underline{L} }A_{\nu}  >   \\
 \notag
 &&  -4   m^{L  \underline{L} } \cdot  <   \derm_{\mu}A_{ \underline{L} } - \derm_{ \underline{L} }A_{\mu}  , \derm_{\nu}A_{L} - \derm_{L}A_{\nu}  >    \\
 \notag
  && -4   m^{e_A e_A} \cdot  <   \derm_{\mu}A_{e_A} - \derm_{e_A }A_{\mu}  ,  \derm_{\nu}A_{e_A} -\derm_{e_A}A_{\nu}  >  \\
 \notag
 && + m_{\mu\nu }  m^{L \underline{L} }  m^{L\underline{L}}    \cdot  < \derm_{L}A_{\underline{L}} - \derm_{\underline{L}}A_{L} , \derm_{\underline{L}}A_{L} - \derm_{L}A_{\underline{L}} > \\
 \notag
&&+ m_{\mu\nu }  m^{L \underline{L} }  m^{\underline{L}L}    \cdot  <  \derm_{\underline{L}}A_{\underline{L}} - \derm_{\underline{L}}A_{\underline{L}} , \derm_{L}A_{L} - \derm_{L}A_{L} >  \\
\notag
&& + m_{\mu\nu }  m^{L \underline{L} }  m^{e_A e_A}    \cdot  < \derm_{e_A}A_{\underline{L}} - \derm_{\underline{L}}A_{e_A} , \derm_{e_A}A_{L} -\derm_{L}A_{e_A} > \\
\notag
&& + m_{\mu\nu }  m^{\underline{L} L}  m^{L\underline{L}}    \cdot  <  \derm_{L}A_{L} - \derm_{L}A_{L} , \derm_{\underline{L}}A_{\underline{L}} - \derm_{\underline{L}}A_{\underline{L}} >  \\
\notag
&& + m_{\mu\nu }  m^{\underline{L} L}  m^{\underline{L}L}    \cdot  <  \derm_{\underline{L}}A_{L} -\derm_{L}A_{\underline{L}} , \derm_{L}A_{\underline{L}} -\derm_{\underline{L}}A_{L} > \\
\notag
&& +  m_{\mu\nu }  m^{\underline{L} L}  m^{e_A e_A}    \cdot  <  \derm_{e_A}A_{L} -\derm_{L}A_{e_A} , \derm_{e_A}A_{\underline{L}} - \derm_{\underline{L}}A_{e_A} > \\
\notag
&& + m_{\mu\nu }  m^{e_A e_A}  m^{L\underline{L}}    \cdot  <  \derm_{L}A_{e_A} - \derm_{e_A}A_{L} , \derm_{\underline{L}}A_{e_A} - \derm_{e_A}A_{\underline{L}} > \\
\notag
&&+  m_{\mu\nu }  m^{e_A e_A}  m^{\underline{L}L}    \cdot  <  \derm_{\underline{L}}A_{e_A} - \derm_{e_A}A_{\underline{L}} , \derm_{L}A_{e_A} - \derm_{e_A}A_{L} > \\
\notag
&& + m_{\mu\nu }  m^{e_A e_A}  m^{e_A e_A}    \cdot  <  \derm_{e_A}A_{e_A} -\derm_{e_A}A_{e_A} , \derm_{e_A}A_{e_A} - \derm_{e_A}A_{e_A} >  \\
  \notag
 &&  -4  m^{L \underline{L}}  \cdot  \big( <   \derm_{\mu}A_{\underline{L}} - \derm_{\underline{L}}A_{\mu}  ,  [A_{\nu},A_{L}] >   + <   [A_{\mu},A_{\underline{L}}] ,  \derm_{\nu}A_{L} - \derm_{L}A_{\nu}  > \big) \\
 \notag
&&   -4   m^{\underline{L}L}  \cdot  \big( <  \derm_{\mu}A_{L} - \derm_{L}A_{\mu}  ,  [A_{\nu},A_{\underline{L}}] >   + <   [A_{\mu},A_{L}] , \derm_{\nu}A_{\underline{L}} - \derm_{\underline{L}}A_{\nu}  > \big)  \\
\notag
 &&  -4  m^{e_A e_A}  \cdot  \big( <  \derm_{\mu}A_{e_A} - \derm_{e_A}A_{\mu}  ,  [A_{\nu},A_{e_A}] >    -4  <   [A_{\mu},A_{e_A}] , \derm_{\nu}A_{e_A} - \derm_{e_A}A_{\nu}  > \big)  \\
\notag
&&  + m_{\mu\nu }  m^{L\underline{L}}  m^{L \underline{L}}    \cdot \big(  <  \derm_{L}A_{\underline{L}} - \derm_{\underline{L}}A_{L} , [A_{\underline{L}},A_{L}] >    +  <  [A_{L},A_{\underline{L}}] , \derm_{\underline{L}}A_{L} -\derm_{L}A_{\underline{L}}  > \big) \\
\notag
&& + m_{\mu\nu }  m^{L\underline{L}}  m^{\underline{L}L}    \cdot \big(  <  \derm_{\underline{L}}A_{\underline{L}} - \derm_{\underline{L}}A_{\underline{L}} , [A_{L},A_{L}] >    +  <  [A_{\underline{L}},A_{\underline{L}}] ,\derm_{L}A_{L} - \derm_{L}A_{L}  > \big) \\
\notag
 && + m_{\mu\nu }  m^{L\underline{L}}  m^{e_A e_A}    \cdot \big(  <  \derm_{e_A }A_{\underline{L}} - \derm_{\underline{L}}A_{e_A} , [A_{e_A},A_{L}] >    +  <  [A_{e_A},A_{\underline{L}}] , \derm_{e_A}A_{L} - \derm_{L}A_{e_A}  > \big) \\
 \notag
 && +   m_{\mu\nu }  m^{\underline{L}L}  m^{L\underline{L}}    \cdot \big(  <  \derm_{L}A_{L} - \derm_{L}A_{L} , [A_{\underline{L}},A_{\underline{L}}] >    +  <  [A_{L},A_{L}] , \derm_{\underline{L}}A_{\underline{L}} - \derm_{\underline{L}}A_{\underline{L}}  > \big) \\
 \notag
 && +   m_{\mu\nu }  m^{\underline{L}L}  m^{\underline{L}L}    \cdot \big(  <  \derm_{\underline{L}}A_{L} - \derm_{L}A_{\underline{L}} , [A_{L},A_{\underline{L}}] >    +  <  [A_{\underline{L}},A_{L}] , \derm_{L}A_{\underline{L}} - \derm_{\underline{L}}A_{L}  > \big) \\
 \notag
 && +  m_{\mu\nu }   m^{\underline{L}L}  m^{e_A e_A}    \cdot \big(  <  \derm_{e_A }A_{L} - \derm_{L}A_{e_A } , [A_{e_A},A_{\underline{L}}] >    +  <  [A_{e_A },A_{L}] , \derm_{e_A}A_{\underline{L}} - \derm_{\underline{L}}A_{e_A}  > \big) \\
 \notag
  && +   m_{\mu\nu }  m^{e_A e_A}  m^{L\underline{L}}    \cdot \big(  <  \derm_{L}A_{e_A} - \derm_{e_A}A_{L} , [A_{\underline{L}},A_{e_A}] >    +  <  [A_{L},A_{e_A}] , \derm_{\underline{L}}A_{e_A} - \derm_{e_A}A_{\underline{L}}  > \big) \\
  \notag
 && +   m_{\mu\nu }  m^{e_A e_A}  m^{\underline{L}L}    \cdot \big(  <  \derm_{\underline{L}}A_{e_A} - \derm_{e_A}A_{\underline{L}} , [A_{L},A_{e_A}] >    +  <  [A_{\underline{L}},A_{e_A}] , \derm_{L}A_{e_A} -\derm_{e_A}A_{L}  > \big) \\
 \notag
 && +   m_{\mu\nu }  m^{e_A e_A}  m^{e_A e_A}    \cdot \big(  <  \derm_{e_A}A_{e_A} - \derm_{e_A}A_{e_A} , [A_{e_A},A_{e_A}] >   \\
 \notag
 &&  +  <  [A_{e_A},A_{e_A}] , \derm_{e_A}A_{e_A} - \derm_{e_A}A_{e_A}  > \big) \\
 \notag
 && -4  m^{L\underline{L}}  \cdot   <   [A_{\mu},A_{\underline{L}}] ,  [A_{\nu},A_{L}] >  -4  m^{\underline{L}L}  \cdot   <   [A_{\mu},A_{L}] ,  [A_{\nu},A_{\underline{L}}] >  -4  m^{e_A e_A}  \cdot   <   [A_{\mu},A_{e_A}] ,  [A_{\nu},A_{e_A}] > \\
 \notag
  &&  + m_{\mu\nu } m^{L\underline{L}}  m^{L\underline{L}}   \cdot   <  [A_{L},A_{\underline{L}}] , [A_{\underline{L}},A_{L}] >  + m_{\mu\nu } m^{L\underline{L}}  m^{e_A e_A}   \cdot   <  [A_{e_A},A_{\underline{L}}] , [A_{e_A},A_{L}] > \\
  \notag
&&  + m_{\mu\nu } m^{\underline{L}L}  m^{\underline{L}L}   \cdot   <  [A_{\underline{L}},A_{L}] , [A_{L},A_{\underline{L}}] >  + m_{\mu\nu } m^{\underline{L}L}  m^{e_A e_A}   \cdot   <  [A_{e_A},A_{L}] , [A_{e_A},A_{\underline{L}}] > \\
\notag
 && + m_{\mu\nu } m^{e_A e_A}  m^{L\underline{L}}   \cdot   <  [A_{L},A_{e_A}] , [A_{\underline{L}},A_{e_A}] > +  m_{\mu\nu } m^{e_A e_A}  m^{\underline{L} L}   \cdot   <  [A_{\underline{L}},A_{e_A}] , [A_{L},A_{e_A}] > \\
     && + O \big(h \cdot  (\pa A)^2 \big)   + O \big(  h  \cdot  A^2 \cdot \pa A \big)     + O \big(  h   \cdot  A^4 \big)  \,  . 
\eeaa
Hence
    \beaa
\notag
&& g^{\alpha\beta} \derm_\alpha \derm_\beta g_{\mu\nu}   \\
\notag
     &=&  P(\pa_\mu h,\pa_\nu h) +  Q_{\mu\nu}(\pa h,\pa h)   + G_{\mu\nu}(h)(\pa h,\pa h)  \\
\notag
 &&  +2  <    \derm_{\mu}A_{L} -  \derm_{L}A_{\mu}  ,  \derm_{\nu}A_{ \underline{L} } -  \derm_{ \underline{L} }A_{\nu}  >   \\
 \notag
 &&  +2 <   \derm_{\mu}A_{ \underline{L} } - \derm_{ \underline{L} }A_{\mu}  ,  \derm_{\nu}A_{L} - \derm_{L}A_{\nu}  >    \\
 \notag
  && -4   <    \derm_{\mu}A_{e_A} -  \derm_{e_A }A_{\mu}  ,   \derm_{\nu}A_{e_A} - \derm_{e_A}A_{\nu}  > \\
  \notag
  &&   + \frac{1}{4} m_{\mu\nu }    \cdot  <   \derm_{L}A_{\underline{L}} - \derm_{\underline{L}}A_{L} ,  \derm_{\underline{L}}A_{L} - \derm_{L}A_{\underline{L}} > \\
\notag
&& -  \frac{1}{2} m_{\mu\nu }     \cdot  <   \derm_{e_A}A_{\underline{L}} -  \derm_{\underline{L}}A_{e_A} ,  \derm_{e_A}A_{L} - \derm_{L}A_{e_A} > \\
\notag
&& +  \frac{1}{4} m_{\mu\nu }     \cdot  <   \derm_{\underline{L}}A_{L} -  \derm_{L}A_{\underline{L}} ,  \derm_{L}A_{\underline{L}} -  \derm_{\underline{L}}A_{L} > \\
\notag
&& -  \frac{1}{2}  m_{\mu\nu }     \cdot  <   \derm_{e_A}A_{L} -  \derm_{L}A_{e_A} ,  \derm_{e_A}A_{\underline{L}} -  \derm_{\underline{L}}A_{e_A} > \\
\notag
&&  -  \frac{1}{2} m_{\mu\nu }     \cdot  <   \derm_{L}A_{e_A} -  \derm_{e_A}A_{L} ,  \derm_{\underline{L}}A_{e_A} -  \derm_{e_A}A_{\underline{L}} > \\
\notag
&& -  \frac{1}{2}  m_{\mu\nu }      \cdot  <   \derm_{\underline{L}}A_{e_A} -  \derm_{e_A}A_{\underline{L}} ,  \derm_{L}A_{e_A} -  \derm_{e_A}A_{L} > \\
\notag
 &&  +2  \big( <   \derm_{\mu}A_{\underline{L}} - \derm_{\underline{L}}A_{\mu}  ,  [A_{\nu},A_{L}] >   + <   [A_{\mu},A_{\underline{L}}] ,   \derm_{\nu}A_{L} -  \derm_{L}A_{\nu}  > \big) \\
 \notag
&&   + 2  \big( <   \derm_{\mu}A_{L} - \derm_{L}A_{\mu}  ,  [A_{\nu},A_{\underline{L}}] >   + <   [A_{\mu},A_{L}] ,  \derm_{\nu}A_{\underline{L}} - \derm_{\underline{L}}A_{\nu}  > \big)  \\
\notag
 &&  -4     \big( <    \derm_{\mu}A_{e_A} - \derm_{e_A}A_{\mu}  ,  [A_{\nu},A_{e_A}] >    -4  <   [A_{\mu},A_{e_A}] ,   \derm_{\nu}A_{e_A} -  \derm_{e_A}A_{\nu}  > \big)  \\
\notag
&&  +  \frac{1}{4} m_{\mu\nu }   \cdot \big(  <  \derm_{L}A_{\underline{L}} -  \derm_{\underline{L}}A_{L} , [A_{\underline{L}},A_{L}] >    +  <  [A_{L},A_{\underline{L}}] ,  \derm_{\underline{L}}A_{L} -  \derm_{L}A_{\underline{L}}  > \big) \\
\notag
 &&  -  \frac{1}{2}  m_{\mu\nu }     \cdot \big(  <   \derm_{e_A }A_{\underline{L}} -  \derm_{\underline{L}}A_{e_A} , [A_{e_A},A_{L}] >    +  <  [A_{e_A},A_{\underline{L}}] ,  \derm_{e_A}A_{L} -  \derm_{L}A_{e_A}  > \big) \\
 \notag
 && +     \frac{1}{4} m_{\mu\nu }      \cdot \big(  <   \derm_{\underline{L}}A_{L} -  \derm_{L}A_{\underline{L}} , [A_{L},A_{\underline{L}}] >    +  <  [A_{\underline{L}},A_{L}] ,  \derm_{L}A_{\underline{L}} -  \derm_{\underline{L}}A_{L}  > \big) \\
 \notag
 &&  -  \frac{1}{2}  m_{\mu\nu }       \cdot \big(  <  \derm_{e_A }A_{L} - \derm_{L}A_{e_A } , [A_{e_A},A_{\underline{L}}] >    +  <  [A_{e_A },A_{L}] ,  \derm_{e_A}A_{\underline{L}} -  \derm_{\underline{L}}A_{e_A}  > \big) \\
 \notag
  &&  -  \frac{1}{2}   m_{\mu\nu }     \cdot \big(  <   \derm_{L}A_{e_A} -  \derm_{e_A}A_{L} , [A_{\underline{L}},A_{e_A}] >    +  <  [A_{L},A_{e_A}] ,  \derm_{\underline{L}}A_{e_A} -  \derm_{e_A}A_{\underline{L}}  > \big) \\
  \notag
 &&    -  \frac{1}{2} m_{\mu\nu }      \cdot \big(  <   \derm_{\underline{L}}A_{e_A} -  \derm_{e_A}A_{\underline{L}} , [A_{L},A_{e_A}] >    +  <  [A_{\underline{L}},A_{e_A}] ,  \derm_{L}A_{e_A} -  \derm_{e_A}A_{L}  > \big) \\
 \notag
 && + 2   <   [A_{\mu},A_{\underline{L}}] ,  [A_{\nu},A_{L}] >  + 2   <   [A_{\mu},A_{L}] ,  [A_{\nu},A_{\underline{L}}] >  -4    <   [A_{\mu},A_{e_A}] ,  [A_{\nu},A_{e_A}] > \\
 \notag
  && +   \frac{1}{4}   m_{\mu\nu }    \cdot   <  [A_{L},A_{\underline{L}}] , [A_{\underline{L}},A_{L}] >   -  \frac{1}{2} m_{\mu\nu }   \cdot   <  [A_{e_A},A_{\underline{L}}] , [A_{e_A},A_{L}] > \\
  \notag
&&  +   \frac{1}{4}  m_{\mu\nu }  \cdot   <  [A_{\underline{L}},A_{L}] , [A_{L},A_{\underline{L}}] >    -  \frac{1}{2} m_{\mu\nu }  \cdot   <  [A_{e_A},A_{L}] , [A_{e_A},A_{\underline{L}}] > \\
\notag
 &&  -  \frac{1}{2} m_{\mu\nu }  \cdot   <  [A_{L},A_{e_A}] , [A_{\underline{L}},A_{e_A}] >   -  \frac{1}{2}  m_{\mu\nu }    \cdot   <  [A_{\underline{L}},A_{e_A}] , [A_{L},A_{e_A}] > \\
     && + O \big(h \cdot  (\pa A)^2 \big)   + O \big(  h  \cdot  A^2 \cdot \pa A \big)     + O \big(  h   \cdot  A^4 \big)  \,  .
\eeaa
Therefore,
  \bea
\notag
&& g^{\alpha\beta} \derm_\alpha \derm_\beta g_{\mu\nu}   \\
\notag
     &=&  P(\pa_\mu h,\pa_\nu h) +  Q_{\mu\nu}(\pa h,\pa h)   + G_{\mu\nu}(h)(\pa h,\pa h)  \\
\notag
 &&  +2  <    \derm_{\mu}A_{L} -  \derm_{L}A_{\mu}  ,  \derm_{\nu}A_{ \underline{L} } -  \derm_{ \underline{L} }A_{\nu}  >   \\
 \notag
 &&  +2 <   \derm_{\mu}A_{ \underline{L} } - \derm_{ \underline{L} }A_{\mu}  ,  \derm_{\nu}A_{L} - \derm_{L}A_{\nu}  >    \\
 \notag
  && -4   <    \derm_{\mu}A_{e_A} -  \derm_{e_A }A_{\mu}  ,   \derm_{\nu}A_{e_A} - \derm_{e_A}A_{\nu}  > \\
  \notag
  &&   + \frac{1}{2} m_{\mu\nu }    \cdot  <   \derm_{L}A_{\underline{L}} - \derm_{\underline{L}}A_{L} ,  \derm_{\underline{L}}A_{L} - \derm_{L}A_{\underline{L}} > \\
\notag
&& -  m_{\mu\nu }     \cdot  <   \derm_{e_A}A_{\underline{L}} -  \derm_{\underline{L}}A_{e_A} ,  \derm_{e_A}A_{L} - \derm_{L}A_{e_A} > \\
\notag
&&  -   m_{\mu\nu }     \cdot  <   \derm_{L}A_{e_A} -  \derm_{e_A}A_{L} ,  \derm_{\underline{L}}A_{e_A} -  \derm_{e_A}A_{\underline{L}} > \\
\notag
 &&  +2  \big( <   \derm_{\mu}A_{\underline{L}} - \derm_{\underline{L}}A_{\mu}  ,  [A_{\nu},A_{L}] >   + <   [A_{\mu},A_{\underline{L}}] ,   \derm_{\nu}A_{L} -  \derm_{L}A_{\nu}  > \big) \\
 \notag
&&   + 2  \big( <   \derm_{\mu}A_{L} - \derm_{L}A_{\mu}  ,  [A_{\nu},A_{\underline{L}}] >   + <   [A_{\mu},A_{L}] ,  \derm_{\nu}A_{\underline{L}} - \derm_{\underline{L}}A_{\nu}  > \big)  \\
\notag
 &&  -4     \big( <    \derm_{\mu}A_{e_A} - \derm_{e_A}A_{\mu}  ,  [A_{\nu},A_{e_A}] >    -4  <   [A_{\mu},A_{e_A}] ,   \derm_{\nu}A_{e_A} -  \derm_{e_A}A_{\nu}  > \big)  \\
\notag
&&  +  \frac{1}{2} m_{\mu\nu }   \cdot   <  \derm_{L}A_{\underline{L}} -  \derm_{\underline{L}}A_{L} , [A_{\underline{L}},A_{L}] >   \\
\notag
 &&  -  \frac{1}{2}  m_{\mu\nu }     \cdot \big(  <   \derm_{e_A }A_{\underline{L}} -  \derm_{\underline{L}}A_{e_A} , [A_{e_A},A_{L}] >    +  <  [A_{e_A},A_{\underline{L}}] ,  \derm_{e_A}A_{L} -  \derm_{L}A_{e_A}  > \big) \\
 \notag
 && +     \frac{1}{2} m_{\mu\nu }      \cdot  <   \derm_{\underline{L}}A_{L} -  \derm_{L}A_{\underline{L}} , [A_{L},A_{\underline{L}}] >     \\
 \notag
 &&  -  \frac{1}{2}  m_{\mu\nu }       \cdot \big(  <  \derm_{e_A }A_{L} - \derm_{L}A_{e_A } , [A_{e_A},A_{\underline{L}}] >    +  <  [A_{e_A },A_{L}] ,  \derm_{e_A}A_{\underline{L}} -  \derm_{\underline{L}}A_{e_A}  > \big) \\
 \notag
  &&  -     m_{\mu\nu }     \cdot \big(  <   \derm_{L}A_{e_A} -  \derm_{e_A}A_{L} , [A_{\underline{L}},A_{e_A}] >    +  <  [A_{L},A_{e_A}] ,  \derm_{\underline{L}}A_{e_A} -  \derm_{e_A}A_{\underline{L}}  > \big) \\
 \notag
 && + 2   <   [A_{\mu},A_{\underline{L}}] ,  [A_{\nu},A_{L}] >  + 2   <   [A_{\mu},A_{L}] ,  [A_{\nu},A_{\underline{L}}] >  -4    <   [A_{\mu},A_{e_A}] ,  [A_{\nu},A_{e_A}] > \\
 \notag
  && +   \frac{1}{2}   m_{\mu\nu }    \cdot   <  [A_{L},A_{\underline{L}}] , [A_{\underline{L}},A_{L}] >   - m_{\mu\nu }   \cdot   <  [A_{e_A},A_{\underline{L}}] , [A_{e_A},A_{L}] > \\
\notag
 &&  -   m_{\mu\nu }  \cdot   <  [A_{L},A_{e_A}] , [A_{\underline{L}},A_{e_A}] >  \\
     && + O \big(h \cdot  (\pa A)^2 \big)   + O \big(  h  \cdot  A^2 \cdot \pa A \big)     + O \big(  h   \cdot  A^4 \big)  \,  ,
\eea
where there are still same terms that add up together.
\end{proof}

\begin{lemma}
We have for all sufficiently smooth function $f$, 
\bea
|\pa  f  |^2 = \frac{1}{2} |\pa_L  f  |^2 +  \frac{1}{2} |\pa_{\underline{L}}  f  |^2 +   |\pa_{e_1}  f  |^2 + |\pa_{e_2}  f  |^2 \; .
\eea
\end{lemma}

\begin{proof}
We compute
\beaa
E ( L,  L ) &=& E (  \pa_{t} - \pa_{r},  \pa_{t} - \pa_{r}) = E (  \pa_{t},  \pa_{t} ) -2 E (  \pa_{t},  \pa_{r})  + E (   \pa_{r}, \pa_{r}) \\
 &=& 1   + E (  \frac{x^{i}}{r} \pa_{i}, \frac{x^{j}}{r} \pa_{j} ) = 1  + \sum_{i=1}^3 \frac{(x^{i})^2}{r^2} \\
  &=& 2 \; ,\\
E ( \underline{L},  \underline{L} ) &=& E (  \pa_{t} + \pa_{r},  \pa_{t} + \pa_{r}) \\
&=& 2 \;,\\
E ( L,  \underline{L} ) &=& E (  \pa_{t} - \pa_{r},  \pa_{t} + \pa_{r}) = E (  \pa_{t},  \pa_{t} ) - E (   \pa_{r}, \pa_{r}) \\
 &=& 1  -1 \\
  &=& 0 \; , 
\eeaa
and we compute for all $A, B \in \{1,2 \}$, 
\beaa
E ( L,  e_A) &=& E ( \underline{L},  e_A) = 0 \;, \\
E ( e_A, e_B ) &=& \delta_{AB} \;.
\eeaa
Therefore, for a sufficiently smooth function $f$, we get 
\beaa
| \pa f |^2 &=& E^{LL} | \pa_{L} f  |^2 + E^{ \underline{L}  \underline{L} } | \pa_{ \underline{L} } f  |^2 + \sum_{A=1}^{2} | \pa_{e_A} f  |^2 \\
&=& \frac{1}{2} | \pa_{L} f  |^2 +\frac{1}{2} | \pa_{ \underline{L} } f  |^2 + \sum_{A=1}^{2} | \pa_{e_A} f  |^2 \; .
\eeaa

\end{proof}

\begin{lemma}\label{equivalencerestrictedpartialderivativeandsumovertangentailframe}
We have for all sufficiently smooth function $f$,
\beaa
 |\rpa f |^2 \sim  |\pa_L f |^2 +  \sum_{A=1}^{2} | \pa_{e_A} f  |^2 \; .
\eeaa

\end{lemma}

\begin{proof}
We have by definition
\beaa
 |\rpa f |^2 :=     |\rpa_0 f |^2 + \sum_{i=1}^3   |\rpa_i f |^2   \; ,
  \eeaa
  where also by definition
\beaa
\rpa_{0} &:=& \frac{ \pa_t +  \pa_r }{2}  = \frac{1}{2} \pa_{L} \, .
\eeaa
We showed, in Lemma \ref{restrictedderivativesintermsofZ}, that
 \beaa
\rpa_i  &=&   \frac{x^j}{r^2}   Z_{ij} \, ,
 \eeaa
 and that we also have
 \bea
e_A = \frac{1}{r}C^{ij}_A Z_{ij} \, ,
\eea
where $C^{ij}_A$ are bounded function on $\SSS^2$.
Thus,
\beaa
 |\rpa f |^2 &=&     |\rpa_0 f |^2 + \sum_{i=1}^{3}   |\rpa_i f |^2   \\
 &=& \frac{1}{4}   |\pa_{L}  f |^2 +  \sum_{i=1}^{3}   | \frac{x^j}{r^2}   Z_{ij}f  |^2 \\
  &=& \frac{1}{4}   |\pa_{L}  f |^2 +  \sum_{i=1}^{3}   | \frac{1}{r} \cdot \frac{x^j}{r}   Z_{ij}f  |^2 \\
  &=& \frac{1}{4}   |\pa_{L}  f |^2 +   \sum_{A=1}^{2}  | \frac{1}{r} \cdot C^{ij}_A Z_{ij} f  |^2 \\
    &=& \frac{1}{4}   |\pa_{L}  f |^2 +  \sum_{A=1}^{2}   | \pa_{e_A} f  |^2 .
  \eeaa

 \end{proof}

\begin{definition}\label{definitionofnormsoncovariantgradiantandrestrictedcovariantgradient}
Let $\Psi$ be a tensor of arbitrary order, for example a one-tensor $\Psi_{\a}$, either a scalar or valued in the Lie algebra $\cal G$. We define
\bea
 | {\der^{(\bf{m})}}  \Psi |^{2} &:=& \sum_{\a, \b \in \cal U}  | {\der^{(\bf{m})}}_\a  \Psi_\b |^{2}  \; , 
 \eea
 and
  \bea
|\rderm  \Psi |^2  &:=& \sum_{\a \in \cal T, \b \in \cal U }  | {\der^{(\bf{m})}}_\a  \Psi_\b |^{2}  \; .
\eea
\end{definition}
\begin{lemma}
Let $\Psi$ be a tensor of arbitrary order, either a scalar or valued in the Lie algebra $\cal G$. We have
\bea
|\pa  \Psi |^2  &\sim&  | {\der^{(\bf{m})}}  \Psi |^{2}  \; , 
\eea
and we have
 \bea
|\rderm  \Psi |^2  &\les& |\pa  \Psi |^2  \; .
\eea

\end{lemma}

\begin{proof}
On one hand, we compute
\beaa
|\pa  \Psi |^2 &=& E^{\a\b} E^{\mu\nu} < {\der^{(\bf{m})}}_\a \Psi_\mu,  {\der^{(\bf{m})}}_\b \Psi_\nu > \\
&=&  \frac{1}{2} E^{\mu\nu} < {\der^{(\bf{m})}}_L \Psi_\mu,  {\der^{(\bf{m})}}_L \Psi_\nu > +  \frac{1}{2} E^{\mu\nu} < {\der^{(\bf{m})}}_{\underline{L}} \Psi_\mu,  {\der^{(\bf{m})}}_{\underline{L} }\Psi_\nu >   \\
&& +  E^{\mu\nu} < {\der^{(\bf{m})}}_{e_1}  \Psi_\mu,  {\der^{(\bf{m})}}_{e_1}  \Psi_\nu > +  E^{\mu\nu} < {\der^{(\bf{m})}}_{e_{2}}  \Psi_\mu,  {\der^{(\bf{m})}}_{e_{2}}  \Psi_\nu > \\
&=&  \frac{1}{4} | {\der^{(\bf{m})}}_L  \Psi_L |^{2} + \frac{1}{4} | {\der^{(\bf{m})}}_{L}  \Psi_{\underline{L}} |^{2} +  \frac{1}{2} | {\der^{(\bf{m})}}_L  \Psi_{e_1} |^{2} +   \frac{1}{2} | {\der^{(\bf{m})}}_L  \Psi_{e_{2}} |^{2}\\
&& +   \frac{1}{4} | {\der^{(\bf{m})}}_{\underline{L}}  \Psi_L |^{2} + \frac{1}{4} | {\der^{(\bf{m})}}_{\underline{L}}  \Psi_{\underline{L}} |^{2} +  \frac{1}{2} | {\der^{(\bf{m})}}_{\underline{L}}  \Psi_{e_1} |^{2} +   \frac{1}{2} | {\der^{(\bf{m})}}_{\underline{L}}  \Psi_{e_{2}} |^{2}\\
&&  +   \frac{1}{2} | {\der^{(\bf{m})}}_{e_1}  \Psi_L |^{2} + \frac{1}{2} | {\der^{(\bf{m})}}_{e_1}  \Psi_{\underline{L}} |^{2} +   | {\der^{(\bf{m})}}_{e_1}  \Psi_{e_1} |^{2} +   | {\der^{(\bf{m})}}_{e_1}  \Psi_{e_{2}} |^{2}\\
&&+ \frac{1}{2} | {\der^{(\bf{m})}}_{e_{2}}  \Psi_L |^{2} + \frac{1}{2} | {\der^{(\bf{m})}}_{e_{2}}  \Psi_{\underline{L}} |^{2} +   | {\der^{(\bf{m})}}_{e_{2}}  \Psi_{e_1} |^{2}  +  | {\der^{(\bf{m})}}_{e_{2}}  \Psi_{e_{2}} |^{2} \\
 &\sim&  | {\der^{(\bf{m})}}  \Psi |^{2}\; .
\eeaa

On the other hand, we have
 \beaa
|\rderm  \Psi |^2  &:=&  \sum_{\a \in \cal T\,,\; \b \in \cal U }  | {\der^{(\bf{m})}}_\a  \Psi_\b |^{2}  \\
&\les&  \sum_{\a \in \cal U\,,\; \b \in \cal U }  | {\der^{(\bf{m})}}_\a  \Psi_\b |^{2}  \\
 &\les& | {\der^{(\bf{m})}}  \Psi |^{2} \\
 &\les& |\pa  \Psi |^2 ,
\eeaa
and similarly, for a tensor of an arbitrary order, we have the same estimates.

\end{proof}

\begin{definition}\label{definitionnormonrestrictedpartialgradient}
Let $\Psi$ be a tensor of arbitrary order, for example a one-tensor $\Psi_{\a}$, either a scalar or valued in the Lie algebra $\cal G$. We define
\bea
|  \rpa \Psi_\a |^2 :=\sum_{\b  \in  \{t, x^1, x^2, x^3 \}}  |  \rpa_{\b} \Psi_\a |^2 \; ,
\eea
and
 \bea
|  \rpa \Psi |^2  &:=&  \sum_{\a  \in  \{t, x^1, x^2, x^3 \}} |  \rpa \Psi_\a |^2  \; .
\eea
\end{definition}

\begin{lemma}\label{equivalenceoftworestrictednormsongradients}
Let $\Psi$ be a tensor of arbitrary order, either a scalar or valued in the Lie algebra $\cal G$. We have
 \bea
|\rderm  \Psi |^2  &\sim&   |  \rpa \Psi |^2  \; .
\eea

\end{lemma}

\begin{proof}

We consider the tensor in $\a$ and in $\b$ that is given by $ E^{\mu\nu} < {\der^{(\bf{m})}}_\a \Psi_\mu,  {\der^{(\bf{m})}}_\b \Psi_\nu > $. Since the trace does not depend on the system of coordinates, we can compute it in wave coordinates for which the Christoffel symbols of the covariant derivative $\derm$ vanish. Thus, by taking $\a, \b \in  \{t, x^1, x^2, x^3 \}$, we get
\beaa
 E^{\mu\nu} < {\der^{(\bf{m})}}_\a \Psi_\mu,  {\der^{(\bf{m})}}_\b \Psi_\nu > &=&  \sum_{\mu, \nu \in  \{t, x^1, x^2, x^3 \}}  E^{\mu\nu} < \pa_\a \Psi_\mu,  \pa_\b \Psi_\nu > \; .
\eeaa
However, since the right hand side is also a tensor in $\a$ and in $\b$, the identity holds true also for the frame $\cal T$, i.e. for any $\a, \b \in  \cal T$, we have
\beaa
 E^{\mu\nu} < {\der^{(\bf{m})}}_\a \Psi_\mu,  {\der^{(\bf{m})}}_\b \Psi_\nu > &=&  \sum_{\mu , \nu \in  \{t, x^1, x^2, x^3 \}}  E^{\mu\nu} < \pa_\a \Psi_\mu,  \pa_\b \Psi_\nu > \; .
\eeaa
Consequently,
\beaa
 \sum_{\ga \in \cal T }  E^{\mu\nu} < {\der^{(\bf{m})}}_\ga \Psi_\mu,  {\der^{(\bf{m})}}_\ga \Psi_\nu > &=&   \sum_{\ga \in \cal T } \sum_{\mu , \nu \in  \{t, x^1, x^2, x^3 \}}  E^{\mu\nu} < \pa_\ga \Psi_\mu,  \pa_\ga \Psi_\nu >  \\
 &=&   \sum_{\ga \in \cal T } \sum_{\mu  \in  \{t, x^1, x^2, x^3 \}}   < \pa_\ga \Psi_\mu,  \pa_\ga \Psi_\mu >  \\
 &=&   \sum_{\ga \in \cal T } \sum_{\mu  \in  \{t, x^1, x^2, x^3 \}}   | \pa_\ga \Psi_{\mu} |^2 .
\eeaa
However, we proved in Lemma \ref{equivalencerestrictedpartialderivativeandsumovertangentailframe}, that $\sum_{\ga \in \cal T } \sum_{\mu  \in  \{t, x^1, x^2, x^3 \}}   | \pa_\ga \Psi_{\mu} |^2 \sim \sum_{\mu  \in  \{t, x^1, x^2, x^3 \}}  | \rpa \Psi_\mu |^2$, hence
\beaa
 \sum_{\ga \in \cal T }  E^{\mu\nu} < {\der^{(\bf{m})}}_\ga \Psi_\mu,  {\der^{(\bf{m})}}_\ga \Psi_\nu > &\sim&    | \rpa \Psi |^2 ,
\eeaa
and since 
 \beaa
 \sum_{\ga \in \cal T }  E^{\mu\nu} < {\der^{(\bf{m})}}_\ga \Psi_\mu,  {\der^{(\bf{m})}}_\ga \Psi_\nu > &\sim&   \sum_{\ga \in \cal T }  \sum_{\mu \in \cal U }  | {\der^{(\bf{m})}}_\ga \Psi_\mu |^2 ,
 \eeaa
 we get the result.
 
\end{proof}

\begin{lemma}\label{estimatesonsourcetermswithusingLorenzgaugeandwavegaugeestimates}
For any $\si \in \cal U $, we have
     \bea
   \notag
 | g^{\la\mu} \derm_{\la}   \derm_{\mu}   A_{\si}  |  &\les& | \derm h | \cdot  |\rderm A |     \\
   \notag
 && +   | \derm  h | \cdot  |\derm  A_{L}  | + | \derm  h_{\cal T L} | \cdot  |\derm A |  \\
   \notag
   && +  \big( | \rderm  h |  +| \derm h_{L \si} | \big) \cdot  \big( |\derm  A_{L}  | +|\rderm  A | \big)  \\
      \notag
         && +   | \rderm  h |    \cdot  |\derm  A  | \\
 \notag
        && +  \big( | \rderm  h |   +| \derm h_{L \si}  |  \big) \cdot  | A  | \cdot  | A_L |   \\
   \notag
           && +   | \rderm  h |  \cdot  | A  |^2   \\
           \notag
                  && +  \big( | \rderm  h |   +| \derm h_{e_{a} \si} |  \big) \cdot    | A_{e_{a}} |  \cdot  | A_L |   \\
          \notag
 &&   +   | A_L  | \cdot    | \derm_\si A_{\underline{L}}  |      +    |  A  | \cdot  | A_\si  | \cdot  | A_L |   \\
    \notag
    &&   +  | A  | \cdot   | \derm_\si A_L  |     \\
       \notag
 &&  + | A_{e_{a}}  | \cdot    | \derm_\si A_{e_{a}} |    +    |  A_{e_{a}}  |^2 \cdot  | A_\si  |   \\
    \notag
  && + O( h \cdot  \pa h \cdot  \pa A) + O( h \cdot  \pa h \cdot  A^2) + O( h \cdot  A \cdot \pa A) + O( h \cdot  A^3) \, . \\
  \eea

For any $\mu, \nu \in \cal U$, we have
    \bea
\notag
 | g^{\alpha\beta} \derm_\alpha \derm_\beta h_{\mu\nu} |     &\les&  | P(\derm_\mu h,\derm_\nu h) | + | Q_{\mu\nu}(\derm h,\derm h)  | + | G_{\mu\nu}(h)(\derm h,\derm h)  | \\
 \notag
    &&+ \big(   | \derm  A_{L}  | +  | \rderm A   | \big) \cdot  |\derm  A  | \\
    \notag
 &&+ \big(  |\derm_\mu  A_{e_{a}}  | +  |\rderm  A  | \big)  \cdot \big(  |\derm_\nu  A_{e_{a}}  | +  |\rderm  A  | \big) \\
   \notag
 && +  | \derm A   |  \cdot | A | \cdot | A_L | + \big(   | \derm A_L  | +  | \rderm A  | \big) \cdot | A |^2   \\
 \notag
&& + \big(   | \derm_\mu  A  | +  | \rderm A   | \big) \cdot | A_\nu | \cdot | A | + \big(   | \derm_\nu  A  | +  | \rderm A   | \big) \cdot | A_\mu | \cdot | A | \\
\notag
&& + \big(   | \derm  A_{L}  | +  | \rderm A   | \big) \cdot | A | \cdot  | A_L |  \\
\notag
&& +  | \derm A   |  \cdot | A | \cdot  | A_L |  +  | \rderm A   |  \cdot | A |^2  \\
 \notag
&& +  | A_\mu | \cdot | A_\nu | \cdot  | A |^2  + | A_L | \cdot  | A |^3 \\
     && + O \big(h \cdot  (\pa A)^2 \big)   + O \big(  h  \cdot  A^2 \cdot \pa A \big)     + O \big(  h   \cdot  A^4 \big)  \,  .
\eea

\end{lemma}

\begin{proof}
For any $\si \in \cal U $, we have

    \beaa
   \notag
&& | g^{\la\mu} \derm_{\la}   \derm_{\mu}   A_{\si}  | \\ 
  \notag
 &\les& | \derm_\si h | \cdot  |\rderm A |     \\
   \notag
 && +   | \derm_\si  h | \cdot  |\derm  A_{L}  | + | \derm  h_{\cal T L} | \cdot  |\derm A |  \\
 \notag
&& +  | \derm_\si h | \cdot  |\rderm A |     \\
\notag
&& +  \big( | \rderm  h |  +| \derm h_{L \si} | \big) \cdot  \big( |\derm  A_{L}  | +|\rderm  A | \big)  \\
 \notag
          && +  | \rderm  h |  \cdot  \big( |\derm  A  | +|\rderm  A | \big)  \\
    \notag
       && +  \big( | \rderm  h |  +| \derm h_{e_{a} \si} | \big) \cdot  |\rderm  A |  \\
 \notag
        && +  \big( | \rderm  h |  +| \derm h_{L \si} | \big) \cdot  | A  | \cdot  | A_L |   \\
   \notag
           && +  | \rderm  h |   \cdot  | A  |  \cdot  | A_{e_{a}} |    \\
   \notag
        && +  \big( | \rderm  h |   +| \derm h_{e_{a} \si} |  \big) \cdot    | A_{e_{a}} |  \cdot  | A_L |   \\
 \notag
 &&   +   | A_L  | \cdot    | \derm_\si A_{\underline{L}}   |      +    |  A  | \cdot  | A_\si  | \cdot  | A_L |   \\
    \notag
    &&   +  | A  | \cdot   | \derm_\si A_L  |  +    |  A  | \cdot  | A_\si  | \cdot  | A_L |     \\
       \notag
 &&  + | A_{e_{a}}  | \cdot    | \derm_\si A_{e_{a}} |    +    |  A_{e_{a}}  |^2 \cdot  | A_\si  |   \\
    \notag
  && + O( h \cdot  \pa h \cdot  \pa A) + O( h \cdot  \pa h \cdot  A^2) + O( h \cdot  A \cdot \pa A) + O( h \cdot  A^3) \; . \\
  \eeaa
  
 This leads to 

    \beaa
   \notag
&& | g^{\la\mu} \derm_{\la}   \derm_{\mu}   A_{\si}  | \\ 
  \notag
 &\les& | \derm h | \cdot  |\rderm A |     \\
   \notag
 && +   | \derm  h | \cdot  |\derm  A_{L}  | + | \derm  h_{\cal T L} | \cdot  |\derm A |  \\
   \notag
   && +  \big( | \rderm  h |  +| \derm h_{L \si} | \big) \cdot  \big( |\derm  A_{L}  | +|\rderm  A | \big)  \\
   \notag
&& +    | \rderm  h |  \cdot  |\derm  A  | \\
      \notag
         && +| \derm h| \cdot  |\rderm  A  | \\
 \notag
        && +  \big( | \rderm  h |   +| \derm h_{L \si}   | \big) \cdot  | A  | \cdot  | A_L |   \\
   \notag
           && +   | \rderm  h |    \cdot  | A  |  \cdot  | A_{e_{a}} |    \\
\notag
        && +  \big( | \rderm  h |   +| \derm h_{e_{a} \si} |  \big) \cdot    | A_{e_{a}} |  \cdot  | A_L |   \\
 \notag
 &&   +   | A_L  | \cdot    | \derm_\si A_{\underline{L}}   |      +    |  A  | \cdot  | A_\si  | \cdot  | A_L |   \\
    \notag
    &&   +  | A  | \cdot   | \derm_\si A_L  |     \\
       \notag
 &&  + | A_{e_{a}}  | \cdot    | \derm_\si A_{e_{a}} |    +    |  A_{e_{a}}  |^2 \cdot  | A_\si  |   \\
    \notag
  && + O( h \cdot  \pa h \cdot  \pa A) + O( h \cdot  \pa h \cdot  A^2) + O( h \cdot  A \cdot \pa A) + O( h \cdot  A^3) \; . \\
  \eeaa

Finally, we write more simply
    \beaa
   \notag
 | g^{\la\mu} \derm_{\la}   \derm_{\mu}   A_{\si}  |  &\les& | \derm h | \cdot  |\rderm A |     \\
   \notag
 && +   | \derm  h | \cdot  |\derm  A_{L}  | + | \derm  h_{\cal T L} | \cdot  |\derm A |  \\
   \notag
   && +  \big( | \rderm  h |  +| \derm h_{L \si} | \big) \cdot  \big( |\derm  A_{L}  | +|\rderm  A | \big)  \\
      \notag
         && +   | \rderm  h |    \cdot  |\derm  A  | \\
 \notag
        && +  \big( | \rderm  h |   +| \derm h_{L \si}  |  \big) \cdot  | A  | \cdot  | A_L |   \\
   \notag
           && +   | \rderm  h |  \cdot  | A  |^2   \\
           \notag
                  && +  \big( | \rderm  h |   +| \derm h_{e_{a} \si} |  \big) \cdot    | A_{e_{a}} |  \cdot  | A_L |   \\
          \notag
 &&   +   | A_L  | \cdot    | \derm_\si A_{\underline{L}}   |      +    |  A  | \cdot  | A_\si  | \cdot  | A_L |   \\
    \notag
    &&   +  | A  | \cdot   | \derm_\si A_L  |     \\
       \notag
 &&  + | A_{e_{a}}  | \cdot    | \derm_\si A_{e_{a}} |    +    |  A_{e_{a}}  |^2 \cdot  | A_\si  |   \\
    \notag
  && + O( h \cdot  \pa h \cdot  \pa A) + O( h \cdot  \pa h \cdot  A^2) + O( h \cdot  A \cdot \pa A) + O( h \cdot  A^3) \; . \\
  \eeaa

To estimate the source term for the wave equation on the metric, we consider that for any $\mu, \nu \in \cal U$, we have
\beaa
| m_{\mu\nu } | \les 1 \, ,
\eeaa
and thus, we have for all $\mu, \nu \in \cal U$
    \bea
\notag
&& | g^{\alpha\beta} \derm_\alpha \derm_\beta h_{\mu\nu} |  \\
\notag
     &\les&  | P(\pa_\mu h,\pa_\nu h) | + | Q_{\mu\nu}(\pa h,\pa h)  | + | G_{\mu\nu}(h)(\pa h,\pa h)  | \\
     \notag
   &&+ \big(   | \derm  A_{L}  | +  | \rderm A   | \big) \cdot  |\derm  A  | \\
 \notag
   &&+  |\derm  A  | \cdot  \big(   | \derm  A_{L}  | +  | \rderm A   | \big) \\
 \notag
  &&+ \big(  |\derm_\mu  A_{e_{a}}  | +  |\rderm  A  | \big)  \cdot \big(  |\derm_\nu  A_{e_{a}}  | +  |\rderm  A  | \big) \\
    \notag
       &&+ \big(   | \derm  A_{L}  | +  | \rderm A   | \big) \cdot  \big(   | \derm  A_{L}  | +  | \rderm A   | \big) \\
\notag
   &&+   \big(   | \rderm  A  | +  | \derm A_L   | \big)  \cdot    |\rderm  A  |  \\
\notag
   &&+       |\rderm  A  | \cdot \big(   | \rderm  A  | +  | \derm A   | \big)   \\
   \notag
   && +  | \derm A   |  \cdot | A | \cdot | A_L | + \big(   | \derm A_L  | +  | \rderm A  | \big) \cdot | A |^2   \\
 \notag
 && + \big(   | \derm  A_{L}  | +  | \rderm A   | \big) \cdot | A |^2 +  | \derm A | \cdot | A | \cdot  | A_L | \\
\notag
&& + \big(   | \derm_\mu  A  | +  | \rderm A   | \big) \cdot | A_\nu | \cdot | A | + \big(   | \derm_\nu  A  | +  | \rderm A   | \big) \cdot | A_\mu | \cdot | A | \\
\notag
&& + \big(   | \derm  A_{L}  | +  | \rderm A   | \big) \cdot | A | \cdot  | A_L |  \\
\notag
&& +  | \derm A   |  \cdot | A | \cdot  | A_L |  +  | \rderm A   |  \cdot | A |^2  \\
 \notag
&& +  | A_\mu | \cdot | A_\nu | \cdot  | A |^2 \\
 \notag
 && + | A_L | \cdot  | A |^3 \\
  \notag
     && + O \big(h \cdot  (\pa A)^2 \big)   + O \big(  h  \cdot  A^2 \cdot \pa A \big)     + O \big(  h   \cdot  A^4 \big)  \,  .
\eea
As a result,
    \bea
\notag
&& | g^{\alpha\beta} \derm_\alpha \derm_\beta h_{\mu\nu} |  \\
\notag
     &\les&  | P(\pa_\mu h,\pa_\nu h) | + | Q_{\mu\nu}(\pa h,\pa h)  | + | G_{\mu\nu}(h)(\pa h,\pa h)  | \\
     \notag
   &&+ \big(   | \derm  A_{L}  | +  | \rderm A   | \big) \cdot  |\derm  A  | \\
 \notag
 &&+ \big(  |\derm_\mu  A_{e_{a}}  | +  |\rderm  A  | \big)  \cdot \big(  |\derm_\nu  A_{e_{a}}  | +  |\rderm  A  | \big) \\
    \notag
       &&+ \big(   | \derm  A_{L}  | +  | \rderm A   | \big) \cdot  \big(   | \derm  A_{L}  | +  | \rderm A   | \big) \\
\notag
   &&+   \big(   | \rderm  A  | +  | \derm A   | \big)  \cdot    |\rderm  A  |  \\
   \notag
 && +  | \derm A   |  \cdot | A | \cdot | A_L | + \big(   | \derm A_L  | +  | \rderm A  | \big) \cdot | A |^2   \\
 \notag
&& + \big(   | \derm_\mu  A  | +  | \rderm A   | \big) \cdot | A_\nu | \cdot | A | + \big(   | \derm_\nu  A  | +  | \rderm A   | \big) \cdot | A_\mu | \cdot | A | \\
\notag
&& + \big(   | \derm  A_{L}  | +  | \rderm A   | \big) \cdot | A | \cdot  | A_L |  \\
\notag
&& +  | \derm A   |  \cdot | A | \cdot  | A_L |  +  | \rderm A   |  \cdot | A |^2  \\
 \notag
&& +  | A_\mu | \cdot | A_\nu | \cdot  | A |^2 \\
 \notag
 && + | A_L | \cdot  | A |^3 \\
  \notag
     && + O \big(h \cdot  (\pa A)^2 \big)   + O \big(  h  \cdot  A^2 \cdot \pa A \big)     + O \big(  h   \cdot  A^4 \big)  \\
     \notag
         &\les&  | P(\pa_\mu h,\pa_\nu h) | + | Q_{\mu\nu}(\pa h,\pa h)  | + | G_{\mu\nu}(h)(\pa h,\pa h)  | \\
 \notag
    &&+ \big(   | \derm  A_{L}  | +  | \rderm A   | \big) \cdot  |\derm  A  | \\
    \notag
 &&+ \big(  |\derm_\mu  A_{e_{a}}  | +  |\rderm  A  | \big)  \cdot \big(  |\derm_\nu  A_{e_{a}}  | +  |\rderm  A  | \big) \\
   \notag
 && +  | \derm A   |  \cdot | A | \cdot | A_L | + \big(   | \derm A_L  | +  | \rderm A  | \big) \cdot | A |^2   \\
 \notag
&& + \big(   | \derm_\mu  A  | +  | \rderm A   | \big) \cdot | A_\nu | \cdot | A | + \big(   | \derm_\nu  A  | +  | \rderm A   | \big) \cdot | A_\mu | \cdot | A | \\
\notag
&& + \big(   | \derm  A_{L}  | +  | \rderm A   | \big) \cdot | A | \cdot  | A_L |  \\
\notag
&& +  | \derm A   |  \cdot | A | \cdot  | A_L |  +  | \rderm A   |  \cdot | A |^2  \\
 \notag
&& +  | A_\mu | \cdot | A_\nu | \cdot  | A |^2  + | A_L | \cdot  | A |^3 \\
  \notag
     && + O \big(h \cdot  (\pa A)^2 \big)   + O \big(  h  \cdot  A^2 \cdot \pa A \big)     + O \big(  h   \cdot  A^4 \big)  \,  .
\eea

\end{proof}

\section{The Lorenz gauge estimate}

\begin{lemma}\label{estimateoncovariantgradientofAL}
We have in the Lorenz gauge that the potential $A$ satisfies the following inequality
\bea
| \derm  A_{L}  | &\les&  | \rderm A  | \, + O \big(  | h  | \cdot  | \derm  A  | \big) \; ,
\eea
where
\bea
 | \derm A_\mu  |^2  :=  E^{\a\b} < \derm_{\b} A_{\mu} \, ,  \derm_{\a} A_{\mu} > \, ,
\eea
and where $L$ and $\rderm$ are constructed using wave coordiantes.
\end{lemma}

\begin{proof}

We showed that in wave coordinates, the Lorenz gauge can be written as 
\beaa
\pa^{\a}  A_{\a}   & =& 0 \\
&=& g^{\a\b}\pa_{\b}  A_{\a} \, .
 \eeaa
 We also showed that
 \beaa
g^{\a\b} = m^{\a\b}-h^{\a\b}  +O^{\a\b}(h^2) \, .
\eeaa
Hence,
\beaa
0 = g^{\a\b}\pa_{\b}  A_{\a}  &=& m^{\a\b}\pa_{\b}  A_{\a}   + \big( -h^{\a\b} +O^{\a\b}(h^2) \big) \cdot \pa_{\b}  A_{\a} \, ,
 \eeaa
and thus, in both the Lorenz gauge and in wave coordinates, we have 
\beaa
 m^{\a\b}\pa_{\b}  A_{\a}  &=&   \big( -h^{\a\b} +O^{\a\b}(h^2) \big) \cdot  \pa_{\b}  A_{\a}  .
 \eeaa
We could then write

\beaa
 m^{\a\b} { \der^{(\bf{m})}}_{\b}  A_{\a}  &=&   \big( -h^{\a\b} +O^{\a\b}(h^2) \big) \cdot \derm_{\b}  A_{\a} \\
&=& O \big( h  \cdot  \pa  A  \big) .
 \eeaa
Note that this expression is coordinate-independent. Decomposing in the null-frame $\cal U$, that is constructed using wave coordinates, we obtain
 \beaa
 m^{\a\b}  { \der^{(\bf{m})}}_{\b}  A_{\a}  &=& m^{\a L}  { \der^{(\bf{m})}}_{L}  A_{\a} + m^{\a\underline{L}}  { \der^{(\bf{m})}}_{\underline{L}}  A_{\a} + m^{\a e_A}  { \der^{(\bf{m})}}_{e_A}  A_{\a}\\
 &=& m^{\underline{L} L}  { \der^{(\bf{m})}}_{L}  A_{\underline{L}} + m^{L\underline{L}}  { \der^{(\bf{m})}}_{\underline{L}}  A_{L} + m^{e_A e_A}  { \der^{(\bf{m})}}_{e_A}  A_{e_A}  \, ,
 \eeaa
and thus,
 \beaa
 m^{\a\b}  { \der^{(\bf{m})}}_{\b}  A_{\a}  &=& - \frac{1}{2}  { \der^{(\bf{m})}}_{L}  A_{\underline{L}} - \frac{1}{2}  { \der^{(\bf{m})}}_{\underline{L}}  A_{L} +  { \der^{(\bf{m})}}_{e_1}  A_{e_1}+   { \der^{(\bf{m})}}_{e_{2}}  A_{e_{2}} \, .
 \eeaa

Hence,
\beaa
- \frac{1}{2}  { \der^{(\bf{m})}}_{L}  A_{\underline{L}} - \frac{1}{2}  { \der^{(\bf{m})}}_{\underline{L}}  A_{L} +  { \der^{(\bf{m})}}_{e_1}  A_{e_1} -  { \der^{(\bf{m})}}_{e_{2}}  A_{e_{2}} &=& O \big( h  \cdot  \pa  A  \big)  .
\eeaa
Consequently 
\bea
\notag
 - \frac{1}{2}  { \der^{(\bf{m})}}_{\underline{L}}  A_{L}   &=& O \big( h  \cdot  \pa  A  \big)  +  \frac{1}{2}  { \der^{(\bf{m})}}_{L}  A_{\underline{L}} -  { \der^{(\bf{m})}}_{e_1}  A_{e_1}-   { \der^{(\bf{m})}}_{e_{2}}  A_{e_{2}}  \, , \\
\eea
which implies that
\beaa
\frac{1}{2} |  { \der^{(\bf{m})}}_{\underline{L}}  A_{L}  | &\leq& O \big(  | h  | \cdot  | \pa  A  | \big) + \frac{1}{2}  |  { \der^{(\bf{m})}}_{L}  A_{\underline{L}}  | +  |  { \der^{(\bf{m})}}_{e_1}  A_{e_1}  | +     |  { \der^{(\bf{m})}}_{e_{2}}  A_{e_{2}}  | \; .
\eeaa

Thus,
\beaa
 | \derm_{\underline{L}}  A_{L}   | &\les& O \big(  | h | \cdot  | \pa  A  | \big) +  | \rderm A  | \; .
\eeaa
Also, since 
\beaa
 | \derm_{L}  A_{L}   |^2  +  | \derm_{e_1}  A_{L}   |^2  + | \derm_{e_2}  A_{L}   |^2  &\les&   | \rderm A  |^2 \, ,
\eeaa
we get
\beaa
 | \derm  A_{L}   | &\les& O \big(  | h  | \cdot  |\pa  A  | \big) +  | \rderm A  | \; ,
\eeaa
where
\beaa
 | \derm A_L  |^2  :=  E^{\a\b} < \derm_{\b} A_L ,  \derm_{\a} A_L > .
\eeaa

\end{proof}

  We remind the following lemma and notation from our previous paper \cite{G4}.
       \begin{lemma}\label{LiederivativeZofMinkwoskimetric}
      
      The Lie derivative in the direction of the Minkowski vector fields of the Minkowski metric, is either null or proportional to the Minkowski metric, that is for any $Z \in {\cal Z}$,
\bea
 \Lie_{Z}  m_{\mu\nu}  =  c_Z \cdot m_{\mu\nu} 
\eea
and 
\bea
 \Lie_{Z}  m^{\mu\nu}  =  - c_Z \cdot m^{\mu\nu} 
\eea
where $c_Z = 0 $ for all $Z \neq S$ and $c_S = 2$.

Thus,
\bea
 \Lie_{Z^I}  m_{\mu\nu}  =  c (I)  \cdot m_{\mu\nu} \; \\
 \Lie_{Z^I}  m^{\mu\nu}  = \hat{c} (I) \cdot m^{\mu\nu} \; ,
\eea
where $c(I)$ and $\hat{c}(I)$ are constants that depend on $Z^I$.

      \end{lemma}

\begin{lemma}\label{LorenzgaugeestimateforgradientofLiederivativesofAL}
We have in the Lorenz gauge that the potential $A$ satisfies the following inequality
\bea
\notag
 | \derm ( \Lie_{Z^J}  A_{L}  ) | &\les& \sum_{|I| \leq |J|} | \rderm ( \Lie_{Z^I} A ) | \, + \sum_{|K| + |M| \leq |J|} O \big(  | \Lie_{Z^K}  h   | \cdot  | \derm  (  \Lie_{Z^M}  A ) | \big)  \;, \\
\eea
where
\bea
 | \derm ( \Lie_{Z^J}  A_\mu )   |^2  :=  E^{\a\b} < \derm_{\b}( \Lie_{Z^J}   A_{\mu} ) \, ,  \derm_{\a} ( \Lie_{Z^J}  A_{\mu} )  > \; ,
\eea
and where $L$ and $\rderm$ are constructed using wave coordinates.
\end{lemma}

\begin{proof}

We already showed in the proof of \ref{estimateoncovariantgradientofAL}, that the Lorenz gauge reads
\beaa
 m^{\a\b} { \der^{(\bf{m})}}_{\b}  A_{\a}  &=&   \big( -h^{\a\b} +O^{\a\b}(h^2) \big) \cdot \derm_{\b}  A_{\a} \\
&=& O \big( h  \cdot  \derm  A  \big) .
 \eeaa
 
 We are going to prove the inequality by induction on $|J|$. For $|J| = 1$, by differentiating, we get
  \beaa
 \Lie_{Z}   (  m^{\a\b} { \der^{(\bf{m})}}_{\b}  A_{\a}  ) &=&   \Lie_{Z}  O \big( h  \cdot  \derm  A  \big) =  \sum_{|K| + |L| \leq 1} O \big(  \Lie_{Z^K} h  \cdot  \Lie_{Z^L} \derm  A  \big) \; .
  \eeaa
 
 Computing directly the left hand side, we get
 \beaa
 \Lie_{Z}   (  m^{\a\b} { \der^{(\bf{m})}}_{\b}  A_{\a}  ) &=& ( \Lie_{Z} m^{\a\b} ) \cdot (    { \der^{(\bf{m})}}_{\b}  A_{\a}  ) +  m^{\a\b} \cdot (  \Lie_{Z}   { \der^{(\bf{m})}}_{\b}  A_{\a}  )   \\
 &=&  - c_Z \cdot m^{\a\b}  { \der^{(\bf{m})}}_{\b}  A_{\a} +  m^{\a\b} \cdot   { \der^{(\bf{m})}}_{\b}  ( \Lie_{Z}   A_{\a} )    \\
&=& O \big( h  \cdot  \pa  A  \big) +  m^{\a\b} \cdot   { \der^{(\bf{m})}}_{\b}  (  \Lie_{Z}   A_{\a} )   .
 \eeaa
And we have
 \beaa
 m^{\a\b}  { \der^{(\bf{m})}}_{\b} (  \Lie_{Z}  A_{\a} )  &=& - \frac{1}{2}  { \der^{(\bf{m})}}_{L}  ( \Lie_{Z}  A_{\underline{L}} ) - \frac{1}{2}  { \der^{(\bf{m})}}_{\underline{L}} ( \Lie_{Z}   A_{L} ) +  { \der^{(\bf{m})}}_{e_1} (  \Lie_{Z}  A_{e_1} ) +   { \der^{(\bf{m})}}_{e_{2}} (   \Lie_{Z}  A_{e_{2}} ) \, .
 \eeaa
 
 Thus, we obtain 
  \beaa
  &&\sum_{|K| + |L| \leq 1} O \big(  \Lie_{Z^K} h  \cdot  \Lie_{Z^L} \derm  A  \big) \\
  &=& O \big( h  \cdot  \pa  A  \big)   - \frac{1}{2}  { \der^{(\bf{m})}}_{L}  ( \Lie_{Z}  A_{\underline{L}} ) - \frac{1}{2}  { \der^{(\bf{m})}}_{\underline{L}} ( \Lie_{Z}   A_{L} ) +  { \der^{(\bf{m})}}_{e_1} (  \Lie_{Z}  A_{e_1} ) +   { \der^{(\bf{m})}}_{e_{2}} (   \Lie_{Z}  A_{e_{2}} ) \\
  \eeaa
  which yields to
    \beaa
     \frac{1}{2}  { \der^{(\bf{m})}}_{\underline{L}} ( \Lie_{Z}   A_{L} )  &=&\sum_{|K| + |L| \leq |J|} O \big(  \Lie_{Z^K} h  \cdot  \Lie_{Z^L} \derm  A  \big)  + O \big( h  \cdot  \derm  A  \big) \\
  &&-  \frac{1}{2}  { \der^{(\bf{m})}}_{L}  ( \Lie_{Z}  A_{\underline{L}} )  +  { \der^{(\bf{m})}}_{e_1} (  \Lie_{Z}  A_{e_1} ) +   { \der^{(\bf{m})}}_{e_{2}} (   \Lie_{Z}  A_{e_{2}} ) \; .
  \eeaa
  Hence,
    \beaa
| { \der^{(\bf{m})}}_{L}  ( \Lie_{Z}  A_{\underline{L}} )  | &\les&\sum_{|K| + |L| \leq |J|} O \big( | \Lie_{Z^K} h  | \cdot | \Lie_{Z^L} \derm  A | \big) + | \rderm  ( \Lie_{Z}   A ) |  \; .
  \eeaa
  
  Assuming now that the inequality is true for $|J| -1$, and using the notation that we defined in such a way that it would allow us to write
 \bea
Z^I(f \cdot g)=\sum_{I_1+I_2=I} (Z^{I_1} f) \cdot (Z^{I_2} g) ,
\eea
we can then write
  \beaa
 \Lie_{Z^J}   (  m^{\a\b} { \der^{(\bf{m})}}_{\b}  A_{\a}  ) &=&   \Lie_{Z^J}  O \big( h  \cdot  \derm  A  \big) =  \sum_{|K| + |L| \leq |J|} O \big(  \Lie_{Z^K} h  \cdot  \Lie_{Z^L} \derm  A  \big) \; .
  \eeaa
 
 Computing directly the left hand side, we get
 \beaa
 \Lie_{Z^J}   (  m^{\a\b} { \der^{(\bf{m})}}_{\b}  A_{\a}  ) &=&\sum_{J_1+J_2= J}  ( \Lie_{Z^{J_1}} m^{\a\b} ) \cdot  (  \Lie_{Z^{J_2}}   { \der^{(\bf{m})}}_{\b}  A_{\a}  )   \\
 &=& \sum_{J_1+J_2= J} \hat{c} (J_1) \cdot m^{\a\b} \cdot   { \der^{(\bf{m})}}_{\b}  ( \Lie_{Z^{J_2}}   A_{\a} )    \; .
 \eeaa
And we have
 \beaa
 m^{\a\b}  { \der^{(\bf{m})}}_{\b} (  \Lie_{Z^{J_2}}  A_{\a} )  &=& - \frac{1}{2}  { \der^{(\bf{m})}}_{L}  ( \Lie_{Z^{J_2}}  A_{\underline{L}} ) - \frac{1}{2}  { \der^{(\bf{m})}}_{\underline{L}} ( \Lie_{Z^{J_2}}   A_{L} ) +  { \der^{(\bf{m})}}_{e_1} (  \Lie_{Z^{J_2}}  A_{e_1} ) +   { \der^{(\bf{m})}}_{e_{2}} (   \Lie_{Z^{J_2}}  A_{e_{2}} ) \, .
 \eeaa
 
 Thus, we obtain 
  \beaa
  &&\sum_{|K| + |L| \leq |J|} O \big(  \Lie_{Z^K} h  \cdot  \Lie_{Z^L} \derm  A  \big) \\
  &=& \sum_{J_1+J_2= J} \hat{c} (J_1) \cdot \Big(  - \frac{1}{2}  { \der^{(\bf{m})}}_{L}  ( \Lie_{Z^{J_2}}  A_{\underline{L}} ) - \frac{1}{2}  { \der^{(\bf{m})}}_{\underline{L}} ( \Lie_{Z^{J_2}}   A_{L} ) +  { \der^{(\bf{m})}}_{e_1} (  \Lie_{Z^{J_2}}  A_{e_1} ) +   { \der^{(\bf{m})}}_{e_{2}} (   \Lie_{Z^{J_2}}  A_{e_{2}} )  \Big) \; ,
  \eeaa
which yields to
    \beaa
    && \frac{1}{2}  { \der^{(\bf{m})}}_{\underline{L}} ( \Lie_{Z^J}   A_{L} )  \\
     &=&\sum_{|K| + |L| \leq |J|} O \big(  \Lie_{Z^K} h  \cdot  \derm  (  \Lie_{Z^L} A )   \big)   \\
&& + \sum_{J_1+J_2= J , \; |J_2| < |J|} \hat{c} (J_1) \cdot \Big(     - \frac{1}{2}  { \der^{(\bf{m})}}_{L}  ( \Lie_{Z^{J_2}}  A_{\underline{L}} ) - \frac{1}{2}  { \der^{(\bf{m})}}_{\underline{L}} ( \Lie_{Z^{J_2}}   A_{L} ) \\
&&+  { \der^{(\bf{m})}}_{e_1} (  \Lie_{Z^{J_2}}  A_{e_1} )   +   { \der^{(\bf{m})}}_{e_{2}} (   \Lie_{Z^{J_2}}  A_{e_{2}} ) \Big) \\
  &&-  \frac{1}{2}  { \der^{(\bf{m})}}_{L}  ( \Lie_{Z^J}  A_{\underline{L}} )  +  { \der^{(\bf{m})}}_{e_1} (  \Lie_{Z^J}  A_{e_1} ) +   { \der^{(\bf{m})}}_{e_{2}} (   \Lie_{Z^J}  A_{e_{2}} ) \; .
  \eeaa
Using the induction hypothesis, we obtain
    \beaa
| { \der^{(\bf{m})}}_{L}  ( \Lie_{Z^J}  A_{\underline{L}} )  | &\les&\sum_{|K| + |L| \leq |J|} O \big( | \Lie_{Z^K} h  | \cdot | \Lie_{Z^L} \derm  A | \big) + \sum_{|I| \leq |J|} | \rderm  ( \Lie_{Z^I}   A ) |  \; .
  \eeaa
   
\end{proof}

\section{The harmonic gauge estimate}

The wave coordinates, or harmonic gauge, reads for $\mu \in \{0, 1, 2, 3 \}$,
\beaa
\Box_{\g} x^\mu = \der^{\a}  \der_{\a}  x^\mu =  0 .
\eeaa
Now, we could like to write $\Box_{\g}$ differently and for this, we will use the following well-known lemma.

\begin{lemma}\label{writingthewaveoperatorinsystemofcoordiantes}
In a given system of coordinates, we have for a sufficiently smooth function $f$, whether it is a scalar or a tensor, the following equality
\bea
\Box_{\g} f = \frac{1}{ \sqrt{ | \det  \g | } } \pa_\b \big( g^{\a\b}  \sqrt{ | \det  \g | } \pa_\a f \big) .
\eea

\end{lemma}

\begin{proof}
We will show that for any sufficiently smooth function $b$ with compact support, we have

\beaa 
\int_{M} ( \Box_{\g} f ) b  \cdot \mbox{dVol} = \int_{M}   \frac{1}{ \sqrt{ | \mbox{det} \, \g | } } \pa_\b \big( g^{\a\b}  \sqrt{ | \mbox{det} \, \g | } \pa_\a f \big) b \cdot  \mbox{dVol}  ,
\eeaa
where $ \mbox{dVol}$ is the volume form on $M$, that is $\mbox{dVol} = \sqrt{ | \mbox{det} \, \g | } \mbox{dx}  $ . We assume here for simplicity that $f$ and $b$ are scalar functions, otherwise we consider the full contraction on the indices of $f $ and $b$ and the same argument by then holds.

Computing, on one hand, the right hand side, we get
\beaa 
 && \int_{M}   \frac{1}{ \sqrt{ | \mbox{det} \, \g | } } \pa_\b \big( g^{\a\b}  \sqrt{ | \mbox{det} \, \g | } \pa_\a f \big) b \cdot  \mbox{dVol} \\
  &=& \int_{M}   \frac{1}{ \sqrt{ | \mbox{det} \, \g | } } \pa_\b \big( g^{\a\b}  \sqrt{ | \mbox{det} \, \g | } \pa_\a f \big) b \cdot   \sqrt{ | \mbox{det} \, \g | } \mbox{dx} \\
 &=&  \int_{M}   \pa_\b \big( g^{\a\b}  \sqrt{ | \mbox{det} \, \g | } \pa_\a f \big) b \cdot  \mbox{dx}  \\
 &=&  \int_{M}  -  \big( g^{\a\b}  \sqrt{ | \mbox{det} \, \g | } \pa_\a f \big) \pa_\b b \cdot  \mbox{dx} ,
\eeaa
where we integrated by parts and used the fact that there are no boundary terms since $b$ is compactly supported. Proceeding we could then write
\beaa 
 \int_{M}   \frac{1}{ \sqrt{ | \mbox{det} \, \g | } } \pa_\b \big( g^{\a\b}  \sqrt{ | \mbox{det} \, \g | } \pa_\a f \big) b \cdot  \mbox{dVol} &=&  \int_{M}  -   g^{\a\b} \sqrt{ | \mbox{det} \, \g | } \pa_\a f  \pa_\b b \cdot  \mbox{dx} \\
 &=&  \int_{M}  -   g^{\a\b} \sqrt{ | \mbox{det} \, \g | } \der_\a f  \der_\b b \cdot  \mbox{dx} \\
 &=& -  \int_{M}   \der^\b f  \der_\b b \cdot \mbox{dVol} .
\eeaa
However, on the other hand, we compute
\beaa 
\int_{M} ( \Box_{\g} f ) b  \cdot \mbox{dVol} &=& \int_{M} (  g^{\a\b} \der_{\b}  \der_{\a}  f ) b  \cdot \mbox{dVol} \\
&=& \int_{M} \der_{\b}   (  g^{\a\b}  \der_{\a}  f ) b  \cdot \mbox{dVol} \\
&& \text{(where we used that $\der g = 0$)} .
\eeaa
By integrating by parts, and since it is a trace, the covariant derivative can then be passed to $b$ also as covariant derivative, and we will have no boundary terms since $b$ is compactly supported, and hence we obtain
\beaa 
\int_{M} ( \Box_{\g} f ) b  \cdot \mbox{dVol} &=&  - \int_{M}  (  g^{\a\b}  \der_{\a}  f )  \der_{\b} b  \cdot \mbox{dVol} .
\eeaa
\end{proof}

\begin{lemma}\label{writingthewavecoordianteconditioninsystemofcoordiantes}
The wave coordinate condition reads
\bea
\pa_\b \big(  h^{\a\b} -  \frac{ 1}{2}  m^{\a\b} h_\la^{\,\, \, \la} + O (|h|^2 )   \big)  &=& 0 ,
\eea
or 
\bea
 \pa_\b \big(  H^{\a\b} -  \frac{ 1}{2}  m^{\a\b} H_\la^{\,\, \, \la} + O (|H|^2 )   \big)   &=& 0,
\eea
where the traces are taken with respect to the background metric $m$.

\end{lemma}

\begin{proof}

Based on what we showed in Lemma \ref{writingthewaveoperatorinsystemofcoordiantes}, we can write the wave coordinate condition as
\beaa
  \frac{1}{ \sqrt{ | \mbox{det} \, \g | } } \pa_\b \big( g^{\a\b}  \sqrt{ | \mbox{det} \, \g | } \pa_\a x^\mu ) =  0 ,
\eeaa
where $x^\mu$ is a scalar and not a tensor. Computing the contraction on $\a$ and $\b$ in the wave system of coordinates, i.e. $\a, \b \in \{x^0, x^1, x^2, x^3 \}$, and since then $\pa_\a x^\mu = \delta_\a^{\mu} \, ,$ we get
 \beaa
  \frac{1}{ \sqrt{ | \mbox{det} \, \g | } } \pa_\b \big( g^{\a\b}  \sqrt{ | \mbox{det} \, \g | } ) &=&  0 ,
\eeaa
which leads to the fact that the wave coordinate condition implies
 \bea
\pa_\b \big( g^{\a\b}  \sqrt{ | \mbox{det} \, \g | } ) &=&  0 .
\eea

Using the formula for the determinant of the sum of two matrices of dimension $4$, we obtain
\beaa
\mbox{det} \, \g &=& \mbox{det} \, ( {\bf m} + {\bf h})  = \mbox{det} \, ( {\bf m} ) +  \mbox{det} \, ( {\bf  m }) \cdot  \mbox{Tr} ( {\bf  m^{-1} h } ) + O (|h|^2 )  \\
&=&  -1 -   \mbox{Tr} ( { \bf m^{-1} h} ) + O (|h|^2 ) ,
\eeaa
where 
\beaa 
 \mbox{Tr} (  { \bf m^{-1} h} ) &=& \mbox{Tr} ( { \bf m^{\a\b}  h_{\a\mu} } ) =  { \bf m^{\a\mu}  h_{\a\mu} } \\
&=& \mbox {Tr} \, h \\
&& \text{(where the trace in the last equality is taken with respect to the background metric $m$)} .
\eeaa
Thus,
\beaa
|\mbox{det} \, \g | &=&  1 +  h_\a^{\,\, \, \a} + O (|h|^2 ) ,
\eeaa
and
\bea
\notag
 \sqrt{ | \mbox{det} \, \g | }  &=&  \big( 1 +  h_\a^{\,\, \, \a} + O (|h|^2 ) \big)^{\frac{1}{2} } \\
 \notag
  &=&   1 +  \frac{1}{2} h_\a^{\,\, \, \a} + O (|h|^2 ) + O ( (h + h^2)^2 )  \\
  &= &1 +  \frac{1}{2} h_\a^{\,\, \, \a} + O (|h|^2 ) .
\eea

On the other hand, we have
\beaa
g^{\a\b} = m^{\a\b} + H^{\a\b} .
\eeaa
Hence,
 \beaa
 g^{\a\b}  \sqrt{ | \mbox{det} \, \g | }  &=&  ( m^{\a\b} + H^{\a\b} ) ( 1 +  \frac{1}{2} h_\la^{\,\, \, \la} + O (|h|^2 ) )  \\
 &=& m^{\a\b} +  \frac{ 1}{2}  m^{\a\b} h_\a^{\,\, \, \a} + O (|h|^2 ) + H^{\a\b} + \frac{1}{2} H^{\a\b}  h_\la^{\,\, \, \la} + O (| H | \cdot |h|^2 ) ) .
\eeaa
Since 
\beaa
H^{\mu\nu} = -h^{\mu\nu}+O^{\mu\nu}(h^2) ,
\eeaa
we obtain 
 \beaa
 g^{\a\b}  \sqrt{ | \mbox{det} \, \g | } &=& m^{\a\b} + H^{\a\b} +  \frac{ 1}{2}  m^{\a\b} h_\la^{\,\, \, \la} + O (|h|^2 )  \\
 &=& m^{\a\b}  - h^{\a\b} +  \frac{ 1}{2}  m^{\a\b} h_\la^{\,\, \, \la} + O (|h|^2 ) \; .
\eeaa
Therefore, we can re-write the wave coordinate condition as 
 \bea
 \notag
 0 &=& \pa_\b \big( m^{\a\b}  - h^{\a\b} +  \frac{ 1}{2}  m^{\a\b} h_\la^{\,\, \, \la} + O (|h|^2 )   \big) \\
 &=& \pa_\b \big(  - h^{\a\b} +  \frac{ 1}{2}  m^{\a\b} h_\la^{\,\, \, \la} + O (|h|^2 )   \big) \\
 \notag
 && \text{(where the trace is taken with respect to the background metric $m$)} .
\eea

Now, using the formula for the determinant of the sum of two matrices of dimension $4$ applied for ${\bf g^{-1}}$, instead of ${\bf g }$, we obtain
\beaa
\mbox{det} \, {\bf g^{-1}} &=& \mbox{det} \, ( {\bf m} + {\bf H})  = \mbox{det} \, ( {\bf m} ) +  \mbox{det} \, ( {\bf  m }) \cdot  \mbox{Tr} ( {\bf  m^{-1} H} ) + O (|H|^2 )  \\
&=&  -1 -   \mbox{Tr} ( { \bf m^{-1} H} ) + O (|H|^2 ) ,
\eeaa
where 
\beaa 
 \mbox{Tr} (  { \bf m^{-1} H} ) &=& \mbox{Tr} ( { \bf m^{\a\b}  H_{\a\mu} } ) =  { \bf m^{\a\mu}  H_{\a\mu} } \\
&=& \mbox {Tr} \, H \\
&& \text{(where the trace is taken with respect to the background metric $m$)} .
\eeaa
Thus,
\beaa
|\mbox{det} \,  {\bf g^{-1}} | &=&  1 +  H_\a^{\,\, \, \a} + O (|H|^2 ) .
\eeaa
However,
\beaa
|\mbox{det} \,  {\bf g} | = |\mbox{det} \,  {\bf g^{-1}} |^{-1} &=& \big(  1 +  H_\a^{\,\, \, \a} + O (|H|^2 ) \big)^{-1} = 1 -  H_\a^{\,\, \, \a} + O (|H|^2 ) ,
\eeaa
and consequently, 
\bea
\notag
 \sqrt{ | \mbox{det} \, \g | }  &=&  \big( 1 - H_\a^{\,\, \, \a} + O (|H|^2 ) \big)^{\frac{1}{2} } \\
 \notag
  &=&   1 -  \frac{1}{2} H_\a^{\,\, \, \a} + O (|H|^2 ) + O ( (H + H^2)^2 )  \\
  &= &1 -  \frac{1}{2} H_\a^{\,\, \, \a} + O (|H|^2 ) .
\eea

Thus,
 \beaa
 g^{\a\b}  \sqrt{ | \mbox{det} \, \g | }  &=&  ( m^{\a\b} + H^{\a\b} ) ( 1 -  \frac{1}{2} H_\la^{\,\, \, \la} + O (|H|^2 ) )  \\
 &=& m^{\a\b} -  \frac{ 1}{2}  m^{\a\b} H_\a^{\,\, \, \a} + O (|H|^2 ) + H^{\a\b} - \frac{1}{2} H^{\a\b}  H_\la^{\,\, \, \la} + O (  |H|^3 ) ) .
\eeaa

We get 
 \beaa
 g^{\a\b}  \sqrt{ | \mbox{det} \, \g | }  &=& m^{\a\b} + H^{\a\b} -  \frac{ 1}{2}  m^{\a\b} H_\la^{\,\, \, \la} + O (|H|^2 )  .
\eeaa
Hence, we can re-write the wave coordinate condition as 
 \bea
 \notag
 0 &=& \pa_\b \big( m^{\a\b} + H^{\a\b} -  \frac{ 1}{2}  m^{\a\b} H_\la^{\,\, \, \la} + O (|H|^2 )   \big) \\
 &=& \pa_\b \big(  H^{\a\b} -  \frac{ 1}{2}  m^{\a\b} H_\la^{\,\, \, \la} + O (|H|^2 )   \big) \\
 \notag
  && \text{(where the trace is taken with respect to the background metric $m$)} .
\eea

\end{proof}

\begin{lemma}\label{waveconditionestimateonzeroLiederivativeofmetric}
We have
\beaa
| \derm  H_{\cal T L} | &\les& | \rderm  H |  + O (|H| \cdot |\derm H| ) \; ,
\eeaa
and
\beaa
| \derm  h_{\cal T L} | &\les& | \rderm  h |  + O (|h| \cdot |\derm h| ) \; .
\eeaa

\end{lemma}

\begin{proof}
In wave coordinates, we have
\beaa
 0 &=& \pa_\b \big(  h^{\a\b} -  \frac{ 1}{2}  m^{\a\b} h_\la^{\,\, \, \la} + O^{\a\b} (|h|^2 )   \big) \\
 &=&  \derm_\b \big(  h^{\a\b} -  \frac{ 1}{2}  m^{\a\b} h_\la^{\,\, \, \la} + O^{\a\b} (|h|^2 )   \big) \\
 &=&  \derm_L \big(  h^{\a L} -  \frac{ 1}{2}  m^{\a L} h_\la^{\,\, \, \la} + O^{\a L} (|h|^2 )   \big)  + \derm_{\underline{L}} \big(  h^{\a \underline{L}} -  \frac{ 1}{2}  m^{\a \underline{L}} h_\la^{\,\, \, \la} + O^{\a \underline{L}} (|h|^2 )   \big) \\
 && + \derm_{e_{A}} \big(  h^{\a e_{A}} -  \frac{ 1}{2}  m^{\a e_{A}} h_\la^{\,\, \, \la} + O^{\a e_{A}} (|h|^2 )   \big) \\
&& \text{(where by $e_A$, we mean a summation on $A \in \{1, 2\}$)} .
\eeaa
Since the above identity is a tensor in $\a$, we can lower the $ \a$-index and write
\beaa
0 &=&  \derm_L \big( m^{L\b } h_{\a \b} -  \frac{ 1}{2} m^{L\b }  m_{\a\b} h_\la^{\,\, \, \la} + m^{L\b } O_{\a \b } (|h|^2 )   \big)  \\
 && + \derm_{\underline{L}} \big( m^{\underline{L} \b }  h_{\a \b} -  \frac{ 1}{2} m^{\underline{L} \b }  m_{\a \b} h_\la^{\,\, \, \la} + m^{\underline{L} \b } O_{\a \b} (|h|^2 )   \big) \\
 && + \derm_{e_{A}} \big(   m^{e_{A} \b } h_{\a \b } -  \frac{ 1}{2}   m^{e_{A} \b } m_{\a e_{A}} h_\la^{\,\, \, \la} +  m^{e_{A} \b } O_{\a \b} (|h|^2 )   \big) \\
  &=&  \derm_L \big( m^{L\underline{L} } h_{\a \underline{L}} -  \frac{ 1}{2} m^{L\underline{L} }  m_{\a \underline{L}} h_\la^{\,\, \, \la} + m^{L\underline{L} } O_{\a \underline{L}} (|h|^2 )   \big)  \\
  && + \derm_{\underline{L}} \big( m^{\underline{L} L }  h_{\a L} -  \frac{ 1}{2} m^{\underline{L} L }  m_{\a L} h_\la^{\,\, \, \la} + m^{\underline{L} L } O_{\a L} (|h|^2 )   \big) \\
 && + \derm_{e_{A}} \big(   m^{e_{A} e_{A} } h_{\a e_{A}} -  \frac{ 1}{2}   m^{e_{A} e_{A} } m_{\a e_{A}} h_\la^{\,\, \, \la} +  m^{e_{A} e_{A} } O_{\a e_{A}} (|h|^2 )   \big) .
\eeaa
We obtain
\beaa
0  &=& - \frac{1}{2} \derm_L \big(  h_{\a \underline{L}} -  \frac{ 1}{2}   m_{\a \underline{L}} h_\la^{\,\, \, \la} +  O_{\a \underline{L}} (|h|^2 )   \big)    - \frac{1}{2} \derm_{\underline{L}} \big(   h_{\a L} -  \frac{ 1}{2}   m_{\a L} h_\la^{\,\, \, \la} +  O_{\a L} (|h|^2 )   \big) \\
 && + \derm_{e_{A}} \big(    h_{\a e_{A}} -  \frac{ 1}{2}    m_{\a e_{A}} h_\la^{\,\, \, \la} +   O_{\a e_{A}} (|h|^2 )   \big) .
\eeaa

Taking $\a \in \cal T$, we get
 \beaa
0  &=& - \frac{1}{2} \derm_L \big(  h_{\cal T \underline{L}} -  \frac{ 1}{2}   m_{\cal T \underline{L}} h_\la^{\,\, \, \la} +  O_{L \underline{L}} (|h|^2 )   \big)    - \frac{1}{2} \derm_{\underline{L}} \big(   h_{\cal T L} -  \frac{ 1}{2}   m_{\cal T L} h_\la^{\,\, \, \la} +  O_{\cal T L} (|h|^2 )   \big) \\
 && + \derm_{e_{A}} \big(    h_{\cal T e_{A}} -  \frac{ 1}{2}    m_{\cal T e_{A}} h_\la^{\,\, \, \la} +   O_{\cal T e_{A}} (|h|^2 )   \big) \\
 &=& - \frac{1}{2} \derm_L \big(  h_{\cal T \underline{L}} -  \frac{ 1}{2}   m_{\cal T \underline{L}} h_\la^{\,\, \, \la} +  O_{\cal T \underline{L}} (|h|^2 )   \big)    - \frac{1}{2} \derm_{\underline{L}} \big(   h_{\cal T L} +  O_{\cal T L} (|h|^2 )   \big) \\
 && + \derm_{e_{A}} \big(    h_{\cal T e_{A}}  -  \frac{ 1}{2}    m_{\cal T e_{A}} h_\la^{\,\, \, \la}+   O_{\cal T e_{A}} (|h|^2 )   \big) .
\eeaa

Finally, we obtain
\bea
\notag
  - \frac{1}{2} \derm_{\underline{L}}    h_{\cal T L} =   \frac{1}{2} \derm_L \big(  h_{\cal T \underline{L}} -  \frac{ 1}{2}   m_{\cal T \underline{L}} h_\la^{\,\, \, \la}  \big)  + \derm_{e_{A}}     ( h_{\cal T e_{A}}  -  \frac{ 1}{2}    m_{\cal T e_{A}} h_\la^{\,\, \, \la})  + O (|h| \cdot |\pa h| ) \; . \\
\eea
Consequently,
\beaa
| \derm_{\underline{L}}    h_{\cal T L} | &\les& | \derm_L \big(  h_{\cal T \underline{L}}  -  m_{\cal T \underline{L}} h_\la^{\,\, \, \la}  \big) | +  | \derm_{e_{A}}  (   h_{\cal T e_{A}}  -  \frac{ 1}{2}    m_{\cal T e_{A}} h_\la^{\,\, \, \la} ) |  + O (|h| \cdot |\pa h| ) .
\eeaa
Thus,
\beaa
| \derm_{\underline{L}}    h_{\cal T L} |  &\les& | \derm_L  h_{\cal T \underline{L}} - \derm_L ( m_{\cal T \underline{L}} h_\la^{\,\, \, \la}  \big) | +  | \derm_{e_{A}}    h_{\cal T e_{A}} -  \frac{ 1}{2} \derm_{e_{A}}  (  m_{\cal T e_{A}} h_\la^{\,\, \, \la} )  |  + O (|h| \cdot |\pa h| ) \; .
\eeaa
Since $\derm m = \derm m^{-1}= 0$, we have
\beaa
  \derm_L ( m_{\cal T \underline{L}} h_\la^{\,\, \, \la} )  = m_{\cal T \underline{L}} \cdot m^{\la\b}  \derm_L h_{\la\b} \; ,
\eeaa
and 
\beaa
  \derm_{e_{A}}  (  m_{\cal T e_{A}} h_\la^{\,\, \, \la} )  = m_{\cal T e_{A}}  \cdot m^{\la\b}  \derm_{e_{A}} h_{\la\b} \; .
\eeaa
Therefore,
\beaa
| \derm_{\underline{L}}    h_{\cal T L} |  &\les& | \derm_L  h_{\cal T \underline{L}} |  + |  m_{\cal T \underline{L}} \cdot m^{\la\b}  \derm_L h_{\la\b}  | +  | \derm_{e_{A}}    h_{\cal T e_{A}} | + | m_{\cal T e_{A}}  \cdot m^{\la\b}  \derm_{e_{A}} h_{\la\b}  |  \\
&& + O (|h| \cdot |\pa h| ) \\
&\les& | \rderm  h |  + O (|h| \cdot |\pa h| )\; ,
\eeaa
and as a result,
\beaa
| \derm  h_{\cal T L} | &\les& | \rderm  h |  + O (|h| \cdot |\pa h| ) \; .
\eeaa
Similarly, using the identity
\beaa
 \derm_\b \big(  H^{\a\b} -  \frac{ 1}{2}  m^{\a\b} H_\la^{\,\, \, \la} + O (|H|^2 )   \big)   &=& 0 \;  ,
\eeaa
we get
\beaa
| \derm  H_{\cal T L} | &\les& | \rderm  H |  + O (|H| \cdot |\pa H| ) \;  .
\eeaa

\end{proof}

\begin{lemma}\label{wavecoordinatesestimateonLiederivativesZonmetric}
We have
\beaa
| \derm ( \Lie_{Z^J}  H_{\cal T L} )  | &\les&\sum_{|K| \leq |J| }  |  \rderm ( \Lie_{Z^K} H ) |  + \sum_{|K|+ |M| \leq |J|}  O (| \Lie_{Z^K} H| \cdot |\derm ( \Lie_{Z^M} H ) | ) \; ,
\eeaa
and
\beaa
| \derm ( \Lie_{Z^J} h_{\cal T L} )  | &\les& \sum_{|K| \leq |J| } | \rderm  ( \Lie_{Z^K}  h)  |  + \sum_{|K|+ |M| \leq |J|}  O (|\Lie_{Z^K} h| \cdot |\derm ( \Lie_{Z^M} h ) | )\; .
\eeaa

\end{lemma}

\begin{proof}
We already showed in Lemma \ref{writingthewavecoordianteconditioninsystemofcoordiantes}, that
\beaa
 0 &=&  \derm_\b \big(  h^{\a\b} -  \frac{ 1}{2}  m^{\a\b} h_\la^{\,\, \, \la} + O^{\a\b} (|h|^2 )   \big) 
 \eeaa
 Differentiating, we obtain
 \beaa
 0 &=& \Lie_{Z^J}  \derm_\b \big(  h^{\a\b} -  \frac{ 1}{2}  m^{\a\b} h_\la^{\,\, \, \la} + O^{\a\b} (|h|^2 )   \big) \\
 &=& \derm_\b \Big( \Lie_{Z^J}  \big(  h^{\a\b} -  \frac{ 1}{2}  m^{\a\b} h_\la^{\,\, \, \la} + O^{\a\b} (|h|^2 )   \big) \Big)\\
  &=& \derm_\b \Big( \Lie_{Z^J}  h^{\a\b} -  \frac{ 1}{2}   \Lie_{Z^J}  ( m^{\a\b}  \cdot h_\la^{\,\, \, \la} )  +  \Lie_{Z^J}  O^{\a\b} (|h|^2 )  \Big) \; .
 \eeaa

Using the notation that we defined in such a way that it would allow us to write
 \bea
Z^I(f \cdot g)=\sum_{I_1+I_2=I} (Z^{I_1} f) \cdot (Z^{I_2} g) ,
\eea
we can then write
 
  \beaa
  \Lie_{Z^J}  (  m^{\a\b} \cdot h_\la^{\,\, \, \la}  ) &=&    \sum_{J_1+J_2=J} (\Lie_{Z^{J_1}} m^{\a \b}  ) \cdot ( \Lie_{Z^{J_2}} h_\la^{\,\, \, \la}  ) \\
  &=&  \sum_{J_1+J_2=J}  \hat{c} ( J_1) \cdot m^{\a \b}   \cdot ( \Lie_{Z^{J_2}} h_\la^{\,\, \, \la} ) \; .
 \eeaa
 
However, we have
 \beaa
  \Lie_{Z^{J_2}}   h_\la^{\,\, \, \la} &=&    \Lie_{Z^{J_2}} (  m^{\la \b} \cdot h_{\la \b}   ) = \sum_{J_3+J_4=J_2} (\Lie_{Z^{J_3}} m^{\la \b}  ) \cdot ( \Lie_{Z^{J_4}} h_{\la \b} ) \\
  &=&  \sum_{J_3+J_4=J_2}  \hat{c} ( J_3) \cdot m^{\la \b}   \cdot ( \Lie_{Z^{J_4}} h_{\la \b} ) \; .
 \eeaa
 Therefore,
   \beaa
  \Lie_{Z^J}  (  m^{\a\b} h_\la^{\,\, \, \la}  ) &=&  \sum_{J_1+J_2=J}  \hat{c} ( J_1) \cdot m^{\a \b}   \cdot \big(  \sum_{J_3+J_4=J_2}  \hat{c} ( J_3) \cdot m^{\la \b}   \cdot ( \Lie_{Z^{J_4}} h_{\la \b} )  \big) \; .
 \eeaa
 
 Hence, we have
  \beaa
 0  &=& \derm_\b \Big( \Lie_{Z^J}  h^{\a\b} -  \frac{ 1}{2}  \cdot  m^{\a \b} \cdot  \sum_{J_1+J_2=J}  \hat{c} ( J_1)    \cdot \big(  \sum_{J_3+J_4=J_2}  \hat{c} ( J_3) \cdot m^{\la \b}   \cdot ( \Lie_{Z^{J_4}} h_{\la \b} )  \big)  \Big)  \\
 && + \derm_\b \sum_{ |K| + |M| \leq |J| }  O^{\a\b} (  \Lie_{Z^K} h \cdot \Lie_{Z^M} h  )  \; .
 \eeaa

Decomposing in the null frame, as we have done earlier in the proof of Lemma \ref{waveconditionestimateonzeroLiederivativeofmetric}, we obtain

 \beaa
0  &=& - \frac{1}{2} \derm_L \big(  \Lie_{Z^J} h_{\cal T \underline{L}} -  \frac{ 1}{2}   m_{\cal T \underline{L}} \cdot  \sum_{J_1+J_2=J}  \hat{c} ( J_1)    \cdot \big(  \sum_{J_3+J_4=J_2}  \hat{c} ( J_3) \cdot m^{\la \b}   \cdot ( \Lie_{Z^{J_4}} h_{\la \b} )  \big)   \big) \\
&&    - \frac{1}{2} \derm_{\underline{L}} \big(  \Lie_{Z^J} h_{\cal T L} -  \frac{ 1}{2}   m_{\cal T L} \cdot  \sum_{J_1+J_2=J}  \hat{c} ( J_1)    \cdot \big(  \sum_{J_3+J_4=J_2}  \hat{c} ( J_3) \cdot m^{\la \b}   \cdot ( \Lie_{Z^{J_4}} h_{\la \b} )  \big)    \big) \\
 && + \derm_{e_{A}} \big(  \Lie_{Z^J}  h_{\cal T e_{A}} -  \frac{ 1}{2}    m_{\cal T e_{A}} \cdot  \sum_{J_1+J_2=J}  \hat{c} ( J_1)    \cdot \big(  \sum_{J_3+J_4=J_2}  \hat{c} ( J_3) \cdot m^{\la \b}   \cdot ( \Lie_{Z^{J_4}} h_{\la \b} )  \big)    \big) \\
&& + \sum_{ |K| + |M| \leq |J|} \Big(  \derm_L O_{L \underline{L}}  (  \Lie_{Z^K} h \cdot \Lie_{Z^M} h  )+ \derm_{\underline{L}} O_{\cal T L}  (  \Lie_{Z^K} h \cdot \Lie_{Z^M} h  ) +  \derm_{e_{A}} O_{\cal T e_{A}}  (  \Lie_{Z^K} h \cdot \Lie_{Z^M} h  )   \Big) \\
 &=& - \frac{1}{2} \derm_L \big( \Lie_{Z^J}  h_{\cal T \underline{L}} -  \frac{ 1}{2}   m_{\cal T \underline{L}} \cdot  \sum_{J_1+J_2=J}  \hat{c} ( J_1)    \cdot \big(  \sum_{J_3+J_4=J_2}  \hat{c} ( J_3) \cdot m^{\la \b}   \cdot ( \Lie_{Z^{J_4}} h_{\la \b} )  \big)   \big)   \\
 && - \frac{1}{2} \derm_{\underline{L}} \big(   \Lie_{Z^J}  h_{\cal T L}    \big)  + \derm_{e_{A}} \big(   \Lie_{Z^J}  h_{\cal T e_{A}}   -  \frac{ 1}{2}    m_{\cal T e_{A}} \cdot  \sum_{J_1+J_2=J}  \hat{c} ( J_1)    \cdot \big(  \sum_{J_3+J_4=J_2}  \hat{c} ( J_3) \cdot m^{\la \b}   \cdot ( \Lie_{Z^{J_4}} h_{\la \b} )  \big)  \big) \\
 && + \sum_{ |K| + |M| \leq |J|} \Big(  \derm_L O_{L \underline{L}}  (  \Lie_{Z^K} h \cdot \Lie_{Z^M} h  )+ \derm_{\underline{L}} O_{\cal T L}  (  \Lie_{Z^K} h \cdot \Lie_{Z^M} h  ) +  \derm_{e_{A}} O_{\cal T e_{A}}  (  \Lie_{Z^K} h \cdot \Lie_{Z^M} h  )   \Big) \\
\eeaa

Using the fact that the covariant derivative $\derm m = 0$, we get

\beaa
\notag
  && - \frac{1}{2} \derm_{\underline{L}} ( \Lie_{Z^J}    h_{\cal T L} ) \\
  &=&   \frac{1}{2} \derm_L \big( \Lie_{Z^J}  h_{\cal T \underline{L}} ) -  \frac{ 1}{2}   m_{\cal T \underline{L}} \cdot  \sum_{J_1+J_2=J}  \hat{c} ( J_1)    \cdot \big(  \sum_{J_3+J_4=J_2}  \hat{c} ( J_3) \cdot m^{\la \b}   \cdot \derm_L ( \Lie_{Z^{J_4}} h_{\la \b} )  \big)     \\
  && + \derm_{e_{A}} ( \Lie_{Z^J}    h_{\cal T e_{A}}  )  -  \frac{ 1}{2}    m_{\cal T e_{A}} \cdot  \sum_{J_1+J_2=J}  \hat{c} ( J_1)    \cdot \big(  \sum_{J_3+J_4=J_2}  \hat{c} ( J_3) \cdot m^{\la \b}   \cdot  \derm_{e_{A}} ( \Lie_{Z^{J_4}} h_{\la \b} )  \big)  \\
 && + \sum_{ |K| + |M| \leq |J|} O \big(  \Lie_{Z^K} h \cdot \derm ( \Lie_{Z^M} h)  \big)  \; .
\eeaa
Consequently,
\beaa
| \derm_{\underline{L}}     ( \Lie_{Z^J}    h_{\cal T L} )  | &\les&  | \rderm  \Lie_{Z^J}  h |  +  \sum_{J_1+J_2=J}  \sum_{J_3+J_4=J_2}   | \rderm  \Lie_{Z^{J_4}}  h |   \\
   &&   + \sum_{|K|+ |L| \leq |J|}   O (| \Lie_{Z^K} h| \cdot |\derm  \Lie_{Z^L} h| ) \; . \\
    &\les& \sum_{|K| \leq |J| }  | \rderm  \Lie_{Z^K}  h |  +  \sum_{J_1+J_2=J}  \sum_{J_3+J_4=J_2}   | \rderm  \Lie_{Z^{J_4}}  h |   \\
   &&   + \sum_{|K|+ |L| \leq |J|}   O (| \Lie_{Z^K} h| \cdot |\derm  \Lie_{Z^L} h| ) \; . \\
\eeaa

Similarly, using the identity
\beaa
 \derm_\b \big(  H^{\a\b} -  \frac{ 1}{2}  m^{\a\b} H_\la^{\,\, \, \la} + O (|H|^2 )   \big)   &=& 0,
\eeaa
and differentiating, using the commutation between $\Lie_{Z^J}$ and $\derm$, we obtain
\beaa
 \derm_\b \big(  \Lie_{Z^J} H^{\a\b} -  \frac{ 1}{2}  \Lie_{Z^J}  ( m^{\a\b} H_\la^{\,\, \, \la} ) + \sum_{|K|+ |L| \leq |J|} O (|\Lie_{Z^K} H| \cdot |\Lie_{Z^L} H|  )   \big)   &=& 0 \; ,
\eeaa
and proceeding similarly as for $h$, we obtain the result.

\end{proof}

\section{Estimating the source terms in the Lorenz and harmonic gauges}

\subsection{Using the gauges estimates}\

\begin{lemma}\label{estimatesonsourcetermsforEinstein-Yang-Millssystem}
In the Lorenz gauge and in wave coordinates, for any $\si \in \cal U $, we have
 \bea
   \notag
 | g^{\la\mu} \derm_{\la}   \derm_{\mu}   A_{\si}  |  
  &\les& | \derm h | \cdot  |\rderm A |       + | \rderm  h | \cdot  |\derm A |  \\
       \notag
           && +  | \rderm  h | \cdot  | A  |^2    \\
           \notag
 &&          +| \derm h_{L \si}   |  \cdot  | A  | \cdot  | A_L |   \\
        \notag
                  && + | \derm h_{e_{a} \si} |  \cdot    | A_{e_{a}} |  \cdot  | A_L |   \\
          \notag
 &&   +   | A_L  | \cdot    | \derm_\si A_{\underline{L}}   |      +    |  A  | \cdot  | A_\si  | \cdot  | A_L |   \\
    \notag
    &&   +  | A  | \cdot   | \rderm A  |     \\
       \notag
 &&  + | A_{e_{a}}  | \cdot    | \derm_\si A_{e_{a}} |    +    |  A_{e_{a}}  |^2 \cdot  | A_\si  |   \\
    \notag
  && + O( h \cdot  \pa h \cdot  \pa A) + O( h \cdot  \pa h \cdot  A^2) + O( h \cdot  A \cdot \pa A) + O( h \cdot  A^3) \, .\\
  \eea

In the Lorenz gauge and in wave coordinates, for any $\mu, \nu \in \cal U$, we have
    \bea
\notag
 | g^{\alpha\beta} \derm_\alpha \derm_\beta h_{\mu\nu} |     &\les&  | P(\derm_\mu h,\derm_\nu h) | + | Q_{\mu\nu}(\derm h,\derm h)  | + | G_{\mu\nu}(h)(\derm h,\derm h)  | \\
 \notag
    && +  | \rderm A   | \cdot  |\derm  A  | \\
    \notag
 &&+ \big(  |\derm_\mu  A_{e_{a}}  | +  |\rderm  A  | \big)  \cdot \big(  |\derm_\nu  A_{e_{a}}  | +  |\rderm  A  | \big) \\
   \notag
 && +  | \derm A   |  \cdot | A | \cdot | A_L | +  | \rderm A  |  \cdot | A |^2    \\
 \notag
&& + \big(   | \derm_\mu  A  | +  | \rderm A   | \big) \cdot | A_\nu | \cdot | A | + \big(   | \derm_\nu  A  | +  | \rderm A   | \big) \cdot | A_\mu | \cdot | A | \\
 \notag
&& +  | A_\mu | \cdot | A_\nu | \cdot  | A |^2  + | A_L | \cdot  | A |^3 \\
     && + O \big(h \cdot  (\pa A)^2 \big)   + O \big(  h  \cdot  A^2 \cdot \pa A \big)     + O \big(  h   \cdot  A^4 \big)  \,  .
\eea

\end{lemma}

\begin{proof}

Using the Lorenz gauge estimate \ref{estimateoncovariantgradientofAL},
\beaa
| \derm  A_{L}  | &\les&  | \rderm A | \, + O \big( | h | \cdot | \pa  A | \big) \; ,
\eeaa
and the wave coordinates estimate \ref{waveconditionestimateonzeroLiederivativeofmetric}, 
\beaa
| \derm  h_{\cal T L} | &\les& | \rderm  h |  + O (|h| \cdot |\pa h| ) \; ,
\eeaa
and plugging this in the estimate on the source term for the wave equation on the Yang-Mills potential, we get
     \beaa
   \notag
 && | g^{\la\mu} \derm_{\la}   \derm_{\mu}   A_{\si}  |  \\
 &\les& | \derm h | \cdot  |\rderm A |     \\
   \notag
 && +   | \derm  h | \cdot  |\derm  A_{L}  | + | \derm  h_{\cal T L} | \cdot  |\derm A |  \\
   \notag
   && +  \big( | \rderm  h |  +| \derm h_{L \si} | \big) \cdot  \big( |\derm  A_{L}  | +|\rderm  A | \big)  \\
      \notag
         && +   | \rderm  h |    \cdot  |\derm  A  | \\
 \notag
        && +  \big( | \rderm  h |   +| \derm h_{L \si}  |  \big) \cdot  | A  | \cdot  | A_L |   \\
   \notag
           && +   | \rderm  h |  \cdot  | A  |^2   \\
           \notag
                  && +  \big( | \rderm  h |   +| \derm h_{e_{a} \si} |  \big) \cdot    | A_{e_{a}} |  \cdot  | A_L |   \\
          \notag
 &&   +   | A_L  | \cdot    | \derm_\si A_{\underline{L}}   |      +    |  A  | \cdot  | A_\si  | \cdot  | A_L |   \\
    \notag
    &&   +  | A  | \cdot   | \derm_\si A_L  |     \\
       \notag
 &&  + | A_{e_{a}}  | \cdot    | \derm_\si A_{e_{a}} |    +    |  A_{e_{a}}  |^2 \cdot  | A_\si  |   \\
    \notag
  && + O( h \cdot  \pa h \cdot  \pa A) + O( h \cdot  \pa h \cdot  A^2) + O( h \cdot  A \cdot \pa A) + O( h \cdot  A^3) \\
     &\les& | \derm h | \cdot  |\rderm A |     \\
   \notag
 && +   | \derm  h | \cdot  |\rderm  A  | + O \big( | h | \cdot |\pa  A | \cdot  | \derm  h | \big) + | \rderm  h | \cdot  |\derm A | + O (|h| \cdot |\pa h| \cdot   |\derm A |  ) \\
   \notag
   && +   | \derm  h |  \cdot  \big( O \big( | h | \cdot | \pa  A | \big) +|\rderm  A | \big)  \\
      \notag
         && +   | \rderm  h |    \cdot  |\derm  A  | \\
 \notag
        && +  \big( | \rderm  h |  +| \derm h_{L \si}   | \big) \cdot  | A  | \cdot  | A_L |   \\
   \notag
           && +   | \rderm  h | \cdot  | A  |^2   \\
             \notag
                  && +  \big( | \rderm  h |   +| \derm h_{e_{a} \si} |  \big) \cdot    | A_{e_{a}} |  \cdot  | A_L |   \\
          \notag
 &&   +   | A_L  | \cdot    | \derm_\si A_{\underline{L}}   |      +    |  A  | \cdot  | A_\si  | \cdot  | A_L |   \\
    \notag
    &&   +  | A  | \cdot  \big( | \rderm A  |  +O( h \cdot   \pa A)   \big)   \\
       \notag
 &&  + | A_{e_{a}}  | \cdot    | \derm_\si A_{e_{a}} |    +    |  A_{e_{a}}  |^2 \cdot  | A_\si  |   \\
    \notag
  && + O( h \cdot  \pa h \cdot  \pa A) + O( h \cdot  \pa h \cdot  A^2) + O( h \cdot  A \cdot \pa A) + O( h \cdot  A^3) \, .\\
  \eeaa
  
  Hence,
    \beaa
   \notag
&& | g^{\la\mu} \derm_{\la}   \derm_{\mu}   A_{\si}  |  \\
 &\les& | \derm h | \cdot  |\rderm A |     \\
   \notag
 && +   | \derm  h | \cdot  |\rderm  A  |  + | \rderm  h | \cdot  |\derm A | + O( h \cdot  \pa h \cdot  \pa A) \\
   \notag
   && +  | \derm  h | \cdot |\rderm  A | + O (|h| \cdot |\pa h| \cdot | \pa  A | ) \big)  \\
      \notag
         && +   | \rderm  h |    \cdot  |\derm  A  | \\
 \notag
        && +  \big( | \rderm  h |  +| \derm h_{L \si}   | \big) \cdot  | A  | \cdot  | A_L |   \\
   \notag
           && +  | \rderm  h | \cdot  | A  |^2      \\
                \notag
                  && +  \big( | \rderm  h |   +| \derm h_{e_{a} \si} |  \big) \cdot    | A_{e_{a}} |  \cdot  | A_L |   \\
          \notag
 &&   +   | A_L  | \cdot    | \derm_\si A_{\underline{L}}   |      +    |  A  | \cdot  | A_\si  | \cdot  | A_L |   \\
    \notag
    &&   +  | A  | \cdot   | \rderm A  |  +O( h \cdot A \cdot   \pa A)     \\
       \notag
 &&  + | A_{e_{a}}  | \cdot    | \derm_\si A_{e_{a}} |    +    |  A_{e_{a}}  |^2 \cdot  | A_\si  |   \\
    \notag
  && + O( h \cdot  \pa h \cdot  \pa A) + O( h \cdot  \pa h \cdot  A^2) + O( h \cdot  A \cdot \pa A) + O( h \cdot  A^3) \\
\notag
    &\les& | \derm h | \cdot  |\rderm A |       + | \rderm  h | \cdot  |\derm A |  \\
       \notag
           && +  | \rderm  h | \cdot  | A  |^2    \\
           \notag
 &&          +| \derm h_{L \si}   |  \cdot  | A  | \cdot  | A_L |   \\
        \notag
                  && + | \derm h_{e_{a} \si} |  \cdot    | A_{e_{a}} |  \cdot  | A_L |   \\
          \notag
 &&   +   | A_L  | \cdot    | \derm_\si A_{\underline{L}}   |      +    |  A  | \cdot  | A_\si  | \cdot  | A_L |   \\
    \notag
    &&   +  | A  | \cdot   | \rderm A  |     \\
       \notag
 &&  + | A_{e_{a}}  | \cdot    | \derm_\si A_{e_{a}} |    +    |  A_{e_{a}}  |^2 \cdot  | A_\si  |   \\
    \notag
  && + O( h \cdot  \pa h \cdot  \pa A) + O( h \cdot  \pa h \cdot  A^2) + O( h \cdot  A \cdot \pa A) + O( h \cdot  A^3) 
  \eeaa

Proceeding by estimating the source term for the wave equation on the metric, and using the Lorenz gauge estimate \ref{estimateoncovariantgradientofAL}, we get

  \beaa
\notag
 && | g^{\alpha\beta} \derm_\alpha \derm_\beta h_{\mu\nu} |   \\
   &\les&  | P(\pa_\mu h,\pa_\nu h) | + | Q_{\mu\nu}(\pa h,\pa h)  | + | G_{\mu\nu}(h)(\pa h,\pa h)  | \\
 \notag
    &&+ \big(   O \big( | h | \cdot | \pa  A | \big) +  | \rderm A   | \big) \cdot  |\derm  A  | \\
    \notag
 &&+ \big(  |\derm_\mu  A_{e_{a}}  | +  |\rderm  A  | \big)  \cdot \big(  |\derm_\nu  A_{e_{a}}  | +  |\rderm  A  | \big) \\
   \notag
 && +  | \derm A   |  \cdot | A | \cdot | A_L | + \big(  O \big( | h | \cdot | \pa  A | \big) +  | \rderm A  | \big) \cdot | A |^2   \\
 \notag
&& + \big(   | \derm_\mu  A  | +  | \rderm A   | \big) \cdot | A_\nu | \cdot | A | + \big(   | \derm_\nu  A  | +  | \rderm A   | \big) \cdot | A_\mu | \cdot | A | \\
\notag
&& + \big(  O \big( | h | \cdot | \pa  A | \big) +  | \rderm A   | \big) \cdot | A | \cdot  | A_L |  \\
\notag
&& +  | \derm A   |  \cdot | A | \cdot  | A_L |  +  | \rderm A   |  \cdot | A |^2  \\
 \notag
&& +  | A_\mu | \cdot | A_\nu | \cdot  | A |^2  + | A_L | \cdot  | A |^3 \\
     && + O \big(h \cdot  (\pa A)^2 \big)   + O \big(  h  \cdot  A^2 \cdot \pa A \big)     + O \big(  h   \cdot  A^4 \big)  \,  \\
       &\les&  | P(\pa_\mu h,\pa_\nu h) | + | Q_{\mu\nu}(\pa h,\pa h)  | + | G_{\mu\nu}(h)(\pa h,\pa h)  | \\
 \notag
    && +  | \rderm A   | \cdot  |\derm  A  | + O \big( | h | \cdot | \pa  A |^2 \big)\\
    \notag
 &&+ \big(  |\derm_\mu  A_{e_{a}}  | +  |\rderm  A  | \big)  \cdot \big(  |\derm_\nu  A_{e_{a}}  | +  |\rderm  A  | \big) \\
   \notag
 && +  | \derm A   |  \cdot | A | \cdot | A_L | +  | \rderm A  |  \cdot | A |^2  + \big(  O \big( | h | \cdot | A |^2 \cdot | \pa  A | \big)   \\
 \notag
&& + \big(   | \derm_\mu  A  | +  | \rderm A   | \big) \cdot | A_\nu | \cdot | A | + \big(   | \derm_\nu  A  | +  | \rderm A   | \big) \cdot | A_\mu | \cdot | A | \\
\notag
&& +  | \rderm A   |  \cdot | A | \cdot  | A_L |  +  O \big( | h | \cdot  | A |^2 \cdot | \pa  A | \big) \\
\notag
&& +  | \derm A   |  \cdot | A | \cdot  | A_L |  +  | \rderm A   |  \cdot | A |^2  \\
 \notag
&& +  | A_\mu | \cdot | A_\nu | \cdot  | A |^2  + | A_L | \cdot  | A |^3 \\
     && + O \big(h \cdot  (\pa A)^2 \big)   + O \big(  h  \cdot  A^2 \cdot \pa A \big)     + O \big(  h   \cdot  A^4 \big)  \,  .
\eeaa
Hence,
  \beaa
\notag
 && | g^{\alpha\beta} \derm_\alpha \derm_\beta h_{\mu\nu} |   \\
       &\les&  | P(\pa_\mu h,\pa_\nu h) | + | Q_{\mu\nu}(\pa h,\pa h)  | + | G_{\mu\nu}(h)(\pa h,\pa h)  | \\
 \notag
    && +  | \rderm A   | \cdot  |\derm  A  | \\
    \notag
 &&+ \big(  |\derm_\mu  A_{e_{a}}  | +  |\rderm  A  | \big)  \cdot \big(  |\derm_\nu  A_{e_{a}}  | +  |\rderm  A  | \big) \\
   \notag
 && +  | \derm A   |  \cdot | A | \cdot | A_L | +  | \rderm A  |  \cdot | A |^2    \\
 \notag
&& + \big(   | \derm_\mu  A  | +  | \rderm A   | \big) \cdot | A_\nu | \cdot | A | + \big(   | \derm_\nu  A  | +  | \rderm A   | \big) \cdot | A_\mu | \cdot | A | \\
\notag
&& +  | \rderm A   |  \cdot | A | \cdot  | A_L |  \\
\notag
&& +  | \derm A   |  \cdot | A | \cdot  | A_L |  +  | \rderm A   |  \cdot | A |^2  \\
 \notag
&& +  | A_\mu | \cdot | A_\nu | \cdot  | A |^2  + | A_L | \cdot  | A |^3 \\
     && + O \big(h \cdot  (\pa A)^2 \big)   + O \big(  h  \cdot  A^2 \cdot \pa A \big)     + O \big(  h   \cdot  A^4 \big)  \,  \\
        &\les&  | P(\pa_\mu h,\pa_\nu h) | + | Q_{\mu\nu}(\pa h,\pa h)  | + | G_{\mu\nu}(h)(\pa h,\pa h)  | \\
 \notag
    && +  | \rderm A   | \cdot  |\derm  A  | \\
    \notag
 &&+ \big(  |\derm_\mu  A_{e_{a}}  | +  |\rderm  A  | \big)  \cdot \big(  |\derm_\nu  A_{e_{a}}  | +  |\rderm  A  | \big) \\
   \notag
 && +  | \derm A   |  \cdot | A | \cdot | A_L | +  | \rderm A  |  \cdot | A |^2    \\
 \notag
&& + \big(   | \derm_\mu  A  | +  | \rderm A   | \big) \cdot | A_\nu | \cdot | A | + \big(   | \derm_\nu  A  | +  | \rderm A   | \big) \cdot | A_\mu | \cdot | A | \\
 \notag
&& +  | A_\mu | \cdot | A_\nu | \cdot  | A |^2  + | A_L | \cdot  | A |^3 \\
     && + O \big(h \cdot  (\pa A)^2 \big)   + O \big(  h  \cdot  A^2 \cdot \pa A \big)     + O \big(  h   \cdot  A^4 \big)  \,  .
\eeaa

\end{proof}

\subsection{Coupled “bad” and “good” wave equations}

\subsubsection{Estimate on “good” components }\

\begin{lemma}\label{estimateonthesourcetermsforgoodcomponentofPoentialAandgoodcompometrich}
In the Lorenz and harmonic gauges, we have the following estimate for the tangential components of the Einstein-Yang-Mills potential,
  \beaa
   \notag
 && | g^{\la\mu} \derm_{\la}   \derm_{\mu}   A_{{\cal T}}   |  \\
 \notag
    &\les& | \derm h | \cdot  |\rderm A |       + | \rderm  h | \cdot  |\derm A |  \\
       \notag
           &&+   | A  | \cdot    | \rderm A  | +  | \derm  h | \cdot  | A  |^2   +  | A  |^3 \\
    \notag
  && + O( h \cdot  \derm h \cdot  \derm A) + O( h \cdot  A \cdot \derm A)  + O( h \cdot  \derm h \cdot  A^2) + O( h \cdot  A^3) \, .
  \eeaa

In the Lorenz and harmonic gauges, we have the following estimate on the components of the metric,
 \bea
\notag
 && | g^{\alpha\beta} \derm_\alpha \derm_\beta h_{ {\cal T} {\cal U}} |   \\
 \notag
    &\les&  | P(\derm_{\cal T} h,\derm_{\cal U} h) | + | Q_{{\cal T} {\cal U}}(\derm h,\derm h)  | + | G_{{\cal T} {\cal U} }(h)(\derm h,\derm h)  | \\
  \notag
    && +  | \rderm A   | \cdot  |\derm  A  |  +    | A |^2 \cdot | \derm A   |  + | A |^4  \\
\notag
     && + O \big(h \cdot  (\derm A)^2 \big)   + O \big(  h  \cdot  A^2 \cdot \derm A \big)     + O \big(  h   \cdot  A^4 \big)  . 
\eea
\end{lemma}

\begin{proof}
Taking $\si \in {\cal T} $ and using the estimate on the Einstein-Yang-Mills potential in the Lorenz and wave gauges, we get

 \beaa
   \notag
 | g^{\la\mu} \derm_{\la}   \derm_{\mu}   A_{{\cal T}}   |  
  &\les& | \derm h | \cdot  |\rderm A |       + | \rderm  h | \cdot  |\derm A |  \\
       \notag
           && +  | \rderm  h | \cdot  | A  |^2    \\
           \notag
 &&          +| \derm h_{L{\cal T}}   |  \cdot  | A  | \cdot  | A_L |   \\
        \notag
                  && + | \derm h_{e_{a} {\cal T}} |  \cdot    | A_{e_{a}} |  \cdot  | A_L |   \\
          \notag
 &&   +   | A_L  | \cdot    | \derm_{\cal T} A_{\underline{L}}   |      +    |  A  | \cdot  | A_{\cal T}  | \cdot  | A_L |   \\
    \notag
    &&   +  | A  | \cdot   | \rderm A  |     \\
       \notag
 &&  + | A_{e_{a}}  | \cdot    | \derm_{\cal T}A_{e_{a}} |    +    |  A_{e_{a}}  |^2 \cdot  | A_{\cal T}  |   \\
    \notag
  && + O( h \cdot  \pa h \cdot  \pa A) + O( h \cdot  \pa h \cdot  A^2) + O( h \cdot  A \cdot \pa A) + O( h \cdot  A^3) \, .\\
  \eeaa
Using again the wave coordinates estimate \ref{waveconditionestimateonzeroLiederivativeofmetric}, we obtain
 \beaa
   \notag
 | g^{\la\mu} \derm_{\la}   \derm_{\mu}   A_{{\cal T}}   |  
  &\les& | \derm h | \cdot  |\rderm A |       + | \rderm  h | \cdot  |\derm A |  \\
       \notag
           && +  | \rderm  h | \cdot  | A  |^2    \\
           \notag
 &&          + | \rderm  h |  \cdot  | A  | \cdot  | A_L |  + O( h \cdot  \pa h \cdot  A^2)  \\
        \notag
                  && + | \derm h_{e_{a} {\cal T}} |  \cdot    | A_{e_{a}} |  \cdot  | A_L |   \\
          \notag
 &&   +   | A_L  | \cdot    | \rderm A_{\underline{L}}   |      +    |  A  | \cdot  | A_{\cal T}  | \cdot  | A_L |   \\
    \notag
    &&   +  | A  | \cdot   | \rderm A  |     \\
       \notag
 &&  + | A_{e_{a}}  | \cdot    | \rderm A_{e_{a}}  |    +    |  A_{e_{a}}  |^2 \cdot  | A_{\cal T}  |   \\
    \notag
  && + O( h \cdot  \pa h \cdot  \pa A) + O( h \cdot  \pa h \cdot  A^2) + O( h \cdot  A \cdot \pa A) + O( h \cdot  A^3)  \\
   &\les& | \derm h | \cdot  |\rderm A |       + | \rderm  h | \cdot  |\derm A |  \\
       \notag
           && +  | \rderm  h | \cdot  | A  |^2    \\
        \notag
                  && + | \derm h_{e_{a} {\cal T}} |  \cdot    | A_{e_{a}} |  \cdot  | A_L |   \\
          \notag
 &&   +   | A  | \cdot    | \rderm A  |      +    |  A  | \cdot  | A_{\cal T}  |^2  \\
    \notag
  && + O( h \cdot  \pa h \cdot  \pa A) + O( h \cdot  \pa h \cdot  A^2) + O( h \cdot  A \cdot \pa A) + O( h \cdot  A^3) \, .\\
  \eeaa
  Hence,
   \beaa
   \notag
 | g^{\la\mu} \derm_{\la}   \derm_{\mu}   A_{{\cal T}}   |  
    &\les& | \derm h | \cdot  |\rderm A |       + | \rderm  h | \cdot  |\derm A |  \\
       \notag
           && +  | \derm  h | \cdot  | A  |^2  +   | A  |^3 \\
          \notag
 &&   +   | A  | \cdot    | \rderm A  |      \\
    \notag
  && + O( h \cdot  \pa h \cdot  \pa A) + O( h \cdot  \pa h \cdot  A^2) + O( h \cdot  A \cdot \pa A) + O( h \cdot  A^3) \, .\\
  \eeaa
  
Concerning the wave equation on the metric in the Lorenz and wave coordinates, for any $\mu \in {\cal T},  \nu \in {\cal U}$, we have
    \beaa
\notag
 | g^{\alpha\beta} \derm_\alpha \derm_\beta h_{ {\cal T} {\cal U}} |     &\les&  | P(\pa_{\cal T} h,\pa_{\cal U} h) | + | Q_{{\cal T} {\cal U}}(\pa h,\pa h)  | + | G_{{\cal T} {\cal U} }(h)(\pa h,\pa h)  | \\
 \notag
    && +  | \rderm A   | \cdot  |\derm  A  | \\
    \notag
 &&+ \big(  |\derm_{\cal T}  A_{e_{a}}  | +  |\rderm  A  | \big)  \cdot \big(  |\derm_{\cal U}  A_{e_{a}}  | +  |\rderm  A  | \big) \\
   \notag
 && +  | \derm A   |  \cdot | A | \cdot | A_L | +  | \rderm A  |  \cdot | A |^2    \\
 \notag
&& + \big(   | \derm_{\cal T}  A  | +  | \rderm A   | \big) \cdot | A_{\cal U} | \cdot | A | + \big(   | \derm_{\cal U}  A  | +  | \rderm A   | \big) \cdot | A_{\cal T} | \cdot | A | \\
 \notag
&& +  | A_{\cal T} | \cdot | A_{\cal U} | \cdot  | A |^2  + | A_L | \cdot  | A |^3 \\
     && + O \big(h \cdot  (\pa A)^2 \big)   + O \big(  h  \cdot  A^2 \cdot \pa A \big)     + O \big(  h   \cdot  A^4 \big)  \\
      &\les&  | P(\pa_{\cal T} h,\pa_{\cal U} h) | + | Q_{{\cal T} {\cal U}}(\pa h,\pa h)  | + | G_{{\cal T} {\cal U} }(h)(\pa h,\pa h)  | \\
 \notag
    && +  | \rderm A   | \cdot  |\derm  A  | \\
    \notag
 &&+   |\rderm  A  |  \cdot   |\derm  A  |  \\
   \notag
 && +  | \derm A   |  \cdot | A | \cdot | A_L | +  | \rderm A  |  \cdot | A |^2    \\
 \notag
&& +    | \rderm A   | \cdot | A |^2 + | \derm  A  |  \cdot | A_{\cal T} | \cdot | A | \\
 \notag
&& +  | A_{\cal T} |  \cdot  | A |^3  + | A_L | \cdot  | A |^3 \\
     && + O \big(h \cdot  (\pa A)^2 \big)   + O \big(  h  \cdot  A^2 \cdot \pa A \big)     + O \big(  h   \cdot  A^4 \big)  \\
           &\les&  | P(\pa_{\cal T} h,\pa_{\cal U} h) | + | Q_{{\cal T} {\cal U}}(\pa h,\pa h)  | + | G_{{\cal T} {\cal U} }(h)(\pa h,\pa h)  | \\
 \notag
    && +  | \rderm A   | \cdot  |\derm  A  | \\
   \notag
 && +  | \derm A   |    \cdot | A_{\cal T} | \cdot | A |   \\
 \notag
&& +    | \rderm A   | \cdot | A |^2 \\
 \notag
&& +  | A_{\cal T} |  \cdot  | A |^3  \\
     && + O \big(h \cdot  (\pa A)^2 \big)   + O \big(  h  \cdot  A^2 \cdot \pa A \big)     + O \big(  h   \cdot  A^4 \big)  .
\eeaa
Hence,
   \beaa
\notag
 | g^{\alpha\beta} \derm_\alpha \derm_\beta h_{ {\cal T} {\cal U}} |      &\les&  | P(\pa_{\cal T} h,\pa_{\cal U} h) | + | Q_{{\cal T} {\cal U}}(\pa h,\pa h)  | + | G_{{\cal T} {\cal U} }(h)(\pa h,\pa h)  | \\
  \notag
    && +  | \rderm A   | \cdot  |\derm  A  | \\
 \notag
&& +    | \derm A   | \cdot | A |^2 \\
 \notag
&& + |A |^4  \\
     && + O \big(h \cdot  (\pa A)^2 \big)   + O \big(  h  \cdot  A^2 \cdot \pa A \big)     + O \big(  h   \cdot  A^4 \big)  .
\eeaa

\end{proof}

\subsubsection{Estimate on “bad” components}

\begin{lemma}\label{structureofthesourcetermsofthewaveoperatoronAandh}
In the Lorenz gauge and in wave coordinates, we have
\beaa
   \notag
 && | g^{\la\mu} \derm_{\la}   \derm_{\mu}   A_{{\underline{L}}}  |   \\
 \notag
    &\les& | \derm h | \cdot  |\rderm A |       + | \rderm  h | \cdot  |\derm A |  \\
       \notag
           &&+   | A  | \cdot    | \rderm A  | +  | \derm  h | \cdot  | A  |^2   +  | A  |^3 \\
       \notag
 && + | A_L  | \cdot    | \derm A   |   + | A_{e_a}  | \cdot    | \derm A_{e_a}  |     \\
    \notag
  && + O( h \cdot  \derm h \cdot  \derm A) + O( h \cdot  A \cdot \derm A)  + O( h \cdot  \derm h \cdot  A^2) + O( h \cdot  A^3) \, .
  \eeaa

In the Lorenz gauge and in wave coordinates, we have
    \bea
\notag
&&  | g^{\alpha\beta} \derm_\alpha \derm_\beta h_{{\underline{L}} {\underline{L}} } |    \\
\notag
            &\les&  | P(\derm_{\underline{L}} h,\derm_{\underline{L}} h) | + | Q_{{\underline{L}}{\underline{L}}}(\derm h,\derm h)  | + | G_{{\underline{L}}{\underline{L}}}(h)(\derm h,\derm h)  | \\
 \notag
    && +  | \rderm A   | \cdot  |\derm  A  | +  | A |^2  \cdot | \derm A   |   +  |A |^4 \\
    \notag
     && + O \big(h \cdot  (\derm A)^2 \big)   + O \big(  h  \cdot  A^2 \cdot \derm A \big)     + O \big(  h   \cdot  A^4 \big)   \\
 &&+   |\derm A_{e_{a}}  |^2 \; .
\eea

\end{lemma}

\begin{proof}
Using the estimate on the Einstein-Yang-Mills potential in the Lorenz gauge and in wave coordinates, \ref{estimatesonsourcetermsforEinstein-Yang-Millssystem}, and taking $\si = {\underline{L}} $, we get
 \beaa
   \notag
 | g^{\la\mu} \derm_{\la}   \derm_{\mu}   A_{{\underline{L}}}  |  
  &\les& | \derm h | \cdot  |\rderm A |       + | \rderm  h | \cdot  |\derm A |  \\
       \notag
           && +  | \rderm  h | \cdot  | A  |^2    \\
           \notag
 &&          +| \derm h_{L {\underline{L}}}   |  \cdot  | A  | \cdot  | A_L |   \\
        \notag
                  && + | \derm h_{e_{a} {\underline{L}}} |  \cdot    | A_{e_{a}} |  \cdot  | A_L |   \\
          \notag
 &&   +   | A_L  | \cdot    | \derm_{\underline{L}} A_{\underline{L}}   |      +    |  A  | \cdot  | A_{\underline{L}}  | \cdot  | A_L |   \\
    \notag
    &&   +  | A  | \cdot   | \rderm A  |     \\
       \notag
 &&  + | A_{e_{a}}  | \cdot    | \derm_{\underline{L}} A_{e_{a}} |    +    |  A_{e_{a}}  |^2 \cdot  | A_{\underline{L}}  |   \\
    \notag
  && + O( h \cdot  \pa h \cdot  \pa A) + O( h \cdot  \pa h \cdot  A^2) + O( h \cdot  A \cdot \pa A) + O( h \cdot  A^3)  \\
    &\les& | \derm h | \cdot  |\rderm A |       + | \rderm  h | \cdot  |\derm A |  \\
       \notag
           && +  | \rderm  h | \cdot  | A  |^2    \\
           \notag
 &&          +| \derm h   |  \cdot  | A  | \cdot  | A_L |   \\
          \notag
 &&   +   | A_L  | \cdot    | \derm A_{\underline{L}}   |      \\
    \notag
    &&   +  | A  | \cdot   | \rderm A  |     \\
       \notag
 &&  + | A_{e_{a}}  | \cdot    | \derm A_{e_{a}} |    +    |  A  |^2 \cdot  | A_{\underline{L}}  |   \\
    \notag
  && + O( h \cdot  \pa h \cdot  \pa A) + O( h \cdot  \pa h \cdot  A^2) + O( h \cdot  A \cdot \pa A) + O( h \cdot  A^3) \, .\\
  \eeaa
Thus,
\beaa
   \notag
 | g^{\la\mu} \derm_{\la}   \derm_{\mu}   A_{{\underline{L}}}  |     &\les& | \derm h | \cdot  |\rderm A |       + | \rderm  h | \cdot  |\derm A |  \\
       \notag
           && +  | \derm  h | \cdot  | A  |^2    \\
    \notag
    &&   +  | A  | \cdot   | \rderm A  |   + | A_L  | \cdot    | \derm A   |   \\
       \notag
 &&  + | A_{e_a}  | \cdot    | \derm A_{e_a}  |    +    |  A  |^3   \\
    \notag
  && + O( h \cdot  \pa h \cdot  \pa A) + O( h \cdot  \pa h \cdot  A^2) + O( h \cdot  A \cdot \pa A) + O( h \cdot  A^3) \, .\\
  \eeaa
  
 Regarding the wave equation for the metric in the Lorenz gauge and in wave coordinates, taking $\mu = \nu = {\underline{L}} $ in \ref{estimatesonsourcetermsforEinstein-Yang-Millssystem}, we get
    \beaa
\notag
 | g^{\alpha\beta} \derm_\alpha \derm_\beta h_{{\underline{L}} {\underline{L}} } |     &\les&  | P(\pa_{\underline{L}} h,\pa_{\underline{L}} h) | + | Q_{{\underline{L}}{\underline{L}}}(\pa h,\pa h)  | + | G_{{\underline{L}}{\underline{L}}}(h)(\pa h,\pa h)  | \\
 \notag
    && +  | \rderm A   | \cdot  |\derm  A  | \\
    \notag
 &&+ \big(  |\derm_{\underline{L}}  A_{e_{a}} | +  |\rderm  A  | \big)  \cdot \big(  |\derm_{\underline{L}}  A_{e_{a}}  | +  |\rderm  A  | \big) \\
   \notag
 && +  | \derm A   |  \cdot | A | \cdot | A_L | +  | \rderm A  |  \cdot | A |^2    \\
 \notag
&& + \big(   | \derm_{\underline{L}}  A  | +  | \rderm A   | \big) \cdot | A_{\underline{L}} | \cdot | A | + \big(   | \derm_{\underline{L}}  A  | +  | \rderm A   | \big) \cdot | A_{\underline{L}} | \cdot | A | \\
 \notag
&& +  | A_{\underline{L}} | \cdot | A_{\underline{L}} | \cdot  | A |^2  + | A_L | \cdot  | A |^3 \\
     && + O \big(h \cdot  (\pa A)^2 \big)   + O \big(  h  \cdot  A^2 \cdot \pa A \big)     + O \big(  h   \cdot  A^4 \big)  \,  \\
        &\les&  | P(\pa_{\underline{L}} h,\pa_{\underline{L}} h) | + | Q_{{\underline{L}}{\underline{L}}}(\pa h,\pa h)  | + | G_{{\underline{L}}{\underline{L}}}(h)(\pa h,\pa h)  | \\
 \notag
    && +  | \rderm A   | \cdot  |\derm  A  | \\
    \notag
 &&+   |\derm A_{e_{a}}  |^2 \\
   \notag
 && +  | \derm A   |  \cdot | A | \cdot | A_L | +  | \rderm A  |  \cdot | A |^2    \\
 \notag
&& +  | \derm A   | \cdot | A_{\underline{L}} | \cdot | A | +  | \derm A   | \cdot | A_{\underline{L}} | \cdot | A | \\
 \notag
&& +  | A_{\underline{L}} | \cdot | A_{\underline{L}} | \cdot  | A |^2  + | A_L | \cdot  | A |^3 \\
     && + O \big(h \cdot  (\pa A)^2 \big)   + O \big(  h  \cdot  A^2 \cdot \pa A \big)     + O \big(  h   \cdot  A^4 \big)  \,  \\
           &\les&  | P(\pa_{\underline{L}} h,\pa_{\underline{L}} h) | + | Q_{{\underline{L}}{\underline{L}}}(\pa h,\pa h)  | + | G_{{\underline{L}}{\underline{L}}}(h)(\pa h,\pa h)  | \\
 \notag
    && +  | \rderm A   | \cdot  |\derm  A  | \\
    \notag
 &&+   |\derm A_{e_{a}}  |^2 \\
   \notag
 && +  | \derm A   |  \cdot | A |^2  +  | \rderm A  |  \cdot | A |^2    \\
 \notag
&& +  | A |^4 \\
     && + O \big(h \cdot  (\pa A)^2 \big)   + O \big(  h  \cdot  A^2 \cdot \pa A \big)     + O \big(  h   \cdot  A^4 \big)  \,  .
\eeaa

\end{proof}

\section{Studying the structure of the Lie derivatives of the source terms in both the Lorenz and harmonic gauges}

First, we remind the following lemma based on what we already showed in \cite{G4}.
  \begin{lemma}\label{LiederivativeZcommuteswithminkowksicovariantderivative}
The Minkowski covariant derivative commutes with the Lie derivative along Minkowski vector fields, that is for any tensor $K$, 
\bea
 \Lie_{Z^I}  \derm  K  =   \derm ( \Lie_{Z^I}   K ) .
\eea
Since  $\derm  K$ is also a tensor, it follows that $\Lie_{Z^I} $ commutes with any product of $\derm$.

Note that the Lie derivatives are not being differentiated in $\derm ( \Lie_{Z^I}   K) $; the differentiation concerns only the tensor $K$.

As a result, we also have 
\bea
 \Lie_{Z^I}  \rderm  K  =   \rderm ( \Lie_{Z^I}   K ) .
\eea

\end{lemma}

\begin{proof}
We already proved in \cite{G4} that
 \bea
 \Lie_{Z^I}  \derm_\a  K_\b  =   \derm_\a ( \Lie_{Z^I}   K_\b ) .
\eea
By taking $\a \in \cal T$, we also get 
 \bea
 \Lie_{Z^I}  \rderm  K_\b  =   \rderm ( \Lie_{Z^I}   K_\b ) .
\eea

The proof works for any tensor $K$ of arbitrary order. 

\end{proof}

\begin{lemma}\label{structureofsourcetermsofLiederivaivesZofthehyperbolicsystemonAandh}

In the Lorenz gauge the Yang-Mills potential satisfies the following equations decomposed in a null-frame, i.e. for any $\si \in \cal U$, we have
     \bea
   \notag
&& \Lie_{Z^I}  g^{\la\mu} \derm_{\la}   \derm_{\mu}   A_{\si}   \\ 
  \notag
 &=&    \sum_{I_1+I_2+I_3+I_4  = I}  \hat{c} (I_1) \cdot   \hat{c} (I_2) \cdot  \Big( \frac{1}{4} (  \derm_{\si}  ( \Lie_{Z^{I_3}}   h_{\underline{L}\underline{L}} ) ) \cdot  \derm_{L} (  \Lie_{Z^{I_4}}  A_{L} )  \\
 \notag
 && + \frac{1}{4} (   \derm_{\si}   (  \Lie_{Z^{I_3}} h_{\underline{L}L} ) )  \cdot \derm_{L}  (  \Lie_{Z^{I_4}} A_{\underline{L}} )     -\frac{1}{2}   (  \derm_{\si}  (  \Lie_{Z^{I_3}} h_{\underline{L}e_{A}} ) ) \cdot \derm_{L}  (  \Lie_{Z^{I_4}} A_{e_{A}}  )  \\
   \notag
 && +  \frac{1}{4} (  \derm_{\si}  (  \Lie_{Z^{I_3}}  h_{L\underline{L}} ) ) \cdot \derm_{\underline{L} } (  \Lie_{Z^{I_4}} A_{L} )  +  \frac{1}{4}  (  \derm_{\si}   (  \Lie_{Z^{I_3}}  h_{LL} ) ) \cdot \derm_{\underline{L} } (  \Lie_{Z^{I_4}} A_{\underline{L}} ) \\
 \notag
 && -  \frac{1}{2} (  \derm_{\si}  (  \Lie_{Z^{I_3}}  h_{Le_{A}} ) ) \cdot  \derm_{\underline{L} } (  \Lie_{Z^{I_4}} A_{e_{A}} )     -  \frac{1}{2} (  \derm_{\si}  (  \Lie_{Z^{I_3}}  h_{e_{A}\underline{L}} ) )  \cdot \derm_{e_{A}} (  \Lie_{Z^{I_4}} A_{L}  )  \\
 \notag
 && -  \frac{1}{2} (   \derm_{\si}  (  \Lie_{Z^{I_3}}  h_{e_{A}L} ) ) \cdot  \derm_{e_{A}} (  \Lie_{Z^{I_4}} A_{\underline{L}} ) +  ( \derm_{\si}   (  \Lie_{Z^{I_3}}  h_{e_{A}e_{A}} ) )  \cdot \derm_{e_{A}} (  \Lie_{Z^{I_4}} A_{e_{A}} ) \Big) \\
   \notag
&& + \frac{1}{4}  \big(    \derm_L   (  \Lie_{Z^{I_3}}  h_{\underline{L}\si} ) -  \derm_{\underline{L}}   (  \Lie_{Z^{I_3}}  h_{L \si} )  \big)   \cdot  \big( \derm_{\underline{L}} (  \Lie_{Z^{I_4}}  A_{L} ) -  \derm_{L} (  \Lie_{Z^{I_4}}  A_{\underline{L}} )   \big) \\
   \notag
&& - \frac{1}{2}    \big(    \derm_L   (  \Lie_{Z^{I_3}}  h_{e_{A}\si} ) -  \derm_{e_{A}}   (  \Lie_{Z^{I_3}}  h_{L \si} )  \big)   \cdot   \big(  \derm_{\underline{L}} (  \Lie_{Z^{I_4}}  A_{e_{A}} ) -  \derm_{e_{A}} (  \Lie_{Z^{I_4}}  A_{\underline{L}} )  \big) \\
     \notag
   && - \frac{1}{2}   \big(    \derm_{\underline{L}} (  \Lie_{Z^{I_3}} h_{e_{A}\si} ) -  \derm_{e_{A}} (  \Lie_{Z^{I_3}} h_{\underline{L}\si} )  \big)   \cdot \big(  \derm_{L}(  \Lie_{Z^{I_4}}  A_{e_{A}} ) -  \derm_{e_{A}} (  \Lie_{Z^{I_4}}  A_{L} )  \big) \\
 \notag
&& +  \sum_{I_1+I_2+I_3+I_4 + I_5 = I} \hat{c} (I_1) \cdot   \hat{c} (I_2) \cdot  \Big( \frac{1}{4}   \big(    \derm_L  ( \Lie_{Z^{I_3}}  h_{\underline{L}\si} ) -  \derm_{\underline{L}}  ( \Lie_{Z^{I_3}}  h_{L \si} )   \big)  \cdot   [ \Lie_{Z^{I_4}}  A_{\underline{L}}, \Lie_{Z^{I_5}}  A_{L}  ]    \\
   \notag
&& - \frac{1}{2}   \big(    \derm_L  ( \Lie_{Z^{I_3}}  h_{e_{A}\si} ) -  \derm_{e_{A}}  ( \Lie_{Z^{I_3}}  h_{L \si} ) \big)   \cdot   [  \Lie_{Z^{I_4}}  A_{\underline{L}},  \Lie_{Z^{I_5}} A_{e_{A}}]   \\
    \notag
   && - \frac{1}{2}   \big(   \derm_{\underline{L}} ( \Lie_{Z^{I_3}} h_{e_{A}\si} ) - \derm_{e_{A}} ( \Lie_{Z^{I_3}} h_{\underline{L}\si} )  \big)   \cdot  [\Lie_{Z^{I_4}}  A_{L}, \Lie_{Z^{I_5}} A_{e_{A}}]    \Big)  \\
   \notag
 && +   \sum_{I_1+I_2+I_3  = I}    \hat{c} (I_1) \cdot \Big( \frac{1}{2}  [  \Lie_{Z^{I_3}} A_{\underline{L}}  , ( \Lie_{Z^{I_4}}  \derm_{L} A_{\si} )  ]  +   \frac{1}{2}  [ \Lie_{Z^{I_3}} A_{L},   \derm_{\underline{L}} ( \Lie_{Z^{I_4}} A_{\si} )  - \derm_{\si} ( \Lie_{Z^{I_4}} A_{\underline{L}} ) ]       \\
    \notag
    &&  +  \frac{1}{2}   [ \Lie_{Z^{I_3}} A_{L},  \derm_{\underline{L}} ( \Lie_{Z^{I_4}}A_{\si} ) ]  + \frac{1}{2}    [ \Lie_{Z^{I_3}} A_{\underline{L}},   \derm_{L} ( \Lie_{Z^{I_4}}  A_{\si} ) - \derm_{\si} ( \Lie_{Z^{I_4}} A_{L} ) ]    \\
       \notag
 &&   -     [\Lie_{Z^{I_3}} A_{e_{a}}, \derm_{e_{a}} ( \Lie_{Z^{I_4}} A_{\si} ) ] -  [\Lie_{Z^{I_3}} A_{e_{a}},  \derm_{e_{a}} ( \Lie_{Z^{I_4}} A_{\si}  ) -  \derm_{\si} ( \Lie_{Z^{I_4}} A_{e_{a}} ) ]     \big)    \\
 \notag
 &&+    \sum_{I_1+I_2+I_3+I_4  = I}  \hat{c} (I_1) \cdot \Big( \frac{1}{2}  [\Lie_{Z^{I_3}} A_{L}, [\Lie_{Z^{I_4}} A_{\underline{L}}, \Lie_{Z^{I_5}} A_{\si}] ]  +    \frac{1}{2}  [\Lie_{Z^{I_3}} A_{\underline{L}}, [\Lie_{Z^{I_4}} A_{L}, \Lie_{Z^{I_5}} A_{\si}] ] \\
 \notag
 &&  -  [\Lie_{Z^{I_3}} A_{e_{a}}, [\Lie_{Z^{I_4}} A_{e_{a}}, \Lie_{Z^{I_5}} A_{\si}] ] \Big)   \\
 \notag
  && +  \Lie_{Z^I}  O( h \cdot  \pa h \cdot  \pa A) +  \Lie_{Z^I}  O( h \cdot  \pa h \cdot  A^2) +  \Lie_{Z^I}  O( h \cdot  A \cdot \pa A) +  \Lie_{Z^I}  O( h \cdot  A^3) \, . \\
  \eea

The perturbations $h$ of the metric $m$, defined to be the Minkowski metric in wave coordinates, solutions to the Einstein-Yang-Mills equations satisfy the following wave equation, i.e. for any $\mu, \nu \in \cal U$, we have
    \bea
\notag
&&  \Lie_{Z^I} g^{\alpha\beta} \derm_\alpha \derm_\beta h_{\mu\nu}   \\
\notag
     &=&   \Lie_{Z^I} P(\pa_\mu h,\pa_\nu h) + \Lie_{Z^I}  Q_{\mu\nu}(\pa h,\pa h)   +  \Lie_{Z^I} G_{\mu\nu}(h)(\pa h,\pa h)  \\
\notag
 &&  + \sum_{I_1+I_2+I_3  = I}   \hat{c} (I_1)   \cdot  \Big( 2  <    \derm_{\mu} ( \Lie_{Z^{I_2}}  A_{L}  ) -  \derm_{L} ( \Lie_{Z^{I_2}}  A_{\mu} ) ,  \derm_{\nu} ( \Lie_{Z^{I_3}} A_{ \underline{L} } ) -  \derm_{ \underline{L} } ( \Lie_{Z^{I_3}} A_{\nu} ) >   \\
 \notag
 &&  +2 <   \derm_{\mu}( \Lie_{Z^{I_2}}  A_{ \underline{L} } )  - \derm_{ \underline{L} } ( \Lie_{Z^{I_2}}  A_{\mu} ) ,  \derm_{\nu} ( \Lie_{Z^{I_3}}  A_{L} ) - \derm_{L} ( \Lie_{Z^{I_3}} A_{\nu}  ) >    \\
 \notag
  && -4   <    \derm_{\mu} ( \Lie_{Z^{I_2}}  A_{e_A} ) -  \derm_{e_A } ( \Lie_{Z^{I_2}}  A_{\mu} ),   \derm_{\nu} ( \Lie_{Z^{I_3}} A_{e_A} ) - \derm_{e_A}( \Lie_{Z^{I_3}} A_{\nu} )  >\Big)  \\
  \notag
  &&   +   \sum_{I_1+I_2+I_3+I_4  + I_5= I}  c (I_1) \cdot   \hat{c} (I_2)  \cdot  \hat{c} (I_3) \\
  \notag
  &&  \cdot \Big( \frac{1}{2} m_{\mu\nu }    \cdot  <   \derm_{L} ( \Lie_{Z^{I_4}}  A_{\underline{L}} ) - \derm_{\underline{L}} ( \Lie_{Z^{I_4}}  A_{L} ) ,  \derm_{\underline{L}}  ( \Lie_{Z^{I_5}}  A_{L} ) - \derm_{L}  ( \Lie_{Z^{I_5}}  A_{\underline{L}} ) > \\
\notag
&& -  2 m_{\mu\nu }     \cdot  <   \derm_{e_A} ( \Lie_{Z^{I_4}}  A_{\underline{L}} ) -  \derm_{\underline{L}} ( \Lie_{Z^{I_4}}  A_{e_A}),  \derm_{e_A} ( \Lie_{Z^{I_5}}  A_{L})  - \derm_{L} ( \Lie_{Z^{I_5}}  A_{e_A} ) >  \Big) \\
\notag
 &&  + \sum_{I_1+I_2+I_3+I_4  = I}  \hat{c} (I_1) \\
 \notag
 && \cdot  \Big( 2  \big( <   \derm_{\mu}  ( \Lie_{Z^{I_2}}  A_{\underline{L}} ) - \derm_{\underline{L}} ( \Lie_{Z^{I_2}}   A_{\mu} ),  [ \Lie_{Z^{I_3}} A_{\nu}, \Lie_{Z^{I_4}}  A_{L}] > \\
 \notag
 &&  + <   [ \Lie_{Z^{I_3}} A_{\mu}, \Lie_{Z^{I_4}}  A_{\underline{L}}] ,   \derm_{\nu} ( \Lie_{Z^{I_2}}   A_{L} ) -  \derm_{L} ( \Lie_{Z^{I_2}}   A_{\nu} )  > \big) \\
 \notag
&&   + 2  \big( <   \derm_{\mu} ( \Lie_{Z^{I_2}}   A_{L}) - \derm_{L} ( \Lie_{Z^{I_2}}   A_{\mu} ) ,  [ \Lie_{Z^{I_3}}  A_{\nu}, \Lie_{Z^{I_4}} A_{\underline{L}}] >  \\
\notag
&& + <   [\Lie_{Z^{I_3}}  A_{\mu}, \Lie_{Z^{I_4}}  A_{L}] ,  \derm_{\nu} ( \Lie_{Z^{I_2}}   A_{\underline{L}})  - \derm_{\underline{L}} ( \Lie_{Z^{I_2}}  A_{\nu} )  > \big)  \\
\notag
 &&  -4     \big( <    \derm_{\mu} ( \Lie_{Z^{I_2}}   A_{e_A} )  - \derm_{e_A} ( \Lie_{Z^{I_2}}  A_{\mu} ) ,  [\Lie_{Z^{I_3}}  A_{\nu}, \Lie_{Z^{I_4}} A_{e_A}] >  \\
 \notag
 &&  -4  <   [\Lie_{Z^{I_3}}  A_{\mu}, \Lie_{Z^{I_4}}  A_{e_A}] ,   \derm_{\nu} ( \Lie_{Z^{I_2}}   A_{e_A} )  -  \derm_{e_A} ( \Lie_{Z^{I_2}}   A_{\nu} ) > \big) \Big)   \\
\notag
&&  +   \sum_{I_1+I_2+I_3+I_4 + I_5 + I_6 = I}   c (I_1) \cdot   \hat{c} (I_2)  \cdot  \hat{c} (I_3) \cdot  m_{\mu\nu }  \Big(  \cdot   <  \derm_{L} ( \Lie_{Z^{I_4}} A_{\underline{L}} )  \\
\notag
&& -  \derm_{\underline{L}} ( \Lie_{Z^{I_4}} A_{L} ) , [\Lie_{Z^{I_5}} A_{\underline{L}}, \Lie_{Z^{I_6}} A_{L}] >   \\
\notag
 &&  -  2 m_{\mu\nu }     \cdot \big(  <   \derm_{e_A }  ( \Lie_{Z^{I_4}} A_{\underline{L}} )  -  \derm_{\underline{L}} ( \Lie_{Z^{I_4}} A_{e_A} ) , [\Lie_{Z^{I_5}} A_{e_A}, \Lie_{Z^{I_6}} A_{L}] >    \\
 \notag
 &&+  <  [A_{e_A},A_{\underline{L}}] ,  \derm_{e_A}A_{L} -  \derm_{L}A_{e_A}  > \big)  \Big)  \\
 \notag
 && +  \sum_{I_1+I_2+I_3+I_4 + I_5  = I}   \hat{c} (I_1) \cdot \Big( 2   <   [  \Lie_{Z^{I_2}} A_{\mu},  \Lie_{Z^{I_3}} A_{\underline{L}}] ,  [  \Lie_{Z^{I_4}} A_{\nu},  \Lie_{Z^{I_5}} A_{L}] > \\
 \notag
 && + 2   <   [  \Lie_{Z^{I_2}} A_{\mu},  \Lie_{Z^{I_3}} A_{L}] ,  [  \Lie_{Z^{I_4}} A_{\nu},  \Lie_{Z^{I_5}} A_{\underline{L}}] >  -4    <   [ \Lie_{Z^{I_2}} A_{\mu},  \Lie_{Z^{I_3}} A_{e_A}] ,  [  \Lie_{Z^{I_4}} A_{\nu},  \Lie_{Z^{I_5}} A_{e_A}] >  \Big) \\
 \notag
  && +   \frac{1}{2}   \sum_{I_1+I_2+I_3+I_4 + I_5 + I_6 + I_7  = I}  c (I_1) \cdot   \hat{c} (I_2)  \cdot  \hat{c} (I_3) \\
  \notag
  && \cdot   \Big( m_{\mu\nu }    \cdot   <  [\Lie_{Z^{I_4}} A_{L}, \Lie_{Z^{I_5}}A_{\underline{L}}] , [\Lie_{Z^{I_6}} A_{\underline{L}}, \Lie_{Z^{I_7}}A_{L}] >   - 2 m_{\mu\nu }   \cdot   <  [\Lie_{Z^{I_4}} A_{e_A}, \Lie_{Z^{I_5}}A_{\underline{L}}] , [ \Lie_{Z^{I_6}} A_{e_A}, \Lie_{Z^{I_7}}A_{L}] > \\
  \notag
     && + \Lie_{Z^I}  O \big(h \cdot  (\derm A)^2 \big)   + \Lie_{Z^I}  O \big(  h  \cdot  A^2 \cdot \derm A \big)     + \Lie_{Z^I}  O \big(  h   \cdot  A^4 \big)  \,  . \\
\eea

    \end{lemma}

    \begin{proof}

We showed -- see Lemma \ref{Einstein-Yang-MillssysteminLorenzwavegauges} -- that in the Lorenz gauge, and for $m$ defined to be the Minkowski metric in wave coordinates, we have
            \beaa
   \notag
g^{\la\mu} \derm_{\la}   \derm_{\mu}   A_{\si}     &=&  m^{\a\ga} \cdot m ^{\mu\la} \cdot ( \derm_{\si}  h_{\ga\la} )  \cdot  \derm_{\a}A_{\mu}     \\
&&  +   \frac{1}{2}  m^{\a\mu} \cdot m^{\b\nu} \cdot  \big(   \derm_\a h_{\b\si} +  \derm_\si h_{\b\a}-  \derm_\b h_{\a\si}  \big)   \cdot  \big( \derm_{\mu}A_{\nu} -  \derm_{\nu}A_{\mu}  \big) \\
 \notag
&& +      \frac{1}{2}  m^{\a\mu} \cdot m^{\b\nu}  \cdot \big(   \derm_\a h_{\b\si} +  \derm_\si h_{\b\a}-  \derm_\b h_{\a\si}  \big)   \cdot   [A_{\mu},A_{\nu}] \\
 \notag
 && -  m^{\a\mu} \cdot \big(  [ A_{\mu}, \derm_{\a} A_{\si} ]  +    [A_{\alpha},  \derm_{\mu}  A_{\si} - \derm_{\si} A_{\mu} ]  \big)  +  m^{\a\mu} \cdot  [A_{\alpha}, [A_{\mu},A_{\si}] ]  \big)  \\
  && + O( h \cdot  \derm h \cdot  \derm A) + O( h \cdot  \derm h \cdot  A^2) + O( h \cdot  A \cdot \derm A) + O( h \cdot  A^3) \, .
  \eeaa
 
 Taking the Lie derivatives of the left hand side, we get
 
             \beaa
   \notag
  && \Lie_{Z^I} g^{\la\mu} \derm_{\la}   \derm_{\mu}   A_{\si}  \\
  \notag
     &=&   \Lie_{Z^{I_1}} ( m^{\a\ga} \cdot  m ^{\mu\la} ) \cdot   \Lie_{Z^{I_2}} (   \derm_{\si}  h_{\ga\la}   \cdot  \derm_{\a}A_{\mu}  )   \\
     \notag
&&  +   \frac{1}{2}  \Lie_{Z^{I_1}}   ( m^{\a\mu}  \cdot  m^{\b\nu} )  \cdot   \Lie_{Z^{I_2}} \Big( \big(   \derm_\a h_{\b\si} +    \derm_\si h_{\b\a}-   \derm_\b h_{\a\si}  \big)   \cdot  \big(   \derm_{\mu} A_{\nu} -    \derm_{\nu} A_{\mu}  \big) \Big) \\
 \notag
&& +      \frac{1}{2}  \Lie_{Z^{I_1}}    ( m^{\a\mu}   \cdot    m^{\b\nu} ) \cdot  \Lie_{Z^{I_2}}  \Big(  \big(      \derm_\a h_{\b\si} +  \derm_\si h_{\b\a}-     \derm_\b h_{\a\si}  \big)   \cdot   [ A_{\mu},  A_{\nu}]   \Big) \\
 \notag
 && -  (  \Lie_{Z^{I_1}}   m^{\a\mu} ) \cdot  \Lie_{Z^{I_2}}  \Big( \big(  [  A_{\mu},  \derm_{\a} A_{\si} ]  +    [ A_{\alpha}, \derm_{\mu}  A_{\si} -  \derm_{\si} A_{\mu} ]  \Big)  \\
 \notag
 &&   +   (  \Lie_{Z^{I_1}}   m^{\a\mu} ) \cdot  \Lie_{Z^{I_2}}  [ A_{\alpha}, [  A_{\mu}, A_{\si}] ]   \\
  && +  \Lie_{Z^I}  O( h \cdot  \derm h \cdot  \derm A) +  \Lie_{Z^I}  O( h \cdot  \derm h \cdot  A^2) +  \Lie_{Z^I}  O( h \cdot  A \cdot \derm A) + \Lie_{Z^I}  O( h \cdot  A^3) \, .
  \eeaa

To lighten the notation, we proceed again on writing the sum here $\sum_{I_1+I_2+I_3+I_4 + \ldots = I}$ on orderings such that the following equalities hold true.

Using the fact that for two tensors $\Phi$ and $\Psi$ valued in the Lie algebra, we have
\beaa
\Lie_{Z^I} [\Phi, \Psi] =\Lie_{Z^I}  ( \Phi \cdot \Psi  -   \Psi \cdot \Phi  ) &=&  \sum_{I_1+I_2 = I}  (  \Lie_{Z^{I_1}} \Phi ) \cdot (  \Lie_{Z^{I_2}} \Psi )   -  ( \Lie_{Z^{I_1}} \Psi ) \cdot (  \Lie_{Z^{I_2}} \Phi )  \\
&=&  \sum_{I_1+I_2 = I} [  \Lie_{Z^{I_1}} \Phi, \Lie_{Z^{I_2}} \Psi] \; ,
\eeaa
we obtain
             \beaa
   \notag
  && \Lie_{Z^I} g^{\la\mu} \derm_{\la}   \derm_{\mu}   A_{\si}  \\
  \notag
     &=&   \sum_{I_1+I_2+I_3+I_4  = I} ( \Lie_{Z^{I_1}} m^{\a\ga} ) \cdot  ( \Lie_{Z^{I_2}} m ^{\mu\la} ) \cdot  (   \Lie_{Z^{I_3}} \derm_{\si}  h_{\ga\la} )  \cdot  (  \Lie_{Z^{I_4}} \derm_{\a}A_{\mu}  )   \\
     \notag
&&  +   \sum_{I_1+I_2+I_3+I_4  = I}  \frac{1}{2}   ( \Lie_{Z^{I_1}} m^{\a\mu} ) \cdot (   \Lie_{Z^{I_2}} m^{\b\nu} ) \cdot   \big(    \Lie_{Z^{I_3}}  \derm_\a h_{\b\si} +   \Lie_{Z^{I_3}}  \derm_\si h_{\b\a}-    \Lie_{Z^{I_3}} \derm_\b h_{\a\si}  \big)   \\
\notag
&& \cdot  \big(   \Lie_{Z^{I_4}} \derm_{\mu}A_{\nu} -    \Lie_{Z^{I_4}} \derm_{\nu}A_{\mu}  \big) \\
 \notag
&& +   \sum_{I_1+I_2+I_3+I_4 + I_5 = I}   \frac{1}{2}    ( \Lie_{Z^{I_1}} m^{\a\mu}   ) \cdot   ( \Lie_{Z^{I_2}} m^{\b\nu} ) \cdot   \big(     \Lie_{Z^{I_3}} \derm_\a h_{\b\si} +  \Lie_{Z^{I_3}}  \derm_\si h_{\b\a}- \Lie_{Z^{I_3}}    \derm_\b h_{\a\si}  \big) \\
\notag
&&   \cdot   [\Lie_{Z^{I_4}} A_{\mu}, \Lie_{Z^{I_5}} A_{\nu}] \\
 \notag
 && -  \sum_{I_1+I_2+I_3  = I}   ( \Lie_{Z^{I_1}} m^{\a\mu} ) \cdot  \big(  [ \Lie_{Z^{I_2}} A_{\mu},  \Lie_{Z^{I_3}} \derm_{\a} A_{\si} ]  +    [ \Lie_{Z^{I_2}} A_{\alpha}, \Lie_{Z^{I_3}} \derm_{\mu}  A_{\si} - \Lie_{Z^{I_3}} \derm_{\si} A_{\mu} ] \big)   \\
 \notag
 &&   +    \sum_{I_1+I_2+I_3+I_4  = I} ( \Lie_{Z^{I_1}} m^{\a\mu} ) \cdot  [\Lie_{Z^{I_2}} A_{\alpha}, [ \Lie_{Z^{I_3}} A_{\mu}, \Lie_{Z^{I_4}}A_{\si}] ]  \\
  && + \Lie_{Z^I}  O( h \cdot  \derm h \cdot  \derm A) + \Lie_{Z^I}  O( h \cdot  \derm h \cdot  A^2) + \Lie_{Z^I}  O( h \cdot  A \cdot \derm A) + \Lie_{Z^I}  O( h \cdot  A^3) \, .
  \eeaa

 Now, using the fact that -- as we had shown in Lemma \ref{LiederivativeZcommuteswithminkowksicovariantderivative} -- the Lie derivative $\Lie_{Z^{I_1}}$ commutes with the Minkowski covariant derivative $\derm_\a$, and the fact that -- see Lemma \ref{LiederivativeZofMinkwoskimetric} --, $\Lie_{Z^{I_1}} m^{\a\b} = \hat{c} (I_1) \cdot m^{\a\b} $, we get
              \beaa
   \notag
  && \Lie_{Z^I} g^{\la\mu} \derm_{\la}   \derm_{\mu}   A_{\si}  \\
  \notag
     &=&   \sum_{I_1+I_2+I_3+I_4  = I}  \hat{c} (I_1) \cdot   \hat{c} (I_2) \cdot m^{\a\ga} \cdot m ^{\mu\la}  \cdot  (  \derm_{\si}  ( \Lie_{Z^{I_3}}   h_{\ga\la} )  ) \cdot  ( \derm_{\a} ( \Lie_{Z^{I_4}}  A_{\mu}  )  ) \\
     \notag
&&  +   \sum_{I_1+I_2+I_3+I_4  = I}  \frac{1}{2}   \hat{c} (I_1) \cdot   \hat{c} (I_2)  \cdot m^{\a\mu}   \cdot  m^{\b\nu}  \cdot   \big(    \derm_\a ( \Lie_{Z^{I_3}}   h_{\b\si} )  +    \derm_\si (  \Lie_{Z^{I_3}} h_{\b\a} ) -    \derm_\b ( \Lie_{Z^{I_3}}   h_{\a\si}  ) \big)   \\
\notag
&& \cdot  \big(   \derm_{\mu} (  \Lie_{Z^{I_4}} A_{\nu} ) -    \derm_{\nu}  \Lie_{Z^{I_4}} A_{\mu} )   \big) \\
 \notag
&& +   \sum_{I_1+I_2+I_3+I_4 + I_5 = I}   \frac{1}{2}     \hat{c} (I_1)  \cdot   \hat{c} (I_2)  \cdot  m^{\a\mu}   \cdot  m^{\b\nu}  \cdot   \big(    \derm_\a ( \Lie_{Z^{I_3}}  h_{\b\si} )  +   \derm_\si  ( \Lie_{Z^{I_3}}  h_{\b\a}) -   \derm_\b ( \Lie_{Z^{I_3}}  h_{\a\si} )  \big) \\
\notag
&&   \cdot   [\Lie_{Z^{I_4}} A_{\mu}, \Lie_{Z^{I_5}} A_{\nu}] \\
 \notag
 && -  \sum_{I_1+I_2+I_3  = I}    \hat{c} (I_1) \cdot  m^{\a\mu}  \cdot  \big(  [ \Lie_{Z^{I_2}} A_{\mu},  \derm_{\a}  (  \Lie_{Z^{I_3}} A_{\si} )  ]  +    [ \Lie_{Z^{I_2}} A_{\alpha}, \derm_{\mu} (\Lie_{Z^{I_3}}  A_{\si})  -  \derm_{\si} ( \Lie_{Z^{I_3}}  A_{\mu} ) ] \big)   \\
 \notag
 &&   +    \sum_{I_1+I_2+I_3+I_4  = I}  \hat{c} (I_1) \cdot  m^{\a\mu}  \cdot  [\Lie_{Z^{I_2}} A_{\alpha}, [ \Lie_{Z^{I_3}} A_{\mu}, \Lie_{Z^{I_4}}A_{\si}] ]  \\
  && + \Lie_{Z^I}  O( h \cdot  \derm h \cdot  \derm A) + \Lie_{Z^I}  O( h \cdot  \derm h \cdot  A^2) + \Lie_{Z^I}  O( h \cdot  A \cdot \derm A) + \Lie_{Z^I}  O( h \cdot  A^3) \, .
  \eeaa
  
 We had already studied the decomposition of such a structure in the null frame ${ \cal U}$, for the zeroth-Lie derivative. However, the structure here for the product of the $I$-Lie derivatives, as we have just shown, is the same. Hence we get the stated result. 

We proceed similarly for $h$: we had already showed (see \eqref{Einstein-Yang-MillssysteminLorenzwavegauges}) that
 
   \beaa
\notag
 && g^{\alpha\beta}\derm_\alpha \derm_\beta h_{\mu\nu} \\
  &=& P(\derm_\mu h,\derm_\nu h)  +  Q_{\mu\nu}(\derm h,\derm h)   + G_{\mu\nu}(h)(\derm h,\derm h)  \\
\notag
 &&   -4   m^{\si\b} \cdot  <   \derm_{\mu}A_{\b} - \derm_{\b}A_{\mu}  ,  \derm_{\nu}A_{\si} - \derm_{\si}A_{\nu}  >    \\
 \notag
 &&   + m_{\mu\nu }  \cdot m^{\si\b}  \cdot  m^{\a\la}    \cdot  <  \derm_{\a}A_{\b} - \derm_{\b}A_{\a} , \derm_{\la}A_{\si} - \derm_{\si}A_{\la} >   \\
 \notag
&&           -4 m^{\si\b}  \cdot  \big( <   \derm_{\mu}A_{\b} - \derm_{\b}A_{\mu}  ,  [A_{\nu},A_{\si}] >   + <   [A_{\mu},A_{\b}] ,  \derm_{\nu}A_{\si} - \derm_{\si}A_{\nu}  > \big)  \\
\notag
&& + m_{\mu\nu }   \cdot m^{\si\b} \cdot  m^{\a\la}    \cdot \big(  <  \derm_{\a}A_{\b} - \derm_{\b}A_{\a} , [A_{\la},A_{\si}] >    +  <  [A_{\a},A_{\b}] , \derm_{\la}A_{\si} - \derm_{\si}A_{\la}  > \big) \\
\notag
 &&  -4 m^{\si\b}  \cdot   <   [A_{\mu},A_{\b}] ,  [A_{\nu},A_{\si}] >      + m_{\mu\nu }  \cdot  m^{\si\b}   \cdot m^{\a\la}   \cdot   <  [A_{\a},A_{\b}] , [A_{\la},A_{\si}] >  \\
     && + O \big(h \cdot  (\derm A)^2 \big)   + O \big(  h  \cdot  A^2 \cdot \derm A \big)     + O \big(  h   \cdot  A^4 \big)  \,  .
\eeaa

Thus, differentiating, and using this time also the fact that $ \Lie_{Z^{I_1}} m_{\mu\nu} = c(I_1 ) \cdot m_{\mu\nu}$, we get

  \beaa
\notag
 && \Lie_{Z^I} g^{\alpha\beta}\derm_\alpha \derm_\beta h_{\mu\nu} \\
  &=&  \Lie_{Z^I} P(\derm_\mu h,\derm_\nu h)  +   \Lie_{Z^I} Q_{\mu\nu}(\derm h,\derm h)   +  \Lie_{Z^I} G_{\mu\nu}(h)(\derm h,\derm h)  \\
\notag
 &&   -4    \sum_{I_1+I_2+I_3  = I}   \hat{c} (I_1)   \cdot m^{\si\b} \cdot  <   \derm_{\mu}  ( \Lie_{Z^{I_2}}A_{\b} )  - \derm_{\b} ( \Lie_{Z^{I_2}}A_{\mu} )   ,  \derm_{\nu} ( \Lie_{Z^{I_3}} A_{\si} ) - \derm_{\si} ( \Lie_{Z^{I_3}} A_{\nu} )  >    \\
 \notag
 &&   +   \sum_{I_1+I_2+I_3+I_4  + I_5= I}  c (I_1) \cdot   \hat{c} (I_2)  \cdot  \hat{c} (I_3) \cdot m_{\mu\nu }  \cdot m^{\si\b}  \cdot  m^{\a\la}  \\
 \notag
 &&  \cdot  <  \derm_{\a} ( \Lie_{Z^{I_4}} A_{\b} )  - \derm_{\b} ( \Lie_{Z^{I_4}}A_{\a} ) , \derm_{\la}  ( \Lie_{Z^{I_5}} A_{\si} ) - \derm_{\si}  ( \Lie_{Z^{I_5}} A_{\la} ) >   \\
 \notag
&&           -4   \sum_{I_1+I_2+I_3+I_4  = I}  \hat{c} (I_1) \cdot m^{\si\b}  \cdot  \big( <   \derm_{\mu} ( \Lie_{Z^{I_2}} A_{\b} ) - \derm_{\b} ( \Lie_{Z^{I_2}} A_{\mu} )  ,  [  \Lie_{Z^{I_3}} A_{\nu}  ,  \Lie_{Z^{I_4}} A_{\si}  ] >   \\
\notag
&& + <   [  \Lie_{Z^{I_2}} A_{\mu}, \Lie_{Z^{I_3}} A_{\b}] ,  \derm_{\nu} ( \Lie_{Z^{I_4}} A_{\si} )  - \derm_{\si} ( \Lie_{Z^{I_4}} A_{\nu} ) > \big)  \\
\notag
&& +  \sum_{I_1+I_2+I_3+I_4 + I_5 + I_6 = I}   c (I_1) \cdot   \hat{c} (I_2)  \cdot  \hat{c} (I_3) \cdot m_{\mu\nu }   \cdot m^{\si\b} \cdot  m^{\a\la}  \\
\notag
&&  \cdot \big(  <  \derm_{\a}  ( \Lie_{Z^{I_4}} A_{\b} )  - \derm_{\b}  ( \Lie_{Z^{I_4}} A_{\a} ) , [  \Lie_{Z^{I_5}} A_{\la},  \Lie_{Z^{I_6}} A_{\si}] >   \\
\notag
&& +  <  [  \Lie_{Z^{I_4}} A_{\a},  \Lie_{Z^{I_5}} A_{\b}] , \derm_{\la}  ( \Lie_{Z^{I_6}} A_{\si} )  - \derm_{\si} ( \Lie_{Z^{I_6}} A_{\la} )  > \big) \\
\notag
 &&  -4 \sum_{I_1+I_2+I_3+I_4 + I_5  = I}   \hat{c} (I_1) \cdot  m^{\si\b}  \cdot   <   [  \Lie_{Z^{I_2}} A_{\mu}, \Lie_{Z^{I_3}}  A_{\b}] ,  [ \Lie_{Z^{I_4}}  A_{\nu},  \Lie_{Z^{I_5}}  A_{\si}] >      \\
 \notag
 && + \sum_{I_1+I_2+I_3+I_4 + I_5 + I_6 + I_7  = I}  c (I_1) \cdot   \hat{c} (I_2)  \cdot  \hat{c} (I_3) \cdot  m_{\mu\nu }  \cdot  m^{\si\b}   \cdot m^{\a\la}   \cdot   <  [ \Lie_{Z^{I_4}}  A_{\a}, \Lie_{Z^{I_5}} A_{\b}] , [ \Lie_{Z^{I_6}}  A_{\la}, \Lie_{Z^{I_7}}  A_{\si}] >  \\
     && +  \Lie_{Z^I} O \big(h \cdot  (\derm A)^2 \big)   +  \Lie_{Z^I}  O \big(  h  \cdot  A^2 \cdot \derm A \big)     + \Lie_{Z^I}   O \big(  h   \cdot  A^4 \big)  \,  .
\eeaa

Decomposing in the null frame $\cal U$ as we did for the zeroth-Lie derivative, we obtain the stated result, since, as we have just exhibited, the structure is the same.

\end{proof}

\begin{lemma}\label{LiederivativesZofsourcestermsforwaveoperatoronAandh}
For any $\si \in \cal U $, we have
     \bea
   \notag
 && |  \Lie_{Z^I}  g^{\la\mu} \derm_{\la}   \derm_{\mu}   A_{\si}  |  \\
 \notag
 &\les& \sum_{|K| + |J| +|M| + |N| \leq |I| } \Big( | \derm ( \Lie_{Z^K} h ) | \cdot  |\rderm ( ( \Lie_{Z^J} A ) |     \\
   \notag
 && +   | \derm ( \Lie_{Z^K}  h )  | \cdot  |\derm  ( \Lie_{Z^J}  A_{L} ) | + | \derm  (  \Lie_{Z^K}  h_{\cal T L} )  | \cdot  |\derm ( \Lie_{Z^J}  A )  |  \\
   \notag
   && +  \big( | \rderm  ( \Lie_{Z^K}  h ) |  +| \derm  ( \Lie_{Z^K}  h_{L \si} )  | \big) \cdot  \big( |\derm  ( \Lie_{Z^J}  A_{L} )  | +|\rderm  ( \Lie_{Z^J}   A ) | \big)  \\
      \notag
         && +   | \rderm  ( \Lie_{Z^K}  h ) |    \cdot  |\derm  ( \Lie_{Z^J}  A ) | \\
 \notag
        && +  \big( | \rderm  ( \Lie_{Z^K}  h ) |   +| \derm  ( \Lie_{Z^K}  h_{L \si} ) |  \big) \cdot  |   \Lie_{Z^J} A  | \cdot  |  \Lie_{Z^M}  A_L |   \\
   \notag
           && +   | \rderm  ( \Lie_{Z^K}   h ) |  \cdot  |   \Lie_{Z^J} A  | \cdot  |   \Lie_{Z^M} A  |   \\
           \notag
                  && +  \big( | \rderm   ( \Lie_{Z^K}  h ) |   +| \derm  ( \Lie_{Z^K}  h_{e_{a} \si} ) |  \big) \cdot    |   \Lie_{Z^J}  A_{e_{a}} |  \cdot  |  \Lie_{Z^M}  A_L |   \\
          \notag
 &&   +   |  \Lie_{Z^K}  A_L  | \cdot    | \derm_\si   ( \Lie_{Z^J}  A_{\underline{L}}  ) |      +    |    \Lie_{Z^K}  A  | \cdot  |   \Lie_{Z^J}  A_\si  | \cdot  |   \Lie_{Z^M}  A_L |   \\
    \notag
    &&   +  |   \Lie_{Z^K}  A  | \cdot   | \derm_\si  ( \Lie_{Z^J} A_L ) |     \\
       \notag
 &&  + |  \Lie_{Z^K}  A_{e_{a}}  | \cdot    | \derm_\si  ( \Lie_{Z^J}  A_{e_{a}} ) |    +    |    \Lie_{Z^K}  A_{e_{a}}  | \cdot  |    \Lie_{Z^J}  A_{e_{a}}  | \cdot  |   \Lie_{Z^M}  A_\si  |   \\
    \notag 
  && +  O( | \Lie_{Z^K}  h | \cdot  | \derm ( \Lie_{Z^J}  h )  | \cdot | \derm ( \Lie_{Z^M}  A) | ) +   O(|  \Lie_{Z^K}  h | \cdot | \derm ( \Lie_{Z^J}  h )| \cdot | \Lie_{Z^M}  A | \cdot | \Lie_{Z^N} A|  ) \\
  \notag
  && +  O( |  \Lie_{Z^K}  h  | \cdot  | \Lie_{Z^J}  A  | \cdot  | \derm (\Lie_{Z^M}   A )  | ) +  O( | \Lie_{Z^K}  h | \cdot | \Lie_{Z^J}  A | \cdot | \Lie_{Z^M}  A  | \cdot  | \Lie_{Z^N} A | ) \Big)  \, . \\
  \eea

For any $\mu, \nu \in \cal U$, we have
    \bea
\notag
 && | \Lie_{Z^I}  g^{\alpha\beta} \derm_\alpha \derm_\beta h_{\mu\nu} |    \\
 \notag
  &\les&  | \Lie_{Z^I}  P(\derm_\mu h,\derm_\nu h) | + | \Lie_{Z^I}  Q_{\mu\nu}(\derm h,\derm h)  | + | \Lie_{Z^I}  G_{\mu\nu}(h)(\derm h,\derm h)  | \\
 \notag
    &&+\sum_{|K| + |J| +|M| + |N|  + |P| \leq |I| } \Big(  \big(   | \derm  ( \Lie_{Z^K}  A_{L} )  | +  | \rderm ( \Lie_{Z^K}  A )  | \big) \cdot  |\derm ( \Lie_{Z^J}  A )  | \\
    \notag
 &&+ \big(  |\derm_\mu ( \Lie_{Z^K}  A_{e_{a}} )  | +  |\rderm ( \Lie_{Z^K}  A ) | \big)  \cdot \big(  |\derm_\nu ( \Lie_{Z^J}  A_{e_{a}} )  | +  |\rderm ( \Lie_{Z^J}  A  ) | \big) \\
   \notag
 && +  | \derm ( \Lie_{Z^K}  A )   |  \cdot |  \Lie_{Z^J}  A | \cdot | \Lie_{Z^M}  A_L | + \big(   | \derm ( \Lie_{Z^K} A_L )  | +  | \rderm ( \Lie_{Z^K}  A ) | \big) \cdot | \Lie_{Z^J}  A |   \cdot | \Lie_{Z^M}  A |   \\
 \notag
&& + \big(   | \derm_\mu  ( \Lie_{Z^K}  A )  | +  | \rderm ( \Lie_{Z^K}  A )   | \big) \cdot |  \Lie_{Z^J}  A_\nu | \cdot |  \Lie_{Z^M}  A | \\
\notag
&& + \big(   | \derm_\nu ( \Lie_{Z^K}  A )  | +  | \rderm ( \Lie_{Z^K}  A )   | \big) \cdot |  \Lie_{Z^J}  A_\mu | \cdot |  \Lie_{Z^M}  A | \\
\notag
&& + \big(   | \derm  ( \Lie_{Z^K}  A_{L} )  | +  | \rderm ( \Lie_{Z^K}  A )   | \big) \cdot | \Lie_{Z^J}  A | \cdot  | \Lie_{Z^M}  A_L |  \\
\notag
&& +  | \derm ( \Lie_{Z^K}  A  )  |  \cdot | \Lie_{Z^J} A | \cdot  |  \Lie_{Z^M}  A_L |  +  | \rderm ( \Lie_{Z^K}  A )  |  \cdot | \Lie_{Z^J}  A | \cdot | \Lie_{Z^M}  A |  \\
 \notag
&& +  |  \Lie_{Z^K} A_\mu | \cdot |  \Lie_{Z^J}  A_\nu | \cdot  |  \Lie_{Z^M}  A | \cdot  |  \Lie_{Z^N}  A |  + | \Lie_{Z^K}  A_L | \cdot  | \Lie_{Z^J}  A | \cdot  | \Lie_{Z^M}  A | \cdot  | \Lie_{Z^N}  A | \\
\notag
     && + O \big( | \Lie_{Z^K}  h | \cdot | \derm ( \Lie_{Z^J} A) | \cdot | \derm ( \Lie_{Z^M} A)|  \big)   + O \big( | \Lie_{Z^K}  h | \cdot |  \Lie_{Z^J} A | \cdot | \Lie_{Z^M} A | \cdot  | \derm ( \Lie_{Z^N}  A )|  \big)   \\
     \notag
     &&  + O \big(  | \Lie_{Z^K}  h |  \cdot  |  \Lie_{Z^J} A | \cdot   | \Lie_{Z^M} A | \cdot  | \Lie_{Z^N} A | \cdot |  \Lie_{Z^P} A | \big)   \Big) \,  .
\eea

\end{lemma}

\begin{proof}
Based on the structure in the null frame $ \cal U$ that we exhibited for $\Lie_{Z^I}   g^{\alpha\beta} \derm_\alpha \derm_\beta h_{\mu\nu} $ and for $ \Lie_{Z^I}   g^{\alpha\beta} \derm_\alpha \derm_\beta A_{\si} $\;, in Lemma \ref{structureofsourcetermsofLiederivaivesZofthehyperbolicsystemonAandh}, and proceeding as we have done for the zeroth-Lie derivative in Lemma \ref{estimatesonsourcetermswithusingLorenzgaugeandwavegaugeestimates}, we get the stated result.
\end{proof}

\begin{lemma}\label{StructureoftheLiederivativesofthesourcetermsofthewaveoperatorforAandh}
In the Lorenz gauge and in wave coordinates, for any $\si \in \cal U $, we have
 \bea
   \notag
&&  |  \Lie_{Z^I}  g^{\la\mu} \derm_{\la}   \derm_{\mu}   A_{\si}  |  \\
 \notag
  &\les&\sum_{|K| + |J| +|M| + |N| \leq |I| } \Big(  | \derm ( \Lie_{Z^K} h ) | \cdot  |\rderm ( \Lie_{Z^J} A ) |       + | \rderm  ( \Lie_{Z^K} h ) | \cdot  |\derm ( \Lie_{Z^J} A ) |  \\
       \notag
           && +  | \rderm  ( \Lie_{Z^K} h )  | \cdot  | \Lie_{Z^J} A  | \cdot    | \Lie_{Z^M} A  |   \\
           \notag
 &&          +| \derm ( \Lie_{Z^K} h_{L \si} )    |  \cdot  | \Lie_{Z^J}  A  | \cdot  | \Lie_{Z^M}  A_L |   \\
        \notag
                  && + | \derm ( \Lie_{Z^K} h_{e_{a} \si} )  |  \cdot    | \Lie_{Z^J}  A_{e_{a}} |  \cdot  | \Lie_{Z^M}  A_L |   \\
          \notag
 &&   +   |  \Lie_{Z^K} A_L  | \cdot    | \derm_\si ( \Lie_{Z^J} A_{\underline{L}}  ) |      +    |  \Lie_{Z^K} A  | \cdot  | \Lie_{Z^J} A_\si  | \cdot  |  \Lie_{Z^M} A_L |   \\
    \notag
    &&   +  | \Lie_{Z^K} A  | \cdot   | \rderm ( \Lie_{Z^J} A )  |     \\
       \notag
 &&  + | \Lie_{Z^K} A_{e_{a}}  | \cdot    | \derm_\si  ( \Lie_{Z^J} A_{e_{a}} ) |    +    |  \Lie_{Z^K} A_{e_{a}}  | \cdot  |  \Lie_{Z^J} A_{e_{a}}  | \cdot  | \Lie_{Z^M} A_\si  |   \\
    \notag 
  && +  O( | \Lie_{Z^K}  h | \cdot  | \derm ( \Lie_{Z^J}  h )  | \cdot | \derm ( \Lie_{Z^M}  A) | ) +   O(|  \Lie_{Z^K}  h | \cdot | \derm ( \Lie_{Z^J}  h )| \cdot | \Lie_{Z^M}  A | \cdot | \Lie_{Z^N} A|  ) \\
  \notag
  && +  O( |  \Lie_{Z^K}  h  | \cdot  | \Lie_{Z^J}  A  | \cdot  | \derm (\Lie_{Z^M}   A )  | ) +  O( | \Lie_{Z^K}  h | \cdot | \Lie_{Z^J}  A | \cdot | \Lie_{Z^M}  A  | \cdot  | \Lie_{Z^N} A | ) \Big)  \, . \\
  \eea

In the Lorenz gauge and in wave coordinates, for any $\mu, \nu \in \cal U$, we have
    \bea
\notag
 && |  \Lie_{Z^I}  g^{\alpha\beta} \derm_\alpha \derm_\beta h_{\mu\nu} |    \\
 \notag
  &\les&   |  \Lie_{Z^I}  P(\derm_\mu h,\derm_\nu h) | + |  \Lie_{Z^I}  Q_{\mu\nu}(\derm h,\derm h)  | + |  \Lie_{Z^I}  G_{\mu\nu}(h)(\derm h,\derm h)  | \\
 \notag
    && +  \sum_{|K| + |J| +|M| + |N|  + |P| \leq |I| }  \Big( | \rderm  ( \Lie_{Z^K}  A )   | \cdot  |\derm   ( \Lie_{Z^J}  A )  | \\
    \notag
 &&+ \big(  |\derm_\mu  ( \Lie_{Z^K}   A_{e_{a}} )   | +  |\rderm  ( \Lie_{Z^K}  A ) | \big)  \cdot \big(  |\derm_\nu   ( \Lie_{Z^J}  A_{e_{a}} )  | +  |\rderm   ( \Lie_{Z^J}  A )  | \big) \\
   \notag
 && +  | \derm  ( \Lie_{Z^K}  A )   |  \cdot |   \Lie_{Z^J}  A | \cdot |   \Lie_{Z^M}  A_L | +  | \rderm  ( \Lie_{Z^K}  A  ) |  \cdot |   \Lie_{Z^J}  A | \cdot   |   \Lie_{Z^M}  A |  \\
 \notag
&& + \big(   | \derm_\mu  ( \Lie_{Z^K} A )  | +  | \rderm ( \Lie_{Z^K}  A  )  | \big) \cdot |  \Lie_{Z^J}  A_\nu | \cdot |  \Lie_{Z^M}   A | \\
\notag
&& + \big(   | \derm_\nu  ( \Lie_{Z^K}  A )   |   +  | \rderm  ( \Lie_{Z^K}  A  )  | \big) \cdot |  \Lie_{Z^J}  A_\mu | \cdot |  \Lie_{Z^M}  A | \\
 \notag
&& +  |  \Lie_{Z^K}  A_\mu | \cdot |  \Lie_{Z^J}  A_\nu | \cdot  |  \Lie_{Z^M}  A | \cdot |  \Lie_{Z^N}  A |  + | \Lie_{Z^K} A_L | \cdot  | \Lie_{Z^J}  A | \cdot  | \Lie_{Z^M}  A | \cdot  | \Lie_{Z^N}  A |  \\
\notag
     && + O \big( | \Lie_{Z^K}  h | \cdot | \derm ( \Lie_{Z^J} A) | \cdot | \derm ( \Lie_{Z^M} A)|  \big)   + O \big( | \Lie_{Z^K}  h | \cdot |  \Lie_{Z^J} A | \cdot | \Lie_{Z^M} A | \cdot  | \derm ( \Lie_{Z^N}  A )|  \big)   \\
     \notag
     &&  + O \big(  | \Lie_{Z^K}  h |  \cdot  |  \Lie_{Z^J} A | \cdot   | \Lie_{Z^M} A | \cdot  | \Lie_{Z^N} A | \cdot |  \Lie_{Z^P} A | \big)   \Big)  \,  .
\eea

\end{lemma}

\begin{proof}
Using the Lorenz gauge estimate that we have shown for Lie derivatives of the Yang-Mills potential in the direction of Minkowski vector fields, in Lemma \ref{LorenzgaugeestimateforgradientofLiederivativesofAL}, 
\bea
\notag
\mid \derm ( \Lie_{Z^J}  A_{L}  ) \mid &\les&  \mid \rderm ( \Lie_{Z^J} A ) \mid \, + \sum_{|K| + |L| \leq |J|} O \big(  | ( \Lie_{Z^K}  h )  | \cdot  | \derm  (  \Lie_{Z^L}  A ) | \big)  , \\
\eea
and injecting in both systems for $\Lie_{Z^I}  g^{\alpha\beta} \derm_\alpha \derm_\beta h_{\mu\nu}  $ and for $\Lie_{Z^I}  g^{\alpha\beta} \derm_\alpha \derm_\beta A $, in what we have shown in Lemma \ref{LiederivativesZofsourcestermsforwaveoperatoronAandh} (as we have done in Lemma \ref{estimatesonsourcetermsforEinstein-Yang-Millssystem} for the zeroth-Lie derivative), we get the desired result.
\end{proof}

\textbf{Notation for the structure of the wave equations for $A$ and $h^1$:}\\

We will use the following notation, in the sense of Lemma \ref{StructureoftheLiederivativesofthesourcetermsofthewaveoperatorforAandh} to illustrate the structure of the coupled wave for the Einstein-Yang-Mills fields:

\subsection{Structure of the wave equation for the “good” components of $A$ and $h^1$}\label{StructureofthewaveequationforthegoodcomponentsofAandhone}\

  \bea
   \notag
 &&   \Lie_{Z^I} \Big( g^{\la\mu} \derm_{\la}   \derm_{\mu}   A_{{\cal T}}   \Big)    \\
 \notag
    &=& \Lie_{Z^I}  \Big[  \derm h  \cdot  \rderm A      +  \rderm  h  \cdot  \derm A   +    A   \cdot   \rderm A  +  \derm  h  \cdot  A^2   +  A^3 \\
    \notag
  && + O( h \cdot  \derm h \cdot  \derm A) + O( h \cdot  A \cdot \derm A)  + O( h \cdot  \derm h \cdot  A^2) + O( h \cdot  A^3)  \Big] \; ,\\
  \eea
  and
 \bea
\notag
 &&  \Lie_{Z^I} \Big(  g^{\alpha\beta} \derm_\alpha \derm_\beta h^1_{ {\cal T} {\cal U}} \Big)    \\
 \notag
    &=& \Lie_{Z^I}  \Big[  P(\derm_{\cal T} h,\derm_{\cal U} h)  + Q_{{\cal T} {\cal U}}(\derm h,\derm h)   +  G_{{\cal T} {\cal U} }(h)(\derm h,\derm h)   \\
    \notag
    &&   +  \rderm A   \cdot  \derm  A   +    A^2 \cdot \derm A     + A^4  \\
\notag
     && + O \big(h \cdot  (\derm A)^2 \big)   + O \big(  h  \cdot  A^2 \cdot \derm A \big)     + O \big(  h   \cdot  A^4 \big) \Big]  \\
  \notag
  && +    \Lie_{Z^I} \Big( g^{\alpha\beta} \derm_\alpha \derm_\beta h^0 \Big) \; ,\\
\eea

where the term $\Lie_{Z^I}  \Big[  P(\derm_{\cal T} h,\derm_{\cal U} h)  + Q_{{\cal T} {\cal U}}(\derm h,\derm h)   +  G_{{\cal T} {\cal U} }(h)(\derm h,\derm h)\Big] $ is estimated in \eqref{estimateonthetermsgeneratedfromtheEinsteinvacuumequationswhichhaveanullstructure}. 

\subsection{Structure of the wave equation for the “bad” components of $A$ and $h^1$}\label{StructureofthewaveequationforthebadcomponentsofAandhone}\

\bea
   \notag
 &&  \Lie_{Z^I} \Big( g^{\la\mu} \derm_{\la}   \derm_{\mu}   A_{{\underline{L}}} \Big)    \\
 \notag
    &=&\Lie_{Z^I}  \Big[   \derm h  \cdot  \rderm A        + \rderm  h \cdot  \derm A  +    A   \cdot     \rderm A   +  \derm  h  \cdot  A^2   +  A^3 \\
    \notag
  && + O( h \cdot  \derm h \cdot  \derm A) + O( h \cdot  A \cdot \derm A)  + O( h \cdot  \derm h \cdot  A^2) + O( h \cdot  A^3) \\
         \notag
 && +  A_L   \cdot   \derm A      +  A_{e_a}  \cdot     \derm A_{e_a}      \Big]   \; ,\\
  \eea
 and
    \bea
\notag
&&   \Lie_{Z^I} \Big( g^{\alpha\beta} \derm_\alpha \derm_\beta h^1_{{\underline{L}} {\underline{L}} } \Big)     \\
\notag
            &=& \Lie_{Z^I}  \Big[   P(\derm_{\cal T} h,\derm_{\cal U} h)  + Q_{{\cal T} {\cal U}}(\derm h,\derm h)   +  G_{{\cal T} {\cal U} }(h)(\derm h,\derm h)   \\
            \notag
         &&    +   \rderm A    \cdot  \derm  A   +  A^2  \cdot  \derm A     +  A^4 \\
    \notag
     && + O \big(h \cdot  (\derm A)^2 \big)   + O \big(  h  \cdot  A^2 \cdot \derm A \big)     + O \big(  h   \cdot  A^4 \big)   \\
         \notag
 &&  +   ( \derm A_{e_{a}} )^2  \Big]  \\
 \notag
 && +     \Lie_{Z^I} \Big(   g^{\alpha\beta} \derm_\alpha \derm_\beta h^0 \Big)  \; , \\
\eea
where the term $\Lie_{Z^I}  \Big[   P(\derm_{\cal T} h,\derm_{\cal U} h)  + Q_{{\cal T} {\cal U}}(\derm h,\derm h)   +  G_{{\cal T} {\cal U} }(h)(\derm h,\derm h)  \Big]$ is estimated in \eqref{estimateonthetermsgeneratedfromtheEinsteinvacuumequationswhichhaveanullstructureplusasquareofderivativeofgoodcomponents}.
 
\section{Pointwise estimates for “good” components of the Einstein-Yang-Mills fields}

We already know, from Lemmas \ref{betterdecayfortangentialderivatives} and \ref{decayrateforfullderivativeintermofZ}, the following estimate for any sufficiently smooth function $f$,
\beaa
(1 + t + |q| ) \cdot |\rpa f | + (1 +  |q| ) \cdot |\pa f | \les \sum_{|I| = 1} |Z^I f | \, .
\eeaa
While this estimate gives a good control on “good” derivatives $|\rpa f |$, we loose control on $|\pa f |$. However, in the region outside $D_t :=\{x\in \R^3;\,
t/2\leq |x|\leq 2t\}$, this estimate does give control on $|\pa f |$, thanks to the following lemma.

\begin{lemma}
Let, $D:=\{(t,x) \in \R \times \R^3;\, t/2\leq |x|\leq 2t\}$, then for any sufficiently smooth function $f$ and for all $(t, x) \nin D$, we have
\beaa
(1 + t + r ) \cdot |\pa f |  \les \sum_{|I| = 1} |Z^I f | \, .
\eeaa

\end{lemma}
\begin{proof}
From Lemma \ref{decayrateforfullderivativeintermofZ}, we have
\beaa
 (1 +  |q| ) \cdot |\pa f | \les \sum_{|I| = 1} |Z^I f | \, .
\eeaa
Now, we look at the cases where $r < \frac{t}{2}$ and where $r > 2t $:

\textbf{Case $r < \frac{t}{2}$:}\\

In that case 
$$ q = r - t < \frac{t}{2} -t = - \frac{t}{2} \leq 0 .$$
Thus, $$q  < 0 .$$
 Since $- r > - \frac{t}{2}$, then $|q| = - q = t - r > t  - \frac{t}{2} = \frac{t}{2} > r$, which leads to
\bea
|q|  &>& \frac{t}{2} \, , \\
|q|  &>& r \, .
\eea

\textbf{Case $r > 2t $:}\\

In that case  $$ q = r - t > 2t -t = t \geq 0 .$$
Thus, $$q  > 0 .$$
 Since $|q| = q = r - t > t $, we have
\bea
|q|  &>& t\, .
\eea
Since $r > 2t $, we have $-2t > -r $, which gives $-t > - \frac{r}{2}$. Thus,
\beaa
|q| = r -t  &>& r  - \frac{r}{2} =\frac{r}{2} .
\eeaa
As a result in all cases, we get for all $(t, x) \nin D$,
\beaa
(1 + t + r )\cdot |\pa f | \les  (1 +  \frac{|q|}{2} + \frac{|q|}{2}  ) \cdot |\pa f | \les \sum_{|I| = 1} |Z^I f | \, .
\eeaa
\end{proof}

In the region $(t, x) \in D$, we lose control on $|\pa f|$. To gain control, we use the special structure satisfied by the fields. In fact, we will use an estimate on the $L^{\infty}$ norm of solutions to wave equations with sources of the form
\bea
 g^{\la\a} \derm_{\la}   \derm_{\a}   \Phi_{\mu\nu}= S_{\mu\nu} \, , 
\eea
where $S_{\mu\nu}$ is the source term, provided that in the region
\bea
\overline{D}_t &:=&D_t \cap \overline{C} = \{x\in\R^3;\, t/2\leq |x|\leq 2t\} \cap \overline{C}  \; ,
\eea
the tensor $H^{\alpha\beta}=g^{\alpha\beta}-m^{\alpha\beta}$ satisfies in $\overline{C}$\;,
\bea
\notag
 |H|\leq \frac{\varepsilon^\prime}{4}\;,\qquad\int_0^{\infty}\!\!
\|\,H(t,\cdot)\|_{L^\infty(\overline{D}_t)}\frac{dt}{1+t} \leq
\frac{\varepsilon^\prime}{4}\;,\qquad |H|_{L{\cal
T}} \leq\frac{\varepsilon^\prime}{4}\, \frac{|q|+1}{1+t+|x|} \; . \\
\eea
This general estimate is derived by Lindblad-Rodnianski in Corollary 7.2 in \cite{LR10} and states that solutions of such wave equations satisfy the following estimate for any $U,V \in \cal U $, and for any $x \in \overline{D}_t=\{x\in\R^3;\, t/2\leq |x|\leq 2t\}$,

\bea\label{LinfinitynormestimateongradientderivedbyLondbladRodnianski}
\notag
&& (1+t+|x|) \cdot |\varpi(q) \cdot   \pa  \Phi_{UV} (t,x)| \\
\notag
&\les& \!\sup_{0\leq \tau\leq t} \sum_{|I|\leq 1}\|\,\varpi(q) \cdot \Lie_{Z^I} \!  \Phi (\tau,\cdot)\|_{L^\infty (\Sigma^{ext}_{\tau} ) }\\
\notag
&& + \int_0^t\Big( \varepsilon^\prime \cdot \alpha \cdot \|\varpi(q) \cdot \pa  \Phi_{UV} (t,\cdot) \|_{L^\infty (\Sigma^{ext}_{\tau} ) } +(1+\tau) \cdot \| \varpi(q) \cdot S_{UV} (\tau,\cdot) \|_{L^\infty(\overline{D}_\tau)} \\
 \notag
&& +\sum_{|I|\leq 2} (1+\tau)^{-1}\cdot \| \varpi(q) \cdot  \Lie_{Z^I}  \Phi (\tau,\cdot)\|_{L^\infty(\overline{D}_\tau)}\Big)\, d\tau \; , \\
\eea
where for $1 + \gamma'\ge 0$, 
\begin{equation}\label{definitionofnewweighttoupgrade}
\varpi :=\varpi(q) :=\begin{cases}
(1+|q|)^{1+\gamma^\prime},\quad\text{when }\quad q>0\; , \\
     1 \,\quad\text{when }\quad
      q<0 \; ,\end{cases}
\end{equation}
where
\bea
\alpha=\max (1+\gamma’,0) = 1+\gamma’ \ge 0  \, ,
\eea
and where for $\Phi_{\a\b}$ scalar field,
\bea
 | \pa \Phi_{UV} |^2 :=   Euc^{\mu\nu}    { \der^{(\bf{m})}}_{\mu} \Phi_{UV} \ \cdot    { \der^{(\bf{m})}}_{\nu} \Phi_{UV} , 
 \eea
 and similarly for a tensor $\Phi_{\a\b}$ valued in the Lie algebra,
 \bea
 | \pa \Phi_{UV} |^2 :=   Euc^{\mu\nu}   <  { \der^{(\bf{m})}}_{\mu} \Phi_{UV} \ ,    { \der^{(\bf{m})}}_{\nu} \Phi_{UV} >.
 \eea
 
Since such an estimate also holds true in the region  $(t, x) \nin D$, the estimate therefore holds true for all  $(t, x)$ in the space-time. We want now apply it to the Einstein-Yang-Mills system in the Lorenz and harmonic gauges. For this, we need to prove first that the assumptions of the theorem are satisfied.

\begin{lemma}\label{aprioridecayestimates}
Let $M\leq \eps$. By choosing $\gamma > \delta > 0 $, we have for all $|I|$, in the entire exterior region $\overline{C} \subset \{ q \geq q_0 \} $, on one hand the following estimates for the full perturbation metric $h$,
   \beaa
 \notag
&&  |  \rderm ( \Lie_{Z^I} h ) (t,x)  |  \\
\notag
 &\leq& \begin{cases} c (\delta) \cdot c (\gamma) \cdot C ( |I| ) \cdot E ( |I| + 3)  \cdot \frac{\eps }{(1+t+|q|)^{2-\delta} (1+|q|)^{\de}},\quad\text{when }\quad q>0,\\
       C ( |I| ) \cdot E ( |I| + 3)  \cdot \eps \cdot \frac{ (1+|q|)^{\frac{1}{2} }}{ (1+t+|q|)^{2-\delta} } \,\quad\text{when }\quad q<0 , \end{cases} \\
      \eeaa
and
              \beaa
 \notag
 |\derm  ( \Lie_{Z^I} h ) (t,x)  |    &\leq& \begin{cases} C ( |I| ) \cdot E ( |I| + 2)  \cdot \frac{\eps }{(1+t+|q|)^{1-\delta} (1+|q|)^{1+\de}},\quad\text{when }\quad q>0,\\
       C ( |I| ) \cdot E ( |I| + 2)  \cdot \frac{\eps  }{(1+t+|q|)^{1-\delta}(1+|q|)^{\frac{1}{2} }}  \,\quad\text{when }\quad q<0 , \end{cases} \\
      \eeaa
      and
                 \beaa
 \notag
 |   \Lie_{Z^I} h (t,x)  | &\leq& \begin{cases} c (\delta) \cdot c (\gamma) \cdot C ( |I| ) \cdot E ( |I| + 2) \cdot  \frac{\eps}{ (1+ t + | q | )^{1-\delta }  (1+| q |   )^{\de}}  ,\quad\text{when }\quad q>0,\\
\notag
    C ( |I| ) \cdot E ( |I| + 2) \cdot  \frac{\eps}{ (1+ t + | q | )^{1-\delta }  } (1+| q |   )^{\frac{1}{2} }  , \,\quad\text{when }\quad q<0 . \end{cases} \\ 
    \eeaa
        On the other hand, we have the following estimates on the potential $A$, 
       \beaa
 \notag
&&  |  \rderm ( \Lie_{Z^I} A ) (t,x)  |  \\
\notag
 &\leq& \begin{cases} c (\gamma) \cdot C ( |I| ) \cdot E ( |I| + 3)  \cdot \frac{\eps }{(1+t+|q|)^{2-\delta} (1+|q|)^{\ga}},\quad\text{when }\quad q>0,\\
       C ( |I| ) \cdot E ( |I| + 3)  \cdot \eps \cdot \frac{ (1+|q|)^{\frac{1}{2} }}{ (1+t+|q|)^{2-\delta} } \,\quad\text{when }\quad q<0 , \end{cases} \\
      \eeaa
and
                  \beaa
 \notag
|\derm ( \Lie_{Z^I} A ) (t,x)    |    &\leq& \begin{cases} C ( |I| ) \cdot E ( |I| + 2)  \cdot \frac{\eps }{(1+t+|q|)^{1-\delta} (1+|q|)^{1+\ga}},\quad\text{when }\quad q>0,\\
       C ( |I| ) \cdot E ( |I| + 2)  \cdot \frac{\eps  }{(1+t+|q|)^{1-\delta}(1+|q|)^{\frac{1}{2} }}  \,\quad\text{when }\quad q<0 , \end{cases} \\
      \eeaa
      and
       \beaa
 \notag
 |   \Lie_{Z^I} A (t,x)  | &\leq& \begin{cases} c (\gamma) \cdot C ( |I| ) \cdot E ( |I| + 2) \cdot  \frac{\eps}{ (1+ t + | q | )^{1-\delta }  (1+| q |   )^{\ga}}  ,\quad\text{when }\quad q>0,\\
\notag
    C ( |I| ) \cdot E ( |I| + 2) \cdot  \frac{\eps}{ (1+ t + | q | )^{1-\delta }  } (1+| q |   )^{\frac{1}{2} }  , \,\quad\text{when }\quad q<0 . \end{cases} \\ 
    \eeaa

In particular, we then have in $\overline{C} \subset \{ q \geq q_0\}$,
   \beaa
   \notag
 |  \rderm (  \Lie_{Z^I} h ) (t,x)  |   &\leq& C(q_0) \cdot c (\gamma) \cdot c (\delta) \cdot C ( |I| ) \cdot E ( |I| + 3)  \cdot \frac{\eps }{(1+t+|q|)^{2-\delta} \cdot (1+|q|)^{\de} } . \\
 \notag
 |  \rderm ( \Lie_{Z^I} A ) (t,x)  |    &\leq& C(q_0) \cdot c (\gamma) \cdot c (\delta) \cdot C ( |I| ) \cdot E ( |I| + 3)  \cdot \frac{\eps }{(1+t+|q|)^{2-\delta} \cdot (1+|q|)^{\ga} } . \\
\notag
 |\derm ( \Lie_{Z^I} h ) (t,x)  |    &\leq& C(q_0) \cdot C ( |I| ) \cdot E ( |I| + 2)  \cdot \frac{\eps }{(1+t+|q|)^{1-\delta} \cdot (1+|q|)^{1+\de}} . \\
\notag
|\derm  ( \Lie_{Z^I} A ) (t,x)  |     &\leq& C(q_0) \cdot C ( |I| ) \cdot E ( |I| + 2)  \cdot \frac{\eps }{(1+t+|q|)^{1-\delta} \cdot (1+|q|)^{1+\ga}} . \\
 \notag
|   \Lie_{Z^I} h (t,x)  |   &\leq& C(q_0) \cdot  c (\gamma)  \cdot c (\delta) \cdot C ( |I| )  \cdot E ( |I| + 2) \frac{\eps }{(1+t+|q|)^{1-\delta}  \cdot (1+|q|)^{\de} } . \\
\notag
|  \Lie_{Z^I} A (t,x)  |   &\leq&   C(q_0)  \cdot  c (\gamma) \cdot C ( |I| ) \cdot E ( |I| +2)  \cdot \frac{\eps }{(1+t+|q|)^{1-\delta} \cdot (1+|q|)^{ \ga}} . \\
      \eeaa
   In particular,   
 \beaa
 |  \Lie_{Z^I} A (t,x)  |^2   &\leq&   C(q_0) \cdot  c (\gamma)  \cdot C ( |I| ) \cdot E ( |I| +2)  \cdot \frac{\eps^2 }{(1+t+|q|)^{2-2\delta} \cdot (1+|q|)^{2 \ga}} , \\
 \notag
|  \Lie_{Z^I} A (t,x)  |^3   &\leq&  C(q_0) \cdot  c (\gamma)  \cdot C ( |I| ) \cdot E ( |I| +2)  \cdot \frac{\eps^3  }{(1+t+|q|)^{3-3\delta}\cdot (1+|q|)^{3 \ga}  } .
      \eeaa

\end{lemma}

\begin{proof}

We have shown in Lemma \ref{estimategoodderivatives}, that for all $I$,       
   \beaa
 \notag
|  \rpa ( \Lie_{Z^I} h^1 ) (t,x)  |   &\leq& \begin{cases} c (\gamma) \cdot C ( |I| ) \cdot E ( |I| + 3)  \cdot \frac{\eps }{(1+t+|q|)^{2-\delta} (1+|q|)^{\gamma}},\quad\text{when }\quad q>0,\\
       C ( |I| ) \cdot E ( |I| + 3)  \cdot \eps \cdot \frac{ (1+|q|)^{\frac{1}{2} }}{ (1+t+|q|)^{2-\delta} } \,\quad\text{when }\quad q<0 , \end{cases} \\
      \eeaa
      
and
 \beaa
 \notag
|  \rpa ( \Lie_{Z^I} h^0 ) (t,x)  |   &\leq& C ( |I|  ) \cdot \frac{\eps }{(1+t+|q|)^{2}  } \\
 \notag
&\leq & \begin{cases}  C ( |I| ) \cdot \frac{\eps }{(1+t+|q|)^{2-\delta} (1+|q|)^{\delta}},\quad\text{when }\quad q>0,\\
 \notag
       C ( |I|  )  \cdot \eps \cdot \frac{ (1+|q|)^{\frac{1}{2} }}{ (1+t+|q|)^{2-\delta} } \,\quad\text{when }\quad q<0 . \end{cases} \\
      \eeaa
Thus, for  $\gamma > \delta > 0 $,
   \bea
 \notag
|  \rpa ( \Lie_{Z^I} h ) (t,x)  |  &\leq& |  \rpa ( \Lie_{Z^I} h^0 ) (t,x)  | + |  \rpa ( \Lie_{Z^I} h^1 ) (t,x)  |  \\
 \notag
& \leq & \begin{cases} c (\de) \cdot c (\gamma) \cdot C ( |I| ) \cdot E ( |I| + 3)  \cdot \frac{\eps }{(1+t+|q|)^{2-\delta} (1+|q|)^{\delta}},\quad\text{when }\quad q>0,\\
 \notag
       C ( |I| ) \cdot E ( |I| + 3)  \cdot \eps \cdot \frac{ (1+|q|)^{\frac{1}{2} }}{ (1+t+|q|)^{2-\delta} } \,\quad\text{when }\quad q<0 , \end{cases} \\
        \notag
  &\leq & C(q_0) \cdot c (\de)  \cdot  c (\gamma) \cdot C ( |I| ) \cdot E ( |I| + 3)  \cdot \frac{\eps }{(1+t+|q|)^{2-\delta} (1+|q|)^{\delta}} \; .\\
      \eea

We also have showed in Lemma \ref{apriordecayestimatesfrombootstrapassumption}, that

       \beaa
 \notag
  |\derm ( \Lie_{Z^I} A ) (t,x)  |  + |\derm ( \Lie_{Z^I} h^1 ) (t,x)  |    &\leq& \begin{cases} C ( |I| ) \cdot E ( |I| + 2)  \cdot \frac{\eps }{(1+t+|q|)^{1-\delta} (1+|q|)^{1+\gamma}},\quad\text{when }\quad q>0,\\
       C ( |I| ) \cdot E ( |I| + 2)  \cdot \frac{\eps  }{(1+t+|q|)^{1-\delta}(1+|q|)^{\frac{1}{2} }}  \,\quad\text{when }\quad q<0 , \end{cases} 
      \eeaa
      and    
 \beaa
 \notag
|\derm  ( \Lie_{Z^I} h^0 ) (t,x)  |   &\leq& C ( |I| )   \cdot \frac{\eps }{(1+t+|q|)^{2} }
      \eeaa
     Hence, for $ \ga > \de$, 
     
          \bea
 \notag
 |\derm (  \Lie_{Z^I} h)  (t,x)  |   &\les &  |\derm (  \Lie_{Z^I} h^0)  (t,x)  |  +  |\derm (  \Lie_{Z^I} h^1)  (t,x)  |  \\
 \notag
  &\leq& \begin{cases} C ( |I| ) \cdot E ( |I| + 2)  \cdot \frac{\eps }{(1+t+|q|)^{1-\delta} (1+|q|)^{1+\delta}},\quad\text{when }\quad q>0,\\
       C ( |I| ) \cdot E ( |I| + 2)  \cdot \frac{\eps  }{(1+t+|q|)^{1-\delta}(1+|q|)^{\frac{1}{2} }}  \,\quad\text{when }\quad q<0 , \end{cases} \\
      \eea

 and therefore,     
             \bea
 \notag
 |\derm  ( \Lie_{Z^I} h ) (t,x)  |    &\leq& C(q_0) \cdot  C ( |I| ) \cdot E ( |I| + 2)  \cdot \frac{\eps }{(1+t+|q|)^{1-\delta}  \cdot (1+|q|)^{1+\delta} } . \\
      \eea

We have shown in Lemmas \ref{spatialdecayoninitialdeataforanapriorestimateontheEinsteinYangMillsfieldswithoutgradiant} and \ref{apriorestimateontheEinsteinYangMillsfieldswithoutgradiant}, that for $\gamma > \max\{0, \delta -1\} $ (which is true when $\gamma > \delta $ ), we have:
  
   \beaa
 \notag
 |  \Lie_{Z^I} A (t,  | x | \cdot \Om)   |  +  |  \Lie_{Z^I} h^1 (t,  | x | \cdot \Om)  |   &\leq& \begin{cases} c (\gamma) \cdot C ( |I| ) \cdot E ( |I| + 2) \cdot  \frac{\eps}{ (1+ t + | q | )^{1-\delta }  (1+| q |   )^{\gamma}}  ,\quad\text{when }\quad q>0,\\
    C ( |I| ) \cdot E ( |I| + 2) \cdot  \frac{\eps}{ (1+ t + | q | )^{1-\delta }  } (1+| q |   )^{\frac{1}{2} }  , \,\quad\text{when }\quad q<0 , \end{cases} \\
      \eeaa
     
     and 
         \beaa
 \notag
|   \Lie_{Z^I} h^0 (t,x)  |   &\leq&  C ( |I| )   \cdot \frac{\eps }{(1+t+|q|) } \; .\\
      \eeaa

Consequently, for $\ga > \de$, 
 \bea
 \notag
|   \Lie_{Z^I} h (t,x)  | &\leq& \begin{cases} c (\delta) \cdot c (\gamma) \cdot C ( |I| ) \cdot E ( |I| + 2) \cdot  \frac{\eps}{ (1+ t + | q | )^{1-\delta }  (1+| q |   )^{\delta}}  ,\quad\text{when }\quad q>0,\\
\notag
    C ( |I| ) \cdot E ( |I| + 2) \cdot  \frac{\eps}{ (1+ t + | q | )^{1-\delta }  } (1+| q |   )^{\frac{1}{2} }  , \,\quad\text{when }\quad q<0 , \end{cases} \\ 
    \eea
    and
 \bea
 \notag
|   \Lie_{Z^I} A (t,x)  | &\leq& \begin{cases}  c (\gamma) \cdot C ( |I| ) \cdot E ( |I| + 2) \cdot  \frac{\eps}{ (1+ t + | q | )^{1-\delta }  (1+| q |   )^{\ga}}  ,\quad\text{when }\quad q>0,\\
\notag
    C ( |I| ) \cdot E ( |I| + 2) \cdot  \frac{\eps}{ (1+ t + | q | )^{1-\delta }  } (1+| q |   )^{\frac{1}{2} }  , \,\quad\text{when }\quad q<0 . \end{cases} \\ 
    \eea
    
Therefore in the exterior $\overline{C}$, we have
    \bea
    \notag
 |   \Lie_{Z^I} h (t,x)  |      &\leq& C(q_0) \cdot  c (\gamma)  \cdot c (\delta) \cdot C ( |I| )  \cdot   E ( |I| + 2) \cdot \frac{\eps}{(1+t+|q|)^{1-\delta} \cdot  (1+| q |   )^{\delta} } \; . \\
      \eea
  \bea
    \notag
 |   \Lie_{Z^I} A (t,x)  |      &\leq& C(q_0) \cdot  c (\gamma)   \cdot C ( |I| )  \cdot   E ( |I| + 2) \cdot  \frac{\eps}{(1+t+|q|)^{1-\delta} \cdot  (1+| q |   )^{\ga} } . \\
      \eea
      
Yet, from Lemma \ref{equivalenceoftworestrictednormsongradients}, we have that
 \beaa
|\rderm  A |^2  &\sim&  |  \rpa A |^2   := \sum_{\a  \in  \{t, x^1, x^2, x^3 \}} |  \rpa A_\a |^2    ,
\eeaa
and therefore,
   \bea
 \notag
&& |  \rderm ( \Lie_{Z^I} A ) (t,x)  |   \\
\notag
 &\leq& \begin{cases}  c (\gamma) \cdot C ( |I| ) \cdot E ( |I| + 3)  \cdot \frac{\eps }{(1+t+|q|)^{2-\delta} (1+|q|)^{\gamma}},\quad\text{when }\quad q>0,\\
       C ( |I| ) \cdot E ( |I| + 3)  \cdot \eps \cdot \frac{ (1+|q|)^{\frac{1}{2} }}{ (1+t+|q|)^{2-\delta} } \,\quad\text{when }\quad q<0 . \end{cases} \\
      \eea
Hence,
   \bea
   \notag
&& |  \rderm ( \Lie_{Z^I} A ) (t,x)  |  \\
   \notag
  &\leq&  C(q_0) \cdot c (\gamma)  \cdot C ( |I| ) \cdot E ( |I| + 3)  \cdot \frac{\eps }{(1+t+|q|)^{2-\delta} \cdot  (1+|q|)^{\gamma}} .\\
      \eea

\end{proof}

\begin{lemma}

Let $M\leq \eps$\,, and let $\gamma >\delta  $ and $\delta \leq \frac{1}{2}$. Then, in wave coordinates, in the exterior region $\overline{C} \subset \{ q \geq q_0 \} $, we have

     \beaa
        \notag
&& | \derm   ( \Lie_{Z^I} h_{\cal T L} ) |  + | \derm   ( \Lie_{Z^I} H_{\cal T L} ) | \\
       &\les& \begin{cases} c (\delta) \cdot c (\gamma) \cdot  C ( |I| ) \cdot E ( |I| + 3)  \cdot \big( \frac{\eps }{(1+t+|q|)^{2-2\delta} (1+|q|)^{2\delta}}  \big),\quad\text{when }\quad q>0 \; ,\\
       \notag
       C ( |I| ) \cdot E ( |I| + 3)  \cdot \big( \frac{ \eps \cdot (1+|q|)^{\frac{1}{2} }}{ (1+t+|q|)^{2-2\delta} } \big) \,\quad\text{when }\quad q<0 \; . \end{cases} 
\eeaa

As a result, more precisely, in $\overline{C} \subset \{ q \geq q_0 \} $, we have 
           \bea
        \notag
&&| \derm   ( \Lie_{Z^I} h_{\cal T L} ) | +  | \derm   ( \Lie_{Z^I} H_{\cal T L} ) |  \\
        \notag
               &\les& C(q_0) \cdot      c (\delta) \cdot c (\gamma) \cdot  C ( |I| ) \cdot E ( |I| + 3)  \cdot \big( \frac{\eps }{(1+t+|q|)^{2-2\delta} (1+|q|)^{2\delta}}  \big) \; , \\
\eea

\end{lemma}

\begin{proof}

In wave coordinates, we showed in Lemma \ref{wavecoordinatesestimateonLiederivativesZonmetric}, that

\beaa
| \derm ( \Lie_{Z^I} h_{\cal T L} )  | &\les& \sum_{|J| \leq |I| } | \rderm  ( \Lie_{Z^J}  h)  |  + \sum_{|K|+ |L| \leq |I|}  O (|\Lie_{Z^K} h| \cdot |\derm ( \Lie_{Z^L} h ) | )
\eeaa

From Lemma \ref{aprioridecayestimates}, we have for $\gamma >\delta  $, for all $|J| \leq |I|$, on one hand
   \beaa
 \notag
&& |  \rderm ( \Lie_{Z^J} h ) (t,x)  |  \\
\notag
 &\leq& \begin{cases} c (\delta) \cdot c (\gamma) \cdot C ( |I| ) \cdot E ( |I| + 3)  \cdot \frac{\eps }{(1+t+|q|)^{2-\delta} (1+|q|)^{\de}},\quad\text{when }\quad q>0  \;,\\
       C ( |I| ) \cdot E ( |I| + 3)  \cdot \eps \cdot \frac{ (1+|q|)^{\frac{1}{2} }}{ (1+t+|q|)^{2-\delta} } \,\quad\text{when }\quad q<0  \; ,\end{cases} \\
      \eeaa
      and on the other hand,
 \bea
 \notag
|   \Lie_{Z^J} h (t,x)  | &\leq& \begin{cases} c (\delta) \cdot c (\gamma) \cdot C ( |I| ) \cdot E ( |I| + 2) \cdot  \frac{\eps}{ (1+ t + | q | )^{1-\delta }  (1+| q |   )^{\de}}  ,\quad\text{when }\quad q>0  \; ,\\
\notag
    C ( |I| ) \cdot E ( |I| + 2) \cdot  \frac{\eps}{ (1+ t + | q | )^{1-\delta }  } (1+| q |   )^{\frac{1}{2} }  , \,\quad\text{when }\quad q<0  \; , \end{cases} \\ 
    \eea
    and
                 \beaa
 \notag
|\derm ( \Lie_{Z^J} h ) (t,x)  |    &\leq& \begin{cases} C ( |I| ) \cdot E ( |I| + 2)  \cdot \frac{\eps }{(1+t+|q|)^{1-\delta} (1+|q|)^{1+\de}},\quad\text{when }\quad q>0  \;,\\
       C ( |I| ) \cdot E ( |I| + 2)  \cdot \frac{\eps  }{(1+t+|q|)^{1-\delta}(1+|q|)^{\frac{1}{2} }}  \,\quad\text{when }\quad q<0 \, . \end{cases} \\
      \eeaa
Therefore, we have for all $|K|+ |L| \leq |I|$, 
                 \bea
 \notag
&& |   \Lie_{Z^K} h (t,x) |  \cdot |\derm   ( \Lie_{Z^L} h ) (t,x)  |  \\
 \notag
  &\leq& \begin{cases} c (\delta) \cdot  c (\gamma) \cdot C ( |I| ) \cdot E ( |I| + 2)  \cdot \frac{\eps^2 }{(1+t+|q|)^{2-2\delta} (1+|q|)^{1+2\de}},\quad\text{when }\quad q>0  \; ,\\
       C ( |I| ) \cdot E ( |I| + 2)  \cdot \frac{\eps^2  }{(1+t+|q|)^{2-2\delta} }  \,\quad\text{when }\quad q<0 \; . \end{cases} \\
      \eea

      As a result,
      \beaa
&& | \derm ( \Lie_{Z^I} h_{\cal T L} )  |  \\
&\les& \begin{cases} c (\delta) \cdot c (\gamma) \cdot C ( |I| ) \cdot E ( |I| + 3)  \cdot \big( \frac{\eps }{(1+t+|q|)^{2-\delta} (1+|q|)^{\de}} +  \frac{\eps^2 }{(1+t+|q|)^{2-2\delta} (1+|q|)^{1+2\de}}\big),\quad\text{when }\quad q>0  \; ,\\
       C ( |I| ) \cdot E ( |I| + 3)  \cdot \big( \frac{ \eps \cdot (1+|q|)^{\frac{1}{2} }}{ (1+t+|q|)^{2-\delta} } + \frac{\eps^2  }{(1+t+|q|)^{2-2\delta} }  \big) \,\quad\text{when }\quad q<0 \; .
       \end{cases}
\eeaa

      We get for $\eps < 1$, 
        \beaa
&& | \derm ( \Lie_{Z^I}  h_{\cal T L} ) |  \\
       &\les& \begin{cases} c (\delta) \cdot c (\gamma) \cdot C ( |I| ) \cdot E ( |I| + 3)  \cdot \big( \frac{\eps }{(1+t+|q|)^{2-2\delta} (1+|q|)^{2\delta}}  \big),\quad\text{when }\quad q>0  \; ,\\
        C ( |I| ) \cdot E ( |I| + 3)  \cdot \big( \frac{ \eps \cdot (1+|q|)^{\frac{1}{2} }}{ (1+t+|q|)^{2-2\delta} } \big) \,\quad\text{when }\quad q<0 \;  . \end{cases}
\eeaa

       Whereas to $H$, we had also showed that 
       \beaa
| \derm ( \Lie_{Z^I} H_{\cal T L} )  | &\les& \sum_{|J| \leq |I| } | \rderm  ( \Lie_{Z^J}  H)  |  + \sum_{|K|+ |L| \leq |I|}  O (|\Lie_{Z^K} H| \cdot |\derm ( \Lie_{Z^L} H ) | )
\eeaa

The same estimates that we used here for $h$, also hold for $H$, and therefore, similarly, we get the same result for $H$ as well.

\end{proof}

\begin{lemma} \label{aprioriestimatesonZLiederivativesofBIGH}

      Let $M\leq \eps$. By choosing $\gamma > \delta $, we have for all $|I|$, $\delta \leq \frac{1}{2}$,\,$\eps \leq 1$,

              \beaa
 \notag
 |\derm   H (t,x)  |    &\les& \begin{cases}  E ( 2)  \cdot \frac{\eps }{(1+t+|q|)^{1-\delta} (1+|q|)^{1+\de}},\quad\text{when }\quad q>0,\\
        E (  2)  \cdot \frac{\eps  }{(1+t+|q|)^{1-\delta}(1+|q|)^{\frac{1}{2} }}  \,\quad\text{when }\quad q<0 , \end{cases} \\
      \eeaa
      and
       \beaa
 \notag
|   H (t,x)  | &\les& \begin{cases} c (\delta) \cdot c (\gamma) \cdot  E (  2) \cdot  \frac{\eps}{ (1+ t + | q | )^{1-\delta }  (1+| q |   )^{\de}}  ,\quad\text{when }\quad q>0,\\
\notag
    E ( 2) \cdot  \frac{\eps}{ (1+ t + | q | )^{1-\delta }  } (1+| q |   )^{\frac{1}{2} }  , \,\quad\text{when }\quad q<0 . \end{cases} \\ 
    \eeaa

\end{lemma}

\begin{proof}

We proved in \cite{G4}, that
  \beaa
H^{\mu\nu}=-h^{\mu\nu}+ O^{\mu\nu}(h^2) .
\eeaa

Since, from Lemma \ref{aprioridecayestimates}, for all $|I|$,

 \bea
 \notag
|   \Lie_{Z^I} h (t,x)  | &\leq& \begin{cases} c (\delta) \cdot c (\gamma) \cdot C ( |I| ) \cdot E ( |I| + 2) \cdot  \frac{\eps}{ (1+ t + | q | )^{1-\delta }  (1+| q |   )^{\de}}  ,\quad\text{when }\quad q>0,\\
\notag
    C ( |I| ) \cdot E ( |I| + 2) \cdot  \frac{\eps}{ (1+ t + | q | )^{1-\delta }  } (1+| q |   )^{\frac{1}{2} }  , \,\quad\text{when }\quad q<0 , \end{cases} 
      \eea
      
       we get
  \beaa
 \notag
&& |   H (t,x)  | \\
  &\les& \begin{cases} c (\delta) \cdot c (\gamma)  \cdot E (  2) \cdot \big( \frac{\eps}{ (1+ t + | q | )^{1-\delta }  (1+| q |   )^{\de}}  + O(  \frac{\eps^2}{ (1+ t + | q | )^{2-2\delta }  (1+| q |   )^{2\de}}  ) \big),\quad\text{when }\quad q>0,\\
\notag
     E (  2) \cdot \big( \frac{\eps}{ (1+ t + | q | )^{1-\delta }  } (1+| q |   )^{\frac{1}{2} }  + O(  \frac{\eps^2}{ (1+ t + | q | )^{2-2\delta }  } (1+| q |   )  ) \big), \,\quad\text{when }\quad q<0 , \end{cases} \\
      &\les& \begin{cases} c (\delta) \cdot c (\gamma) \cdot E (  2) \cdot  \frac{\eps}{ (1+ t + | q | )^{1-\delta }  (1+| q |   )^{\de}}  ,\quad\text{when }\quad q>0,\\
\notag
    E ( 2) \cdot \big( \frac{\eps}{ (1+ t + | q | )^{1-\delta }  } (1+| q |   )^{\frac{1}{2} }  + O(  \frac{\eps}{ (1+ t + | q | )^{1-\delta }  } (1+| q |   )^{\frac{1}{2} }     \cdot  \frac{\eps \cdot (1+| q |   )^{\frac{1}{2} } }{ (1+ t + | q | )^{1-\delta }  }   \big), \,\quad\text{when }\quad q<0 .\end{cases} 
      \eeaa
      
      Thus,
        \bea
 \notag
&& |   H (t,x)  | \\
      &\les& \begin{cases} c (\delta) \cdot c (\gamma)  \cdot E (  2) \cdot  \frac{\eps}{ (1+ t + | q | )^{1-\delta }  (1+| q |   )^{\de}}  ,\quad\text{when }\quad q>0,\\
\notag
 E (  2) \cdot \frac{\eps}{ (1+ t + | q | )^{1-\delta }  } (1+| q |   )^{\frac{1}{2} }  , \,\quad\text{when }\quad q<0 .\end{cases} \\
      \eea
However, given the fact that in the expression      
  \beaa
H^{\mu\nu}=-h^{\mu\nu}+ O^{\mu\nu}(h^2) ,
\eeaa
here the $O^{\mu\nu}(h^2)$ happen to be a product of tensors of $m$ with $h^2$, we then also have that
  \bea
\derm_\a H^{\mu\nu}=- \derm_\a h^{\mu\nu}+ O_\a^{\, \,\,\, \mu\nu} (h \cdot \derm h ). 
\eea      

Since for all $|I|$,
               \bea\label{apriorestimimateonproducthandgradientofh}
 \notag
&& |   \Lie_{Z^I} h (t,x) |  \cdot |\derm  ( \Lie_{Z^I} h ) (t,x)  |  \\
 \notag
  &\leq& \begin{cases} c (\delta) \cdot  c (\gamma) \cdot C ( |I| ) \cdot E ( |I| + 2)  \cdot \frac{\eps^2 }{(1+t+|q|)^{2-2\delta} (1+|q|)^{1+2\de}},\quad\text{when }\quad q>0,\\
       C ( |I| ) \cdot E ( |I| + 2)  \cdot \frac{\eps^2  }{(1+t+|q|)^{2-2\delta} }  \,\quad\text{when }\quad q<0 . \end{cases} \\
      \eea
      
we obtain,
              \beaa
 \notag
 && |\derm ( \Lie_{Z^I}  h ) (t,x)  | +  |   \Lie_{Z^I} h (t,x) |  \cdot |\derm  ( \Lie_{Z^I} h ) (t,x)  |  \\
  &\leq& \begin{cases} c (\delta) \cdot  c (\gamma) \cdot C ( |I| ) \cdot E ( |I| + 2)  \cdot \big( \frac{\eps }{(1+t+|q|)^{1-\delta} (1+|q|)^{1+\de}} + \frac{\eps^2 }{(1+t+|q|)^{2-2\delta} (1+|q|)^{1+2\de}}\big),\quad\text{when }\quad q>0,\\
       C ( |I| ) \cdot E ( |I| + 2)  \cdot \big( \frac{\eps  }{(1+t+|q|)^{1-\delta}(1+|q|)^{\frac{1}{2} }} +  \frac{\eps^2  }{(1+t+|q|)^{2-2\delta} } \big) \,\quad\text{when }\quad q<0 , \end{cases} \\
       \eeaa
       Thus, if $\delta \leq \frac{1}{2}$ and $\eps \leq 1$, we get 
       \beaa
       \notag
 && |\derm  ( \Lie_{Z^I} h ) (t,x)  | +  |   \Lie_{Z^I} h (t,x) |  \cdot |\derm  ( \Lie_{Z^I} h ) (t,x)  |  \\
 \notag
         &\leq& \begin{cases} c (\delta) \cdot  c (\gamma) \cdot C ( |I| ) \cdot E ( |I| + 2)  \cdot \frac{\eps }{(1+t+|q|)^{1-\delta} (1+|q|)^{1+\de}} ,\quad\text{when }\quad q>0,\\
       C ( |I| ) \cdot E ( |I| + 2)  \cdot  \frac{\eps  }{(1+t+|q|)^{1-\delta}(1+|q|)^{\frac{1}{2} }}   \,\quad\text{when }\quad q<0 , \end{cases} \\
      \eeaa
which gives the result for  $|\derm   H (t,x)  |$.

\end{proof}

      \begin{lemma}\label{aprioriestimatesonZLiederivativesonbigH}
      Let $M\leq \eps$. By choosing $\gamma > \delta  $, we have for all $|I|$, $\delta \leq \frac{1}{2}$,\,$\eps \leq 1$,

              \beaa
 \notag
 |\derm ( \Lie_{Z^I} H ) (t,x)  |    &\leq& \begin{cases} C ( |I| ) \cdot E ( |I| + 2)  \cdot \frac{\eps }{(1+t+|q|)^{1-\delta} (1+|q|)^{1+\de}} \;,\quad\text{when }\quad q>0 \;,\\
       C ( |I| ) \cdot E ( |I| + 2)  \cdot \frac{\eps  }{(1+t+|q|)^{1-\delta} (1+|q|)^{\frac{1}{2} }}  \; , \quad\text{when }\quad q<0 \; , \end{cases} \\
      \eeaa
      and
       \beaa
 \notag
|  \Lie_{Z^I} H (t,x)  | &\leq& \begin{cases} c (\delta) \cdot c (\gamma) \cdot C ( |I| ) \cdot E ( |I| + 2) \cdot  \frac{\eps}{ (1+ t + | q | )^{1-\delta }  (1+| q |   )^{\de}} \;  ,\quad\text{when }\quad q>0 \;  ,\\
\notag
    C ( |I| ) \cdot E ( |I| + 2) \cdot  \frac{\eps}{ (1+ t + | q | )^{1-\delta }  } (1+| q |   )^{\frac{1}{2} }  \;  , \,\quad\text{when }\quad q<0  \;  . \end{cases} \\ 
    \eeaa
Consequently,

   \beaa
 \notag |  \rderm ( \Lie_{Z^I} H ) (t,x)  | 
 &\leq& \begin{cases} c (\delta) \cdot c (\gamma) \cdot C ( |I| ) \cdot E ( |I| + 3)  \cdot \frac{\eps }{(1+t+|q|)^{2-\delta} (1+|q|)^{\de}} \;  ,\quad\text{when }\quad q>0 \;  ,\\
       C ( |I| ) \cdot E ( |I| + 3)  \cdot \eps \cdot \frac{ (1+|q|)^{\frac{1}{2} }}{ (1+t+|q|)^{2-\delta} } \; , \quad\text{when }\quad q<0 \;  . \end{cases} \\
      \eeaa
      
\end{lemma}

\begin{proof}
We have already showed in \cite{G4} that 
  \beaa
H^{\mu\nu}=-h^{\mu\nu}+ O^{\mu\nu}(h^2) .
\eeaa

Using again that here that $O^{\mu\nu}(h^2)$ and $O_\a^{\, \,\,\, \mu\nu} (h \cdot \pa h )$ are in fact product of constant tensors with $h$ and $\derm h$ and using the Leibniz rule for Lie derivatives, we obtain that for all $Z \in \cal Z$, 
  \beaa
 \Lie_{Z}  H^{\mu\nu} &=& -  \Lie_{Z}  h^{\mu\nu}+ O^{\mu\nu}( h\cdot  \Lie_{Z} h ) \\
 \derm_\a  ( \Lie_{Z}    H)^{\mu\nu} &=& -  \derm_\a  ( \Lie_{Z}  h)^{\mu\nu}+ O_\a^{\, \,\,\, \mu\nu} (  \derm  h \cdot  \Lie_{Z}     h ) + O_\a^{\, \,\,\, \mu\nu} (    h \cdot   \derm (\Lie_{Z}     h)  ). 
\eeaa
Since $|h|$ and  $ | \Lie_{Z^I} h|$ obey the same estimate, and since $|\derm h|$ and  $| \derm ( \Lie_{Z^I} h ) |$ obey the same estimate, we then derive the same estimate for $|\Lie_{Z} H^{\mu\nu}|$ as for  $| H^{\mu\nu}|$ and the same estimate for $|\derm ( \Lie_{Z} H )^{\mu\nu}|$ as for  $| \derm H^{\mu\nu}|$. By induction , we get the result for all $|I|$.

Now, since we have already controlled the $| \Lie_{Z^I} H |$, we could use
\bea
\notag
(1 + t + |q| ) \cdot |\rderm ( \Lie_{Z^I} H ) | \les \sum_{|J| \leq |I| + 1}  | \Lie_{Z^J} H |    
\eea
to derive the estimate for $|\rderm ( \Lie_{Z^I} H ) |$.

\end{proof}

\section{Estimates on the zero-covariant derivative of the special components of the Einstein-Yang-Mills fields}

\begin{lemma}\label{estimateonpartialderivativeofALcomponent}
Under the bootstrap assumption \eqref{bootstrap} holding for all $|J| \leq  \lfloor \frac{|I|}{2} \rfloor  $\;, we have the following estimate for the Einstein-Yang-Mills potential, 
\bea\label{partialderivativeofLeiderivativeofALcomponent}
\notag
  && | \pa  \Lie_{Z^I}  A_{L}  |  \\
  \notag
  &\les& E(  \lfloor \frac{|I|}{2} \rfloor   )  \cdot \Big(  \sum_{|J|\leq |I|} |   \rderm  \Lie_{Z^J} A | + \sum_{|J|\leq |I| -1} |   \derm  \Lie_{Z^J} A |  + \sum_{|K| + |M| \leq |I|} O \big(  | ( \Lie_{Z^K}  h )  | \cdot  | \derm  (  \Lie_{Z^M}  A ) | \big) \Big)  \; ,\\
\eea
where 
\beaa
| \pa  \Lie_{Z^I}  A_{L} |^2 = | \pa_t ( \Lie_{Z^I} A )_{L} |^2 + \sum_{i=1}^3 | \pa_i ( \Lie_{Z^I} A)_{L} |^2 \;. 
\eeaa
Under the bootstrap assumption \eqref{bootstrap} holding for all $|J| \leq  \lfloor \frac{|I|}{2} \rfloor  $\;, we have the following estimate for the Einstein-Yang-Mills metric, 
\bea\label{partialderivativeofLiederivativeofHTangentialLcompo}
  \notag
  && | \pa  \Lie_{Z^I}  H_{\cal T L}  | \\
    \notag
   &\les& E(  \lfloor \frac{|I|}{2} \rfloor   )  \cdot \Big(   \sum_{|J|\leq |I|} |   \rderm  \Lie_{Z^J} H | + \sum_{|J|\leq |I| -1} |   \derm  \Lie_{Z^J} H |  + \sum_{|K|+ |M| \leq |I|}  O (|\Lie_{Z^K} H| \cdot |\derm ( \Lie_{Z^M} H )| )   \Big)  \; ,\\
\eea
\bea\label{partialderivativeofLiederivativeofHLLcomp}
  \notag
  && | \pa  \Lie_{Z^I}  H_{ L L}  |  \\
    \notag
    &\les& E(  \lfloor \frac{|I|}{2} \rfloor   )  \cdot \Big(  \sum_{|J|\leq |I|} |   \rderm  \Lie_{Z^J} H | + \sum_{|J|\leq |I| -2} |   \derm  \Lie_{Z^J} H |  + \sum_{|K|+ |M| \leq |I|}  O (|\Lie_{Z^K} H| \cdot |\derm ( \Lie_{Z^M} H )| )   \Big)  \; ,\\
\eea
and 
\bea\label{partialderivativeofLiederivativeofsmallhTangentialLcompo}
  \notag
  && | \pa  \Lie_{Z^I}  h_{\cal T L}  | \\
    \notag
   &\les& E(  \lfloor \frac{|I|}{2} \rfloor   )  \cdot \Big(  \sum_{|J|\leq |I|} |   \rderm  \Lie_{Z^J} h | + \sum_{|J|\leq |I| -1} |   \derm  \Lie_{Z^J} h |  + \sum_{|K|+ |M| \leq |I|}  O (|\Lie_{Z^K} h| \cdot |\derm ( \Lie_{Z^M} h )| )   \Big)  \; ,\\
\eea
\bea\label{partialderivativeofLiederivativeofsmallhLLcomp}
  \notag
  && | \pa  \Lie_{Z^I}  h_{ L L}  |  \\
  \notag
  &\les& E(  \lfloor \frac{|I|}{2} \rfloor   )  \cdot \Big(  \sum_{|J|\leq |I|} |   \rderm  \Lie_{Z^J} h | + \sum_{|J|\leq |I| -2} |   \derm  \Lie_{Z^J} h |  + \sum_{|K|+ |M| \leq |I|}  O (|\Lie_{Z^K} h| \cdot |\derm ( \Lie_{Z^M} h )| )   \Big)  \; ,\\
\eea
where 
\beaa
| \pa  \Lie_{Z^I}  H_{\cal T L}  |^2 = | \pa_t ( \Lie_{Z^I} H  )_{\cal T L} |^2 + \sum_{i=1}^3 | \pa_i ( \Lie_{Z^I} H  )_{\cal T L} |^2 \;. 
\eeaa
In these estimates the partial differentiation actually concerns a partial derivative of the stated components of the Lie derivatives of the tensor (as opposed to the covariant derivative), and where the $\sum_{|J|\leq  -1}$ is understood to be a vanishing sum.

\end{lemma}

\begin{proof}
Based on the estimate derived by Lindblad-Rodnianski in \cite{LR10}, we have that if for $|Z^J H|\le C$ for all $|J| \leq  \lfloor \frac{|I|}{2} \rfloor  $, which is a condition satisfied by our bootstrap assumption, we have
\beaa
&&|\pa Z^I H_{ L\cal T} | \\
&\les&  E(  \lfloor \frac{|I|}{2} \rfloor   )  \cdot  \Big (\sum_{|J|\leq |I|}|\rpa Z^J
H|+\!\!\!\!\sum_{|J|\leq |I|-1}\!\!\!|\pa Z^J H|\,\, +
 \sum_{\,\,\,\,|K|+|M| \leq |I|}|Z^{K} H| \cdot |\pa Z^{M} H|\Big )\; ,
\eeaa
\beaa
&&|\pa Z^I H_{ L L} |  \\
&\les&  E(  \lfloor \frac{|I|}{2} \rfloor   )  \cdot  \Big( \sum_{|J|\leq |I|}|\rpa Z^J
H|+\!\!\!\!\sum_{|J|\leq |I|-2}\!\!\!|\pa Z^J H|\,\, +
 \sum_{\,\,\,\,|K|+|M| \leq |I|}|Z^{K} H| \cdot |\pa Z^{M} H| \Big)\; .
\eeaa

However,  the commutation of two vector fields in $\cal Z$ is a linear combination of vector fields in $\cal Z$, and that a commutation of a vector field in $\cal Z$ and of a vector $\pa_\mu$, $\mu \in  \{t, x^1, x^2, x^3 \}$, gives a linear combination of vectors of the form $\pa_\mu$. Using that fact, and the fact that we are summing over all $Z \in {\cal Z}$, we get
\beaa
&&\sum_{|J|\leq |I|}|\rpa Z^J
H|+\!\!\!\!\sum_{|J|\leq |I|-1}\!\!\!|\pa \Lie_{Z^J} H |\,\, + \sum_{\,\,\,\,|I_1|+||I_2| \leq |I|}|\Lie_{Z^{I_{2}}} H| \cdot |\derm \Lie_{Z^{I_1}} H| \\
&\les& \sum_{|J|\leq |I|}|\rderm \Lie_{Z^J} H|+\!\!\!\!\sum_{|J|\leq |I|-1}\!\!\!|\derm \Lie_{Z^J} H|\,\, + \sum_{\,\,\,\,|K|+|M| \leq |I|}| \Lie_{Z^{K}} H| \cdot |\derm \Lie_{Z^{M}} H| \; . \\
\eeaa
Using the fact that the Minkowski covariant derivative commutes with the Lie derivative in direction of Minkowski vector fields, we obtain the result for \eqref{partialderivativeofLiederivativeofHTangentialLcompo} and \eqref{partialderivativeofLiederivativeofHLLcomp}.

Based on the equations that we showed in the proof of Lemma \label{LorenzgaugeestimateforLiederivatives}, where we established the Lorenz gauge estimate for the Lie derivatives of the Einstein-Yang-Mils potential, and adapting exactly the same argument of Lindblad-Rodnianski as in Appendix D of \cite{LR10} for the Einstein vacuum metric, we obtain
\beaa
&&  | \pa  \Lie_{Z^I}  A_{L}  |  \\
  &\les& E(  \lfloor \frac{|I|}{2} \rfloor   )  \cdot  \Big(  \sum_{|J|\leq |I|} |   \rpa Z^J A | + \sum_{|J|\leq |I| -1} |   \pa  Z^J A |  + \sum_{|K| + |M| \leq |J|} O \big(  | ( Z^K  h )  | \cdot  | \pa  (  Z^M  A ) | \big) \Big) \; ,
\eeaa
and then similarly, using the special structure of the Minkowski vector fields, we obtain the stated result in \eqref{partialderivativeofLeiderivativeofALcomponent}.

\end{proof}

\begin{definition}\label{definitionforintegralalongthenullcoordinateplusboundaryterm}
For a function $f$, we define $\int\limits_{s,\, \Om=const} f  $ as the integral at a fixed $\Om \in \SSS^2$, from $(t, | x | \cdot \Om)$ along the line $(\tau, r \cdot \Om)$ such that $r+\tau =  | x | +t $ (i.e. along a fixed null coordinate $s:=\tau+r $) till we reach the hyperplane $\tau=0$, to which we also add the generated boundary term at the hyperplane prescribed by $\tau=0$. In other words,
\bea\label{definitionequationforintegralalongthenullcoordinateplusboundaryterm}
\notag
\int\limits_{s,\, \Om=const} |f |  (t, | x | \cdot \Om)  &=&   \int_{t+| x | }^{| x |  } \pa_r | f (t + | x | - r,  r  \cdot \Om ) | dr +  | f \big(0, ( t + | x |) \cdot \Om \big) |  \; .\\
\eea
\end{definition}

\begin{remark}\label{remarkabouttheasymptoticbehaviouroftheboundarytermonhyperplaneprescribedtqual0}
Under the bootstrap assumption, we showed in \cite{G4} that the fields $ \Lie_{Z^I}  A (0,x)$\, $\Lie_{Z^I}  h^1 (0,x)$ and $ \derm ( \Lie_{Z^I}  A )  (0,x)$\, $\derm (\Lie_{Z^I}  h^1 ) (0,x)$\,, have a certain asymptotic behaviour on the hyperplane prescribed by $t=0$, which will be always used to estimate the boundary term in \eqref{definitionequationforintegralalongthenullcoordinateplusboundaryterm}, either for $|  \Lie_{Z^I}  A \big(0, ( t + | x |) \cdot \Om \big) |$ and $| \Lie_{Z^I}  h^1 \big(0, ( t + | x |) \cdot \Om \big) |$\, or for $| \derm ( \Lie_{Z^I}  A ) \big(0, ( t + | x |) \cdot \Om \big) |$ and $| \derm ( \Lie_{Z^I}   h^1 ) \big(0, ( t + | x |) \cdot \Om \big) |$\,.
\end{remark}

\begin{lemma}\label{estimategoodcomponentspotentialandmetric}
Let $M\leq \eps$\,. Under the bootstrap assumption holding for all $|J| \leq  \lfloor \frac{|I|}{2} \rfloor  $\,, we have
        \bea
        \notag
 |   \Lie_{Z^I}  A_{L}  |  (t, | x | \cdot \Om)  &\leq& \int\limits_{s,\,\Om=const} \sum_{|J|\leq |I| -1} |   \derm (  \Lie_{Z^J} A )  | \\
 \notag
     && + \begin{cases} c (\delta) \cdot c (\gamma) \cdot C ( |I| ) \cdot E ( |I| + 3)  \cdot \big( \frac{ \eps   }{ (1+t+|q|)^{2-2\delta} \cdot  (1+|q|)^{\ga - 1} } \big),\quad\text{when }\quad q>0 \; ,\\
       C ( |I| ) \cdot E ( |I| + 3)  \cdot \big( \frac{ \eps \cdot (1+|q|)^{\frac{3}{2} }   }{ (1+t+|q|)^{2-2\delta} } \big) \,\quad\text{when }\quad q<0 \; , \end{cases} 
       \eea
        where $\int\limits_{s,\, \Om=const}  |   \derm  ( \Lie_{Z^J} A )  |   $ is the integral defined as in Definition \ref{definitionforintegralalongthenullcoordinateplusboundaryterm}.\\

Also, under the bootstrap assumption holding for all $|J| \leq  \lfloor \frac{|I|}{2} \rfloor  $,
        \beaa
        \notag
|  \Lie_{Z^I} h_{\cal T L} |  + |  \Lie_{Z^I} H_{\cal T L} |   &\les&  \int\limits_{s,\,\Om=const} \sum_{|J|\leq |I| -1} \big( |\derm  ( \Lie_{Z^J} h) | + |\derm  ( \Lie_{Z^J} H )  | \big) \\
&& + \begin{cases} c (\delta) \cdot c (\gamma) \cdot C ( |I| ) \cdot E ( |I| + 3)  \cdot \frac{ \eps   }{ (1+t+|q|) } ,\quad\text{when }\quad q>0 \; ,\\
        C ( |I| ) \cdot E ( |I| + 3) \cdot  \frac{ \eps \cdot (1+|q|)^{\frac{1}{2} + 2 \delta } }{ (1+t+|q|) }  \,\quad\text{when }\quad q<0 \; . \end{cases} \\
       \eeaa
   and
           \beaa
        \notag
|  \Lie_{Z^I} h_{L L} |  + |  \Lie_{Z^I} H_{ L L} |   &\les&  \int\limits_{s,\,\Om=const} \sum_{|J|\leq |I| -2} \big( |\derm  ( \Lie_{Z^J} h ) | + |\derm  ( \Lie_{Z^J} H ) | \big) \\
&& + \begin{cases} c (\delta) \cdot c (\gamma) \cdot C ( |I| ) \cdot E ( |I| + 3)  \cdot \frac{ \eps   }{ (1+t+|q|) } ,\quad\text{when }\quad q>0 \; ,\\
        C ( |I| ) \cdot E ( |I| + 3) \cdot  \frac{ \eps \cdot (1+|q|)^{\frac{1}{2} + 2 \delta } }{ (1+t+|q|) }  \,\quad\text{when }\quad q<0 \; . \end{cases} \\
       \eeaa    
  
             As a result, more precisely, in $\overline{C} \subset \{ q \geq q_0 \} $, we have 
                    \bea
        \notag
|  \Lie_{Z^I} h_{\cal T L} |    + |  \Lie_{Z^I} H_{\cal T L} |                    &\les&   \int\limits_{s,\,\Om=const} \sum_{|J|\leq |I| -1} \big( |\derm  ( \Lie_{Z^J} h) | + |\derm (  \Lie_{Z^J} H  ) | \big) \\
\notag
&& + C(q_0) \cdot   c (\delta) \cdot c (\gamma) \cdot C ( |I| ) \cdot E ( |I| + 3)  \cdot \frac{ \eps   }{ (1+t+|q|) } \; , \\
       \eea
and
                \bea
        \notag
|  \Lie_{Z^I} h_{L L} |    + |  \Lie_{Z^I} H_{L L} |                    &\les&   \int\limits_{s,\,\Om=const} \sum_{|J|\leq |I| -2} \big( |\derm  ( \Lie_{Z^J} h) | + |\derm  ( \Lie_{Z^J} H )  | \big) \\
\notag
&& + C(q_0) \cdot   c (\delta) \cdot c (\gamma) \cdot C ( |I| ) \cdot E ( |I| + 3)  \cdot \frac{ \eps   }{ (1+t+|q|) } \; . \\
       \eea

\end{lemma}

            \begin{proof}
            
            In the Lorenz gauge, we showed in Lemma \ref{estimateonpartialderivativeofALcomponent}, that

\beaa
  | \pa  \Lie_{Z^I}  A_{L}  |  \les \sum_{|J|\leq |I|} |   \rderm  ( \Lie_{Z^J} A ) | + \sum_{|J|\leq |I| -1} |   \derm  ( \Lie_{Z^J} A ) |  + \sum_{|K| + |M| \leq |I|} O \big(  | ( \Lie_{Z^K}  h )  | \cdot  | \derm  (  \Lie_{Z^M}  A ) | \big) \; ,
\eeaa

For $\gamma >\delta  $, we have
   \beaa
 \notag
&& |  \rderm ( \Lie_{Z^I} A ) (t,x)  |   \\
\notag
 &\leq& \begin{cases} c (\delta) \cdot c (\gamma) \cdot C ( |I| ) \cdot E ( |I| + 3)  \cdot \frac{\eps }{(1+t+|q|)^{2-\delta} (1+|q|)^{\ga}},\quad\text{when }\quad q>0 \; ,\\
       C ( |I| ) \cdot E ( |I| + 3)  \cdot \eps \cdot \frac{ (1+|q|)^{\frac{1}{2} }}{ (1+t+|q|)^{2-\delta} } \,\quad\text{when }\quad q<0 \;, \end{cases} \\
      \eeaa
and for all $|K| + |M| \leq |I|$, we have
 \bea
 \notag
|   \Lie_{Z^K} h (t,x)  | &\leq& \begin{cases} c (\delta) \cdot c (\gamma) \cdot C ( |I| ) \cdot E ( |I| + 2) \cdot  \frac{\eps}{ (1+ t + | q | )^{1-\delta }  (1+| q |   )^{\de}}  ,\quad\text{when }\quad q>0 \; ,\\
\notag
    C ( |I| ) \cdot E ( |I| + 2) \cdot  \frac{\eps}{ (1+ t + | q | )^{1-\delta }  } (1+| q |   )^{\frac{1}{2} }  , \,\quad\text{when }\quad q<0 \; , \end{cases} \\ 
    \eea
    and
                 \beaa
 \notag
|\derm ( \Lie_{Z^M} A ) (t,x)  |     &\leq& \begin{cases} C ( |I| ) \cdot E ( |I| + 2)  \cdot \frac{\eps }{(1+t+|q|)^{1-\delta} (1+|q|)^{1+\ga}},\quad\text{when }\quad q>0 \; ,\\
       C ( |I| ) \cdot E ( |I| + 2)  \cdot \frac{\eps  }{(1+t+|q|)^{1-\delta}(1+|q|)^{\frac{1}{2} }}  \,\quad\text{when }\quad q<0 \; . \end{cases} \\
      \eeaa
Thus,
                 \bea
 \notag
&& \sum_{|K| + |M| \leq |I|}  |   \Lie_{Z^K} h (t,x) |  \cdot |\derm (  \Lie_{Z^M} A ) (t,x)  |   \\
 \notag
  &\leq& \begin{cases} c (\delta) \cdot  c (\gamma) \cdot C ( |I| ) \cdot E ( |I| + 2)  \cdot \frac{\eps^2 }{(1+t+|q|)^{2-2\delta} (1+|q|)^{1+\ga+\de}},\quad\text{when }\quad q>0 \; ,\\
       C ( |I| ) \cdot E ( |I| + 2)  \cdot \frac{\eps^2  }{(1+t+|q|)^{2-2\delta} }  \,\quad\text{when }\quad q<0 \; . \end{cases} \\
      \eea

      As a result,
      \beaa
&& | \pa ( \Lie_{Z^I}  A_{L} )  |  \\
&\les&  \sum_{|J|\leq |I| -1} |   \derm  ( \Lie_{Z^J} A )  |  \\
&& + \begin{cases} c (\delta) \cdot c (\gamma) \cdot C ( |I| ) \cdot E ( |I| + 3) \cdot \big( \frac{\eps }{(1+t+|q|)^{2-\delta} (1+|q|)^{\ga}} +  \frac{\eps^2 }{(1+t+|q|)^{2-2\delta} (1+|q|)^{1+\ga+\de}}\big),\quad\text{when }\quad q>0 \; ,\\
       C ( |I| ) \cdot E ( |I| + 3)  \cdot \big( \frac{ \eps \cdot (1+|q|)^{\frac{1}{2} }}{ (1+t+|q|)^{2-\delta} } + \frac{\eps^2  }{(1+t+|q|)^{2-2\delta} }  \big) \,\quad\text{when }\quad q<0 \;, \end{cases} 
\eeaa

      We get for $\eps < 1$, 
        \beaa
&& | \pa  ( \Lie_{Z^I}  A_{L} )  |  \\
&\les&  \sum_{|J|\leq |I| -1} |   \derm   ( \Lie_{Z^J} A )  |  \\
       && + \begin{cases} c (\delta) \cdot c (\gamma) \cdot C ( |I| ) \cdot E ( |I| + 3)  \cdot \big( \frac{\eps }{(1+t+|q|)^{2-2\delta} (1+|q|)^{\ga}}  \big),\quad\text{when }\quad q>0 \; ,\\
        C ( |I| ) \cdot E ( |I| + 3)  \cdot \big( \frac{ \eps \cdot (1+|q|)^{\frac{1}{2} }}{ (1+t+|q|)^{2-2\delta} } \big) \,\quad\text{when }\quad q<0 \; . \end{cases}
\eeaa

Since $\Sigma_{t=0}$ is diffeomorphic to $\R^3$, for each $x \in \Sigma_{0}$, there exists $\Om \in \SSS^3$, such that $x =r \cdot \Om $. As previously, applying the fundamental theorem of calculus by integrating at a fixed $\Om$, from $(t, | x | \cdot \Om)\in \overline{C}$, along the line $(\tau, r \cdot \Om)$ such that $r+\tau =  | x | +t :=s $ till we reach the hyperplane $\tau=0$, and considering that such a line will be totally contained in $\overline{C}$, we obtain, as before, that

        \beaa
 |   \Lie_{Z^I} A_{L}  |  &\les&\int\limits_{s,\Om=const}  \sum_{|J|\leq |I| -1} |   \derm  ( \Lie_{Z^J} A )  |     \\
 && + \begin{cases} c (\delta) \cdot c (\gamma) \cdot C ( |I| ) \cdot E ( |I| + 3)  \cdot \big( \frac{\eps \cdot (1+|q|) }{(1+t+|q|)^{2-2\delta} (1+|q|)^{\ga}}  \big),\quad\text{when }\quad q>0 \; ,\\
 C ( |I| ) \cdot E ( |I| + 3)  \cdot \big( \frac{ \eps \cdot (1+|q|)^{ \frac{3}{2} }}{ (1+t+|q|)^{2-2\delta} } \big) \,\quad\text{when }\quad q<0 \; , \end{cases} 
 \eeaa
  where the sum over $|I|-1$ is absent if $|I| = 0$\,.  \\

Concerning $H$, we have shown in Lemma \ref{estimateonpartialderivativeofALcomponent}, that
\beaa
  | \pa  \Lie_{Z^I}  H_{\cal T L}  |  \les \sum_{|J|\leq |I|} |   \rderm  \Lie_{Z^J} H | + \sum_{|J|\leq |I| -1} |   \derm  \Lie_{Z^J} H |  + \sum_{|K|+ |M| \leq |I|}  O (|\Lie_{Z^K} H| \cdot |\derm ( \Lie_{Z^M} H )| ) \; ,
\eeaa

Integrating as before, we get

        \beaa
 |  \Lie_{Z^I} H_{\cal T L}  |  &\les&    \int\limits_{s,\, \Om=const} \sum_{|J|\leq |I| -1} |\derm \Lie_{Z^J} H| \\
&& + \begin{cases} c (\delta) \cdot c (\gamma) \cdot C ( |I| ) \cdot E ( |I| + 3) \cdot \big( \frac{\eps \cdot (1+|q|) }{(1+t+|q|)^{2-2\delta} (1+|q|)^{2\delta}}  \big),\quad\text{when }\quad q>0  \; ,\\
       C ( |I| ) \cdot E ( |I| + 3)  \cdot \big( \frac{ \eps \cdot (1+|q|)^{ \frac{3}{2} }}{ (1+t+|q|)^{2-2\delta} } \big) \,\quad\text{when }\quad q<0 \; . \end{cases} \\
       &\les& \int\limits_{s,\,\Om=const} \sum_{|J|\leq |I| -1} |\derm \Lie_{Z^J} H | \\ 
       && + \begin{cases} c (\delta) \cdot c (\gamma) \cdot  C ( |I| ) \cdot E ( |I| + 3)  \cdot \big( \frac{ \eps  \cdot (1+|q|)^{ 1- 2\delta } \cdot  (1+|q|)^{ 2\delta } }{ (1+t+|q|) \cdot (1+t+|q|)^{1-2\delta} \cdot (1+|q|)^{2\delta} } \big),\quad\text{when }\quad q>0  \; ,\\
        C ( |I| ) \cdot E ( |I| + 3) \cdot \big( \frac{ \eps \cdot (1+|q|)^{\frac{1}{2} } \cdot (1+|q|)^{ 1- 2\delta } \cdot  (1+|q|)^{ 2\delta } }{ (1+t+|q|) \cdot (1+t+|q|)^{1-2\delta} } \big) \,\quad\text{when }\quad q<0 \; . \end{cases} \\
            &\les& \int\limits_{s,\,\Om=const} \sum_{|J|\leq |I| -1} |\derm \Lie_{Z^K} H | \\ 
            && + \begin{cases} c (\delta) \cdot c (\gamma) \cdot C ( |I| ) \cdot E ( |I| + 3)  \cdot \big( \frac{ \eps  \cdot (1+|q|)^{ 1- 2\delta } }{ (1+t+|q|) \cdot (1+t+|q|)^{1-2\delta} } \big),\quad\text{when }\quad q>0  \; ,\\
      C ( |I| ) \cdot E ( |I| + 3)  \cdot \big( \frac{ \eps \cdot (1+|q|)^{\frac{1}{2} + 2 \delta } \cdot (1+|q|)^{ 1- 2\delta }  }{ (1+t+|q|) \cdot (1+t+|q|)^{1-2\delta} } \big) \,\quad\text{when }\quad q<0  \; . \end{cases} 
       \eeaa
Hence,
        \beaa
        \notag
 |   \Lie_{Z^I} H_{\cal T L}  |   &\leq& \int\limits_{s,\, \Om=const} \sum_{|J|\leq |I| -1} |\derm \Lie_{Z^J} H | \\ 
    \notag
 && + \begin{cases} c (\delta) \cdot c (\gamma) \cdot C ( |I| ) \cdot E ( |I| + 3)  \cdot \frac{ \eps   }{ (1+t+|q|) } ,\quad\text{when }\quad q>0  \; ,\\
        C ( |I| ) \cdot E ( |I| + 3)  \cdot  \frac{ \eps \cdot (1+|q|)^{\frac{1}{2} + 2 \delta } }{ (1+t+|q|) }  \,\quad\text{when }\quad q<0 \; . \end{cases} \\
       \eeaa
             
Similarly,
        \beaa
        \notag
 |   \Lie_{Z^I} H_{L L}  |   &\leq& \int\limits_{s,\,\Om=const} \sum_{|J|\leq |I| -2} |\derm \Lie_{Z^J} H | \\ 
       \notag
 && + \begin{cases} c (\delta) \cdot c (\gamma) \cdot C ( |I| ) \cdot E ( |I| + 3)  \cdot \frac{ \eps   }{ (1+t+|q|) } ,\quad\text{when }\quad q>0  \; ,\\
        C ( |I| ) \cdot E ( |I| + 3)  \cdot  \frac{ \eps \cdot (1+|q|)^{\frac{1}{2} + 2 \delta } }{ (1+t+|q|) }  \,\quad\text{when }\quad q<0 \; . \end{cases} \\
       \eeaa
       The same estimates hold for $h$ as well.

            \end{proof}
            
\section{Improved pointwise estimates on the gradient of “good” components of the Einstein-Yang-Mills fields}

We going first start by verifying the assumptions for \eqref{LinfinitynormestimateongradientderivedbyLondbladRodnianski}, so that we could apply it.
\begin{lemma}

Let $M\leq \eps$\,. With our conditions on the initial data, generated from the bootstrap assumption, we have
\bea
\notag
 |H|\leq \frac{\varepsilon^\prime}{4}\;,\qquad\int_0^{\infty}\!\!
\|\,H(t,\cdot)\|_{L^\infty(\overline{D}_t)} \cdot \frac{dt}{(1+t)} \leq
\frac{\varepsilon^\prime}{4} \;,\qquad |H_{L{\cal
T}} | \leq\frac{\varepsilon^\prime}{4}\cdot \frac{(|q|+1)}{(1+t+|x|)} \; ,\\
\eea
with
\bea\label{definitionofepsprime}
\varepsilon^\prime = 4 C  \cdot  c (\gamma)  \cdot c (\delta)  \cdot E ( 3)  \cdot \eps \; ,
\eea
and where  $C$ is some constant.

\end{lemma}

\begin{proof}

  We have already showed in Lemma \ref{aprioriestimatesonZLiederivativesonbigH}, that
        \beaa
 \notag
&& |   \Lie_{Z^I} H (t,x)  | \\
      &\leq& \begin{cases} c (\delta) \cdot c (\gamma) \cdot C ( |I| ) \cdot E ( |I| + 2) \cdot  \frac{\eps}{ (1+ t + | q | )^{1-\delta }  (1+| q |   )^{\de}}  ,\quad\text{when }\quad q>0,\\
\notag
    C ( |I| ) \cdot E ( |I| + 2) \cdot \frac{\eps}{ (1+ t + | q | )^{1-\delta }  } (1+| q |   )^{\frac{1}{2} }  , \,\quad\text{when }\quad q<0 .\end{cases} \\
      \eeaa
      
Consequently,
      
 \beaa
 \notag
|   \Lie_{Z^I} H (t,x)  |   &\leq&  c (\gamma)  \cdot c (\delta) \cdot C ( |I| ) \cdot E ( |I| + 2)  \cdot \frac{\eps}{(1+t+|q|)^{\frac{1}{2}-\delta}  } \; .
      \eeaa

      Hence,
     \bea
|   H (t,x)  |   &\les& c (\gamma)  \cdot c (\delta) \cdot E (  2) \cdot \eps \, . 
      \eea
We also have

   \beaa
 \notag
  \frac{|   H (t,x)  | }{1+t}  &\les& c (\gamma)  \cdot c (\delta) \cdot E (  2)  \cdot \frac{\eps}{(1+t)^{\frac{3}{2}-\delta}  }   ,
      \eeaa
     and therefore  

   \beaa
 \notag
\int_0^{\infty}  \frac{|   H (t,x)  | }{1+t} dt &\les&  c (\gamma)  \cdot c (\delta) \cdot E (  2)  [ \frac{\eps}{(1+t)^{\frac{1}{2}-\delta}  } ]_{\infty}^{0}  \\
   &\les& c (\gamma)  \cdot c (\delta)  \cdot E ( 2) \cdot \eps  
      \eeaa
 and thus,
 
\bea
\int_0^{\infty}\!\!  |\,H(t, x) | \frac{dt}{1+t} \les c (\gamma)  \cdot c (\delta) \cdot E ( 2) \cdot \eps   \, .
\eea

We also have,

 \bea
 \notag
|   \Lie_{Z^I} h (t,x)  | &\leq& \begin{cases} c (\delta) \cdot c (\gamma) \cdot C ( |I| ) \cdot E ( |I| + 2) \cdot  \frac{\eps}{ (1+ t + | q | )^{1-\delta }  (1+| q |   )^{\de}}  ,\quad\text{when }\quad q>0,\\
\notag
    C ( |I| ) \cdot E ( |I| + 2) \cdot  \frac{\eps}{ (1+ t + | q | )^{1-\delta }  } (1+| q |   )^{\frac{1}{2} }  , \,\quad\text{when }\quad q<0 , \end{cases} \\ 
    \notag
      &\leq&  c (\gamma)  \cdot c (\delta) \cdot C ( |I| )  \cdot   E ( |I| + 2) \frac{\eps}{(1+t+|q|)^{\frac{1}{2}-\delta}  } . \\
      \eea

We showed in Lemma \eqref{estimategoodcomponentspotentialandmetric}, that

       \beaa
        \notag
 |   H_{\cal T L} |   &\les& \begin{cases} c (\delta) \cdot c (\gamma) \cdot E (  3)  \cdot \frac{ \eps   }{ (1+t+|q|) } ,\quad\text{when }\quad q>0,\\
      E ( 3)  \cdot  \frac{ \eps \cdot (1+|q|)^{\frac{1}{2} + 2 \delta } }{ (1+t+|q|) }  \,\quad\text{when }\quad q<0 . \end{cases} 
       \eeaa
Thus,

  \bea
|   H_{\cal T L}  (t,x)  |    &\les&  c (\gamma)  \cdot c (\delta) \cdot  E ( 3)  \cdot \frac{\eps \cdot (1+|q|)}{(1+t+|q|)  }  \; .
      \eea

By choosing $\varepsilon^\prime$ as in \eqref{definitionofepsprime}, we obtain the result.

\end{proof}

Now that the assumptions for the estimate \ref{LinfinitynormestimateongradientderivedbyLondbladRodnianski} on solutions to the wave equation with sources are satisfied, we would like to apply it to the Einstein-Yang-Mills system.   
For this, we always start by estimating the source terms.

\subsection{Upgrading the estimates for the “good” components of the potential} \

  \begin{lemma}\label{decayestimateonthesourcetermforgoodpoentialAcomponents}
Let $M\leq \eps$\,. We have for $\ga > \de $\,, and $0 \leq \de \leq 1$\,, 

                                        \beaa
 |  \Lie_{Z^I} \Big( g^{\la\mu} \derm_{\la}   \derm_{\mu}   A_{{\cal T}}  \Big)  |   &\leq& \begin{cases} c (\delta) \cdot c (\gamma) \cdot C ( |I| ) \cdot E ( |I| + 3)   \cdot \frac{\eps }{(1+t+|q|)^{3-3\delta} (1+|q|)^{1+\ga}},\quad\text{when }\quad q>0\; ,\\
       C ( |I| ) \cdot E ( |I| + 3)  \cdot \frac{\eps   \cdot (1+| q |   )^{\frac{3}{2} }  }{(1+t+|q|)^{3-3\delta} }  \,\quad\text{when }\quad q<0 \; . \end{cases} \\
           \notag
       \eeaa

    Thus, in the exterior region $\overline{C}$, we have
                                                \bea
                                                \notag
&& (1+ t) \cdot \| \,  \Lie_{Z^I} \Big( g^{\la\mu} \derm_{\la}   \derm_{\mu}   A_{{\cal T}}  (t,\cdot) \Big)   \|_{L^\infty(\overline{D}_t)} \\
\notag
& \leq& C(q_0) \cdot c (\delta) \cdot c (\gamma) \cdot C ( |I| ) \cdot E ( |I| + 3)   \cdot \frac{\eps }{(1+t+|q|)^{2-3\delta} (1+|q|)^{1+\ga}}  \; . \\
       \eea

\end{lemma}

\begin{proof}
We have already established in Lemma \ref{estimateonthesourcetermsforgoodcomponentofPoentialAandgoodcompometrich}, that
  \bea\label{waveequationforAtaugoodcompoenents}
   \notag
 && | g^{\la\mu} \derm_{\la}   \derm_{\mu}   A_{{\cal T}}   |  \\
 \notag
    &\les& | \derm h | \cdot  |\rderm A |       + | \rderm  h | \cdot  |\derm A |  \\
       \notag
           &&+   | A  | \cdot    | \rderm A  | +  | \derm  h | \cdot  | A  |^2   +  | A  |^3 \\
    \notag
  && + O( h \cdot  \derm h \cdot  \derm A) + O( h \cdot  A \cdot \derm A)  + O( h \cdot  \derm h \cdot  A^2) + O( h \cdot  A^3) \, , \\
  \eea
and we showed in Lemma \ref{StructureoftheLiederivativesofthesourcetermsofthewaveoperatorforAandh} that this structure is preserved for the Lie derivatives in the direction of Minkowski vector fields, which we summed up as being
  \bea
   \notag
 &&   \Lie_{Z^I} \Big( g^{\la\mu} \derm_{\la}   \derm_{\mu}   A_{{\cal T}}   \Big)    \\
 \notag
    &=& \Lie_{Z^I}  \Big[  \derm h  \cdot  \rderm A      +  \rderm  h  \cdot  \derm A   +    A   \cdot   \rderm A  +  \derm  h  \cdot  A^2   +  A^3 \\
    \notag
  && + O( h \cdot  \derm h \cdot  \derm A) + O( h \cdot  A \cdot \derm A)  + O( h \cdot  \derm h \cdot  A^2) + O( h \cdot  A^3)  \Big] \; .\\
  \eea

Now, we would like to estimate each of the terms.

\textbf{The term $  \Lie_{Z^I} \big(   \derm h \cdot  \rderm A   \big)      +  \Lie_{Z^I} \big(     \rderm  h  \cdot  \derm A \big)  $:}\\

We have
      \beaa
 \notag
&&  |  \rderm ( \Lie_{Z^I} A ) (t,x)  |  \\
\notag
 &\leq& \begin{cases} c (\gamma) \cdot C ( |I| ) \cdot E ( |I| + 3)  \cdot \frac{\eps }{(1+t+|q|)^{2-\delta} (1+|q|)^{\ga}},\quad\text{when }\quad q>0 \; ,\\
       C ( |I| ) \cdot E ( |I| + 3)  \cdot \eps \cdot \frac{ (1+|q|)^{\frac{1}{2} }}{ (1+t+|q|)^{2-\delta} } \,\quad\text{when }\quad q<0 \; , \end{cases} \\
      \eeaa
and for $\ga > \de$,
              \beaa
 \notag
 |\derm  ( \Lie_{Z^I} h ) (t,x)  |    &\leq& \begin{cases} C ( |I| ) \cdot E ( |I| + 2)  \cdot \frac{\eps }{(1+t+|q|)^{1-\delta} (1+|q|)^{1+\de}},\quad\text{when }\quad q>0 \;,\\
       C ( |I| ) \cdot E ( |I| + 2)  \cdot \frac{\eps  }{(1+t+|q|)^{1-\delta}(1+|q|)^{\frac{1}{2} }}  \,\quad\text{when }\quad q<0 \; , \end{cases} \\
      \eeaa
      hence,
           \beaa
 \notag
&&  |  \Lie_{Z^I} \big(   \derm h \cdot  \rderm A   \big)   |   \\
\notag
&\leq& \begin{cases}  c (\gamma) \cdot C ( |I| ) \cdot E ( |I| + 3)    \cdot \frac{\eps }{(1+t+|q|)^{3-2\delta} (1+|q|)^{1+\ga+\de}},\quad\text{when }\quad q>0,\\
       C ( |I| ) \cdot E ( |I| + 3)  \cdot \frac{\eps  }{(1+t+|q|)^{3-2\delta} }  \,\quad\text{when }\quad q<0 \; .\end{cases} \\
      \eeaa
      
        On the other hand,

   \beaa
 \notag
&&  |  \rderm ( \Lie_{Z^I} h ) (t,x)  |  \\
\notag
 &\leq& \begin{cases} c (\delta) \cdot c (\gamma) \cdot C ( |I| ) \cdot E ( |I| + 3)  \cdot \frac{\eps }{(1+t+|q|)^{2-\delta} (1+|q|)^{\de}},\quad\text{when }\quad q>0 \; ,\\
       C ( |I| ) \cdot E ( |I| + 3)  \cdot \eps \cdot \frac{ (1+|q|)^{\frac{1}{2} }}{ (1+t+|q|)^{2-\delta} } \,\quad\text{when }\quad q<0 \;, \end{cases} \\
      \eeaa
and
                 \beaa
 \notag
|\derm ( \Lie_{Z^I} A ) (t,x)    |    &\leq& \begin{cases} C ( |I| ) \cdot E ( |I| + 2)  \cdot \frac{\eps }{(1+t+|q|)^{1-\delta} (1+|q|)^{1+\ga}},\quad\text{when }\quad q>0\; ,\\
       C ( |I| ) \cdot E ( |I| + 2)  \cdot \frac{\eps  }{(1+t+|q|)^{1-\delta}(1+|q|)^{\frac{1}{2} }}  \,\quad\text{when }\quad q<0 \; , \end{cases} \\
      \eeaa
    thus,
                     \beaa
 \notag
 | \Lie_{Z^I} \big(     \rderm  h  \cdot  \derm A \big)  |    &\leq& \begin{cases} c (\delta) \cdot c (\gamma) \cdot C ( |I| ) \cdot E ( |I| + 3)   \cdot \frac{\eps }{(1+t+|q|)^{3-2\delta} (1+|q|)^{1+\ga+\de}},\quad\text{when }\quad q>0\; ,\\
       C ( |I| ) \cdot E ( |I| + 3)  \cdot \frac{\eps  }{(1+t+|q|)^{3-2\delta} }  \,\quad\text{when }\quad q<0 \; . \end{cases} \\
      \eeaa
      
      As a result,
      
                           \bea
 \notag
&& | \Lie_{Z^I} \big(   \derm h \cdot  \rderm A   \big)      +  \Lie_{Z^I} \big(     \rderm  h  \cdot  \derm A \big)   |   \\
\notag
  &\leq& \begin{cases} c (\delta) \cdot c (\gamma) \cdot C ( |I| ) \cdot E ( |I| + 3)   \cdot \frac{\eps }{(1+t+|q|)^{3-2\delta} (1+|q|)^{1+\ga+\de}},\quad\text{when }\quad q>0\; ,\\
       C ( |I| ) \cdot E ( |I| + 3)  \cdot \frac{\eps  }{(1+t+|q|)^{3-2\delta} }  \,\quad\text{when }\quad q<0 \; . \end{cases} \\
      \eea

          \textbf{The term  $  \Lie_{Z^I} \big(   A  \cdot     \rderm A   \big)  $:}\\
          
 We have
           \beaa
 \notag |  \rderm ( \Lie_{Z^I} A ) (t,x)  | 
 &\leq& \begin{cases} c (\delta) \cdot c (\gamma) \cdot C ( |I| ) \cdot E ( |I| + 3)  \cdot \frac{\eps }{(1+t+|q|)^{2-\delta} (1+|q|)^{\gamma}},\quad\text{when }\quad q>0,\\
       C ( |I| ) \cdot E ( |I| + 3)  \cdot \eps \cdot \frac{ (1+|q|)^{\frac{1}{2} }}{ (1+t+|q|)^{2-\delta} } \,\quad\text{when }\quad q<0 , \end{cases} \\
          \eeaa
  and 
        \bea
 \notag
| \Lie_{Z^I}  A  (t,x)  |        &\leq& \begin{cases} c (\gamma) \cdot  C ( |I| ) \cdot E ( |I| +  2)  \cdot \frac{\eps }{(1+t+|q|)^{1-\delta} (1+|q|)^{\ga}},\quad\text{when }\quad q>0,\\
       C ( |I| ) \cdot E ( |I| + 2 )  \cdot \frac{\eps \cdot (1+| q |   )^{\frac{1}{2}} }{(1+t+|q|)^{1-\delta} }  \,\quad\text{when }\quad q<0 . \end{cases} \\
      \notag
      \eea
      Therefore, for $\eps \leq 1$, we get
            \bea
 \notag 
  | \Lie_{Z^I} \big(   A  \cdot     \rderm A   \big)  |   &\leq&  \begin{cases} c (\delta) \cdot c (\gamma) \cdot C ( |I| ) \cdot E ( |I| + 3)  \cdot \frac{\eps }{(1+t+|q|)^{3-2\delta} (1+|q|)^{2\gamma}},\quad\text{when }\quad q>0,\\
       C ( |I| ) \cdot E ( |I| + 3)  \cdot \eps \cdot \frac{ (1+|q|)}{ (1+t+|q|)^{3-2\delta} } \,\quad\text{when }\quad q<0 . \end{cases} \\
          \eea

\textbf{The term $ \Lie_{Z^I} \big(   \derm  h  \cdot   A^2 \big) $:}\\

      We have shown in Lemma \eqref{aprioridecayestimates}, that for $\ga > \de$,
                    \beaa
 \notag
 |\derm (  \Lie_{Z^I} h ) (t,x)  |   &\leq& \begin{cases} C ( |I| ) \cdot E ( |I| + 2)  \cdot \frac{\eps }{(1+t+|q|)^{1-\delta} (1+|q|)^{1+\de}},\quad\text{when }\quad q>0 \; ,\\
       C ( |I| ) \cdot E ( |I| + 2)  \cdot \frac{\eps  }{(1+t+|q|)^{1-\delta}(1+|q|)^{\frac{1}{2} }}  \,\quad\text{when }\quad q<0 \; , \end{cases} 
       \eeaa
   and
        \bea
 \notag
| \Lie_{Z^I}  A  (t,x)  |^2        &\leq& \begin{cases} c (\gamma) \cdot  C ( |I| ) \cdot E ( |I| +  2)  \cdot \frac{\eps }{(1+t+|q|)^{2-2\delta} (1+|q|)^{2\ga}},\quad\text{when }\quad q>0 \; ,\\
       C ( |I| ) \cdot E ( |I| + 2 )  \cdot \frac{\eps \cdot (1+| q |   ) }{(1+t+|q|)^{2-2\delta} }  \,\quad\text{when }\quad q<0 \; . \end{cases} \\
      \notag
      \eea
       Thus, 
           \bea
           \notag
 | \Lie_{Z^I} \big(   \derm  h  \cdot   A^2 \big) |      &\leq& \begin{cases} c (\gamma) \cdot  C ( |I| ) \cdot E ( |I| +  2)  \cdot \frac{\eps }{(1+t+|q|)^{3-3\delta} (1+|q|)^{1+\de+2\ga}},\quad\text{when }\quad q>0 \; ,\\
       C ( |I| ) \cdot E ( |I| + 2 )  \cdot \frac{\eps \cdot (1+| q |   )^{\frac{1}{2} } }{(1+t+|q|)^{3-3\delta}  }  \,\quad\text{when }\quad q<0  \;. \end{cases} \\
      \eea

        \textbf{The term $ \Lie_{Z^I} \big(   A^3  \big) $:}\\
      
  We have
           \bea
 \notag
| \Lie_{Z^I}  A  (t,x)  |^3        &\leq& \begin{cases} c (\gamma) \cdot  C ( |I| ) \cdot E ( |I| +  2)  \cdot \frac{\eps }{(1+t+|q|)^{3-3\delta} (1+|q|)^{3\ga}},\quad\text{when }\quad q>0 \;,\\
       C ( |I| ) \cdot E ( |I| + 2 )  \cdot \frac{\eps \cdot (1+| q |   )^{\frac{3}{2} } }{(1+t+|q|)^{3-3\delta} }  \,\quad\text{when }\quad q<0\; . \end{cases} \\
      \eea

Now, we want to study the terms  $ \Lie_{Z^I}  O( h \cdot  \derm h \cdot  \derm A)$\,, $  \Lie_{Z^I}  O( h \cdot  A \cdot \derm A)$\,, $ \Lie_{Z^I}  O( h \cdot  \derm h \cdot  A^2)$\,, $ \Lie_{Z^I}  O( h \cdot  A^3)$\,.

We recall that we obtained in \eqref{apriorestimimateonproducthandgradientofh}, that we have
                     \bea
 \notag
&& |   \Lie_{Z^I} h (t,x) |  \cdot |\derm ( \Lie_{Z^I} h ) (t,x)  |  \\
 \notag
  &\leq& \begin{cases} c (\delta) \cdot  c (\gamma) \cdot C ( |I| ) \cdot E ( |I| + 2)  \cdot \frac{\eps^2 }{(1+t+|q|)^{2-2\delta} (1+|q|)^{1+2\de}},\quad\text{when }\quad q>0,\\
       C ( |I| ) \cdot E ( |I| + 2)  \cdot \frac{\eps^2  }{(1+t+|q|)^{2-2\delta} }  \,\quad\text{when }\quad q<0 , \end{cases} \\
      \eea
      and from Lemma \ref{aprioridecayestimates}, we have
                   \bea
 \notag
&& |   \Lie_{Z^I} h (t,x) |  \cdot |\derm ( \Lie_{Z^I} A ) (t,x)  |  \\
 \notag
  &\leq& \begin{cases} c (\delta) \cdot  c (\gamma) \cdot C ( |I| ) \cdot E ( |I| + 2)  \cdot \frac{\eps^2 }{(1+t+|q|)^{2-2\delta} (1+|q|)^{1+\de+\ga}},\quad\text{when }\quad q>0,\\
       C ( |I| ) \cdot E ( |I| + 2)  \cdot \frac{\eps^2  }{(1+t+|q|)^{2-2\delta} }  \,\quad\text{when }\quad q<0 . \end{cases} \\
      \eea

      \textbf{The term $ \Lie_{Z^I} O( h \cdot  \derm h \cdot  \derm A)$:}\\ 
      We have
                           \bea
 \notag
&& |   \Lie_{Z^I} h (t,x) |  \cdot |\derm ( \Lie_{Z^I} h ) (t,x)  | \cdot  |\derm ( \Lie_{Z^I} A ) (t,x)  | \\
 \notag
  &\leq& \begin{cases} c (\delta) \cdot  c (\gamma) \cdot C ( |I| ) \cdot E ( |I| + 2)  \cdot \frac{\eps^3 }{(1+t+|q|)^{3-3\delta} (1+|q|)^{2+2\de+\ga}},\quad\text{when }\quad q>0,\\
       C ( |I| ) \cdot E ( |I| + 2)  \cdot \frac{\eps^3  }{(1+t+|q|)^{3-3\delta} \cdot (1+|q|)^{\frac{1}{2}}}  \,\quad\text{when }\quad q<0 . \end{cases} \\
      \eea

      \textbf{The term $ \Lie_{Z^I} O( h \cdot  A \cdot \derm A)$:}\\ 
      
      We have
                      \beaa
 \notag
 |   \Lie_{Z^I} h (t,x)  | &\leq& \begin{cases} c (\delta) \cdot c (\gamma) \cdot C ( |I| ) \cdot E ( |I| + 2) \cdot  \frac{\eps}{ (1+ t + | q | )^{1-\delta }  (1+| q |   )^{\de}}  ,\quad\text{when }\quad q>0,\\
\notag
    C ( |I| ) \cdot E ( |I| + 2) \cdot  \frac{\eps}{ (1+ t + | q | )^{1-\delta }  } (1+| q |   )^{\frac{1}{2} }  , \,\quad\text{when }\quad q<0 , \end{cases} \\ 
    \eeaa
    and therefore,
              \bea
 \notag
&& |   \Lie_{Z^I} h  |  \cdot | \Lie_{Z^I}  A    |  \cdot |\derm ( \Lie_{Z^I} A )   |     \\
\notag
 &\leq& \begin{cases}c (\delta) \cdot   c (\gamma) \cdot  C ( |I| ) \cdot E ( |I| +  2)  \cdot \frac{\eps }{(1+t+|q|)^{3-3\delta} (1+|q|)^{1+\de+2\ga}},\quad\text{when }\quad q>0,\\
       C ( |I| ) \cdot E ( |I| + 2 )  \cdot \frac{\eps \cdot (1+|q|)^{\frac{1}{2}} }{(1+t+|q|)^{3-3\delta} }  \,\quad\text{when }\quad q<0 . \end{cases} \\
      \eea
      
           \textbf{The term $\Lie_{Z^I} O( h \cdot  \derm h \cdot   A^2)$:}\\ 
      
      Since under the bootstrap assumption  $|   \Lie_{Z^I} h  | \leq  c (\gamma)  \cdot c (\delta) \cdot C ( |I| )  \cdot E ( |I| + 2)$, we have
                   \bea
 \notag
&&  |   \Lie_{Z^I} h (t,x) |  \cdot |\derm ( \Lie_{Z^I} h ) (t,x)  | \cdot  |   \Lie_{Z^I} A (t,x) |^2   \\
 \notag
 &\leq &c (\delta) \cdot  | \derm  h |  \cdot  | A  |^2     \\
\notag
 &\leq& \begin{cases} c (\delta) \cdot  c (\gamma) \cdot  C ( |I| ) \cdot E ( |I| +  2)  \cdot \frac{\eps }{(1+t+|q|)^{3-3\delta} (1+|q|)^{1+\de+2\ga}},\quad\text{when }\quad q>0 \; ,\\
       C ( |I| ) \cdot E ( |I| + 2 )  \cdot \frac{\eps \cdot (1+| q |   )^{\frac{1}{2} } }{(1+t+|q|)^{3-3\delta}  }  \,\quad\text{when }\quad q<0  \;. \end{cases} \\
      \eea

            \textbf{The term $ \Lie_{Z^I}   O( h \cdot  A^3) $:}\\ 
            Since under the bootstrap assumption  $|   \Lie_{Z^I} h   | \leq c (\gamma)  \cdot c (\delta) \cdot C ( |I| )  \cdot E ( |I| + 2)$, we have
                      \bea
 \notag
 |   \Lie_{Z^I} h   | \cdot   |   \Lie_{Z^I} A  |^3  &\leq& c (\delta) \cdot   | \Lie_{Z^I}  A    |^3   \\
&\leq& \begin{cases} c (\delta) \cdot   c (\gamma) \cdot  C ( |I| ) \cdot E ( |I| +  2)  \cdot \frac{\eps }{(1+t+|q|)^{3-3\delta} (1+|q|)^{3\ga}},\quad\text{when }\quad q>0 \;,\\
       C ( |I| ) \cdot E ( |I| + 2 )  \cdot \frac{\eps \cdot (1+| q |   )^{\frac{3}{2} } }{(1+t+|q|)^{3-3\delta} }  \,\quad\text{when }\quad q<0\; . \end{cases} \\
      \eea
      
              \textbf{The final term  $| \Lie_{Z^I}   \big(  g^{\la\mu} \derm_{\la}   \derm_{\mu}   A_{{\cal T}}    \big)   $:}\\

  Putting all the terms together, respectively in the order of the terms in the right hand side of \eqref{waveequationforAtaugoodcompoenents}, and for interest, we write more precise estimates that what is needed for $ \Lie_{Z^I} O( h \cdot  \derm h \cdot   A^2)$ and $  \Lie_{Z^I} O( h \cdot  A^3) $\,, we obtain
                                            \beaa
  && | \Lie_{Z^I}  \big(  g^{\la\mu} \derm_{\la}   \derm_{\mu}   A_{{\cal T}} \big)   | \\
 \notag
  \leq&& \begin{cases} c (\delta) \cdot c (\gamma) \cdot C ( |I| ) \cdot E ( |I| + 3)   \cdot \frac{\eps }{(1+t+|q|)^{3-2\delta} (1+|q|)^{1+\ga+\de}},\quad\text{when }\quad q>0\; ,\\
       C ( |I| ) \cdot E ( |I| + 3)  \cdot \frac{\eps  }{(1+t+|q|)^{3-2\delta} }  \,\quad\text{when }\quad q<0 \; . \end{cases} \\
           \notag
  &+&  \begin{cases} c (\delta) \cdot c (\gamma) \cdot C ( |I| ) \cdot E ( |I| + 3)  \cdot \frac{\eps }{(1+t+|q|)^{3-2\delta} (1+|q|)^{2\gamma}},\quad\text{when }\quad q>0,\\
       C ( |I| ) \cdot E ( |I| + 3)  \cdot \eps \cdot \frac{ (1+|q|)}{ (1+t+|q|)^{3-2\delta} } \,\quad\text{when }\quad q<0 . \end{cases}  \\
                   \notag 
    &+& \begin{cases} c (\gamma) \cdot  C ( |I| ) \cdot E ( |I| +  2)  \cdot \frac{\eps }{(1+t+|q|)^{3-3\delta} (1+|q|)^{1+\de+2\ga}},\quad\text{when }\quad q>0 \; ,\\
       C ( |I| ) \cdot E ( |I| + 2 )  \cdot \frac{\eps \cdot (1+| q |   )^{\frac{1}{2} } }{(1+t+|q|)^{3-3\delta}  }  \,\quad\text{when }\quad q<0  \;. \end{cases} \\
 \notag
       &+& \begin{cases} c (\gamma) \cdot  C ( |I| ) \cdot E ( |I| +  2)  \cdot \frac{\eps }{(1+t+|q|)^{3-3\delta} (1+|q|)^{3\ga}},\quad\text{when }\quad q>0 \;,\\
       C ( |I| ) \cdot E ( |I| + 2 )  \cdot \frac{\eps \cdot (1+| q |   )^{\frac{3}{2} } }{(1+t+|q|)^{3-3\delta} }  \,\quad\text{when }\quad q<0\; , \end{cases} 
       \eeaa 
       \beaa
\notag
  &+& \begin{cases} c (\delta) \cdot  c (\gamma) \cdot C ( |I| ) \cdot E ( |I| + 2)  \cdot \frac{\eps^3 }{(1+t+|q|)^{3-3\delta} (1+|q|)^{2+2\de+\ga}},\quad\text{when }\quad q>0,\\
       C ( |I| ) \cdot E ( |I| + 2)  \cdot \frac{\eps^3  }{(1+t+|q|)^{3-3\delta} \cdot (1+|q|)^{\frac{1}{2}}}  \,\quad\text{when }\quad q<0 . \end{cases} \\
 \notag
 &+& \begin{cases}c (\delta) \cdot   c (\gamma) \cdot  C ( |I| ) \cdot E ( |I| +  2)  \cdot \frac{\eps }{(1+t+|q|)^{3-3\delta} (1+|q|)^{1+\de+2\ga}},\quad\text{when }\quad q>0,\\
       C ( |I| ) \cdot E ( |I| + 2 )  \cdot \frac{\eps \cdot (1+|q|)^{\frac{1}{2}} }{(1+t+|q|)^{3-3\delta} }  \,\quad\text{when }\quad q<0 . \end{cases} \\
        \notag
  &+& \begin{cases} c (\delta) \cdot  c (\gamma) \cdot C ( |I| ) \cdot E ( |I| + 2)  \cdot \frac{\eps^4 }{(1+t+|q|)^{4-4\delta} (1+|q|)^{1+2\de+2\ga}},\quad\text{when }\quad q>0,\\
       C ( |I| ) \cdot E ( |I| + 2)  \cdot \frac{\eps^4 \cdot (1+|q|)  }{(1+t+|q|)^{4-4\delta} }  \,\quad\text{when }\quad q<0  . \end{cases} \\
 \notag
&+& \begin{cases} c (\delta) \cdot c (\gamma) \cdot C ( |I| ) \cdot E ( |I| + 2) \cdot  \frac{\eps}{ (1+ t + | q | )^{4-4\delta }  (1+| q |   )^{\de+3\ga}}  ,\quad\text{when }\quad q>0,\\
\notag
    C ( |I| ) \cdot E ( |I| + 2) \cdot  \frac{\eps \cdot (1+| q |   )^2 }{ (1+ t + | q | )^{4-4\delta }  }   , \,\quad\text{when }\quad q<0 . \end{cases} \\ 
    \eeaa
    
    Consequently, considering that $\ga > 1 $, and therefore $2\ga > 1 +\ga $, we get for $0 \leq \de \leq 1$, the following estimate under the bootstrap assumption
                                            \bea\label{decayestimateonthetermsinthesourcetermsforwaveeqonAtau}
                                            \notag
  && |  \Lie_{Z^I}  \big( g^{\la\mu} \derm_{\la}   \derm_{\mu}   A_{{\cal T}}   \big) | \\
 \notag
  \leq&& \begin{cases} c (\delta) \cdot c (\gamma) \cdot C ( |I| ) \cdot E ( |I| + 3)   \cdot \frac{\eps }{(1+t+|q|)^{3-3\delta} (1+|q|)^{1+\ga}},\quad\text{when }\quad q>0\; ,\\
       C ( |I| ) \cdot E ( |I| + 3)  \cdot \frac{\eps   \cdot (1+| q |   )^{\frac{3}{2} }  }{(1+t+|q|)^{3-3\delta} }  \,\quad\text{when }\quad q<0 \; . \end{cases} \\
       \eea
       
      Since for $ q \geq q_0$, we have for $q \leq 0$ that $|q| \leq |q_0|$, we therefore obtain that for $q \leq 0$, we have
      \beaa
       \frac{(1+| q |   )^{\frac{3}{2} }  }{(1+t+|q|)^{3-3\delta} } \leq    \frac{(1+| q_0 |   )^{\frac{3}{2} }  }{(1+t+|q|)^{3-3\delta} } \; .
       \eeaa
       Also, since for $ q_0 \leq q \leq 0$, we have
       \beaa
   \frac{1}{   (1+|q_0|)^{1+\ga} }  \leq   \frac{1}{ (1+|q|)^{1+\ga}} \;,
       \eeaa
        
  we obtain that in the exterior region $\overline{C}$, for all $q$, we have
                                                \beaa
  && | \Lie_{Z^I}  \big(  g^{\la\mu} \derm_{\la}   \derm_{\mu}   A_{{\cal T}} \big)    | \\
 \notag
  \leq&& C(q_0) \cdot c (\delta) \cdot c (\gamma) \cdot C ( |I| ) \cdot E ( |I| + 3)   \cdot \frac{\eps }{(1+t+|q|)^{3-3\delta} (1+|q|)^{1+\ga}} 
       \eeaa 
    
      \end{proof}

Now, let us point out the following lemma:
 
  \begin{lemma}
  We point out the follwoing equivalence,
  \beaa
  1+t+|x| \sim 1+t+|q|  \;.
  \eeaa
   
  \end{lemma}
  
  \begin{proof}

When $q \geq 0$, then, $|q| = q = r - t $. Hence, $1+t+|q| =  1+t+ r - t = 1+ |x| $. Therefore,
   \beaa
  1+t+|x| &=&   1+2t+|q| \leq 2( 1+t+|q| ) \; ,\\
   1+t+|x|  &\geq& 1+|x| = 1+t+|q| \;.
  \eeaa
  Hence, for  $q \geq 0$, we have
  \beaa
  1+t+|x|  \sim 1+t+|q|  \;.
  \eeaa
  When $q \leq 0$, then, $|q| = - q = t - r$. Hence, $1+t+|q| =  1+t+  t -r = 1+ 2t - |x| $. Therefore,
   \beaa
1+t+|q|  &\leq& 1+2t \leq  2( 1+t+|x| )  \;,\\
1+t+|x| &=&  1+ 2t - |q|  \leq  2( 1+t+|q| )\; .
  \eeaa
Thus, for  $q \leq 0$, we have $$1+t+|x|  \sim 1+t+|q|  .$$

 \end{proof}

  \begin{lemma}\label{estimatefortangentialcomponentspotential}
In the Lorenz and harmonic gauges, the Einstein-Yang-Mills fields satisfy in the exterior region $\overline{C}$\,, for $M\leq \eps$\,, for \,$\ga > \de $\,, and $0 \leq \delta \leq \frac{1}{4}$\,, 
 \beaa
 \,|\derm A_{{\cal T}}  (t,x)| &\les&   C(q_0) \cdot  c (\gamma)  \cdot c (\delta)   \cdot E (4) \cdot \frac{\eps }{ (1+t+|q|)} \; . 
\eeaa
\end{lemma}

\begin{proof}

We have shown in Lemma \ref{estimateonthesourcetermsforgoodcomponentofPoentialAandgoodcompometrich} that in the Lorenz and the harmonic gauges, we have the following estimate for wave equations for the tangential components of the Einstein-Yang-Mills potential,
  \bea
   \notag
 && | g^{\la\mu} \derm_{\la}   \derm_{\mu}   A_{{\cal T}}   |  \\
 \notag
    &\les& | \derm h | \cdot  |\rderm A |       + | \rderm  h | \cdot  |\derm A |  \\
       \notag
           && +   | A  | \cdot    | \rderm A  |    +  | \derm  h | \cdot  | A  |^2  +   | A  |^3 \\
    \notag
  && + O( h \cdot  \derm h \cdot  \derm A) + O( h \cdot  A \cdot \derm A) + O( h \cdot  \derm h \cdot  A^2)  + O( h \cdot  A^3) \, .
  \eea

Using the following estimate for all $U,V\in\{L,\Lb,A,B\}$,
\beaa
\notag
&& (1+t+|x|)  \cdot   |\varpi(q)\pa  \Phi_{UV} (t,x)| \\
&\les& \!\sup_{0\leq \tau\leq t} \sum_{|I|\leq 1}\|\,\varpi(q) \Lie_{Z^I} \! \phi(\tau,\cdot)\|_{L^\infty  (\Sigma^{ext}_{\tau} ) }\\
\notag
&& + \int_0^t\Big( \varepsilon^\prime \cdot  ( 1+\gamma’) \cdot  \|\varpi(q) \pa  \Phi_{UV} (\tau,\cdot) \|_{L^\infty  (\Sigma^{ext}_{\tau} )} +(1+\tau) \cdot  \| \varpi(q)  S_{UV} (\tau,\cdot) \|_{L^\infty(\overline{D}_\tau)} \\
 \notag
&& +\sum_{|I|\leq 2} (1+\tau)^{-1}  \cdot   \| \varpi(q) \Lie_{Z^I}  \Phi_{UV} (\tau,\cdot)\|_{L^\infty(\overline{D}_\tau)}\Big)\, d\tau \; , \\
\eeaa
with $1 + \gamma' = 0$, and therefore
\beaa
\varpi &:=& \varpi(q) :=\begin{cases}
(1+|q|)^{1+\gamma^\prime} = 1,\quad\text{when }\quad q>0\\
      1 \,\quad\text{when }\quad
      q<0\end{cases},\\
      &=& 1 \, 
\eeaa
leads to
\begin{multline}
(1+t+|x|) \cdot \,|\derm A_{{\cal T}} (t,x)|  \les\!\sup_{0\leq
\tau\leq t}
\sum_{|I|\leq 1}\|\, \Lie_{Z^I} \! A (\tau,\cdot)\|_{L^\infty (\Sigma^{ext}_{\tau} ) }\\
+ \int_0^t\Big( (1+\tau)  \cdot   \|
\, g^{\la\mu} \derm_{\la}   \derm_{\mu}   A_{{\cal T}}  (\tau,\cdot)   \|_{L^\infty(\overline{D}_\tau)} +\sum_{|I|\leq 2} (1+\tau)^{-1}  \cdot   \| \Lie_{Z^I} A(\tau,\cdot)\|_{L^\infty(\overline{D}_\tau)}\Big)\, d\tau.
\end{multline}

  Since
  \beaa
  |   \Lie_{Z^I} h (t,x)  | + |   \Lie_{Z^I} A (t,x)  |   &\leq&  c (\gamma)  \cdot c (\delta) \cdot C ( |I| )  \cdot E ( |I| + 2)  \cdot \frac{\eps}{(1+t+|q|)^{\frac{1}{2}-\delta}  } ,\\
 &\leq&    c (\gamma)  \cdot c (\delta) \cdot C ( |I| )  \cdot E ( |I| + 2) \cdot \eps  \,,
     \eeaa
we have
\bea
\!\sup_{0\leq \tau\leq t} \sum_{|I|\leq 1}\|\,  \Lie_{Z^I} \! A (\tau,\cdot)\|_{L^\infty (\Sigma^{ext}_{\tau} ) } \leq  c (\gamma)  \cdot c (\delta)   \cdot E ( 3) \cdot \eps
\eea
and 
\bea
\notag
 (1+\tau)^{-1}  \cdot   \|  \Lie_{Z^I} A(\tau,\cdot)\|_{L^\infty(\overline{D}_\tau)} &\leq&    c (\gamma)  \cdot c (\delta) \cdot C ( |I| )  \cdot E ( |I| + 2) \cdot \frac{\eps}{(1+t)^{\frac{3}{2}-\delta}  } .\\
\eea
We have already proved in Lemma \eqref{decayestimateonthesourcetermforgoodpoentialAcomponents}, that
\bea
\notag
&& (1+ t) \cdot \| \, g^{\la\mu} \derm_{\la}   \derm_{\mu}   A_{{\cal T}}  (t,\cdot)  \|_{L^\infty(\overline{D}_t)} \\
\notag
& \les& C(q_0) \cdot c (\delta) \cdot c (\gamma)  \cdot E ( 3)   \cdot \frac{\eps }{(1+t+|q|)^{2-3\delta} (1+|q|)^{1+\ga}}   \\
\notag
 &\les& C(q_0) \cdot  c (\delta) \cdot  c (\gamma) \cdot E (  3)   \cdot \frac{\eps  }{(1+t)^{2 -3\delta  } } \; .\\
\eea

Finally, we get
\beaa
&& (1+t+|x|) \cdot \,|\derm A_{{\cal T}} (t,x)| \\
 &\les& \!\sup_{0\leq \tau\leq t} \sum_{|I|\leq 1}\|\,  \Lie_{Z^I} \! A (\tau,\cdot)\|_{L^\infty (\Sigma^{ext}_{\tau} ) }\\
&& + \int_0^t\Big( (1+\tau)\| \, g^{\la\mu} \derm_{\la}   \derm_{\mu}   A_{{\cal T}}  (\tau,\cdot)  \|_{L^\infty(\overline{D}_\tau)} +\sum_{|I|\leq 2} (1+\tau)^{-1} \|  \Lie_{Z^I} A(\tau,\cdot)\|_{L^\infty(\overline{D}_\tau)}\Big)\, d\tau \\
&\les&   c (\gamma)  \cdot c (\delta) \cdot  E (3) \cdot \eps  +  \int_0^t  C(q_0) \cdot  c (\delta) \cdot  c (\gamma)  \cdot E (  3)   \cdot \frac{\eps  }{(1+\tau)^{2 -3\delta  } } d\tau \\
&& +   \int_0^t  c (\gamma)  \cdot c (\delta)  \cdot E ( 4) \cdot \frac{\eps}{(1+\tau)^{\frac{3}{2}-\delta }  }   d\tau \\
&\les&   c (\gamma)  \cdot c (\delta) \cdot  E (3) \cdot \eps  +  \int_0^t  C(q_0) \cdot  c (\delta) \cdot  c (\gamma)  \cdot E (  3)   \cdot \frac{\eps  }{(1+\tau)^{\frac{3}{2}-\delta  } \cdot (1+\tau)^{\frac{1}{2}-2\delta  }} d\tau \\
&& +   \int_0^t  c (\gamma)  \cdot c (\delta)  \cdot E ( 4) \cdot \frac{\eps}{(1+\tau)^{\frac{3}{2}-\delta }  }   d\tau \\
&& \text{(we choose here $\delta \leq \frac{1}{4} $)} \\
&\les&   C(q_0) \cdot  c (\gamma)  \cdot c (\delta)   \cdot E (4) \cdot \eps \big[  \frac{1}{(1+\tau)^{\frac{1}{2} -\delta } } \big]^{\infty}_{0} \\
&& \text{(which is integrable for $\delta \leq \frac{1}{4} <  \frac{1}{2}$)}.
\eeaa
Hence, for $\delta \leq \frac{1}{4}$, we obtain
\beaa
 (1+t+|x|) \cdot \,|\derm A_{{\cal T}} (t,x)|  &\les&  C(q_0) \cdot   c (\gamma)  \cdot c (\delta)   \cdot E (4) \cdot \eps \; ,
\eeaa
 which leads to 
 \bea
 \,|\derm A_{{\cal T}} (t,x)|  &\les&  C(q_0) \cdot   c (\gamma)  \cdot c (\delta)   \cdot E (4) \cdot \frac{\eps }{ (1+t+|x|)} \;  . 
\eea

\end{proof}

\subsection{Upgrading the estimates for the “good” components of the metric}\

\begin{lemma}\label{decayestimateonsourcetermsforgoodcomponentofmetrichtauU}
We have in the exterior region $\overline{C}$, for $M\leq \eps$\,, for $\ga > \de $\,, and $0 \leq \de \leq 1$\,,

  \bea
\notag
&& (1+t ) \cdot | \Lie_{Z^I} \Big( g^{\alpha\beta} \derm_\alpha \derm_\beta h_{ {\cal T} {\cal U}}  \Big) |  \\
&\les&  \begin{cases}  c (\delta) \cdot c (\gamma) \cdot C ( |I| ) \cdot E ( |I| + 3)  \cdot \frac{\eps }{(1+t+|q|)^{2-3\delta} (1+|q|)^{1+2\de}},\quad\text{when }\quad q>0  \;  ,\\
       \notag
       C ( |I| ) \cdot E ( |I| + 3)  \cdot \frac{\eps  \cdot (1+| q |   )^2 }{(1+t+|q|)^{2-3\delta} }  \,\quad\text{when }\quad q<0 \; .\end{cases} \\
        \eea

        Thus, in the exterior region $\overline{C}$, we have
          \bea
\notag
&&(1+t ) \cdot |  \Lie_{Z^I} \Big( g^{\alpha\beta} \derm_\alpha \derm_\beta h_{ {\cal T} {\cal U}} \Big) |_{L^\infty(\overline{D}_\tau)} \\
 &\les&
       C(q_0) \cdot  c (\delta) \cdot c (\gamma) \cdot E ( |I| +3)  \cdot \frac{\eps }{(1+t )^{2-3\delta} (1+|q|)^{1+2\de} } \; . 
    \eea
\end{lemma}

\begin{proof}

We showed in Lemma \ref{estimateonthesourcetermsforgoodcomponentofPoentialAandgoodcompometrich} that in the Lorenz and harmonic gauges, we have the following wave equations on the components of the metric,
 \bea
\notag
 && | g^{\alpha\beta} \derm_\alpha \derm_\beta h_{ {\cal T} {\cal U}} |   \\
 \notag
    &\les&  | P(\derm_{\cal T} h,\derm_{\cal U} h) | + | Q_{{\cal T} {\cal U}}(\derm h,\derm h)  | + | G_{{\cal T} {\cal U} }(h)(\derm h,\derm h)  | \\
  \notag
    && +  | \rderm A   | \cdot  |\derm  A  |  +    | A |^2 \cdot | \derm A   |   + |A |^4  \\
\notag
     && + O \big(h \cdot  (\derm A)^2 \big)   + O \big(  h  \cdot  A^2 \cdot \derm A \big)     + O \big(  h   \cdot  A^4 \big)  . 
\eea
we showed in Lemma \ref{StructureoftheLiederivativesofthesourcetermsofthewaveoperatorforAandh} that this structure is preserved for the Lie derivatives in the direction of Minkowski vector fields, which we wrote as
 \bea
\notag
 &&  \Lie_{Z^I} \Big(  g^{\alpha\beta} \derm_\alpha \derm_\beta h^1_{ {\cal T} {\cal U}} \Big)    \\
 \notag
    &=&\Lie_{Z^I} \Big[  P(\derm_{\cal T} h,\derm_{\cal U} h)  + Q_{{\cal T} {\cal U}}(\derm h,\derm h)   +  G_{{\cal T} {\cal U} }(h)(\derm h,\derm h)   \\
    \notag
    &&   +  \rderm A   \cdot  \derm  A   +    A^2 \cdot \derm A     + A^4  \\
\notag
     && + O \big(h \cdot  (\derm A)^2 \big)   + O \big(  h  \cdot  A^2 \cdot \derm A \big)     + O \big(  h   \cdot  A^4 \big) \Big]  \\
  \notag
  && +    \Lie_{Z^I} \Big( g^{\alpha\beta} \derm_\alpha \derm_\beta h^0 \Big) \; .\\
\eea

\textbf{The term $\Lie_{Z^I} \Big[  P(\derm_{\cal T} h,\derm_{\cal U} h)  + Q_{{\cal T} {\cal U}}(\derm h,\derm h)   +  G_{{\cal T} {\cal U} }(h)(\derm h,\derm h) \Big] $:}\\

Based on Lemma 4.2 and Corollary 9.7 in \cite{LR10}, we have
 \bea\label{estimateonthetermsgeneratedfromtheEinsteinvacuumequationswhichhaveanullstructure}
\notag
&& | \Lie_{Z^I}  P(\pa_{\cal T} h,\pa_{\cal U} h) | + | \Lie_{Z^I}  Q_{{\cal T} {\cal U}}(\pa h,\pa h)  | + | \Lie_{Z^I}  G_{{\cal T} {\cal U} }(h)(\pa h,\pa h)  | \\
\notag
    &\les& | \Lie_{Z^I} \Big[  \rderm h \cdot \derm h + h \cdot (\derm h)^2  \Big] | \; . \\
\eea

However, we already showed in Lemma \ref{aprioridecayestimates}, that $\ga > \de$,

   \beaa
 \notag
&&  |  \rderm ( \Lie_{Z^I} h ) (t,x)  |  \\
\notag
 &\leq& \begin{cases} c (\delta) \cdot c (\gamma) \cdot C ( |I| ) \cdot E ( |I| + 3)  \cdot \frac{\eps }{(1+t+|q|)^{2-\delta} (1+|q|)^{\de}},\quad\text{when }\quad q>0,\\
       C ( |I| ) \cdot E ( |I| + 3)  \cdot \eps \cdot \frac{ (1+|q|)^{\frac{1}{2} }}{ (1+t+|q|)^{2-\delta} } \,\quad\text{when }\quad q<0 , \end{cases} \\
      \eeaa
and
              \beaa
 \notag
 |\derm  ( \Lie_{Z^I} h ) (t,x)  |    &\leq& \begin{cases} C ( |I| ) \cdot E ( |I| + 2)  \cdot \frac{\eps }{(1+t+|q|)^{1-\delta} (1+|q|)^{1+\de}},\quad\text{when }\quad q>0 \;,\\
       C ( |I| ) \cdot E ( |I| + 2)  \cdot \frac{\eps  }{(1+t+|q|)^{1-\delta}(1+|q|)^{\frac{1}{2} }}  \,\quad\text{when }\quad q<0 \; . \end{cases} \\
      \eeaa
 Thus, for $\ga > \de$, we have  
           \beaa
 \notag
&& \Lie_{Z^I} \Big( \rderm h \cdot \derm h \Big)  \\
\notag
&\leq& \begin{cases}  c (\delta) \cdot c (\gamma) \cdot C ( |I| ) \cdot E ( |I| + 3)  \cdot \frac{\eps }{(1+t+|q|)^{3-2\delta} (1+|q|)^{1+2\de}},\quad\text{when }\quad q>0,\\
       C ( |I| ) \cdot E ( |I| + 3)  \cdot \frac{\eps  }{(1+t+|q|)^{3-2\delta} }  \,\quad\text{when }\quad q<0 \; .\end{cases} \\
      \eeaa

Whereas to the term $  \Lie_{Z^I} \Big(  h \cdot ( \derm h )^2 \Big) $\,, we have
              \beaa
 \notag
 |\derm  ( \Lie_{Z^I} h )  |^2    &\leq& \begin{cases} C ( |I| ) \cdot E ( |I| + 2)  \cdot \frac{\eps^2 }{(1+t+|q|)^{2-2\delta} (1+|q|)^{2+2\de}},\quad\text{when }\quad q>0,\\
       C ( |I| ) \cdot E ( |I| + 2)  \cdot \frac{\eps^2  }{(1+t+|q|)^{2-2\delta}(1+|q|)}  \,\quad\text{when }\quad q<0 , \end{cases} \\
      \eeaa
      and
                 \beaa
 \notag
 |   \Lie_{Z^I} h (t,x)  | &\leq& \begin{cases} c (\delta) \cdot c (\gamma) \cdot C ( |I| ) \cdot E ( |I| + 2) \cdot  \frac{\eps}{ (1+ t + | q | )^{1-\delta }  (1+| q |   )^{\de}}  ,\quad\text{when }\quad q>0,\\
\notag
    C ( |I| ) \cdot E ( |I| + 2) \cdot  \frac{\eps}{ (1+ t + | q | )^{1-\delta }  } (1+| q |   )^{\frac{1}{2} }  , \,\quad\text{when }\quad q<0 . \end{cases} \\ 
    \eeaa
    Hence,

 \beaa
\notag
 | \Lie_{Z^I} \Big(  h \cdot ( \derm h )^2 \Big)|  &\leq& \begin{cases}  c (\delta) \cdot c (\gamma) \cdot C ( |I| ) \cdot E ( |I| + 2) \cdot \frac{\eps^3 }{(1+t+|q|)^{3-3\delta} (1+|q|)^{2+3\de}},\quad\text{when }\quad q>0,\\
       C ( |I| ) \cdot E ( |I| + 2)  \cdot \frac{\eps^3  }{(1+t+|q|)^{3-3\delta}(1+|q|)^{\frac{1}{2} }  }  \,\quad\text{when }\quad q<0 , \end{cases} \\
      \eeaa
Therefore, for $\de \leq 1$, we have for $\eps \leq 1$,

 \bea
\notag
&& | \Lie_{Z^I}  P(\derm_{\cal T} h,\derm_{\cal U} h) | + | \Lie_{Z^I}  Q_{{\cal T} {\cal U}}(\derm h,\derm h)  | + | \Lie_{Z^I}  G_{{\cal T} {\cal U} }(h)(\derm h,\derm h)  | \\
\notag
&\leq& \begin{cases}  c (\delta) \cdot c (\gamma) \cdot C ( |I| ) \cdot E ( |I| + 3)  \cdot \frac{\eps }{(1+t+|q|)^{3-2\delta} (1+|q|)^{1+2\de}},\quad\text{when }\quad q>0  \;  ,\\
\notag
       C ( |I| ) \cdot E ( |I| + 3)  \cdot \frac{\eps  }{(1+t+|q|)^{3-2\delta} }  \,\quad\text{when }\quad q<0 \; ,\end{cases} \\
       \notag
       &&+  \begin{cases}  c (\delta) \cdot c (\gamma) \cdot C ( |I| ) \cdot E ( |I| + 2) \cdot \frac{\eps^3 }{(1+t+|q|)^{3-3\delta} (1+|q|)^{2+3\de}},\quad\text{when }\quad q>0  \;  ,\\
       \notag
       C ( |I| ) \cdot E ( |I| + 2)  \cdot \frac{\eps^3  }{(1+t+|q|)^{3-3\delta}(1+|q|)^{\frac{1}{2} }  }  \,\quad\text{when }\quad q<0 \; . \end{cases} \\
       \notag
       &\leq& \begin{cases}  c (\delta) \cdot c (\gamma) \cdot C ( |I| ) \cdot E ( |I| + 3)  \cdot \frac{\eps }{(1+t+|q|)^{3-3\delta} (1+|q|)^{1+2\de}},\quad\text{when }\quad q>0  \;  ,\\
       \notag
       C ( |I| ) \cdot E ( |I| + 3)  \cdot \frac{\eps  }{(1+t+|q|)^{3-2\delta} }  \,\quad\text{when }\quad q<0 \; .\end{cases} \\
      \eea

\textbf{The terms $ \Lie_{Z^I}  \Big[ \rderm A    \cdot  \derm  A    +  O \big(h \cdot  (\derm A)^2 \big) \Big]   $:}\\

Since $\ga > \de$, the estimates that we obtained from the bootstrap assumption for $ | \Lie_{Z^I}  \Big( \rderm h \cdot \derm h \Big) |   $ and $| \Lie_{Z^I}  \Big( h \cdot ( \derm h )^2 \Big) |$ are also true for $|\Lie_{Z^I}  \Big( \rderm A|\cdot \derm A \Big) |$ and $ | \Lie_{Z^I}  \Big( A \cdot (\derm A)^2 \Big) |$, respectively.
Hence,
        \bea
 \notag
 && |  \Lie_{Z^I}  \Big[ \rderm A    \cdot  \derm  A    +  O \big(h \cdot  (\derm A)^2 \big) \Big]   | \\
 \notag
      &\les& \begin{cases}  c (\delta) \cdot c (\gamma) \cdot C ( |I| ) \cdot E ( |I| + 3)  \cdot \frac{\eps }{(1+t+|q|)^{3-3\delta} (1+|q|)^{1+2\de}},\quad\text{when }\quad q>0  \;  ,\\
       \notag
       C ( |I| ) \cdot E ( |I| + 3)  \cdot \frac{\eps  }{(1+t+|q|)^{3-2\delta} }  \,\quad\text{when }\quad q<0 \; .\end{cases} \\
      \eea
      
      \textbf{The terms $ \Lie_{Z^I}  \Big[   A^2   \cdot  \derm A    + O \big(  h  \cdot  A^2 \cdot \derm A \big)  \Big]  $:}\\
      
      We have
                 \bea
           \notag
&&  | \Lie_{Z^I}  \Big(   A^2   \cdot  \derm A  \Big) |  \\
  \notag
     &\leq& \begin{cases} c (\gamma) \cdot  C ( |I| ) \cdot E ( |I| +  2)  \cdot \frac{\eps }{(1+t+|q|)^{3-3\delta} (1+|q|)^{1+3\ga}},\quad\text{when }\quad q>0 \; ,\\
       C ( |I| ) \cdot E ( |I| + 2 )  \cdot \frac{\eps \cdot (1+| q |   )^{\frac{1}{2} } }{(1+t+|q|)^{3-3\delta}  }  \,\quad\text{when }\quad q<0  \;. \end{cases} \\
      \eea
      
Also, under the bootstrap assumption, we have $|   \Lie_{Z^I} h   | \leq c (\gamma)  \cdot c (\delta) \cdot C ( |I| )  \cdot E ( |I| + 2)$. Therefore,
        \bea
 \notag
 &&   | \Lie_{Z^I}  \Big(   A^2   \cdot  \derm A  \Big) | +  | \Lie_{Z^I}  O \big(  h  \cdot  A^2 \cdot \derm A \big) |   \\
 \notag
  &\leq& \begin{cases} c (\delta) \cdot c (\gamma) \cdot  C ( |I| ) \cdot E ( |I| +  2)  \cdot \frac{\eps }{(1+t+|q|)^{3-3\delta} (1+|q|)^{1+3\ga}},\quad\text{when }\quad q>0 \; ,\\
       C ( |I| ) \cdot E ( |I| + 2 )  \cdot \frac{\eps \cdot (1+| q |   )^{\frac{1}{2} } }{(1+t+|q|)^{3-3\delta}  }  \,\quad\text{when }\quad q<0  \;. \end{cases} \\
      \eea

      \textbf{The terms $ \Lie_{Z^I}  \Big[   A^4  + O \big(  h  \cdot  A^4  \big) \Big]     $:}\\

We have
 \bea
 \notag
|   \Lie_{Z^I} A   |^4 &\leq& \begin{cases}  c (\gamma) \cdot C ( |I| ) \cdot E ( |I| + 2) \cdot  \frac{\eps}{ (1+ t + | q | )^{4-4\delta }  (1+| q |   )^{4\ga}}  ,\quad\text{when }\quad q>0,\\
\notag
    C ( |I| ) \cdot E ( |I| + 2) \cdot  \frac{\eps \cdot (1+| q |   )^2}{ (1+ t + | q | )^{4-4\delta }  }   , \,\quad\text{when }\quad q<0 . \end{cases} \\ 
    \eea
     Under the bootstrap assumption, we have $|   \Lie_{Z^I} h   | \leq c (\gamma)  \cdot c (\delta) \cdot C ( |I| )  \cdot E ( |I| + 2)$, therefore,
     
 \bea
 \notag
&&|  \Lie_{Z^I}  \Big[   A^4  + O \big(  h  \cdot  A^4  \big) \Big]   | \\
\notag
 &\leq& \begin{cases} c (\delta) \cdot   c (\gamma) \cdot C ( |I| ) \cdot E ( |I| + 2) \cdot  \frac{\eps}{ (1+ t + | q | )^{4-4\delta }  (1+| q |   )^{4\ga}}  ,\quad\text{when }\quad q>0,\\
\notag
    C ( |I| ) \cdot E ( |I| + 2) \cdot  \frac{\eps \cdot (1+| q |   )^2}{ (1+ t + | q | )^{4-4\delta }  }   , \,\quad\text{when }\quad q<0 . \end{cases} \\ 
    \eea
    
       \textbf{The final term $| \Lie_{Z^I}  \Big(  g^{\alpha\beta} \derm_\alpha \derm_\beta h_{ {\cal T} {\cal U}}  \Big) |$ :}\\
       
       Putting all together, we obtain
        \beaa
\notag
 && | \Lie_{Z^I}  \Big(  g^{\alpha\beta} \derm_\alpha \derm_\beta h_{ {\cal T} {\cal U}} \Big)  |   \\
 &\leq&  \begin{cases}  c (\delta) \cdot c (\gamma) \cdot C ( |I| ) \cdot E ( |I| + 3)  \cdot \frac{\eps }{(1+t+|q|)^{3-3\delta} (1+|q|)^{1+2\de}},\quad\text{when }\quad q>0  \;  ,\\
       \notag
       C ( |I| ) \cdot E ( |I| + 3)  \cdot \frac{\eps  }{(1+t+|q|)^{3-2\delta} }  \,\quad\text{when }\quad q<0 \; .\end{cases} \\
      && + \begin{cases}  c (\delta) \cdot c (\gamma) \cdot C ( |I| ) \cdot E ( |I| + 3)  \cdot \frac{\eps }{(1+t+|q|)^{3-3\delta} (1+|q|)^{1+2\de}},\quad\text{when }\quad q>0  \;  ,\\
       \notag
       C ( |I| ) \cdot E ( |I| + 3)  \cdot \frac{\eps  }{(1+t+|q|)^{3-2\delta} }  \,\quad\text{when }\quad q<0 \; .\end{cases} \\
  && + \begin{cases} c (\delta) \cdot c (\gamma) \cdot  C ( |I| ) \cdot E ( |I| +  2)  \cdot \frac{\eps }{(1+t+|q|)^{3-3\delta} (1+|q|)^{1+3\ga}},\quad\text{when }\quad q>0 \; ,\\
       C ( |I| ) \cdot E ( |I| + 2 )  \cdot \frac{\eps \cdot (1+| q |   )^{\frac{1}{2} } }{(1+t+|q|)^{3-3\delta}  }  \,\quad\text{when }\quad q<0  \;. \end{cases} \\
       \notag
 &&+ \begin{cases} c (\delta) \cdot   c (\gamma) \cdot C ( |I| ) \cdot E ( |I| + 2) \cdot  \frac{\eps}{ (1+ t + | q | )^{4-4\delta }  (1+| q |   )^{4\ga}}  ,\quad\text{when }\quad q>0 \; ,\\
\notag
    C ( |I| ) \cdot E ( |I| + 2) \cdot  \frac{\eps \cdot (1+| q |   )^2}{ (1+ t + | q | )^{4-4\delta }  }   , \,\quad\text{when }\quad q<0 \; . \end{cases} \\ 
    \eeaa
    However, examining the last term in the exterior, we have
    \beaa
 && c (\delta) \cdot   c (\gamma) \cdot C ( |I| ) \cdot E ( |I| + 2) \cdot  \frac{\eps}{ (1+ t + | q | )^{4-4\delta }  (1+| q |   )^{4\ga}} \\
 &\leq& c (\delta) \cdot   c (\gamma) \cdot C ( |I| ) \cdot E ( |I| + 2) \cdot  \frac{\eps}{ (1+ t + | q | )^{3-3\delta }  (1+| q |   )^{1+4\ga - \de}} \\
 &\leq& c (\delta) \cdot   c (\gamma) \cdot C ( |I| ) \cdot E ( |I| + 2) \cdot  \frac{\eps}{ (1+ t + | q | )^{3-3\delta }  (1+| q |   )^{1+(\ga - \de) + 3\ga}} \\
  &\leq& c (\delta) \cdot   c (\gamma) \cdot C ( |I| ) \cdot E ( |I| + 2) \cdot  \frac{\eps}{ (1+ t + | q | )^{3-3\delta }  (1+| q |   )^{1+2\de}} \\
  && \text{(given that $\ga > \de$).}
    \eeaa

Thus, for $\ga > \de$, $0 \leq \de \leq 1$,
        \bea\label{aprioridecayestimateontauUcompoforfullmetrich}
\notag
 && | \Lie_{Z^I}  \Big( g^{\alpha\beta} \derm_\alpha \derm_\beta h_{ {\cal T} {\cal U}} \Big)  |   \\
 &\leq&  \begin{cases}  c (\delta) \cdot c (\gamma) \cdot C ( |I| ) \cdot E ( |I| + 3)  \cdot \frac{\eps }{(1+t+|q|)^{3-3\delta} (1+|q|)^{1+2\de}},\quad\text{when }\quad q>0  \;  ,\\
       \notag
       C ( |I| ) \cdot E ( |I| + 3)  \cdot \frac{\eps  \cdot (1+| q |   )^2 }{(1+t+|q|)^{3-3\delta} }  \,\quad\text{when }\quad q<0 \; .\end{cases} \\
        \eea
Thus, we obtain
  \bea
\notag
&& (1+t ) \cdot | \Lie_{Z^I} \Big( g^{\alpha\beta} \derm_\alpha \derm_\beta h_{ {\cal T} {\cal U}} \Big) |  \\
\notag
&\les&  \begin{cases}  c (\delta) \cdot c (\gamma) \cdot C ( |I| ) \cdot E ( |I| + 3)  \cdot \frac{\eps }{(1+t+|q|)^{2-3\delta} (1+|q|)^{1+2\de}},\quad\text{when }\quad q>0  \;  ,\\
       \notag
       C ( |I| ) \cdot E ( |I| + 3)  \cdot \frac{\eps \cdot (1+| q |   )^2  }{(1+t+|q|)^{2-3\delta} }  \,\quad\text{when }\quad q<0 \; .\end{cases} \\
        \eea

\end{proof}

\begin{lemma}\label{upgradedestimatesongoodcomponnentforh1andh0}
In the Lorenz and harmonic gauges, the Einstein-Yang-Mills fields satisfy in the exterior region $\overline{C}$, for $M\leq \eps$\,, for $\ga > \de $\,, and $0 \leq \delta \leq \frac{1}{4}$\,, 
 \bea
 \,|\derm h^1_{ {\cal T} {\cal U}} (t,x)|  &\les&  C(q_0) \cdot   c (\gamma)  \cdot c (\delta)   \cdot E (4) \cdot \frac{\eps }{ (1+t+|q|)} \; ,
\eea
we have
 \bea\label{theupgradedestimateongoodcomponenth}
 \,|\derm h_{ {\cal T} {\cal U}} (t,x)|  &\les&  C(q_0) \cdot   c (\gamma)  \cdot c (\delta)   \cdot E (4) \cdot \frac{\eps }{ (1+t+|q|)} \; . 
\eea
\end{lemma}

\begin{proof}

Given the linearity of the wave operator, we have
\beaa
 g^{\la\mu} \derm_{\la}   \derm_{\mu}    h^1_{ {\cal T} {\cal U}} =  g^{\la\mu} \derm_{\la}   \derm_{\mu}    h_{ {\cal T} {\cal U}} - g^{\la\mu} \derm_{\la}   \derm_{\mu}    h^0_{ {\cal T} {\cal U}} \; .
\eeaa
However, we have from Lemma \ref{estimateonthesourcetermsforhzerothesphericallsymmtrpart},
\beaa
(1+t) \cdot |  g^{\la\mu} \derm_{\la}   \derm_{\mu}    h^0  | &\les&c (\gamma) \cdot  E (  2)  \cdot \frac{\eps \cdot (1+t)  }{(1+t+|q|)^{3} } \\
&\les&c (\gamma) \cdot  E (  2)  \cdot \frac{\eps  }{(1+t+|q|)^{2} }\; ,
\eeaa
and in the exterior, we have already estimated in Lemma \ref{decayestimateonsourcetermsforgoodcomponentofmetrichtauU}, that for $\eps \leq 1$, and $0 \leq \de \leq 1$,
      \bea
\notag
&& (1+t ) \cdot | g^{\alpha\beta} \derm_\alpha \derm_\beta h_{ {\cal T} {\cal U}} |_{L^\infty(\overline{D}_\tau)}  \\
\notag
&\les& C(q_0) \cdot  c (\delta) \cdot c (\gamma) \cdot E (  3)  \cdot \frac{\eps }{(1+t )^{2-3\delta} (1+|q|)^{1+2\de} } \\
       \notag
        &\les& C(q_0) \cdot  c (\delta) \cdot c (\gamma) \cdot E (  3)  \cdot \frac{\eps }{(1+t )^{2-3\delta} } \; . \\
    \eea
We therefore get
          \bea
\notag
(1+t ) \cdot | g^{\alpha\beta} \derm_\alpha \derm_\beta h^1_{ {\cal T} {\cal U}} |_{L^\infty(\overline{D}_\tau)}  &\les&
       C(q_0) \cdot  c (\delta) \cdot c (\gamma) \cdot E (  3)  \cdot \frac{\eps }{(1+t )^{2-3\delta}  } \; . \\
    \eea

    Finally, we get
\beaa
&& (1+t+|x|) \cdot \,|\derm h^1_{ {\cal T} {\cal U}} (t,x)| \\
 &\les& \!\sup_{0\leq \tau\leq t} \sum_{|I|\leq 1}\|\, \Lie_{Z^I}\! h^1_{ {\cal T} {\cal U}} (\tau,\cdot)\|_{L^\infty (\Sigma^{ext}_{\tau} ) }\\
&& + \int_0^t\Big( (1+\tau) \cdot \| \, g^{\la\mu} \derm_{\la}   \derm_{\mu}   h^1_{ {\cal T} {\cal U}} (\tau,\cdot)  \|_{L^\infty(\overline{D}_\tau)} +\sum_{|I|\leq 2} (1+\tau)^{-1}  \cdot  \|  \Lie_{Z^I} h^1_{ {\cal T} {\cal U}} (\tau,\cdot)\|_{L^\infty(\overline{D}_\tau)}\Big)\, d\tau \\
&\les&   c (\gamma)  \cdot c (\delta) \cdot  E (3) \cdot \eps  +  \int_0^t  C(q_0) \cdot c (\delta) \cdot c (\gamma) \cdot E (  3)  \cdot \frac{\eps }{(1+\tau )^{2-3\delta}  } d\tau \\
&& +   \int_0^t   c (\gamma)  \cdot c (\delta)  \cdot E ( 4) \cdot \frac{\eps}{(1+\tau)^{\frac{3}{2}-\delta}  }   d\tau \\
&\les&   c (\gamma)  \cdot c (\delta) \cdot  E (3) \cdot \eps  +  \int_0^t  C(q_0) \cdot c (\delta) \cdot c (\gamma) \cdot E (  3)  \cdot \frac{\eps }{ (1+\tau )^{\frac{3}{2}-\delta}  \cdot (1+\tau )^{\frac{1}{2}-2\delta}  } d\tau \\
&& +   \int_0^t   c (\gamma)  \cdot c (\delta)  \cdot E ( 4) \cdot \frac{\eps}{(1+\tau)^{\frac{3}{2}-\delta}  }   d\tau \\
&& \text{(we choose here again $\delta < \frac{1}{4} $)}  \\
&\les&   C(q_0) \cdot c (\gamma)  \cdot c (\delta)   \cdot E (4) \cdot \eps \big[  \frac{1}{(1+\tau)^{\frac{1}{2} -\delta} } \big]^{\infty}_{0} \\
&& \text{(for $\delta \leq \frac{1}{4} < \frac{1}{2}$)}  .
\eeaa
Hence, for $\delta \leq \frac{1}{4}$, we obtain
\beaa
 (1+t+|x|) \cdot \,|\derm h^1_{ {\cal T} {\cal U}} (t,x)|  &\les&  C(q_0) \cdot  c (\gamma)  \cdot c (\delta)   \cdot E (4) \cdot \eps \; ,
\eeaa
 which leads to 
 \bea
 \,|\derm h^1_{ {\cal T} {\cal U}} (t,x)| &\les&   C(q_0) \cdot c (\gamma)  \cdot c (\delta)   \cdot E (4) \cdot \frac{\eps }{ (1+t+|x|)} \;  . 
\eea

Since from Lemma \ref{Liederivativesofsphericalsymmetricpart}, we get that for $M \leq \eps$\,, 
\beaa
 \,|\derm h_{ {\cal T} {\cal U}} |  &\les& |\derm h^0_{ {\cal T} {\cal U}} |  + |\derm h^1_{ {\cal T} {\cal U}} | \\
 &\leq&\frac{\eps }{(1+t+|q|)^{2} } +  C(q_0) \cdot   c (\gamma)  \cdot c (\delta)   \cdot E (4) \cdot \frac{\eps }{ (1+t+|q|)} \\
 &\leq& C(q_0) \cdot   c (\gamma)  \cdot c (\delta)   \cdot E (4) \cdot \frac{\eps }{ (1+t+|q|)} \; ,
\eeaa
we obtain the result \eqref{theupgradedestimateongoodcomponenth}.
\end{proof}

\section{Pointwise estimates for “bad” components of the metric}

  \begin{lemma}\label{decayestimateinexterioronsourcetermforhbadcomponent}
  
  We have in the exterior region $\overline{C}$\,, for $M\leq \eps$\,, for $\ga > \de $\,, and $0 \leq \de \leq \frac{1}{4}$\,,
  \bea
\notag
(1+t ) \cdot | g^{\alpha\beta} \derm_\alpha \derm_\beta h_{{\underline{L}}{\underline{L}}} |  &\les&C(q_0) \cdot  c (\gamma)  \cdot c (\delta)   \cdot E (4) \cdot \frac{\eps }{ (1+t+|q|) } \; .
\eea

\end{lemma}
\begin{proof}
  
In the Lorenz gauge and in wave coordinates, from Lemma \ref{structureofthesourcetermsofthewaveoperatoronAandh}, we have
    \bea
\notag
&&  | g^{\alpha\beta} \derm_\alpha \derm_\beta h_{{\underline{L}} {\underline{L}} } |    \\
\notag
            &\les&  | P(\derm_{\underline{L}} h,\derm_{\underline{L}} h) | + | Q_{{\underline{L}}{\underline{L}}}(\derm h,\derm h)  | + | G_{{\underline{L}}{\underline{L}}}(h)(\derm h,\derm h)  | \\
 \notag
    && +  | \rderm A   | \cdot  |\derm  A  | +  | A |^2  \cdot | \derm A   |   +  |A |^4 \\
    \notag
     && + O \big(h \cdot  (\derm A)^2 \big)   + O \big(  h  \cdot  A^2 \cdot \derm A \big)     + O \big(  h   \cdot  A^4 \big)   \\
 &&+   |\derm A_{e_{a}}  |^2 \; .
\eea

\textbf{The term $  | P(\derm_{\underline{L}} h,\derm_{\underline{L}} h) | + | Q_{{\underline{L}}{\underline{L}}}(\derm h,\derm h)  | + | G_{{\underline{L}}{\underline{L}}}(h)(\derm h,\derm h)  |  $:}\\

  Based on Corollary 7.2 in \cite{LR10}, we have
  
  \bea\label{estimateonthetermsgeneratedfromtheEinsteinvacuumequationswhichhaveanullstructureplusasquareofderivativeofgoodcomponents}
  \notag
&& | P(\derm_{\underline{L}} h,\derm_{\underline{L}} h) | + | Q_{{\underline{L}}{\underline{L}}}(\pa h,\pa h)  | + | G_{{\underline{L}}{\underline{L}}}(h)(\pa h,\pa h)  | \\
\notag
 &\les& |\derm h_{\cal TU} |^2+ |\rderm h| \cdot |\derm h| + |h| \cdot |\derm h|^2 \;. \\
\eea

Thus, we have the same estimate we made on $  | P(\derm_{\cal T} h,\derm_{\cal U} h) | + | Q_{{\cal T} {\cal U}}(\derm h,\derm h)  | + | G_{{\cal T} {\cal U} }(h)(\derm h,\derm h)  |$, except with the new term $|\derm h_{\cal TU} |^2$. Consequently, for the term $| g^{\alpha\beta} \derm_\alpha \derm_\beta h_{{\underline{L}} {\underline{L}} } |$, we have the same estimate as for  $| g^{\alpha\beta} \derm_\alpha \derm_\beta h_{ {\cal T} {\cal U}} |$ plus the additional $|\derm h_{ {\cal T} {\cal U}} |^2 + |\derm A_{{\cal T}} |^2$. Thus, in the exterior, for $\ga > \de $, and $0 \leq \de \leq \frac{1}{4}$, we have
          \bea\label{decayestimateongoodtermsplusthebadtermwithoutdecayinsourceforh}
\notag
 &&  | g^{\alpha\beta} \derm_\alpha \derm_\beta h_{{\underline{L}} {\underline{L}} } |    \\
 \notag
 &\leq&  \begin{cases}  c (\delta) \cdot c (\gamma) \cdot E (  3)  \cdot \frac{\eps }{(1+t+|q|)^{3-3\delta} (1+|q|)^{1+2\de}},\quad\text{when }\quad q>0  \;  ,\\
       \notag
      E ( 3)  \cdot \frac{\eps  \cdot (1+| q |   )^2 }{(1+t+|q|)^{3-3\delta} }  \,\quad\text{when }\quad q<0 \; \end{cases} \\
       \notag
       && + |\derm h_{ {\cal T} {\cal U}} |^2 + |\derm A_{{\cal T}} |^2 \; .\\
       \eea
       However, we have already estimated, in Lemma \ref{upgradedestimatesongoodcomponnentforh1andh0}, the following term, for $0 \leq \de < \frac{1}{4}$, by
 \beaa
 \notag
 \,|\derm h_{ {\cal T} {\cal U}} |^2 &\les&   C(q_0) \cdot c (\gamma)  \cdot c (\delta)   \cdot E (4) \cdot \frac{\eps }{ (1+t+|q|)^2} \; .
\eeaa

\textbf{The term $   |\derm A_{e_{a}}  |^2 $:}\\

We have already estimated, in Lemma \ref{estimatefortangentialcomponentspotential}, for $0 \leq \de < \frac{1}{4}$,
 \bea
 \,|\derm A_{{\cal T}} |^2  &\les&  C(q_0) \cdot  c (\gamma)  \cdot c (\delta)   \cdot E (4) \cdot \frac{\eps }{ (1+t+|q|)^2} \; .
\eea

\textbf{The final term $   | g^{\alpha\beta} \derm_\alpha \derm_\beta h_{{\underline{L}} {\underline{L}} } |$:}\\

Finally, the term $| g^{\alpha\beta} \derm_\alpha \derm_\beta h_{{\underline{L}} {\underline{L}} } |$ has the same estimate as for  $| g^{\alpha\beta} \derm_\alpha \derm_\beta h_{ {\cal T} {\cal U}} |$ plus the additional  $C(q_0) \cdot  c (\gamma)  \cdot c (\delta)   \cdot E (4) \cdot \frac{\eps }{ (1+t+|q|)^2}$. Thus, in the exterior, for $\ga > \de $, and $0 \leq \de \leq \frac{1}{4}$, we have

            \bea
\notag
 &&  | g^{\alpha\beta} \derm_\alpha \derm_\beta h_{{\underline{L}} {\underline{L}} } |    \\
 \notag     
        &\leq&  \begin{cases}  c (\delta) \cdot c (\gamma) \cdot E (  3)  \cdot \frac{\eps }{(1+t+|q|)^{3-3\delta} (1+|q|)^{1+2\de}},\quad\text{when }\quad q>0  \;  ,\\
       \notag
      E (  3)  \cdot \frac{\eps  \cdot (1+| q |   )^2 }{(1+t+|q|)^{3-3\delta} }  \,\quad\text{when }\quad q<0 \; \end{cases} \\
       \notag
       && + C(q_0) \cdot  c (\gamma)  \cdot c (\delta)   \cdot E (4) \cdot \frac{\eps }{ (1+t + |q| )^2} \\
       \notag
 &\leq& C(q_0) \cdot  c (\gamma)  \cdot c (\delta)   \cdot E (4) \cdot \frac{\eps }{ (1+t +|q|)^2} \; .
        \eea

\end{proof}

\begin{lemma}
In the Lorenz and harmonic gauges, the Einstein-Yang-Mills fields satisfy in the exterior region $\overline{C}$\,, for $M\leq \eps$\,, for $\ga > \de $\,, and $\delta \leq \frac{1}{4}$\,, 
\beaa
 \,|\derm h^1_{{\underline{L}}{\underline{L}}} (t,x)|  &\les&  C(q_0) \cdot   c (\gamma)  \cdot c (\delta)   \cdot E (4) \cdot \eps \cdot \frac{ ( 1 +  \log (1+t) ) }{ (1+t)} .
\eeaa

\end{lemma}
\begin{proof}
Since 
\beaa
 | g^{\alpha\beta} \derm_\alpha \derm_\beta h^1_{{\underline{L}}{\underline{L}}} |  \leq  | g^{\alpha\beta} \derm_\alpha \derm_\beta h_{{\underline{L}}{\underline{L}}} |  +  | g^{\alpha\beta} \derm_\alpha \derm_\beta h^0_{{\underline{L}}{\underline{L}}} | \; ,
\eeaa
and 
\beaa
(1+t) \cdot |  g^{\la\mu} \derm_{\la}   \derm_{\mu}    h^0  | &\les&c (\gamma) \cdot  E (  2)  \cdot \frac{\eps  }{(1+t+|q|)^{2} }\; ,
\eeaa
and since we have shown, in Lemma \ref{decayestimateinexterioronsourcetermforhbadcomponent}, that in the exterior 
 \bea
\notag
(1+t ) \cdot | g^{\alpha\beta} \derm_\alpha \derm_\beta h_{{\underline{L}}{\underline{L}}} |  &\les&C(q_0) \cdot  c (\gamma)  \cdot c (\delta)   \cdot E (4) \cdot \frac{\eps }{ (1+t + |q|) } \; ,
\eea
we get that
  \bea
\notag
(1+t ) \cdot | g^{\alpha\beta} \derm_\alpha \derm_\beta h^1_{{\underline{L}}{\underline{L}}} |  &\les&C(q_0) \cdot  c (\gamma)  \cdot c (\delta)   \cdot E (4) \cdot \frac{\eps }{ (1+t + |q|) }  \\
\notag
 &\les&C(q_0) \cdot  c (\gamma)  \cdot c (\delta)   \cdot E (4) \cdot \frac{\eps }{ (1+t ) } \; . \\
\eea

We can now apply the estimate to get

\beaa
&& (1+t+|x|) \cdot \,|\derm h^1_{{\underline{L}}{\underline{L}}} (t,x)| \\
 &\les& \!\sup_{0\leq \tau\leq t} \sum_{|I|\leq 1}\|\, \Lie_{Z^I} \! h^1 (\tau,\cdot)\|_{L^\infty (\Sigma^{ext}_{\tau} ) }\\
&& + \int_0^t\Big( (1+\tau)  \cdot \| \, g^{\la\mu} \derm_{\la}   \derm_{\mu}   h^1_{{\underline{L}}{\underline{L}}} (\tau,\cdot)  \|_{L^\infty(\overline{D}_\tau)} +\sum_{|I|\leq 2} (1+\tau)^{-1}  \cdot \|  \Lie_{Z^I} h^1(\tau,\cdot)\|_{L^\infty(\overline{D}_\tau)}\Big)\, d\tau \\
&\les&   c (\gamma)  \cdot c (\delta) \cdot  E (3) \cdot \eps   +  \int_0^t C(q_0) \cdot   c (\delta) \cdot c (\gamma) \cdot E (  4)  \cdot \frac{\eps }{(1+\tau )  } d\tau  \\
&& +   \int_0^t   c (\gamma)  \cdot c (\delta)  \cdot E ( 4) \cdot \frac{\eps}{(1+\tau)^{\frac{3}{2}-\delta}  }   d\tau \\
&\les&   c (\gamma)  \cdot c (\delta)   \cdot E (4) \cdot \eps \Big[  \frac{1}{(1+\tau)^{\frac{1}{2} -\delta} } \Big]^{t}_{0} + C(q_0)  \cdot c (\gamma)  \cdot c (\delta)   \cdot E (4) \cdot \big[  \log (1+\tau)  \big]^{t}_{0} \\
&& \text{(where $\delta \leq \frac{1}{4} < \frac{1}{2}$)} \\
&\les&   C(q_0) \cdot  c (\gamma)  \cdot c (\delta)   \cdot E (4) \cdot \eps \cdot \Big( 1 +  \log (1+t)   \Big) 
\eeaa
Hence, we obtain
\beaa
 (1+t+|x|) \cdot \,|\derm h^1_{{\underline{L}}{\underline{L}}} (t,x)|  &\les& C(q_0) \cdot   c (\gamma)  \cdot c (\delta)   \cdot E (4)  \cdot \Big( 1 +  \log (1+t)   \Big) \, ,
\eeaa
 which leads to 
 \bea
 \notag
 \,|\derm h^1_{{\underline{L}}{\underline{L}}} (t,x)| &\les&  C(q_0) \cdot   c (\gamma)  \cdot c (\delta)   \cdot E (4) \cdot \eps \cdot \frac{ ( 1 +  \log (1+t) )}{ (1+t+|x|)} \\
 \notag
 &\les&   C(q_0) \cdot  c (\gamma)  \cdot c (\delta)   \cdot E (4) \cdot \eps \cdot \frac{ ( 1 +  \log (1+t) ) }{ (1+t)} \; .\\
\eea
\end{proof}

\section{Improved pointwise estimates with better decay in spatial infinity}
        
    From now on, we are in fact going often to make estimates for $\de < \frac{1}{4}$.

\subsection{Estimate for the “bad” component of the Einstein-Yang-Mills potential}\

For this, we will apply again the estimate of Lindblad-Rodnianski in Corollary 7.2 in \cite{LR10} (see \eqref{LinfinitynormestimateongradientderivedbyLondbladRodnianski}), for $ g^{\la\mu} \derm_{\la}   \derm_{\mu}    \Lie_{Z^J} A_{\cal T} $, but this time with a non-zero weight $\varpi(q) $, i.e. where $1 + \gamma' >  0$ and non-zero, and with
\bea
\varpi :=\varpi(q) :=\begin{cases}
(1+|q|)^{1+\gamma^\prime},\quad\text{when }\quad q>0 \; ,\\
     1 \,\quad\text{when }\quad q<0 \; .\end{cases}
\eea

We will use the weighted version of the estimate of Lindblad-Rodnianski to improve the estimate for  $ | \derm A_{e_{a}} | $, and then, we will use this upgraded estimate to get an estimate on $ | \derm A_{\underline{L}}  | $ .

The weighted version of the Lindblad-Rodnianski (see \eqref{LinfinitynormestimateongradientderivedbyLondbladRodnianski}) gives that in the exterior, for all $V\in\{L,\Lb,A,B\}$,
\beaa
\notag
&& (1+t+|x|) \cdot |\varpi(q)\derm ( \Lie_{Z^J} A)_V (t,x)| \\
 &\les& \!\sup_{0\leq \tau\leq t} \sum_{|I|\leq 1}\|\,\varpi(q) ( \Lie_{Z^I} \! \Lie_{Z^J} A) (\tau,\cdot)\|_{L^\infty (\Sigma^{ext}_{\tau} ) }\\
\notag
&& + \int_0^t\Big( \varepsilon^\prime \cdot (1 + \gamma’)  \cdot \|\varpi(q) \derm  ( \Lie_{Z^J} A)_V (\tau,\cdot) \|_{L^\infty (\Sigma^{ext}_{\tau} ) }  \\
\notag
&& +(1+\tau) \cdot \| \varpi(q)
 g^{\la\mu} \derm_{\la}   \derm_{\mu}  ( \Lie_{Z^J} A)_V (\tau,\cdot) \|_{L^\infty(\overline{D}_\tau)} \\
 \notag
&& +\sum_{|I|\leq 2} (1+\tau)^{-1} \cdot \| \varpi(q) \Lie_{Z^I}  \Lie_{Z^J} A (\tau,\cdot)\|_{L^\infty(\overline{D}_\tau)}\Big)\, d\tau. \\
\eeaa

Now, we are going to study the terms in the right hand side of this estimate one by one.

\begin{lemma}
For $\gamma^\prime$ such that  $-1 \leq \gamma^\prime \leq \gamma - \delta$, we have
 \bea
 \notag
 \!\sup_{0\leq \tau\leq t} \sum_{|I|\leq 1}\|\,\varpi(q) \Lie_{Z^I} \! \Lie_{Z^J} A (\tau,\cdot)\|_{L^\infty (\Sigma^{ext}_{\tau} ) } &\leq&  C ( |J| ) \cdot E ( |J| + 3)  \cdot \eps \; ,\\ 
      \eea
      and
       \bea
 \notag
 \!\sup_{0\leq \tau\leq t} \sum_{|I|\leq 1}\|\,\varpi(q) \Lie_{Z^I} \! \Lie_{Z^J} h^1 (\tau,\cdot)\|_{L^\infty (\Sigma^{ext}_{\tau} ) } &\leq&  C ( |J| ) \cdot E ( |J| + 3)  \cdot \eps \; .\\ 
      \eea
\end{lemma}

\begin{proof}
In fact, the term
\beaa
 \!\sup_{0\leq \tau\leq t} \sum_{|I|\leq 1}\|\,\varpi(q) \Lie_{Z^I} \! \Lie_{Z^J} A (\tau,\cdot)\|_{L^\infty (\Sigma^{ext}_{\tau} ) } &=&  \!\sup_{0\leq \tau\leq t} \sum_{|I|\leq |J| + 1 }\|\,\varpi(q) \Lie_{Z^I}  A (\tau,\cdot)\|_{L^\infty (\Sigma^{ext}_{\tau} ) } \, .
\eeaa
However, by bootstrap assumption, we have
                 \bea
 \notag
|   \Lie_{Z^I} A (t,x)  | &\leq& \begin{cases} c (\delta) \cdot c (\gamma) \cdot C ( |I| ) \cdot E ( |I| + 2) \cdot  \frac{\eps}{ (1+ t + | q | )^{1-\delta }  (1+| q |   )^{\gamma}}  ,\quad\text{when }\quad q>0 \; ,\\
\notag
    C ( |I| ) \cdot E ( |I| + 2) \cdot  \frac{\eps \cdot (1+|q|)^{\frac{1}{2}} }{ (1+ t + | q | )^{1-\delta }  }   , \,\quad\text{when }\quad q<0 \; .  \end{cases} \\ 
    \eea
     Thus,
  \beaa
\varpi(q) \cdot |  \Lie_{Z^I} A (t,x)  |    &\leq& \begin{cases} c (\delta) \cdot c (\gamma) \cdot C ( |I| ) \cdot E ( |I| + 2) \cdot  \frac{\eps \cdot (1+|q|)^{1+\gamma^\prime} }{ (1+ t + | q | )^{1-\delta }  (1+| q |   )^{\gamma}}  ,\quad\text{when }\quad q>0 \; ,\\
\notag
    C ( |I| ) \cdot E ( |I| + 2) \cdot  \frac{\eps \cdot (1+|q|)^{\frac{1}{2}}  }{ (1+ t + | q | )^{1-\delta }  }   , \,\quad\text{when }\quad q<0 \;.  \end{cases} \\ 
      &\leq& \begin{cases} c (\delta) \cdot c (\gamma) \cdot C ( |I| ) \cdot E ( |I| + 2) \cdot  \frac{\eps \cdot (1+|q|)^{1+\gamma^\prime - \ga} }{ (1+ t + | q | )^{1-\delta } }  ,\quad\text{when }\quad q>0 \;,\\
\notag
    C ( |I| ) \cdot E ( |I| + 2) \cdot  \frac{\eps  }{ (1+ t + | q | )^{\frac{1}{2} -\delta }  }   , \,\quad\text{when }\quad q<0 \; .  \end{cases} \\ 
    \eeaa
     Hence, if $1+\gamma^\prime - \gamma \leq 1 - \delta$ or said differently, if $-1 \leq \gamma^\prime \leq \gamma - \delta$, then 
        \beaa
\|\,\varpi(q) \Lie_{Z^I}  A (\tau,\cdot)\|_{L^\infty (\Sigma^{ext}_{\tau} ) }    &\leq&   c (\delta) \cdot c (\gamma) \cdot C ( |I| ) \cdot E ( |I| + 2)  \; .
      \eeaa

 Thus,
   \bea
   \notag
 \!\sup_{0\leq \tau\leq t} \sum_{|I|\leq 1}\|\,\varpi(q) \Lie_{Z^I} \! \Lie_{Z^J} A (\tau,\cdot)\|_{L^\infty (\Sigma^{ext}_{\tau} ) } &\leq&  C ( |J| ) \cdot E ( |J| + 3)  \cdot \eps  \; . \\
      \eea
 
Since these estimates for $A$ hold also for $h^1$, we get the same conclusion for $h^1$.
 \end{proof}

\begin{lemma}
We have
  \bea
 \notag
  &&\int_0^t  \varepsilon^\prime \cdot   (1+\gamma' )  \cdot  \|\varpi(q)  \cdot  \derm ( \Lie_{Z^J} A)_V (\tau,\cdot) \|_{L^\infty (\Sigma^{ext}_{\tau} ) } d \tau \\
  \notag
          &\les&    c (\gamma^\prime)  \cdot c (\gamma)  \cdot c (\delta)  \cdot E ( 3)  \cdot \eps  \cdot  \int_0^t  \frac{1}{(1+\tau)} \cdot (1+\tau+|x|) \cdot \|\varpi(q)  \cdot \derm  \Lie_{Z^J} A_V (\tau,\cdot) \|_{L^\infty (\Sigma^{ext}_{\tau} ) } d \tau  \; , \\
 \eea
 and
   \bea
 \notag
  &&\int_0^t  \varepsilon^\prime \cdot   (1+\gamma' )  \cdot  \|\varpi(q)  \cdot  \derm ( \Lie_{Z^J} h^1)_{UV} (\tau,\cdot) \|_{L^\infty (\Sigma^{ext}_{\tau} ) } d \tau \\
  \notag
          &\les&    c (\gamma^\prime)  \cdot c (\gamma)  \cdot c (\delta)  \cdot E ( 3)  \cdot \eps  \cdot  \int_0^t  \frac{1}{(1+\tau)} \cdot (1+\tau+|x|) \cdot \|\varpi(q)  \cdot \derm  \Lie_{Z^J} h^1_{UV} (\tau,\cdot) \|_{L^\infty (\Sigma^{ext}_{\tau} ) } d \tau  \; , \\
 \eea
where $(1+\gamma' ) > 0 \,$.
 \end{lemma}
 
 \begin{proof}
 Recall that $\varepsilon^\prime$ was the constant that appears on the conditions on $H$ so that the pointwise estimate on the gradient would hold true. We showed that for some constant $C$, we have
 \bea
\varepsilon^\prime \cdot  (1+\gamma' ) = 4 C  \cdot  c (\gamma)  \cdot c (\delta)  \cdot E ( 3)  \cdot \eps \cdot (1+\gamma' ) \; .
\eea

Thus,
 \bea
 \notag
  &&\int_0^t  \varepsilon^\prime \cdot (1+\gamma' )  \cdot  \|\varpi(q)  \cdot  \derm ( \Lie_{Z^J} A)_V (\tau,\cdot) \|_{L^\infty(\Sigma^{ext}_{\tau} ) } d \tau \\
   \notag
   &\les&   \int_0^t    c (\gamma)  \cdot c (\delta)  \cdot E ( 3)  \cdot \eps \cdot (1+\gamma' )  \cdot  \|\varpi(q)  \cdot  \derm (  \Lie_{Z^J} A)_V (\tau,\cdot) \|_{L^\infty (\Sigma^{ext}_{\tau} ) } d \tau \\
    \notag
      &\les&   \int_0^t  \frac{(1+t+|x|)}{(1+t)} \cdot  c (\gamma)  \cdot c (\delta)  \cdot E ( 3)  \cdot \eps \cdot (1+\gamma' )  \cdot  \|\varpi(q)  \cdot  \derm ( \Lie_{Z^J} A)_V (\tau,\cdot) \|_{L^\infty (\Sigma^{ext}_{\tau} ) } d \tau \\
       \notag
         &\les&   c (\gamma)  \cdot c (\delta)  \cdot E ( 3)  \cdot \eps \cdot (1+\gamma' )  \cdot  \int_0^t  \frac{1}{(1+\tau)} \cdot (1+\tau+|x|) \cdot \|\varpi(q)  \cdot  \derm ( \Lie_{Z^J} A)_V (\tau,\cdot) \|_{L^\infty (\Sigma^{ext}_{\tau} ) } d \tau  \; .\\
 \eea
 These same estimates on $A$ hold for $h^1$.
 \end{proof}

\begin{lemma}
For $\gamma^\prime$ such that $-1 \leq \gamma^\prime < \gamma - \delta$, and $ \delta <    1/2 $, we have
\bea
\notag
&&  \int_0^t \sum_{|I|\leq 2} (1+\tau)^{-1} \cdot \| \varpi(q)  \cdot  \Lie_{Z^I}  \Lie_{Z^J} A (\tau,\cdot)\|_{L^\infty(\overline{D}_\tau)} \, d\tau \\
\notag
  &\leq&  c (\gamma^\prime) \cdot  c (\delta) \cdot c (\gamma) \cdot C ( |J| ) \cdot E ( |J| + 4) \cdot \eps \, , \\
\eea
and
\bea
\notag
&&  \int_0^t \sum_{|I|\leq 2} (1+\tau)^{-1} \cdot \| \varpi(q)  \cdot  \Lie_{Z^I}  \Lie_{Z^J} h^1 (\tau,\cdot)\|_{L^\infty(\overline{D}_\tau)} \, d\tau \\
\notag
  &\leq&  c (\gamma^\prime) \cdot  c (\delta) \cdot c (\gamma) \cdot C ( |J| ) \cdot E ( |J| + 4) \cdot \eps \, . \\
\eea
\end{lemma}
  
  \begin{proof}
  
  Regarding the term $  \int_0^t \sum_{|I|\leq 2} (1+\tau)^{-1} \cdot \| \varpi(q) \cdot \Lie_{Z^I}  \Lie_{Z^J} A (\tau,\cdot)\|_{L^\infty(\overline{D}_\tau)} \, d\tau $, we now choose $-1 \leq \gamma^\prime < \gamma - \delta$ and $ \delta <    1/2 $, then 
 
  \beaa
\varpi(q) \cdot |  \Lie_{Z^I} A (t,x)  |    &\leq& \begin{cases} c (\delta) \cdot c (\gamma) \cdot C ( |I| ) \cdot E ( |I| + 2) \cdot  \frac{\eps \cdot (1+|q|)^{1+\gamma^\prime} }{ (1+ t + | q | )^{1-\delta }  (1+| q |   )^{\gamma}}  ,\quad\text{when }\quad q>0,\\
\notag
    C ( |I| ) \cdot E ( |I| + 2) \cdot  \frac{\eps \cdot (1+|q|)^{\frac{1}{2}}  }{ (1+ t + | q | )^{1-\delta }  }   , \,\quad\text{when }\quad q<0 .  \end{cases} \\ 
     &\leq& \begin{cases} c (\delta) \cdot c (\gamma) \cdot C ( |I| ) \cdot E ( |I| + 2) \cdot  \frac{\eps \cdot (1+t +|q|)^{1+\gamma^\prime - \gamma} }{ (1+ t + | q | )^{1-\delta } }  ,\quad\text{when }\quad q>0,\\
\notag
    C ( |I| ) \cdot E ( |I| + 2) \cdot  \frac{\eps   }{ (1+ t + | q | )^{\frac{1}{2}-\delta }  }   , \,\quad\text{when }\quad q<0 .  \end{cases} 
    \eeaa
    Hence,
    \beaa
   \varpi(q) \cdot |  \Lie_{Z^I} A (t,x)  |         &\leq& \begin{cases} c (\delta) \cdot c (\gamma) \cdot C ( |I| ) \cdot E ( |I| + 2) \cdot  \frac{\eps }{ (1+ t + | q | )^{\gamma-\gamma^\prime-\delta } }  ,\quad\text{when }\quad q>0,\\
\notag
    C ( |I| ) \cdot E ( |I| + 2) \cdot  \frac{\eps }{ (1+ t + | q | )^{ \frac{1}{2} -\delta }  }   , \,\quad\text{when }\quad q<0 .  \end{cases} \\ 
      &\leq& c (\delta) \cdot c (\gamma) \cdot C ( |I| ) \cdot E ( |I| + 2) \cdot \Big( \frac{\eps }{ (1+ t + | q | )^{\gamma-\gamma^\prime-\delta } }  + \frac{\eps }{ (1+ t + | q | )^{ \frac{1}{2} -\delta }  }  \Big)  \; .
    \eeaa
    Consequently,
        \bea
        \notag
 && \int_0^t  \frac{1}{(1+\tau)}\cdot   \| \varpi(q)\cdot  \Lie_{Z^I}  \Lie_{Z^J} A (\tau,\cdot)\|_{L^\infty(\overline{D}_\tau)} \, d\tau \\
         \notag
  &\leq&  \int_0^t   c (\delta) \cdot c (\gamma) \cdot C ( |I| +|J| ) \cdot E ( |I| +|J| + 2) \cdot \Big( \frac{\eps }{ (1+ \tau )^{1+ \gamma-\gamma^\prime-\delta } }  +  \frac{\eps  }{ (1+ \tau  )^{\frac{3}{2} -\delta }   } \Big) d\tau \\
          \notag
    &\leq&  c (\gamma^\prime) \cdot  c (\delta) \cdot c (\gamma) \cdot C ( |I| +|J|) \cdot E ( |I| +|J| + 2) \cdot \Big[ \frac{\eps }{ (1+ \tau )^{ \gamma-\gamma^\prime-\delta } }  +  \frac{\eps  }{ (1+ \tau  )^{ \frac{1}{2} -\delta  }  } \Big]^{t}_{0}  \\
    \notag
    &\leq&  c (\gamma^\prime) \cdot  c (\delta) \cdot c (\gamma) \cdot C ( |I| +|J|) \cdot E ( |I| +|J|+ 2) \cdot \eps \, .
    \eea
    Summing over $|I|\leq 2$, we obtain the result.
    Similarly, we get the same estimates for $h^1$.
      \end{proof}

  \begin{lemma}\label{GronwallinequalitiesoncomponentsofgradientofAandgradientofh1withweight}
Let $M \leq \eps$\,. For $\gamma^\prime$ such that $-1 \leq \gamma^\prime < \gamma - \delta$\,, and $ \delta <    1/2 $\,, we have for all $U, V\in \cal U := \{L,\Lb,A,B\}$\,,
     \bea
   \notag
&& (1+t+|x|) \cdot |\varpi(q) \cdot \derm ( \Lie_{Z^J} A)_V (t,x)| \\
   \notag
 &\les&   c (\gamma^\prime) \cdot  c (\delta) \cdot c (\gamma) \cdot C ( |J|  ) \cdot E ( |J| + 4) \cdot \eps \, \\
\notag
&& + c (\gamma^\prime) \cdot c (\gamma)  \cdot c (\delta)  \cdot E ( 3)  \cdot \eps  \cdot  \int_0^t  \frac{1}{(1+\tau)} \cdot (1+\tau+|x|) \cdot \|\varpi(q)  \cdot  \derm  ( \Lie_{Z^J} A)_V (\tau,\cdot) \|_{L^\infty (\Sigma^{ext}_{\tau} )} d \tau \\
       \notag
&& + \int_0^t  (1+\tau) \cdot \| \varpi(q)  \cdot g^{\la\mu} \derm_{\la}   \derm_{\mu}  (  \Lie_{Z^J} A)_V (\tau,\cdot) \|_{L^\infty(\overline{D}_\tau)} d\tau \; ,\\
\eea
and 
     \bea
   \notag
&& (1+t+|x|) \cdot |\varpi(q) \cdot \derm ( \Lie_{Z^J} h^1)_{UV} (t,x)| \\
   \notag
 &\les&   c (\gamma^\prime) \cdot  c (\delta) \cdot c (\gamma) \cdot C ( |J|  ) \cdot E ( |J| + 4) \cdot \eps \, \\
\notag
&& + c (\gamma^\prime) \cdot c (\gamma)  \cdot c (\delta)  \cdot E ( 3)  \cdot \eps  \cdot  \int_0^t  \frac{1}{(1+\tau)} \cdot (1+\tau+|x|) \cdot \|\varpi(q)  \cdot  \derm  ( \Lie_{Z^J} h^1)_{UV} (\tau,\cdot) \|_{L^\infty (\Sigma^{ext}_{\tau} )} d \tau \\
       \notag
&& + \int_0^t  (1+\tau) \cdot \| \varpi(q)  \cdot g^{\la\mu} \derm_{\la}   \derm_{\mu}  (  \Lie_{Z^J} h^1)_{UV} (\tau,\cdot) \|_{L^\infty(\overline{D}_\tau)} d\tau  \; .\\
\eea

  \end{lemma}
  
  \begin{proof}
  
Plugging all the previous estimates together, we get for all $V\in\{L,\Lb,A,B\}$,
\beaa
\notag
&& (1+t+|x|) \cdot |\varpi(q) \cdot \derm ( \Lie_{Z^J} A)_V (t,x)| \\
 &\les& \!\sup_{0\leq \tau\leq t} \sum_{|I|\leq 1}\|\,\varpi(q) \Lie_{Z^I} \! \Lie_{Z^J} A(\tau,\cdot)\|_{L^\infty (\Sigma^{ext}_{\tau} )}\\
\notag
&& + \int_0^t\Big( \varepsilon^\prime \cdot (1+ \gamma^\prime) \cdot \|\varpi(q) \cdot \derm ( \Lie_{Z^J} A)_V (\tau,\cdot) \|_{L^\infty (\Sigma^{ext}_{\tau} ) } \\
\notag
&& +(1+\tau) \cdot \| \varpi(q) \cdot  g^{\la\mu} \derm_{\la}   \derm_{\mu}  ( \Lie_{Z^J} A)_V (\tau,\cdot) \|_{L^\infty(\overline{D}_\tau)} \\
 \notag
&& +\sum_{|I|\leq 2} (1+\tau)^{-1} \cdot \| \varpi(q) \cdot \Lie_{Z^I}  \Lie_{Z^J} A (\tau,\cdot)\|_{L^\infty(\overline{D}_\tau)}\Big)\, d\tau \\
\notag
&\leq&  C ( |J| ) \cdot E ( |J| + 3)  \cdot \eps \\
\notag
&& + c (\gamma^\prime) \cdot c (\gamma)  \cdot c (\delta)  \cdot E ( 3)  \cdot \eps  \cdot  \int_0^t  \frac{1}{(1+\tau)} \cdot (1+\tau+|x|) \cdot \|\varpi(q) \cdot \derm  ( \Lie_{Z^J} A)_V (\tau,\cdot) \|_{L^\infty (\Sigma^{ext}_{\tau} ) } d \tau \\
       \notag
&& + \int_0^t  (1+\tau)\| \varpi(q) \cdot g^{\la\mu} \derm_{\la}   \derm_{\mu} (  \Lie_{Z^J} A)_V (\tau,\cdot) \|_{L^\infty(\overline{D}_\tau)} \\
\notag
 && +  c (\gamma^\prime)  \cdot  c (\delta) \cdot c (\gamma) \cdot C ( |J| ) \cdot E ( |J| + 4) \cdot \eps \, .
\eeaa

We obtain, the same estimate for $h^1$ as for $A$.
\end{proof}

Now, we would like to improve the estimate for  $ | \derm A_{e_{a}} | $ and for this, we will use the weighted version of the estimate of Lindblad-Rodnianski. Then, we will use this upgraded estimate to get an estimate on $ | \derm A_{\underline{L}}  | $ .

  \begin{lemma}\label{weightedestimatefortangentialcomponentspotential}
 
  We have in the exterior region $\overline{C}$\,, for $M\leq \eps$\,, for $\ga > \de $\,, $0 \leq \de < \frac{1}{2}$\, and for
   \bea
-1 &\leq& \gamma^\prime <  \ga - \delta \; ,
  \eea
the following estimate
              \bea
   \notag
&& |\derm \ A_{{\cal T}} (t,x)| \\
   \notag
   &\leq& C(q_0) \cdot c (\gamma^\prime)  \cdot  c (\delta) \cdot c (\gamma)  \cdot E (  4) \cdot \eps \cdot    \frac{ 1}{ (1+t+|q|)^{1-c (\gamma^\prime) \cdot c (\gamma)  \cdot c (\delta)  \cdot E ( 3)  \cdot \eps }  \cdot (1+|q|)^{1+\gamma^\prime} } \; .
\eea

\end{lemma}

  \begin{proof}
  In fact, we have already showed that if $-1 \leq \gamma^\prime < \gamma - \delta$, and $ \delta <    1/2 $, then for all $V\in\{L,\Lb,A,B\}$, we have
     \beaa
   \notag
&& (1+t+|x|) \cdot |\varpi(q) \cdot  \derm \ A_V (t,x)| \\
 &\les&   c (\gamma^\prime)  \cdot  c (\delta) \cdot c (\gamma)  \cdot E (  4) \cdot \eps \, \\
\notag
&& + c (\gamma^\prime)  \cdot c (\gamma)  \cdot c (\delta)  \cdot E ( 3)  \cdot \eps  \cdot  \int_0^t  \frac{1}{(1+\tau)} \cdot (1+\tau+|x|) \cdot \|\varpi(q) \cdot  \derm   A_V (\tau,\cdot) \|_{L^\infty (\Sigma^{ext}_{\tau} )} d \tau \\
       \notag
&& + \int_0^t  (1+\tau) \cdot \| \varpi(q)  \cdot  g^{\la\mu} \derm_{\la}   \derm_{\mu}  A_V (\tau,\cdot) \|_{L^\infty(\overline{D}_\tau)}  d \tau \; ,
\eeaa
which we will in particular apply for $V \in \cal T : = \{L,A,B\} \subset \cal U $. However, we have already established that for $\ga > \de $, and $0 \leq \de \leq 1\; $,
                                 \beaa
 | g^{\la\mu} \derm_{\la}   \derm_{\mu}   A_{{\cal T}}   |   &\leq& \begin{cases} c (\delta) \cdot c (\gamma) \cdot  E (  3)   \cdot \frac{\eps }{(1+t+|q|)^{3-3\delta} (1+|q|)^{1+\ga}},\quad\text{when }\quad q>0\; ,\\
       E (  3)  \cdot \frac{\eps   \cdot (1+| q |   )^{\frac{3}{2} }  }{(1+t+|q|)^{3-3\delta} }  \,\quad\text{when }\quad q<0 \; . \end{cases} \\
           \notag
       \eeaa 
       
  We, therefore, have established that
     \beaa
   \notag
&& (1+t) \cdot |\varpi(q) \cdot  g^{\la\mu} \derm_{\la}   \derm_{\mu}   A_{{\cal T}}  (t,\cdot)  \|_{L^\infty(\overline{D}_t)}  \\
   \notag
&\les&  c (\delta) \cdot  c (\gamma) \cdot   E ( 3)   \cdot  \begin{cases}    \frac{\eps \cdot \varpi(q) }{(1+t+|q|)^{2-3\delta} (1+|q|)^{1+\ga}},\quad\text{when }\quad q>0\; ,\\
\notag
 \frac{\eps  \cdot (1+| q |   )^{\frac{3}{2} } \cdot\varpi(q)  }{(1+t+|q|)^{2-3\delta} }  \,\quad\text{when }\quad q<0 \; . \end{cases}  \\
 \notag
 &\leq&   c (\delta) \cdot  c (\gamma) \cdot   E (  3) \cdot \eps \cdot  \begin{cases}  \frac{ (1+|q|)^{1+\gamma^\prime}}{(1+t+|q|)^{2-3\delta  } (1+|q|)^{1+\ga} },\quad\text{when }\quad q>0,\\
    \notag
   \frac{      (1+| q |   )^{\frac{3}{2} }  }{ (1+t+|q|)^{2-3\delta  } } \,\quad\text{when }\quad q<0 . \end{cases} \\
       \notag
        &\leq& C(q_0) \cdot   c (\delta) \cdot  c (\gamma) \cdot   E ( 3) \cdot  \frac{ \eps }{(1+t+|q|)^{2-2\delta  } (1+|q|)^{\ga -\de -\gamma^\prime  } } \; .
\eeaa

However, for f $-1 \leq \gamma^\prime < \gamma - \delta$, we have $\ga -\de -\gamma^\prime > 0$, and thus,
     \beaa
   \notag
&& (1+t) \cdot |\varpi(q) \cdot  g^{\la\mu} \derm_{\la}   \derm_{\mu}   A_{{\cal T}}  (t,\cdot)  \|_{L^\infty(\overline{D}_t)}  \\
   \notag
        &\leq& C(q_0) \cdot   c (\delta) \cdot  c (\gamma) \cdot   E (  3) \cdot  \frac{ \eps }{(1+t+|q|)^{2-2\delta  } } \; .
\eeaa

Now, taking $\de < \frac{1}{2}$, we thereby obtain that 
  \beaa
   \int_0^t  (1+\tau)\| \varpi(q)  g^{\la\mu} \derm_{\la}   \derm_{\mu}  A_{{\cal T}} (\tau,\cdot) \|_{L^\infty(\overline{D}_\tau)} & \leq& C(q_0) \cdot c (\gamma^\prime) \cdot c (\delta) \cdot  c (\gamma) \cdot E (  3)  \cdot \eps \, .
   \eeaa

   As a result, we get
        \bea
   \notag
&& (1+t+|x|) \cdot |\varpi(q)\derm \ A_{{\cal T}} (t,x)| \\
   \notag
 &\les& C(q_0) \cdot  c (\gamma^\prime) \cdot  c (\delta) \cdot c (\gamma)  \cdot E (  4) \cdot \eps \, \\
\notag
&& + c (\gamma^\prime) \cdot c (\gamma)  \cdot c (\delta)  \cdot E ( 3)  \cdot \eps  \cdot  \int_0^t  \frac{1}{(1+\tau)} \cdot (1+\tau+|x|) \cdot \|\varpi(q) \derm   A_{{\cal T}} (\tau,\cdot) \|_{L^\infty (\Sigma^{ext}_{\tau} )} d \tau .
\eea
Using Grönwall lemma, we get

          \bea
   \notag
&& (1+t+|x|) \cdot |\varpi(q)\derm \ A_{{\cal T}} (t,x)| \\
   \notag
 &\les&  C(q_0) \cdot c (\gamma^\prime) \cdot  c (\delta) \cdot c (\gamma)  \cdot E (  4) \cdot \eps  \cdot \exp \big( c (\gamma^\prime) \cdot c (\gamma)  \cdot c (\delta)  \cdot E ( 3)  \cdot \eps  \cdot  \int_0^t  \frac{1}{(1+\tau)} d\tau  \big) \\
 \notag
  &\les&  C(q_0) \cdot c (\gamma^\prime) \cdot  c (\delta) \cdot c (\gamma)  \cdot E (  4) \cdot \eps  \cdot   (1+t)^{c (\gamma^\prime)  \cdot c (\gamma)  \cdot c (\delta)  \cdot E ( 3)  \cdot \eps }  \; .
\eea

  Thus, using $(1+t+|x|) \sim  (1+t+|q|)$, we get in the exterior region $\overline{C}$, the following estimate
            \bea
   \notag
&& |\derm \ A_{{\cal T}} (t,x)| \\
   \notag
  &\les& C(q_0) \cdot  c (\gamma^\prime) \cdot  c (\delta) \cdot c (\gamma)  \cdot E (  4) \cdot \eps  \cdot \frac{1}{  (1+t+|q|)^{1-c (\gamma^\prime) \cdot c (\gamma)  \cdot c (\delta)  \cdot E ( 3)  \cdot \eps }  \cdot \varpi(q) } \\
   &\leq& C(q_0) \cdot c (\gamma^\prime)  \cdot  c (\delta) \cdot c (\gamma)  \cdot E (  4) \cdot \eps \cdot  \begin{cases}  \frac{ 1}{ (1+t+|q|)^{1-c (\gamma^\prime) \cdot c (\gamma)  \cdot c (\delta)  \cdot E ( 3)  \cdot \eps }  \cdot (1+|q|)^{1+\gamma^\prime} },\quad\text{when }\quad q>0,\\
    \notag
   \frac{     1 }{  (1+t+|q|)^{1-c (\gamma^\prime) \cdot c (\gamma)  \cdot c (\delta)  \cdot E ( 3)  \cdot \eps }  } \, ,\quad\text{when }\quad q<0  \end{cases}  \\
   \notag
   &\leq& C(q_0) \cdot c (\gamma^\prime)  \cdot  c (\delta) \cdot c (\gamma)  \cdot E (  4) \cdot \eps \cdot    \frac{ 1}{ (1+t+|q|)^{1-c (\gamma^\prime) \cdot c (\gamma)  \cdot c (\delta)  \cdot E ( 3)  \cdot \eps }  \cdot (1+|q|)^{1+\gamma^\prime} } \; .
\eea

  Integrating, we obtain in the exterior region, 
              \bea
   \notag
&& | \ A_{{\cal T}} (t,x)| \\
   \notag
   &\leq& C(q_0) \cdot c (\gamma^\prime)  \cdot  c (\delta) \cdot c (\gamma)  \cdot E (  4) \cdot \eps \cdot   \frac{ 1}{ (1+t+|q|)^{1-c (\gamma^\prime) \cdot c (\gamma)  \cdot c (\delta)  \cdot E ( 3)  \cdot \eps }  \cdot (1+|q|)^{\gamma^\prime} } \; , \\ 
\eea
  where $-1 \leq \gamma^\prime <  \ga -\delta \; $. 
  \end{proof}

\subsection{The terms   $| A_L  | \cdot    | \derm A   |$  and $ | A_{e_{a}}  | \cdot    | \derm A_{e_{a}} |\; $}\

\begin{lemma}
We have, in the exterior $\overline{C}$,
        \beaa
| A_L  | \cdot    | \derm A   | &\les&  C(q_0) \cdot  c (\delta) \cdot c (\gamma) \cdot E (  3)   \cdot \frac{\eps }{(1+t+|q|)^{3-3\delta} (1+|q|)^{2\ga}} \; .\\
      \eeaa

\end{lemma}

\begin{proof}

We have shown in Lemma \ref{estimategoodcomponentspotentialandmetric}, that
              \bea
        \notag
 |   \Lie_{Z^I}  A_{L}  |   &\leq& \int\limits_{s,\,\Om=const} \sum_{|J|\leq |I| -1} |   \derm  \Lie_{Z^J} A | \\
      \notag
     && + \begin{cases} c (\delta) \cdot c (\gamma) \cdot C ( |I| ) \cdot E ( |I| + 3)  \cdot \big( \frac{ \eps   }{ (1+t+|q|)^{2-2\delta} \cdot  (1+|q|)^{\ga - 1} } \big),\quad\text{when }\quad q>0 \; ,\\
       C ( |I| ) \cdot E ( |I| + 3)  \cdot \big( \frac{ \eps \cdot (1+|q|)^{\frac{3}{2} }   }{ (1+t+|q|)^{2-2\delta} } \big) \,\quad\text{when }\quad q<0 \; . \end{cases} 
       \eea
         Thus, for $|I| = 0$, we get
                    \bea
        \notag
 |     A_{L}  |   &\les&  \begin{cases} c (\delta) \cdot c (\gamma) \cdot E (  3)  \cdot \frac{ \eps   }{ (1+t+|q|)^{2-2\delta} \cdot  (1+|q|)^{\ga - 1}  }  ,\quad\text{when }\quad q>0 \; ,\\
         E (  3) \cdot  \frac{ \eps \cdot (1+|q|)^{\frac{3}{2}  } }{ (1+t+|q|)^{2-2\delta} }  \,\quad\text{when }\quad q<0 \; . \end{cases} \\
       \eea
         On the other hand, from the bootstrap assumption, we have
                        \beaa
 \notag
|\derm A      |    &\leq& \begin{cases}  E ( 2)  \cdot \frac{\eps }{(1+t+|q|)^{1-\delta} (1+|q|)^{1+\ga}},\quad\text{when }\quad q>0 \; ,\\
        E ( 2)  \cdot \frac{\eps  }{(1+t+|q|)^{1-\delta}(1+|q|)^{\frac{1}{2} }}  \,\quad\text{when }\quad q<0 \; . \end{cases} \\
      \eeaa
     Hence,
     \beaa
| A_L  | \cdot    | \derm A   | \les \begin{cases}  c (\delta) \cdot c (\gamma) \cdot E (  3)  \cdot \frac{\eps }{(1+t+|q|)^{3-3\delta} (1+|q|)^{2\ga}},\quad\text{when }\quad q>0\; ,\\
        E ( 3)  \cdot \frac{\eps  \cdot (1+|q|) }{(1+t+|q|)^{3-3\delta}}  \,\quad\text{when }\quad q<0 \;  . \end{cases} \\
      \eeaa
     Consequently, in the exterior $\overline{C}$, we have
         \beaa
| A_L  | \cdot    | \derm A   | &\les&  C(q_0) \cdot  c (\delta) \cdot c (\gamma) \cdot E (  3)   \cdot \frac{\eps }{(1+t+|q|)^{3-3\delta} (1+|q|)^{2\ga}} \; .\\
      \eeaa
      
      \end{proof}
      
      \begin{lemma}
Let $M\leq \eps$\,, and let $\ga > \de $\,, $0 \leq \de < \frac{1}{2}$ and $-1 \leq \gamma^\prime <  \ga - \delta $\,. Then, we have
     \beaa
      |  A_{e_a}   |       &\leq& C(q_0) \cdot c (\gamma^\prime)  \cdot  c (\delta) \cdot c (\gamma)  \cdot E (  4) \cdot \eps \cdot    \frac{ 1}{ (1+t+|q|)^{1-c (\gamma^\prime) \cdot c (\gamma)  \cdot c (\delta)  \cdot E ( 3)  \cdot \eps }  \cdot (1+|q|)^{\gamma^\prime} } \; ,
\eeaa 
and as a result,
     \beaa
  &&    | A_{e_a}   |  \cdot | \derm   A_{e_a}  |    \\
       &\leq& C(q_0) \cdot c (\gamma^\prime)  \cdot  c (\delta) \cdot c (\gamma)  \cdot E (  4) \cdot \eps \cdot    \frac{ 1}{ (1+t+|q|)^{2-c (\gamma^\prime) \cdot c (\gamma)  \cdot c (\delta)  \cdot E ( 3)  \cdot \eps }  \cdot (1+|q|)^{1+2\gamma^\prime} } \; .
\eeaa 
\end{lemma}

\begin{proof}

      In order to estimate the term $ | A_{e_{a}}  |$, we would like to estimate $ | \pa A_{e_{a}}  |$ and then integrate along $s= constant$ and $\Om = constant$ as before. However, for this, we will use the special fact that it is an $ A_{e_{a}}$ component and that we have also estimated $\derm  A_{e_{a}}  $. In fact,
      \beaa
    \pa_r  \Lie_{Z^I} A_{e_{a}} = \derm_r  ( \Lie_{Z^I} A_{e_{a}} )  +  \Lie_{Z^I} A (\derm_{r} e_a) \; .
     \eeaa

      However, as we already computed in the Appendix of \cite{G2}, we have $\derm_{r} e_a = 0 $ (the computation was carried out for a Schwarzschild metric, yet Minkowski is a special case where the Schwarzschild mass in null). Consequently, we get that
      \bea \label{rderiofAacomp}
    \pa_r  \Lie_{Z^I} A_{e_{a}} = \derm_r (  \Lie_{Z^I} A_{e_{a}} ) 
     \eea 
   and as a result, integrating as we have already shown in \cite{G4}, 
       \beaa
\notag
| \Lie_{Z^I} A_{e_{a}}  (t, | x | \cdot \Om) | &=& | \Lie_{Z^I} A_{e_{a}} \big(0, ( t + | x |) \cdot \Om \big) |  + \int_{| x | }^{t + | x |  } \pa_r | (\Lie_{Z^I} A_{e_{a}}  (t + | x | - r,  r  \cdot \Om ) ) | dr \\
 &\leq& | \Lie_{Z^I} A_{e_{a}}  \big(0, ( t + | x | ) \cdot \Om \big) |  + \int_{| x | }^{t + | x |  } | \pa_r  (\Lie_{Z^I} A_{e_{a}}  (t + | x | - r,  r  \cdot \Om ) ) | dr \\
  &\leq& | \Lie_{Z^I} A_{e_{a}}  \big(0, ( t + | x | ) \cdot \Om \big)|  + \int_{| x | }^{t + | x |  } | \derm_r  (\Lie_{Z^I} A_{e_a} (t + | x | - r,  r  \cdot \Om ) ) | dr \; . \\
\eeaa
     However, we have shown that in the exterior, for $\ga > \de $, $0 \leq \de < \frac{1}{2}$ and $-1 \leq \gamma^\prime <  \ga - \delta $, we have

     \beaa
      | \derm_r   A_{e_a}   | &\les&   | \derm  A_{e_a}  |  \leq  | \derm  A_{\cal T}  |  \\
      &\leq& C(q_0) \cdot c (\gamma^\prime)  \cdot  c (\delta) \cdot c (\gamma)  \cdot E (  4) \cdot \eps \cdot    \frac{ 1}{ (1+t+|q|)^{1-c (\gamma^\prime) \cdot c (\gamma)  \cdot c (\delta)  \cdot E ( 3)  \cdot \eps }  \cdot (1+|q|)^{1+\gamma^\prime} } \; .
\eeaa     
Hence, proceeding the integration, as in \cite{G4}, we get that in the exterior $\overline{C}$,
     \beaa
      |  A_{e_a}   |       &\leq& C(q_0) \cdot c (\gamma^\prime)  \cdot  c (\delta) \cdot c (\gamma)  \cdot E (  4) \cdot \eps \cdot    \frac{ 1}{ (1+t+|q|)^{1-c (\gamma^\prime) \cdot c (\gamma)  \cdot c (\delta)  \cdot E ( 3)  \cdot \eps }  \cdot (1+|q|)^{\gamma^\prime} } \; .
\eeaa 
As a result, we have
     \beaa
  &&    | A_{e_a}   |  \cdot | \derm   A_{e_a}  |    \\
       &\leq& C(q_0) \cdot c (\gamma^\prime)  \cdot  c (\delta) \cdot c (\gamma)  \cdot E (  4) \cdot \eps \cdot    \frac{ 1}{ (1+t+|q|)^{2-c (\gamma^\prime) \cdot c (\gamma)  \cdot c (\delta)  \cdot E ( 3)  \cdot \eps }  \cdot (1+|q|)^{1+2\gamma^\prime} } \; .
\eeaa 

\end{proof}

\subsection{Estimate on all the components of the Yang-Mills potential}\

\begin{lemma}\label{upgradedestimateonthegradientofthefullcompoenentsofpoentialA}

In the Lorenz gauge and in wave coordinates, we have for $M\leq \eps$\,, for \,$\ga > \de $\,, and $0 \leq \de < \frac{1}{2}$\,, and for  
  \beaa
\de &\leq& \gamma^\prime <  \ga - \delta \; ,
  \eeaa
the following estimate
   \bea
 \notag
 \,|\derm A_{{\cal U}} (t,x)|  &\les&    C(q_0) \cdot c (\gamma^\prime)  \cdot  c (\delta) \cdot c (\gamma)  \cdot E (  4)  \frac{\eps }{ (1+t+|q|)^{1-c (\gamma^\prime)\cdot c (\gamma)  \cdot c (\delta)  \cdot E ( 4)  \cdot \eps} \cdot (1 +|q|)^{1+\gamma^\prime}} \, . \\
\eea

\end{lemma}

\begin{proof}

We already showed, in Lemma \ref{structureofthesourcetermsofthewaveoperatoronAandh}, that
\beaa
   \notag
 && | g^{\la\mu} \derm_{\la}   \derm_{\mu}   A_{{\underline{L}}}  |   \\
   &\les& | \derm h | \cdot  |\rderm A |       + | \rderm  h | \cdot  |\derm A |  \\
   \notag
   && +           | A  | \cdot    | \rderm A  | +  | \derm  h | \cdot  | A  |^2   +  | A  |^3 \\
       \notag
           && + | A_L  | \cdot    | \derm A   |    + | A_{e_{a}}  | \cdot    | \derm A_{e_{a}} |      \\
    \notag
 && + O( h \cdot  \derm h \cdot  \derm A) + O( h \cdot  A \cdot \derm A)  + O( h \cdot  \derm h \cdot  A^2) + O( h \cdot  A^3) \, .
   \eeaa

    We have already shown, in \eqref{decayestimateonthetermsinthesourcetermsforwaveeqonAtau}, that in the exterior region $\overline{C}$, for $\ga > \de $, and $0 \leq \de \leq 1$,

                                            \bea
                                            \notag
&&| \derm h | \cdot  |\rderm A |    + | \rderm  h | \cdot  |\derm A |   +           | A  | \cdot    | \rderm A  | +  | \derm  h | \cdot  | A  |^2   +  | A  |^3  \\
\notag
 && + O( h \cdot  \derm h \cdot  \derm A) + O( h \cdot  A \cdot \derm A)  + O( h \cdot  \derm h \cdot  A^2) + O( h \cdot  A^3) \\
 \notag
&\les& \begin{cases} c (\delta) \cdot c (\gamma) \cdot E (  3)   \cdot \frac{\eps }{(1+t+|q|)^{3-3\delta} (1+|q|)^{1+\ga}},\quad\text{when }\quad q>0\; ,\\
\notag
        E (  3)  \cdot \frac{\eps  }{(1+t+|q|)^{\frac{3}{2}-3\delta} }  \,\quad\text{when }\quad q<0 \; , \end{cases} \\
 &\les&       C(q_0) \cdot    c (\delta) \cdot c (\gamma) \cdot E (  3)   \cdot \frac{\eps }{(1+t+|q|)^{3-3\delta} (1+|q|)^{1+\ga}} \; .
       \eea
        
    Thus, the only new terms to estimate are  $| A_L  | \cdot    | \derm A   |$  and $ | A_{e_{a}}  | \cdot    | \derm A_{e_{a}} |\; $.
    
    However, we have shown that in the exterior $\overline{C}$, for $\ga > \de $, $0 \leq \de < \frac{1}{2}$ and $-1 \leq \gamma^\prime <  \ga - \delta $,
                            
        \beaa
&& | A_L  | \cdot    | \derm A   | + | A_{e_a}   |  \cdot | \derm   A_{e_a}  |  \\
 &\les&  C(q_0) \cdot  c (\delta) \cdot c (\gamma) \cdot E (  3)   \cdot \frac{\eps }{(1+t+|q|)^{3-3\delta} (1+|q|)^{2\ga}} \\
 \notag
   && +  C(q_0) \cdot c (\gamma^\prime)  \cdot  c (\delta) \cdot c (\gamma)  \cdot E (  4) \cdot \eps \cdot    \frac{ 1}{ (1+t+|q|)^{2-c (\gamma^\prime) \cdot c (\gamma)  \cdot c (\delta)  \cdot E ( 3)  \cdot \eps }  \cdot (1+|q|)^{1+2\gamma^\prime} } \\
    &\les&  C(q_0) \cdot  c (\delta) \cdot c (\gamma) \cdot E (  3)   \cdot \frac{\eps }{(1+t+|q|)^{2}  (1+t+|q|)^{1-3\delta}  (1+|q|)^{2\ga} } \\
 \notag
   && +  C(q_0) \cdot c (\gamma^\prime)  \cdot  c (\delta) \cdot c (\gamma)  \cdot E (  4) \cdot \eps \cdot    \frac{ 1}{ (1+t+|q|)^{2-c (\gamma^\prime) \cdot c (\gamma)  \cdot c (\delta)  \cdot E ( 3)  \cdot \eps }  \cdot (1+|q|)^{1+2\gamma^\prime} } \\
      &\les&  C(q_0) \cdot  c (\delta) \cdot c (\gamma) \cdot E (  3)   \cdot \frac{\eps }{(1+t+|q|)^{2}   (1+|q|)^{1+ 2 (\ga - \de) - \de} } \\
 \notag
   && +  C(q_0) \cdot c (\gamma^\prime)  \cdot  c (\delta) \cdot c (\gamma)  \cdot E (  4) \cdot \eps \cdot    \frac{ 1}{ (1+t+|q|)^{2-c (\gamma^\prime) \cdot c (\gamma)  \cdot c (\delta)  \cdot E ( 3)  \cdot \eps }  \cdot (1+|q|)^{1+2\gamma^\prime} } \; .\\
   \eeaa
   
   Now, given that $\gamma^\prime <  \ga - \delta$,
   
           \beaa
&& | A_L  | \cdot    | \derm A   | + | A_{e_a}   |  \cdot | \derm   A_{e_a}  |  \\
      &\les&  C(q_0) \cdot  c (\delta) \cdot c (\gamma) \cdot E (  3)   \cdot \frac{\eps }{(1+t+|q|)^{2}   (1+|q|)^{1+ 2 \ga^\prime - \de} } \\
 \notag
   && +  C(q_0) \cdot c (\gamma^\prime)  \cdot  c (\delta) \cdot c (\gamma)  \cdot E (  4) \cdot \eps \cdot    \frac{ 1}{ (1+t+|q|)^{2-c (\gamma^\prime) \cdot c (\gamma)  \cdot c (\delta)  \cdot E ( 3)  \cdot \eps }  \cdot (1+|q|)^{1+2\gamma^\prime} } \\
     &\les&  C(q_0) \cdot c (\gamma^\prime)  \cdot  c (\delta) \cdot c (\gamma)  \cdot E (  4)  \cdot    \frac{ \eps}{ (1+t+|q|)^{2-c (\gamma^\prime) \cdot c (\gamma)  \cdot c (\delta)  \cdot E ( 3)  \cdot \eps }  \cdot (1+|q|)^{1+2\gamma^\prime -\de} } \\
 && \text{(given that $E(3) \leq E(4)$)} .\\
   \eeaa
   
Therefore,
\beaa
   \notag
 && | g^{\la\mu} \derm_{\la}   \derm_{\mu}   A_{{\underline{L}}}  |   \\
&\les&C(q_0) \cdot    c (\delta) \cdot c (\gamma) \cdot E (  3)   \cdot \frac{\eps }{(1+t+|q|)^{3-3\delta} (1+|q|)^{1+\ga}} \\
&& + C(q_0) \cdot c (\gamma^\prime)  \cdot  c (\delta) \cdot c (\gamma)  \cdot E (  4) \cdot \eps \cdot    \frac{ 1}{ (1+t+|q|)^{2-c (\gamma^\prime) \cdot c (\gamma)  \cdot c (\delta)  \cdot E ( 3)  \cdot \eps }  \cdot (1+|q|)^{1+2\gamma^\prime - \de} } \\
&\les&   c (\delta) \cdot c (\gamma) \cdot E (  3)   \cdot \frac{\eps }{(1+t+|q|)^{2} (1+|q|)^{2+\ga - 3\de}} \\
&& + C(q_0) \cdot c (\gamma^\prime)  \cdot  c (\delta) \cdot c (\gamma)  \cdot E (  4) \cdot \eps \cdot    \frac{ 1}{ (1+t+|q|)^{2-c (\gamma^\prime) \cdot c (\gamma)  \cdot c (\delta)  \cdot E ( 3)  \cdot \eps }  \cdot (1+|q|)^{1+2\gamma^\prime - \de} } \\
&\les&   c (\delta) \cdot c (\gamma) \cdot E (  3)   \cdot \frac{\eps }{(1+t+|q|)^{2} (1+|q|)^{1+ (\ga - \de ) + 1- 2\de}} \\
&& + C(q_0) \cdot c (\gamma^\prime)  \cdot  c (\delta) \cdot c (\gamma)  \cdot E (  4) \cdot \eps \cdot    \frac{ 1}{ (1+t+|q|)^{2-c (\gamma^\prime) \cdot c (\gamma)  \cdot c (\delta)  \cdot E ( 3)  \cdot \eps }  \cdot (1+|q|)^{1+2\gamma^\prime - \de} } \; . \\
\eeaa

Now, given that $\gamma^\prime < \gamma - \delta$, and that $E(3) \leq E(4)$, we get
\beaa
   \notag
 && | g^{\la\mu} \derm_{\la}   \derm_{\mu}   A_{{\underline{L}}}  |   \\
&\les&    C(q_0) \cdot c (\gamma^\prime)  \cdot  c (\delta) \cdot c (\gamma)  \cdot E (  4) \cdot \frac{\eps }{(1+t+|q|)^{2} (1+|q|)^{1+ \ga^\prime + 1- 2\de}} \\
&& + C(q_0) \cdot c (\gamma^\prime)  \cdot  c (\delta) \cdot c (\gamma)  \cdot E (  4) \cdot    \frac{ \eps}{ (1+t+|q|)^{2-c (\gamma^\prime) \cdot c (\gamma)  \cdot c (\delta)  \cdot E ( 3)  \cdot \eps }  \cdot (1+|q|)^{1+2\gamma^\prime - \de} } \; . \\
\eeaa

Applying the Lindblad-Rodnianski estimate again (see \eqref{LinfinitynormestimateongradientderivedbyLondbladRodnianski}) for $V \in \cal U$, since we also assumed that $\gamma^\prime$ is such that $-1 \leq \gamma^\prime < \gamma - \delta$, and that $ \delta <    1/2 $, we finally obtain that for all $V\in  \{L,\Lb,A,B\}$, given what we had already proved, that
     \beaa
   \notag
&& (1+t+|x|) \cdot |\varpi(q) \cdot \derm  A_{\cal U} (t,x)| \\
   \notag
 &\les&   c (\gamma^\prime) \cdot  c (\delta) \cdot c (\gamma) \cdot E (  4) \cdot \eps \, \\
\notag
&& + c (\gamma^\prime) \cdot c (\gamma)  \cdot c (\delta)  \cdot E ( 3)  \cdot \eps  \cdot  \int_0^t  \frac{1}{(1+\tau)} \cdot (1+\tau+|x|) \cdot \|\varpi(q)  \cdot  \derm   A_{\cal U} (\tau,\cdot) \|_{L^\infty (\Sigma^{ext}_{\tau} )} d \tau \\
       \notag
&& + \int_0^t  (1+\tau) \cdot \| \varpi(q)  \cdot g^{\la\mu} \derm_{\la}   \derm_{\mu}  A_{\cal U} (\tau,\cdot) \|_{L^\infty(\overline{D}_\tau)} d \tau \\
\notag
&\les&  c (\gamma^\prime) \cdot  c (\delta) \cdot c (\gamma) \cdot E (  4) \cdot \eps   \\
&& + c (\gamma^\prime) \cdot c (\gamma)  \cdot c (\delta)  \cdot E ( 3)  \cdot \eps  \cdot  \int_0^t  \frac{1}{(1+\tau)} \cdot (1+\tau+|x|) \cdot \|\varpi(q)  \cdot  \derm   A_{\cal U} (\tau,\cdot) \|_{L^\infty (\Sigma^{ext}_{\tau} )} d \tau \\
&& +   \int_0^t    C(q_0) \cdot c (\gamma^\prime)  \cdot  c (\delta) \cdot c (\gamma)  \cdot E (  4) \cdot \eps \cdot    \frac{ (1+\tau) \cdot (1+|q|)^{1+\gamma^{\prime}}}{ (1+\tau+|q|)^{2 }  \cdot (1+|q|)^{1+\gamma^\prime + 1 -2\de} }   d\tau \\
&& +   \int_0^t    C(q_0) \cdot c (\gamma^\prime)  \cdot  c (\delta) \cdot c (\gamma)  \cdot E (  4) \cdot \eps \cdot    \frac{ (1+\tau) \cdot (1+|q|)^{1+\gamma^{\prime}}}{ (1+\tau+|q|)^{2-c (\gamma^\prime) \cdot c (\gamma)  \cdot c (\delta)  \cdot E ( 3)  \cdot \eps }  \cdot (1+|q|)^{1+2\gamma^\prime -\de} }   d\tau \\
&\les&  c (\gamma^\prime) \cdot  c (\delta) \cdot c (\gamma) \cdot E (  4) \cdot \eps   \\
&& + c (\gamma^\prime) \cdot c (\gamma)  \cdot c (\delta)  \cdot E ( 3)  \cdot \eps  \cdot  \int_0^t  \frac{1}{(1+\tau)} \cdot (1+\tau+|x|) \cdot \|\varpi(q)  \cdot  \derm   A_{\cal U} (\tau,\cdot) \|_{L^\infty (\Sigma^{ext}_{\tau} )} d \tau \\
&& +   \int_0^t    C(q_0) \cdot c (\gamma^\prime)  \cdot  c (\delta) \cdot c (\gamma)  \cdot E (  4) \cdot     \frac{ \eps}{ (1+\tau+|q|)^{1-c (\gamma^\prime) \cdot c (\gamma)  \cdot c (\delta)  \cdot E ( 4)  \cdot \eps }  \cdot (1+|q|)^{1-2\de} }   d\tau \\
&& +   \int_0^t    C(q_0) \cdot c (\gamma^\prime)  \cdot  c (\delta) \cdot c (\gamma)  \cdot E (  4) \cdot     \frac{ \eps}{ (1+\tau+|q|)^{1-c (\gamma^\prime) \cdot c (\gamma)  \cdot c (\delta)  \cdot E ( 4)  \cdot \eps }  \cdot (1+|q|)^{\gamma^\prime -\de} }   d\tau \\
&& \text{(given that $E(3) \leq E(4)$) }. \\
\eeaa
Now, we choose $\ga^\prime \geq \de$ and as a result, we obtain 
     \beaa
   \notag
&& (1+t+|x|) \cdot |\varpi(q) \cdot \derm  A_{\cal U} (t,x)| \\
   \notag
&\les&  c (\gamma^\prime) \cdot  c (\delta) \cdot c (\gamma) \cdot E (  4) \cdot \eps  +  C(q_0) \cdot c (\gamma^\prime)  \cdot  c (\delta) \cdot c (\gamma)  \cdot E (  4) \cdot \eps \cdot \big[   (1+t)^{c (\gamma^\prime) \cdot c (\gamma)  \cdot c (\delta)  \cdot E ( 4)  \cdot \eps }   \big]^{t}_{0} \\
 && + c (\gamma^\prime) \cdot c (\gamma)  \cdot c (\delta)  \cdot E ( 3)  \cdot \eps  \cdot  \int_0^t  \frac{1}{(1+\tau)} \cdot (1+\tau+|x|) \cdot \|\varpi(q)  \cdot  \derm   A_{\cal U} (\tau,\cdot) \|_{L^\infty (\Sigma^{ext}_{\tau} )} d \tau \\
 &\les&  C(q_0) \cdot c (\gamma^\prime)  \cdot  c (\delta) \cdot c (\gamma)  \cdot E (  4) \cdot \eps \cdot    (1+t)^{c (\gamma^\prime) \cdot c (\gamma)  \cdot c (\delta)  \cdot E ( 4)  \cdot \eps }    \\
 && + c (\gamma^\prime) \cdot c (\gamma)  \cdot c (\delta)  \cdot E ( 3)  \cdot \eps  \cdot  \int_0^t  \frac{1}{(1+\tau)} \cdot (1+\tau+|x|) \cdot \|\varpi(q)  \cdot  \derm   A_{\cal U} (\tau,\cdot) \|_{L^\infty (\Sigma^{ext}_{\tau} )} d \tau \; .
\eeaa
Using Grönwall inequality, we get
    \beaa
   \notag
&& (1+t+|x|) \cdot |\varpi(q) \cdot \derm  A_{\cal U} (t,x)| \\
   \notag
   &\les &  C(q_0) \cdot c (\gamma^\prime)  \cdot  c (\delta) \cdot c (\gamma)  \cdot E (  4) \cdot \eps \cdot    (1+t)^{c (\gamma^\prime) \cdot c (\gamma)  \cdot c (\delta)  \cdot E ( 4)  \cdot \eps }  \cdot (1+t)^{c (\gamma^\prime) \cdot c (\gamma)  \cdot c (\delta)  \cdot E ( 3)  \cdot \eps} \; .
   \eeaa
   
Hence, for $\delta < \frac{1}{2}$, we obtain in the exterior,
\beaa
 (1+t+|x|)| \cdot |\varpi(q) \cdot |\derm A_{{\cal U}} (t,x)|  &\les&    C(q_0) \cdot c (\gamma^\prime)  \cdot  c (\delta) \cdot c (\gamma)  \cdot E (  4) \cdot \eps  \cdot (1+ t )^{c (\gamma^\prime)  \cdot c (\gamma)  \cdot c (\delta)  \cdot E ( 4)  \cdot \eps } \, ,
\eeaa
 which leads to 
 \bea
 \notag
 \,|\derm A_{{\cal U}} (t,x)|  &\les&    C(q_0) \cdot c (\gamma^\prime)  \cdot  c (\delta) \cdot c (\gamma)  \cdot E (  4)  \frac{\eps }{ (1+t+|x|)^{1-c (\gamma^\prime)\cdot c (\gamma)  \cdot c (\delta)  \cdot E ( 4)  \cdot \eps} \cdot (1 +|q|)^{1+\gamma^\prime}} \, . \\
\eea

  \end{proof}

\subsection{Improved estimate for the bad component of the Yang-Mills metric}\

We will start by improving first the decay estimate on the good components of the Einstein-Yang-Mils metric $h_{\cal T \cal U}$\,, in \eqref{upgradedestimatesongoodcomponnentforh1andh0}, so as to have more decay in $|q|$. 
\begin{lemma}
We have for $M\leq \eps$\,, for $\ga > \de $, $\ga > 1 $, and $0 \leq \delta \leq \frac{1}{4}$, and for
\beaa
 - 1 \leq \gamma^\prime \leq 0 \; ,
\eeaa
that
 \beaa
 \,|\derm h^1_{ {\cal T} {\cal U}} (t,x)|     &\leq& C(q_0) \cdot c (\gamma^\prime)  \cdot  c (\delta) \cdot c (\gamma)  \cdot E (  4) \cdot    \frac{ \eps  }{ (1+t+|q|)^{1-c (\gamma^\prime) \cdot c (\gamma)  \cdot c (\delta)  \cdot E ( 3)  \cdot \eps }  \cdot (1+|q|)^{1+\gamma^\prime} } \; . \\ 
\eeaa
\end{lemma}

\begin{proof}

We have shown in Lemma \ref{GronwallinequalitiesoncomponentsofgradientofAandgradientofh1withweight}, that for $\gamma^\prime$ such that $-1 \leq \gamma^\prime < \gamma - \delta$, and $ \delta <    1/2 $, we have
     \bea
   \notag
&& (1+t+|q|) \cdot |\varpi(q) \cdot \derm ( \Lie_{Z^J} h^1)_{\cal T \cal U} (t,x)| \\
   \notag
 &\les&   c (\gamma^\prime) \cdot  c (\delta) \cdot c (\gamma) \cdot C ( |J|  ) \cdot E ( |J| + 4) \cdot \eps \, \\
\notag
&& + c (\gamma^\prime) \cdot c (\gamma)  \cdot c (\delta)  \cdot E ( 3)  \cdot \eps  \cdot  \int_0^t  \frac{1}{(1+\tau)} \cdot (1+\tau+|x|) \cdot \|\varpi(q)  \cdot  \derm  ( \Lie_{Z^J} h^1)_{\cal T \cal U} (\tau,\cdot) \|_{L^\infty (\Sigma^{ext}_{\tau} )} d \tau \\
       \notag
&& + \int_0^t  (1+\tau) \cdot \| \varpi(q)  \cdot g^{\la\mu} \derm_{\la}   \derm_{\mu}  (  \Lie_{Z^J} h^1)_{\cal T \cal U} (\tau,\cdot) \|_{L^\infty(\overline{D}_\tau)}  d \tau \; .\\
\eea
On the other hand, we have shown in Lemma \ref{decayestimateonsourcetermsforgoodcomponentofmetrichtauU}, that for $\ga > \de $\,, $\ga > 1 $\,, and $0 \leq \delta \leq \frac{1}{4} < \frac{1}{2}$\,, 
      \bea
\notag
(1+t ) \cdot | g^{\alpha\beta} \derm_\alpha \derm_\beta h_{ {\cal T} {\cal U}} |_{L^\infty(\overline{D}_\tau)}  &\les&
       C(q_0) \cdot  c (\delta) \cdot c (\gamma) \cdot E (  3)  \cdot \frac{\eps }{(1+t )^{2-3\delta} (1+|q|)^{1+2\de} }  \; . \\
    \eea
    Thus, in the exterior,
          \beaa
\notag
(1+t ) \cdot \varpi(q) \cdot | g^{\alpha\beta} \derm_\alpha \derm_\beta h_{ {\cal T} {\cal U}} |_{L^\infty(\overline{D}_\tau)}  &\les&   C(q_0) \cdot  c (\delta) \cdot c (\gamma) \cdot E (  3)  \cdot \frac{\eps \cdot (1+|q|)^{1+\gamma^\prime}}{(1+t )^{2-3\delta} (1+|q|)^{1+2\de} }  \\
&\les&   C(q_0) \cdot  c (\delta) \cdot c (\gamma) \cdot E (  3)  \cdot \frac{\eps }{(1+t )^{2-3\delta} (1+|q|)^{2\de-\gamma^\prime} }\; . \\
    \eeaa
    Hence, if  $ - 1 \leq \gamma^\prime \leq 2\de$, we get
          \bea
\notag
(1+t )\cdot \varpi(q) \cdot | g^{\alpha\beta} \derm_\alpha \derm_\beta h_{ {\cal T} {\cal U}} |_{L^\infty(\overline{D}_\tau)}  &\les&
       C(q_0) \cdot  c (\delta) \cdot c (\gamma) \cdot E (  3)  \cdot \frac{\eps }{(1+t )^{2-3\delta}  } \; . \\
    \eea

Therefore, using the estimate on $h^0$ that we obtained in Lemma \ref{estimateonthesourcetermsforhzerothesphericallsymmtrpart}, we get
          \beaa
\notag
&& (1+t )\cdot \varpi(q) \cdot | g^{\alpha\beta} \derm_\alpha \derm_\beta h^1_{ {\cal T} {\cal U}} |_{L^\infty(\overline{D}_\tau)} \\
 &\les& (1+t )\cdot \varpi(q) \cdot | g^{\alpha\beta} \derm_\alpha \derm_\beta h^0_{ {\cal T} {\cal U}} |_{L^\infty(\overline{D}_\tau)} \\
&& +  C(q_0) \cdot  c (\delta) \cdot c (\gamma) \cdot E (  3)  \cdot \frac{\eps }{(1+t )^{2-3\delta}  }  \\
 &\les&   c (\gamma) \cdot  E ( 2)  \cdot \frac{\eps  \cdot (1+|q|)^{1+\gamma^\prime} }{(1+t+|q|)^{2} } \\
&& +  C(q_0) \cdot  c (\delta) \cdot c (\gamma) \cdot E (  3)  \cdot \frac{\eps }{(1+t )^{2-3\delta}  } \;. \\
    \eeaa
    If  $ - 1 \leq \gamma^\prime \leq 0$, then 
            \bea
\notag
&& (1+t )\cdot \varpi(q) \cdot | g^{\alpha\beta} \derm_\alpha \derm_\beta h^1_{ {\cal T} {\cal U}} |_{L^\infty(\overline{D}_\tau)} \\
\notag
 &\les&     c (\gamma) \cdot  E ( 2)  \cdot \frac{\eps }{(1+t+|q|) } + C(q_0) \cdot  c (\delta)  \cdot c (\gamma) \cdot E (  3)  \cdot \frac{\eps \cdot M }{(1+t+|q|) }  \;. \\
    \eea

    Proceeding exactly as earlier, we obtain
         \beaa
   \notag
&& (1+t+|q|) \cdot |\varpi(q) \cdot \derm h^1_{\cal T \cal U} (t,x)| \\
   \notag
 &\les&   c (\gamma^\prime) \cdot  c (\delta) \cdot c (\gamma)  \cdot E ( 4) \cdot \eps \, \\
\notag
&& + c (\gamma^\prime) \cdot c (\gamma)  \cdot c (\delta)  \cdot E ( 3)  \cdot \eps  \cdot  \int_0^t  \frac{1}{(1+\tau)} \cdot (1+\tau+|x|) \cdot \|\varpi(q)  \cdot  \derm   h^1_{\cal T \cal U} (\tau,\cdot) \|_{L^\infty (\Sigma^{ext}_{\tau} )} d \tau \\
       \notag
&& +   \int_0^t   c (\gamma) \cdot  E ( 2)  \cdot \frac{\eps }{(1+\tau) }  d \tau + \int_0^t         C(q_0) \cdot  c (\delta) \cdot c (\gamma) \cdot E (  3)    \cdot \frac{\eps }{(1+t )^{2-3\delta}  }  d \tau \\
&& \text{(we choose here again $\delta < \frac{1}{4} $)}  \\
&\les&    C(q_0) \cdot  c (\gamma^\prime) \cdot  c (\delta) \cdot c (\gamma) \cdot E (  4) \cdot \eps \cdot  \big[  \frac{1}{(1+\tau)^{1 -3\delta} } \big]^{\infty}_{0} +  c (\gamma) \cdot  E ( 2)  \cdot \eps  \cdot \ln(1+ t)  \\
 && + c (\gamma^\prime) \cdot c (\gamma)  \cdot c (\delta)  \cdot E ( 3)  \cdot \eps  \cdot  \int_0^t  \frac{1}{(1+\tau)} \cdot (1+\tau+|x|) \cdot \|\varpi(q)  \cdot  \derm   h^1_{\cal T \cal U} (\tau,\cdot) \|_{L^\infty (\Sigma^{ext}_{\tau} )} d \tau \\
&& \text{(for $\delta \leq \frac{1}{4} < \frac{1}{3}$)} \\
&\les&    C(q_0) \cdot  c (\gamma^\prime) \cdot  c (\delta) \cdot c (\gamma) \cdot E (  4) \cdot \eps \cdot  \ln(1+ t)    \\
 && + c (\gamma^\prime) \cdot c (\gamma)  \cdot c (\delta)  \cdot E ( 3)  \cdot \eps  \cdot  \int_0^t  \frac{1}{(1+\tau)} \cdot (1+\tau+|x|) \cdot \|\varpi(q)  \cdot  \derm   h^1_{\cal T \cal U} (\tau,\cdot) \|_{L^\infty (\Sigma^{ext}_{\tau} )} d \tau  \;. . \\
\eeaa
Hence, for $\delta \leq \frac{1}{4}$, we obtain using Grönwall lemma,
\beaa
 && (1+t+|q|) \cdot  \varpi(q) \,|\derm h^1_{ {\cal T} {\cal U}} (t,x)|  \\
 &\les& C(q_0) \cdot c (\gamma^\prime) \cdot  c (\delta) \cdot c (\gamma)  \cdot E (  4) \cdot \eps  \cdot \ln(1+ t) \cdot   (1+t)^{c (\gamma^\prime)  \cdot c (\gamma)  \cdot c (\delta)  \cdot E ( 3)  \cdot \eps }   \; .
\eeaa
Using the fact that $\ln(1+ t)$ grows slower than any power, we get
 \beaa
 \,|\derm h^1_{ {\cal T} {\cal U}} (t,x)|     &\leq& C(q_0) \cdot c (\gamma^\prime)  \cdot  c (\delta) \cdot c (\gamma)  \cdot E (  4)  \cdot  \frac{ \eps  }{ (1+t+|q|)^{1-c (\gamma^\prime) \cdot c (\gamma)  \cdot c (\delta)  \cdot E ( 3)  \cdot \eps }  \cdot (1+|q|)^{1+\gamma^\prime} } \; . \\ 
\eeaa
\end{proof}

Now, we could use the upgraded estimate on the “good” components of the metric to upgrade the estimate on the “bad” components of the metric.

\begin{lemma}
For $M \leq \eps$\,, for $\ga > \de $\,, $0 \leq \de < \frac{1}{4}$ and for $  -\frac{1}{2} \leq \gamma^\prime \leq 0 $\,, we have
 \bea
 \notag
 \,|\derm h^1_{\cal U \cal U} |  &\les&    C(q_0) \cdot c (\gamma^\prime)  \cdot  c (\delta) \cdot c (\gamma)  \cdot E (  4)  \cdot \frac{ \eps }{ (1+t+|x|)^{1-c (\gamma^\prime)\cdot c (\gamma)  \cdot c (\delta)  \cdot E ( 4)  \cdot \eps} \cdot (1 +|q|)^{1+\gamma^\prime}} \, . \\
\eea

\end{lemma}

\begin{proof}

We have shown in \eqref{decayestimateongoodtermsplusthebadtermwithoutdecayinsourceforh}, that in the exterior, for $\ga > \de $, and $0 \leq \de \leq \frac{1}{4}$, 
          \bea
\notag
 &&  | g^{\alpha\beta} \derm_\alpha \derm_\beta h_{{\underline{L}} {\underline{L}} } |    \\
 \notag
 &\leq&  \begin{cases}  c (\delta) \cdot c (\gamma) \cdot E (  3)  \cdot \frac{\eps }{(1+t+|q|)^{3-3\delta} (1+|q|)^{1+2\de}},\quad\text{when }\quad q>0  \;  ,\\
       \notag
      E ( 3)  \cdot \frac{\eps  \cdot (1+| q |   )^2 }{(1+t+|q|)^{3-3\delta} }  \,\quad\text{when }\quad q<0 \; \end{cases} \\
       \notag
       && + |\derm h_{ {\cal T} {\cal U}} |^2 + |\derm A_{{\cal T}} |^2 \; . \\
       \eea
       We have just upgraded the estimate on  $|\derm h_{ {\cal T} {\cal U}} |^2 $, and obtained for
\beaa
 - 1 \leq \gamma^\prime \leq 0 \; ,
\eeaa
that
 \beaa
 \,|\derm h^1_{ {\cal T} {\cal U}} |^2     &\leq& C(q_0) \cdot c (\gamma^\prime)  \cdot  c (\delta) \cdot c (\gamma)  \cdot E (  4) \cdot   \frac{ \eps }{ (1+t+|q|)^{2-c (\gamma^\prime) \cdot c (\gamma)  \cdot c (\delta)  \cdot E ( 3)  \cdot \eps }  \cdot (1+|q|)^{2+2\gamma^\prime} } \; . \\ 
\eeaa
Also, we had upgraded the estimate on $|\derm A_{{\cal T}} |^2 $ and obtained that for $\ga > \de $, $0 \leq \de < \frac{1}{2}$ and for
   \bea
-1 &\leq& \gamma^\prime <  \ga - \delta \; ,
  \eea
the following estimate
              \bea
   \notag
&& |\derm \ A_{{\cal T}} |^2 \\
   \notag
   &\leq& C(q_0) \cdot c (\gamma^\prime)  \cdot  c (\delta) \cdot c (\gamma)  \cdot E (  4) \cdot     \frac{ \eps }{ (1+t+|q|)^{2-c (\gamma^\prime) \cdot c (\gamma)  \cdot c (\delta)  \cdot E ( 3)  \cdot \eps }  \cdot (1+|q|)^{2+2\gamma^\prime} } \; ,
\eea
Hence, for
  \bea
-1 &\leq& \gamma^\prime <  \ga - \delta, \quad \quad \text{ and} \quad \quad  \gamma^\prime \leq 0  \; ,
  \eea
we have in the exterior $\overline{C}$,
       \bea
\notag
 &&  | g^{\alpha\beta} \derm_\alpha \derm_\beta h_{{\underline{L}} {\underline{L}} } |    \\
 \notag
 &\les&  \begin{cases}  c (\delta) \cdot c (\gamma)  \cdot E (  3)  \cdot \frac{\eps }{(1+t+|q|)^{3-3\delta} (1+|q|)^{1+2\de}},\quad\text{when }\quad q>0  \;  ,\\
       \notag
      E ( 3)  \cdot \frac{\eps  \cdot (1+| q |   )^2 }{(1+t+|q|)^{3-3\delta} }  \,\quad\text{when }\quad q<0 \; \end{cases} \\
       \notag
       && +C(q_0) \cdot c (\gamma^\prime)  \cdot  c (\delta) \cdot c (\gamma)  \cdot E (  4) \cdot \eps \cdot    \frac{ 1}{ (1+t+|q|)^{2-c (\gamma^\prime) \cdot c (\gamma)  \cdot c (\delta)  \cdot E ( 3)  \cdot \eps }  \cdot (1+|q|)^{2+2\gamma^\prime} }\\
        &\les&  C(q_0) \cdot  c (\delta) \cdot c (\gamma)  \cdot E (  3)  \cdot \frac{\eps }{(1+t+|q|)^{3-3\delta} (1+|q|)^{1+2\de}}  \; \\
       \notag
       && +C(q_0) \cdot c (\gamma^\prime)  \cdot  c (\delta) \cdot c (\gamma)  \cdot E (  4) \cdot    \frac{ \eps}{ (1+t+|q|)^{2-c (\gamma^\prime) \cdot c (\gamma)  \cdot c (\delta)  \cdot E ( 3)  \cdot \eps }  \cdot (1+|q|)^{2+2\gamma^\prime} } \; .\\
       \eea
   
       Hence,
       \bea
       \notag
 &&  | g^{\alpha\beta} \derm_\alpha \derm_\beta h_{{\underline{L}} {\underline{L}} } |    \\
 \notag
           &\les&  C(q_0) \cdot  c (\delta) \cdot c (\gamma)  \cdot E (  3)  \cdot \frac{\eps }{(1+t+|q|)^{2} (1+|q|)^{2-\de}}  \; \\
       \notag
       && +C(q_0) \cdot c (\gamma^\prime)  \cdot  c (\delta) \cdot c (\gamma)  \cdot E (  4) \cdot  \frac{ \eps}{ (1+t+|q|)^{2-c (\gamma^\prime) \cdot c (\gamma)  \cdot c (\delta)  \cdot E ( 3)  \cdot \eps }  \cdot (1+|q|)^{2+2\gamma^\prime} }\\
            &\les&   +C(q_0) \cdot c (\gamma^\prime)  \cdot  c (\delta) \cdot c (\gamma)  \cdot E (  4) \cdot   \frac{\eps }{(1+t+|q|)^{2} (1+|q|)^{1 + (1-\de) }}  \; \\
       \notag
       && +C(q_0) \cdot c (\gamma^\prime)  \cdot  c (\delta) \cdot c (\gamma)  \cdot E (  4) \cdot   \frac{  \eps }{ (1+t+|q|)^{2-c (\gamma^\prime) \cdot c (\gamma)  \cdot c (\delta)  \cdot E ( 3)  \cdot \eps }  \cdot (1+|q|)^{1 + (1 + 2\gamma^\prime) } } \; .\\
       \eea
       
       Thus, we choose $2\ga^{\prime} \geq -1$, or differently speaking if $\ga^{\prime} \geq -\frac{1}{2}$, we get
              \bea
       \notag
 &&  | g^{\alpha\beta} \derm_\alpha \derm_\beta h_{{\underline{L}} {\underline{L}} } |    \\
 \notag
           &\les&  C(q_0) \cdot c (\gamma^\prime)  \cdot  c (\delta) \cdot c (\gamma)  \cdot E (  4) \cdot \eps \cdot    \frac{  \eps }{ (1+t+|q|)^{2-c (\gamma^\prime) \cdot c (\gamma)  \cdot c (\delta)  \cdot E ( 3)  \cdot \eps }  \cdot (1+|q|) } \; .\\
       \eea
   Now,  given that 
       \beaa
(1+t) \cdot |  g^{\la\mu} \derm_{\la}   \derm_{\mu}    h^0  | &\les&c (\gamma) \cdot  E (  2)  \cdot \frac{\eps  }{(1+t+|q|)^{2} }\; ,
\eeaa
we get
              \beaa
       \notag
 &&  (1+t) \cdot  | g^{\alpha\beta} \derm_\alpha \derm_\beta h^1_{{\underline{L}} {\underline{L}} } |    \\
 \notag
           &\les&  C(q_0) \cdot c (\gamma^\prime)  \cdot  c (\delta) \cdot c (\gamma)  \cdot E (  4) \cdot \eps \cdot    \frac{  \eps }{ (1+t+|q|)^{1-c (\gamma^\prime) \cdot c (\gamma)  \cdot c (\delta)  \cdot E ( 3)  \cdot \eps }  \cdot (1+|q|) }\\
      && +     c (\gamma) \cdot  E (  2)  \cdot \frac{\eps  }{(1+t+|q|)^{2} } \\
      \notag
      &\leq&  C(q_0) \cdot c (\gamma^\prime)  \cdot  c (\delta) \cdot c (\gamma)  \cdot E (  4)  \cdot \eps \cdot    \frac{  \eps }{ (1+t+|q|)^{1-c (\gamma^\prime) \cdot c (\gamma)  \cdot c (\delta)  \cdot E ( 3)  \cdot \eps }  \cdot (1+|q|)} \; .\\
       \eeaa

   As a result, in the exterior
          \beaa
\notag
&&(1+t ) \cdot \varpi(q) \cdot | g^{\alpha\beta} \derm_\alpha \derm_\beta h^1_{{\underline{L}} {\underline{L}} } |_{L^\infty(\overline{D}_\tau)} \\
 &\les&    C(q_0) \cdot c (\gamma^\prime)  \cdot  c (\delta) \cdot c (\gamma)  \cdot E (  4) \cdot \frac{ \eps \cdot (1+|q|)^{1+\gamma^\prime}}{(1+t+ |q| )^{1-c (\gamma^\prime) \cdot c (\gamma)  \cdot c (\delta)  \cdot E ( 3)  \cdot \eps }  (1+|q|) }  \\
&\les&    C(q_0) \cdot c (\gamma^\prime)  \cdot  c (\delta) \cdot c (\gamma)  \cdot E (  4)   \cdot \frac{ \eps  }{(1+t + |q| )^{1-c (\gamma^\prime) \cdot c (\gamma)  \cdot c (\delta)  \cdot E ( 3)  \cdot \eps }  (1+|q|)^{-\gamma^\prime} } \; . \\
    \eeaa
    Hence, if  $  -\frac{1}{2} \leq \gamma^\prime \leq 0$, then $ - \gamma^\prime \geq 0$ and we have
          \beaa
\notag
&& (1+t )\cdot \varpi(q) \cdot | g^{\alpha\beta} \derm_\alpha \derm_\beta h^1_{{\underline{L}} {\underline{L}} } |_{L^\infty(\overline{D}_\tau)}  \\
&\les&  C(q_0) \cdot  c (\delta) \cdot c (\gamma) \cdot E (  3)  \cdot \frac{ \eps}{(1+t+ |q| )^{1-c (\gamma^\prime) \cdot c (\gamma)  \cdot c (\delta)  \cdot E ( 3)  \cdot \eps }   } \; . \\
    \eeaa
    
Consequently, applying the estimate in Corollary 7.2 of \cite{LR10} (see \eqref{LinfinitynormestimateongradientderivedbyLondbladRodnianski}), we get
     \beaa
   \notag
&& (1+t+|q|) \cdot \varpi(q) \cdot | \derm  h^1 | \\
   \notag
 &\les&   c (\gamma^\prime) \cdot  c (\delta) \cdot c (\gamma) \cdot E (  4) \cdot \eps \, \\
\notag
&& + c (\gamma^\prime) \cdot c (\gamma)  \cdot c (\delta)  \cdot E ( 3)  \cdot \eps  \cdot  \int_0^t  \frac{1}{(1+\tau)} \cdot (1+\tau+|x|) \cdot \|\varpi(q)  \cdot  \derm   h^1 (\tau,\cdot) \|_{L^\infty (\Sigma^{ext}_{\tau} )} d \tau \\
       \notag
&& + \int_0^t  (1+\tau) \cdot \| \varpi(q)  \cdot g^{\la\mu} \derm_{\la}   \derm_{\mu}  h^1 (\tau,\cdot) \|_{L^\infty(\overline{D}_\tau)}  d \tau  \\
\notag
&\les&  c (\gamma^\prime) \cdot  c (\delta) \cdot c (\gamma) \cdot E (  4) \cdot \eps   \\
&& + c (\gamma^\prime) \cdot c (\gamma)  \cdot c (\delta)  \cdot E ( 3)  \cdot \eps  \cdot  \int_0^t  \frac{1}{(1+\tau)} \cdot (1+\tau+|x|) \cdot \|\varpi(q)  \cdot  \derm   h^1  (\tau,\cdot) \|_{L^\infty (\Sigma^{ext}_{\tau} )} d \tau \\
&& +   \int_0^t        C(q_0) \cdot c (\gamma^\prime)  \cdot  c (\delta) \cdot c (\gamma)  \cdot E (  4) \cdot  \frac{ \eps }{(1+\tau+ |q| )^{1-c (\gamma^\prime) \cdot c (\gamma)  \cdot c (\delta)  \cdot E ( 3)  \cdot \eps }   }   d\tau \\
&\les&  c (\gamma^\prime) \cdot  c (\delta) \cdot c (\gamma) \cdot E (  4) \cdot \eps   \\
&& + c (\gamma^\prime) \cdot c (\gamma)  \cdot c (\delta)  \cdot E ( 3)  \cdot \eps  \cdot  \int_0^t  \frac{1}{(1+\tau)} \cdot (1+\tau+|x|) \cdot \|\varpi(q)  \cdot  \derm   h^1 (\tau,\cdot) \|_{L^\infty (\Sigma^{ext}_{\tau} )} d \tau \\
&& +   \int_0^t    C(q_0) \cdot c (\gamma^\prime)  \cdot  c (\delta) \cdot c (\gamma)  \cdot E (  4) \cdot     \frac{  \eps }{ (1+\tau)^{1-c (\gamma^\prime) \cdot c (\gamma)  \cdot c (\delta)  \cdot E ( 4)  \cdot \eps }  }   d\tau \\
&& \text{(given that $E(3) \leq E(4)$) } \\
&\les&  c (\gamma^\prime) \cdot  c (\delta) \cdot c (\gamma) \cdot E (  4) \cdot \eps  +  C(q_0) \cdot c (\gamma^\prime)  \cdot  c (\delta) \cdot c (\gamma)  \cdot E (  4) \cdot  \eps \cdot ( 1+M) \cdot \big[   (1+t)^{c (\gamma^\prime) \cdot c (\gamma)  \cdot c (\delta)  \cdot E ( 4)  \cdot \eps }   \big]^{t}_{0} \\
&& + c (\gamma^\prime) \cdot c (\gamma)  \cdot c (\delta)  \cdot E ( 3)  \cdot \eps  \cdot  \int_0^t  \frac{1}{(1+\tau)} \cdot (1+\tau+|x|) \cdot \|\varpi(q)  \cdot  \derm   h^1 (\tau,\cdot) \|_{L^\infty (\Sigma^{ext}_{\tau} )} d \tau  \; .\\
\eeaa
Hence,

     \beaa
   \notag
&& (1+t+|q|) \cdot \varpi(q) \cdot | \derm  h^1 | \\
   \notag
 &\les&  C(q_0) \cdot c (\gamma^\prime)  \cdot  c (\delta) \cdot c (\gamma)  \cdot E (  4) \cdot  \eps  \cdot    (1+t)^{c (\gamma^\prime) \cdot c (\gamma)  \cdot c (\delta)  \cdot E ( 4)  \cdot \eps }    \\
 && + c (\gamma^\prime) \cdot c (\gamma)  \cdot c (\delta)  \cdot E ( 3)  \cdot \eps  \cdot  \int_0^t  \frac{1}{(1+\tau)} \cdot (1+\tau+|x|) \cdot \|\varpi(q)  \cdot  \derm   h  (\tau,\cdot) \|_{L^\infty (\Sigma^{ext}_{\tau} )} d \tau \; .
\eeaa
Using Grönwall inequality, we get
    \beaa
   \notag
&& (1+t+|q|) \cdot \varpi(q) \cdot | \derm  h^1 | \\
   \notag
   &\les &  C(q_0) \cdot c (\gamma^\prime)  \cdot  c (\delta) \cdot c (\gamma)  \cdot E (  4) \cdot  \eps \cdot    (1+t)^{c (\gamma^\prime) \cdot c (\gamma)  \cdot c (\delta)  \cdot E ( 4)  \cdot \eps }  \cdot (1+t)^{c (\gamma^\prime) \cdot c (\gamma)  \cdot c (\delta)  \cdot E ( 3)  \cdot \eps} \; .
   \eeaa
   
Therefore, we obtain in the exterior,
\beaa
 && (1+t+|q|)| \cdot \varpi(q) \cdot |\derm h^1 | \\
  &\les&    C(q_0) \cdot c (\gamma^\prime)  \cdot  c (\delta) \cdot c (\gamma)  \cdot E (  4) \cdot  \eps  \cdot (1+ t )^{c (\gamma^\prime)  \cdot c (\gamma)  \cdot c (\delta)  \cdot E ( 4)  \cdot \eps } \, .
\eeaa
Finally, for  $\ga > \de $, $0 \leq \de \leq \frac{1}{4}$ and for $  -\frac{1}{2} \leq \gamma^\prime \leq 0 $, which satisfies the condition that
  \bea
-1 &\leq& \gamma^\prime <  \ga - \delta, \quad \quad \text{ and} \quad \quad  \gamma^\prime \leq 0  \; ,
  \eea
 we obtain
 \bea
 \notag
 \,|\derm h^1 |  &\les&    C(q_0) \cdot c (\gamma^\prime)  \cdot  c (\delta) \cdot c (\gamma)  \cdot E (  4)  \cdot \frac{ \eps }{ (1+t+|x|)^{1-c (\gamma^\prime)\cdot c (\gamma)  \cdot c (\delta)  \cdot E ( 4)  \cdot \eps} \cdot (1 +|q|)^{1+\gamma^\prime}} \, . \\
\eea

\end{proof}

\section{Upgraded estimates for the Lie derivatives of the Einstein-Yang-Mills fields}\label{UpgradedestimatesfortheLiederivativesoftheEinstein-Yang-Mills fields}

In this section, we would like to upgrade the decay estimates for the Lie derivatives of the Einstein-Yang-Mills fields. To achieve this, we wish to adapt the argument of Lindblad-Rodnianski in their proof of Proposition 10.2 in \cite{LR10}, however, we are confronted to new challenges (see Subsection \ref{dealingwithLiederivativesupgrade}). In fact, our argument needs an estimate on the commutator term that deals with $A_{\cal T}$ separately (in some sense that we explain in Subsection \ref{dealingwithLiederivativesupgrade}) without referring to the other components, unlike the case for the Einstein vacuum equations and also unlike the case for the Einstein-Maxwell system. We are going to establish the “correct” estimate on the commutator term that allows us to close our argument.

We are going to prove the upgraded estimates on the Lie derivatives by induction: we shall show an upgraded estimate for the zeroth Lie derivative -- based on what we have already proven --\,, will then assume that the following upgraded estimate holds for all $K$ such that $|K| \leq |J| - 1 = k -1 \in \N $ and then, we prove the estimate for $|K| = |J| = k$, and thereby we would have proven the following lemma for all $K$ (by induction on $|K|$):

\begin{lemma}\label{upgradedestimateonLiederivativesoffields}
Let $M \leq \eps$\,. In the Lorenz and harmonic gauges, the Einstein-Yang-Mills fields satisfy for $\ga \geq 3 \de $\,, and $0 < \de \leq \frac{1}{4}$\,, and for any $|K|  \in \N $\,,
            \bea\label{formulaforinductionhypothesisonA}
 \notag
&& |\derm  ( \Lie_{Z^K} A) (t,x)  |    \\
\notag
&\leq&   C(q_0)   \cdot  c (\delta) \cdot c (\gamma) \cdot C(|K|) \cdot E ( |K| + 4)  \cdot \frac{\eps }{(1+t+|q|)^{1-      c (\gamma)  \cdot c (\delta)  \cdot c(|K|) \cdot E ( |K|+ 4) \cdot \eps } \cdot (1+|q|)^{1+\gamma - 2\de }}    \; ,\\
      \eea
      
            \bea\label{formulaforinductionhypothesisonh^1}
 \notag
&& |\derm  ( \Lie_{Z^K} h^1) (t,x)  |   \\
\notag
 &\leq&   C(q_0)   \cdot  c (\delta) \cdot c (\gamma) \cdot C(|K|) \cdot E (  |K|  +4)  \cdot \frac{\eps   }{(1+t+|q|)^{1-   c (\gamma)  \cdot c (\delta)  \cdot c(|K|) \cdot E ( |K|+ 4)\cdot  \eps } \cdot (1+|q|)}     \; , \\
      \eea
      
      where $c (\gamma)$, $c (\delta)$, $c(|K|)$ are constants that depend on $\gamma$, on $\de$, and on $|K|$, respectively.
 \end{lemma}
  
\begin{proof}
The way we are going to prove this estimate for $|K| = |J| = k$, is by using what we have already established in Lemma \eqref{GronwallinequalitiesoncomponentsofgradientofAandgradientofh1withweight}: that is for $\gamma^\prime$ such that $-1 \leq \gamma^\prime < \gamma - \delta$, and $ \delta <    1/2 $, we have for all $U, V\in \cal U$,
     \bea
   \notag
&& (1+t+|x|) \cdot |\varpi(q) \cdot \derm ( \Lie_{Z^J} A)_V (t,x)| \\
   \notag
 &\les&   c (\gamma^\prime) \cdot  c (\delta) \cdot c (\gamma) \cdot C ( |J|  ) \cdot E ( |J| + 4) \cdot \eps \, \\
\notag
&& + c (\gamma^\prime) \cdot c (\gamma)  \cdot c (\delta)  \cdot E ( 3)  \cdot \eps  \cdot  \int_0^t  \frac{1}{(1+\tau)} \cdot (1+\tau+|x|) \cdot \|\varpi(q)  \cdot  \derm  ( \Lie_{Z^J} A)_V (\tau,\cdot) \|_{L^\infty (\Sigma^{ext}_{\tau} )} d \tau \\
       \notag
&& + \int_0^t  (1+\tau) \cdot \| \varpi(q)  \cdot g^{\la\mu} \derm_{\la}   \derm_{\mu}  (  \Lie_{Z^J} A)_V (\tau,\cdot) \|_{L^\infty(\overline{D}_\tau)} d\tau \; ,\\
\eea
and 
     \bea
   \notag
&& (1+t+|x|) \cdot |\varpi(q) \cdot \derm ( \Lie_{Z^J} h^1)_{UV} (t,x)| \\
   \notag
 &\les&   c (\gamma^\prime) \cdot  c (\delta) \cdot c (\gamma) \cdot C ( |J|  ) \cdot E ( |J| + 4) \cdot \eps \, \\
\notag
&& + c (\gamma^\prime) \cdot c (\gamma)  \cdot c (\delta)  \cdot E ( 3)  \cdot \eps  \cdot  \int_0^t  \frac{1}{(1+\tau)} \cdot (1+\tau+|x|) \cdot \|\varpi(q)  \cdot  \derm  ( \Lie_{Z^J} h^1)_{UV} (\tau,\cdot) \|_{L^\infty (\Sigma^{ext}_{\tau} )} d \tau \\
       \notag
&& + \int_0^t  (1+\tau) \cdot \| \varpi(q)  \cdot g^{\la\mu} \derm_{\la}   \derm_{\mu}  (  \Lie_{Z^J} h^1)_{UV} (\tau,\cdot) \|_{L^\infty(\overline{D}_\tau)} d\tau  \; .\\
\eea

However, for this, we need to estimate the commutator term for the Lie derivatives and the source terms for the Lie derivatives. We start first by showing that the initial step for the induction is satisfied. We will assume from now on that $M \leq \eps$\;, as stated in Lemma \ref{upgradedestimateonLiederivativesoffields} that we want to prove.

 \subsection{The initialisation for the Induction for the Einstein-Yang-Mills fields}\

\begin{lemma}\label{Step1fortheinductiononbothAandh1}
We have for $\ga \geq 3 \de $, and $0 < \de \leq \frac{1}{4}$,
   \bea
 \notag
 \,|\derm A |  &\les&    C(q_0)   \cdot  c (\delta) \cdot c (\gamma)  \cdot E (  4) \cdot  \frac{\eps }{ (1+t+|q|)^{1-c (\gamma)  \cdot c (\delta)  \cdot E ( 4)  \cdot \eps} \cdot (1 +|q|)^{1+\ga -2\de}} \, , \\
\eea
 \bea
 \notag
 \,|\derm h^1 |  &\les&    C(q_0)   \cdot  c (\delta) \cdot c (\gamma)  \cdot E (  4)  \cdot \frac{\eps  }{ (1+t+|q|)^{1- c (\gamma)  \cdot c (\delta)  \cdot E ( 4)  \cdot \eps} \cdot (1 +|q|)} \, . \\
\eea

\end{lemma}

 \begin{proof}
 
We showed in Lemma \ref{upgradedestimateonthegradientofthefullcompoenentsofpoentialA}, that in the Lorenz gauge and in wave coordinates, we have for $\ga > \de $, and $0 \leq \de \leq \frac{1}{4}$, and for  
  \beaa
\de &\leq&  \gamma^\prime <  \ga - \delta \; ,
  \eeaa
  the following estimate
   \bea
 \notag
 \,|\derm A |  &\les&    C(q_0) \cdot c (\gamma^\prime)  \cdot  c (\delta) \cdot c (\gamma)  \cdot E (  4)  \frac{\eps }{ (1+t+|q|)^{1-c (\gamma^\prime)\cdot c (\gamma)  \cdot c (\delta)  \cdot E ( 4)  \cdot \eps} \cdot (1 +|q|)^{1+\gamma^\prime}} \, . \\
\eea
Thus, by taking $ \gamma^\prime =  \ga - 2\delta$ where $0 < \de $, with $\ga \geq 3\de$, we satisfy the conditions and we obtain
   \bea
 \notag
 \,|\derm A |  &\les&    C(q_0)   \cdot  c (\delta) \cdot c (\gamma)  \cdot E (  4)  \frac{\eps }{ (1+t+|q|)^{1-c (\gamma)  \cdot c (\delta)  \cdot E ( 4)  \cdot \eps} \cdot (1 +|q|)^{1+\ga -2\de}} \, . \\
\eea

Whereas to the metric, we showed that for $\ga > \de $, $0 \leq \de \leq \frac{1}{4}$ and for $  -\frac{1}{2} \leq \gamma^\prime \leq 0 $, we have
 \bea
 \notag
 \,|\derm h^1 |  &\les&    C(q_0) \cdot c (\gamma^\prime)  \cdot  c (\delta) \cdot c (\gamma)  \cdot E (  4)  \cdot \frac{\eps  }{ (1+t+|q|)^{1-c (\gamma^\prime)\cdot c (\gamma)  \cdot c (\delta)  \cdot E ( 4)  \cdot \eps} \cdot (1 +|q|)^{1+\gamma^\prime}} \, . \\
\eea
Thus, by taking another $\gamma^\prime$\,, different than the one used for $A$\,, where this time we choose $ \ga^\prime = 0$\,, we get
 \bea
 \notag
 \,|\derm h^1 |  &\les&    C(q_0)   \cdot  c (\delta) \cdot c (\gamma)  \cdot E (  4)  \cdot \frac{\eps }{ (1+t+|q|)^{1- c (\gamma)  \cdot c (\delta)  \cdot E ( 4)  \cdot \eps} \cdot (1 +|q|)} \, . \\
\eea
 
 \end{proof}

 \subsection{The commutator term for the Einstein-Yang-Mills fields}\

 Now, we are left to study the terms  $\int_0^t (1+\tau) \cdot \| \varpi(q)   g^{\la\mu} \derm_{\la}   \derm_{\mu}   \Lie_{Z^J} A (\tau,\cdot) \|_{L^\infty(\overline{D}_\tau)} d\tau $ and $\int_0^t (1+\tau) \cdot \| \varpi(q)   g^{\la\mu} \derm_{\la}   \derm_{\mu}   \Lie_{Z^J} h^1 (\tau,\cdot) \|_{L^\infty(\overline{D}_\tau)} d\tau $. For this, we actually need a slightly more refined commutation formula than the one derived by Lindblad-Rodnianski in Proposition 5.3 in \cite{LR10}, in order to be able to upgrade the estimates for the derivatives of $A$, which require that we deal first with $A_{\cal T}$ separately, in order to then transfer the estimate to $A_{\cal U}$ (unlike the case for the metric $h^1$ where the induction hypothesis suffices to upgrade the estimate for $h^1_{\cal U \cal U} $ directly).

\begin{lemma}\label{theexactcommutatortermwithdependanceoncomponents}
We have for all $I$\;, or any $U \in {\cal U}$\;,
\bea
\notag
&& \Lie_{Z^I}  ( g^{\la\mu} \derm_{\la}   \derm_{\mu}     \Phi_{U} ) - g^{\la\mu}    \derm_{\la}   \derm_{\mu}  (  \Lie_{Z^I} \Phi_{U}  )  \\
\notag
&=&  \sum_{I_1 + I_2 = I, \; I_2 \neq I} \hat{c}(I_1) \cdot m^{\la\mu} \cdot \derm_{\la}   \derm_{\mu} (  \Lie_{Z^{I_2}}  \Phi_{U} ) \\
\notag
&& +  \sum_{I_2 + I_4 + I_5 + I_6  = I, \; I_2 \neq I }   \hat{c}(I_6) \cdot   \hat{c}(I_5)   \cdot \Big(   \frac{1}{4}   ( \Lie_{Z^{I_4}}   H)_{L  L}    \cdot  \derm_{\underline{L}}   \derm_{ \underline{L}}   (  \Lie_{Z^{I_2}}  \Phi_{U} ) \\
\notag
&&  + \frac{1}{4}  ( \Lie_{Z^{I_4}}   H)_{L  \underline{L}}    \cdot  \derm_{\underline{L}}   \derm_{L}   (  \Lie_{Z^{I_2}}  \Phi_{U} )    - \frac{1}{2}   ( \Lie_{Z^{I_4}}   H)_{L e_A }    \cdot  \derm_{\underline{L}}   \derm_{e_A}    (  \Lie_{Z^{I_2}}  \Phi_{U} )    \\
   \notag
 && - \frac{1}{2} m^{\mu \b}   ( \Lie_{Z^{I_4}}   H)_{\underline{L}\b}  \cdot    \derm_{L}   \derm_{\mu}     (  \Lie_{Z^{I_2}}  \Phi_{U} )   +    m^{\mu \b}   ( \Lie_{Z^{I_4}}   H)_{e_A \b}   \cdot   \derm_{e_A}   \derm_{\mu}    (  \Lie_{Z^{I_2}}  \Phi_{U} )   \Big)  \; .\\
 \eea
 
 \end{lemma}
 
 \begin{proof}
 
We compute for any $U \in {\cal U}$, 
\beaa
&& \Lie_{Z^{I}}  ( g^{\la\mu} \derm_{\la}   \derm_{\mu}     \Phi_{U} )  \\
&=&  \Lie_{Z^{I}}  \Big( (m^{\la\mu} + H^{\la\mu} ) \cdot \derm_{\la}   \derm_{\mu}     \Phi_{U} \Big) \\
&=& \sum_{I_1 + I_2 = I} \Big( \Lie_{Z^{I_1}} m^{\la\mu} + \Lie_{Z^{I_1}}  H^{\la\mu} \Big)  \cdot  \Lie_{Z^{I_2}}  \derm_{\la}   \derm_{\mu}     \Phi_{U} \\
&=&\sum_{I_1 + I_2 = I} \hat{c}(I_1) \cdot m^{\la\mu} \cdot \Lie_{Z^{I_2}}  \derm_{\la}   \derm_{\mu}  \Phi_{U}  + \sum_{I_1 + I_2 = I} ( \Lie_{Z^{I_1}}  H^{\la\mu} ) \cdot   \Lie_{Z^{I_2}}  \derm_{\la}   \derm_{\mu}  \Phi_{U} \; .
\eeaa

We note that 
\beaa
&& \Lie_{Z^{I_1}}  H^{\la\mu} \\
 &=&  \Lie_{Z^{I_1}}(  m^{\la\a} \cdot m^{\mu\b} \cdot  H_{\a\b}  ) \\
  &=&\sum_{I_3 + I_4   = I_1 }  \Lie_{Z^{I_3}} (  m^{\la\a}  \cdot m^{\mu\b} ) \cdot  \Lie_{Z^{I_4}} (  H_{\a\b}  ) = \sum_{I_3 + I_4   = I_1 } \Big( \sum_{I_5 + I_6   = I_3 } (  \Lie_{Z^{I_5}}  m^{\la\a}  \cdot  \Lie_{Z^{I_6}}  m^{\mu\b} ) \cdot  \Lie_{Z^{I_4}} (  H_{\a\b}  ) \Big)    \\
  &=&  \sum_{ I_4 + I_5 + I_6  = I_1 } (   \hat{c}(I_5) \cdot  m^{\la\a} )  \cdot (  \hat{c}(I_6) \cdot m^{\mu\b} )  \cdot (   \Lie_{Z^{I_4}}    H_{\a\b}  )  \; .
\eeaa
With the notation $ m^{\la\a}  \cdot m^{\mu\b}  \cdot ( \Lie_{Z^I}   H_{\a\b}  )  =  ( \Lie_{Z^I}   H)^{\la\mu}  $, and using the fact that the Lie derivatives in the direction of Minkowski vector fields commute with $\derm$, we can then write
\beaa
&& \Lie_{Z^I}  ( g^{\la\mu} \derm_{\la}   \derm_{\mu}     \Phi_{U} )  \\
&=& \sum_{I_1 + I_2 = I} \hat{c}(I_1) \cdot m^{\la\mu} \cdot \derm_{\la}   \derm_{\mu} (  \Lie_{Z^{I_2}}  \Phi_{U} ) \\
&& +  \sum_{I_2 + I_4 + I_5 + I_6  = I } (   \hat{c}(I_5) \cdot  m^{\la\a} )  \cdot (  \hat{c}(I_6) \cdot m^{\mu\b} )  \cdot (   \Lie_{Z^{I_4}}    H_{\a\b}  ) \cdot  \derm_{\la}   \derm_{\mu} (  \Lie_{Z^{I_2}}  \Phi_{U} )   \\
&=& \sum_{I_1 + I_2 = I} \hat{c}(I_1) \cdot m^{\la\mu} \cdot \derm_{\la}   \derm_{\mu} (  \Lie_{Z^{I_2}}  \Phi_{U} ) \\
&& +  \sum_{I_2 + I_4 + I_5 + I_6  = I }   \hat{c}(I_5) \cdot   \hat{c}(I_6)   \cdot (   \Lie_{Z^{I_4}}     H)^{\la\mu} \cdot  \derm_{\la}   \derm_{\mu} (  \Lie_{Z^{I_2}}  \Phi_{U} )    \\
&=&    m^{\la\mu} \cdot \derm_{\la}   \derm_{\mu} (  \Lie_{Z^{I}}  \Phi_{U} ) +       H^{\la\mu} \cdot  \derm_{\la}   \derm_{\mu} (  \Lie_{Z^{I}}  \Phi_{U} )    \\
&& + \sum_{I_1 + I_2 = I, \; I_2 \neq I} \hat{c}(I_1) \cdot m^{\la\mu} \cdot \derm_{\la}   \derm_{\mu} (  \Lie_{Z^{I_2}}  \Phi_{U} ) \\
&& +  \sum_{I_2 + I_4 + I_5 + I_6  = I, \; I_2 \neq I }   \hat{c}(I_5) \cdot   \hat{c}(I_6)   \cdot (   \Lie_{Z^{I_4}}     H)^{\la\mu} \cdot  \derm_{\la}   \derm_{\mu} (  \Lie_{Z^{I_2}}  \Phi_{U} )    \; .
\eeaa

Therefore, 
\bea
\notag
&& \Lie_{Z^I}  ( g^{\la\mu} \derm_{\la}   \derm_{\mu}     \Phi_{U} ) - g^{\la\mu}    \derm_{\la}   \derm_{\mu}  (  \Lie_{Z^I} \Phi_{U}  )  \\
\notag
&=&  \sum_{I_1 + I_2 = I, \; I_2 \neq I} \hat{c}(I_1) \cdot m^{\la\mu} \cdot \derm_{\la}   \derm_{\mu} (  \Lie_{Z^{I_2}}  \Phi_{U} ) \\
\notag
&& +  \sum_{I_2 + I_4 + I_5 + I_6  = I, \; I_2 \neq I }   \hat{c}(I_5) \cdot   \hat{c}(I_6)   \cdot (   \Lie_{Z^{I_4}}     H)^{\la\mu} \cdot  \derm_{\la}   \derm_{\mu} (  \Lie_{Z^{I_2}}  \Phi_{U} )   \; . \\
\eea

Now, to estimate the term with  $H$, we decompose the contractions in the null frame $\cal U$, 
\beaa
 && ( \Lie_{Z^{I_4}}   H)^{\la\mu}  \cdot  \derm_{\la}   \derm_{\mu}    (  \Lie_{Z^{I_2}}  \Phi_{U} )  \\
 & =&  ( \Lie_{Z^{I_4}}   H)^{\underline{L}\mu}  \cdot   \derm_{\underline{L}}   \derm_{\mu}   (  \Lie_{Z^{I_2}}  \Phi_{U} )  +  ( \Lie_{Z^{I_4}}   H)^{L \mu}   \cdot  \derm_{L}   \derm_{\mu}     (  \Lie_{Z^{I_2}}  \Phi_{U} )  \\
 && +  ( \Lie_{Z^{I_4}}   H)^{e_A \mu}   \cdot  \derm_{e_A}   \derm_{\mu}   (  \Lie_{Z^{I_2}}  \Phi_{U} )  \\
 & =&  m^{\underline{L} \a} m^{\mu \b}  ( \Lie_{Z^{I_4}}   H)_{\a\b}   \cdot  \derm_{\underline{L}}   \derm_{\mu}   (  \Lie_{Z^{I_2}}  \Phi_{U} )   \\
 && +m^{L \a} m^{\mu \b}   ( \Lie_{Z^{I_4}}   H)_{\a\b}   \cdot  \derm_{L}   \derm_{\mu}  (  \Lie_{Z^{I_2}}  \Phi_{U} )  +   m^{e_A \a} m^{\mu \b}   ( \Lie_{Z^{I_4}}   H)_{\a\b}   \cdot  \derm_{e_A}   \derm_{\mu}  (  \Lie_{Z^{I_2}}  \Phi_{U} )   \\
 & =&  m^{\underline{L} L} m^{\mu \b}  ( \Lie_{Z^{I_4}}   H)_{L\b}    \cdot  \derm_{\underline{L}}   \derm_{\mu}   (  \Lie_{Z^{I_2}}  \Phi_{U} )  \\
 && +m^{L \underline{L}} m^{\mu \b}   ( \Lie_{Z^{I_4}}   H)_{\underline{L}\b}   \cdot  \derm_{L}   \derm_{\mu} (  \Lie_{Z^{I_2}}  \Phi_{U} )   +   m^{e_A e_A} m^{\mu \b}   ( \Lie_{Z^{I_4}}   H)_{e_A \b}     \cdot \derm_{e_A}   \derm_{\mu}   (  \Lie_{Z^{I_2}}  \Phi_{U} )  \; .
 \eeaa
 
 Using the fact that  $m^{L \underline{L}} = - \frac{1}{2} $ and $m^{e_A e_B}  =  \delta_{AB} $, we obtain
\beaa
 && ( \Lie_{Z^{I_4}}   H)^{\la\mu}  \cdot  \derm_{\la}   \derm_{\mu}   (  \Lie_{Z^{I_2}}  \Phi_{U} )  \\
 & =&   - \frac{1}{2} m^{\mu \b}  ( \Lie_{Z^{I_4}}   H)_{L\b}    \cdot  \derm_{\underline{L}}   \derm_{\mu}   (  \Lie_{Z^{I_2}}  \Phi_{U} )   \\
 && - \frac{1}{2} m^{\mu \b}   ( \Lie_{Z^{I_4}}   H)_{\underline{L}\b}  \cdot    \derm_{L}   \derm_{\mu}     (  \Lie_{Z^{I_2}}  \Phi_{U} )  +    m^{\mu \b}   ( \Lie_{Z^{I_4}}   H)_{e_A \b}   \cdot   \derm_{e_A}   \derm_{\mu}   (  \Lie_{Z^{I_2}}  \Phi_{U} )  \\
  & =&   - \frac{1}{2} m^{\mu L }  ( \Lie_{Z^{I_4}}   H)_{L  L}    \cdot  \derm_{\underline{L}}   \derm_{\mu}   (  \Lie_{Z^{I_2}}  \Phi_{U} )   - \frac{1}{2} m^{\mu  \underline{L}}  ( \Lie_{Z^{I_4}}   H)_{L  \underline{L}}    \cdot  \derm_{\underline{L}}   \derm_{\mu}   (  \Lie_{Z^{I_2}}  \Phi_{U} )   \\
  && - \frac{1}{2} m^{\mu e_A}  ( \Lie_{Z^{I_4}}   H)_{L e_A }    \cdot  \derm_{\underline{L}}   \derm_{\mu}   (  \Lie_{Z^{I_2}}  \Phi_{U} )   \\
 && - \frac{1}{2} m^{\mu \b}   ( \Lie_{Z^{I_4}}   H)_{\underline{L}\b}  \cdot    \derm_{L}   \derm_{\mu}    (  \Lie_{Z^{I_2}}  \Phi_{U} )   +    m^{\mu \b}   ( \Lie_{Z^{I_4}}   H)_{e_A \b}   \cdot   \derm_{e_A}   \derm_{\mu}   (  \Lie_{Z^{I_2}}  \Phi_{U} )   \\
  & =&   - \frac{1}{2} m^{ \underline{L} L }  ( \Lie_{Z^{I_4}}   H)_{L  L}    \cdot  \derm_{\underline{L}}   \derm_{ \underline{L}}   (  \Lie_{Z^{I_2}}  \Phi_{U} )   - \frac{1}{2} m^{L  \underline{L}}  ( \Lie_{Z^{I_4}}   H)_{L  \underline{L}}    \cdot  \derm_{\underline{L}}   \derm_{L}   (  \Lie_{Z^{I_2}}  \Phi_{U} )   \\
  && - \frac{1}{2} m^{e_A e_A}  ( \Lie_{Z^{I_4}}   H)_{L e_A }    \cdot  \derm_{\underline{L}}   \derm_{e_A} (  \Lie_{Z^{I_2}}  \Phi_{U} )    \\
 && - \frac{1}{2} m^{\mu \b}   ( \Lie_{Z^{I_4}}   H)_{\underline{L}\b}  \cdot    \derm_{L}   \derm_{\mu}   (  \Lie_{Z^{I_2}}  \Phi_{U} )  +    m^{\mu \b}   ( \Lie_{Z^{I_4}}   H)_{e_A \b}   \cdot   \derm_{e_A}   \derm_{\mu}    (  \Lie_{Z^{I_2}}  \Phi_{U} )  \; .
 \eeaa
 Thus,
 \bea
 \notag
 && ( \Lie_{Z^{I_4}}   H)^{\la\mu}  \cdot  \derm_{\la}   \derm_{\mu}    (  \Lie_{Z^{I_2}}  \Phi_{U} )   \\
  \notag
  & =&   \frac{1}{4}   ( \Lie_{Z^{I_4}}   H)_{L  L}    \cdot  \derm_{\underline{L}}   \derm_{ \underline{L}}   (  \Lie_{Z^{I_2}}  \Phi_{U} )   + \frac{1}{4}  ( \Lie_{Z^{I_4}}   H)_{L  \underline{L}}    \cdot  \derm_{\underline{L}}   \derm_{L}   (  \Lie_{Z^{I_2}}  \Phi_{U} )   \\
   \notag
  && - \frac{1}{2}   ( \Lie_{Z^{I_4}}   H)_{L e_A }    \cdot  \derm_{\underline{L}}   \derm_{e_A}    (  \Lie_{Z^{I_2}}  \Phi_{U} )    \\
   \notag
 && - \frac{1}{2} m^{\mu \b}   ( \Lie_{Z^{I_4}}   H)_{\underline{L}\b}  \cdot    \derm_{L}   \derm_{\mu}     (  \Lie_{Z^{I_2}}  \Phi_{U} )   +    m^{\mu \b}   ( \Lie_{Z^{I_4}}   H)_{e_A \b}   \cdot   \derm_{e_A}   \derm_{\mu}    (  \Lie_{Z^{I_2}}  \Phi_{U} )   \; . \\
 \eea
  We get then the desired result.
  
 \end{proof}
 
 \begin{lemma}\label{commutationformaulamoreprecisetoconservegpodcomponentsstructure}
 Let  $\Phi_{\mu}$ be a one-tensor valued in the Lie algebra or a scalar. Then, we have for all $I$, and for any $V \in \cal T$, 
 \bea\label{morerefinedcommutationformularusingonlytangentialcomponennts}
\notag
&&| \Lie_{Z^I}  ( g^{\la\mu} \derm_{\la}   \derm_{\mu}     \Phi_{V} ) - g^{\la\mu}    \derm_{\la}   \derm_{\mu}  (  \Lie_{Z^I} \Phi_{V}  ) |  \\
\notag
&\les&  \sum_{|K| < |I| }  | g^{\la\mu} \cdot \derm_{\la}   \derm_{\mu} (  \Lie_{Z^{K}}  \Phi_{V} ) | \\
\notag
&& +  \sum_{|J| + |K| \leq |I|, \; |K| < |I| }  \Big(   |   ( \Lie_{Z^{J}}   H)_{L  L} |   \cdot    \frac{1}{(1+t+|q|)} \cdot  \sum_{|M| \leq |K|+1}  | \derm ( \Lie_{Z^M}  \Phi )  |  \\
   \notag
 && +    |   ( \Lie_{Z^{J}}   H)_{L  L}   | \cdot   \frac{1}{(1+|q|)} \cdot   \sum_{|M| \leq |K|+1}    \sum_{ V^\prime \in \cal T } | \derm   ( \Lie_{Z^M}  \Phi _{V^\prime} ) |   \\
\notag
 && +   | ( \Lie_{Z^{J}}   H)_{L  \underline{L}} |   \cdot | \derm_{L}    \derm_{\underline{L}}  (  \Lie_{Z^{K}}  \Phi_{V} )  | \\
\notag
    &&  +    \frac{1}{(1 + t + |q|)  }  \cdot   |  ( \Lie_{Z^{J}}   H)_{L e_A }   | \cdot | \sum_{|M| \leq |K|+1} | \derm (  \Lie_{Z^{M}}  \Phi )| \\
    \notag
 && +       | m^{\mu \b}   ( \Lie_{Z^{J}}   H)_{\underline{L}\b}  \cdot    \derm_{L}   \derm_{\mu}     (  \Lie_{Z^{K}}  \Phi_{V} )  | + |   m^{\mu \b}   ( \Lie_{Z^{J}}   H)_{e_A \b}   \cdot   \derm_{e_A}   \derm_{\mu}    (  \Lie_{Z^{K}}  \Phi_{V} )   | \Big) \; .\\
  \eea
Furthermore, let $\Phi_{\mu\nu }$ be a two-tensor (or a one-tensor $\Phi_{\mu}$) valued either in the Lie algebra or a scalar. Then, we have for any $U,  V \in \cal U$,
 \bea\label{commutationformaulaforallcomponentsasknowninlitter}
\notag
&&| \Lie_{Z^I}  ( g^{\la\mu} \derm_{\la}   \derm_{\mu}     \Phi_{UV} ) - g^{\la\mu}    \derm_{\la}   \derm_{\mu}  (  \Lie_{Z^I} \Phi_{UV}  ) |  \\
\notag
&\les&  \sum_{|K| < |I| }  | g^{\la\mu} \cdot \derm_{\la}   \derm_{\mu} (  \Lie_{Z^{K}}  \Phi_{UV} ) | \\
\notag
&& +  \sum_{|J| + |K| \leq |I|, \; |K| < |I| }  \Big(      |   ( \Lie_{Z^{J}}   H)_{L  L}   | \cdot   \frac{1}{(1+|q|)} \cdot   \sum_{|M| \leq |K|+1}    | \derm   ( \Lie_{Z^M}  \Phi) |   \\
\notag
 && +   | ( \Lie_{Z^{J}}   H)_{L  \underline{L}} |   \cdot | \derm_{L}    \derm_{\underline{L}}  (  \Lie_{Z^{K}}  \Phi_{UV} )  | \\
\notag
    &&  +    \frac{1}{(1 + t + |q|)  }  \cdot   |  ( \Lie_{Z^{J}}   H)_{L e_A }   | \cdot | \sum_{|M| \leq |K|+1} | \derm (  \Lie_{Z^{M}}  \Phi )| \\
    \notag
 && +       | m^{\mu \b}   ( \Lie_{Z^{J}}   H)_{\underline{L}\b}  \cdot    \derm_{L}   \derm_{\mu}     (  \Lie_{Z^{K}}  \Phi_{UV} )  | + |   m^{\mu \b}   ( \Lie_{Z^{J}}   H)_{e_A \b}   \cdot   \derm_{e_A}   \derm_{\mu}    (  \Lie_{Z^{K}}  \Phi_{UV} )   | \Big) \; .\\
  \eea

 \end{lemma}
 
 \begin{proof}
 
 From Lemma \ref{theexactcommutatortermwithdependanceoncomponents}, we immediately get hat for any $U \in \cal U$,
 \bea
\notag
&&| \Lie_{Z^I}  ( g^{\la\mu} \derm_{\la}   \derm_{\mu}     \Phi_{U} ) - g^{\la\mu}    \derm_{\la}   \derm_{\mu}  (  \Lie_{Z^I} \Phi_{U}  ) |  \\
\notag
&\les&  \sum_{ |K| < |I| }  | m^{\la\mu} \cdot \derm_{\la}   \derm_{\mu} (  \Lie_{Z^{K}}  \Phi_{U} ) | \\
\notag
&& +  \sum_{J, \, K, \,  |J| + |K| \leq |I|, \; |K| < |I| }  \Big( |   ( \Lie_{Z^{J}}   H)_{L  L}    \cdot  \derm_{\underline{L}}   \derm_{ \underline{L}}   (  \Lie_{Z^{K}}  \Phi_{U} ) | \\
\notag
&&  + | ( \Lie_{Z^{J}}   H)_{L  \underline{L}}    \cdot  \derm_{\underline{L}}   \derm_{L}   (  \Lie_{Z^{K}}  \Phi_{U} )  |  + |  ( \Lie_{Z^{J}}   H)_{L e_A }    \cdot  \derm_{\underline{L}}   \derm_{e_A}    (  \Lie_{Z^{K}}  \Phi_{U} )  |   \\
   \notag
 && + | m^{\mu \b}   ( \Lie_{Z^{J}}   H)_{\underline{L}\b}  \cdot    \derm_{L}   \derm_{\mu}     (  \Lie_{Z^{K}}  \Phi_{U} )  | + |   m^{\mu \b}   ( \Lie_{Z^{J}}   H)_{e_A \b}   \cdot   \derm_{e_A}   \derm_{\mu}    (  \Lie_{Z^{K}}  \Phi_{U} )   | \Big)  \; ,\\
 \eea
and also, from the proof of Lemma \ref{theexactcommutatortermwithdependanceoncomponents}, we get, in addition, that for any $U \in \cal U$\;,
  
   \bea
\notag
&&| \Lie_{Z^I}  ( g^{\la\mu} \derm_{\la}   \derm_{\mu}     \Phi_{U} ) - g^{\la\mu}    \derm_{\la}   \derm_{\mu}  (  \Lie_{Z^I} \Phi_{U}  ) |  \\
\notag
&\les&  \sum_{ |K| < |I| }  | g^{\la\mu} \cdot \derm_{\la}   \derm_{\mu} (  \Lie_{Z^{K}}  \Phi_{U} ) | \\
\notag
&& +  \sum_{J, \, K, \,  |J| + |K| \leq |I|, \; |K| < |I| }  \Big( |   ( \Lie_{Z^{J}}   H)_{L  L}    \cdot  \derm_{\underline{L}}   \derm_{ \underline{L}}   (  \Lie_{Z^{K}}  \Phi_{U} ) | \\
\notag
&&  + | ( \Lie_{Z^{J}}   H)_{L  \underline{L}}    \cdot  \derm_{\underline{L}}   \derm_{L}   (  \Lie_{Z^{K}}  \Phi_{U} )  |  + |  ( \Lie_{Z^{J}}   H)_{L e_A }    \cdot  \derm_{\underline{L}}   \derm_{e_A}    (  \Lie_{Z^{K}}  \Phi_{U} )  |   \\
   \notag
 && + | m^{\mu \b}   ( \Lie_{Z^{J}}   H)_{\underline{L}\b}  \cdot    \derm_{L}   \derm_{\mu}     (  \Lie_{Z^{K}}  \Phi_{U} )  | + |   m^{\mu \b}   ( \Lie_{Z^{J}}   H)_{e_A \b}   \cdot   \derm_{e_A}   \derm_{\mu}    (  \Lie_{Z^{K}}  \Phi_{U} )   | \Big)  \; .\\
 \eea
 
  Now, we are going to estimate the terms one by one.
  
 \textbf{The term $   ( \Lie_{Z^{I_4}}   H)_{L  L}    \cdot  \derm_{\underline{L}}   \derm_{ \underline{L}}   (  \Lie_{Z^{I_2}}  \Phi_{U} ) $}:
 
 We have, see Lemma 3.3 in \cite{BFJST1}, that for any tensor $\Psi_{UV}$ and for any $U \in \cal U$ and for any $V \in \cal T$,
 \bea\label{estimateforgradientoftensorestaimetedbyLiederivativesoftensorwithrefinedcomponentsforthebadfactor}
 \notag
 |\derm \Psi_{UV} | \les \sum_{|I| \leq 1}  \frac{1}{(1+t+|q|)} \cdot | \Lie_{Z^I} \Psi |  +  \sum_{U^\prime \in \cal U,  V^\prime \in \cal T } \sum_{|I| \leq 1}  \frac{1}{(1+|q|)} \cdot  | \Lie_{Z^I} \Psi_{U^\prime V^\prime} | \; . \\
 \eea
 We can then apply the above for $\Psi_{UV} = \derm_{U} (  \Lie_{Z^{K}}   \Phi_{V} ) $, and we get
  \beaa
 && |\derm \derm_{U} (  \Lie_{Z^{K}}   \Phi_{V} ) | \\
 &\les& \sum_{|M| \leq 1}  \frac{1}{(1+t+|q|)} \cdot | \Lie_{Z^M} \derm (  \Lie_{Z^{K}}   \Phi)  |  +    \sum_{U^\prime \in \cal U,  V^\prime \in \cal T } \sum_{|M| \leq 1}  \frac{1}{(1+|q|)} \cdot  | \Lie_{Z^M} \derm_{U^\prime } (  \Lie_{Z^{K}}   \Phi_{V^\prime } ) | \; .
 \eeaa
 Using the commutation of the Lie derivative in the direction of Minkowski vector fields with the covariant derivative of Minkowski $\derm$, we obtain 
  \bea
  \notag
 && | \derm_{\underline{L}}   \derm_{ \underline{L}}   (  \Lie_{Z^{K}}  \Phi_{V} )   | \\
   \notag
 &\les &\sum_{|M| \leq |K|+1}  \frac{1}{(1+t+|q|)} \cdot | \derm ( \Lie_{Z^M}  \Phi )  |  +  \sum_{ V^\prime \in \cal T }   \sum_{|M| \leq |K|+1}  \frac{1}{(1+|q|)} \cdot  | \derm   ( \Lie_{Z^M}  \Phi _{V^\prime} ) | \; .\\
 \eea
 Therefore, for all $V \in \cal T$, 
 \bea
 \notag
 && \sum_{  |J| + |K| \leq |I|, \; |K| < |I| }  |   ( \Lie_{Z^{J}}   H)_{L  L}    \cdot  \derm_{\underline{L}}   \derm_{ \underline{L}}   (  \Lie_{Z^{K}}  \Phi_{V} ) | \\
  \notag
  &\les&  \sum_{  |J| + |K| \leq |I|, \; |K| < |I| }    |   ( \Lie_{Z^{J}}   H)_{L  L} |   \cdot    \sum_{|M| \leq |K|+1}  \frac{1}{(1+t+|q|)} \cdot | \derm ( \Lie_{Z^M}  \Phi )  |  \\
   \notag
 && + \sum_{  |J| + |K| \leq |I|, \; |K| < |I| }   |   ( \Lie_{Z^{J}}   H)_{L  L}   | \cdot   \sum_{|M| \leq |K|+1}  \frac{1}{(1+|q|)} \cdot   \big( \sum_{ V^\prime \in \cal T } | \derm   ( \Lie_{Z^M}  \Phi _{V^\prime} ) | \big)  \;. \\
 \eea
 
 Since, as we have already shown in Lemma \ref{decayrateforfullderivativeintermofZ}, 
  \beaa
 |\derm \Psi_{UV} | \les \sum_{|I| \leq 1}  \frac{1}{(1+|q|)} \cdot  | \Lie_{Z^I} \Psi | \; ,
 \eeaa
 we also get that for all $U \in \cal U$,
  \bea
   \notag
 && \sum_{  |J| + |K| \leq |I|, \; |K| < |I| }  |   ( \Lie_{Z^{J}}   H)_{L  L}    \cdot  \derm_{\underline{L}}   \derm_{ \underline{L}}   (  \Lie_{Z^{K}}  \Phi_{U} ) | \\
  \notag
  &\les&  \sum_{  |J| + |K| \leq |I|, \; |K| < |I| }   |   ( \Lie_{Z^{J}}   H)_{L  L}   | \cdot   \sum_{|M| \leq |K|+1}  \frac{1}{(1+|q|)} \cdot    | \derm   ( \Lie_{Z^M}  \Phi  ) |  \;. \\
 \eea
 Obviously, the above estimate is also true for a two-tensor.
   
 \textbf{The term $ \derm_{\underline{L}}   \derm_{L}   (  \Lie_{Z^{K}}  \Phi_{U} )$}:
 
We have
\beaa
 && \derm_{\underline{L}}   \derm_{L}   (  \Lie_{Z^{K}}  \Phi_{U} ) \\
  &=&   \pa_{\underline{L}}  \Big(  \derm_{L}   (  \Lie_{Z^{K}}  \Phi_{U} ) \Big) -    \derm_{ \derm_{\underline{L}} L }   (  \Lie_{Z^{K}}  \Phi_{U} ) -    \derm_{L }   (  \Lie_{Z^{K}}  \Phi_{ \derm_{\underline{L}}  U} )  \; .
  \eeaa

We have shown in the Appendix of \cite{G2} (by taking in \cite{G2}, a zero Schwarzschild mass, which gives a Minkowski metric), that for all $U \in {\cal U}$,
  \beaa
 \derm_{\underline{L}} U  =  \derm_{L } U = 0 \; .
\eeaa
Thus,
\beaa
  \derm_{\underline{L}}   \derm_{L}   (  \Lie_{Z^{K}}  \Phi_{U} ) &=&   \pa_{\underline{L}}  \Big(  \pa_{L}   (  \Lie_{Z^{K}}  \Phi_{U} ) \Big)   \; .
  \eeaa
As a result, we could write,
\bea
\notag
  \derm_{\underline{L}}   \derm_{L}   (  \Lie_{Z^{K}}  \Phi_{U} ) &=&    \pa_{L}  \Big(  \pa_{\underline{L}}  (  \Lie_{Z^{K}}  \Phi_{U} ) \Big)   \\
  \notag
  &=&    \pa_{L}  \Big(  \derm_{\underline{L}}  (  \Lie_{Z^{K}}  \Phi_{U} )  \Big)   \\
  &=&    \derm_{L}    \derm_{\underline{L}}  (  \Lie_{Z^{K}}  \Phi_{U} )     \; .
  \eea
Hence,
 \bea
 \notag
 && \sum_{  |J| + |K| \leq |I|, \; |K| < |I| }  | ( \Lie_{Z^{J}}   H)_{L  \underline{L}}    \cdot  \derm_{\underline{L}}   \derm_{L}   (  \Lie_{Z^{K}}  \Phi_{U} )  |  \\
  \notag
 &\les&  \sum_{  |J| + |K| \leq |I|, \; |K| < |I| }  | ( \Lie_{Z^{J}}   H)_{L  \underline{L}} |   \cdot | \derm_{L}    \derm_{\underline{L}}  (  \Lie_{Z^{K}}  \Phi_{U} )  | \; .\\
 \eea
 
\textbf{The term $\derm_{\underline{L}}   \derm_{e_A}   (  \Lie_{Z^{K}}  \Phi_{U} )$}:

We look at the region $t \geq 1 $ or $r \geq 1$ and we distinguish the cases $q > 0$ and $q \leq 0$:
  
\textbf{In the case where $q = r-t > 0$}:

We showed in Lemma \ref{restrictedderivativesintermsofZ}, that
 \beaa
\rpa_i  &=&   \frac{x^j}{r^2}   Z_{ij} \, ,
 \eeaa
 and hence, we can write
 \bea
e_A = \frac{1}{r}C^{ij}_A Z_{ij} \, ,
\eea
where $C^{ij}_A$ are in fact linear combinations of  $\frac{x^j}{r}$. However, 

\beaa
\pa_{\underline{L}} \frac{x^j}{r} &=&(  \pa_{t} - \pa_{r} ) \frac{x^j}{r}   =  - \pa_{r} ( \frac{x^j}{r} )  =   -  \frac{x^{i}}{r} \pa_{i} ( \frac{x^j}{r}  )   \\
&=& -  \frac{x^{i}}{r} ( \frac{\de_{i}^{\;\; j} }{r}  ) +   \frac{x^{i}}{r}   \frac{x^j}{r^2} ( \pa_{i} r  ) \; .
\eeaa
We have
\beaa
 \pa_{i} r &=&  \pa_{i} (  \sqrt{(x^1)^2 + (x^2)^2 +(x^3)^2 } ) = \frac{1}{2} \frac{2 x_i}{  \sqrt{(x^1)^2 + (x^2)^2 +(x^3)^2 } } \\
 &=& \frac{x_i}{r} \; .
\eeaa
Thus,
\beaa
 \pa_{r} ( \frac{x^j}{r} ) &=&  \frac{x^{i}}{r} ( \frac{\de_{i}^{\;\; j} }{r}  ) -    \frac{x^{i}}{r}   \frac{x^j}{r^2} ( \pa_{i} r  )  =  \frac{x^{i}}{r} ( \frac{\de_{i}^{\;\; j} }{r}  ) -  \frac{x^{i} \cdot x_i}{r^2}   \frac{x^j}{r^2}  =   \frac{x^{j}}{r^2}  -  \frac{r^2}{r^2} \cdot   \frac{x^j}{r^2}  \\
  &=&  0 \;.
\eeaa
As a result
\bea\label{Lbarderivativeofxioverrisinfactzero}
\pa_{\underline{L}} \frac{x^j}{r} = 0 \; ,
\eea
and therefore, 
 \beaa
\pa_{\underline{L}} C^{ij}_A = 0 \, .
\eeaa
We also have 
 \beaa
\pa_{\underline{L}} \frac{1}{r} &=& - \pa_{r}  ( \frac{1}{r}  ) =   \frac{1}{r^2}   \, .
\eeaa
Thus,
\bea
\pa_{\underline{L}} ( \frac{1}{r}  C^{ij}_A ) =  \frac{1}{r^2}   C^{ij}_A \; .
\eea
Since
\beaa 
  \derm_{e_A}   (  \Lie_{Z^{K}}  \Phi_{U} ) &=&  \frac{1}{r} C^{ij}_A \derm_{Z_{ij} }  (  \Lie_{Z^{K}}  \Phi_{U} ) \; ,
\eeaa
therefore,
  \bea
  \notag
  \derm_{\underline{L}}   \derm_{e_A}   (  \Lie_{Z^{K}}  \Phi_{U} ) &=&  \pa_{\underline{L}} ( \frac{1}{r}  C^{ij}_A ) \cdot \derm_{Z_{ij} }  (  \Lie_{Z^{K}}  \Phi_{U} )  + \frac{1}{r} C^{ij}_A   \derm_{\underline{L}}  \derm_{Z_{ij} }  (  \Lie_{Z^{K}}  \Phi_{U} ) \\
    \notag
   &=&   \frac{1}{r^2}  C^{ij}_A \cdot \derm_{Z_{ij} }  (  \Lie_{Z^{K}}  \Phi_{U} )  + \frac{1}{r} C^{ij}_A   \derm_{\underline{L}}  \derm_{Z_{ij} }  (  \Lie_{Z^{K}}  \Phi_{U} ) \; . \\
  \eea

However, $ Z_{ij} = x_{j} \pa_{i} - x_{i} \pa_{j} $, and thus,
\beaa
  \derm_{Z_{ij} }  (  \Lie_{Z^{K}}  \Phi ) = x_{j} \derm_{ \frac{\pa}{\pa x_i}}  (  \Lie_{Z^{K}}  \Phi )   - x_{i} \derm_{ \frac{\pa}{\pa x_j}}    (  \Lie_{Z^{K}}  \Phi )   \, .
\eeaa
Therefore,
\beaa
 | \frac{1}{r^2}  C^{ij}_A \cdot \derm_{Z_{ij} }  (  \Lie_{Z^{K}}  \Phi_{U} )  | &\leq&   | \frac{x_{j}}{r^2}  C^{ij}_A  |  \cdot | \derm_{ \frac{\pa}{\pa x_i}}  (  \Lie_{Z^{K}}  \Phi )| +    | \frac{x_{i}}{r^2}  C^{ij}_A  | \cdot | \derm_{ \frac{\pa}{\pa x_j}}    (  \Lie_{Z^{K}}  \Phi )  | \; .
 \eeaa
 Given that $C^{ij}_A$ and $\frac{x_{i}}{r}$ are bounded, we get 
 \bea
 \notag
 | \frac{1}{r^2}  C^{ij}_A \cdot \derm_{Z_{ij} }  (  \Lie_{Z^{K}}  \Phi_{U} )  | &\leq&    \frac{1}{r}   \cdot | \derm_{ \frac{\pa}{\pa x_i}}  (  \Lie_{Z^{K}}  \Phi )| +   \frac{1}{r} \cdot | \derm_{ \frac{\pa}{\pa x_j}}    (  \Lie_{Z^{K}}  \Phi )  | \\
 &\leq&    \frac{1}{r}   \cdot | \derm (  \Lie_{Z^{K}}  \Phi )|  \; .
 \eea
 
 Whereas to the term $\frac{1}{r} C^{ij}_A   \derm_{\underline{L}}  \derm_{Z_{ij} }  (  \Lie_{Z^{K}}  \Phi_{U} ) $, we first note that since it is a tensor, we have
 
 \beaa
 | \frac{1}{r} C^{ij}_A   \derm_{\underline{L}}  \derm_{Z_{ij} }  (  \Lie_{Z^{K}}  \Phi_{U} ) | \les  | \frac{1}{r} C^{ij}_A   \derm_{\underline{L}}  \derm_{Z_{ij} }  (  \Lie_{Z^{K}}  \Phi ) | \; .
 \eeaa
 
 Thus, we can compute in wave coordinates the right side of the above inequality in order to make an estimate. First, we compute for $\mu \in \{0, 1, 2, 3 \}$, 
 \beaa
  &&  \derm_{Z_{ij} }  (  \Lie_{Z^{K}}  \Phi_{\mu} ) \\
      &=&   x_{j} \derm_{ \frac{\pa}{\pa x_i}}  (  \Lie_{Z^{K}}  \Phi_{\mu} )   - x_{i} \derm_{ \frac{\pa}{\pa x_j}}    (  \Lie_{Z^{K}}  \Phi_{\mu} )   \\
  &=&   x_{j} \pa_{i}  (  \Lie_{Z^{K}}  \Phi_{\mu} )   - x_{i} \pa_{j}    (  \Lie_{Z^{I_2}}  \Phi_{\mu} )   \\ 
    &=&  \pa_{Z_{ij}}  (  \Lie_{Z^{K}}  \Phi_{\mu} )   \\ 
    &=&  \Lie_{Z_{ij}}  (   \Lie_{Z^{K}}  \Phi_{\mu} )  +   \Lie_{Z^{K}}  \Phi ( [Z_{ij},  \frac{\pa}{\pa x_{\mu}} ] ) \\
     &=&  \Lie_{Z_{ij}}  \Lie_{Z^{K}}  \Phi_{\mu}  + \sum_{m} \Lie_{Z^{J_{m_1}} } \Lie_{ [Z_{ij}, Z^{m}] }  \Lie_{Z^{J_{m_2}} } \Phi_{\mu}   +  \Lie_{Z^{K}}  \Phi ( [Z_{ij},  \frac{\pa}{\pa x_{\mu}} ] ) \; ,
  \eeaa
  where by $Z^{m}$, we mean the $m$-th Minkowski vector field in the product $Z^K$, and the rest in the product is $Z^{J_{m_1}}$ (what is before $Z^m$) and $Z^{J_{m_2}}$ (what is after $Z^m$).

  However,  the commutation of two vector fields in $\cal Z$ is a linear combination of vector fields in $\cal Z$, and the commutation of a vector field in $\cal Z$ and of a vector $\pa_\mu$, $\mu \in  \{t, x^1, x^2, x^3 \}$, gives a linear combination of vectors of the form $\pa_\mu$. Using that fact, we get that for all $ \mu \in (t, x^1, x^2, x^3)$,
  \beaa
\derm_{Z_{ij} }  (  \Lie_{Z^{K}}  \Phi_{\mu} ) &=&  \sum_{i} a_i \Lie_{Z^{I_i}} \Phi_{\mu} + \sum_{i} b_i \cdot \Lie_{Z^{K}} \Phi_{\mu^i}  \; ,
\eeaa
for some $\mu^i \in  \{t, x^1, x^2, x^3 \}$, and $a_i, b_i$ are constants, and $|I_i| \leq |K| + 1$.

Thus,
\beaa
 \derm_{\underline{L}}   \derm_{Z_{ij} }  (  \Lie_{Z^{K}}  \Phi_{\mu} ) &=&  \sum_{i} a_i    \derm_{\underline{L}}   ( \Lie_{Z^{I_i}} \Phi_{\mu})  + \sum_{i} b_i \cdot   \derm_{\underline{L}}   ( \Lie_{Z^{K}} \Phi_{\mu^i} ) 
\eeaa
and therefore, for all $\mu \in  \{t, x^1, x^2, x^3 \}$
 \beaa
 |   \derm_{\underline{L}}  \derm_{Z_{ij} }  (  \Lie_{Z^{K}}  \Phi_{\mu} ) | \les \sum_{|M| \leq |K|+1} | \derm_{\underline{L}}   ( \Lie_{Z^{M}} \Phi )  | \; ,
 \eeaa
and as a result
 \bea
 \notag
 | \frac{1}{r} C^{ij}_A   \derm_{\underline{L}}  \derm_{Z_{ij} }  (  \Lie_{Z^{K}}  \Phi_{U} ) |   \les \sum_{|M| \leq |K|+1}  \frac{1}{r} \cdot  | \derm  ( \Lie_{Z^{M}} \Phi )  | \; .
 \eea

 Finally, we obtain for $q > 0$,
   \bea
  \notag
  | \derm_{\underline{L}}   \derm_{e_A}   (  \Lie_{Z^{K}}  \Phi_{U} ) | &\les&  | \frac{1}{r^2}  C^{ij}_A \cdot \derm_{Z_{ij} }  (  \Lie_{Z^{K}}  \Phi_{U} )  |  +  | \frac{1}{r} C^{ij}_A   \derm_{\underline{L}}  \derm_{Z_{ij} }  (  \Lie_{Z^{K}}  \Phi_{U} ) |  \\
  &\les&    \sum_{|M| \leq |K|+1} \frac{1}{r}   \cdot | \derm (  \Lie_{Z^{M}}  \Phi )| \; .
  \eea

Since we look in the region $t \geq 1 $ or $r \geq 1$ and we are in the case where $q  > 0$, this imposes that $r \geq 1 $ and therefore, in that region, we have 
\bea
\frac{1}{r} \leq \frac{1}{1+r} \; ,
\eea
and as we have shown in Lemma \ref{decayfactorsforregions}, for $q \geq 0$, we have
\beaa
 \frac{1}{1+r} &= & \frac{1}{1 + t + |q|  }  \;.
\eeaa
Consequently,
   \bea
  \notag
  | \derm_{\underline{L}}   \derm_{e_A}   (  \Lie_{Z^{K}}  \Phi_{U} ) | &\les&    \sum_{|M| \leq |K|+1}  \frac{1}{1+r}   \cdot | \derm (  \Lie_{Z^{M}}  \Phi )|  \\
  \notag
  &\les&  \frac{1}{(1 + t + |q|)  }  \cdot   \sum_{|M| \leq |K|+1} | \derm (  \Lie_{Z^{M}}  \Phi )| \; . \\
  \eea
  
  Thus, in the region $t \geq 1 $ or $r \geq 1$, for $q  > 0$, we have
     \bea
  \notag
  && \sum_{  |J| + |K| \leq |I|, \; |K| < |I| }   |  ( \Lie_{Z^{J}}   H)_{L e_A }    \cdot  \derm_{\underline{L}}   \derm_{e_A}    (  \Lie_{Z^{K}}  \Phi_{U} )  | \\
    \notag
    &\les&    \sum_{  |J| + |K| \leq |I|, \; |K| < |I| }    \frac{1}{(1 + t + |q|)  }  \cdot   |  ( \Lie_{Z^{J}}   H)_{L e_A }   | \cdot | \sum_{|M| \leq |K|+1} | \derm (  \Lie_{Z^{M}}  \Phi )| \; .
  \eea
  
\textbf{ In the case where $q = r-t < 0$}:
 
 We then recall that we showed in Lemma \ref{restrictedderivativesintermsofZ}, that we have also a different presentation for $\rpa_{i}$\,, which is more suitable to make estimates in the region $q  < 0$, that is
 \beaa
\rpa_{i} &=& \frac{  - \frac{x_i}{r} \frac{x^j}{r}   Z_{0j}  +  Z_{0i}    }{  t  }  \; , \\
 \eeaa
 and therefore, we can see $e_A$ as linear combinations of the form 
 \bea
e_A = \frac{1}{t} G_{A}^{j} Z_{0j} +  \frac{1}{t} B_{A}^{i} Z_{0i} \, ,
\eea
where $G_{A}^{j}$ is a linear combination of $\frac{x_i}{r} \frac{x^j}{r}$ and $B_{A}^{i}$ are constants. Based on what we have shown in \eqref{Lbarderivativeofxioverrisinfactzero}, we have then
 \beaa
\pa_{\underline{L}} G_{A}^{j}  = \pa_{\underline{L}} B_{A}^{j}  =  0 \, .
\eeaa
Thus, we can write
 \bea
e_A = \frac{1}{t} C_{A}^{j} Z_{0j} \, ,
\eea
with $C_{A}^{j} $ bounded and furthermore,
 \beaa
\pa_{\underline{L}} C_{A}^{j}  = 0 \, .
\eeaa
We also have 
 \beaa
\pa_{\underline{L}} \frac{1}{t} &=&  \pa_{t}  ( \frac{1}{t}  ) =   - \frac{1}{t^2}   \, .
\eeaa
Thus,
\beaa
\pa_{\underline{L}} ( \frac{1}{t}   C_{A}^{j}   ) &=&  - \frac{1}{t^2}   C^{j}_A \; .\\
\eeaa

Since
\beaa 
  \derm_{e_A}   (  \Lie_{Z^{K}}  \Phi_{U} ) &=& \frac{1}{t} C_{A}^{j} \derm_{Z_{0j} }  (  \Lie_{Z^{K}}  \Phi_{U} )   \; ,
\eeaa
therefore,
  \bea\label{Lbareacovarderiderivativeofliederivativephi}
  \notag
  \derm_{\underline{L}}   \derm_{e_A}   (  \Lie_{Z^{K}}  \Phi_{U} ) &=&  \pa_{\underline{L}} ( \frac{1}{t} C_{A}^{j}  ) \cdot \derm_{Z_{0j}}  (  \Lie_{Z^{K}}  \Phi_{U} )  + \frac{1}{t} C_{A}^{j}   \derm_{\underline{L}}  \derm_{Z_{0j} }  (  \Lie_{Z^{K}}  \Phi_{U} ) \\
    \notag
   &=& -  \frac{1}{t^2} C_{A}^{j}  \cdot \derm_{Z_{0j} }  (  \Lie_{Z^{K}}  \Phi_{U} )  + \frac{1}{t} C_{A}^{j}    \derm_{\underline{L}}  \derm_{Z_{0j} }  (  \Lie_{Z^{K}}  \Phi_{U} ) \; . \\
  \eea

However, $Z_{0j} =   x_{j} \pa_{t} + t \pa_{j}  $, and thus,
\beaa
  \derm_{Z_{0j} }  (  \Lie_{Z^{K}}  \Phi ) = x_{j} \derm_{ \frac{\pa}{\pa t}}  (  \Lie_{Z^{K}}  \Phi )   + t  \derm_{ \frac{\pa}{\pa x_j}}    (  \Lie_{Z^{K}}  \Phi )   \, .
\eeaa
Therefore, on one hand
\beaa
 | \frac{1}{t^2}  C^{j}_A \cdot \derm_{Z_{0j} }  (  \Lie_{Z^{K}}  \Phi_{U} )  | &\leq&   | \frac{x_{j}}{t^2}  C^{j}_A  |  \cdot | \derm_{ \frac{\pa}{\pa t}}  (  \Lie_{Z^{K}}  \Phi )| +    | \frac{t}{t^2}  C^{j}_A  | \cdot | \derm_{ \frac{\pa}{\pa x_j}}    (  \Lie_{Z^{K}}  \Phi )  | \; .
 \eeaa
 Now, we recall that we are in the region $q = r-t < 0$ and therefore in that region, 
 \beaa
 \frac{r}{t} \leq 1 \; ,
 \eeaa
 and therefore, since also $C^{j}_A$ is bounded, we get 
 \bea
 \notag
 | \frac{1}{t^2}  C^{j}_A \cdot \derm_{Z_{0j} }  (  \Lie_{Z^{K}}  \Phi_{U} )  |  &\leq&    \frac{1}{t}   \cdot | \derm_{ \frac{\pa}{\pa t}}  (  \Lie_{Z^{K}}  \Phi )| +   \frac{1}{t} \cdot | \derm_{ \frac{\pa}{\pa x_j}}    (  \Lie_{Z^{K}}  \Phi )  | \\
 &\leq&    \frac{1}{t}   \cdot | \derm (  \Lie_{Z^{K}}  \Phi )|  \; .
 \eea
 
On the other hand, for the other term in \eqref{Lbareacovarderiderivativeofliederivativephi}, proceeding as earlier, we obtain
 
 \beaa
 | \frac{1}{t} C_{A}^{j}    \derm_{\underline{L}}  \derm_{Z_{0j} }  (  \Lie_{Z^{K}}  \Phi_{U} ) |&\les&   | \frac{1}{t} C_{A}^{j}    \derm_{\underline{L}}  \derm_{Z_{0j} }  (  \Lie_{Z^{K}}  \Phi) | \\
&\les& \sum_{|M| \leq |K|+1}  \frac{1}{t} \cdot  | \derm  ( \Lie_{Z^{M}} \Phi )  | \; .
 \eeaa

 Finally, we obtain for $q < 0$,
   \bea
  \notag
  |   \derm_{\underline{L}}   \derm_{e_A}   (  \Lie_{Z^{K}}  \Phi_{U} ) |  &\les&    \sum_{|M| \leq |K|+1} \frac{1}{t}   \cdot | \derm (  \Lie_{Z^{M}}  \Phi )| \; .
  \eea

Since we look in the region $t \geq 1 $ or $r \geq 1$ and we are in the case where $q  < 0$, this imposes, as we have already shown in Lemma \ref{decayfactorsforregions}, that
\bea
 \frac{1}{t} &\les&  \frac{1}{1+t} \les  \frac{1}{1 + t + |q|  } \; ,
  \eea
  and hence
   \bea
  \notag
  | \derm_{\underline{L}}   \derm_{e_A}   (  \Lie_{Z^{K}}  \Phi_{U} ) |  &\les&  \frac{1}{(1 + t + |q|)  }  \cdot   \sum_{|M| \leq |K|+1} | \derm (  \Lie_{Z^{M}}  \Phi )| \; . \\
  \eea
  
  Thus, in the region $t \geq 1 $ or $r \geq 1$, for $q  < 0$, we have
     \bea
  \notag
  && \sum_{  |J| + |K| \leq |I|, \; |K| < |I| }   |  ( \Lie_{Z^{J}}   H)_{L e_A }    \cdot  \derm_{\underline{L}}   \derm_{e_A}    (  \Lie_{Z^{K}}  \Phi_{U} )  | \\
    \notag
    &\les&    \sum_{  |J| + |K| \leq |I|, \; |K| < |I| }    \frac{1}{(1 + t + |q|)  }  \cdot   |  ( \Lie_{Z^{J}}   H)_{L e_A }   | \cdot | \sum_{|M| \leq |K|+1} | \derm (  \Lie_{Z^{M}}  \Phi )| \; .
  \eea

Finally, putting all together we obtain the result. We note the the calculations are also true for a two-tensor.
 \end{proof}

\begin{lemma}\label{Thecommutationformulawithpossibilityofsperationoftangentialcomponentsaswell}

Let  $\Phi_{\mu}$ be a tensor valued either in the Lie algebra or a scalar, satisfying the following tensorial wave equation
\beaa
 g^{\la\a} \derm_{\la}   \derm_{\a}   \Phi_{\mu}= S_{\mu} \, , 
\eeaa
where $S_{\mu}$ is the source term. Then, we have for any $V \in \cal T$,
 \beaa
\notag
&&| \Lie_{Z^I}  ( g^{\la\mu} \derm_{\la}   \derm_{\mu}     \Phi_{V} ) - g^{\la\mu}    \derm_{\la}   \derm_{\mu}  (  \Lie_{Z^I} \Phi_{V}  ) |  \\
  \notag
   &\les&  \sum_{|K| < |I| }  | g^{\la\mu} \cdot \derm_{\la}   \derm_{\mu} (  \Lie_{Z^{K}}  \Phi_{V} ) | \\
   \notag
&&+  \frac{1}{(1+t+|q|)}  \cdot \sum_{|K|\leq |I|,}\,\, \sum_{|J|+(|K|-1)_+\le |I|} \,\,\, | \Lie_{Z^{J}} H |\, \cdot | \derm ( \Lie_{Z^K}  \Phi )  | \\
   \notag
&& +   \frac{1}{(1+|q|)}  \cdot \sum_{|K|\leq |I|,}\,\, \sum_{|J|+(|K|-1)_+\le |I|} \,\,\, | \Lie_{Z^{J}} H_{L  L} |\, \cdot  \big( \sum_{ V^\prime \in \cal T } | \derm ( \Lie_{Z^K}  \Phi _{V^\prime} )  | \:  \big) \; ,
\eeaa
  where $(|K|-1)_+=|K|-1$ if $|K|\geq 1$ and $(|K|-1)_+=0$ if $|K|=0$.\\
  Furthermore, let $\Phi_{\mu\nu }$ be a two-tensor (or a one-tensor $\Phi_{\mu}$) valued either in the Lie algebra or a scalar. Then, we have for any $U,  V \in \cal U$,
 \beaa
\notag
&&| \Lie_{Z^I}  ( g^{\la\mu} \derm_{\la}   \derm_{\mu}     \Phi ) - g^{\la\mu}    \derm_{\la}   \derm_{\mu}  (  \Lie_{Z^I} \Phi  ) |  \\
  \notag
   &\les&  \sum_{|K| < |I| }  | g^{\la\mu} \cdot \derm_{\la}   \derm_{\mu} (  \Lie_{Z^{K}}  \Phi_{V} ) | \\
   \notag
&&+  \frac{1}{(1+t+|q|)}  \cdot \sum_{|K|\leq |I|,}\,\, \sum_{|J|+(|K|-1)_+\le |I|} \,\,\, | \Lie_{Z^{J}} H |\, \cdot | \derm ( \Lie_{Z^K}  \Phi )  | \\
   \notag
&& +   \frac{1}{(1+|q|)}  \cdot \sum_{|K|\leq |I|,}\,\, \sum_{|J|+(|K|-1)_+\le |I|} \,\,\, | \Lie_{Z^{J}} H_{L  L} |\, \cdot   | \derm ( \Lie_{Z^K}  \Phi)  |  \; .
\eeaa
 \end{lemma}
 
 \begin{proof}
 We estimate
 \beaa
  &&    | ( \Lie_{Z^{J}}   H)_{L  \underline{L}} |   \cdot | \derm_{L}    \derm_{\underline{L}}  (  \Lie_{Z^{K}}  \Phi_{U} )  | \\
\notag
 && +       | m^{\mu \b}   ( \Lie_{Z^{J}}   H)_{\underline{L}\b}  \cdot    \derm_{L}   \derm_{\mu}     (  \Lie_{Z^{K}}  \Phi_{U} )  | + |   m^{\mu \b}   ( \Lie_{Z^{J}}   H)_{e_A \b}   \cdot   \derm_{e_A}   \derm_{\mu}    (  \Lie_{Z^{K}}  \Phi_{U} )   | \Big) \\
 \notag
 &\les &   |  \Lie_{Z^{J}}   H |   \cdot | \derm_{L}    \derm  (  \Lie_{Z^{K}}  \Phi )  |  + |  \Lie_{Z^{J}}   H |  \cdot  |  \derm_{L}   \derm     (  \Lie_{Z^{K}}  \Phi )  | \\
 \notag
 && + | \Lie_{Z^{J}}   H |     \cdot  | \derm_{e_A}   \derm    (  \Lie_{Z^{K}}  \Phi )   | \Big) \\
 &\les & | \Lie_{Z^{J}}   H |     \cdot  | \rderm   \derm    (  \Lie_{Z^{K}}  \Phi )   |
  \eeaa
  In wave coordinates, for fixed $\mu, \nu \in \{t, x^1, x^2, x^3 \}$, we can look at each component and estimate, as we have shown in Lemma \ref{betterdecayfortangentialderivatives}, that
  \beaa
   | \rderm   \derm_{\mu}    (  \Lie_{Z^{K}}  \Phi_{\nu} )   | &\leq&    | \rpa   \derm_{\mu}    (  \Lie_{Z^{K}}  \Phi_{\nu} )   | \\
   &\leq& \frac{1}{(1 + t + |q| )}  \cdot \sum_{|I| = 1} |Z^I  \derm_{\mu}    (  \Lie_{Z^{K}}  \Phi_{\nu} )   |  .
\eeaa
Using the fact that commutation of a Minkowski vector field $Z$ and a wave coordinate vector field $\frac{\pa}{\pa x_{\mu}}$, gives a linear combination of wave coordinates vector fields, we get
\beaa
 \sum_{|I| = 1}  |Z^I  \derm_{\mu}    (  \Lie_{Z^{K}}  \Phi_{\nu} )   | \les \sum_{|I| \leq 1}  |\Lie_{Z^I}   \derm (  \Lie_{Z^{K}}  \Phi )   |  \; .
\eeaa
Summing over all indices $\mu, \nu \in \{t, x^1, x^2, x^3 \}$, we obtain
\beaa
  | \rderm   \derm    (  \Lie_{Z^{K}}  \Phi )   |  \les \frac{1}{(1 + t + |q| )}  \cdot  \sum_{|I| \leq 1}  |\Lie_{Z^I}   \derm (  \Lie_{Z^{K}}  \Phi )   |  \;.
\eeaa
As a result, using the fact that the Lie derivative $\Lie_{Z^{I}}$ commutes with $\derm$, we get
 \beaa
  &&    | ( \Lie_{Z^{J}}   H)_{L  \underline{L}} |   \cdot | \derm_{L}    \derm_{\underline{L}}  (  \Lie_{Z^{K}}  \Phi_{U} )  | \\
\notag
 && +       | m^{\mu \b}   ( \Lie_{Z^{J}}   H)_{\underline{L}\b}  \cdot    \derm_{L}   \derm_{\mu}     (  \Lie_{Z^{K}}  \Phi_{U} )  | + |   m^{\mu \b}   ( \Lie_{Z^{J}}   H)_{e_A \b}   \cdot   \derm_{e_A}   \derm_{\mu}    (  \Lie_{Z^{K}}  \Phi_{U} )   | \Big) \\
 \notag
 &\les & | \Lie_{Z^{J}}   H |     \cdot  \frac{1}{(1 + t + |q| )}  \cdot   \sum_{|M| \leq |K| +1}  |   \derm (  \Lie_{Z^{M}}  \Phi )   | \; .
 \eeaa

Thus, for any $V \in \cal T$,
 \beaa
\notag
&&| \Lie_{Z^I}  ( g^{\la\mu} \derm_{\la}   \derm_{\mu}     \Phi_{V} ) - g^{\la\mu}    \derm_{\la}   \derm_{\mu}  (  \Lie_{Z^I} \Phi_{V}  ) |  \\
\notag
&\les&  \sum_{|K| < |I| }  | g^{\la\mu} \cdot \derm_{\la}   \derm_{\mu} (  \Lie_{Z^{K}}  \Phi_{V} ) | \\
\notag
&& +  \sum_{|J| + |K| \leq |I|, \; |K| < |I| }  \Big(   |   ( \Lie_{Z^{J}}   H)_{L  L} |   \cdot    \frac{1}{(1+t+|q|)} \cdot  \sum_{|M| \leq |K|+1}  | \derm ( \Lie_{Z^M}  \Phi )  |  \\
   \notag
 && +    |   ( \Lie_{Z^{J}}   H)_{L  L}   | \cdot   \frac{1}{(1+|q|)} \cdot   \sum_{|M| \leq |K|+1}    \sum_{ V^\prime \in \cal T } | \derm   ( \Lie_{Z^M}  \Phi _{V^\prime} ) |   \\
\notag
    &&  +    \frac{1}{(1 + t + |q|)  }  \cdot   |  ( \Lie_{Z^{J}}   H)_{L e_A }   | \cdot | \sum_{|M| \leq |K|+1} | \derm (  \Lie_{Z^{M}}  \Phi )| \\
 && +  | \Lie_{Z^{J}}   H |     \cdot  \frac{1}{(1 + t + |q| )}  \cdot \sum_{|M| \leq |K| +1}  |   \derm (  \Lie_{Z^{M}}  \Phi )   |  \\
 &\les&  \sum_{|K| < |I| }  | g^{\la\mu} \cdot \derm_{\la}   \derm_{\mu} (  \Lie_{Z^{K}}  \Phi_{V} ) | \\
\notag
&& +  \sum_{|J| + |K| \leq |I|, \; |K| < |I| }  \Big(   |   \Lie_{Z^{J}}   H  |   \cdot    \frac{1}{(1+t+|q|)} \cdot  \sum_{|M| \leq |K|+1}  | \derm ( \Lie_{Z^M}  \Phi )  |  \\
   \notag
 && +    |   ( \Lie_{Z^{J}}   H)_{L  L}   | \cdot   \frac{1}{(1+|q|)} \cdot   \sum_{|M| \leq |K|+1}    \sum_{ V^\prime \in \cal T } | \derm   ( \Lie_{Z^M}  \Phi _{V^\prime} ) | \Big)  \\
  \eeaa

Consequently, for all $V \in \cal T$, 
 \beaa
\notag
&&| \Lie_{Z^I}  ( g^{\la\mu} \derm_{\la}   \derm_{\mu}     \Phi_{V} ) - g^{\la\mu}    \derm_{\la}   \derm_{\mu}  (  \Lie_{Z^I} \Phi_{V}  ) |  \\
  \notag
   &\les&  \sum_{|K| < |I| }  | g^{\la\mu} \cdot \derm_{\la}   \derm_{\mu} (  \Lie_{Z^{K}}  \Phi_{V} ) | \\
   \notag
&&+  \frac{1}{(1+t+|q|)}  \cdot \sum_{|K|\leq |I|,}\,\, \sum_{|J|+(|K|-1)_+\le |I|} \,\,\, | \Lie_{Z^{J}} H |\, \cdot | \derm ( \Lie_{Z^K}  \Phi )  | \\
   \notag
&& +   \frac{1}{(1+|q|)}  \cdot \sum_{|K|\leq |I|,}\,\, \sum_{|J|+(|K|-1)_+\le |I|} \,\,\, |( \Lie_{Z^{J}} H)_{L  L} |\, \cdot  \big( \sum_{ V^\prime \in \cal T } | \derm ( \Lie_{Z^K}  \Phi _{V^\prime} )  | \: , \big) \; .
\eeaa

Similarly, we get the result also for a 1-tensor, or 2-tensor $\Phi_{UV}$, and where $U, V \in \cal U$, yet the sum then on the right hand side should run for $V^\prime \in \cal U $.
\end{proof}

 \begin{lemma}\label{estimatethatallowsupgradeincorporatingtermsfromthecommutationformula}
 
We have for $\gamma^\prime$ such that $-1 \leq \gamma^\prime < \gamma - \delta$, and for $ \delta <    1/2 $, and for all $U, V\in  \{L,\Lb,A,B\}$,

\bea
   \notag
&& (1+t+|q|) \cdot |\varpi(q) \cdot \derm ( \Lie_{Z^J} A)_{V} (t,x)| \\
   \notag
 &\les&   c (\gamma^\prime) \cdot  c (\delta) \cdot c (\gamma) \cdot C ( |J|  ) \cdot E ( |J| + 4) \cdot \eps \, \\
\notag
&& + c (\gamma^\prime) \cdot c (\gamma)  \cdot c (\delta)  \cdot E ( 3)  \cdot \eps  \cdot  \int_0^t  \frac{1}{(1+\tau)} \cdot (1+\tau+|q|) \cdot \|\varpi(q)  \cdot  \derm  ( \Lie_{Z^J} A_{V}) (\tau,\cdot) \|_{L^\infty (\Sigma^{ext}_{\tau} )} d \tau \\
       \notag
    && + \sum_{|K| \leq |J|} \int_0^t (1+\tau) \cdot  \varpi(q) \cdot  \|    \Lie_{Z^K}  g^{\la\mu} \derm_{\la}   \derm_{\mu}  A_V (\tau,\cdot) \|_{L^\infty(\overline{D}_\tau)} d\tau \\
  \notag
& &+ \int_0^t   \frac {(1+\tau) \cdot  \varpi(q)}{(1+\tau+|q|)} \cdot    \sum_{|K|\leq |J|} \Big( \sum_{|J^{\prime}|+(|K|-1)_+\le |J|} \,\,\,
|\Lie_{Z^{J^{\prime}}} H|\cdot {|\derm ( \Lie_{Z^{K} } A) |}  \Big) d\tau \\
\notag
& &+  \int_0^t  \frac {(1+\tau) \cdot  \varpi(q)}{(1+|q|)} \cdot  \sum_{|K|\leq |J|} \Big( \sum_{|J^{\prime}|+(|K|-1)_+\leq |J|} \!\!\!\!\!| \Lie_{Z^{J^{\prime}} }H_{LL}|    \cdot \sum_{X \in \cal V} |\derm ( \Lie_{Z^{K}} A )_{X} | \Big) d\tau \; ,
\eea
where
\bea
\cal V :=  \begin{cases}  \cal T \;  ,\quad\text{if }\quad V \in \cal T \; ,\\
   \cal U \; , \,\quad\text{if }\quad V \in \cal U \; . \end{cases}   
\eea
Furthermore,
\bea
   \notag
&& (1+t+|q|) \cdot |\varpi(q) \cdot \derm ( \Lie_{Z^J} h^1 )_{UV} (t,x)| \\
   \notag
 &\les&   c (\gamma^\prime) \cdot  c (\delta) \cdot c (\gamma) \cdot C ( |J|  ) \cdot E ( |J| + 4) \cdot \eps \, \\
\notag
&& + c (\gamma^\prime) \cdot c (\gamma)  \cdot c (\delta)  \cdot E ( 3)  \cdot \eps  \cdot  \int_0^t  \frac{1}{(1+\tau)} \cdot (1+\tau+|q|) \cdot \|\varpi(q)  \cdot  \derm  ( \Lie_{Z^J} h^1_{UV} ) (\tau,\cdot) \|_{L^\infty (\Sigma^{ext}_{\tau} )} d \tau \\
       \notag
    && + \sum_{|K| \leq |J|} \int_0^t (1+\tau) \cdot  \varpi(q) \cdot  \|    \Lie_{Z^K}  g^{\la\mu} \derm_{\la}   \derm_{\mu}  h^1_{UV} (\tau,\cdot) \|_{L^\infty(\overline{D}_\tau)} d\tau \\
  \notag
& &+ \int_0^t   \frac {(1+\tau) \cdot  \varpi(q)}{(1+\tau+|q|)} \cdot    \sum_{|K|\leq |J|} \Big( \sum_{|J^{\prime}|+(|K|-1)_+\le |J|} \,\,\,
|\Lie_{Z^{J^{\prime}}} H|\cdot {|\derm ( \Lie_{Z^{K} } h^1) |}  \Big) d\tau \\
\notag
& &+  \int_0^t  \frac {(1+\tau) \cdot  \varpi(q)}{(1+|q|)} \cdot   \sum_{|K|\leq |J|} \Big( \sum_{|J^{\prime}|+(|K|-1)_+\leq |J|} \!\!\!\!\!| \Lie_{Z^{J^{\prime}} }H_{LL}|    \cdot {|\derm ( \Lie_{Z^{K}} h^1 ) |}  d\tau \Big) \; .
\eea

 \end{lemma}

\begin{proof}
Based on the commutation formula that we established in Lemma \ref{Thecommutationformulawithpossibilityofsperationoftangentialcomponentsaswell}, we have
\beaa
\notag
&& \int_0^t (1+\tau)\| \varpi(q)   g^{\la\mu} \derm_{\la}   \derm_{\mu}  ( \Lie_{Z^J} A)_V (\tau,\cdot) \|_{L^\infty(D_\tau)} d\tau \\
\notag
  \les  && \int_0^t (1+\tau) \cdot  \varpi(q) \cdot \|    \Lie_{Z^J}  g^{\la\mu} \derm_{\la}   \derm_{\mu}  A_V (\tau,\cdot) \|_{L^\infty(D_\tau)} d\tau \\
  \notag
& +& \int_0^t   \frac {(1+\tau) \cdot \varpi(q)}{1+\tau+|q|} \cdot    \sum_{|K|\leq |J|} \Big( \sum_{|J^{\prime}|+(|K|-1)_+\le |J|} \,\,\,
|\Lie_{Z^{J^{\prime}}} H|\cdot {|\derm ( \Lie_{Z^{K} } A) |}  \Big) d\tau \\
\notag
& +&  \int_0^t  \frac {(1+\tau) \cdot \varpi(q)}{1+|q|} \cdot   \sum_{|K|\leq |J|} \Big( \sum_{|J^{\prime}|+(|K|-1)_+\leq |J|} \!\!\!\!\!|  \Lie_{Z^{J^{\prime}} }H_{LL}|   \cdot \sum_{X \in \cal V} |\derm ( \Lie_{Z^{K}} A )_{X} |  d\tau  \Big) \\
\notag
 && + \sum_{|K| \leq |J|} \int_0^t (1+\tau) \cdot \varpi(q) \cdot \|    \Lie_{Z^K}  g^{\la\mu} \derm_{\la}   \derm_{\mu}  A_V (\tau,\cdot) \|_{L^\infty(D_\tau)} d\tau . \\
\eeaa
As a result, injecting in Lemma \ref{GronwallinequalitiesoncomponentsofgradientofAandgradientofh1withweight}, we get
     \bea
   \notag
&& (1+t+|q|) \cdot |\varpi(q)\derm ( \Lie_{Z^J} A)_V (t,x)| \\
\notag
 &\les&   c (\gamma^\prime) \cdot  c (\delta) \cdot c (\gamma) \cdot C ( |J| ) \cdot E ( |J| + 4) \cdot \eps \, \\
\notag
&& + c (\gamma^\prime)  \cdot c (\gamma)  \cdot c (\delta)  \cdot E ( 3)  \cdot \eps  \cdot  \int_0^t  \frac{1}{(1+\tau)} \cdot (1+\tau+|q|) \cdot \|\varpi(q) \derm  ( \Lie_{Z^J} A)_V (\tau,\cdot) \|_{L^\infty} d \tau \\
       \notag
&& + \int_0^t  (1+\tau) \cdot \| \varpi(q)  g^{\la\mu} \derm_{\la}   \derm_{\mu}  ( \Lie_{Z^J} A)_V (\tau,\cdot) \|_{L^\infty(D_\tau)}  d\tau \\
  &\les &  c (\gamma^\prime)  \cdot  c (\delta) \cdot c (\gamma) \cdot C ( |J| ) \cdot E ( |J| + 4) \cdot \eps \, \\
\notag
&& + c (\gamma^\prime) \cdot c (\gamma)  \cdot c (\delta)  \cdot E ( 3)  \cdot \eps  \cdot  \int_0^t  \frac{1}{(1+\tau)} \cdot (1+\tau+|q|) \cdot \|\varpi(q) \derm  ( \Lie_{Z^J} A)_V (\tau,\cdot) \|_{L^\infty} d \tau \\
       \notag
    && + \sum_{|K| \leq |J|} \int_0^t (1+\tau) \cdot  \varpi(q)  \cdot \|    \Lie_{Z^K}  g^{\la\mu} \derm_{\la}   \derm_{\mu}  A_V (\tau,\cdot) \|_{L^\infty(D_\tau)} d\tau \\
  \notag
& &+ \int_0^t   \frac {(1+\tau) \cdot  \varpi(q)}{1+\tau+|q|} \cdot    \sum_{|K|\leq |J|} \Big( \sum_{|J^{\prime}|+(|K|-1)_+\le |J|} \,\,\,
|\Lie_{Z^{J^{\prime}}} H|\cdot {|\derm ( \Lie_{Z^{K} } A) |}  \Big) d\tau \\
\notag
& &+  \int_0^t  \frac {(1+\tau) \cdot  \varpi(q)}{1+|q|} \cdot   \sum_{|K|\leq |J|} \Big( \sum_{|J^{\prime}|+(|K|-1)_+\leq |J|} \!\!\!\!\!|  \Lie_{Z^{J^{\prime}} }H_{LL}| \cdot\sum_{X \in \cal V} |\derm ( \Lie_{Z^{K}} A )_{X} | d\tau  \Big)   \; .
\eea
Similarly for $h^1$, we get the desired estimate.
\end{proof}

\subsection{Establishing a Grönwall type inequality for the Lie derivatives of the Einstein-Yang-Mills fields}\

\begin{lemma}\label{estimateonthegooddecayingpartofthecommutatorterm}
For $\gamma^\prime$ such that $-1 \leq \gamma^\prime < \gamma - \delta$, and $ \delta <    1/2 $, we have
  \beaa
   &&    \int_0^t   \frac {(1+\tau) \cdot  \varpi(q)}{(1+\tau+|q|)} \cdot    \sum_{|K|\leq |J|} \Big( \sum_{|J^{\prime}|+(|K|-1)_+\le |J|} \,\,\,
|\Lie_{Z^{J^{\prime}}} H|\cdot {|\derm ( \Lie_{Z^{K} } A) |}  \Big) d\tau \\
&\les&  c (\delta) \cdot c (\gamma) \cdot C ( |J| ) \cdot E ( |J| + 2)  \cdot \eps^2 \; 
\notag
\eeaa
and
  \beaa
   &&    \int_0^t   \frac {(1+\tau) \cdot  \varpi(q)}{(1+\tau+|q|)} \cdot   \sum_{|K|\leq |J|} \Big( \sum_{|J^{\prime}|+(|K|-1)_+\le |J|} \,\,\,
|\Lie_{Z^{J^{\prime}}} H|\cdot {|\derm ( \Lie_{Z^{K} } h^1) |}  \Big) d\tau \\
&\les&  c (\delta) \cdot c (\gamma) \cdot C ( |J| ) \cdot E ( |J| + 2)  \cdot \eps^2 \; .
\notag
\eeaa

\end{lemma}
\begin{proof}
We showed in Lemmas \eqref{apriordecayestimatesfrombootstrapassumption} and \ref{aprioriestimatesonZLiederivativesofBIGH}, that in the exterior,
 \beaa
 \notag
|\derm ( \Lie_{Z^J} A ) (t,x)  | + |\derm ( \Lie_{Z^J} h^1 ) (t,x)  |   &\leq& \begin{cases} C ( |J| ) \cdot E ( |J| + 2)  \cdot \frac{\eps }{(1+t+|q|)^{1-\delta} (1+|q|)^{1+\ga}},\quad\text{when }\quad q>0,\\
       C ( |J| ) \cdot E ( |J| + 2)  \cdot \frac{\eps  }{(1+t+|q|)^{1-\delta}(1+|q|)^{\frac{1}{2} }}  \,\quad\text{when }\quad q<0 , \end{cases} \\
       \eeaa
       and 
              \beaa
 \notag
|  \Lie_{Z^I} H (t,x)  | &\leq& \begin{cases} c (\delta) \cdot c (\gamma) \cdot C ( |I| ) \cdot E ( |I| + 2) \cdot  \frac{\eps}{ (1+ t + | q | )^{1-\delta }  (1+| q |   )^{\de}}  ,\quad\text{when }\quad q>0,\\
\notag
    C ( |I| ) \cdot E ( |I| + 2) \cdot  \frac{\eps}{ (1+ t + | q | )^{1-\delta }  } (1+| q |   )^{\frac{1}{2} }  , \,\quad\text{when }\quad q<0 . \end{cases} \\ 
    \eeaa
    Thus, for all $|J^{\prime}|, |K|\leq |J|$,
    \beaa
   |\Lie_{Z^{J^{\prime}}} H|\cdot {|\derm ( \Lie_{Z^{K} } A) |}   &\leq& \begin{cases}  c (\delta) \cdot c (\gamma) \cdot C ( |J| ) \cdot E ( |J| + 2)  \cdot \frac{\eps^2 }{(1+t+|q|)^{2-2\delta} (1+|q|)^{1+\ga + \de}},\quad\text{when }\quad q>0,\\
       C ( |J| ) \cdot E ( |J| + 2)  \cdot \frac{\eps^2  }{(1+t+|q|)^{2-2\delta}}  \,\quad\text{when }\quad q<0 . \end{cases} \\
       \eeaa
Hence, in the exterior, we have
           \beaa
  && \varpi(q) \cdot  |\Lie_{Z^{J^{\prime}}} H|\cdot {|\derm ( \Lie_{Z^{K} } A) |}   \\
  \notag
  &\leq& \begin{cases}   c (\delta) \cdot c (\gamma) \cdot C ( |J| ) \cdot E ( |J| + 2)  \cdot \frac{\eps^2 \cdot  (1+|q|)^{1+\gamma^\prime} } {(1+t+|q|)^{2-2\delta} (1+|q|)^{1+\ga+\de}},\quad\text{when }\quad q>0,\\
       C ( |J| ) \cdot E ( |J| + 2)  \cdot \frac{\eps^2  }{(1+t+|q|)^{2-2\delta}}  \,\quad\text{when }\quad q<0 , \end{cases} \\
         &\leq&  \begin{cases}     c (\delta) \cdot c (\gamma)  \cdot C ( |J| ) \cdot E ( |J| + 2)  \cdot   \frac{\eps^2 }{(1+t+|q|)^{2-2\delta} (1+|q|)^{-\ga^\prime+\ga+\de} }, \quad\text{when }\quad q>0,\\
           C ( |J| ) \cdot E ( |J| + 2)  \cdot \frac{\eps^2  }{(1+t+|q|)^{2-2\delta}}  \,\quad\text{when }\quad q<0 , \end{cases} \\
            &\leq&   c (\delta) \cdot c (\gamma)  \cdot C ( |J| ) \cdot E ( |J| + 2)  \cdot   \frac{\eps^2 }{(1+t+|q|)^{2-2\delta}  } \; .
       \eeaa
       Thus, for $\gamma^\prime$ such that $-1 \leq \gamma^\prime < \gamma - \delta$,
       \beaa
    && \frac {(1+ t)\cdot  \varpi(q)}{1+t+|q|}   \cdot  |\Lie_{Z^{J^{\prime}}} H|\cdot {|\derm ( \Lie_{Z^{K} } A) |}   \\
    &\leq&   c (\delta) \cdot c (\gamma)  \cdot C ( |J| ) \cdot E ( |J| + 2)  \cdot   \frac{\eps^2 }{(1+t+|q|)^{2-2\delta}  } \; .
       \eeaa
       This means that if $ \delta <    1/2 $, then the right hand side is integrable and we get
      \beaa
      \notag
    &&   \int_0^t   \frac {(1+\tau) \cdot \varpi(q)}{1+\tau+|q|} \cdot  \Big(   \sum_{|K|\leq |J|,}\,\, \sum_{|J^{\prime}|+(|K|-1)_+\le |J|} \,\,\,
|\Lie_{Z^{J^{\prime}}} H|\cdot {|\derm ( \Lie_{Z^{K} } A) |}  \Big) d\tau \\
&\les&  c (\delta) \cdot c (\gamma) \cdot C ( |J| ) \cdot E ( |J| + 2)  \cdot \eps^2 \;  .\\
\notag
\eeaa
Since the same estimates used here for $A$ hold for $h^1$, we get the result for $h^1$ as well.

     \end{proof}

We now examine the other term in the commutator estimate. For that term, we are in fact going to use the induction hypothesis. For this, we start first by translating an estimate on $h$ into an estimate on $H$ so that we could use either the bootstrap assumption on $h^1$ or the induction hypothesis on $h^1$.

\begin{lemma}
We have
\beaa
&&  \sum_{|I|\leq \, |J|}| \Lie_{Z^I} H_{LL}|  \\
 &\les&  \sum_{|I|\leq \, |J|}| \Lie_{Z^I} h_{LL} |    \\
&&  +    c (\delta) \cdot c (\gamma) \cdot C ( |J| ) \cdot E ( |J| + 2) \cdot  \begin{cases}   \frac{\eps^2}{ (1+ t + | q | )^{2-2\delta }  (1+| q |   )^{2\de}}  ,\quad\text{when }\quad q>0 \; ,\\
\notag
     \frac{\eps^2}{ (1+ t + | q | )^{2-2\delta }  } (1+| q |   )  , \,\quad\text{when }\quad q<0\; . \end{cases} 
\eeaa
 Consequently, using the harmonic gauge estimate, for $\eps \leq 1$, and for $\de \leq \frac{1}{4}$, we have
       \beaa
       &&  \sum_{|I|\leq \, |J|}| \Lie_{Z^I} H_{LL}|    \\
      &\les& \sum_{|I|\leq \, |J|-2} \,\,\,\int\limits_{s,\Om=const} |\derm \Lie_{Z^I} h| \\
&&  +     c (\delta) \cdot c (\gamma) \cdot C ( |J| ) \cdot E ( |J| + 3) \cdot  \begin{cases}   \frac{\eps}{ (1+ t + | q | ) }  ,\quad\text{when }\quad q>0 \; ,\\
\notag
     \frac{\eps}{ (1+ t + | q | )} (1+| q |   )^{1/2+2\delta}  , \,\quad\text{when }\quad q<0 \; , \end{cases} 
     \eeaa
where the $\sum_{|I|\leq \, |J|-2} \,\,\,\int\limits_{s,\Om=const} |\derm \Lie_{Z^I} h| $ is vanishing when $|J| \leq 1$.
\end{lemma}

\begin{proof}

We showed in \cite{G4}, that  $H_{\mu\nu}=-h_{\mu\nu}+ O_{\mu\nu}(h^2) $, thus,
  \beaa
\Lie_{Z^I}H_{\mu\nu}=- \Lie_{Z^I} h_{\mu\nu}+ \sum_{|J| + |K| \leq |I|} O_{\mu\nu}(\Lie_{Z^J} h \cdot \Lie_{Z^K} h ) \; .
\eeaa
Hence,
\beaa
| \Lie_{Z^I}H_{\mu\nu} | &\les& | \Lie_{Z^I} h_{\mu\nu}| + \sum_{|J| + |K| \leq |I|}  |\Lie_{Z^J} h | \cdot |\Lie_{Z^K} h | \\
&\les & | \Lie_{Z^I} h_{\mu\nu}| +      \begin{cases} c (\delta) \cdot c (\gamma) \cdot C ( |I| ) \cdot E ( |I| + 2) \cdot  \frac{\eps^2}{ (1+ t + | q | )^{2-2\delta }  (1+| q |   )^{2\de}}  ,\quad\text{when }\quad q>0 \; ,\\
\notag
    C ( |I| ) \cdot E ( |I| + 2) \cdot  \frac{\eps^2}{ (1+ t + | q | )^{2-2\delta }  } (1+| q |   )  , \,\quad\text{when }\quad q<0\; , \end{cases} 
     \eeaa

Consequently,
\bea\label{linkboigHLLandsnmallhLL}
&&  \sum_{|I|\leq \, |J|}| \Lie_{Z^I} H_{LL}| (t, | x | \cdot \Om) \\
\notag
&\les&  \sum_{|I|\leq \, |J|}| \Lie_{Z^I} h_{LL} | (t, | x | \cdot \Om)  \\
\notag
&&  +     c (\delta) \cdot c (\gamma) \cdot C ( |J| ) \cdot E ( |J| + 2) \cdot \begin{cases}   \frac{\eps^2}{ (1+ t + | q | ) }  ,\quad\text{when }\quad q>0 \; ,\\
\notag
     \frac{\eps^2}{ (1+ t + | q | )} (1+| q |   )^{1/2+\delta}  , \,\quad\text{when }\quad q<0 \; . \end{cases} 
     \eea

By integration, and taking into consideration Definition \ref{definitionforintegralalongthenullcoordinateplusboundaryterm} which absorbs the boundary term, we have
 \beaa
 | \Lie_{Z^I} h | &\leq& \int\limits_{s,\, \Om=const} |\derm \Lie_{Z^I} h| \; .
 \eeaa
 Consequently, given the estimate in Lemma \ref{estimateonpartialderivativeofALcomponent}, derived from the harmonic gauge under the bootstrap assumption holding true for all $|J| \leq  \lfloor \frac{|I|}{2} \rfloor  $, we get that

           \beaa
        \notag
 \sum_{|J|\leq |I| }  |  \Lie_{Z^I} h_{ L L} |    &\les&  \int\limits_{s,\, \Om=const} \sum_{|J|\leq |I| -2}  |\derm \Lie_{Z^J} h |  \\
&& + \begin{cases} c (\delta) \cdot c (\gamma) \cdot C ( |I| ) \cdot E ( |I| + 3)  \cdot \frac{ \eps   }{ (1+t+|q|) } ,\quad\text{when }\quad q>0 \; ,\\
        C ( |I| ) \cdot E ( |I| + 3) \cdot  \frac{ \eps \cdot (1+|q|)^{\frac{1}{2} + 2 \delta } }{ (1+t+|q|) }  \,\quad\text{when }\quad q<0 \; . \end{cases} \\
       \eeaa    
  
       Thus, using the above and injecting in \eqref{linkboigHLLandsnmallhLL}, we get for $\eps \leq 1$\,, and for $\de \leq \frac{1}{4}$\,,
       \beaa
       &&  \sum_{|I|\leq \, |J|}| \Lie_{Z^I} H_{LL}|   \\
 &\les&  \sum_{|I|\leq \, |J|}| \Lie_{Z^I} h_{LL} |    \\
&& +    c (\delta) \cdot c (\gamma) \cdot C ( |J| ) \cdot E ( |J| + 2) \cdot  \begin{cases}   \frac{\eps^2}{ (1+ t + | q | )^{2-2\delta }  (1+| q |   )^{2\de}}  ,\quad\text{when }\quad q>0 \; ,\\
\notag
     \frac{\eps^2}{ (1+ t + | q | )^{2-2\delta }  } (1+| q |   )  , \,\quad\text{when }\quad q<0\; , \end{cases} \\
     \notag
 &\les& \sum_{|I|\leq \, |J|-2} \,\,\,\int\limits_{s,\, \Om=const}
|\derm \Lie_{Z^I} h| \\
&&  +     c (\delta) \cdot c (\gamma) \cdot C ( |J| ) \cdot E ( |J| + 3) \cdot  \begin{cases}   \frac{\eps}{ (1+ t + | q | ) }  ,\quad\text{when }\quad q>0 \; ,\\
\notag
     \frac{\eps}{ (1+ t + | q | )} (1+| q |   )^{1/2+2\delta}  , \,\quad\text{when }\quad q<0 \; , \end{cases} \\ 
&&      +    c (\delta) \cdot c (\gamma) \cdot C ( |J| ) \cdot E ( |J| + 2) \cdot  \begin{cases}   \frac{\eps^2}{ (1+ t + | q | ) \cdot (1+ t + | q | )^{1-2\delta }  (1+| q |   )^{2\de}}  ,\quad\text{when }\quad q>0 \; ,\\
\notag
     \frac{\eps^2 \cdot (1+| q |   )^{\frac{1}{2}} }{ (1+ t + | q | )  } \cdot \frac{(1+| q |   )^{\frac{1}{2}}}{(1+ t + | q | )^{1-2\delta }} , \,\quad\text{when }\quad q<0\; , \end{cases} \\
      &\les& \sum_{|I|\leq \, |J|-2} \,\,\,\int\limits_{s,\, \Om=const}
|\derm \Lie_{Z^I} h| \\
&&  +     c (\delta) \cdot c (\gamma) \cdot C ( |J| ) \cdot E ( |J| + 3) \cdot  \begin{cases}   \frac{\eps}{ (1+ t + | q | ) }  ,\quad\text{when }\quad q>0 \; ,\\
\notag
     \frac{\eps}{ (1+ t + | q | )} (1+| q |   )^{1/2+2\delta}  , \,\quad\text{when }\quad q<0 \; . \end{cases} 
     \eeaa
Similarly, we get the estimate for $\sum_{|I|\leq \, |J|-2} | \Lie_{Z^I}H |$. 
     \end{proof}

Now, we would like to use the induction hypothesis. First, let us point out the following.

\begin{lemma}\label{estimateonzeroderivativeofhunderinducthypothesi}
For $M \leq \eps$, the induction hypothesis on the Einstein-Yang-Mills metric $h^1$, gives that in that exterior, for all $|K|\leq |J| -1$,
       \bea
 \notag
&& \int\limits_{s,\, \Om=const} \sum_{|I|\leq \, |J|-2}   |\derm  ( \Lie_{Z^I} h)   |    \\
\notag
 &\leq&   C(q_0)   \cdot  c (\delta) \cdot c (\gamma) \cdot C(|J|) \cdot E (  |J|  +2)  \cdot \frac{\eps }{(1+t+|q|)^{1-   c (\gamma)  \cdot c (\delta)  \cdot c(|J|) \cdot E ( |J|+ 2)\cdot  \eps } }     \; . \\
      \eea
      
\end{lemma}

     \begin{proof}
In fact, the induction hypothesis is on $h^1$, that in the exterior, for $|K|\leq |J| -1  $, 
            \bea\label{formulaforinductionhypothesisonh^1withKlesJminues1}
 \notag
&& |\derm  ( \Lie_{Z^K} h^1)   |   \\
\notag
 &\leq&   C(q_0)   \cdot  c (\delta) \cdot c (\gamma) \cdot C(|K|) \cdot E (  |K|  +4)  \cdot \frac{\eps }{(1+t+|q|)^{1-   c (\gamma)  \cdot c (\delta)  \cdot c(|K|) \cdot E ( |K|+ 4)\cdot  \eps } \cdot (1+|q|)}     \; . \\
      \eea
Yet, we know that  
 \beaa
 \notag
|\derm (  \Lie_{Z^K} h^0 )   |   &\leq&  C ( |K| )   \cdot \frac{\eps }{(1+t+|q|)^{2} } \; .
      \eeaa
      Thus,
                  \bea
 \notag
&& |\derm  ( \Lie_{Z^K} h)   |  \\
 \notag
 &\leq& |\derm  ( \Lie_{Z^K} h^1)   | +  |\derm  ( \Lie_{Z^K} h^0)   | \\
\notag
 &\leq&   C(q_0)   \cdot  c (\delta) \cdot c (\gamma) \cdot C(|K|) \cdot E (  |K|  +4)  \cdot \frac{\eps }{(1+t+|q|)^{1-   c (\gamma)  \cdot c (\delta)  \cdot c(|K|) \cdot E ( |K|+ 4)\cdot  \eps } \cdot (1+|q|)}     \; . \\
      \eea
            Hence, by integrating, we obtain
            \bea
 \notag
&& \int\limits_{s,\Om=const} \sum_{|I|\leq \, |J|-2}   |\derm  ( \Lie_{Z^I} h) (t,x)  |    \\
\notag
 &\leq&   C(q_0)   \cdot  c (\delta) \cdot c (\gamma) \cdot C(|J|) \cdot E (  |J|  +2)  \cdot \frac{\eps }{(1+t+|q|)^{1-   c (\gamma)  \cdot c (\delta)  \cdot c(|J|) \cdot E ( |J|+ 2)\cdot  \eps } }     \; .\\
      \eea
      
           \end{proof}
\begin{lemma}\label{estimateonthebadpartofthecommutatortermtobtainaGronwallforcomponents}

For $M \leq \eps$, under the induction hypothesis on both $A$ and on $h^1$, for  $|K|\leq |J| -1$, we have for $\gamma^\prime$ such that $-1 \leq \gamma^\prime < \gamma - \delta$, and $ \delta \leq  \frac{1}{4} $, and for $\cal V \in \{\cal T\; , \cal U\}$,
\beaa
& &  \int_0^t  \frac {(1+\tau) \cdot \varpi(q)}{1+|q|} \cdot   \sum_{|K|\leq |J|} \Big( \sum_{|J^{\prime}|+(|K|-1)_+\leq |J|} \!\!\!\!\!| \Lie_{Z^{J^{\prime}} }H_{LL}|   \cdot\sum_{X \in \cal V} |\derm ( \Lie_{Z^{K}} A )_{X} | d\tau  \Big)\; . \\
&\leq &   \int_0^t  C(q_0)   \cdot  c (\delta) \cdot c (\gamma) \cdot C(|J|) \cdot E ( |J| + 3)  \cdot \frac{\eps \cdot  \varpi(q) }{(1+t+|q|)^{1-      c (\gamma)  \cdot c (\delta)  \cdot c(|J|) \cdot E ( |J|+ 3) \cdot \eps } \cdot (1+|q|)^{2+\gamma - 2 \de }}  d\tau \\
&&+    \int_0^t    c (\delta) \cdot c (\gamma) \cdot E (  4 )  \cdot \frac{\eps }{(1+t+|q|) }     \cdot   (1+t+|q|) \cdot \varpi(q) \cdot \sum_{|K| = |J|} \sum_{X \in \cal V} |\derm ( \Lie_{Z^{K}} A )_{X} |  d\tau \; ,
 \eeaa
 and
 \beaa
& &  \int_0^t  \frac {(1+\tau) \cdot \varpi(q)}{1+|q|} \cdot \sum_{|K|\leq |J|} \Big(  \sum_{|J^{\prime}|+(|K|-1)_+\leq |J|} \!\!\!\!\!|  \Lie_{Z^{J^{\prime}} }H_{LL}|   \cdot {|\derm ( \Lie_{Z^{K}} h^1 ) |}  d\tau  \Big)  \; . \\
&\leq &   \int_0^t  C(q_0)   \cdot  c (\delta) \cdot c (\gamma) \cdot C(|J|) \cdot E ( |J| + 3)  \cdot \frac{\eps \cdot  \varpi(q) }{(1+t+|q|)^{1-      c (\gamma)  \cdot c (\delta)  \cdot c(|J|) \cdot E ( |J|+ 3) \cdot \eps } \cdot (1+|q|)^2}  d\tau \\
&&+    \int_0^t    c (\delta) \cdot c (\gamma) \cdot E (  4 )  \cdot \frac{\eps }{(1+t+|q|) }     \cdot   (1+t+|q|) \cdot \varpi(q) \cdot  \sum_{|K| = |J|} |\derm ( \Lie_{Z^{K}} h^1 ) |  d\tau \; .
 \eeaa

\end{lemma}

\begin{proof}

We decompose the following sum as  
\beaa
& &  \sum_{|K|\leq |J|}\Big(\sum_{|J^{\prime}|+(|K|-1)_+\leq |J|} \!\!\!\!\!|  \Lie_{Z^{J^{\prime}} }H_{LL}|   \cdot  \sum_{X \in \cal V} |\derm ( \Lie_{Z^{K}} A )_{X} |  \Big)    \\
&\les &  \sum_{|K| \leq  |J| -1}\Big(\sum_{|J^{\prime}|+(|K|-1)_+\leq |J|} \!\!\!\!\!|  \Lie_{Z^{J^{\prime}} }H_{LL}|   \cdot \sum_{X \in \cal V} |\derm ( \Lie_{Z^{K}} A )_{X} |   \Big)   \\
& &  + \sum_{|K| = |J|}\Big(\sum_{|J^{\prime}|+(|K|-1)_+\leq |J|} \!\!\!\!\!|  \Lie_{Z^{J^{\prime}} }H_{LL}|   \cdot  \sum_{X \in \cal V} |\derm ( \Lie_{Z^{K}} A )_{X} |  \Big)   \\
&\les &  \sum_{|K| \leq  |J| -1}\Big(\sum_{|I| \leq |J|} \!\!\!\!\ |  \Lie_{Z^{I} }H_{LL}|  \cdot  \sum_{X \in \cal V} |\derm ( \Lie_{Z^{K}} A )_{X} |  \Big)    \\
& &  + \sum_{|K| = |J|} \Big( \sum_{|I| \leq 1} \!\!\!\!\ |  \Lie_{Z^{I} }H_{LL}|    \cdot  \sum_{X \in \cal V} |\derm ( \Lie_{Z^{K}} A )_{X} |  \Big)   \;.
\eeaa
We also have the same estimate for $h^1$ instead of $A$.\\

\textbf{The sum for $|K| \leq |J| - 1$ :}

For $\eps \leq 1$, and for $\de \leq \frac{1}{4}$, we have
       \beaa
       &&  \sum_{|I|\leq \, |J|}| \Lie_{Z^I} H_{LL}|  \\
      &\les& \sum_{|I|\leq \, |J|-2} \,\,\,\int\limits_{s,\Om=const} |\derm \Lie_{Z^I} h| \\
&&  +     c (\delta) \cdot c (\gamma) \cdot C ( |J| ) \cdot E ( |J| + 3) \cdot  \begin{cases}   \frac{\eps}{ (1+ t + | q | ) }  ,\quad\text{when }\quad q>0 \; ,\\
\notag
     \frac{\eps}{ (1+ t + | q | )} (1+| q |   )^{1/2+2\delta}  , \,\quad\text{when }\quad q<0 \; . \end{cases} 
     \eeaa

  Thus, for $M\leq \eps$, using what we have shown in Lemma \ref{estimateonzeroderivativeofhunderinducthypothesi} under the induction hypothesis on $h^1$, and given that $E ( |J|+ 2) \leq E ( |J|+ 3)$, 
         \beaa
       &&  \sum_{|I|\leq \, |J|}| \Lie_{Z^I} H_{LL}|   \\
      &\les&    C(q_0)   \cdot  c (\delta) \cdot c (\gamma) \cdot C(|J|) \cdot E (  |J|  +2)  \cdot \frac{\eps }{(1+t+|q|)^{1-   c (\gamma)  \cdot c (\delta)  \cdot c(|J|) \cdot E ( |J|+ 2)\cdot  \eps } }  \\
&&  +     c (\delta) \cdot c (\gamma) \cdot C ( |J| ) \cdot E ( |J| + 3) \cdot  \begin{cases}   \frac{\eps}{ (1+ t + | q | ) }  ,\quad\text{when }\quad q>0 \; ,\\
\notag
     \frac{\eps}{ (1+ t + | q | )} (1+| q |   )^{1/2+2\delta}  , \,\quad\text{when }\quad q<0 \; , \end{cases} \\
\notag
 &\leq&   C(q_0)   \cdot  c (\delta) \cdot c (\gamma) \cdot C(|J|) \cdot E (  |J|  +3)  \cdot \frac{\eps }{(1+t+|q|)^{1-   c (\gamma)  \cdot c (\delta)  \cdot c(|J|) \cdot E ( |J|+ 2)\cdot  \eps } }     \; .\\
      \eeaa

Whereas the induction hypothesis on the Einstein-Yang-Mills potential $A$ gives for all $|K|\leq |J| -1  $,
      \beaa
 \notag
&& |\derm  ( \Lie_{Z^K} A)   |    \\
\notag
&\leq&   C(q_0)   \cdot  c (\delta) \cdot c (\gamma) \cdot C(|K|) \cdot E ( |K| + 4)  \cdot \frac{\eps }{(1+t+|q|)^{1-      c (\gamma)  \cdot c (\delta)  \cdot c(|K|) \cdot E ( |K|+ 4) \cdot \eps } \cdot (1+|q|)^{1+\gamma - 2 \de }}    \; \\
&\leq&   C(q_0)   \cdot  c (\delta) \cdot c (\gamma) \cdot C(|J|) \cdot E ( |J| + 3)  \cdot \frac{\eps }{(1+t+|q|)^{1-      c (\gamma)  \cdot c (\delta)  \cdot c(|J|) \cdot E ( |J|+ 3) \cdot \eps } \cdot (1+|q|)^{1+\gamma - 2 \de }}    \; .\\
      \eeaa
Therefore, 
\beaa
&&  \sum_{|K| \leq  |J| -1}\Big(\sum_{|I| \leq |J|} \!\!\!\!\ | \Lie_{Z^{I} }H_{LL}|     \cdot \sum_{X \in \cal V} |\derm ( \Lie_{Z^{K}} A )_{X} | \Big)   \\
&\les&   C(q_0)   \cdot  c (\delta) \cdot c (\gamma) \cdot C(|J|) \cdot E ( |J| + 3)  \cdot \frac{\eps }{(1+t+|q|)^{2-      c (\gamma)  \cdot c (\delta)  \cdot c(|J|) \cdot E ( |J|+ 3) \cdot \eps } \cdot (1+|q|)^{1+\gamma - 2 \de }}    \; .\\
      \eeaa
Thus, for $\eps \leq 1$, and for $\de \leq \frac{1}{4}$,
\bea
\notag
&& \frac {(1+t) \cdot \varpi(q)}{1+|q|} \cdot  \sum_{|K| \leq  |J| -1}\Big(\sum_{|I| \leq |J|} \!\!\!\!\ | \Lie_{Z^{I} }H_{LL}|    \cdot\sum_{X \in \cal V} |\derm ( \Lie_{Z^{K}} A )_{X} |  \Big) \\
\notag
&\les&   C(q_0)   \cdot  c (\delta) \cdot c (\gamma) \cdot C(|J|) \cdot E ( |J| + 3)  \cdot \frac{\eps \cdot  \varpi(q) }{(1+t+|q|)^{1-      c (\gamma)  \cdot c (\delta)  \cdot c(|J|) \cdot E ( |J|+ 3) \cdot \eps } \cdot (1+|q|)^{2+\gamma - 2 \de }}    \; . \\
      \eea
Whereas to $h^1$, since the induction hypothesis gives that for all $|K|\leq |J| -1  $,
      \beaa
 \notag
&& |\derm  ( \Lie_{Z^K} h^1)   |    \\
\notag
&\leq&   C(q_0)   \cdot  c (\delta) \cdot c (\gamma) \cdot C(|K|) \cdot E ( |K| + 4)  \cdot \frac{\eps }{(1+t+|q|)^{1-      c (\gamma)  \cdot c (\delta)  \cdot c(|K|) \cdot E ( |K|+ 4) \cdot \eps } \cdot (1+|q|)}    \; \\
&\leq&   C(q_0)   \cdot  c (\delta) \cdot c (\gamma) \cdot C(|J|) \cdot E ( |J| + 3)  \cdot \frac{\eps }{(1+t+|q|)^{1-      c (\gamma)  \cdot c (\delta)  \cdot c(|J|) \cdot E ( |J|+ 3) \cdot \eps } \cdot (1+|q|)}    \; ,\\
      \eeaa
    similarly, we obtain
    \bea
\notag
&& \frac {(1+t) \cdot \varpi(q)}{1+|q|} \cdot  \sum_{|K| \leq  |J| -1}\Big(\sum_{|I| \leq |J|} \!\!\!\!\ | \Lie_{Z^{I} }H_{LL}|   \cdot {|\derm ( \Lie_{Z^{K}} h^1 ) |}  \Big)  \\
\notag
&\les&   C(q_0)   \cdot  c (\delta) \cdot c (\gamma) \cdot C(|J|) \cdot E ( |J| + 3)  \cdot \frac{\eps \cdot  \varpi(q) }{(1+t+|q|)^{1-      c (\gamma)  \cdot c (\delta)  \cdot c(|J|) \cdot E ( |J|+ 3) \cdot \eps } \cdot (1+|q|)^2}    \; . \\
      \eea
      
\textbf{The sum for $|K| = |J| $ :}

The harmonic gauge estimate gives
        \beaa
       &&  \sum_{|I|\leq \, 1}| \Lie_{Z^I} H_{LL}|    \\
      &\les&     c (\delta) \cdot c (\gamma)  \cdot E ( 4 ) \cdot  \begin{cases}   \frac{\eps}{ (1+ t + | q | ) }  ,\quad\text{when }\quad q>0 \; ,\\
\notag
     \frac{\eps}{ (1+ t + | q | )} (1+| q |   )^{1/2+2\delta}  , \,\quad\text{when }\quad q<0 \; . \end{cases} 
     \eeaa
 Thus,
   \beaa
       && \frac{1}{(1+|q|)} \sum_{|I|\leq \, 1}| \Lie_{Z^I} H_{LL}|   \\
      &\les&     c (\delta) \cdot c (\gamma)  \cdot E ( 4 ) \cdot  \begin{cases}   \frac{\eps}{ (1+ t + | q | ) \cdot (1+|q|)  }  ,\quad\text{when }\quad q>0 \; ,\\
\notag
     \frac{\eps}{ (1+ t + | q | )} (1+| q |   )^{-1/2+2\delta}  , \,\quad\text{when }\quad q<0 \; , \end{cases} \\
 &\leq&    c (\delta) \cdot c (\gamma) \cdot E (  4 )  \cdot \frac{\eps }{(1+t+|q|)}     \\
 && \text{(since $\de \leq \frac{1}{4} \leq \frac{1}{2} $) .}
      \eeaa

Thus,
\beaa
 &&   \frac{1}{(1+|q|)}  \sum_{|K| = |J|} \cdot \Big( \sum_{|I| \leq 1} \!\!\!\!\ | \Lie_{Z^{I} }H_{LL}|   \cdot \sum_{X \in \cal V} |\derm ( \Lie_{Z^{K}} A )_{X} | \Big)  \\
\notag
 &\leq&   c (\delta) \cdot c (\gamma) \cdot E (  4 )  \cdot \frac{\eps }{(1+t+|q|)}     \cdot \sum_{|K| = |J|}\sum_{X \in \cal V} |\derm ( \Lie_{Z^{K}} A )_{X} | \; . \\
      \eeaa
      Hence,
\bea
\notag
 &&  \frac {(1+t) \cdot \varpi(q)}{1+|q|} \cdot  \sum_{|K| = |J|} \Big( \sum_{|I| \leq 1} \!\!\!\!\ | \Lie_{Z^{I} }H_{LL}|   \cdot \sum_{X \in \cal V} |\derm ( \Lie_{Z^{K}} A )_{X} | \Big)  \\
\notag
 &\leq&    c (\delta) \cdot c (\gamma) \cdot E (  4 )  \cdot \frac{\eps }{(1+t+|q|)  }     \cdot   (1+t+|q|) \cdot \varpi(q) \cdot \sum_{|K| = |J|} \sum_{X \in \cal V} |\derm ( \Lie_{Z^{K}} A )_{X} |  \; . \\
      \eea
And, similarly for $h^1$,
\bea
\notag
 &&  \frac {(1+t) \cdot \varpi(q)}{1+|q|} \cdot  \sum_{|K| = |J|} \Big( \sum_{|I| \leq 1} \!\!\!\!\ | \Lie_{Z^{I} }H_{LL}|  \cdot {|\derm ( \Lie_{Z^{K}} h^1 ) |}   \Big) \\
\notag
 &\leq&    c (\delta) \cdot c (\gamma) \cdot E (  4 )  \cdot \frac{\eps }{(1+t+|q|)  }     \cdot   (1+t+|q|) \cdot \varpi(q) \cdot \sum_{|K| = |J|} |\derm ( \Lie_{Z^{K}} h^1 ) |   \; . \\
      \eea
      
      Finally, put the two terms together, we obtain the result.

\end{proof}

\subsection{The commutator term and the Grönwall inequality on the Lie derivatives of the Einstein-Yang-Mills fields}\

      \begin{lemma}\label{almostGrnwallinequalitongradientofAandh1afterestimatedthecommutator}
       For $M \leq \eps \leq 1$, under the induction hypothesis holding true for both $A$ and on $h^1$, for all $|K|\leq |J| -1$, and for $\gamma^\prime$ such that $-1 \leq \gamma^\prime < \gamma - \delta$, and  $ \delta \leq  \frac{1}{4} $, we have for any $U, V \in \cal U$,
and for
 \bea
\cal V :=  \begin{cases}  \cal T \;  ,\quad\text{if }\quad V \in \cal T \; ,\\
   \cal U \; , \,\quad\text{if }\quad V \in \cal U \; , \end{cases}   
\eea
the following estimates
      \beaa
   \notag
&& (1+t+|q|) \cdot |\varpi(q) \cdot \derm ( \Lie_{Z^J} A)_V | \\
   \notag
   &\les&   c (\gamma^\prime) \cdot  c (\delta) \cdot c (\gamma) \cdot C ( |J|  ) \cdot E ( |J| + 4) \cdot \eps \, \\
\notag
&& + c (\gamma^\prime)   \cdot  c (\delta) \cdot c (\gamma) \cdot E (  4 )  \cdot \eps  \cdot  \int_0^t  \frac{1}{(1+\tau)} \cdot (1+\tau+|q|) \cdot  \sum_{X \in \cal V}  \|\varpi(q)  \cdot  \derm  ( \Lie_{Z^J} A)_{X} (\tau,\cdot) \|_{L^\infty (\Sigma^{ext}_{\tau} )} d \tau \\
       \notag
    && + \sum_{|K| \leq |J|} \int_0^t (1+\tau) \cdot  \varpi(q) \cdot  \|    \Lie_{Z^K}  g^{\la\mu} \derm_{\la}   \derm_{\mu}  A_V (\tau,\cdot) \|_{L^\infty(\overline{D}_\tau)} d\tau \\
  \notag
& & +  C(q_0)   \cdot  c (\delta) \cdot c (\gamma) \cdot C(|J|) \cdot E ( |J| + 3)  \cdot  \int_0^t  \frac{\eps \cdot  \varpi(q) }{(1+t+|q|)^{1- c (\gamma)  \cdot c (\delta)  \cdot c(|J|) \cdot E ( |J|+ 3) \cdot \eps } \cdot (1+|q|)^{2+\gamma - 2 \de }}  d\tau \; ,
 \eeaa
and
 \beaa
   \notag
&& (1+t+|q|) \cdot |\varpi(q) \cdot \derm ( \Lie_{Z^J} h^1)_{UV} | \\
   \notag
   &\les&   c (\gamma^\prime) \cdot  c (\delta) \cdot c (\gamma) \cdot C ( |J|  ) \cdot E ( |J| + 4) \cdot \eps \, \\
\notag
&& + c (\gamma^\prime)   \cdot  c (\delta) \cdot c (\gamma) \cdot E (  4 )  \cdot \eps  \cdot  \int_0^t  \frac{1}{(1+\tau)} \cdot (1+\tau+|q|) \cdot \|\varpi(q)  \cdot  \derm  ( \Lie_{Z^J} h^1) (\tau,\cdot) \|_{L^\infty (\Sigma^{ext}_{\tau} )} d \tau \\
       \notag
    && + \sum_{|K| \leq |J|} \int_0^t (1+\tau) \cdot  \varpi(q) \cdot  \|    \Lie_{Z^K}  g^{\la\mu} \derm_{\la}   \derm_{\mu}  h^1_{UV} (\tau,\cdot) \|_{L^\infty(\overline{D}_\tau)} d\tau \\
  \notag
& & +   C(q_0)   \cdot  c (\delta) \cdot c (\gamma) \cdot C(|J|) \cdot E ( |J| + 3)  \cdot  \int_0^t \frac{\eps \cdot  \varpi(q) }{(1+t+|q|)^{1-      c (\gamma)  \cdot c (\delta)  \cdot c(|J|) \cdot E ( |J|+ 3) \cdot \eps } \cdot (1+|q|)^2}  d\tau \; .
 \eeaa
      \end{lemma}
      
       \begin{proof}
       
      Based on what we proved in Lemmas \ref{estimateonthegooddecayingpartofthecommutatorterm} and \ref{estimateonthebadpartofthecommutatortermtobtainaGronwallforcomponents} under the assumptions of this lemma, and injecting in the estimates in Lemma \ref{estimatethatallowsupgradeincorporatingtermsfromthecommutationformula}, we get that for any $V \in \cal U$,
\beaa
   \notag
&& (1+t+|q|) \cdot |\varpi(q) \cdot \derm ( \Lie_{Z^J} A)_V (t,x)| \\
   \notag
 &\les&   c (\gamma^\prime) \cdot  c (\delta) \cdot c (\gamma) \cdot C ( |J|  ) \cdot E ( |J| + 4) \cdot \eps \, \\
\notag
&& + c (\gamma^\prime) \cdot c (\gamma)  \cdot c (\delta)  \cdot E ( 3)  \cdot \eps  \cdot  \int_0^t  \frac{1}{(1+\tau)} \cdot (1+\tau+|q|) \cdot \|\varpi(q)  \cdot  \derm  ( \Lie_{Z^J} A)_{V} (\tau,\cdot) \|_{L^\infty (\Sigma^{ext}_{\tau} )} d \tau \\
       \notag
    && + \sum_{|K| \leq |J|} \int_0^t (1+\tau) \cdot  \varpi(q) \cdot  \|    \Lie_{Z^K}  g^{\la\mu} \derm_{\la}   \derm_{\mu}  A_V (\tau,\cdot) \|_{L^\infty(\overline{D}_\tau)} d\tau \\
  \notag
& &+ \int_0^t   \frac {(1+\tau) \cdot  \varpi(q)}{(1+\tau+|q|)} \cdot  \Big(   \sum_{|K|\leq |J|,}\,\, \sum_{|J^{\prime}|+(|K|-1)_+\le |J|} \,\,\,
|\Lie_{Z^{J^{\prime}}} H|\cdot {|\derm ( \Lie_{Z^{K} } A) |}  \Big) d\tau \\
\notag
& &+  \int_0^t  \frac {(1+\tau) \cdot  \varpi(q)}{(1+|q|)} \cdot
 \sum_{|K|\leq |J|}\Big(\sum_{|J^{\prime}|+(|K|-1)_+\leq |J|} \!\!\!\!\!| \Lie_{Z^{J^{\prime}} }H_{LL}|    \cdot \sum_{X \in \cal V} |\derm ( \Lie_{Z^{K}} A )_{X} |  \Big)  d\tau \; 
\eeaa
\beaa
   &\les&   c (\gamma^\prime) \cdot  c (\delta) \cdot c (\gamma) \cdot C ( |J|  ) \cdot E ( |J| + 4) \cdot \eps \, \\
\notag
&& + c (\gamma^\prime) \cdot c (\gamma)  \cdot c (\delta)  \cdot E ( 3)  \cdot \eps  \cdot  \int_0^t  \frac{1}{(1+\tau)} \cdot (1+\tau+|q|) \cdot \|\varpi(q)  \cdot  \derm  ( \Lie_{Z^J} A)_{V} (\tau,\cdot) \|_{L^\infty (\Sigma^{ext}_{\tau} )} d \tau \\
       \notag
    && + \sum_{|K| \leq |J|} \int_0^t (1+\tau) \cdot  \varpi(q) \cdot  \|    \Lie_{Z^K}  g^{\la\mu} \derm_{\la}   \derm_{\mu}  A_V (\tau,\cdot) \|_{L^\infty(\overline{D}_\tau)} d\tau \\
  \notag
  && + c (\delta) \cdot c (\gamma) \cdot C ( |J| ) \cdot E ( |J| + 2)  \cdot \eps^2 \\
  \notag
& & +  \int_0^t  C(q_0)   \cdot  c (\delta) \cdot c (\gamma) \cdot C(|J|) \cdot E ( |J| + 3)  \cdot \frac{\eps \cdot  \varpi(q) }{(1+t+|q|)^{1-      c (\gamma)  \cdot c (\delta)  \cdot c(|J|) \cdot E ( |J|+ 3) \cdot \eps } \cdot (1+|q|)^{2+\gamma - 2 \de }}  d\tau \\
&&+    \int_0^t    c (\delta) \cdot c (\gamma) \cdot E (  4 )  \cdot \frac{\eps }{(1+t+|q|) }     \cdot   (1+t+|q|) \cdot \varpi(q) \cdot   \sum_{X \in \cal V} |\derm ( \Lie_{Z^{K}} A )_{X} |  d\tau \; .
 \eeaa
 Hence, for $\eps \leq 1$, we get
\beaa
   \notag
&& (1+t+|q|) \cdot |\varpi(q) \cdot \derm ( \Lie_{Z^J} A)_V | \\
   \notag
   &\les&   c (\gamma^\prime) \cdot  c (\delta) \cdot c (\gamma) \cdot C ( |J|  ) \cdot E ( |J| + 4) \cdot \eps \, \\
\notag
&& + c (\gamma^\prime)    \cdot  c (\delta) \cdot c (\gamma) \cdot E (  4 )  \cdot \eps  \cdot  \int_0^t  \frac{1}{(1+\tau)} \cdot (1+\tau+|q|) \cdot  \sum_{X \in \cal V}  \|\varpi(q)  \cdot  \derm  ( \Lie_{Z^J} A)_{X} (\tau,\cdot) \|_{L^\infty (\Sigma^{ext}_{\tau} )} d \tau \\
       \notag
    && + \sum_{|K| \leq |J|} \int_0^t (1+\tau) \cdot  \varpi(q) \cdot  \|    \Lie_{Z^K}  g^{\la\mu} \derm_{\la}   \derm_{\mu}  A_V (\tau,\cdot) \|_{L^\infty(\overline{D}_\tau)} d\tau \\
  \notag
& & +  \int_0^t  C(q_0)   \cdot  c (\delta) \cdot c (\gamma) \cdot C(|J|) \cdot E ( |J| + 3)  \cdot \frac{\eps \cdot  \varpi(q) }{(1+t+|q|)^{1-      c (\gamma)  \cdot c (\delta)  \cdot c(|J|) \cdot E ( |J|+ 3) \cdot \eps } \cdot (1+|q|)^{2+\gamma - 2 \de }}  d\tau \; .
 \eeaa
 Similarly, we get the result for $h^1$, where we note that the last term has a different decaying factor in $|q|$ under the integral.
 
 \end{proof}

      \begin{lemma}\label{GronwallinequalitongradientofAandh1withliederivativesofsourceterms}

    Let $ 0\leq \delta \leq  \frac{1}{4} $, and $\ga > \de$, and $M \leq \eps \leq 1$. We assume the induction hypothesis holding true for both $A$ and on $h^1$, for all $|K|\leq |J| -1$. Let
\beaa
\varpi(q) :=\begin{cases}
(1+|q|)^{1+\gamma^\prime},\quad\text{when }\quad q>0\; , \\
     1 \,\quad\text{when }\quad   q<0 \;, \end{cases} 
\eeaa
and let
 \beaa
\cal V :=  \begin{cases}  \cal T \;  ,\quad\text{if }\quad V \in \cal T \; ,\\
   \cal U \; , \,\quad\text{if }\quad V \in \cal U \; . \end{cases}   
\eeaa
We have for $\ga^\prime =  \gamma - 2 \de$, the following estimate in the exterior for $A$, 
    
            \beaa
   \notag
&& (1+t+|q|) \cdot |\varpi(q) \cdot \derm ( \Lie_{Z^J} A)_V | \\
   \notag
   &\les&   C(q_0)  \cdot  c (\delta) \cdot c (\gamma) \cdot C ( |J|  ) \cdot E ( |J| + 4) \cdot \eps \cdot (1+t)^{     c (\gamma)  \cdot c (\delta)  \cdot c(|J|) \cdot E ( |J|+ 4) \cdot \eps }  \\
\notag
&& +  c (\delta) \cdot c (\gamma) \cdot E (  4 )  \cdot \eps  \cdot  \int_0^t  \frac{1}{(1+\tau)} \cdot (1+\tau+|q|) \cdot  \sum_{X \in \cal V}  \|\varpi(q)  \cdot  \derm  ( \Lie_{Z^J} A)_{X} (\tau,\cdot) \|_{L^\infty (\Sigma^{ext}_{\tau} )} d \tau \\
       \notag
    && + \sum_{|K| \leq |J|} \int_0^t (1+\tau) \cdot  \varpi(q) \cdot  \|    \Lie_{Z^K}  g^{\la\mu} \derm_{\la}   \derm_{\mu}  A_V (\tau,\cdot) \|_{L^\infty(\overline{D}_\tau)} d\tau \; .
 \eeaa
Also, we have for $\ga^\prime = 0$, the following estimate in the exterior for $h^1$, 
             \beaa
   \notag
&& (1+t+|q|) \cdot |\varpi(q) \cdot \derm ( \Lie_{Z^J} h^1)_{UV} | \\
   \notag
   &\les&   C(q_0)  \cdot  c (\delta) \cdot c (\gamma) \cdot C ( |J|  ) \cdot E ( |J| + 4) \cdot \eps \cdot (1+t)^{     c (\gamma)  \cdot c (\delta)  \cdot c(|J|) \cdot E ( |J|+ 4) \cdot \eps }  \\
\notag
&& + c (\delta) \cdot c (\gamma) \cdot E (  4 )  \cdot \eps  \cdot  \int_0^t  \frac{1}{(1+\tau)} \cdot (1+\tau+|q|) \cdot \|\varpi(q)  \cdot  \derm  ( \Lie_{Z^J} h^1) (\tau,\cdot) \|_{L^\infty (\Sigma^{ext}_{\tau} )} d \tau \\
       \notag
    && + \sum_{|K| \leq |J|} \int_0^t (1+\tau) \cdot  \varpi(q) \cdot  \|    \Lie_{Z^K}  g^{\la\mu} \derm_{\la}   \derm_{\mu}  h^1_{UV} (\tau,\cdot) \|_{L^\infty(\overline{D}_\tau)} d\tau \; .
 \eeaa
 \end{lemma}
 
 \begin{proof}
 
        We Choose $ 0\leq \delta \leq  \frac{1}{4} $, and $\gamma > \de $, and $M \leq \eps \leq 1$. We assume the induction hypothesis holding true for both $A$ and on $h^1$, for all $|K|\leq |J| -1$\,.
             
    Choosing for $A$, an $\varpi(q)$ such that $\gamma^\prime = \gamma - 2 \de <  \gamma -  \de$, for our choice of $\de$, we could then apply what we have already shown in Lemma \eqref{almostGrnwallinequalitongradientofAandh1afterestimatedthecommutator},
 and get
            \beaa
   \notag
&& (1+t+|q|) \cdot |\varpi(q) \cdot \derm ( \Lie_{Z^J} A)_V | \\
   \notag
   &\les&   c (\gamma^\prime) \cdot  c (\delta) \cdot c (\gamma) \cdot C ( |J|  ) \cdot E ( |J| + 4) \cdot \eps \, \\
\notag
&& + c (\gamma^\prime)    \cdot  c (\delta) \cdot c (\gamma) \cdot E (  4 )  \cdot \eps  \cdot  \int_0^t  \frac{1}{(1+\tau)} \cdot (1+\tau+|q|) \cdot  \sum_{X \in \cal V}  \|\varpi(q)  \cdot  \derm  ( \Lie_{Z^J} A)_{X} (\tau,\cdot) \|_{L^\infty (\Sigma^{ext}_{\tau} )} d \tau \\
       \notag
    && + \sum_{|K| \leq |J|} \int_0^t (1+\tau) \cdot  \varpi(q) \cdot  \|    \Lie_{Z^K}  g^{\la\mu} \derm_{\la}   \derm_{\mu}  A_V (\tau,\cdot) \|_{L^\infty(\overline{D}_\tau)} d\tau \\
  \notag
& & +  \int_0^t  C(q_0)   \cdot  c (\delta) \cdot c (\gamma) \cdot C(|J|) \cdot E ( |J| + 3)  \cdot \frac{\eps \cdot  \varpi(q) }{(1+t+|q|)^{1-      c (\gamma)  \cdot c (\delta)  \cdot c(|J|) \cdot E ( |J|+ 3) \cdot \eps } \cdot (1+|q|)^{2+\gamma - 2 \de }}  d\tau \\
   &\les&   c (\gamma^\prime) \cdot  c (\delta) \cdot c (\gamma) \cdot C ( |J|  ) \cdot E ( |J| + 4) \cdot \eps \, \\
\notag
&& +   c (\delta) \cdot c (\gamma) \cdot E (  4 )  \cdot \eps  \cdot  \int_0^t  \frac{1}{(1+\tau)} \cdot (1+\tau+|q|) \cdot  \sum_{X \in \cal V}  \|\varpi(q)  \cdot  \derm  ( \Lie_{Z^J} A)_{X} (\tau,\cdot) \|_{L^\infty (\Sigma^{ext}_{\tau} )} d \tau \\
       \notag
    && + \sum_{|K| \leq |J|} \int_0^t (1+\tau) \cdot  \varpi(q) \cdot  \|    \Lie_{Z^K}  g^{\la\mu} \derm_{\la}   \derm_{\mu}  A_V (\tau,\cdot) \|_{L^\infty(\overline{D}_\tau)} d\tau \\
  \notag
& & +  \int_0^t  C(q_0)   \cdot  c (\delta) \cdot c (\gamma) \cdot C(|J|) \cdot E ( |J| + 3)  \cdot \frac{\eps  }{(1+t)^{1-      c (\gamma)  \cdot c (\delta)  \cdot c(|J|) \cdot E ( |J|+ 3) \cdot \eps } \cdot (1+|q|)}  d\tau \\
   &\les&   c (\gamma^\prime) \cdot  c (\delta) \cdot c (\gamma) \cdot C ( |J|  ) \cdot E ( |J| + 4) \cdot \eps \, \\
\notag
&& + c (\delta) \cdot c (\gamma) \cdot E (  4 )  \cdot \eps  \cdot  \int_0^t  \frac{1}{(1+\tau)} \cdot (1+\tau+|q|) \cdot  \sum_{X \in \cal V}  \|\varpi(q)  \cdot  \derm  ( \Lie_{Z^J} A)_{X} (\tau,\cdot) \|_{L^\infty (\Sigma^{ext}_{\tau} )} d \tau \\
       \notag
    && + \sum_{|K| \leq |J|} \int_0^t (1+\tau) \cdot  \varpi(q) \cdot  \|    \Lie_{Z^K}  g^{\la\mu} \derm_{\la}   \derm_{\mu}  A_V (\tau,\cdot) \|_{L^\infty(\overline{D}_\tau)} d\tau \\
  \notag
& & + C(q_0)   \cdot  c (\delta) \cdot c (\gamma) \cdot C(|J|) \cdot E ( |J| + 3)  \cdot \eps \cdot (1+t)^{     c (\gamma)  \cdot c (\delta)  \cdot c(|J|) \cdot E ( |J|+ 3) \cdot \eps }   \; ,
 \eeaa
 Given that $E ( |J|+ 3) \leq E ( |J|+ 4)$, we get the result. 

Now,  choosing for $h^1$, an $\varpi(q)$ such that $\gamma^\prime = 0 <  \gamma -  \de$, for our choice of $\ga > \de$, we could then apply what we have already shown in Lemma \eqref{almostGrnwallinequalitongradientofAandh1afterestimatedthecommutator} for $h^1$ and get similarly the following result,

 \beaa
   \notag
&& (1+t+|q|) \cdot |\varpi(q) \cdot \derm ( \Lie_{Z^J} h^1)_V | \\
   \notag
   &\les&   c (\gamma^\prime) \cdot  c (\delta) \cdot c (\gamma) \cdot C ( |J|  ) \cdot E ( |J| + 4) \cdot \eps \, \\
\notag
&& + c (\gamma^\prime)   \cdot  c (\delta) \cdot c (\gamma) \cdot E (  4 )  \cdot \eps  \cdot  \int_0^t  \frac{1}{(1+\tau)} \cdot (1+\tau+|q|) \cdot \|\varpi(q)  \cdot  \derm  ( \Lie_{Z^J} h^1) (\tau,\cdot) \|_{L^\infty (\Sigma^{ext}_{\tau} )} d \tau \\
       \notag
    && + \sum_{|K| \leq |J|} \int_0^t (1+\tau) \cdot  \varpi(q) \cdot  \|    \Lie_{Z^K}  g^{\la\mu} \derm_{\la}   \derm_{\mu}  h^1_{UV} (\tau,\cdot) \|_{L^\infty(\overline{D}_\tau)} d\tau \\
  \notag
& & +  \int_0^t  C(q_0)   \cdot  c (\delta) \cdot c (\gamma) \cdot C(|J|) \cdot E ( |J| + 3)  \cdot \frac{\eps \cdot  \varpi(q) }{(1+t+|q|)^{1-      c (\gamma)  \cdot c (\delta)  \cdot c(|J|) \cdot E ( |J|+ 3) \cdot \eps } \cdot (1+|q|)^2}  d\tau \; .
 \eeaa
 
  \end{proof}

  \subsection{Upgrading the induction hypothesis on the good components of the Yang-Mills potential}\
  
     \begin{lemma}

     Let $ 0\leq \delta \leq  \frac{1}{4} $, and $\ga > \de$, and $M \leq \eps \leq 1$. We assume the induction hypothesis holding true for both $A$ and on $h^1$, for all $|K|\leq |J| -1$. 
   
  Then, in the exterior region $\overline{C}$, we have the following estimate on the tangential components of the Einstein-Yang-Mills potential,
      \beaa
   \notag
&& | \derm ( \Lie_{Z^J} A)_{\cal T} | \\
   \notag
        &\les&   C(q_0)  \cdot  c (\delta) \cdot c (\gamma) \cdot C ( |J|  ) \cdot E ( |J| + 4) \cdot \eps \cdot \frac{1}{ (1+t)^{  1-   c (\gamma)  \cdot c (\delta)  \cdot c(|J|) \cdot E ( |J|+ 4) \cdot \eps }\cdot  (1+|q|)^{1+\ga - 2 \de}  }  \; .
 \eeaa
   
   \end{lemma}
   
   \begin{proof}
   
   Let $\varpi(q)$ defined as in \eqref{definitionofnewweighttoupgrade}, that is 
\beaa
\varpi(q) :=\begin{cases}
(1+|q|)^{1+\gamma^\prime},\quad\text{when }\quad q>0 \; ,\\
     1 \,\quad\text{when }\quad q<0 \; .\end{cases}
\eeaa

Based on Lemma \ref{GronwallinequalitongradientofAandh1withliederivativesofsourceterms}, we have for $\ga^\prime =  \gamma - 2 \de$\,, the following estimate in the exterior for $A$\,, 
    
            \beaa
   \notag
&& (1+t+|q|) \cdot \sum_{X \in \cal T} |\varpi(q) \cdot \derm ( \Lie_{Z^J} A)_{X}  | \\
   \notag
   &\les&   C(q_0)  \cdot  c (\delta) \cdot c (\gamma) \cdot C ( |J|  ) \cdot E ( |J| + 4) \cdot \eps \cdot (1+t)^{     c (\gamma)  \cdot c (\delta)  \cdot c(|J|) \cdot E ( |J|+ 4) \cdot \eps }  \\
\notag
&& +   c (\delta) \cdot c (\gamma) \cdot E (  4 )  \cdot \eps  \cdot  \int_0^t  \frac{1}{(1+\tau)} \cdot (1+\tau+|q|) \cdot  \sum_{X \in \cal T}  \|\varpi(q)  \cdot  \derm  ( \Lie_{Z^J} A)_{X}  (\tau,\cdot) \|_{L^\infty (\Sigma^{ext}_{\tau} )} d \tau \\
       \notag
    && + \sum_{|K| \leq |J|} \int_0^t (1+\tau) \cdot  \varpi(q) \cdot  \sum_{X \in \cal T}   \|    \Lie_{Z^K}  g^{\la\mu} \derm_{\la}   \derm_{\mu}  A_{X}  (\tau,\cdot) \|_{L^\infty(\overline{D}_\tau)} d\tau \; .
 \eeaa

       We have from Lemma \ref{decayestimateonthesourcetermforgoodpoentialAcomponents}, that for $\ga > \de $\,, and $0 \leq \de \leq 1$\,, 

                                        \beaa
 | \Lie_{Z^I} g^{\la\mu} \derm_{\la}   \derm_{\mu}   A_{{\cal T}}   |   &\leq& \begin{cases} c (\delta) \cdot c (\gamma) \cdot C ( |I| ) \cdot E ( |I| + 3)   \cdot \frac{\eps }{(1+t+|q|)^{3-3\delta} (1+|q|)^{1+\ga}},\quad\text{when }\quad q>0\; ,\\
       C ( |I| ) \cdot E ( |I| + 3)  \cdot \frac{\eps   \cdot (1+| q |   )^{\frac{3}{2} }  }{(1+t+|q|)^{3-3\delta} }  \,\quad\text{when }\quad q<0 \; . \end{cases} \\
           \notag
       \eeaa

    Thus, in the exterior region $\overline{C}$, we have for all $|K | \leq |I|$,
                                                \beaa
                                                \notag
&& (1+ t) \cdot  \varpi(q) \cdot  \| \, \Lie_{Z^K} g^{\la\mu} \derm_{\la}   \derm_{\mu}   A_{{\cal T}}  (t,\cdot)  \|_{L^\infty(\overline{D}_t)} \\
\notag
& \leq& C(q_0) \cdot c (\delta) \cdot c (\gamma) \cdot C ( |I| ) \cdot E ( |I| + 3)   \cdot \frac{\eps \cdot \varpi(q) }{(1+t+|q|)^{2-3\delta} (1+|q|)^{1+\ga}}  \\
& \leq& C(q_0) \cdot c (\delta) \cdot c (\gamma) \cdot C ( |I| ) \cdot E ( |I| + 3)   \cdot \frac{\eps \cdot (1+|q|)^{1+\ga^{\prime}} }{(1+t+|q|)^{2-3\delta} (1+|q|)^{1+\ga}}  \\
& \leq& C(q_0) \cdot c (\delta) \cdot c (\gamma) \cdot C ( |I| ) \cdot E ( |I| + 3)   \cdot \frac{\eps \cdot (1+|q|)^{1+\ga - 2 \de} }{(1+t+|q|)^{2-3\delta} (1+|q|)^{1+\ga}}  \\
& \leq& C(q_0) \cdot c (\delta) \cdot c (\gamma) \cdot C ( |I| ) \cdot E ( |I| + 3)   \cdot \frac{\eps  }{(1+t+|q|)^{2-3\delta} (1+|q|)^{1+2\de}}  \\
& \leq& C(q_0) \cdot c (\delta) \cdot c (\gamma) \cdot C ( |I| ) \cdot E ( |I| + 3)   \cdot \frac{\eps  }{(1+t)^{2-3\delta} }  \; .\\
       \eeaa
       Hence, for the tangential components ${\cal T}$, we obtain
       \beaa
     &&  \int_0^t (1+\tau) \cdot  \varpi(q) \cdot  \|    \Lie_{Z^K}  g^{\la\mu} \derm_{\la}   \derm_{\mu}  A_{\cal T} (\tau,\cdot) \|_{L^\infty(\overline{D}_\tau)} d\tau \\
   & \leq&      \int_0^t   C(q_0) \cdot c (\delta) \cdot c (\gamma) \cdot C ( |I| ) \cdot E ( |I| + 3)   \cdot \frac{\eps  }{(1+\tau )^{2-3\delta} }d\tau \\
      & \leq&   C(q_0) \cdot c (\delta) \cdot c (\gamma) \cdot C ( |I| ) \cdot E ( |I| + 3)   \cdot \eps \cdot \Big[ \frac{1  }{(1+\tau )^{1-3\delta} } \Big]^{\infty}_{0} \\
        & \leq&   C(q_0) \cdot c (\delta) \cdot c (\gamma) \cdot C ( |I| ) \cdot E ( |I| + 3)   \cdot \eps \\
        && \text{(since $\de  \leq  \frac{1}{4} <  \frac{1}{3}$). }
 \eeaa

 Thus,
       \beaa
   \notag
&& (1+t+|q|) \cdot \sum_{X \in \cal T} |\varpi(q) \cdot \derm ( \Lie_{Z^J} A)_{X} | \\
   \notag
   &\les&   C(q_0)  \cdot  c (\delta) \cdot c (\gamma) \cdot C ( |J|  ) \cdot E ( |J| + 4) \cdot \eps \cdot (1+t)^{     c (\gamma)  \cdot c (\delta)  \cdot c(|J|) \cdot E ( |J|+ 4) \cdot \eps }  \\
\notag
&& +  c (\delta) \cdot c (\gamma) \cdot E (  4 )  \cdot \eps  \cdot  \int_0^t  \frac{1}{(1+\tau)} \cdot (1+\tau+|q|) \cdot \sum_{X \in \cal T} \|\varpi(q)  \cdot  \derm  ( \Lie_{Z^J} A)_{X} (\tau,\cdot) \|_{L^\infty (\Sigma^{ext}_{\tau} )} d \tau \\
       \notag
    && + C(q_0) \cdot c (\delta) \cdot c (\gamma) \cdot C ( |I| ) \cdot E ( |I| + 3)   \cdot \eps  \\
     &\les&   C(q_0)  \cdot  c (\delta) \cdot c (\gamma) \cdot C ( |J|  ) \cdot E ( |J| + 4) \cdot \eps \cdot (1+t)^{     c (\gamma)  \cdot c (\delta)  \cdot c(|J|) \cdot E ( |J|+ 4) \cdot \eps }  \\
\notag
&& + c (\delta) \cdot c (\gamma) \cdot E (  4 )  \cdot \eps  \cdot  \int_0^t  \frac{1}{(1+\tau)} \cdot (1+\tau+|q|) \cdot \sum_{X \in \cal T} \|\varpi(q)  \cdot  \derm  ( \Lie_{Z^J} A)_{X} (\tau,\cdot) \|_{L^\infty (\Sigma^{ext}_{\tau} )} d \tau \; .\\
 \eeaa
 Using Grönwall lemma, we get
      \beaa
   \notag
&& (1+t+|q|) \cdot\sum_{X \in \cal T}  |\varpi(q) \cdot \derm ( \Lie_{Z^J} A)_{X} | \\
   \notag
     &\les&   C(q_0)  \cdot  c (\delta) \cdot c (\gamma) \cdot C ( |J|  ) \cdot E ( |J| + 4) \cdot \eps \cdot (1+t)^{     c (\gamma)  \cdot c (\delta)  \cdot c(|J|) \cdot E ( |J|+ 4) \cdot \eps }  \cdot (1+t)^{ c (\delta) \cdot c (\gamma) \cdot E (  4 )}   \\
      &\les&   C(q_0)  \cdot  c (\delta) \cdot c (\gamma) \cdot C ( |J|  ) \cdot E ( |J| + 4) \cdot \eps \cdot (1+t)^{     c (\gamma)  \cdot c (\delta)  \cdot c(|J|) \cdot E ( |J|+ 4) \cdot \eps }   \; , \\
 \eeaa
 which gives that in the exterior,
      \beaa
   \notag
&& | \derm ( \Lie_{Z^J} A)_{\cal T} | \\
   \notag
        &\les&   C(q_0)  \cdot  c (\delta) \cdot c (\gamma) \cdot C ( |J|  ) \cdot E ( |J| + 4) \cdot \eps \cdot \frac{1}{ (1+t)^{  1-   c (\gamma)  \cdot c (\delta)  \cdot c(|J|) \cdot E ( |J|+ 4) \cdot \eps }\cdot  (1+|q|)^{1+\ga - 2 \de}  }  \; .
 \eeaa
 
 \end{proof}

\subsection{The terms $\Lie_{Z^I}  \big(  A_L  \cdot    \derm ( \Lie_{Z^I} A  )  \big) $  and $\Lie_{Z^I}    \big(   A_{e_{a}} \cdot    \derm ( \Lie_{Z^I} A_{e_{a}} ) \big) \; $}\

\begin{lemma}\label{upgradeddecayestimateonthetermALgradientA}
     Let $ 0\leq \delta \leq  \frac{1}{4} $, and $\ga > \de$, and $M \leq \eps \leq 1$. We assume the induction hypothesis holding true for both $A$ and on $h^1$, for all $|K|\leq |J| -1$. 
   
  Then, in the exterior region $\overline{C}$, we have the following estimate,
        \beaa
 && \sum_{|K| + |I| \leq |J| }  | \Lie_{Z^K} A_L  | \cdot    | \derm  ( \Lie_{Z^I}  A  ) |  \\
    &\les & C(q_0)   \cdot  c (\delta) \cdot c (\gamma) \cdot C(|J|) \cdot E ( |J| + 3)  \cdot \frac{\eps }{(1+t+|q|)^{2- c (\gamma)  \cdot c (\delta)  \cdot c(|J|) \cdot E ( |J|+ 3) \cdot \eps } \cdot (1+|q|)^{1+2\gamma -4\de }} \; .
  \eeaa

\end{lemma}

\begin{proof}

We now decompose the sum in the following term,
\beaa
&& \sum_{|K| + |I| \leq |J| }  | \Lie_{Z^K} A_L  | \cdot    | \derm  ( \Lie_{Z^I}  A  ) | \\
&\leq &  \sum_{|I| = |J| }  | A_L  | \cdot    | \derm  ( \Lie_{Z^I}  A  ) | +   \sum_{|K| + |I| \leq |J|\;, |I| < |J| }  | \Lie_{Z^K} A_L  |  \cdot    | \derm  ( \Lie_{Z^I}  A  ) | \; .
\eeaa

\textbf{The term $  \sum_{|I| = |J| }  | A_L  | \cdot    | \derm  ( \Lie_{Z^I}  A  ) |$:}\\

We have shown that
                   \bea
        \notag
 |     A_{L}  |   &\les&  \begin{cases} c (\delta) \cdot c (\gamma) \cdot E (  3)  \cdot \frac{ \eps   }{ (1+t+|q|)^{2-2\delta} \cdot  (1+|q|)^{\ga - 1}  }  ,\quad\text{when }\quad q>0 \; ,\\
         E (  3) \cdot  \frac{ \eps \cdot (1+|q|)^{\frac{3}{2}  } }{ (1+t+|q|)^{2-2\delta} }  \,\quad\text{when }\quad q<0 \; , \end{cases} \\
       \eea      
and for all $|I| \leq |J|$, under the bootstrap assumption, we have
                        \beaa
 \notag
|\derm  ( \Lie_{Z^I} A )      |    &\leq& \begin{cases}   C ( |J| )  \cdot E ( |J| + 2)  \cdot \frac{\eps }{(1+t+|q|)^{1-\delta} (1+|q|)^{1+\ga}},\quad\text{when }\quad q>0 \; ,\\
       C ( |J| ) \cdot    E (  |J| + 2)  \cdot \frac{\eps  }{(1+t+|q|)^{1-\delta}(1+|q|)^{\frac{1}{2} }}  \,\quad\text{when }\quad q<0 \; . \end{cases} \\
      \eeaa
Consequently, for all  $|I| \leq |J|$,
     \beaa
| A_L  | \cdot    | \derm   ( \Lie_{Z^I}  A )    | \les \begin{cases}   C ( |J| )  \cdot E ( |J| + 2)  \cdot  c (\delta) \cdot c (\gamma) \cdot E (  3)  \cdot \frac{\eps }{(1+t+|q|)^{3-3\delta} (1+|q|)^{2\ga}},\quad\text{when }\quad q>0\; ,\\
        C ( |J| )  \cdot E ( |J| + 2)  \cdot   E ( 3)  \cdot \frac{\eps  \cdot (1+|q|) }{(1+t+|q|)^{3-3\delta}}  \,\quad\text{when }\quad q<0 \;  . \end{cases} \\
      \eeaa
Thus, in the exterior $\overline{C}$, we have
         \beaa
 \sum_{|I| = |J| }  | A_L  | \cdot    | \derm ( \Lie_{Z^I} A )   | &\les&  C(q_0) \cdot  C ( |J| )  \cdot E ( |J| + 2)  \cdot  c (\delta) \cdot c (\gamma) \cdot E (  3)   \cdot \frac{\eps }{(1+t+|q|)^{3-3\delta} (1+|q|)^{2\ga}} \; .\\
      \eeaa

\textbf{The term $   \sum_{|K| + |I| \leq |J|,\; |I| < |J| }  | \Lie_{Z^K} A_L  |  \cdot    | \derm  ( \Lie_{Z^I}  A  ) | $:}\\

We have shown in Lemma \ref{estimategoodcomponentspotentialandmetric}, that for all $|K| \leq |J|$,
              \bea
        \notag
 |  \Lie_{Z^K}  A_{L}  |   &\leq& \int\limits_{s,\Om=const} \sum_{|M|\leq |J| -1} |   \derm (  \Lie_{Z^M} A ) | \\
      \notag
     && + \begin{cases} c (\delta) \cdot c (\gamma) \cdot C ( |J| ) \cdot E ( |J| + 3)  \cdot \big( \frac{ \eps   }{ (1+t+|q|)^{2-2\delta} \cdot  (1+|q|)^{\ga - 1} } \big),\quad\text{when }\quad q>0 \; ,\\
       C ( |J| ) \cdot E ( |J| + 3)  \cdot \big( \frac{ \eps \cdot (1+|q|)^{\frac{3}{2} }   }{ (1+t+|q|)^{2-2\delta} } \big) \,\quad\text{when }\quad q<0 \; . \end{cases} 
       \eea
       Hence,
              \beaa
     &&  \sum_{|K| + |I| \leq |J|,\; |I| < |J| }  | \Lie_{Z^K} A_L  | \cdot    |  \derm ( \Lie_{Z^I} A  ) | \\
       &\leq& \sum_{|K| + |I| \leq |J|,\; |I| < |J| }  \Big(  \int\limits_{s,\Om=const} \sum_{|M|\leq |J| -1} |   \derm  \Lie_{Z^M} A |  \Big)  \cdot  |  \derm ( \Lie_{Z^I} A  ) |  \\
      \notag
     && +    \sum_{|K| + |I| \leq |J|,\; |I| < |J| }  \Big( \begin{cases} c (\delta) \cdot c (\gamma) \cdot C ( |J| ) \cdot E ( |J| + 3)  \cdot \big( \frac{ \eps   }{ (1+t+|q|)^{2-2\delta} \cdot  (1+|q|)^{\ga - 1} } \big),\quad\text{when }\quad q>0 \; ,\\
       C ( |J| ) \cdot E ( |J| + 3)  \cdot \Big( \frac{ \eps \cdot (1+|q|)^{\frac{3}{2} }   }{ (1+t+|q|)^{2-2\delta} } \big) \,\quad\text{when }\quad q<0 \; , \end{cases}  \Big) \\
   && \quad \quad \quad  \quad \quad \quad \quad \quad   \cdot |  \derm ( \Lie_{Z^I} A  ) | \; .
   \eeaa
   Therefore, in the exterior, we get
          \bea \label{sumlessAL}
          \notag
     &&  \sum_{|K| + |I| \leq |J|,\; |I| < |J| }  | \Lie_{Z^K} A_L  | \cdot    |  \derm ( \Lie_{Z^I} A  ) | \\
     \notag
    &\leq& \sum_{|K| + |I| \leq |J|,\; |I| < |J| }  \Big(  \int\limits_{s,\Om=const} \sum_{|M|\leq |J| -1} |   \derm  \Lie_{Z^M} A |  \Big)  \cdot  |  \derm ( \Lie_{Z^I} A  ) |  \\
      \notag
     && +    \sum_{|K| + |I| \leq |J|,\; |I| < |J| } \Big(  C(q_0) \cdot  c (\delta) \cdot c (\gamma) \cdot C ( |J| ) \cdot E ( |J| + 3)  \cdot  \frac{ \eps   }{ (1+t+|q|)^{2-2\delta} \cdot  (1+|q|)^{\ga - 1} }  \Big) \\
 && \quad \quad \quad  \quad \quad \quad \quad \quad   \cdot |  \derm ( \Lie_{Z^I} A  ) | \; .
       \eea

However, under the induction hypothesis, we have
                   \bea
 \notag
&&\sum_{|M|\leq |J| -1} |   \derm  ( \Lie_{Z^M} A ) |    \\
\notag
&\leq&   C(q_0)   \cdot  c (\delta) \cdot c (\gamma) \cdot C(|J|) \cdot E ( |J| + 3)  \cdot \frac{\eps }{(1+t+|q|)^{1-      c (\gamma)  \cdot c (\delta)  \cdot c(|J|) \cdot E ( |J|+ 3) \cdot \eps } \cdot (1+|q|)^{1+\gamma - 2\de }}    \; ,\\
      \eea
   and consequently, by integrating, we obtain 
                         \bea
 \notag
&& \int\limits_{s,\Om=const} \sum_{|M|\leq |J| -1} |   \derm ( \Lie_{Z^M} A)  |    \\
\notag
&\leq&   C(q_0)   \cdot  c (\delta) \cdot c (\gamma) \cdot C(|J|) \cdot E ( |J| + 3)  \cdot \frac{\eps }{(1+t+|q|)^{1-      c (\gamma)  \cdot c (\delta)  \cdot c(|J|) \cdot E ( |J|+ 3) \cdot \eps } \cdot (1+|q|)^{\gamma - 2\de }}    \; .\\
      \eea
      
On one hand, as a result, we get in the exterior, for all $|I| \leq |J| - 1$,

\beaa
&& \sum_{|K| + |I| \leq |J|,\; |I| < |J| }  \Big(  \int\limits_{s,\Om=const} \sum_{|M|\leq |J| -1} |   \derm  \Lie_{Z^M} A | \Big)  \cdot  |  \derm ( \Lie_{Z^I} A  ) |  \\
&\leq & C(q_0)   \cdot  c (\delta) \cdot c (\gamma) \cdot C(|J|) \cdot E ( |J| + 3)  \cdot \frac{\eps }{(1+t+|q|)^{2-      c (\gamma)  \cdot c (\delta)  \cdot c(|I|) \cdot E ( |I|+ 3) \cdot \eps } \cdot (1+|q|)^{1+2\gamma - 4\de }} \; .
\eeaa

On the other hand, from the bootstrap assumption, we have in the exterior region, for all $|I| \leq |J|$,       
     
                              \beaa
 \notag
      |\derm  ( \Lie_{Z^I} A )   |     &\leq& C(q_0) \cdot C ( |J| ) \cdot E ( |J| + 2)  \cdot \frac{\eps }{(1+t+|q|)^{1-\delta} \cdot (1+|q|)^{1+\ga}} \; ,
 \eeaa
            
      and as a result,

      \beaa
            \notag
     &&     \sum_{|K| + |I| \leq |J|,\; |I| < |J| } \Big(  C(q_0) \cdot  c (\delta) \cdot c (\gamma) \cdot C ( |J| ) \cdot E ( |J| + 3)  \cdot  \frac{ \eps   }{ (1+t+|q|)^{2-2\delta} \cdot  (1+|q|)^{\ga - 1} }  \Big) \\
 && \quad \quad \quad  \quad \quad \quad \quad \quad   \cdot |  \derm ( \Lie_{Z^I} A  ) | \\
&\leq&   C(q_0)   \cdot  c (\delta) \cdot c (\gamma) \cdot C(|J|) \cdot E ( |J| + 3)  \cdot \frac{\eps }{(1+t+|q|)^{3-3\de } \cdot (1+|q|)^{2\gamma }}    \; ,\\
      \eeaa
      
      Putting together in \eqref{sumlessAL}, we obtain,
      
                \beaa
     &&  \sum_{|K| + |I| \leq |J|,\; |I| < |J| }  | \Lie_{Z^K} A_L  | \cdot    |  \derm ( \Lie_{Z^I} A  ) | \\
  &\leq & C(q_0)   \cdot  c (\delta) \cdot c (\gamma) \cdot C(|J|) \cdot E ( |J| + 3)  \cdot \frac{\eps }{(1+t+|q|)^{2-      c (\gamma)  \cdot c (\delta)  \cdot c(|I|) \cdot E ( |I|+ 3) \cdot \eps } \cdot (1+|q|)^{1+2\gamma - 4\de }} \\
  && + C(q_0)   \cdot  c (\delta) \cdot c (\gamma) \cdot C(|J|) \cdot E ( |J| + 3)  \cdot \frac{\eps }{(1+t+|q|)^{3-3\de } \cdot (1+|q|)^{2\gamma }} \\
    &\leq & C(q_0)   \cdot  c (\delta) \cdot c (\gamma) \cdot C(|J|) \cdot E ( |J| + 3)  \cdot \frac{\eps }{(1+t+|q|)^{2-      c (\gamma)  \cdot c (\delta)  \cdot c(|I|) \cdot E ( |I|+ 3) \cdot \eps } \cdot (1+|q|)^{1+2\gamma - 4\de }} \\
  && + C(q_0)   \cdot  c (\delta) \cdot c (\gamma) \cdot C(|J|) \cdot E ( |J| + 3)  \cdot \frac{\eps }{(1+t+|q|)^{2 } \cdot (1+|q|)^{1+2\gamma -3\de }} \\
    &\leq & C(q_0)   \cdot  c (\delta) \cdot c (\gamma) \cdot C(|J|) \cdot E ( |J| + 3)  \cdot \frac{\eps }{(1+t+|q|)^{2- c (\gamma)  \cdot c (\delta)  \cdot c(|J|) \cdot E ( |J|+ 3) \cdot \eps } \cdot (1+|q|)^{1+2\gamma -4\de}} \\
  &&  \text{(for $\de \leq \frac{1}{4} \leq 1 $)}.
  \eeaa

\textbf{Finally:}\\

We obtain
\beaa
&& \sum_{|K| + |I| \leq |J| }  | \Lie_{Z^K} A_L  | \cdot    | \derm  ( \Lie_{Z^I}  A  ) | \\
&\leq &  \sum_{|I| = |J| }  | A_L  | \cdot    | \derm  ( \Lie_{Z^I}  A  ) | +   \sum_{|K| + |I| \leq |J|\;, |I| < |J| }  | \Lie_{Z^K} A_L  |  \cdot    | \derm  ( \Lie_{Z^I}  A  ) | \\
 &\les&  C(q_0) \cdot  C ( |J| )  \cdot E ( |J| + 2)  \cdot  c (\delta) \cdot c (\gamma) \cdot E (  3)   \cdot \frac{\eps }{(1+t+|q|)^{3-3\delta} (1+|q|)^{2\ga }} \\
 && + C(q_0)   \cdot  c (\delta) \cdot c (\gamma) \cdot C(|J|) \cdot E ( |J| + 3)  \cdot \frac{\eps }{(1+t+|q|)^{2- c (\gamma)  \cdot c (\delta)  \cdot c(|J|) \cdot E ( |J|+ 3) \cdot \eps } \cdot (1+|q|)^{1+2\gamma -4\de}}  \; ,
\eeaa
and hence we get the result.

       \end{proof}
                  
\begin{lemma}\label{upgradeddecayestimateonthetermAe_agradAe_a}

   Let $ 0\leq \delta \leq  \frac{1}{4} $, and $\ga > \de$, and $M \leq \eps \leq 1$. We assume the induction hypothesis holding true for both $A$ and on $h^1$, for all $|K|\leq |J| -1$. Then, in the exterior region $\overline{C}$, we have the following estimate
     \beaa
  &&   \sum_{|K| + |I| \leq |J| }   | \Lie_{Z^K}  A_{e_a}   |  \cdot | \derm   ( \Lie_{Z^I} A_{e_a} )  |    \\
       &\leq&  C(q_0)  \cdot  c (\delta) \cdot c (\gamma) \cdot C ( |J|  ) \cdot E ( |J| + 4) \cdot     \frac{\eps}{ (1+t)^{  2-   c (\gamma)  \cdot c (\delta)  \cdot c(|J|) \cdot E ( |J|+ 4) \cdot \eps }\cdot  (1+|q|)^{1+2\ga - 4 \de}  }  \; .
\eeaa

\end{lemma}

\begin{proof}

   We have shown in  \eqref{rderiofAacomp} that
      \beaa
    \pa_r \Lie_{Z^I} A_{e_{a}} = \derm_r  ( \Lie_{Z^I}  A_{e_{a}} )
     \eeaa  
   and thus,
       \beaa
\notag
| \Lie_{Z^I} A_{e_{a}}  (t, | x | \cdot \Om) |   &\leq& | \Lie_{Z^I} A_{e_{a}}  \big(0, ( t + | x | ) \cdot \Om \big)|  + \int_{| x | }^{t + | x |  } | \derm_r  (\Lie_{Z^I} A_{e_a} (t + | x | - r,  r  \cdot \Om ) ) | dr \; .\\
\eeaa

Yet, we showed already that for $ 0\leq \delta \leq  \frac{1}{4} $, and $\ga > \de$, and $M \leq \eps \leq 1$, under the induction hypothesis holding true for both $A$ and on $h^1$ for all $|K|\leq |J| -1$, then we have the following upgraded estimate for the tangential components,
      \beaa
   \notag
&& | \derm ( \Lie_{Z^K} A)_{\cal T} | \\
   \notag
        &\les&   C(q_0)  \cdot  c (\delta) \cdot c (\gamma) \cdot C ( |J|  ) \cdot E ( |J| + 4) \cdot \eps \cdot \frac{1}{ (1+t)^{  1-   c (\gamma)  \cdot c (\delta)  \cdot c(|J|) \cdot E ( |J|+ 4) \cdot \eps }\cdot  (1+|q|)^{1+\ga - 2 \de}  }  \; .
 \eeaa
Thus, under these assumptions, we get 
     \beaa
     && | \derm_r  (\Lie_{Z^K}  A_{e_a} )   | \\
      &\les&   | \derm  ( \Lie_{Z^K}  A_{e_a} )   |  \leq  | \derm  ( \Lie_{Z^K}  A_{\cal T} )  |  \\
        &\les&   C(q_0)  \cdot  c (\delta) \cdot c (\gamma) \cdot C ( |J|  ) \cdot E ( |J| + 4) \cdot \eps \cdot \frac{1}{ (1+t)^{  1-   c (\gamma)  \cdot c (\delta)  \cdot c(|J|) \cdot E ( |J|+ 4) \cdot \eps }\cdot  (1+|q|)^{1+\ga - 2 \de}  }  \; .
 \eeaa
 
Hence, proceeding the integration, as earlier (see \cite{G4}), we get
     \beaa
      |  \Lie_{Z^K}  A_{e_a}   |              &\les&   C(q_0)  \cdot  c (\delta) \cdot c (\gamma) \cdot C ( |J|  ) \cdot E ( |J| + 4) \cdot \eps \cdot \frac{1}{ (1+t)^{  1-   c (\gamma)  \cdot c (\delta)  \cdot c(|J|) \cdot E ( |J|+ 4) \cdot \eps }\cdot  (1+|q|)^{\ga - 2 \de}  }  \; .
\eeaa 

As a result, we have
     \beaa
  &&   \sum_{|K| + |I| \leq |J| }   | \Lie_{Z^K}  A_{e_a}   |  \cdot | \derm   ( \Lie_{Z^I} A_{e_a} )  |    \\
       &\leq&  C(q_0)  \cdot  c (\delta) \cdot c (\gamma) \cdot C ( |J|  ) \cdot E ( |J| + 4) \cdot \eps \cdot    \frac{1}{ (1+t)^{  2-   c (\gamma)  \cdot c (\delta)  \cdot c(|J|) \cdot E ( |J|+ 4) \cdot \eps }\cdot  (1+|q|)^{1+2\ga - 4 \de}  }  \; .
\eeaa 

\end{proof}

      \begin{lemma} \label{SteptwoforinductonforA}
   
     Let $ 0 < \delta \leq  \frac{1}{4} $, and $\ga \geq 3 \de$, and $M \leq \eps \leq 1$. We assume the induction hypothesis holding true for both $A$ and on $h^1$, for all $|K|\leq |J| -1$. 
   
  Then, in the exterior region $\overline{C}$, we have the following estimate on all the components of the Einstein-Yang-Mills potential,
      \beaa
   \notag
&& | \derm ( \Lie_{Z^J} A) | \\
   \notag
        &\les&   C(q_0)  \cdot  c (\delta) \cdot c (\gamma) \cdot C ( |J|  ) \cdot E ( |J| + 4) \cdot \eps \cdot \frac{1}{ (1+t)^{  1-   c (\gamma)  \cdot c (\delta)  \cdot c(|J|) \cdot E ( |J|+ 4) \cdot \eps }\cdot  (1+|q|)^{1+\ga - 2 \de}  }  \; .
 \eeaa
    \end{lemma}
    
      \begin{proof}
      We have shown in Lemma \ref{decayestimateonthesourcetermforgoodpoentialAcomponents}, that for $\ga > \de $, and $0 \leq \de \leq 1$, in the exterior region $\overline{C}$, we have
                                                \bea
                                                \notag
&& \| \, \Lie_{Z^I} g^{\la\mu} \derm_{\la}   \derm_{\mu}   A_{{\cal T}}  (t,\cdot)  \|_{L^\infty(\overline{D}_t)} \\
\notag
& \leq& C(q_0) \cdot c (\delta) \cdot c (\gamma) \cdot C ( |I| ) \cdot E ( |I| + 3)   \cdot \frac{\eps }{(1+t+|q|)^{3-3\delta} (1+|q|)^{1+\ga}}  \; . \\
       \eea

However, to estimate $ \| \, \Lie_{Z^I} g^{\la\mu} \derm_{\la}   \derm_{\mu}   A_{{\underline{L}}} \|_{L^\infty(\overline{D}_t)}$, the only new terms are the following, which we already estimated in Lemmas \ref{upgradeddecayestimateonthetermAe_agradAe_a} and \eqref{upgradeddecayestimateonthetermALgradientA}, as follows,
             \beaa
 && \sum_{|K| + |I| \leq |J| } \Big(  | \Lie_{Z^K} A_L  | \cdot    | \derm  ( \Lie_{Z^I}  A  ) | +   | \Lie_{Z^K}  A_{e_a}   |  \cdot | \derm   ( \Lie_{Z^I} A_{e_a} )  |  \Big) \\
    &\les & C(q_0)   \cdot  c (\delta) \cdot c (\gamma) \cdot C(|J|) \cdot E ( |J| + 3)  \cdot \frac{\eps }{(1+t+|q|)^{2- c (\gamma)  \cdot c (\delta)  \cdot c(|J|) \cdot E ( |J|+ 3) \cdot \eps } \cdot (1+|q|)^{1+2\gamma - 4\de }} \\
       && + C(q_0)  \cdot  c (\delta) \cdot c (\gamma) \cdot C ( |J|  ) \cdot E ( |J| + 4) \cdot \eps \cdot    \frac{1}{ (1+t)^{  2-   c (\gamma)  \cdot c (\delta)  \cdot c(|J|) \cdot E ( |J|+ 4) \cdot \eps }\cdot  (1+|q|)^{1+2\ga - 4 \de}  }  \\
          &\les & C(q_0)  \cdot  c (\delta) \cdot c (\gamma) \cdot C ( |J|  ) \cdot E ( |J| + 4) \cdot  \frac{\eps }{(1+t+|q|)^{2- c (\gamma)  \cdot c (\delta)  \cdot c(|J|) \cdot E ( |J|+ 4) \cdot \eps } \cdot (1+|q|)^{1+2\gamma -4\de }} \\
       && \text{(using the fact that $E ( |J| + 4)  \geq E ( |J| + 3)$)}\; .
\eeaa 
Therefore,
\bea
   \notag
 && | \Lie_{Z^I } g^{\la\mu} \derm_{\la}   \derm_{\mu}   A_{{\underline{L}}}  |   \\
    \notag
&\les&   C(q_0) \cdot c (\delta) \cdot c (\gamma) \cdot C ( |I| ) \cdot E ( |I| + 3)   \cdot \frac{\eps }{(1+t+|q|)^{2} \cdot (1+|q|)^{2+\ga -3\de}} \\
   \notag
&& +    C(q_0)  \cdot  c (\delta) \cdot c (\gamma) \cdot C ( |I|  ) \cdot E ( |I| + 4) \cdot  \frac{\eps }{(1+t+|q|)^{2- c (\gamma)  \cdot c (\delta)  \cdot c(|J|) \cdot E ( |J|+ 4) \cdot \eps } \cdot (1+|q|)^{1+2\gamma -4\de }} \; . \\
\eea
   Let $\varpi(q)$ defined as in \eqref{definitionofnewweighttoupgrade}, then, we have  
    
                                           \beaa
                                                \notag
&& (1+ t) \cdot  \varpi(q) \cdot  \| \, \Lie_{Z^I} g^{\la\mu} \derm_{\la}   \derm_{\mu}   A_{{\underline{L}}}    \|_{L^\infty(\overline{D}_t)} \\
\notag
&\les&   C(q_0) \cdot c (\delta) \cdot c (\gamma) \cdot C ( |I| ) \cdot E ( |I| + 3)   \cdot \frac{\eps  \cdot (1+|q|)^{1+\ga^{\prime}} }{(1+t+|q|) \cdot (1+|q|)^{2+\ga -3\de}} \\
       && +  C(q_0)  \cdot  c (\delta) \cdot c (\gamma) \cdot C ( |I|  ) \cdot E ( |I| + 4) \cdot   \frac{\eps  \cdot (1+|q|)^{1+\ga^{\prime}}  }{ (1+t)^{  1-   c (\gamma)  \cdot c (\delta)  \cdot c(|I|) \cdot E ( |I|+ 4) \cdot \eps }\cdot  (1+|q|)^{1+2 ga - 4 \de  }  }  \\
       &\les&   C(q_0) \cdot c (\delta) \cdot c (\gamma) \cdot C ( |I| ) \cdot E ( |I| + 3)   \cdot \frac{\eps   }{(1+t+|q|) \cdot (1+|q|)^{1+\ga -3\de -\ga^\prime}} \\
       && +  C(q_0)  \cdot  c (\delta) \cdot c (\gamma) \cdot C ( |I|  ) \cdot E ( |I| + 4) \cdot   \frac{\eps   }{ (1+t)^{  1-   c (\gamma)  \cdot c (\delta)  \cdot c(|I|) \cdot E ( |I|+ 4) \cdot \eps }\cdot  (1+|q|)^{2\ga - 4 \de  - \ga^\prime }  }  \; .\\
\eeaa 

We choose $\ga^\prime =  \gamma - 2 \de$, with $\ga \geq 3\de$, then we have the following estimate in the exterior for $A$, for all $|I| \leq |J|$,
                                           \beaa
                                                \notag
&& (1+ t) \cdot  \varpi(q) \cdot  \| \, \Lie_{Z^I} g^{\la\mu} \derm_{\la}   \derm_{\mu}   A_{{\underline{L}}}    \|_{L^\infty(\overline{D}_t)} \\
\notag
       &\les&   C(q_0) \cdot c (\delta) \cdot c (\gamma) \cdot C ( |I| ) \cdot E ( |I| + 3)   \cdot \frac{\eps   }{(1+t+|q|) \cdot (1+|q|)^{1 - \de }} \\
       && +  C(q_0)  \cdot  c (\delta) \cdot c (\gamma) \cdot C ( |I|  ) \cdot E ( |I| + 4) \cdot   \frac{\eps   }{ (1+t)^{  1-   c (\gamma)  \cdot c (\delta)  \cdot c(|I|) \cdot E ( |I|+ 4) \cdot \eps }\cdot  (1+|q|)^{\ga - 2 \de  }  }  \\
          &\les&   C(q_0) \cdot c (\delta) \cdot c (\gamma) \cdot C ( |J| ) \cdot E ( |J| + 3)   \cdot \frac{\eps   }{(1+t)^{  1-   c (\gamma)  \cdot c (\delta)  \cdot c(|J|) \cdot E ( |J|+ 4) \cdot \eps } } \\
       && \text{(since $\de \leq 1$ and  $\ga \geq 3\de$ and $\de \geq 0$)}.
\eeaa

Now, for $0< \delta \leq  \frac{1}{4} $, we get $\ga \geq 3 \de > \de $. For $M \leq \eps \leq 1$ and under the induction hypothesis holding true for both $A$ and on $h^1$, for all $|K|\leq |J| -1$, we have showed in Lemma \ref{GronwallinequalitongradientofAandh1withliederivativesofsourceterms}, the following estimate in the exterior for $A$, for $\ga^\prime =  \gamma - 2 \de$,
            \beaa
   \notag
&& (1+t+|q|) \cdot |\varpi(q) \cdot \derm ( \Lie_{Z^J} A)_V | \\
   \notag
   &\les&   C(q_0)  \cdot  c (\delta) \cdot c (\gamma) \cdot C ( |J|  ) \cdot E ( |J| + 4) \cdot \eps \cdot (1+t)^{     c (\gamma)  \cdot c (\delta)  \cdot c(|J|) \cdot E ( |J|+ 4) \cdot \eps }  \\
\notag
&& +  c (\delta) \cdot c (\gamma) \cdot E (  4 )  \cdot \eps  \cdot  \int_0^t  \frac{1}{(1+\tau)} \cdot (1+\tau+|q|) \cdot \|\varpi(q)  \cdot  \derm  ( \Lie_{Z^J} A)_{V} (\tau,\cdot) \|_{L^\infty (\Sigma^{ext}_{\tau} )} d \tau \\
       \notag
    && + \sum_{|K| \leq |J|} \int_0^t (1+\tau) \cdot  \varpi(q) \cdot  \|    \Lie_{Z^K}  g^{\la\mu} \derm_{\la}   \derm_{\mu}  A_V (\tau,\cdot) \|_{L^\infty(\overline{D}_\tau)} d\tau \; .
 \eeaa

    Hence,
            \beaa
   \notag
&& (1+t+|q|) \cdot |\varpi(q) \cdot \derm ( \Lie_{Z^J} A)_{{\cal U}}  | \\
   \notag
   &\les&   C(q_0)  \cdot  c (\delta) \cdot c (\gamma) \cdot C ( |J|  ) \cdot E ( |J| + 4) \cdot \eps \cdot (1+t)^{     c (\gamma)  \cdot c (\delta)  \cdot c(|J|) \cdot E ( |J|+ 4) \cdot \eps }  \\
\notag
&& +   c (\delta) \cdot c (\gamma) \cdot E (  4 )  \cdot \eps  \cdot  \int_0^t  \frac{1}{(1+\tau)} \cdot (1+\tau+|q|) \cdot \|\varpi(q)  \cdot  \derm  ( \Lie_{Z^J} A)_{{\cal U}}  (\tau,\cdot) \|_{L^\infty (\Sigma^{ext}_{\tau} )} d \tau \\
       \notag
    && + \sum_{|K| \leq |J|} \int_0^t   C(q_0) \cdot c (\delta) \cdot c (\gamma) \cdot C ( |J| ) \cdot E ( |J| + 3)   \cdot \frac{\eps   }{(1+t)^{  1-   c (\gamma)  \cdot c (\delta)  \cdot c(|J|) \cdot E ( |J|+ 4) \cdot \eps } }  d\tau \\
     &\les&   C(q_0)  \cdot  c (\delta) \cdot c (\gamma) \cdot C ( |J|  ) \cdot E ( |J| + 4) \cdot \eps \cdot (1+t)^{     c (\gamma)  \cdot c (\delta)  \cdot c(|J|) \cdot E ( |J|+ 4) \cdot \eps }  \\
\notag
&& + c (\delta) \cdot c (\gamma) \cdot E (  4 )  \cdot \eps  \cdot  \int_0^t  \frac{1}{(1+\tau)} \cdot (1+\tau+|q|) \cdot \|\varpi(q)  \cdot  \derm  ( \Lie_{Z^J} A) (\tau,\cdot) \|_{L^\infty (\Sigma^{ext}_{\tau} )} d \tau \; . \\
 \eeaa
 Using Grönwall lemma, we get
      \beaa
   \notag
&& (1+t+|q|) \cdot |\varpi(q) \cdot \derm ( \Lie_{Z^J} A)_{\cal U} | \\
   \notag
     &\les&   C(q_0)  \cdot  c (\delta) \cdot c (\gamma) \cdot C ( |J|  ) \cdot E ( |J| + 4) \cdot \eps \cdot (1+t)^{     c (\gamma)  \cdot c (\delta)  \cdot c(|J|) \cdot E ( |J|+ 4) \cdot \eps }  \cdot (1+t)^{ c (\delta) \cdot c (\gamma) \cdot E (  4 )}   \\
      &\les&   C(q_0)  \cdot  c (\delta) \cdot c (\gamma) \cdot C ( |J|  ) \cdot E ( |J| + 4) \cdot \eps \cdot (1+t)^{     c (\gamma)  \cdot c (\delta)  \cdot c(|J|) \cdot E ( |J|+ 4) \cdot \eps }  \; ,  \\
 \eeaa
and thus, we obtain in the exterior,
      \beaa
   \notag
&& | \derm ( \Lie_{Z^J} A)_{\cal U} | \\
   \notag
        &\les&   C(q_0)  \cdot  c (\delta) \cdot c (\gamma) \cdot C ( |J|  ) \cdot E ( |J| + 4) \cdot \eps \cdot \frac{1}{ (1+t)^{  1-   c (\gamma)  \cdot c (\delta)  \cdot c(|J|) \cdot E ( |J|+ 4) \cdot \eps }\cdot  (1+|q|)^{1+\ga - 2 \de}  }  \; .
 \eeaa
 
 \end{proof}

\subsection{Upgrading the hypothesis assumption on the metric}\

            \begin{lemma}\label{decayestimateonsourcetermsforfullcompoofhusingonlyinductionhyposothattoupgradetoJ}
                   We have in the exterior region $\overline{C}$, for $\ga \geq 3\de $, and $0 < \delta \leq \frac{1}{4}$, under the induction hypothesis holding true on $A$ and $h$ for all $|K| \leq |J| -1 $, the following estimate
                      \beaa
\notag
 &&\sum_{|K| = |J|}  | \Lie_{Z^K}g^{\alpha\beta} \derm_\alpha \derm_\beta h_{{\underline{L}} {\underline{L}} } |    \\
 \notag
            &\les& C(q_0) \cdot   c (\gamma)  \cdot c (\delta)   \cdot E (4) \cdot  C(|J|) \cdot \frac{\eps }{ (1+t+|q|)} \cdot \sum_{|K| = |J|} |\derm ( \Lie_{Z^K} h_{ {\cal T} {\cal U}} ) | \\
       && +  C(q_0)  \cdot  c (\delta) \cdot c (\gamma) \cdot C ( |J|  ) \cdot E ( |J| + 4) \cdot  \frac{\eps}{ (1+t)^{  2-   c (\gamma)  \cdot c (\delta)  \cdot c(|J|) \cdot E ( |J|+ 4) \cdot \eps }\cdot  (1+|q|)^{2-\de}  } \; .
             \eeaa

             \end{lemma}
             \begin{proof}
             
Based on Lemma \ref{StructureoftheLiederivativesofthesourcetermsofthewaveoperatorforAandh}
and on Lemma \ref{decayestimateonsourcetermsforgoodcomponentofmetrichtauU}, in particular \eqref{aprioridecayestimateontauUcompoforfullmetrich} in the proof, we get that in the exterior, for $\ga > \de $\,, and $0 \leq \de \leq \frac{1}{4}$\,, 
          \beaa
\notag
 &&  | \Lie_{Z^J}g^{\alpha\beta} \derm_\alpha \derm_\beta h_{{\underline{L}} {\underline{L}} } |    \\
 \notag
 &\leq&  \begin{cases}  c (\delta) \cdot c (\gamma) \cdot E ( |J| + 3)  \cdot \frac{\eps }{(1+t+|q|)^{3-3\delta} (1+|q|)^{1+2\de}},\quad\text{when }\quad q>0  \;  ,\\
       \notag
      E ( |J| + 3)  \cdot \frac{\eps  \cdot (1+| q |   )^2 }{(1+t+|q|)^{3-3\delta} }  \,\quad\text{when }\quad q<0 \; \end{cases} \\
       \notag
       && + \sum_{|K| + |I| \leq |J|} \Big(  |\derm ( \Lie_{Z^K} h_{ {\cal T} {\cal U}} ) | \cdot |\derm ( \Lie_{Z^I} h_{ {\cal T} {\cal U}} ) | + |\derm (  \Lie_{Z^K} A_{{\cal T}} ) | \cdot |\derm (  \Lie_{Z^I} A_{{\cal T}} ) |  \Big)  \\
              \notag
            &\les&  C(q_0) \cdot  c (\delta) \cdot c (\gamma)  \cdot E ( |J| + 3)  \cdot \frac{\eps }{(1+t+|q|)^{3-3\delta} (1+|q|)^{1+2\de}}  \; \\
       \notag
         && + \sum_{|K| + |I| \leq |J|}  |\derm (  \Lie_{Z^K} A_{{\cal T}} ) | \cdot |\derm (  \Lie_{Z^I} A_{{\cal T}} ) |    \\
                \notag
             && + \sum_{|K| + |I| \leq |J|}   |\derm ( \Lie_{Z^K} h_{ {\cal T} {\cal U}} ) | \cdot |\derm ( \Lie_{Z^I} h_{ {\cal T} {\cal U}} ) |   \; .
       \eeaa
       
       Yet, we have already upgraded the estimate on $ |\derm (  \Lie_{Z^J} A_{{\cal T}} ) |$, for $\ga \geq 3\de > 0$ and thus, we have
       \beaa
         &&  \sum_{|K| + |I| \leq |J|}  |\derm (  \Lie_{Z^K} A_{{\cal T}} ) | \cdot |\derm (  \Lie_{Z^I} A_{{\cal T}} ) |    \\
         &\les&   C(q_0)  \cdot  c (\delta) \cdot c (\gamma) \cdot C ( |J|  ) \cdot E ( |J| + 4) \cdot \eps \cdot \frac{1}{ (1+t)^{  2-   c (\gamma)  \cdot c (\delta)  \cdot c(|J|) \cdot E ( |J|+ 4) \cdot \eps }\cdot  (1+|q|)^{2+2\ga - 4 \de}  } \; .
\eeaa
 Yet, for $\ga \geq 3\de$, we have $2\gamma - 4\de \geq 6\de - 4\de = 2\de \geq 0$. As a result,
           \beaa
\notag
 &&  | \Lie_{Z^J}g^{\alpha\beta} \derm_\alpha \derm_\beta h_{{\underline{L}} {\underline{L}} } |    \\
 \notag
            &\les&  C(q_0) \cdot  c (\delta) \cdot c (\gamma)  \cdot E ( |J| + 3)  \cdot \frac{\eps }{(1+t+|q|)^{2} (1+|q|)^{2 - \de}}  \; \\
       \notag
         && + C(q_0)  \cdot  c (\delta) \cdot c (\gamma) \cdot C ( |J|  ) \cdot E ( |J| + 4) \cdot  \frac{\eps}{ (1+t)^{  2-   c (\gamma)  \cdot c (\delta)  \cdot c(|J|) \cdot E ( |J|+ 4) \cdot \eps }\cdot  (1+|q|)^{2}  } \\
                         \notag
             && + \sum_{|K| + |I| \leq |J|}   |\derm ( \Lie_{Z^K} h_{ {\cal T} {\cal U}} ) | \cdot |\derm ( \Lie_{Z^I} h_{ {\cal T} {\cal U}} ) |   \\
    &\les&  C(q_0)  \cdot  c (\delta) \cdot c (\gamma) \cdot C ( |J|  ) \cdot E ( |J| + 4) \cdot  \frac{\eps}{ (1+t)^{  2-   c (\gamma)  \cdot c (\delta)  \cdot c(|J|) \cdot E ( |J|+ 4) \cdot \eps }\cdot  (1+|q|)^{2-\de}  } \\
                         \notag
             && + \sum_{|K| + |I| \leq |J|}   |\derm ( \Lie_{Z^K} h_{ {\cal T} {\cal U}} ) | \cdot |\derm ( \Lie_{Z^I} h_{ {\cal T} {\cal U}} ) |   \; .
       \eeaa

  However, we can decompose the following sum as
      \beaa
       && \sum_{|K| + |I| \leq |J|}   |\derm ( \Lie_{Z^K} h_{ {\cal T} {\cal U}} ) | \cdot |\derm ( \Lie_{Z^I} h_{ {\cal T} {\cal U}} ) | \\
       &=& \sum_{|K| = |J|}  2\cdot  |\derm  h_{ {\cal T} {\cal U}}  | \cdot |\derm ( \Lie_{Z^J} h_{ {\cal T} {\cal U}} ) |  +  \sum_{|K| + |I| \leq |J|-1}   |\derm ( \Lie_{Z^K} h_{ {\cal T} {\cal U}} ) | \cdot |\derm ( \Lie_{Z^I} h_{ {\cal T} {\cal U}} ) |    \; .
      \eeaa

\textbf{The terms $  \sum_{|K| + |I| \leq |J|-1}    |\derm ( \Lie_{Z^K} h_{ {\cal T} {\cal U}} ) | \cdot |\derm ( \Lie_{Z^I} h_{ {\cal T} {\cal U}} ) | 
   $}:

      Since we also took $\ga \geq 3 \de > 0$, we can use the induction hypothesis, and we get for all $|K| \leq |J| -1 $,      
            \bea
 \notag
&& |\derm  ( \Lie_{Z^K} h)   |  \\
  \notag
 &\leq & |\derm  ( \Lie_{Z^K} h^1)   | +  |\derm  ( \Lie_{Z^K} h^0)   |   \\
\notag
 &\leq&   C(q_0)   \cdot  c (\delta) \cdot c (\gamma) \cdot C(|J|) \cdot E (  |J|  +3)  \cdot \frac{\eps }{(1+t+|q|)^{1-   c (\gamma)  \cdot c (\delta)  \cdot c(|J|) \cdot E ( |J|+ 3)\cdot  \eps } \cdot (1+|q|)}     \\
  && + |\derm (  \Lie_{Z^K} h^0 )  | \; .
      \eea
   However, we have shown in Lemma \ref{Liederivativesofsphericalsymmetricpart}, that for $M\leq \eps$, we have for all $|K|\leq |J| -1 $,
   \beaa
 \notag
|\derm (  \Lie_{Z^K} h^0 )  |   &\leq&  C ( |J| )   \cdot \frac{\eps }{(1+t+|q|)^{2} }  \; .
      \eeaa   
      Thus, for all $|K| \leq |J| -1 $,      
            \beaa
 \notag
&& |\derm  ( \Lie_{Z^K} h)   |  \\
\notag
 &\leq&   C(q_0)   \cdot  c (\delta) \cdot c (\gamma) \cdot C(|J|) \cdot E (  |J|  +3)  \cdot \frac{\eps }{(1+t+|q|)^{1-   c (\gamma)  \cdot c (\delta)  \cdot c(|J|) \cdot E ( |J|+ 3)\cdot  \eps } \cdot (1+|q|)}     \\
  && + C ( |J| )   \cdot \frac{\eps }{(1+t+|q|)^{2} } \\
   &\leq&   C(q_0)   \cdot  c (\delta) \cdot c (\gamma) \cdot C(|J|) \cdot E (  |J|  +3)  \cdot \frac{\eps }{(1+t+|q|)^{1-   c (\gamma)  \cdot c (\delta)  \cdot c(|J|) \cdot E ( |J|+ 3)\cdot  \eps } \cdot (1+|q|)}      \; .
      \eeaa
Consequently,
     \beaa
     &&  \sum_{|K| + |I| \leq |J|-1}   |\derm ( \Lie_{Z^K} h_{ {\cal T} {\cal U}} ) | \cdot |\derm ( \Lie_{Z^I} h_{ {\cal T} {\cal U}} ) |  \\
     &\leq &   C(q_0)   \cdot  c (\delta) \cdot c (\gamma) \cdot C(|J|) \cdot E (  |J|  +3)  \cdot \frac{\eps }{(1+t+|q|)^{2-   c (\gamma)  \cdot c (\delta)  \cdot c(|J|) \cdot E ( |J|+ 3)\cdot  \eps } \cdot (1+|q|)^2}     \; .
     \eeaa

  \textbf{The term $   \sum_{|K| = |J|} |\derm  h_{ {\cal T} {\cal U}}  | \cdot |\derm ( \Lie_{Z^K} h_{ {\cal T} {\cal U}} ) |$}:

        We have shown in Lemma \ref{upgradedestimatesongoodcomponnentforh1andh0}, that in the exterior region $\overline{C}$\,, for $\ga > \de $\,, and $0 \leq \delta \leq \frac{1}{4}$\,, 
 \beaa
 \,|\derm h^1_{ {\cal T} {\cal U}} (t,x)|  &\les&  C(q_0) \cdot   c (\gamma)  \cdot c (\delta)   \cdot E (4) \cdot \frac{\eps }{ (1+t+|q|)} \; . 
\eeaa
Thus, for $M \leq \eps$,
\beaa
 \,|\derm h_{ {\cal T} {\cal U}} |  &\les&   \,|\derm h^1_{ {\cal T} {\cal U}} |   +  \,|\derm h^0_{ {\cal T} {\cal U}} |  \leq C(q_0) \cdot   c (\gamma)  \cdot c (\delta)   \cdot E (4) \cdot \frac{\eps }{ (1+t+|q|)} +    \frac{\eps }{(1+t+|q|)^{2} } \\
 &\les&   C(q_0) \cdot   c (\gamma)  \cdot c (\delta)   \cdot E (4) \cdot \frac{\eps }{ (1+t+|q|)} \; .
 \eeaa
Consequently,
\beaa
 &&  |\derm  h_{ {\cal T} {\cal U}}  | \cdot  \sum_{|K| = |J|} |\derm ( \Lie_{Z^K} h_{ {\cal T} {\cal U}} ) |   \\
   &\les&  C(q_0) \cdot   c (\gamma)  \cdot c (\delta)   \cdot E (4) \cdot  \frac{\eps }{ (1+t+|q|)} \cdot \sum_{|K| = |J|}  |\derm ( \Lie_{Z^K} h_{ {\cal T} {\cal U}} ) |  \; .
   \eeaa
   
   \textbf{The whole term}:
   
Therefore, we obtain
      \beaa
       && \sum_{|K| + |I| \leq |J|} |\derm ( \Lie_{Z^K} h_{ {\cal T} {\cal U}} ) | \cdot  |\derm ( \Lie_{Z^I} h_{ {\cal T} {\cal U}} ) | \\
       &=&   C(q_0) \cdot   c (\gamma)  \cdot c (\delta)   \cdot E (4) \cdot  \frac{\eps }{ (1+t+|q|)} \cdot   \sum_{|K| = |J|} |\derm ( \Lie_{Z^K} h_{ {\cal T} {\cal U}} ) | \\
       && + C(q_0)   \cdot  c (\delta) \cdot c (\gamma) \cdot C(|J|) \cdot E ( |J| + 3)  \cdot \frac{\eps }{(1+t+|q|)^{2-      c (\gamma)  \cdot c (\delta)  \cdot c(|J|) \cdot E ( |J|+ 3) \cdot \eps }\cdot (1+|q|)^2 } \; .
             \eeaa
             
 Finally, we get            
                       \beaa
\notag
 &&  | \Lie_{Z^J}g^{\alpha\beta} \derm_\alpha \derm_\beta h_{{\underline{L}} {\underline{L}} } |    \\
 \notag
            &\les& C(q_0)  \cdot  c (\delta) \cdot c (\gamma) \cdot C ( |J|  ) \cdot E ( |J| + 4) \cdot  \frac{\eps}{ (1+t)^{  2-   c (\gamma)  \cdot c (\delta)  \cdot c(|J|) \cdot E ( |J|+ 4) \cdot \eps }\cdot  (1+|q|)^{2-\de}  } \\
            && + C(q_0) \cdot   c (\gamma)  \cdot c (\delta)   \cdot E (4) \cdot  \frac{\eps }{ (1+t+|q|)} \cdot  \sum_{|K| = |J|} |\derm ( \Lie_{Z^K} h_{ {\cal T} {\cal U}} ) | \\
       && + C(q_0)   \cdot  c (\delta) \cdot c (\gamma) \cdot C(|J|) \cdot E ( |J| + 3)  \cdot \frac{\eps }{(1+t+|q|)^{2-      c (\gamma)  \cdot c (\delta)  \cdot c(|J|) \cdot E ( |J|+ 3) \cdot \eps } \cdot (1+|q|)^{2} } \; .
             \eeaa
       Summing over all $J$, with fixed length $|J|$, we get the desired result.
                    \end{proof}

\begin{lemma}\label{upgradedestimateforhwithaconstantCqo}

      Let $ 0 < \delta \leq  \frac{1}{4} $, and $\ga \geq 3 \de$, and $M \leq \eps \leq 1$. We assume the induction hypothesis holding true for both $A$ and on $h^1$, for all $|K|\leq |J| -1$. 
   
  Then, in the exterior region $\overline{C}$, we have the following estimate on all the components of the metric,
              \beaa
 \notag
&& |\derm  ( \Lie_{Z^K} h^1) (t,x)  |   \\
\notag
 &\leq&   C(q_0)   \cdot  c (\delta) \cdot c (\gamma) \cdot C(|J|) \cdot E (  |J|  +4)  \cdot \frac{\eps }{(1+t+|q|)^{1-    C(q_0)  \cdot    c (\gamma)  \cdot c (\delta)   \cdot E ( |J|+ 4) \cdot C ( |J|  )  \cdot  \eps } \cdot (1+|q|)}     \; . 
      \eeaa

\end{lemma}

       \begin{proof}
Based on what we have shown in Lemma \ref{decayestimateonsourcetermsforfullcompoofhusingonlyinductionhyposothattoupgradetoJ}, we get
          \beaa
\notag
&&(1+t+| q |  ) \cdot \varpi(q) \cdot \sum_{|K| \leq |J|}  | \Lie_{Z^K} g^{\alpha\beta} \derm_\alpha \derm_\beta h^1_{{\underline{L}} {\underline{L}} } |_{L^\infty(\overline{D}_\tau)} \\
      &\les& (1+t+| q |  ) \cdot \varpi(q) \cdot \sum_{|K| \leq |J|}  | \Lie_{Z^K} g^{\alpha\beta} \derm_\alpha \derm_\beta h^0_{{\underline{L}} {\underline{L}} } |_{L^\infty(\overline{D}_\tau)} \\ \\
      && + C(q_0) \cdot   c (\gamma)  \cdot c (\delta)   \cdot E (4)\cdot C(|J|)  \cdot  \frac{\eps }{ (1+t+|q|)} \cdot(1+t+| q |  )  \cdot \varpi(q) \cdot \sum_{|K| \leq |J|}  |\derm ( \Lie_{Z^K} h_{ {\cal T} {\cal U}} ) | \\
       && +  C(q_0)  \cdot  c (\delta) \cdot c (\gamma) \cdot C ( |J|  ) \cdot E ( |J| + 4) \cdot  \frac{\eps  \cdot (1+|q|)^{1+\gamma^\prime}}{ (1+t)^{  1-   c (\gamma)  \cdot c (\delta)  \cdot c(|J|) \cdot E ( |J|+ 4) \cdot \eps }\cdot  (1+|q|)^{2-\de}  } \; .
             \eeaa

Taking $\ga^\prime = 0$, based on Lemma \ref{estimateonthesourcetermsforhzerothesphericallsymmtrpart}, we have for $\de \leq 1$,
\beaa
&& (1+t+| q |  ) \cdot \varpi(q) \cdot  | \Lie_{ Z^J}  ( g^{\la\mu} \derm_{\la}   \derm_{\mu}    h^0 ) | \\
 &\les& C ( |J| ) \cdot c (\gamma) \cdot  E ( |J| + 2)  \cdot \frac{\eps  \cdot \varpi(q)}{(1+t+|q|)^{2} } \\
 &\les& C ( |J| ) \cdot c (\gamma) \cdot  E ( |J| + 2)  \cdot \frac{\eps }{(1+t) } \; .
\eeaa

Since $\de \leq 1$, we obtain

          \beaa
\notag
 &&  \sum_{|K| \leq |J|}  \int_0^t   (1+t+| q |  ) \cdot \varpi(q) \cdot | \Lie_{Z^K} g^{\alpha\beta} \derm_\alpha \derm_\beta h^1_{{\underline{L}} {\underline{L}} } |_{L^\infty(\overline{D}_\tau)}  d \tau  \\
      &\les&  \int_0^t C(q_0) \cdot   c (\gamma)  \cdot c (\delta)   \cdot E (4) \cdot C(|J|)  \cdot  \frac{\eps }{ (1+t+|q|)} \cdot(1+t+| q |  ) \cdot \varpi(q) \cdot \sum_{|K| \leq |J|}  |\derm ( \Lie_{Z^J} h_{ {\cal T} {\cal U}} ) | d \tau   \\
        && +   C ( |J| ) \cdot c (\gamma) \cdot  E ( |J| + 2)  \cdot \eps \cdot \ln(1+t) \\
       && +  C(q_0)   \cdot  c (\delta) \cdot c (\gamma) \cdot C(|J|) \cdot E ( |J| + 4)  \cdot \eps \cdot (1+t)^{    c (\gamma)  \cdot c (\delta)  \cdot c(|J|) \cdot E ( |J|+ 4) \cdot \eps } \\
        &\les&  \int_0^t C(q_0) \cdot   c (\gamma)  \cdot c (\delta)   \cdot E (4) \cdot C(|J|)  \cdot  \frac{\eps }{ (1+t+|q|)} \cdot(1+t+| q |  ) \cdot \varpi(q) \cdot \sum_{|K| \leq |J|}  |\derm ( \Lie_{Z^J} h_{ {\cal T} {\cal U}} ) | d \tau   \\
              && +  C(q_0)   \cdot  c (\delta) \cdot c (\gamma) \cdot C(|J|) \cdot E ( |J| + 4)  \cdot \eps   \cdot (1+t)^{    c (\gamma)  \cdot c (\delta)  \cdot c(|J|) \cdot E ( |J|+ 4) \cdot \eps } \; .
             \eeaa

However, we have shown in Lemma \ref{GronwallinequalitongradientofAandh1withliederivativesofsourceterms}, that for $ 0\leq \delta \leq  \frac{1}{4} $\,, and $\ga > \de$\,, and $M \leq \eps \leq 1$\,, under the induction hypothesis holding true for both $A$ and on $h^1$\,, for all $|K|\leq |J| -1$\,, we have for $\ga^\prime = 0$\,, the following estimate in the exterior for $h^1$\,, 
             \beaa
   \notag
&& (1+t+|q|) \cdot |\varpi(q) \cdot \derm ( \Lie_{Z^J} h^1)_{UV} | \\
   \notag
   &\les&   C(q_0)  \cdot  c (\delta) \cdot c (\gamma) \cdot C ( |J|  ) \cdot E ( |J| + 4) \cdot \eps \cdot (1+t)^{     c (\gamma)  \cdot c (\delta)  \cdot c(|J|) \cdot E ( |J|+ 4) \cdot \eps }  \\
\notag
&& + c (\delta) \cdot c (\gamma) \cdot E (  4 )  \cdot \eps  \cdot  \int_0^t  \frac{1}{(1+\tau)} \cdot (1+\tau+|q|) \cdot \|\varpi(q)  \cdot  \derm  ( \Lie_{Z^J} h^1) (\tau,\cdot) \|_{L^\infty (\Sigma^{ext}_{\tau} )} d \tau \\
       \notag
    && + \sum_{|K| \leq |J|} \int_0^t (1+\tau) \cdot  \varpi(q) \cdot  \|    \Lie_{Z^K}  g^{\la\mu} \derm_{\la}   \derm_{\mu}  h^1_{UV} (\tau,\cdot) \|_{L^\infty(\overline{D}_\tau)} d\tau \; .
 \eeaa

Injecting, we obtain,
             \beaa
   \notag
&& (1+t+|q|) \cdot |\varpi(q) \cdot \derm ( \Lie_{Z^J} h^1) | \\
   \notag
   &\les&   C(q_0)  \cdot  c (\delta) \cdot c (\gamma) \cdot C ( |J|  ) \cdot E ( |J| + 4) \cdot \eps \cdot (1+t)^{     c (\gamma)  \cdot c (\delta)  \cdot c(|J|) \cdot E ( |J|+ 4) \cdot \eps }  \\
\notag
&& + c (\delta) \cdot c (\gamma) \cdot E (  4 )  \cdot \eps  \cdot  \int_0^t  \frac{1}{(1+\tau)} \cdot (1+\tau+|q|) \cdot \|\varpi(q)  \cdot  \derm  ( \Lie_{Z^J} h^1) (\tau,\cdot) \|_{L^\infty (\Sigma^{ext}_{\tau} )} d \tau \\
       \notag
&& +   \sum_{|K| \leq |J|}    \int_0^t C(q_0) \cdot   c (\gamma)  \cdot c (\delta)   \cdot E (4)\cdot C ( |J|  )   \cdot  \frac{\eps }{ (1+t+|q|)} \cdot(1+t+| q |  ) \cdot \varpi(q) \cdot |\derm ( \Lie_{Z^K} h_{ {\cal T} {\cal U}} ) | d \tau   \\
 &\les&   C(q_0)  \cdot  c (\delta) \cdot c (\gamma) \cdot C ( |J|  ) \cdot E ( |J| + 4) \cdot \eps \cdot (1+t)^{     c (\gamma)  \cdot c (\delta)  \cdot c(|J|) \cdot E ( |J|+ 4) \cdot \eps }  \\
\notag
&& + C(q_0)  \cdot  c (\delta) \cdot c (\gamma) \cdot E (  4 )  \cdot C(|J|)  \cdot \eps \\
&&  \times  \int_0^t  \frac{1}{(1+\tau)} \cdot (1+\tau+|q|) \cdot \sum_{|K| \leq |J|} \|\varpi(q)  \cdot  \derm  ( \Lie_{Z^K} h^1) (\tau,\cdot) \|_{L^\infty (\Sigma^{ext}_{\tau} )} d \tau \; .
  \eeaa
Summing over all $J$ with length $|J|$\,, we obtain the following Grönwall type inequality,
             \beaa
   \notag
&& (1+t+|q|) \cdot  \sum_{|K| \leq |J|}  |\varpi(q) \cdot \derm ( \Lie_{Z^K} h^1) | \\
   \notag
 &\les&   C(q_0)  \cdot  c (\delta) \cdot c (\gamma) \cdot C ( |J|  ) \cdot E ( |J| + 4) \cdot \eps \cdot (1+t)^{     c (\gamma)  \cdot c (\delta)  \cdot c(|J|) \cdot E ( |J|+ 4) \cdot \eps }  \\
\notag
&& + C(q_0)  \cdot  c (\delta) \cdot c (\gamma) \cdot E (  4 )  \cdot C(|J|)  \cdot \eps \\
&&  \times  \int_0^t  \frac{1}{(1+\tau)} \cdot (1+\tau+|q|) \cdot \sum_{|K| \leq |J|} \|\varpi(q)  \cdot  \derm  ( \Lie_{Z^K} h^1) (\tau,\cdot) \|_{L^\infty (\Sigma^{ext}_{\tau} )} d \tau \; .
  \eeaa
 Using Grönwall lemma, we obtain
    \beaa
   \notag
&& (1+t+|q|) \cdot \sum_{|K| \leq |J|} |\varpi(q) \cdot \derm ( \Lie_{Z^K} h^1) | \\
   \notag
    &\les &    C(q_0)  \cdot  c (\delta) \cdot c (\gamma) \cdot C ( |J|  ) \cdot E ( |J| + 4) \cdot \eps \cdot (1+t)^{     c (\gamma)  \cdot c (\delta)  \cdot c(|J|) \cdot E ( |J|+ 4) \cdot \eps } \\
    && \times  (1+t)^{C(q_0)  \cdot c (\gamma)  \cdot c (\delta)  \cdot E ( 4)  \cdot C ( |J|  ) \cdot \eps }  \\
    &\les &      C(q_0)  \cdot  c (\delta) \cdot c (\gamma) \cdot C ( |J|  ) \cdot E ( |J| + 4) \cdot \eps\cdot (1+t)^{ C(q_0)  \cdot    c (\gamma)  \cdot c (\delta)   \cdot E ( |J|+ 4) \cdot C ( |J|  ) \cdot \eps } \; .
       \eeaa
   
Therefore, given that $\ga^\prime = 0$, we get
            \beaa
 \notag
&& |\derm  ( \Lie_{Z^K} h^1) (t,x)  |   \\
\notag
 &\leq&   C(q_0)   \cdot  c (\delta) \cdot c (\gamma) \cdot C(|J|) \cdot E (  |J|  +4)  \cdot \frac{\eps  }{(1+t+|q|)^{1-   C(q_0)  \cdot    c (\gamma)  \cdot c (\delta)   \cdot E ( |J|+ 4) \cdot C ( |J|  )  \cdot \eps \cdot (1+M) } \cdot (1+|q|)}     \; . 
      \eeaa

       \end{proof}   
          
  \begin{remark}
  We notice that in Lemma \ref{upgradedestimateforhwithaconstantCqo}, we upgraded the estimate \eqref{formulaforinductionhypothesisonh^1} for $\Lie_{Z^J} h$, without having to deal first with $\Lie_{Z^J} h_{ {\cal T} {\cal U} }$ at the expense that we obtained a constant $C(q_0)$ in the exponent for $t$. However, one can get rid of that $C(q_0)$ by upgrading \eqref{formulaforinductionhypothesisonh^1} for $\Lie_{Z^J} h_{ {\cal T} {\cal U} }$ first and then upgrading for  $\Lie_{Z^J} h$, as we did for $A$ by upgrading first the “good” components, using \eqref{morerefinedcommutationformularusingonlytangentialcomponennts} in Lemma \ref{commutationformaulamoreprecisetoconservegpodcomponentsstructure} and then using this upgraded estimate to upgrade for the “bad” components through \label{commutationformaulaforallcomponentsasknowninlitter}, and thus for the full components. Doing so for $h$\,, we obtain the following slightly improved Lemma \ref{SteptwoforinductonforhwithoutCqzero} (improved with respect to Lemma \ref{upgradedestimateforhwithaconstantCqo}).
  \end{remark}
  \begin{lemma}\label{SteptwoforinductonforhwithoutCqzero}

   Let $ 0 < \delta \leq  \frac{1}{4} $, and $\ga \geq 3 \de$. We assume the induction hypothesis holding true for both $A$ and on $h^1$, for all $|K|\leq |J| -1$. 
   
  Then, in the exterior region $\overline{C}$, we have the following estimate on all the components of the metric,
              \beaa
 \notag
&& |\derm  ( \Lie_{Z^K} h^1) (t,x)  |   \\
\notag
 &\leq&   C(q_0)   \cdot  c (\delta) \cdot c (\gamma) \cdot C(|J|) \cdot E (  |J|  +4)  \cdot \frac{\eps  }{(1+t+|q|)^{1-       c (\gamma)  \cdot c (\delta)   \cdot E ( |J|+ 4) \cdot C ( |J|  )  \cdot  \eps } \cdot (1+|q|)}     \; . 
      \eeaa
\end{lemma}

Lemmas \ref{SteptwoforinductonforA} and \ref{SteptwoforinductonforhwithoutCqzero}, along with Lemma \ref{Step1fortheinductiononbothAandh1}, close the induction argument to prove Lemma \ref{upgradedestimateonLiederivativesoffields}.
\end{proof}

\section{Other upgraded decay estimates for the Einstein-Yang-Mills fields}

We will show upgraded estimates, other than the one we proved in Lemma \ref{upgradedestimateonLiederivativesoffields}. We assume from now on that $M \leq \eps$\,, so that we do not need to re-mention this assumption each time.

\begin{lemma}\label{upgradedestimateonLiederivativesoffieldswithoutgradiant}
In the Lorenz and harmonic gauges, the Einstein-Yang-Mills fields satisfy for $\ga \geq 3 \de $, and $0 < \de \leq \frac{1}{4}$, and for any $|K|  \in \N $,
            \bea
 \notag
&& | \Lie_{Z^K} A (t,x)  |    \\
\notag
&\leq&   C(q_0)   \cdot  c (\delta) \cdot c (\gamma) \cdot C(|K|) \cdot E ( |K| + 4)  \cdot \frac{\eps }{(1+t+|q|)^{1-      c (\gamma)  \cdot c (\delta)  \cdot c(|K|) \cdot E ( |K|+ 4) \cdot \eps } \cdot (1+|q|)^{\gamma - 2\de }}    \; ,
      \eea
      and
            \bea
 \notag
&& |\Lie_{Z^K} h^1 (t,x)  |   \\
\notag
 &\leq&   C(q_0)   \cdot  c (\delta) \cdot c (\gamma) \cdot C(|K|) \cdot E (  |K|  +4)  \cdot \frac{\eps   }{(1+t+|q|)^{1-   c (\gamma)  \cdot c (\delta)  \cdot c(|K|) \cdot E ( |K|+ 4)\cdot  \eps } }     \; ,
      \eea
      
      where $c (\gamma)$, $c (\delta)$, $c(|K|)$ are constants that depend respectively on $\gamma$, on $\de$, and on $|K|$.
 \end{lemma}
  
  \begin{proof}
By integrating the estimates in Lemma \ref{upgradedestimateonLiederivativesoffields}, at a fixed $\Om \in \SSS^2$ along the line $(\tau, r \cdot \Om)$, as in \cite{G4}, we obtain
            \bea
 \notag
&& | \Lie_{Z^K} A (t,x)  |    \\
\notag
&\leq&   C(q_0)   \cdot  c (\delta) \cdot c (\gamma) \cdot C(|K|) \cdot E ( |K| + 4)  \cdot \frac{\eps }{(1+t+|q|)^{1-      c (\gamma)  \cdot c (\delta)  \cdot c(|K|) \cdot E ( |K|+ 4) \cdot \eps } \cdot (1+|q|)^{\gamma - 2\de }}    \; ,\\
      \eea
      and
                  \bea
 \notag
&& |\Lie_{Z^K} h^1 (t,x)  |   \\
\notag
 &\leq&   C(q_0)   \cdot  c (\delta) \cdot c (\gamma) \cdot C(|K|) \cdot E (  |K|  +4)  \cdot \frac{\eps  \cdot \ln(1+|q|)}{(1+t+|q|)^{1-   c (\gamma)  \cdot c (\delta)  \cdot c(|K|) \cdot E ( |K|+ 4)\cdot  \eps } }     \; . \\
      \eea
      Thus, for any $\eps^\prime$, we have
                        \bea
 \notag
&& |\Lie_{Z^K} h^1 (t,x)  |   \\
\notag
 &\leq&   C(q_0)   \cdot  c (\delta) \cdot c (\gamma) \cdot C(|K|) \cdot E (  |K|  +4)  \cdot \frac{\eps }{(1+t+|q|)^{1-   c (\gamma)  \cdot c (\delta)  \cdot c(|K|) \cdot E ( |K|+ 4)\cdot  \eps - \eps^{\prime} } }     \; , \\
      \eea
      Choosing $\eps^\prime = c (\gamma)  \cdot c (\delta)  \cdot c(|K|) \cdot E ( |K|+ 4)\cdot  \eps$, we get the result, with the understanding that the constants are not the same.
      
  \end{proof}

\begin{lemma}\label{improvedtangentialderivatives}
For $\ga \geq 3 \de $, and $0 < \de \leq \frac{1}{4}$, we have in the exterior region $\overline{C} \subset \{ q \geq q_0 \} $, for all $ I $, 
             \bea
 \notag
 |\rderm ( \Lie_{Z^I}  A ) |    &\leq&   C(q_0)   \cdot  c (\delta) \cdot c (\gamma) \cdot C(|I|) \cdot E ( |I| + 5)  \cdot \frac{\eps }{(1+t+|q|)^{2-      c (\gamma)  \cdot c (\delta)  \cdot c(|I|) \cdot E ( |I|+ 5) \cdot \eps } \cdot (1+|q|)^{\gamma - 2\de }}    \; ,
      \eea      
 and  
             \bea
 \notag
 |\rderm ( \Lie_{Z^I}  h^1 ) |    &\leq&   C(q_0)   \cdot  c (\delta) \cdot c (\gamma) \cdot C(|I|) \cdot E (  |I|  +5)  \cdot \frac{\eps }{(1+t+|q|)^{2-   c (\gamma)  \cdot c (\delta)  \cdot c(|I|) \cdot E ( |I|+ 5)\cdot  \eps } }     \; .
      \eea
      
\end{lemma}

\begin{proof}

   Using Lemma \ref{betterdecayfortangentialderivatives}, and Lemma \ref{upgradedestimateonLiederivativesoffieldswithoutgradiant}, we obtain the result.

       \end{proof}

\begin{lemma}\label{upgradedestimatesonh}
For $\ga \geq 3 \de $, and $0 < \de \leq \frac{1}{4}$, we have in the exterior region $\overline{C} \subset \{ (t, x) \; | \;  q \geq q_0 \} $, for all $|I|$, 
   \beaa
 \notag
&&  |  \rderm ( \Lie_{Z^I} h ) (t,x)  |  \\
\notag
&\leq&   C(q_0)   \cdot  c (\delta) \cdot c (\gamma) \cdot C(|I|) \cdot E (  |I|  +5)  \cdot \frac{\eps }{(1+t+|q|)^{2-   c (\gamma)  \cdot c (\delta)  \cdot c(|I|) \cdot E ( |I|+ 5)\cdot  \eps } }     \; , 
\eeaa
              \beaa
 \notag
 && |\derm  ( \Lie_{Z^I} h ) (t,x)  | \\
  &\leq&   C(q_0)   \cdot  c (\delta) \cdot c (\gamma) \cdot C(|I|) \cdot E (  |I|  +4)  \cdot \frac{\eps }{(1+t+|q|)^{1-   c (\gamma)  \cdot c (\delta)  \cdot c(|I|) \cdot E ( |I|+ 4)\cdot  \eps } \cdot (1+|q|)}     \; , 
      \eeaa
            and
                 \beaa
 \notag
 |   \Lie_{Z^I} h (t,x)  | &\leq&   C(q_0)   \cdot  c (\delta) \cdot c (\gamma) \cdot C(|I|) \cdot E (  |I|  +4)  \cdot \frac{\eps }{(1+t+|q|)^{1-   c (\gamma)  \cdot c (\delta)  \cdot c(|I|) \cdot E ( |I|+ 4)\cdot  \eps } }     \; . \\
      \eeaa

\end{lemma}

\begin{proof}

We have shown in Lemma \ref{improvedtangentialderivatives}, that for $\ga \geq 3 \de $, and $0 < \de \leq \frac{1}{4}$, for all $I$,       
             \beaa
 \notag
 |\rderm ( \Lie_{Z^I}  h^1 ) |    &\leq&   C(q_0)   \cdot  c (\delta) \cdot c (\gamma) \cdot C(|I|) \cdot E (  |I|  +5)  \cdot \frac{\eps  }{(1+t+|q|)^{2-   c (\gamma)  \cdot c (\delta)  \cdot c(|I|) \cdot E ( |I|+ 5)\cdot  \eps } }     \; , \\
      \eeaa
      
and in Lemma \ref{tangentialderivativesphericalsymmetricpart},
 \beaa
 \notag
|  \rderm ( \Lie_{Z^I} h^0 ) (t,x)  |   &\leq& C ( |I|  ) \cdot \frac{\eps }{(1+t+|q|)^{2}  } \; .
      \eeaa
Consequently,
   \bea
 \notag
|  \rderm ( \Lie_{Z^I} h ) (t,x)  |  &\leq& |  \rderm ( \Lie_{Z^I} h^0 ) (t,x)  | + |  \rderm ( \Lie_{Z^I} h^1 ) (t,x)  |  \\
 \notag
 &\leq&   C(q_0)   \cdot  c (\delta) \cdot c (\gamma) \cdot C(|I|) \cdot E (  |I|  +5)  \cdot \frac{\eps  }{(1+t+|q|)^{2-   c (\gamma)  \cdot c (\delta)  \cdot c(|I|) \cdot E ( |I|+ 5)\cdot  \eps } }     \; . \\
      \eea

  Also, from Lemma \ref{upgradedestimateonLiederivativesoffields}, we have for $\ga \geq 3 \de $, and $0 < \de \leq \frac{1}{4}$,
            \beaa
 \notag
&& |\derm  ( \Lie_{Z^K} h^1) (t,x)  |   \\
\notag
 &\leq&   C(q_0)   \cdot  c (\delta) \cdot c (\gamma) \cdot C(|K|) \cdot E (  |K|  +4)  \cdot \frac{\eps   }{(1+t+|q|)^{1-   c (\gamma)  \cdot c (\delta)  \cdot c(|K|) \cdot E ( |K|+ 4)\cdot  \eps } \cdot (1+|q|)}     \; , \\
      \eeaa
and since from Lemma \ref{Liederivativesofsphericalsymmetricpart}, we have 
 \beaa
 \notag
|\derm  ( \Lie_{Z^I} h^0 ) (t,x)  |   &\leq& C ( |I| )   \cdot \frac{\eps }{(1+t+|q|)^{2} } \; ,
      \eeaa
we then get,
           \bea
 \notag
 && |\derm  ( \Lie_{Z^K} h)  | \\
 \notag
  &\leq& |\derm  ( \Lie_{Z^K} h^1 )  | +  |\derm  ( \Lie_{Z^I} h^0 )  |   \\
\notag
 &\leq&   C(q_0)   \cdot  c (\delta) \cdot c (\gamma) \cdot C(|K|) \cdot E (  |K|  +4)  \cdot \frac{\eps}{(1+t+|q|)^{1-   c (\gamma)  \cdot c (\delta)  \cdot c(|K|) \cdot E ( |K|+ 4)\cdot  \eps } \cdot (1+|q|)}     \; . \\
      \eea
By integrating, as in \cite{G4}, we obtain that for $\ga \geq 3 \de $, and $0 < \de \leq \frac{1}{4}$,

 \bea
 \notag
&& | \Lie_{Z^K} h  |  \\
\notag
 &\leq&   C(q_0)   \cdot  c (\delta) \cdot c (\gamma) \cdot C(|K|) \cdot E (  |K|  +4)  \cdot \frac{\eps   }{(1+t+|q|)^{1-   c (\gamma)  \cdot c (\delta)  \cdot c(|K|) \cdot E ( |K|+ 4)\cdot  \eps } }     \; . \\
      \eea

\end{proof}

\begin{lemma}\label{upgradedestimategoodcomponentspotentialandmetric}
Under the bootstrap assumption holding for all $|J| \leq  \lfloor \frac{|I|}{2} \rfloor  $, we have for $\ga \geq 3 \de $, and $0 < \de \leq \frac{1}{4}$, in the exterior region,

      \beaa
&& | \pa ( \Lie_{Z^I}  A_{L} )  |  \\
   &\leq&    C(q_0)   \cdot  c (\delta) \cdot c (\gamma) \cdot C(|I|) \cdot E ( |I| + 5)  \cdot \frac{\eps }{(1+t+|q|)^{1-      c (\gamma)  \cdot c (\delta)  \cdot c(|I|) \cdot E ( |I|+ 5) \cdot \eps } \cdot (1+|q|)^{1+\gamma - 2\de }} \; ,
      \eeaa   
and
              \beaa
        \notag
 |  \Lie_{Z^I}  A_{L}  |    &\leq& C(q_0)   \cdot  c (\delta) \cdot c (\gamma) \cdot C(|I|) \cdot E ( |I| + 5)  \cdot \frac{\eps }{(1+t+|q|)^{1-      c (\gamma)  \cdot c (\delta)  \cdot c(|I|) \cdot E ( |I|+ 5) \cdot \eps } \cdot (1+|q|)^{\gamma - 2\de }} \; .
      \eeaa

\end{lemma}

            \begin{proof}
            
            In the Lorenz gauge, we showed in Lemma \ref{estimateonpartialderivativeofALcomponent}, that

\beaa
  | \pa  \Lie_{Z^I}  A_{L}  |  \les \sum_{|J|\leq |I|} |   \rderm  ( \Lie_{Z^J} A ) | + \sum_{|J|\leq |I| -1} |   \derm  ( \Lie_{Z^J} A ) |  + \sum_{|K| + |M| \leq |I|} O \big(  | ( \Lie_{Z^K}  h )  | \cdot  | \derm  (  \Lie_{Z^M}  A ) | \big) \; .
\eeaa

From Lemma \ref{improvedtangentialderivatives}, we have for $\ga \geq 3 \de $, and $0 < \de \leq \frac{1}{4}$,
             \bea
 \notag
&&  \sum_{|J|\leq |I|}    |\rderm ( \Lie_{Z^J}  A ) |   \\
 \notag &\leq&   C(q_0)   \cdot  c (\delta) \cdot c (\gamma) \cdot C(|I|) \cdot E ( |I| + 5)  \cdot \frac{\eps }{(1+t+|q|)^{2-      c (\gamma)  \cdot c (\delta)  \cdot c(|I|) \cdot E ( |I|+ 5) \cdot \eps } \cdot (1+|q|)^{\gamma - 2\de }}    \; ,
      \eea     
 
Also, based on Lemma \ref{upgradedestimatesonh}, we have for all $|K| \leq |I|$,
                 \beaa
 \notag
 |   \Lie_{Z^K} h  | &\leq&   C(q_0)   \cdot  c (\delta) \cdot c (\gamma) \cdot C(|I|) \cdot E (  |I|  +4)  \cdot \frac{\eps }{(1+t+|q|)^{1-   c (\gamma)  \cdot c (\delta)  \cdot c(|I|) \cdot E ( |I|+ 4)\cdot  \eps } }     \; ,
      \eeaa
      
    and from Lemma \ref{upgradedestimateonLiederivativesoffields}, we have for all $|M| \leq |I|$, 
            \beaa
 \notag
&& |\derm  ( \Lie_{Z^M} A)  |    \\
\notag
&\leq&   C(q_0)   \cdot  c (\delta) \cdot c (\gamma) \cdot C(|I|) \cdot E ( |I| + 4)  \cdot \frac{\eps }{(1+t+|q|)^{1-      c (\gamma)  \cdot c (\delta)  \cdot c(|I|) \cdot E ( |I|+ 4) \cdot \eps } \cdot (1+|q|)^{1+\gamma - 2\de }}    \; .
      \eeaa

Thus,
                 \bea
 \notag
&& \sum_{|K| + |M| \leq |I|}  |   \Lie_{Z^K} h  |  \cdot |\derm (  \Lie_{Z^M} A )   |   \\
 \notag
&\leq&   C(q_0)   \cdot  c (\delta) \cdot c (\gamma) \cdot C(|I|) \cdot E ( |I| + 4)  \cdot \frac{\eps }{(1+t+|q|)^{2-      c (\gamma)  \cdot c (\delta)  \cdot c(|I|) \cdot E ( |I|+ 4) \cdot \eps } \cdot (1+|q|)^{1+\gamma - 2\de }}    \; . 
      \eea

    Also,
    \beaa
 && \sum_{|J|\leq |I| -1} |   \derm  ( \Lie_{Z^J} A ) |  \\
 &\leq&    C(q_0)   \cdot  c (\delta) \cdot c (\gamma) \cdot C(|I|) \cdot E ( |I| + 3)  \cdot \frac{\eps }{(1+t+|q|)^{1-      c (\gamma)  \cdot c (\delta)  \cdot c(|I|) \cdot E ( |I|+ 3) \cdot \eps } \cdot (1+|q|)^{1+\gamma - 2\de }}    \; .
      \eeaa
      
      As a result,
      \beaa
&& | \pa ( \Lie_{Z^I}  A_{L} )  |  \\
  &\leq&  C(q_0)   \cdot  c (\delta) \cdot c (\gamma) \cdot C(|I|) \cdot E ( |I| + 5)  \cdot \frac{\eps }{(1+t+|q|)^{2-      c (\gamma)  \cdot c (\delta)  \cdot c(|I|) \cdot E ( |I|+ 5) \cdot \eps } \cdot (1+|q|)^{\gamma - 2\de }}    \\
 && +  C(q_0)   \cdot  c (\delta) \cdot c (\gamma) \cdot C(|I|) \cdot E ( |I| + 3)  \cdot \frac{\eps }{(1+t+|q|)^{1-      c (\gamma)  \cdot c (\delta)  \cdot c(|I|) \cdot E ( |I|+ 3) \cdot \eps } \cdot (1+|q|)^{1+\gamma - 2\de }} \\
 && +  C(q_0)   \cdot  c (\delta) \cdot c (\gamma) \cdot C(|I|) \cdot E ( |I| + 4)  \cdot \frac{\eps }{(1+t+|q|)^{2-      c (\gamma)  \cdot c (\delta)  \cdot c(|I|) \cdot E ( |I|+ 4) \cdot \eps } \cdot (1+|q|)^{1+\gamma - 2\de }} \\
   &\leq&    C(q_0)   \cdot  c (\delta) \cdot c (\gamma) \cdot C(|I|) \cdot E ( |I| + 5)  \cdot \frac{\eps }{(1+t+|q|)^{1-      c (\gamma)  \cdot c (\delta)  \cdot c(|I|) \cdot E ( |I|+ 5) \cdot \eps } \cdot (1+|q|)^{1+\gamma - 2\de }} \; .
      \eeaa

By integrating along the line $(\tau, r \cdot \Om)$ such that $r+\tau =  | x | +t $ till we reach the hyperplane $\tau=0$, i.e. along the null coordinate $s=\tau+r$, we obtain

              \beaa
        \notag
 |  \Lie_{Z^I}  A_{L}  |    &\leq& C(q_0)   \cdot  c (\delta) \cdot c (\gamma) \cdot C(|I|) \cdot E ( |I| + 5)  \cdot \frac{\eps }{(1+t+|q|)^{1-      c (\gamma)  \cdot c (\delta)  \cdot c(|I|) \cdot E ( |I|+ 5) \cdot \eps } \cdot (1+|q|)^{\gamma - 2\de }} \; .
      \eeaa   
            \end{proof}

\begin{lemma}\label{upgradedestimatesonproducts}
We have the following estimates for $|J|\,, |K| \leq |I|$\,  
\beaa
 && |\derm  ( \Lie_{Z^J} h )   | \cdot  | \Lie_{Z^K} A   |  \\
    &\les& C(q_0)   \cdot  c (\delta) \cdot c (\gamma) \cdot C(|I|) \cdot E ( |I| + 4)  \cdot \frac{\eps  }{(1+t+|q|)^{2-      c (\gamma)  \cdot c (\delta)  \cdot c(|I|) \cdot E ( |I|+ 4) \cdot \eps } \cdot (1+|q|)^{1+\gamma - 2\de }}  \; ,
 \eeaa
 
            \beaa
 \notag
 && |\derm  ( \Lie_{Z^J} h )   | \cdot  |   \Lie_{Z^K} h   |   \\
  &\leq&   C(q_0)   \cdot  c (\delta) \cdot c (\gamma) \cdot C(|I|) \cdot E (  |I|  +4)  \cdot \frac{\eps  }{(1+t+|q|)^{2-   c (\gamma)  \cdot c (\delta)  \cdot c(|I|) \cdot E ( |I|+ 4)\cdot  \eps } \cdot (1+|q|)} \; ,
  \eeaa

                \beaa
 \notag
 && |\derm  ( \Lie_{Z^J} A )   | \cdot  |   \Lie_{Z^K} A   |  \\
    &\leq&  C(q_0)   \cdot  c (\delta) \cdot c (\gamma) \cdot C(|I|) \cdot E ( |I| + 4)  \cdot \frac{\eps }{(1+t+|q|)^{2-      c (\gamma)  \cdot c (\delta)  \cdot c(|I|) \cdot E ( |I|+ 4) \cdot \eps } \cdot (1+|q|)^{1+2\gamma - 4\de }}   \; ,
      \eeaa

             \beaa
 \notag
 && |\derm  ( \Lie_{Z^J} A )   | \cdot  |   \Lie_{Z^K} h   |   \\
    &\leq&  C(q_0)   \cdot  c (\delta) \cdot c (\gamma) \cdot C(|I|) \cdot E ( |I| + 4)  \cdot \frac{\eps  }{(1+t+|q|)^{2-      c (\gamma)  \cdot c (\delta)  \cdot c(|I|) \cdot E ( |I|+ 4) \cdot \eps } \cdot (1+|q|)^{1+\gamma - 2\de }}  \; ,
    \eeaa

                  \beaa
 \notag
 &&  |\derm  ( \Lie_{Z^J} A )   | \cdot   |\derm  ( \Lie_{Z^K} h )   |  \\
      &\leq&  C(q_0)   \cdot  c (\delta) \cdot c (\gamma) \cdot C(|I|) \cdot E ( |I| + 4)  \cdot \frac{\eps  }{(1+t+|q|)^{2-      c (\gamma)  \cdot c (\delta)  \cdot c(|I|) \cdot E ( |K|+ 4) \cdot \eps } \cdot (1+|q|)^{2+\gamma - 2\de }}  \; ,
        \eeaa

                      \beaa
 \notag
&& | \Lie_{Z^J} A    | \cdot   | \Lie_{Z^K} h    |   \\
    &\leq&  C(q_0)   \cdot  c (\delta) \cdot c (\gamma) \cdot C(|I|) \cdot E ( |I| + 4)  \cdot \frac{\eps  }{(1+t+|q|)^{2-      c (\gamma)  \cdot c (\delta)  \cdot c(|I|) \cdot E ( |I|+ 4) \cdot \eps } \cdot (1+|q|)^{\gamma - 2\de }}   \; .
    \eeaa
    
\end{lemma}
 
 \begin{proof}
 Based on the upgraded estimates that we have just shown in Lemma \ref{upgradedestimateonLiederivativesoffields}, Lemma \ref{upgradedestimateonLiederivativesoffieldswithoutgradiant}, and Lemma \ref{upgradedestimatesonh}, we obtain in the exterior region
 \beaa
 && |\derm  ( \Lie_{Z^J} h )   | \cdot  | \Lie_{Z^K} A   |  \\
 &\les&C(q_0)   \cdot  c (\delta) \cdot c (\gamma) \cdot C(|I|) \cdot E (  |I|  +4)  \cdot \frac{\eps  }{(1+t+|q|)^{1-   c (\gamma)  \cdot c (\delta)  \cdot c(|I|) \cdot E ( |I|+ 4)\cdot  \eps } \cdot (1+|q|)} \\
        && \times  C(q_0)   \cdot  c (\delta) \cdot c (\gamma) \cdot C(|I|) \cdot E ( |I| + 4)  \cdot \frac{\eps }{(1+t+|q|)^{1-      c (\gamma)  \cdot c (\delta)  \cdot c(|I|) \cdot E ( |I|+ 4) \cdot \eps } \cdot (1+|q|)^{\gamma - 2\de }}  \\
    &\les& C(q_0)   \cdot  c (\delta) \cdot c (\gamma) \cdot C(|I|) \cdot E ( |I| + 4)  \cdot \frac{\eps }{(1+t+|q|)^{2-      c (\gamma)  \cdot c (\delta)  \cdot c(|I|) \cdot E ( |I|+ 4) \cdot \eps } \cdot (1+|q|)^{1+\gamma - 2\de }}  \; .
 \eeaa
 
Also,
            \beaa
 \notag
 && |\derm  ( \Lie_{Z^J} h )   | \cdot  |   \Lie_{Z^K} h   |   \\
 \notag
  &\leq&   C(q_0)   \cdot  c (\delta) \cdot c (\gamma) \cdot C(|I|) \cdot E (  |I|  +4)  \cdot \frac{\eps}{(1+t+|q|)^{1-   c (\gamma)  \cdot c (\delta)  \cdot c(|I|) \cdot E ( |I|+ 4)\cdot  \eps } \cdot (1+|q|)} \\
 \notag
  &&\times   C(q_0)   \cdot  c (\delta) \cdot c (\gamma) \cdot C(|I|) \cdot E (  |I|  +4)  \cdot \frac{\eps }{(1+t+|q|)^{1-   c (\gamma)  \cdot c (\delta)  \cdot c(|I|) \cdot E ( |I|+ 4)\cdot  \eps } }   \\
  \notag
  &\leq&   C(q_0)   \cdot  c (\delta) \cdot c (\gamma) \cdot C(|I|) \cdot E (  |I|  +4)  \cdot \frac{\eps  }{(1+t+|q|)^{2-   c (\gamma)  \cdot c (\delta)  \cdot c(|I|) \cdot E ( |I|+ 4)\cdot  \eps } \cdot (1+|q|)} \; .
  \eeaa

And,
                \beaa
 \notag
 && |\derm  ( \Lie_{Z^J} A )   | \cdot  |   \Lie_{Z^K} A   |  \\
   &\leq& C(q_0)   \cdot  c (\delta) \cdot c (\gamma) \cdot C(|I|) \cdot E ( |I| + 4)  \cdot \frac{\eps }{(1+t+|q|)^{1-      c (\gamma)  \cdot c (\delta)  \cdot c(|I|) \cdot E ( |I|+ 4) \cdot \eps } \cdot (1+|q|)^{1+\gamma - 2\de }}  \\
  \notag
 &&\times    C(q_0)   \cdot  c (\delta) \cdot c (\gamma) \cdot C(|I|) \cdot E ( |I| + 4)  \cdot \frac{\eps }{(1+t+|q|)^{1-      c (\gamma)  \cdot c (\delta)  \cdot c(|I|) \cdot E ( |I|+ 4) \cdot \eps } \cdot (1+|q|)^{\gamma - 2\de }} \\
 \notag 
    &\leq&  C(q_0)   \cdot  c (\delta) \cdot c (\gamma) \cdot C(|I|) \cdot E ( |I| + 4)  \cdot \frac{\eps }{(1+t+|q|)^{2-      c (\gamma)  \cdot c (\delta)  \cdot c(|I|) \cdot E ( |I|+ 4) \cdot \eps } \cdot (1+|q|)^{1+2\gamma - 4\de }} 
        \eeaa

We estimate for the following product,
             \beaa
 \notag
 && |\derm  ( \Lie_{Z^J} A )   | \cdot  |   \Lie_{Z^K} h   |   \\
  &\leq& C(q_0)   \cdot  c (\delta) \cdot c (\gamma) \cdot C(|I|) \cdot E ( |I| + 4)  \cdot \frac{\eps }{(1+t+|q|)^{1-      c (\gamma)  \cdot c (\delta)  \cdot c(|I|) \cdot E ( |K|+ 4) \cdot \eps } \cdot (1+|q|)^{1+\gamma - 2\de }}  \\
  \notag
  && \times  C(q_0)   \cdot  c (\delta) \cdot c (\gamma) \cdot C(|I|) \cdot E (  |I|  +4)  \cdot \frac{\eps }{(1+t+|q|)^{1-   c (\gamma)  \cdot c (\delta)  \cdot c(|I|) \cdot E ( |I|+ 4)\cdot  \eps } }   \\
    &\leq&  C(q_0)   \cdot  c (\delta) \cdot c (\gamma) \cdot C(|I|) \cdot E ( |I| + 4)  \cdot \frac{\eps  }{(1+t+|q|)^{2-      c (\gamma)  \cdot c (\delta)  \cdot c(|I|) \cdot E ( |I|+ 4) \cdot \eps } \cdot (1+|q|)^{1+\gamma - 2\de }}  \;.
    \eeaa

We also have
                  \beaa
 \notag
 &&  |\derm  ( \Lie_{Z^J} A )   | \cdot   |\derm  ( \Lie_{Z^K} h )   |  \\
  &\leq&    C(q_0)   \cdot  c (\delta) \cdot c (\gamma) \cdot C(|I|) \cdot E ( |I| + 4)  \cdot \frac{\eps }{(1+t+|q|)^{1-      c (\gamma)  \cdot c (\delta)  \cdot c(|I|) \cdot E ( |K|+ 4) \cdot \eps } \cdot (1+|q|)^{1+\gamma - 2\de }}  \\
 \notag
  && \times  C(q_0)   \cdot  c (\delta) \cdot c (\gamma) \cdot C(|I|) \cdot E (  |I|  +4)  \cdot \frac{\eps  }{(1+t+|q|)^{1-   c (\gamma)  \cdot c (\delta)  \cdot c(|I|) \cdot E ( |I|+ 4)\cdot  \eps } \cdot (1+|q|)} \\
 \notag
    &\leq&  C(q_0)   \cdot  c (\delta) \cdot c (\gamma) \cdot C(|I|) \cdot E ( |I| + 4)  \cdot \frac{\eps   }{(1+t+|q|)^{2-      c (\gamma)  \cdot c (\delta)  \cdot c(|I|) \cdot E ( |K|+ 4) \cdot \eps } \cdot (1+|q|)^{2+\gamma - 2\de }}  \; .
        \eeaa
Whereas to the following term, we estimate
                      \beaa
 \notag
&& | \Lie_{Z^J} A    | \cdot   | \Lie_{Z^K} h    |   \\
 &\leq&  C(q_0)   \cdot  c (\delta) \cdot c (\gamma) \cdot C(|I|) \cdot E ( |I| + 4)  \cdot \frac{\eps }{(1+t+|q|)^{1-      c (\gamma)  \cdot c (\delta)  \cdot c(|I|) \cdot E ( |I|+ 4) \cdot \eps } \cdot (1+|q|)^{\gamma - 2\de }}   \\
      &&\times   C(q_0)   \cdot  c (\delta) \cdot c (\gamma) \cdot C(|I|) \cdot E (  |I|  +4)  \cdot \frac{\eps   }{(1+t+|q|)^{1-   c (\gamma)  \cdot c (\delta)  \cdot c(|I|) \cdot E ( |I|+ 4)\cdot  \eps } }   \\ 
     \notag
    &\leq&  C(q_0)   \cdot  c (\delta) \cdot c (\gamma) \cdot C(|I|) \cdot E ( |I| + 4)  \cdot \frac{\eps   }{(1+t+|q|)^{2-      c (\gamma)  \cdot c (\delta)  \cdot c(|I|) \cdot E ( |I|+ 4) \cdot \eps } \cdot (1+|q|)^{\gamma - 2\de }}   \; .
    \eeaa
    
     \end{proof}
     
\section{Structure of the source terms for the wave operator on the Einstein-Yang-Mills potential}

\subsection{Structure of the “good” source terms for the Einstein-Yang-Mills potential}\
 
  \begin{lemma}\label{TheactualusefulstrzuctureofthesourcetermsforthewaveequationonpoentialAusingbootstrap}

For $\ga \geq 3 \de $\;, for $0 < \de \leq \frac{1}{4}$ in the bootstrap assumption, and for $\eps$ small enough depending on $\ga$\;, $\de$\;, and $|I|$\;, we have the following estimate for all $I$\;, and for all $V \in \cal{U}$\;,
                                                       \beaa
   \notag
&&  |  \Lie_{Z^I}  g^{\la\mu} \derm_{\la}   \derm_{\mu}   A_{V}  |  \\
   &\les&           \sum_{|K| \leq |I |}  \Big[   \;       O \big(   C(q_0)   \cdot  c (\delta) \cdot c (\gamma) \cdot C(|I|) \cdot E ( \lfloor \frac{|I|}{2} \rfloor + 5)  \cdot \frac{\eps   \cdot  | \derm ( \Lie_{Z^K} h^1 ) | }{(1+t+|q|)^{2-      c (\gamma)  \cdot c (\delta)  \cdot c(|I|) \cdot E ( \lfloor \frac{|I|}{2} \rfloor + 5) \cdot \eps } \cdot (1+|q|)^{\gamma - 2\de }}   \big) \\
   && + O \big(    C(q_0)   \cdot  c (\delta) \cdot c (\gamma) \cdot C(|I|) \cdot E ( \lfloor \frac{|I|}{2} \rfloor   +5)  \cdot \frac{\eps  \cdot  |\derm ( \Lie_{Z^K} A ) | }{(1+t+|q|)^{2-   c (\gamma)  \cdot c (\delta)  \cdot c(|I|) \cdot E ( \lfloor \frac{|I|}{2} \rfloor + 5)\cdot  \eps } }  \big) \\
   &&  +  C(q_0)   \cdot  c (\delta) \cdot c (\gamma) \cdot C(|I|) \cdot E ( \lfloor \frac{|I|}{2} \rfloor + 4)  \cdot \frac{\eps   \cdot    | \rderm ( \Lie_{Z^K} A )  | }{(1+t+|q|)^{1-      c (\gamma)  \cdot c (\delta)  \cdot c(|I|) \cdot E (\lfloor \frac{|I|}{2} \rfloor + 4) \cdot \eps } \cdot ( 1 + |q| )^{ \frac{1}{2}}  }    \\
      \notag
        &&    +     C(q_0)   \cdot  c (\delta) \cdot c (\gamma) \cdot C(|I|) \cdot E ( \lfloor \frac{|I|}{2} \rfloor + 4)  \cdot \frac{\eps  \cdot   | \rderm  ( \Lie_{Z^K} h^1 ) | }{(1+t+|q|)^{1-      c (\gamma)  \cdot c (\delta)  \cdot c(|I|) \cdot E ( \lfloor \frac{|I|}{2} \rfloor  + 4) \cdot \eps } \cdot (1+|q|)^{1+\gamma - 2\de }}    \\
      &&    +       O \big(     C(q_0)   \cdot  c (\delta) \cdot c (\gamma) \cdot C(|I|) \cdot E ( \lfloor \frac{|I|}{2} \rfloor  + 5)  \cdot \frac{\eps \cdot |  \Lie_{Z^K} A  | }{(1+t+|q|)^{2-      c (\gamma)  \cdot c (\delta)  \cdot c(|I|) \cdot E ( \lfloor \frac{|I|}{2} \rfloor + 5) \cdot \eps } \cdot (1+|q|)^{\gamma - 2\de }}   \\
         &&    +         O \big( C(q_0)   \cdot  c (\delta) \cdot c (\gamma) \cdot C(|I|) \cdot E ( \lfloor \frac{|I|}{2} \rfloor  + 4)  \cdot \frac{\eps \cdot | \Lie_{Z^K} h^1 | }{(1+t+|q|)^{2-      c (\gamma)  \cdot c (\delta)  \cdot c(|I|) \cdot E ( \lfloor \frac{|I|}{2} \rfloor + 4) \cdot \eps } \cdot (1+|q|)^{1+2\gamma - 4\de }}     \big) \; \Big] \\
         \notag
         && +     |   \Lie_{Z^I}    \big(  A_L   \cdot     \derm A    \big)  |   +  |  \Lie_{Z^I} \big( A_{e_a}  \cdot     \derm A_{e_a} \big)  |  \;  \\
                 &&    +     C(q_0)   \cdot  c (\delta) \cdot c (\gamma) \cdot C(|I|) \cdot E ( \lfloor \frac{|I|}{2} \rfloor + 5)  \cdot \frac{\eps^2 }{(1+t+|q|)^{3-      c (\gamma)  \cdot c (\delta)  \cdot c(|I|) \cdot E ( \lfloor \frac{|I|}{2} \rfloor  + 5) \cdot \eps } \cdot (1+|q|)^{1+\gamma - 2\de }}    \\
          \notag
\eeaa
For all $V \in \cal{T}$, we have the same estimate without the last term $   \Lie_{Z^I}    \big(  A_L   \cdot     \derm A    \big)  +  \Lie_{Z^I} \big( A_{e_a}  \cdot     \derm A_{e_a} \big)$.

      \begin{remark}\label{sufficientdecayinqinfactorfortangentialderivativeofA}
      We need sufficient decay in $q$ in the factor, in order to control the tangential derivatives $\rderm ( \Lie_{Z^I} A ) $ using the energy estimate established in Lemma \ref{Theveryfinal } and see \eqref{estimateonthetermthatcontainstangentialderivativeofAinsourcetermsforwaveopertaorinAwithdecayinqfactor}.
      
      \end{remark}
      
\end{lemma}
      
\begin{proof}

First, we consider the fact that for a function $f$\;, the sum could be decomposed as 
  \bea\label{generalinequalityaboutsplittingasumthatisconstrainedtoadduptolenth}
  \notag
  \sum_{|K| + |J| +|M| + |N| \leq |I| } |f| &\leq& \sum_{|J|\,, |M|\,, |N| \leq  \lfloor \frac{|I|}{2} \rfloor\;,\;  |K| \leq|I|  } |f| + \sum_{|K|\,, |M|\,, |N| \leq  \lfloor \frac{|I|}{2} \rfloor\;,\;  |J| \leq |I| } 
|f| \\
\notag
&& + \sum_{|K|\,, |J|\,, |N| \leq  \lfloor \frac{|I|}{2} \rfloor\;,\;  |M| \leq |I| } |f| + \sum_{|K|\,, |J|\,, |M| \leq \lfloor \frac{|I|}{2} \rfloor\;,\;  |N| \leq  |I| }  |f| \; . \\
  \eea

Now, Looking at the structure of the source term for the $A$, that we showed in \eqref{StructureoftheLiederivativesofthesourcetermsofthewaveoperatorforAandh} -- which preserves the structure of \eqref{structureofthesourcetermsofthewaveoperatoronAandh} --, and then, using Lemmas \ref{upgradedestimateonLiederivativesoffields}, \ref{upgradedestimateonLiederivativesoffieldswithoutgradiant}, \ref{improvedtangentialderivatives}, \ref{upgradedestimatesonh}, \ref{upgradedestimategoodcomponentspotentialandmetric}, and \ref{upgradedestimatesonproducts}, we obtain
  
\beaa
   \notag
&&  |  \Lie_{Z^I}  g^{\la\mu} \derm_{\la}   \derm_{\mu}   A_{\si}  |  \\
  \notag
  &\les&  \sum_{|K| \leq |I |}  \Big[   \;   C(q_0)   \cdot  c (\delta) \cdot c (\gamma) \cdot C(|I|) \cdot E (  \lfloor \frac{|I|}{2} \rfloor + 5)  \cdot \frac{\eps \cdot  | \derm ( \Lie_{Z^K} h ) | }{(1+t+|q|)^{2-      c (\gamma)  \cdot c (\delta)  \cdot c(|I|) \cdot E (\lfloor \frac{|I|}{2} \rfloor + 5) \cdot \eps } \cdot (1+|q|)^{\gamma - 2\de }}   \\
  && +  C(q_0)   \cdot  c (\delta) \cdot c (\gamma) \cdot C(|I|) \cdot E (  \lfloor \frac{|I|}{2} \rfloor   +4)  \cdot \frac{\eps\cdot  |\rderm ( \Lie_{Z^K} A ) |}{(1+t+|q|)^{1-   c (\gamma)  \cdot c (\delta)  \cdot c(|I|) \cdot E ( \lfloor \frac{|I|}{2} \rfloor + 4)\cdot  \eps } \cdot (1+|q|)}       \\
  &&    +     C(q_0)   \cdot  c (\delta) \cdot c (\gamma) \cdot C(|I|) \cdot E (\lfloor \frac{|I|}{2} \rfloor  + 4)  \cdot \frac{\eps  \cdot   | \rderm  ( \Lie_{Z^K} h ) | }{(1+t+|q|)^{1-      c (\gamma)  \cdot c (\delta)  \cdot c(|I|) \cdot E (\lfloor \frac{|I|}{2} \rfloor + 4) \cdot \eps } \cdot (1+|q|)^{1+\gamma - 2\de }}    \\
  && +    C(q_0)   \cdot  c (\delta) \cdot c (\gamma) \cdot C(|I|) \cdot E (  \lfloor \frac{|I|}{2} \rfloor   +5)  \cdot \frac{\eps \cdot  |\derm ( \Lie_{Z^K} A ) | }{(1+t+|q|)^{2-   c (\gamma)  \cdot c (\delta)  \cdot c(|I|) \cdot E (\lfloor \frac{|I|}{2} \rfloor + 5)\cdot  \eps } }  \\
   \notag
           && +     C(q_0)   \cdot  c (\delta) \cdot c (\gamma) \cdot C(|I|) \cdot E (\lfloor \frac{|I|}{2} \rfloor  + 5)  \cdot \frac{\eps  \cdot |  \Lie_{Z^K} A  | }{(1+t+|q|)^{2-      c (\gamma)  \cdot c (\delta)  \cdot c(|I|) \cdot E (\lfloor \frac{|I|}{2} \rfloor + 5) \cdot \eps } \cdot (1+|q|)^{\gamma - 2\de }}    \\
           && +     C(q_0)   \cdot  c (\delta) \cdot c (\gamma) \cdot C(|I|) \cdot E (\lfloor \frac{|I|}{2} \rfloor + 4)  \cdot \frac{\eps \cdot    | \rderm ( \Lie_{Z^K} A )  | }{(1+t+|q|)^{1-      c (\gamma)  \cdot c (\delta)  \cdot c(|I|) \cdot E ( \lfloor \frac{|I|}{2} \rfloor + 4) \cdot \eps } \cdot (1+|q|)^{\gamma - 2\de }}     \\
                    \notag
 &&          +    C(q_0)   \cdot  c (\delta) \cdot c (\gamma) \cdot C(| I |) \cdot E ( \lfloor \frac{|I|}{2} \rfloor  + 4)  \cdot \frac{\eps^2 \cdot   | \derm ( \Lie_{Z^K} h )    |  }{(1+t+|q|)^{2-      c (\gamma)  \cdot c (\delta)  \cdot c(|I|) \cdot E ( \lfloor \frac{|I|}{2} \rfloor + 4) \cdot \eps } \cdot (1+|q|)^{2\gamma - 4\de }}    \\
 &&  + C(q_0)   \cdot  c (\delta) \cdot c (\gamma) \cdot C(|I|) \cdot E ( \lfloor \frac{|I|}{2} \rfloor  + 4)  \cdot \frac{\eps^2  \cdot  | \Lie_{Z^K}  A | }{(1+t+|q|)^{2-      c (\gamma)  \cdot c (\delta)  \cdot c(|I|) \cdot E (\lfloor \frac{|I|}{2} \rfloor  + 4) \cdot \eps } \cdot (1+|q|)^{1+\gamma - 2\de }}    \\
               \notag
                &&+      C(q_0)   \cdot  c (\delta) \cdot c (\gamma) \cdot C(|I|) \cdot E ( \lfloor \frac{|I|}{2} \rfloor + 4)  \cdot \frac{\eps^2 \cdot  |  \Lie_{Z^K} A  | }{(1+t+|q|)^{2-      c (\gamma)  \cdot c (\delta)  \cdot c(|I|) \cdot E (\lfloor \frac{|I|}{2} \rfloor  + 4) \cdot \eps } \cdot (1+|q|)^{2\gamma - 4\de }}     \\
                \notag
                  && +  O \big( C(q_0)   \cdot  c (\delta) \cdot c (\gamma) \cdot C(|I|) \cdot E ( \lfloor \frac{|I|}{2} \rfloor  + 4)  \cdot \frac{\eps^2 \cdot  | \Lie_{Z^K}  h | }{(1+t+|q|)^{2-      c (\gamma)  \cdot c (\delta)  \cdot c(|I|) \cdot E ( \lfloor \frac{|I|}{2} \rfloor + 4) \cdot \eps } \cdot (1+|q|)^{2+\gamma - 2\de }}  \big)   \\
                  && +  O \big(    C(q_0)   \cdot  c (\delta) \cdot c (\gamma) \cdot C(|I|) \cdot E ( \lfloor \frac{|I|}{2} \rfloor  + 4)  \cdot \frac{\eps^2 \cdot  | \derm ( \Lie_{Z^K}  h )  |   }{(1+t+|q|)^{2-      c (\gamma)  \cdot c (\delta)  \cdot c(|I|) \cdot E ( \lfloor \frac{|I|}{2} \rfloor + 4) \cdot \eps } \cdot (1+|q|)^{1+\gamma - 2\de }} \big) \\
                  && +  O \big(  C(q_0)   \cdot  c (\delta) \cdot c (\gamma) \cdot C(|I|) \cdot E (  \lfloor \frac{|I|}{2} \rfloor   +4)  \cdot \frac{\eps^2  \cdot | \derm ( \Lie_{Z^K}  A) |  }{(1+t+|q|)^{2-   c (\gamma)  \cdot c (\delta)  \cdot c(|I|) \cdot E ( \lfloor \frac{|I|}{2} \rfloor + 4)\cdot  \eps } \cdot (1+|q|)}  \big) \\
                  \notag
           &&       + O \big( C(q_0)   \cdot  c (\delta) \cdot c (\gamma) \cdot C(|I|) \cdot E (\lfloor \frac{|I|}{2} \rfloor + 4)  \cdot \frac{\eps^2 \cdot | \Lie_{Z^K} h | }{(1+t+|q|)^{2-      c (\gamma)  \cdot c (\delta)  \cdot c(|I|) \cdot E (\lfloor \frac{|I|}{2} \rfloor + 4) \cdot \eps } \cdot (1+|q|)^{1+2\gamma - 4\de }}     \big)  \\
           && + O \big(  C(q_0)   \cdot  c (\delta) \cdot c (\gamma) \cdot C(|I|) \cdot E ( \lfloor \frac{|I|}{2} \rfloor + 4)  \cdot \frac{\eps^2  \cdot | \Lie_{Z^K} A|}{(1+t+|q|)^{2-      c (\gamma)  \cdot c (\delta)  \cdot c(|I|) \cdot E ( \lfloor \frac{|I|}{2} \rfloor + 4) \cdot \eps } \cdot (1+|q|)^{1+\gamma - 2\de }}    \big)      \\
           && + O  \big(   C(q_0)   \cdot  c (\delta) \cdot c (\gamma) \cdot C(|I|) \cdot E (\lfloor \frac{|I|}{2} \rfloor  + 4)  \cdot \frac{\eps^2  \cdot |\derm ( \Lie_{Z^K} A ) | }{(1+t+|q|)^{2-      c (\gamma)  \cdot c (\delta)  \cdot c(|I|) \cdot E ( \lfloor \frac{|I|}{2} \rfloor + 4) \cdot \eps } \cdot (1+|q|)^{\gamma - 2\de }}  \big) \; \Big]  \\
\notag
&& +      |  \Lie_{Z^I}    \big(  A_L   \cdot     \derm A    \big)  | +|  \Lie_{Z^I} \big( A_{e_a}  \cdot     \derm A_{e_a} \big)|
                \eeaa

Grouping the terms together, and using the fact that since $\ga \geq 3 \de$\,, we have $\ga -2\de \geq \de \geq 0$\,, we get the following factors for each term.

    \textbf{The factor for $  | \derm ( \Lie_{Z^I} h ) | $:}\\
    
    We estimate
\beaa
  \notag
&&  C(q_0)   \cdot  c (\delta) \cdot c (\gamma) \cdot C(|I|) \cdot E ( \lfloor \frac{|I|}{2} \rfloor + 5)  \cdot \frac{\eps \cdot  | \derm ( \Lie_{Z^I} h ) | }{(1+t+|q|)^{2-      c (\gamma)  \cdot c (\delta)  \cdot c(|I|) \cdot E ( \lfloor \frac{|I|}{2} \rfloor + 5) \cdot \eps } \cdot (1+|q|)^{\gamma - 2\de }}   \\
  &&          +    C(q_0)   \cdot  c (\delta) \cdot c (\gamma) \cdot C(| I |) \cdot E (\lfloor \frac{|I|}{2} \rfloor  + 4)  \cdot \frac{\eps^2 \cdot   | \derm ( \Lie_{Z^I} h )    |  }{(1+t+|q|)^{2-      c (\gamma)  \cdot c (\delta)  \cdot c(|I|) \cdot E ( \lfloor \frac{|I|}{2} \rfloor  + 4) \cdot \eps } \cdot (1+|q|)^{2\gamma - 4\de }}    \\
                    && +  O \big(    C(q_0)   \cdot  c (\delta) \cdot c (\gamma) \cdot C(|I|) \cdot E ( \lfloor \frac{|I|}{2} \rfloor  + 4)  \cdot \frac{\eps^2  \cdot  | \derm ( \Lie_{Z^I}  h )  |   }{(1+t+|q|)^{2-      c (\gamma)  \cdot c (\delta)  \cdot c(|I|) \cdot E ( \lfloor \frac{|I|}{2} \rfloor + 4) \cdot \eps } \cdot (1+|q|)^{1+\gamma - 2\de }} \big) \\
     &\les&               O \big(   C(q_0)   \cdot  c (\delta) \cdot c (\gamma) \cdot C(|I|) \cdot E ( \lfloor \frac{|I|}{2} \rfloor  + 5)  \cdot \frac{\eps    \cdot  | \derm ( \Lie_{Z^I} h ) | }{(1+t+|q|)^{2-      c (\gamma)  \cdot c (\delta)  \cdot c(|I|) \cdot E ( \lfloor \frac{|I|}{2} \rfloor + 5) \cdot \eps } \cdot (1+|q|)^{\gamma - 2\de }}   \big)\; .
  \eeaa

      \textbf{The factor for $ |\derm ( \Lie_{Z^I} A ) | $:}\\

       We compute
\beaa
   \notag
    &&     C(q_0)   \cdot  c (\delta) \cdot c (\gamma) \cdot C(|I|) \cdot E ( \lfloor \frac{|I|}{2} \rfloor   +5)  \cdot \frac{\eps    \cdot  |\derm ( \Lie_{Z^I} A ) | }{(1+t+|q|)^{2-   c (\gamma)  \cdot c (\delta)  \cdot c(|I|) \cdot E ( \lfloor \frac{|I|}{2} \rfloor + 5)\cdot  \eps } }  \\
       && +  O \big(  C(q_0)   \cdot  c (\delta) \cdot c (\gamma) \cdot C(|I|) \cdot E (  \lfloor \frac{|I|}{2} \rfloor  +4)  \cdot \frac{\eps^2  \cdot | \derm ( \Lie_{Z^I}  A) |  }{(1+t+|q|)^{2-   c (\gamma)  \cdot c (\delta)  \cdot c(|I|) \cdot E (\lfloor \frac{|I|}{2} \rfloor + 4)\cdot  \eps } \cdot (1+|q|)}  \big) \\
       && + O  \big(   C(q_0)   \cdot  c (\delta) \cdot c (\gamma) \cdot C(|I|) \cdot E (\lfloor \frac{|I|}{2} \rfloor  + 4)  \cdot \frac{\eps^2 \cdot |\derm ( \Lie_{Z^I} A ) | }{(1+t+|q|)^{2-      c (\gamma)  \cdot c (\delta)  \cdot c(|I|) \cdot E ( \lfloor \frac{|I|}{2} \rfloor + 4) \cdot \eps } \cdot (1+|q|)^{\gamma - 2\de }}  \big)  \\
\notag
    &\les&  O \big(    C(q_0)   \cdot  c (\delta) \cdot c (\gamma) \cdot C(|I|) \cdot E (  \lfloor \frac{|I|}{2} \rfloor   +5)  \cdot \frac{\eps  \cdot  |\derm ( \Lie_{Z^I} A ) | }{(1+t+|q|)^{2-   c (\gamma)  \cdot c (\delta)  \cdot c(|I|) \cdot E ( \lfloor \frac{|I|}{2} \rfloor + 5)\cdot  \eps } }  \big)  \; .  
    \eeaa

           \textbf{The factor for $ |\rderm ( \Lie_{Z^I} A ) |$:}\\
We have for $\ga \geq 3\ga $\;, and therefore for $\ga - 2 \de \geq    \de$\;, that 
  \bea
   \notag
    &&  C(q_0)   \cdot  c (\delta) \cdot c (\gamma) \cdot C(|I|) \cdot E (  \lfloor \frac{|I|}{2} \rfloor   +4)  \cdot \frac{\eps \cdot  |\rderm ( \Lie_{Z^I} A ) |}{(1+t+|q|)^{1-   c (\gamma)  \cdot c (\delta)  \cdot c(|I|) \cdot E (\lfloor \frac{|I|}{2} \rfloor + 4)\cdot  \eps } \cdot (1+|q|)}       \\
       \notag
               && +     C(q_0)   \cdot  c (\delta) \cdot c (\gamma) \cdot C(|I|) \cdot E (\lfloor \frac{|I|}{2} \rfloor + 4)  \cdot \frac{\eps \cdot    | \rderm ( \Lie_{Z^I} A )  | }{(1+t+|q|)^{1-      c (\gamma)  \cdot c (\delta)  \cdot c(|I|) \cdot E (\lfloor \frac{|I|}{2} \rfloor + 4) \cdot \eps } \cdot (1+|q|)^{\gamma - 2\de }}     \\
                  \notag
         &\les&         C(q_0)   \cdot  c (\delta) \cdot c (\gamma) \cdot C(|I|) \cdot E ( \lfloor \frac{|I|}{2} \rfloor + 4)  \cdot \frac{\eps \cdot    | \rderm ( \Lie_{Z^I} A )  | }{(1+t+|q|)  \cdot (1+|q|)^{ \de -      c (\gamma)  \cdot c (\delta)  \cdot c(|I|) \cdot E ( \lfloor \frac{|I|}{2} \rfloor + 4) \cdot \eps}  }     \\
         \notag
         &\les&         C(q_0)   \cdot  c (\delta) \cdot c (\gamma) \cdot C(|I|) \cdot E ( \lfloor \frac{|I|}{2} \rfloor + 4)  \cdot \frac{\eps \cdot    | \rderm ( \Lie_{Z^I} A )  | }{(1+t+|q|)  \cdot (1+|q|)^{ \de -      c (\gamma)  \cdot c (\delta)  \cdot c(|I|) \ \cdot \eps}  }     \\
         \notag
         && \text{(since $E ( \lfloor \frac{|I|}{2} \rfloor + 4) \leq 1$ ).} \; \\
      \eea 
      Thus, for $\eps$ small enough depending on $\ga$\;, $\de$\;, and $|I|$\;, we have
      \bea
      \de -      c (\gamma)  \cdot c (\delta)  \cdot c(|I|) \ \cdot \eps &\geq& 0 \; ,
      \eea
       and therefore,
        \bea\label{estimateonthetermthatcontainstangentialderivativeofAinsourcetermsforwaveopertaorinAwithdecayinqfactor}
   \notag
    &&  C(q_0)   \cdot  c (\delta) \cdot c (\gamma) \cdot C(|I|) \cdot E (  \lfloor \frac{|I|}{2} \rfloor   +4)  \cdot \frac{\eps \cdot  |\rderm ( \Lie_{Z^I} A ) |}{(1+t+|q|)^{1-   c (\gamma)  \cdot c (\delta)  \cdot c(|I|) \cdot E (\lfloor \frac{|I|}{2} \rfloor + 4)\cdot  \eps } \cdot (1+|q|)}       \\
       \notag
               && +     C(q_0)   \cdot  c (\delta) \cdot c (\gamma) \cdot C(|I|) \cdot E (\lfloor \frac{|I|}{2} \rfloor + 4)  \cdot \frac{\eps \cdot    | \rderm ( \Lie_{Z^I} A )  | }{(1+t+|q|)^{1-      c (\gamma)  \cdot c (\delta)  \cdot c(|I|) \cdot E (\lfloor \frac{|I|}{2} \rfloor + 4) \cdot \eps } \cdot (1+|q|)^{\gamma - 2\de }}     \\
                  \notag
         &\les&         C(q_0)   \cdot  c (\delta) \cdot c (\gamma) \cdot C(|I|) \cdot E ( \lfloor \frac{|I|}{2} \rfloor + 4)  \cdot \frac{\eps \cdot    | \rderm ( \Lie_{Z^I} A )  | }{(1+t+|q|)    }     \\
      \eea

         \textbf{The factor for $ | \rderm  ( \Lie_{Z^I} h ) |$:}\\
         
         It is
        \beaa
   \notag
        &&         C(q_0)   \cdot  c (\delta) \cdot c (\gamma) \cdot C(|I|) \cdot E (\lfloor \frac{|I|}{2} \rfloor  + 4)  \cdot \frac{\eps  \cdot   | \rderm  ( \Lie_{Z^I} h ) | }{(1+t+|q|)^{1-      c (\gamma)  \cdot c (\delta)  \cdot c(|I|) \cdot E (\lfloor \frac{|I|}{2} \rfloor + 4) \cdot \eps } \cdot (1+|q|)^{1+\gamma - 2\de }}    \; .
       \eeaa
  
     \textbf{The factors for $| \Lie_{Z^I} A |$:}\\
     
     We estimate,
          \beaa
   \notag
             &&     C(q_0)   \cdot  c (\delta) \cdot c (\gamma) \cdot C(|I|) \cdot E ( \lfloor \frac{|I|}{2} \rfloor  + 5)  \cdot \frac{\eps  \cdot |  \Lie_{Z^I} A  | }{(1+t+|q|)^{2-      c (\gamma)  \cdot c (\delta)  \cdot c(|I|) \cdot E ( \lfloor \frac{|I|}{2} \rfloor  + 5) \cdot \eps } \cdot (1+|q|)^{\gamma - 2\de }}    \\
             \notag
              &&  + C(q_0)   \cdot  c (\delta) \cdot c (\gamma) \cdot C(|I|) \cdot E (\lfloor \frac{|I|}{2} \rfloor  + 4)  \cdot \frac{\eps^2  \cdot  | \Lie_{Z^I}  A | }{(1+t+|q|)^{2-      c (\gamma)  \cdot c (\delta)  \cdot c(|I|) \cdot E (\lfloor \frac{|I|}{2} \rfloor  + 4) \cdot \eps } \cdot (1+|q|)^{1+\gamma - 2\de }}    \\
               \notag
                &&+      C(q_0)   \cdot  c (\delta) \cdot c (\gamma) \cdot C(|I|) \cdot E ( \lfloor \frac{|I|}{2} \rfloor  + 4)  \cdot \frac{\eps^2 \cdot  |  \Lie_{Z^I} A  | }{(1+t+|q|)^{2-      c (\gamma)  \cdot c (\delta)  \cdot c(|I|) \cdot E ( \lfloor \frac{|I|}{2} \rfloor  + 4) \cdot \eps } \cdot (1+|q|)^{2\gamma - 4\de }}     \\
                \notag
        && + O \big(  C(q_0)   \cdot  c (\delta) \cdot c (\gamma) \cdot C(|I|) \cdot E ( \lfloor \frac{|I|}{2} \rfloor  + 4)  \cdot \frac{\eps^2  \cdot | \Lie_{Z^I} A|}{(1+t+|q|)^{2-      c (\gamma)  \cdot c (\delta)  \cdot c(|I|) \cdot E ( \lfloor \frac{|I|}{2} \rfloor + 4) \cdot \eps } \cdot (1+|q|)^{1+\gamma - 2\de }}    \big)      \\
 \notag
             &\les&  O \big(     C(q_0)   \cdot  c (\delta) \cdot c (\gamma) \cdot C(|I|) \cdot E ( \lfloor \frac{|I|}{2} \rfloor + 5)  \cdot \frac{\eps \cdot |  \Lie_{Z^I} A  | }{(1+t+|q|)^{2-      c (\gamma)  \cdot c (\delta)  \cdot c(|I|) \cdot E ( \lfloor \frac{|I|}{2} \rfloor + 5) \cdot \eps } \cdot (1+|q|)^{\gamma - 2\de }}   \big)  \;.
                    \eeaa

        \textbf{The factors for $|\Lie_{Z^I} h |$:}\\
We have
                         \beaa
                                      &&   O \big( C(q_0)   \cdot  c (\delta) \cdot c (\gamma) \cdot C(|I|) \cdot E ( \lfloor \frac{|I|}{2} \rfloor  + 4)  \cdot \frac{\eps^2  \cdot  | \Lie_{Z^I}  h | }{(1+t+|q|)^{2-      c (\gamma)  \cdot c (\delta)  \cdot c(|I|) \cdot E ( \lfloor \frac{|I|}{2} \rfloor  + 4) \cdot \eps } \cdot (1+|q|)^{2+\gamma - 2\de }}  \big)   \\
                                      \notag
                                                 &&       + O \big( C(q_0)   \cdot  c (\delta) \cdot c (\gamma) \cdot C(|I|) \cdot E (\lfloor \frac{|I|}{2} \rfloor  + 4)  \cdot \frac{\eps^2 \cdot | \Lie_{Z^I} h | }{(1+t+|q|)^{2-      c (\gamma)  \cdot c (\delta)  \cdot c(|I|) \cdot E ( \lfloor \frac{|I|}{2} \rfloor + 4) \cdot \eps } \cdot (1+|q|)^{1+2\gamma - 4\de }}     \big)  \\
                                                 &\les&O \big( C(q_0)   \cdot  c (\delta) \cdot c (\gamma) \cdot C(|I|) \cdot E ( \lfloor \frac{|I|}{2} \rfloor  + 4)  \cdot \frac{\eps  \cdot | \Lie_{Z^I} h | }{(1+t+|q|)^{2-      c (\gamma)  \cdot c (\delta)  \cdot c(|I|) \cdot E (\lfloor \frac{|I|}{2} \rfloor + 4) \cdot \eps } \cdot (1+|q|)^{1+2\gamma - 4\de }}     \big) \; .
                                                      \eeaa

                \textbf{The result:}\\
                                                         
Putting all together, we get that

                                                       \beaa
   \notag
&&  |  \Lie_{Z^I}  g^{\la\mu} \derm_{\la}   \derm_{\mu}   A_{V}  |  \\
   &\les&           \sum_{|K| \leq |I |}  \Big[   \;       O \big(   C(q_0)   \cdot  c (\delta) \cdot c (\gamma) \cdot C(|I|) \cdot E ( \lfloor \frac{|I|}{2} \rfloor + 5)  \cdot \frac{\eps   \cdot  | \derm ( \Lie_{Z^K} h ) | }{(1+t+|q|)^{2-      c (\gamma)  \cdot c (\delta)  \cdot c(|I|) \cdot E ( \lfloor \frac{|I|}{2} \rfloor + 5) \cdot \eps } \cdot (1+|q|)^{\gamma - 2\de }}   \big) \\
   && + O \big(    C(q_0)   \cdot  c (\delta) \cdot c (\gamma) \cdot C(|I|) \cdot E ( \lfloor \frac{|I|}{2} \rfloor   +5)  \cdot \frac{\eps  \cdot  |\derm ( \Lie_{Z^K} A ) | }{(1+t+|q|)^{2-   c (\gamma)  \cdot c (\delta)  \cdot c(|I|) \cdot E ( \lfloor \frac{|I|}{2} \rfloor + 5)\cdot  \eps } }  \big) \\
   &&  +  C(q_0)   \cdot  c (\delta) \cdot c (\gamma) \cdot C(|I|) \cdot E ( \lfloor \frac{|I|}{2} \rfloor + 4)  \cdot \frac{\eps   \cdot    | \rderm ( \Lie_{Z^K} A )  | }{(1+t+|q|)^{1-      c (\gamma)  \cdot c (\delta)  \cdot c(|I|) \cdot E (\lfloor \frac{|I|}{2} \rfloor + 4) \cdot \eps } \cdot ( 1 + |q| )^{ \frac{1}{2}}  }    \\
      \notag
        &&    +     C(q_0)   \cdot  c (\delta) \cdot c (\gamma) \cdot C(|I|) \cdot E ( \lfloor \frac{|I|}{2} \rfloor + 4)  \cdot \frac{\eps  \cdot   | \rderm  ( \Lie_{Z^K} h ) | }{(1+t+|q|)^{1-      c (\gamma)  \cdot c (\delta)  \cdot c(|I|) \cdot E ( \lfloor \frac{|I|}{2} \rfloor  + 4) \cdot \eps } \cdot (1+|q|)^{1+\gamma - 2\de }}    \\
      &&    +       O \big(     C(q_0)   \cdot  c (\delta) \cdot c (\gamma) \cdot C(|I|) \cdot E ( \lfloor \frac{|I|}{2} \rfloor  + 5)  \cdot \frac{\eps \cdot |  \Lie_{Z^K} A  | }{(1+t+|q|)^{2-      c (\gamma)  \cdot c (\delta)  \cdot c(|I|) \cdot E ( \lfloor \frac{|I|}{2} \rfloor + 5) \cdot \eps } \cdot (1+|q|)^{\gamma - 2\de }}   \\
         &&    +         O \big( C(q_0)   \cdot  c (\delta) \cdot c (\gamma) \cdot C(|I|) \cdot E ( \lfloor \frac{|I|}{2} \rfloor  + 4)  \cdot \frac{\eps \cdot | \Lie_{Z^K} h | }{(1+t+|q|)^{2-      c (\gamma)  \cdot c (\delta)  \cdot c(|I|) \cdot E ( \lfloor \frac{|I|}{2} \rfloor + 4) \cdot \eps } \cdot (1+|q|)^{1+2\gamma - 4\de }}     \big) \; \Big] \\
         \notag
         && +     |   \Lie_{Z^I}    \big(  A_L   \cdot     \derm A    \big)  |   +  |  \Lie_{Z^I} \big( A_{e_a}  \cdot     \derm A_{e_a} \big)  |  \; ,
\eeaa

Now, since the energy that we want to control is on $h^1 = h - h^0$, and not on $h$ (we recall that we showed that the energy for $h^0$ is infinite), we want to obtain estimates involving only $h^1$\,. For this, we estimate
\beaa
| \derm ( \Lie_{Z^I} h )  | &\leq &  | \derm ( \Lie_{Z^I} h^1 )  | + | \derm ( \Lie_{Z^I} h^0 )  |
 \eeaa

Now, using Lemma \ref{Liederivativesofsphericalsymmetricpart}, we have for all $I$\,, the following estimates on $h^0$ which do not involve the bootstrap assumptions and therefore, theses estimates on $h^0$ can be used for any number of Lie derivatives on $h^0$ without ruining our bootstrap argument, 
 \beaa
 \notag
|   \Lie_{Z^I} h^0  |   &\leq& c (\delta) \cdot C ( |I| )   \cdot \frac{M }{(1+t+|q|) }  \; , \\
 \notag
|\derm (  \Lie_{Z^I} h^0 )   |   &\leq&  C ( |I| )   \cdot \frac{M }{(1+t+|q|)^{2} } \; .
      \eeaa
   Hence, for $M \leq \eps $\,, 
    \beaa
 \notag
|   \Lie_{Z^I} h  |   &\leq& |   \Lie_{Z^I} h^1  | + c (\delta) \cdot C ( |I| )   \cdot \frac{\eps }{(1+t+|q|) }  \; , \\
 \notag
|\derm (  \Lie_{Z^I} h)  |   &\leq& |\derm (  \Lie_{Z^I} h^1)  |  +  C ( |I| )   \cdot \frac{\eps }{(1+t+|q|)^{2} } \; .
      \eeaa
      
   Injecting, we obtain
                                                          \beaa
   \notag
&&  |  \Lie_{Z^I}  g^{\la\mu} \derm_{\la}   \derm_{\mu}   A_{V}  |  \\
   &\les&           \sum_{|K| \leq |I |}  \Big[   \;       O \big(   C(q_0)   \cdot  c (\delta) \cdot c (\gamma) \cdot C(|I|) \cdot E ( \lfloor \frac{|I|}{2} \rfloor + 5)  \cdot \frac{\eps   \cdot  | \derm ( \Lie_{Z^K} h^1 ) | }{(1+t+|q|)^{2-      c (\gamma)  \cdot c (\delta)  \cdot c(|I|) \cdot E ( \lfloor \frac{|I|}{2} \rfloor + 5) \cdot \eps } \cdot (1+|q|)^{\gamma - 2\de }}   \big) \\
   && + O \big(    C(q_0)   \cdot  c (\delta) \cdot c (\gamma) \cdot C(|I|) \cdot E ( \lfloor \frac{|I|}{2} \rfloor   +5)  \cdot \frac{\eps  \cdot  |\derm ( \Lie_{Z^K} A ) | }{(1+t+|q|)^{2-   c (\gamma)  \cdot c (\delta)  \cdot c(|I|) \cdot E ( \lfloor \frac{|I|}{2} \rfloor + 5)\cdot  \eps } }  \big) \\
   &&  +  C(q_0)   \cdot  c (\delta) \cdot c (\gamma) \cdot C(|I|) \cdot E ( \lfloor \frac{|I|}{2} \rfloor + 4)  \cdot \frac{\eps   \cdot    | \rderm ( \Lie_{Z^K} A )  | }{(1+t+|q|)^{1-      c (\gamma)  \cdot c (\delta)  \cdot c(|I|) \cdot E (\lfloor \frac{|I|}{2} \rfloor + 4) \cdot \eps } \cdot ( 1 + |q| )^{ \frac{1}{2}}  }    \\
      \notag
        &&    +     C(q_0)   \cdot  c (\delta) \cdot c (\gamma) \cdot C(|I|) \cdot E ( \lfloor \frac{|I|}{2} \rfloor + 4)  \cdot \frac{\eps  \cdot   | \rderm  ( \Lie_{Z^K} h^1 ) | }{(1+t+|q|)^{1-      c (\gamma)  \cdot c (\delta)  \cdot c(|I|) \cdot E ( \lfloor \frac{|I|}{2} \rfloor  + 4) \cdot \eps } \cdot (1+|q|)^{1+\gamma - 2\de }}    \\
      &&    +       O \big(     C(q_0)   \cdot  c (\delta) \cdot c (\gamma) \cdot C(|I|) \cdot E ( \lfloor \frac{|I|}{2} \rfloor  + 5)  \cdot \frac{\eps \cdot |  \Lie_{Z^K} A  | }{(1+t+|q|)^{2-      c (\gamma)  \cdot c (\delta)  \cdot c(|I|) \cdot E ( \lfloor \frac{|I|}{2} \rfloor + 5) \cdot \eps } \cdot (1+|q|)^{\gamma - 2\de }}   \\
         &&    +         O \big( C(q_0)   \cdot  c (\delta) \cdot c (\gamma) \cdot C(|I|) \cdot E ( \lfloor \frac{|I|}{2} \rfloor  + 4)  \cdot \frac{\eps \cdot | \Lie_{Z^K} h^1 | }{(1+t+|q|)^{2-      c (\gamma)  \cdot c (\delta)  \cdot c(|I|) \cdot E ( \lfloor \frac{|I|}{2} \rfloor + 4) \cdot \eps } \cdot (1+|q|)^{1+2\gamma - 4\de }}     \big) \; \Big] \\
         \notag
         && +     |   \Lie_{Z^I}    \big(  A_L   \cdot     \derm A    \big)  |   +  |  \Lie_{Z^I} \big( A_{e_a}  \cdot     \derm A_{e_a} \big)  |  \; 
\eeaa
\beaa
\notag
&&  +  O \big(   C(q_0)   \cdot  c (\delta) \cdot c (\gamma) \cdot C(|I|) \cdot E ( \lfloor \frac{|I|}{2} \rfloor + 5)  \cdot \frac{\eps^2 }{(1+t+|q|)^{4-      c (\gamma)  \cdot c (\delta)  \cdot c(|I|) \cdot E ( \lfloor \frac{|I|}{2} \rfloor + 5) \cdot \eps } \cdot (1+|q|)^{\gamma - 2\de }}   \big)  \\
      \notag
          &&    +     C(q_0)   \cdot  c (\delta) \cdot c (\gamma) \cdot C(|I|) \cdot E ( \lfloor \frac{|I|}{2} \rfloor + 4)  \cdot \frac{\eps^2 }{(1+t+|q|)^{3-      c (\gamma)  \cdot c (\delta)  \cdot c(|I|) \cdot E ( \lfloor \frac{|I|}{2} \rfloor  + 4) \cdot \eps } \cdot (1+|q|)^{1+\gamma - 2\de }}    \\
          \notag
             &&    +         O \big( C(q_0)   \cdot  c (\delta) \cdot c (\gamma) \cdot C(|I|) \cdot E ( \lfloor \frac{|I|}{2} \rfloor  + 4)  \cdot \frac{\eps^2  }{(1+t+|q|)^{3-      c (\gamma)  \cdot c (\delta)  \cdot c(|I|) \cdot E ( \lfloor \frac{|I|}{2} \rfloor + 4) \cdot \eps } \cdot (1+|q|)^{1+2\gamma - 4\de }}     \big) \; .
      \eeaa
      
      Thus, we obtain the desired result for $\ga \geq 3 \de $\;, for $0 < \de \leq \frac{1}{4}$\;, and for $\eps$ small enough depending on $\ga$\;, $\de$\;, and $|I|$\;.

\end{proof}

\subsection{Structure of the “bad” term $\Lie_{Z^I}    \big(  A_L   \cdot     \derm A    \big) $ for the Einstein-Yang-Mills potential}\

Now, we would like to study the behaviour of one of the remaining so-called “bad” term  $\Lie_{Z^I}    \big(  A_L   \cdot     \derm A    \big) $\,.

First, as we have always done earlier, we estimate, using also Lemmas \ref{LiederivativeZofMinkwoskimetric} and \ref{LiederivativeZcommuteswithminkowksicovariantderivative},
\beaa
&& | \Lie_{Z^I}    \big(  A_L   \cdot     \derm A      \big)  | \\
&\les& \sum_{|J| + |K| \leq |I| }  | \Lie_{Z^J}  A_L |  \cdot     |  \derm ( \Lie_{Z^K} A )  |  
\eeaa

Now, since it is a “bad” term, we will decompose the following sum in way that we can extract from it a term that is not so “bad” or in fact good as it will turn out, and we will deal with the rest later. We decompose the sum as follows,
\bea\label{decompositionofsumfortroublesometermsforsourcesonA}
\notag
&& | \Lie_{Z^I}    \big(  A_L   \cdot     \derm A    \big)    | \\
\notag
&\les& \sum_{|K| = |I| }  \Big(  |  A_L |  \cdot     |  \derm ( \Lie_{Z^K} A )  |  +  |  \Lie_{Z^K} A_L |  \cdot     |  \derm A   |   \Big) \\
\notag
&& +  \sum_{|J| + |K| \leq |I| -1}    | \Lie_{Z^J}  A_L |  \cdot     |  \derm ( \Lie_{Z^K} A )  |     \; . \\
\eea

We notice that there are terms that are not that “bad” after all, since the Lie derivatives for these special components behave worse than the zeroth Lie derivative. Exploiting the very \textit{special} structure of these components, we shall show the the terms that involve zeroth Lie derivative are in fact “good” terms. We shall explain with the following lemmas.

\begin{lemma}\label{simpleestimateonlyonthe ALcomponenetwithOUTLiederivative}
We have in the exterior region,
          \beaa
        \notag
 |  A_{L}  |   &\les&  C(q_0)   \cdot  c (\delta) \cdot c (\gamma)  \cdot E ( 3) \cdot \frac{ \eps   }{ (1+t+|q|)^{2-2\delta}  \cdot  (1+|q|)^{\ga - 1}  } \; .\\
       \eeaa
        \end{lemma}
  
   \begin{proof}  
By taking $|I| = 0$ in Lemma \ref{estimategoodcomponentspotentialandmetric}, and restricting the estimate to the exterior region, we get the stated result.
\end{proof}

\begin{corollary} \label{TheleeadingtermfortroublesomecomponentsintehsourcesforA}
We have for $\ga \geq 3 \de $\,, and $0 < \de \leq \frac{1}{4}$,
\beaa
&& \sum_{|K| = |I| }  \Big(  |  A_L |  \cdot     |  \derm ( \Lie_{Z^K} A )  |  +  |  \Lie_{Z^K} A_L |  \cdot     |  \derm A   |   \Big) \\
&\les& \sum_{|K| = |I| }  \Big( C(q_0)   \cdot  c (\delta) \cdot c (\gamma)  \cdot E ( 3) \cdot \frac{ \eps  \cdot     |  \derm ( \Lie_{Z^K} A )  |   }{ (1+t+|q|)^{2-2\delta}  \cdot  (1+|q|)^{\ga - 1}  } \\
&&  +    C(q_0)   \cdot  c (\delta) \cdot c (\gamma)  \cdot E ( 4)  \cdot \frac{\eps \cdot   |  \Lie_{Z^K} A_L |   }{(1+t+|q|)^{1-      c (\gamma)  \cdot c (\delta)   \cdot E ( 4) \cdot \eps } \cdot (1+|q|)^{1+\gamma - 2\de }}     \Big) 
\eeaa

\end{corollary}
  \begin{proof}
As a result of Lemma \ref{simpleestimateonlyonthe ALcomponenetwithOUTLiederivative}, and injecting in the “good” sum in \eqref{decompositionofsumfortroublesometermsforsourcesonA}, we get

\beaa
&& | \Lie_{Z^I}    \big(  A_L   \cdot     \derm A    \big)    | \\
&\les& \sum_{|K| = |I| }  \Big(  |  A_L |  \cdot     |  \derm ( \Lie_{Z^K} A )  |  +  |  \Lie_{Z^K} A_L |  \cdot     |  \derm A   |   \Big) \\
&\les& \sum_{|K| = |I| }  \Big(C(q_0)   \cdot  c (\delta) \cdot c (\gamma)  \cdot E ( 3) \cdot \frac{ \eps  \cdot     |  \derm ( \Lie_{Z^K} A )  |   }{ (1+t+|q|)^{2-2\delta}  \cdot  (1+|q|)^{\ga - 1}  } +  |  \Lie_{Z^K} A_L |  \cdot     |  \derm A   |    \Big) 
\eeaa

There still remains a term  $|  \Lie_{Z^K} A_L |  \cdot     |  \derm A   | $\,; we estimate it using \eqref{formulaforinductionhypothesisonA} in Lemma \ref{upgradedestimateonLiederivativesoffields}, without forgetting that it is an $A_L$ component in the product. We obtain for $\ga \geq 3 \de $\,, and $0 < \de \leq \frac{1}{4}$\,,
           
           \beaa
 && \sum_{|K| = |I| }             |  \Lie_{Z^K} A_L |  \cdot     |  \derm A   | \\
  &\les& \sum_{|K| = |I| }             C(q_0)   \cdot  c (\delta) \cdot c (\gamma)  \cdot E ( 4)  \cdot \frac{\eps \cdot   |  \Lie_{Z^K} A_L |   }{(1+t+|q|)^{1-      c (\gamma)  \cdot c (\delta)   \cdot E ( 4) \cdot \eps } \cdot (1+|q|)^{1+\gamma - 2\de }}  \; .
           \eeaa
           \end{proof}

\begin{lemma}  \label{estimateonlowerordertermsforthepotentialAinthesourcetermsforenergy}
We have for $\ga \geq 3 \de $\,, and $0 < \de \leq \frac{1}{4}$\,,
               \beaa
&& \sum_{|J| + |K| \leq |I| -1}    | \Lie_{Z^J}  A_L |  \cdot     |  \derm ( \Lie_{Z^K} A )  |     \\
&\les& C(q_0)   \cdot  c (\delta) \cdot c (\gamma) \cdot C( \lfloor \frac{|I|-1}{2} \rfloor) \cdot E ( \lfloor \frac{|I|-1}{2} \rfloor + 4) \\
&& \times \sum_{ |K| \leq |I| -1}    \Big(    \frac{\eps  \cdot     |  \derm ( \Lie_{Z^K} A )  | }{(1+t+|q|)^{1-      c (\gamma)  \cdot c (\delta)  \cdot c( \lfloor \frac{|I|-1}{2} \rfloor) \cdot E ( \lfloor \frac{|I|-1}{2} \rfloor+ 4) \cdot \eps } \cdot (1+|q|)^{\gamma - 2\de }}    \\
&& + \frac{\eps \cdot  | \Lie_{Z^K}  A |  }{(1+t+|q|)^{1-      c (\gamma)  \cdot c (\delta)  \cdot c(\lfloor \frac{|I|-1}{2} \rfloor ) \cdot E (\lfloor \frac{|I|-1}{2} \rfloor + 4) \cdot \eps } \cdot (1+|q|)^{1+\gamma - 2\de }}   \Big) \;. 
\eeaa       

\end{lemma}  
\begin{proof}

We decompose the sum as follows,

\beaa
&& \sum_{|J| + |K| \leq |I| -1}    | \Lie_{Z^J}  A_L |  \cdot     |  \derm ( \Lie_{Z^K} A )  |     \\
&\leq& \sum_{ |J| \leq   \lfloor \frac{|I|-1}{2} \rfloor   ,\; |K| \leq |I| -1}   | \Lie_{Z^J}  A_L |  \cdot     |  \derm ( \Lie_{Z^K} A )  |     \\
&& +  \sum_{ |K| \leq   \lfloor \frac{|I|-1}{2} \rfloor    ,\; |J| \leq |I| -1}    | \Lie_{Z^J}  A_L |  \cdot     |  \derm ( \Lie_{Z^K} A )  |  
\eeaa
However, we know from Lemma \ref{formulaforinductionhypothesisonA} and Lemma \ref{upgradedestimateonLiederivativesoffieldswithoutgradiant}, that for $\ga \geq 3 \de $\,, and $0 < \de \leq \frac{1}{4}$\,, and for any $|K| \leq   \lfloor \frac{|I|-1}{2} \rfloor $, we have

            \beaa
 \notag
&& |\derm  ( \Lie_{Z^K} A) (t,x)  |    \\
\notag
&\leq&   C(q_0)   \cdot  c (\delta) \cdot c (\gamma) \cdot C(\lfloor \frac{|I|-1}{2} \rfloor ) \cdot E ( \lfloor \frac{|I|-1}{2} \rfloor  + 4)  \\
\notag
&& \times  \frac{\eps }{(1+t+|q|)^{1-      c (\gamma)  \cdot c (\delta)  \cdot c(\lfloor \frac{|I|-1}{2} \rfloor ) \cdot E (\lfloor \frac{|I|-1}{2} \rfloor + 4) \cdot \eps } \cdot (1+|q|)^{1+\gamma - 2\de }}    \; ,
      \eeaa
  and
            \beaa
 \notag
&& | \Lie_{Z^K} A (t,x)  |    \\
\notag
&\leq&   C(q_0)   \cdot  c (\delta) \cdot c (\gamma) \cdot C( \lfloor \frac{|I|-1}{2} \rfloor) \cdot E ( \lfloor \frac{|I|-1}{2} \rfloor + 4)  \\
\notag
&& \times \frac{\eps }{(1+t+|q|)^{1-      c (\gamma)  \cdot c (\delta)  \cdot c( \lfloor \frac{|I|-1}{2} \rfloor) \cdot E ( \lfloor \frac{|I|-1}{2} \rfloor+ 4) \cdot \eps } \cdot (1+|q|)^{\gamma - 2\de }}    \; .
      \eeaa
      
           Hence, we obtain
           
           \beaa
&& \sum_{|J| + |K| \leq |I| -1}    | \Lie_{Z^J}  A_L |  \cdot     |  \derm ( \Lie_{Z^K} A )  |    \\
&\leq&  \sum_{ |J| \leq   \lfloor \frac{|I|-1}{2} \rfloor   ,\; |K| \leq |I| -1}  \Big(   | \Lie_{Z^J}  A_L |  \cdot     |  \derm ( \Lie_{Z^K} A )  | +  | \Lie_{Z^K}  A_L |  \cdot     |  \derm ( \Lie_{Z^J} A )  |     \Big) \\
&\leq& C(q_0)   \cdot  c (\delta) \cdot c (\gamma) \cdot C( \lfloor \frac{|I|-1}{2} \rfloor) \cdot E ( \lfloor \frac{|I|-1}{2} \rfloor + 4) \\
&& \times \sum_{ |K| \leq |I| -1}    \Big(    \frac{\eps  \cdot     |  \derm ( \Lie_{Z^K} A )  | }{(1+t+|q|)^{1-      c (\gamma)  \cdot c (\delta)  \cdot c( \lfloor \frac{|I|-1}{2} \rfloor) \cdot E ( \lfloor \frac{|I|-1}{2} \rfloor+ 4) \cdot \eps } \cdot (1+|q|)^{\gamma - 2\de }}    \\
&& + \frac{\eps \cdot  | \Lie_{Z^K}  A_L |  }{(1+t+|q|)^{1-      c (\gamma)  \cdot c (\delta)  \cdot c(\lfloor \frac{|I|-1}{2} \rfloor ) \cdot E (\lfloor \frac{|I|-1}{2} \rfloor + 4) \cdot \eps } \cdot (1+|q|)^{1+\gamma - 2\de }}  \Big) \;. 
\eeaa
           
           Thus, we obtain the stated result.
     \end{proof}

               \begin{corollary}\label{thestructureofthebadtermALdermAusingboostrapassumptionanddecompistionofthesumandlowerordertermsexhibitedtodealwithAL}
               We have for $\ga \geq 3 \de $\,, and $0 < \de \leq \frac{1}{4}$\,,
               
                                         \beaa
\notag
&& | \Lie_{Z^I}    \big(  A_L   \cdot     \derm A    \big)    | \\
&\les&  C(q_0)   \cdot  c (\delta) \cdot c (\gamma)  \cdot E ( 4)  \\
&& \times \sum_{|K| = |I| }  \Big(  \frac{ \eps  \cdot     |  \derm ( \Lie_{Z^K} A )  |   }{ (1+t+|q|)  \cdot  (1+|q|)^{\gamma - 2\de}  } \\
&&  +    \frac{\eps \cdot   |  \Lie_{Z^K} A_L |   }{(1+t+|q|)^{1-      c (\gamma)  \cdot c (\delta)   \cdot E ( 4) \cdot \eps } \cdot (1+|q|)^{1+\gamma - 2\de }}    \Big) \\
&& +  C(q_0)   \cdot  c (\delta) \cdot c (\gamma) \cdot C( \lfloor \frac{|I|-1}{2} \rfloor) \cdot E ( \lfloor \frac{|I|-1}{2} \rfloor + 4) \\
&& \times \sum_{ |K| \leq |I| -1}    \Big(    \frac{\eps  \cdot     |  \derm ( \Lie_{Z^K} A )  | }{(1+t+|q|)^{1-      c (\gamma)  \cdot c (\delta)  \cdot c( \lfloor \frac{|I|-1}{2} \rfloor) \cdot E ( \lfloor \frac{|I|-1}{2} \rfloor+ 4) \cdot \eps } \cdot (1+|q|)^{\gamma - 2\de }}    \\
&& + \frac{\eps \cdot  | \Lie_{Z^K}  A |  }{(1+t+|q|)^{1-      c (\gamma)  \cdot c (\delta)  \cdot c(\lfloor \frac{|I|-1}{2} \rfloor ) \cdot E (\lfloor \frac{|I|-1}{2} \rfloor + 4) \cdot \eps } \cdot (1+|q|)^{1+\gamma - 2\de }}   \Big) \;. 
\eeaa

               \end{corollary}

               \begin{proof}
               Putting together Corollary \ref{TheleeadingtermfortroublesomecomponentsintehsourcesforA} and Lemma \ref{estimateonlowerordertermsforthepotentialAinthesourcetermsforenergy}, and injecting in \eqref{decompositionofsumfortroublesometermsforsourcesonA}, we obtain
               
                              \beaa
\notag
&& | \Lie_{Z^I}    \big(  A_L   \cdot     \derm A    \big)    | \\
&\les& \sum_{|K| = |I| }  \Big( C(q_0)   \cdot  c (\delta) \cdot c (\gamma)  \cdot E ( 3) \cdot \frac{ \eps  \cdot     |  \derm ( \Lie_{Z^K} A )  |   }{ (1+t+|q|)^{2-2\delta}  \cdot  (1+|q|)^{\ga - 1}  } \\
&&  +    C(q_0)   \cdot  c (\delta) \cdot c (\gamma)  \cdot E ( 4)  \cdot \frac{\eps \cdot   |  \Lie_{Z^K} A_L |   }{(1+t+|q|)^{1-      c (\gamma)  \cdot c (\delta)   \cdot E ( 4) \cdot \eps } \cdot (1+|q|)^{1+\gamma - 2\de }}      \Big) \\
&& +  C(q_0)   \cdot  c (\delta) \cdot c (\gamma) \cdot C( \lfloor \frac{|I|-1}{2} \rfloor) \cdot E ( \lfloor \frac{|I|-1}{2} \rfloor + 4) \\
&& \times \sum_{ |K| \leq |I| -1}    \Big(    \frac{\eps  \cdot     |  \derm ( \Lie_{Z^K} A )  | }{(1+t+|q|)^{1-      c (\gamma)  \cdot c (\delta)  \cdot c( \lfloor \frac{|I|-1}{2} \rfloor) \cdot E ( \lfloor \frac{|I|-1}{2} \rfloor+ 4) \cdot \eps } \cdot (1+|q|)^{\gamma - 2\de }}    \\
&& + \frac{\eps \cdot  | \Lie_{Z^K}  A_L |  }{(1+t+|q|)^{1-      c (\gamma)  \cdot c (\delta)  \cdot c(\lfloor \frac{|I|-1}{2} \rfloor ) \cdot E (\lfloor \frac{|I|-1}{2} \rfloor + 4) \cdot \eps } \cdot (1+|q|)^{1+\gamma - 2\de }}   \Big) \;. 
\eeaa     

    Using the fact that
\beaa
 && C(q_0)   \cdot c (\delta) \cdot c (\gamma)  \cdot E ( 3) \cdot \frac{ \eps  \cdot     |  \derm ( \Lie_{Z^K} A )  |   }{ (1+t+|q|)^{2-2\delta}  \cdot  (1+|q|)^{\ga - 1}  }  \\
 &\les & C(q_0)   \cdot  c (\delta) \cdot c (\gamma)  \cdot E ( 3) \cdot \frac{ \eps  \cdot     |  \derm ( \Lie_{Z^K} A )  |   }{ (1+t+|q|)  \cdot  (1+|q|)^{1-2\delta + \ga - 1}  }  \\
  &\les &   C(q_0)   \cdot  c (\delta) \cdot c (\gamma)  \cdot E ( 4) \cdot \frac{ \eps  \cdot     |  \derm ( \Lie_{Z^K} A )  |   }{ (1+t+|q|)  \cdot  (1+|q|)^{\ga-2\delta }  }  \; ,
 \eeaa
we get the desired result.

               \end{proof}

\section{Structure of the source terms for the wave operator on the Einstein-Yang-Mills metric}

\subsection{Structure of the “good” source terms for the Einstein-Yang-Mills metric}\

\begin{lemma}\label{actualusefulstrzuctureofthesourcetermsforthewaveequationontheMetricsmallhoneusingbootstrap}
For $M \leq \eps$\;, for $\ga \geq 3 \de $\,, and $0 < \de \leq \frac{1}{4}$\,, we have in the exterior region, for all $U\,,V \in \cal{U}$\,,
\beaa
   \notag
&&  |  \Lie_{Z^I}  g^{\la\mu} \derm_{\la}   \derm_{\mu}  h^1_{UV}  |  \\
  \notag
  &\les& \sum_{|K| \leq |I |}  \Big[ C(q_0)   \cdot  c (\delta) \cdot c (\gamma) \cdot C(|I|) \cdot E (\lfloor \frac{|I|}{2} \rfloor + 5)  \cdot \frac{\eps \cdot  | \derm ( \Lie_{Z^K} h^1 ) | }{(1+t+|q|)^{2-      c (\gamma)  \cdot c (\delta)  \cdot c(|I|) \cdot E ( \lfloor \frac{|I|}{2} \rfloor+ 5) \cdot \eps }} \\
  && + C(q_0)   \cdot  c (\delta) \cdot c (\gamma) \cdot C(|I|) \cdot E (\lfloor \frac{|I|}{2} \rfloor+ 5)  \cdot \frac{\eps \cdot | \derm ( \Lie_{Z^K} A ) | }{(1+t+|q|)^{2-      c (\gamma)  \cdot c (\delta)  \cdot c(|I|) \cdot E ( \lfloor \frac{|I|}{2} \rfloor+ 5) \cdot \eps } \cdot (1+|q|)^{\gamma - 2\de }}  \\  
  && + C(q_0)   \cdot  c (\delta) \cdot c (\gamma) \cdot C(|I|) \cdot E ( \lfloor \frac{|I|}{2} \rfloor+ 4)  \cdot \frac{\eps \cdot | \rderm ( \Lie_{Z^K} A ) | }{(1+t+|q|)^{1-      c (\gamma)  \cdot c (\delta)  \cdot c(|I|) \cdot E ( \lfloor \frac{|I|}{2} \rfloor+ 4) \cdot \eps } \cdot (1+|q|)^{1+\gamma - 2\de }}  \\
 && +  C(q_0)   \cdot  c (\delta) \cdot c (\gamma) \cdot C(|I|) \cdot E ( \lfloor \frac{|I|}{2} \rfloor  +4)  \cdot \frac{\eps\cdot  |\rderm ( \Lie_{Z^K} h^1 ) |}{(1+t+|q|)^{1-   c (\gamma)  \cdot c (\delta)  \cdot c(|I|) \cdot E (\lfloor \frac{|I|}{2} \rfloor+ 4)\cdot  \eps } \cdot (1+|q|)}       \\
 && +  C(q_0)   \cdot  c (\delta) \cdot c (\gamma) \cdot C(|I|) \cdot E ( \lfloor \frac{|I|}{2} \rfloor+ 4)  \cdot \frac{\eps^2  \cdot  |\Lie_{Z^K}  A | }{(1+t+|q|)^{2-      c (\gamma)  \cdot c (\delta)  \cdot c(|I|) \cdot E ( \lfloor \frac{|I|}{2} \rfloor+ 4) \cdot \eps } \cdot (1+|q|)^{1+2\gamma - 4\de }} \\
    && +   C(q_0)   \cdot  c (\delta) \cdot c (\gamma) \cdot C(|I|) \cdot E ( \lfloor \frac{|I|}{2} \rfloor  +4)  \cdot \frac{\eps^2 \cdot  |\Lie_{Z^K}  h^1   |}{(1+t+|q|)^{2-   c (\gamma)  \cdot c (\delta)  \cdot c(|I|) \cdot E ( \lfloor \frac{|I|}{2} \rfloor+ 4)\cdot  \eps } \cdot (1+|q|)^2}    \Big]        \\
         &&+  | \Lie_{Z^I}  ( \derm h_{\cal TU} ) ^2   +   \Lie_{Z^I}   ( \derm A_{e_{a}}  )^2  | \\
 \notag
 && +   | g^{\alpha\beta} \derm_\alpha \derm_\beta h^0 |  \; \\
    && +  C(q_0)   \cdot  c (\delta) \cdot c (\gamma) \cdot C(|I|) \cdot E ( \lfloor \frac{|I|}{2} \rfloor  +5)  \cdot \frac{\eps^2 }{(1+t+|q|)^{3-   c (\gamma)  \cdot c (\delta)  \cdot c(|I|) \cdot E (\lfloor \frac{|I|}{2} \rfloor+ 5)\cdot  \eps } \cdot (1+|q|)}       \; .
  \eeaa

  For all $V \in \cal{T}$\,, we have the same estimate without the last term $  \Lie_{Z^I}  ( \derm h_{\cal TU} ) ^2  + \Lie_{Z^I}   ( \derm A_{e_{a}}  )^2  $\,.

  \end{lemma}
  
  \begin{proof}
  
Examining the structure of the source terms for the wave operator applied to $h^1$, exhibited in \eqref{StructureoftheLiederivativesofthesourcetermsofthewaveoperatorforAandh} -- that preserves for the Lie derivatives of the source terms, the structure exhibited in \eqref{structureofthesourcetermsofthewaveoperatoronAandh} --, we get using Lemmas \ref{upgradedestimateonLiederivativesoffields}, \ref{upgradedestimateonLiederivativesoffieldswithoutgradiant}, \ref{improvedtangentialderivatives}, \ref{upgradedestimatesonh}, \ref{upgradedestimategoodcomponentspotentialandmetric}, and \ref{upgradedestimatesonproducts}, and by using \eqref{generalinequalityaboutsplittingasumthatisconstrainedtoadduptolenth} to split the sum, that

\beaa
   \notag
&&  |  \Lie_{Z^I}  g^{\la\mu} \derm_{\la}   \derm_{\mu}  h^1  |  \\
  \notag
  &\les&   \sum_{|K| \leq |I |}  \Big[   \; C(q_0)   \cdot  c (\delta) \cdot c (\gamma) \cdot C(|I|) \cdot E (  \lfloor \frac{|I|}{2} \rfloor + 5)  \cdot \frac{\eps \cdot  | \derm ( \Lie_{Z^K} h ) | }{(1+t+|q|)^{2-      c (\gamma)  \cdot c (\delta)  \cdot c(|I|) \cdot E (  \lfloor \frac{|I|}{2} \rfloor+ 5) \cdot \eps }}   \\
  && +  C(q_0)   \cdot  c (\delta) \cdot c (\gamma) \cdot C(|I|) \cdot E (  \lfloor \frac{|I|}{2} \rfloor  +4)  \cdot \frac{\eps\cdot  |\rderm ( \Lie_{Z^K} h ) |}{(1+t+|q|)^{1-   c (\gamma)  \cdot c (\delta)  \cdot c(|I|) \cdot E ( \lfloor \frac{|I|}{2} \rfloor+ 4)\cdot  \eps } \cdot (1+|q|)}       \\
&& +     C(q_0)   \cdot  c (\delta) \cdot c (\gamma) \cdot C(|I|) \cdot E (   \lfloor \frac{|I|}{2} \rfloor  +4)  \cdot \frac{\eps^2 \cdot  |\Lie_{Z^K}  h   |}{(1+t+|q|)^{2-   c (\gamma)  \cdot c (\delta)  \cdot c(|I|) \cdot E ( \lfloor \frac{|I|}{2} \rfloor+ 4)\cdot  \eps } \cdot (1+|q|)^2} \\
&& +  C(q_0)   \cdot  c (\delta) \cdot c (\gamma) \cdot C(|I|) \cdot E (  \lfloor \frac{|I|}{2} \rfloor  +4)  \cdot \frac{\eps^2 \cdot | \derm ( \Lie_{Z^K} h )    | }{(1+t+|q|)^{2-   c (\gamma)  \cdot c (\delta)  \cdot c(|I|) \cdot E (  \lfloor \frac{|I|}{2} \rfloor+ 4)\cdot  \eps } \cdot (1+|q|)} \\
 && +  C(q_0)   \cdot  c (\delta) \cdot c (\gamma) \cdot C(|I|) \cdot E (  \lfloor \frac{|I|}{2} \rfloor+ 5)  \cdot \frac{\eps \cdot | \derm ( \Lie_{Z^K} A ) | }{(1+t+|q|)^{2-      c (\gamma)  \cdot c (\delta)  \cdot c(|I|) \cdot E (  \lfloor \frac{|I|}{2} \rfloor+ 5) \cdot \eps } \cdot (1+|q|)^{\gamma - 2\de }}\\
 && + C(q_0)   \cdot  c (\delta) \cdot c (\gamma) \cdot C(|I|) \cdot E (  \lfloor \frac{|I|}{2} \rfloor + 4)  \cdot \frac{\eps \cdot | \rderm ( \Lie_{Z^K} A ) | }{(1+t+|q|)^{1-      c (\gamma)  \cdot c (\delta)  \cdot c(|I|) \cdot E (  \lfloor \frac{|I|}{2} \rfloor+ 4) \cdot \eps } \cdot (1+|q|)^{1+\gamma - 2\de }}  \\
  && + C(q_0)   \cdot  c (\delta) \cdot c (\gamma) \cdot C(|I|) \cdot E (  \lfloor \frac{|I|}{2} \rfloor + 4)  \cdot \frac{\eps^2  \cdot | \derm ( \Lie_{Z^K} A ) | }{(1+t+|q|)^{2-      c (\gamma)  \cdot c (\delta)  \cdot c(|I|) \cdot E (  \lfloor \frac{|I|}{2} \rfloor+ 4) \cdot \eps } \cdot (1+|q|)^{2\gamma - 4\de }} \\
   && + C(q_0)   \cdot  c (\delta) \cdot c (\gamma) \cdot C(|I|) \cdot E (  \lfloor \frac{|I|}{2} \rfloor + 4)  \cdot \frac{\eps^2  \cdot  |\Lie_{Z^K}  A | }{(1+t+|q|)^{2-      c (\gamma)  \cdot c (\delta)  \cdot c(|I|) \cdot E (  \lfloor \frac{|I|}{2} \rfloor+ 4) \cdot \eps } \cdot (1+|q|)^{1+2\gamma - 4\de }}  \\
        && +  C(q_0)   \cdot  c (\delta) \cdot c (\gamma) \cdot C(|I|) \cdot E (  \lfloor \frac{|I|}{2} \rfloor+ 4)  \cdot \frac{\eps^3 \cdot  |\Lie_{Z^K}  A | }{(1+t+|q|)^{3-      c (\gamma)  \cdot c (\delta)  \cdot c(|I|) \cdot E ( | \lfloor \frac{|I|}{2} \rfloor+ 4) \cdot \eps } \cdot (1+|q|)^{3\gamma - 6\de }}   \\
  && +  O \big(      C(q_0)   \cdot  c (\delta) \cdot c (\gamma) \cdot C(|I|) \cdot E (  \lfloor \frac{|I|}{2} \rfloor + 4)  \cdot \frac{\eps \cdot  |\Lie_{Z^K}  h  | }{(1+t+|q|)^{2-      c (\gamma)  \cdot c (\delta)  \cdot c(|I|) \cdot E (  \lfloor \frac{|I|}{2} \rfloor+ 4) \cdot \eps } \cdot (1+|q|)^{2+2\gamma - 4\de }} \big)  \\
  && +   O \big(  C(q_0)   \cdot  c (\delta) \cdot c (\gamma) \cdot C(|I|) \cdot E ( \lfloor \frac{|I|}{2} \rfloor + 4)  \cdot \frac{\eps  \cdot | \derm ( \Lie_{Z^K} A ) | }{(1+t+|q|)^{2-      c (\gamma)  \cdot c (\delta)  \cdot c(|I|) \cdot E (  \lfloor \frac{|I|}{2} \rfloor+ 4) \cdot \eps } \cdot (1+|q|)^{1+\gamma - 2\de }} \big) \; \Big]  \\
   &&+ |\derm h_{\cal TU} |^2  +   |\derm A_{e_{a}}  |^2  \\
 \notag
 && +   | g^{\alpha\beta} \derm_\alpha \derm_\beta h^0 |  
  \eeaa

Putting the terms together and using the fact that $\ga \geq 3 \de  \geq 2\de  $\,, we get the following factors for each term,

    \textbf{The factor for $  | \derm ( \Lie_{Z^I} h ) | $:}\\
 We have
 \beaa   
&&C(q_0)   \cdot  c (\delta) \cdot c (\gamma) \cdot C(|I|) \cdot E (\lfloor \frac{|I|}{2} \rfloor+ 5)  \cdot \frac{\eps \cdot  | \derm ( \Lie_{Z^I} h ) | }{(1+t+|q|)^{2-      c (\gamma)  \cdot c (\delta)  \cdot c(|I|) \cdot E ( \lfloor \frac{|I|}{2} \rfloor+ 5) \cdot \eps }}   \\
&& +  C(q_0)   \cdot  c (\delta) \cdot c (\gamma) \cdot C(|I|) \cdot E ( \lfloor \frac{|I|}{2} \rfloor  +4)  \cdot \frac{\eps^2 \cdot | \derm ( \Lie_{Z^I} h )    | }{(1+t+|q|)^{2-   c (\gamma)  \cdot c (\delta)  \cdot c(|I|) \cdot E ( \lfloor \frac{|I|}{2} \rfloor+ 4)\cdot  \eps } \cdot (1+|q|)} \\
&\les&C(q_0)   \cdot  c (\delta) \cdot c (\gamma) \cdot C(|I|) \cdot E ( \lfloor \frac{|I|}{2} \rfloor+ 5)  \cdot \frac{\eps \cdot  | \derm ( \Lie_{Z^I} h ) | }{(1+t+|q|)^{2-      c (\gamma)  \cdot c (\delta)  \cdot c(|I|) \cdot E ( \lfloor \frac{|I|}{2} \rfloor+ 5) \cdot \eps }}  \; .
\eeaa

      \textbf{The factor for $ |\derm ( \Lie_{Z^I} A ) | $:}\\

We get
\beaa
 &&   C(q_0)   \cdot  c (\delta) \cdot c (\gamma) \cdot C(|I|) \cdot E ( \lfloor \frac{|I|}{2} \rfloor + 5)  \cdot \frac{\eps \cdot | \derm ( \Lie_{Z^I} A ) | }{(1+t+|q|)^{2-      c (\gamma)  \cdot c (\delta)  \cdot c(|I|) \cdot E ( \lfloor \frac{|I|}{2} \rfloor+ 5) \cdot \eps } \cdot (1+|q|)^{\gamma - 2\de }}\\
    && + C(q_0)   \cdot  c (\delta) \cdot c (\gamma) \cdot C(|I|) \cdot E (\lfloor \frac{|I|}{2} \rfloor + 4)  \cdot \frac{\eps^2  \cdot | \derm ( \Lie_{Z^I} A ) | }{(1+t+|q|)^{2-      c (\gamma)  \cdot c (\delta)  \cdot c(|I|) \cdot E ( \lfloor \frac{|I|}{2} \rfloor+ 4) \cdot \eps } \cdot (1+|q|)^{2\gamma - 4\de }} \\
&& +   O \big(  C(q_0)   \cdot  c (\delta) \cdot c (\gamma) \cdot C(|I|) \cdot E ( \lfloor \frac{|I|}{2} \rfloor+ 4)  \cdot \frac{\eps  \cdot | \derm ( \Lie_{Z^I} A ) | }{(1+t+|q|)^{2-      c (\gamma)  \cdot c (\delta)  \cdot c(|I|) \cdot E (\lfloor \frac{|I|}{2} \rfloor+ 4) \cdot \eps } \cdot (1+|q|)^{1+\gamma - 2\de }} \big) \\
&\les &  C(q_0)   \cdot  c (\delta) \cdot c (\gamma) \cdot C(|I|) \cdot E ( \lfloor \frac{|I|}{2} \rfloor+ 5)  \cdot \frac{\eps \cdot | \derm ( \Lie_{Z^I} A ) | }{(1+t+|q|)^{2-      c (\gamma)  \cdot c (\delta)  \cdot c(|I|) \cdot E ( \lfloor \frac{|I|}{2} \rfloor+ 5) \cdot \eps } \cdot (1+|q|)^{\gamma - 2\de }}  \; .
\eeaa

           \textbf{The factor for $ |\rderm ( \Lie_{Z^I} A ) |$:}\\
           
It is 
\beaa
 && C(q_0)   \cdot  c (\delta) \cdot c (\gamma) \cdot C(|K|) \cdot E ( \lfloor \frac{|I|}{2} \rfloor + 4)  \cdot \frac{\eps \cdot | \rderm ( \Lie_{Z^I} A ) | }{(1+t+|q|)^{1-      c (\gamma)  \cdot c (\delta)  \cdot c(|K|) \cdot E ( \lfloor \frac{|I|}{2} \rfloor|+ 4) \cdot \eps } \cdot (1+|q|)^{1+\gamma - 2\de }}  \\
   \eeaa
      
         \textbf{The factor for $ | \rderm  ( \Lie_{Z^I} h ) |$:}\\
         It is 
         \beaa
           &&   C(q_0)   \cdot  c (\delta) \cdot c (\gamma) \cdot C(|I|) \cdot E ( \lfloor \frac{|I|}{2} \rfloor +4)  \cdot \frac{\eps\cdot  |\rderm ( \Lie_{Z^I} h ) |}{(1+t+|q|)^{1-   c (\gamma)  \cdot c (\delta)  \cdot c(|I|) \cdot E (\lfloor \frac{|I|}{2} \rfloor+ 4)\cdot  \eps } \cdot (1+|q|)}       \\
     \eeaa
     
     \textbf{The factors for $|\Lie_{Z^I} A|$:}\\
     
     We obtain
   \beaa
      &&  C(q_0)   \cdot  c (\delta) \cdot c (\gamma) \cdot C(|I|) \cdot E (\lfloor \frac{|I|}{2} \rfloor+ 4)  \cdot \frac{\eps^2  \cdot  |\Lie_{Z^I}  A | }{(1+t+|q|)^{2-      c (\gamma)  \cdot c (\delta)  \cdot c(|I|) \cdot E (\lfloor \frac{|I|}{2} \rfloor+ 4) \cdot \eps } \cdot (1+|q|)^{1+2\gamma - 4\de }}  \\
        && +  C(q_0)   \cdot  c (\delta) \cdot c (\gamma) \cdot C(|K|) \cdot E (\lfloor \frac{|I|}{2} \rfloor+ 4)  \cdot \frac{\eps^3 \cdot  |\Lie_{Z^I}  A | }{(1+t+|q|)^{3-      c (\gamma)  \cdot c (\delta)  \cdot c(|I|) \cdot E ( |\lfloor \frac{|I|}{2} \rfloor+ 4) \cdot \eps } \cdot (1+|q|)^{3\gamma - 6\de }}   \\
        &\les& C(q_0)   \cdot  c (\delta) \cdot c (\gamma) \cdot C(|I|) \cdot E (\lfloor \frac{|I|}{2} \rfloor+ 4)  \cdot \frac{\eps^2  \cdot  |\Lie_{Z^I}  A | }{(1+t+|q|)^{2-      c (\gamma)  \cdot c (\delta)  \cdot c(|I|) \cdot E ( \lfloor \frac{|I|}{2} \rfloor+ 4) \cdot \eps } \cdot (1+|q|)^{1+2\gamma - 4\de }}  \;.
   \eeaa
                    
        \textbf{The factors for $|\Lie_{Z^I} h|$:}\\
 We estimate
              \beaa
                &&     C(q_0)   \cdot  c (\delta) \cdot c (\gamma) \cdot C(|I|) \cdot E (  \lfloor \frac{|I|}{2} \rfloor +4)  \cdot \frac{\eps^2 \cdot  |\Lie_{Z^I}  h   |}{(1+t+|q|)^{2-   c (\gamma)  \cdot c (\delta)  \cdot c(|I|) \cdot E ( \lfloor \frac{|I|}{2} \rfloor+ 4)\cdot  \eps } \cdot (1+|q|)^2} \\
                  && +  O \big(      C(q_0)   \cdot  c (\delta) \cdot c (\gamma) \cdot C(|I|) \cdot E (\lfloor \frac{|I|}{2} \rfloor+ 4)  \cdot \frac{\eps \cdot  |\Lie_{Z^I}  h  | }{(1+t+|q|)^{2-      c (\gamma)  \cdot c (\delta)  \cdot c(|I|) \cdot E ( \lfloor \frac{|I|}{2} \rfloor+ 4) \cdot \eps } \cdot (1+|q|)^{2+2\gamma - 4\de }} \big)  \\
                  &\les&  C(q_0)   \cdot  c (\delta) \cdot c (\gamma) \cdot C(|I|) \cdot E (  \lfloor \frac{|I|}{2} \rfloor +4)  \cdot \frac{\eps^2 \cdot  |\Lie_{Z^I}  h   |}{(1+t+|q|)^{2-   c (\gamma)  \cdot c (\delta)  \cdot c(|I|) \cdot E (\lfloor \frac{|I|}{2} \rfloor+ 4)\cdot  \eps } \cdot (1+|q|)^2} \; .
        \eeaa
                                                     
Putting all together, we obtain
\beaa
   \notag
&&  |  \Lie_{Z^I}  g^{\la\mu} \derm_{\la}   \derm_{\mu}  h^1_{UV}  |  \\
  \notag
  &\les& \sum_{|K| \leq |I |}  \Big[ C(q_0)   \cdot  c (\delta) \cdot c (\gamma) \cdot C(|I|) \cdot E (\lfloor \frac{|I|}{2} \rfloor + 5)  \cdot \frac{\eps \cdot  | \derm ( \Lie_{Z^K} h ) | }{(1+t+|q|)^{2-      c (\gamma)  \cdot c (\delta)  \cdot c(|I|) \cdot E ( \lfloor \frac{|I|}{2} \rfloor+ 5) \cdot \eps }} \\
  && + C(q_0)   \cdot  c (\delta) \cdot c (\gamma) \cdot C(|I|) \cdot E (\lfloor \frac{|I|}{2} \rfloor+ 5)  \cdot \frac{\eps \cdot | \derm ( \Lie_{Z^K} A ) | }{(1+t+|q|)^{2-      c (\gamma)  \cdot c (\delta)  \cdot c(|I|) \cdot E ( \lfloor \frac{|I|}{2} \rfloor+ 5) \cdot \eps } \cdot (1+|q|)^{\gamma - 2\de }}  \\  
  && + C(q_0)   \cdot  c (\delta) \cdot c (\gamma) \cdot C(|I|) \cdot E ( \lfloor \frac{|I|}{2} \rfloor+ 4)  \cdot \frac{\eps \cdot | \rderm ( \Lie_{Z^K} A ) | }{(1+t+|q|)^{1-      c (\gamma)  \cdot c (\delta)  \cdot c(|I|) \cdot E ( \lfloor \frac{|I|}{2} \rfloor+ 4) \cdot \eps } \cdot (1+|q|)^{1+\gamma - 2\de }}  \\
 && +  C(q_0)   \cdot  c (\delta) \cdot c (\gamma) \cdot C(|I|) \cdot E ( \lfloor \frac{|I|}{2} \rfloor  +4)  \cdot \frac{\eps\cdot  |\rderm ( \Lie_{Z^K} h ) |}{(1+t+|q|)^{1-   c (\gamma)  \cdot c (\delta)  \cdot c(|I|) \cdot E (\lfloor \frac{|I|}{2} \rfloor+ 4)\cdot  \eps } \cdot (1+|q|)}       \\
 && +  C(q_0)   \cdot  c (\delta) \cdot c (\gamma) \cdot C(|I|) \cdot E ( \lfloor \frac{|I|}{2} \rfloor+ 4)  \cdot \frac{\eps^2  \cdot  |\Lie_{Z^K}  A | }{(1+t+|q|)^{2-      c (\gamma)  \cdot c (\delta)  \cdot c(|I|) \cdot E ( \lfloor \frac{|I|}{2} \rfloor+ 4) \cdot \eps } \cdot (1+|q|)^{1+2\gamma - 4\de }} \\
    && +   C(q_0)   \cdot  c (\delta) \cdot c (\gamma) \cdot C(|I|) \cdot E ( \lfloor \frac{|I|}{2} \rfloor  +4)  \cdot \frac{\eps^2 \cdot  |\Lie_{Z^K}  h   |}{(1+t+|q|)^{2-   c (\gamma)  \cdot c (\delta)  \cdot c(|I|) \cdot E ( \lfloor \frac{|I|}{2} \rfloor+ 4)\cdot  \eps } \cdot (1+|q|)^2}    \Big]        \\
         &&+  | \Lie_{Z^I}  ( \derm h_{\cal TU} ) ^2   +   \Lie_{Z^I}   ( \derm A_{e_{a}}  )^2  | \\
 \notag
 && +   | g^{\alpha\beta} \derm_\alpha \derm_\beta h^0 |  \; .
  \eeaa
  
  Now, again, since we want estimates involving $h^1$, for which the expression of the energy is defined, and not $h$, we use Lemma \ref{Liederivativesofsphericalsymmetricpart}, to estimate for $M \leq \eps $\,, 
    \beaa
 \notag
|   \Lie_{Z^I} h  |   &\leq& |   \Lie_{Z^I} h^1  | + c (\delta) \cdot C ( |I| )   \cdot \frac{\eps }{(1+t+|q|) }  \; , \\
 \notag
|\derm (  \Lie_{Z^I} h)  |   &\leq& |\derm (  \Lie_{Z^I} h^1)  |  +  C ( |I| )   \cdot \frac{\eps }{(1+t+|q|)^{2} } \; .
      \eeaa
Thus, injecting, we get
\beaa
   \notag
&&  |  \Lie_{Z^I}  g^{\la\mu} \derm_{\la}   \derm_{\mu}  h^1_{UV}  |  \\
  \notag
  &\les& \sum_{|K| \leq |I |}  \Big[ C(q_0)   \cdot  c (\delta) \cdot c (\gamma) \cdot C(|I|) \cdot E (\lfloor \frac{|I|}{2} \rfloor + 5)  \cdot \frac{\eps \cdot  | \derm ( \Lie_{Z^K} h^1 ) | }{(1+t+|q|)^{2-      c (\gamma)  \cdot c (\delta)  \cdot c(|I|) \cdot E ( \lfloor \frac{|I|}{2} \rfloor+ 5) \cdot \eps }} \\
  && + C(q_0)   \cdot  c (\delta) \cdot c (\gamma) \cdot C(|I|) \cdot E (\lfloor \frac{|I|}{2} \rfloor+ 5)  \cdot \frac{\eps \cdot | \derm ( \Lie_{Z^K} A ) | }{(1+t+|q|)^{2-      c (\gamma)  \cdot c (\delta)  \cdot c(|I|) \cdot E ( \lfloor \frac{|I|}{2} \rfloor+ 5) \cdot \eps } \cdot (1+|q|)^{\gamma - 2\de }}  \\  
  && + C(q_0)   \cdot  c (\delta) \cdot c (\gamma) \cdot C(|I|) \cdot E ( \lfloor \frac{|I|}{2} \rfloor+ 4)  \cdot \frac{\eps \cdot | \rderm ( \Lie_{Z^K} A ) | }{(1+t+|q|)^{1-      c (\gamma)  \cdot c (\delta)  \cdot c(|I|) \cdot E ( \lfloor \frac{|I|}{2} \rfloor+ 4) \cdot \eps } \cdot (1+|q|)^{1+\gamma - 2\de }}  \\
 && +  C(q_0)   \cdot  c (\delta) \cdot c (\gamma) \cdot C(|I|) \cdot E ( \lfloor \frac{|I|}{2} \rfloor  +4)  \cdot \frac{\eps\cdot  |\rderm ( \Lie_{Z^K} h^1 ) |}{(1+t+|q|)^{1-   c (\gamma)  \cdot c (\delta)  \cdot c(|I|) \cdot E (\lfloor \frac{|I|}{2} \rfloor+ 4)\cdot  \eps } \cdot (1+|q|)}       \\
 && +  C(q_0)   \cdot  c (\delta) \cdot c (\gamma) \cdot C(|I|) \cdot E ( \lfloor \frac{|I|}{2} \rfloor+ 4)  \cdot \frac{\eps^2  \cdot  |\Lie_{Z^K}  A | }{(1+t+|q|)^{2-      c (\gamma)  \cdot c (\delta)  \cdot c(|I|) \cdot E ( \lfloor \frac{|I|}{2} \rfloor+ 4) \cdot \eps } \cdot (1+|q|)^{1+2\gamma - 4\de }} \\
    && +   C(q_0)   \cdot  c (\delta) \cdot c (\gamma) \cdot C(|I|) \cdot E ( \lfloor \frac{|I|}{2} \rfloor  +4)  \cdot \frac{\eps^2 \cdot  |\Lie_{Z^K}  h^1   |}{(1+t+|q|)^{2-   c (\gamma)  \cdot c (\delta)  \cdot c(|I|) \cdot E ( \lfloor \frac{|I|}{2} \rfloor+ 4)\cdot  \eps } \cdot (1+|q|)^2}    \Big]        \\
         &&+  | \Lie_{Z^I}  ( \derm h_{\cal TU} ) ^2   +   \Lie_{Z^I}   ( \derm A_{e_{a}}  )^2  | \\
 \notag
 && +   | g^{\alpha\beta} \derm_\alpha \derm_\beta h^0 |    \; 
  \eeaa
  \beaa
  \notag
  && + C(q_0)   \cdot  c (\delta) \cdot c (\gamma) \cdot C(|I|) \cdot E (\lfloor \frac{|I|}{2} \rfloor + 5)  \cdot \frac{\eps^2  }{(1+t+|q|)^{4-      c (\gamma)  \cdot c (\delta)  \cdot c(|I|) \cdot E ( \lfloor \frac{|I|}{2} \rfloor+ 5) \cdot \eps }} \\
  \notag
   && +  C(q_0)   \cdot  c (\delta) \cdot c (\gamma) \cdot C(|I|) \cdot E ( \lfloor \frac{|I|}{2} \rfloor  +4)  \cdot \frac{\eps^2 }{(1+t+|q|)^{3-   c (\gamma)  \cdot c (\delta)  \cdot c(|I|) \cdot E (\lfloor \frac{|I|}{2} \rfloor+ 4)\cdot  \eps } \cdot (1+|q|)}       \\
   \notag
       && +   C(q_0)   \cdot  c (\delta) \cdot c (\gamma) \cdot C(|I|) \cdot E ( \lfloor \frac{|I|}{2} \rfloor  +4)  \cdot \frac{\eps^3 }{(1+t+|q|)^{3-   c (\gamma)  \cdot c (\delta)  \cdot c(|I|) \cdot E ( \lfloor \frac{|I|}{2} \rfloor+ 4)\cdot  \eps } \cdot (1+|q|)^2}       \\
  \eeaa
  
  \end{proof}

  \subsection{Structure of the “bad” source terms for the Einstein-Yang-Mills metric}\ \label{Structurebadsourcetermsmetricinssubsection}\

We examine now, the seemingly troublesome terms $\Lie_{Z^I}  ( \derm h_{\cal TU} ) ^2  + \Lie_{Z^I}   ( \derm A_{e_{a}}  )^2 $, which we had left over in our estimates. In fact, these terms, by a naive estimate, using Lemmas \ref{LiederivativeZofMinkwoskimetric} and \ref{LiederivativeZcommuteswithminkowksicovariantderivative}, are evaluated as 
\bea\label{straightforwardderivativeoftroubletermsforhinthestructure}
\notag
&& | \Lie_{Z^I}  ( \derm h_{\cal TU} ) ^2  + \Lie_{Z^I}   ( \derm A_{e_{a}}  )^2 \big)  | \\
\notag
&\les& \sum_{|J| + |K| \leq |I| }  | \derm  ( \Lie_{Z^J} h_{\cal TU} ) |  \cdot     |  \derm  ( \Lie_{Z^K} h_{\cal TU}  )  |  +  \sum_{|J| + |K| \leq |I| }  |\derm ( \Lie_{Z^J} A_{e_a} )  |  \cdot  | \derm ( \Lie_{Z^K} A_{e_a} ) | \, .\\
\eea
As we have done for the Einstein-Yang-Mills potential, we are going to decompose the sum in way where we extract a “good” term, and a “bad” term of lower order. 

\begin{lemma}\label{Thebadtermforthemetricisestimatehereforetheclosureofenergyestimate}
 We have $\ga \geq 3 \de $\,, and $0 < \de \leq \frac{1}{4}$\,,
\beaa
\notag
&& \sum_{|J| + |K| \leq |I| }  | \derm  ( \Lie_{Z^J} h_{\cal TU} ) |  \cdot     |  \derm  ( \Lie_{Z^K} h_{\cal TU}  )  |  +  \sum_{|J| + |K| \leq |I| }  |\derm ( \Lie_{Z^J} A_{e_a} )  |  \cdot  | \derm ( \Lie_{Z^K} A_{e_a} ) | \\
&\les& \sum_{|K| = |I| }  C(q_0) \cdot   c (\gamma)  \cdot c (\delta)   \cdot E (4) \cdot  \Big(  \frac{\eps \cdot     | \derm  ( \Lie_{Z^J} h_{\cal TU} ) | }{ (1+t+|q|)} +  \frac{\eps \cdot |\derm ( \Lie_{Z^J} A_{e_a} )  | }{ (1+t+|q|)} \Big) \\
&& + C(q_0)   \cdot  c (\delta) \cdot c (\gamma) \cdot C( \lfloor \frac{|I|-1}{2} \rfloor) \cdot E (  \lfloor \frac{|I|-1}{2} \rfloor+ 4)  \\
&& \times \sum_{ |K| \leq |I| -1}    \Big(  \frac{\eps   \cdot | \derm  ( \Lie_{Z^K} h_{\cal TU} ) | }{(1+t+|q|)^{1-   c (\gamma)  \cdot c (\delta)  \cdot c( \lfloor \frac{|I|-1}{2} \rfloor) \cdot E (  \lfloor \frac{|I|-1}{2} \rfloor+ 4)\cdot  \eps } \cdot (1+|q|)}  \\
&& +\frac{\eps  \cdot     |\derm ( \Lie_{Z^K} A_{e_a} )  | }{(1+t+|q|)^{1-      c (\gamma)  \cdot c (\delta)  \cdot c(|K|) \cdot E (  \lfloor \frac{|I|-1}{2} \rfloor+ 4) \cdot \eps } \cdot (1+|q|)^{1+\gamma - 2\de }}  \Big) \;. 
\eeaa

\end{lemma}
\begin{proof}

We split the sum in \eqref{straightforwardderivativeoftroubletermsforhinthestructure} as follows

\beaa
&& \sum_{|J| + |K| \leq |I| } \Big( | \derm  ( \Lie_{Z^J} h_{\cal TU} ) |  \cdot     |  \derm  ( \Lie_{Z^K} h_{\cal TU}  )  |  +  |\derm ( \Lie_{Z^J} A_{e_a} )  |  \cdot  | \derm ( \Lie_{Z^K} A_{e_a} ) |  \Big) \\
\notag
&\les& \sum_{|K| = |I| }  \Big(   | \derm   h_{\cal TU}  |  \cdot     | \derm  ( \Lie_{Z^J} h_{\cal TU} ) |  +  |\derm  A_{e_a}   | \cdot     |\derm ( \Lie_{Z^J} A_{e_a} )  |   \Big) \\
\notag
&& +  \sum_{|J| + |K| \leq |I| -1}   \Big( | \derm  ( \Lie_{Z^J} h_{\cal TU} ) | \cdot    | \derm  ( \Lie_{Z^K} h_{\cal TU} ) | +  |\derm ( \Lie_{Z^J} A_{e_a} )  | \cdot     |\derm ( \Lie_{Z^K} A_{e_a} )  |   \Big) \; . \\
\eeaa
We then see that the leading terms have a “good” structure, thanks to estimates in Lemma \ref{estimatefortangentialcomponentspotential} and in Lemma \ref{upgradedestimatesongoodcomponnentforh1andh0}, in particular \eqref{theupgradedestimateongoodcomponenth}, and we obtain that for \,$\ga > \de $\,, and $0 \leq \delta \leq \frac{1}{4}$\,,

\beaa
&& \sum_{|K| = |I| }  \Big(   | \derm   h_{\cal TU}  |  \cdot     | \derm  ( \Lie_{Z^J} h_{\cal TU} ) |  +  |\derm  A_{e_a}   | \cdot     |\derm ( \Lie_{Z^J} A_{e_a} )  |   \Big) \\
&\les&  \sum_{|K| = |I| }  \Big(  C(q_0) \cdot   c (\gamma)  \cdot c (\delta)   \cdot E (4) \cdot \frac{\eps \cdot     | \derm  ( \Lie_{Z^J} h_{\cal TU} ) | }{ (1+t+|q|)} \\
&& + C(q_0) \cdot  c (\gamma)  \cdot c (\delta)   \cdot E (4) \cdot \frac{\eps \cdot |\derm ( \Lie_{Z^J} A_{e_a} )  | }{ (1+t+|q|)} \Big) \; .
\eeaa

Whereas to the lower order terms, we write them as 
\beaa
&&  \sum_{|J| + |K| \leq |I| -1}   \Big( | \derm  ( \Lie_{Z^J} h_{\cal TU} ) | \cdot    | \derm  ( \Lie_{Z^K} h_{\cal TU} ) | +  |\derm ( \Lie_{Z^J} A_{e_a} )  | \cdot     |\derm ( \Lie_{Z^K} A_{e_a} )  |   \Big)  \\
&\leq & \sum_{ |J| \leq   \lfloor \frac{|I|-1}{2} \rfloor   ,\; |K| \leq |I| -1}    \Big( | \derm  ( \Lie_{Z^J} h_{\cal TU} ) | \cdot    | \derm  ( \Lie_{Z^K} h_{\cal TU} ) | +  |\derm ( \Lie_{Z^J} A_{e_a} )  | \cdot     |\derm ( \Lie_{Z^K} A_{e_a} )  |   \Big)
\eeaa
Using Lemma \ref{upgradedestimateonLiederivativesoffields},  we have $\ga \geq 3 \de $\,, and $0 < \de \leq \frac{1}{4}$\,,
\beaa
&& \sum_{ |J| \leq   \lfloor \frac{|I|-1}{2} \rfloor   ,\; |K| \leq |I| -1}    \Big( | \derm  ( \Lie_{Z^J} h_{\cal TU} ) | \cdot    | \derm  ( \Lie_{Z^K} h_{\cal TU} ) | +  |\derm ( \Lie_{Z^J} A_{e_a} )  | \cdot     |\derm ( \Lie_{Z^K} A_{e_a} )  | \\
&\les& C(q_0)   \cdot  c (\delta) \cdot c (\gamma) \cdot C( \lfloor \frac{|I|-1}{2} \rfloor) \cdot E (  \lfloor \frac{|I|-1}{2} \rfloor+ 4)  \\
&& \times \sum_{ |K| \leq |I| -1}    \Big(  \frac{\eps   \cdot | \derm  ( \Lie_{Z^K} h_{\cal TU} ) | }{(1+t+|q|)^{1-   c (\gamma)  \cdot c (\delta)  \cdot c( \lfloor \frac{|I|-1}{2} \rfloor) \cdot E (  \lfloor \frac{|I|-1}{2} \rfloor+ 4)\cdot  \eps } \cdot (1+|q|)}  \\
&& +\frac{\eps  \cdot     |\derm ( \Lie_{Z^K} A_{e_a} )  | }{(1+t+|q|)^{1-      c (\gamma)  \cdot c (\delta)  \cdot c(|K|) \cdot E (  \lfloor \frac{|I|-1}{2} \rfloor+ 4) \cdot \eps } \cdot (1+|q|)^{1+\gamma - 2\de }}  \Big) \;. 
\eeaa

Putting the estimate on the leading term and the one on the lower order term together, we get the desired result.
\end{proof}

 \section{The energy estimate}
 
 We recapitulate the following definitions from \cite{G5}.

 \begin{definition}\label{defofthestreessenergymomentumtensorforwaveequationhere}
 Let either $\Phi$ be a tensor either valued in the Lie algebra $\cal G$ or a scalar. In particular, we will have in this paper either $\Phi_{V} = A_{V}$ or $\Phi_{UV} = h_{UV}$\;.
 We consider the following non-symmetric tensor for wave equations:
 when $\Phi =A$\;, we define
\bea\label{definitionofthewavestreessenergymomentumtensorfortensorphiVwhereVisanyvector}
\notag
T^{(\bf{g}) \; \mu}_{\;\;\;\;\;\;\;\;\;\;\nu} (\Phi_{V} )    =  g^{\mu\a}< \derm_\a \Phi_{V} ,   \derm_\nu \Phi_{V} > - \frac{1}{2} m^{\mu}_{\;\;\;\nu}  \cdot g^{\a\b}  < \derm_\a \Phi_{V} ,   \derm_\b \Phi_{V} > \;, \\
 \eea
 and when $\Phi = h$\;, we define
 \bea\label{definitionofthewavestreessenergymomentumtensorfortensorphiVwhereUVareanyvector}
\notag
T^{(\bf{g}) \; \mu}_{\;\;\;\;\;\;\;\;\;\;\nu} (\Phi_{UV} )    =  g^{\mu\a}< \derm_\a \Phi_{UV}  ,   \derm_\nu \Phi_{UV}  > - \frac{1}{2} m^{\mu}_{\;\;\;\nu}  \cdot g^{\a\b}  < \derm_\a \Phi_{UV}  ,   \derm_\b \Phi_{UV}  > \;, \\
 \eea
where we raise index with respect to the Minkowski metric $m$\,, defined to be Minkowski in wave coordinates. We consider $\Phi$ to be a field decaying fast enough at spatial infinity, so that there is no contribution at null infinity.

\end{definition}

\begin{definition}\label{definitionofLtilde}
We define $ \widetilde{L}^{\nu} $ as a vector at a point in the space-time, such that  $T_{\widetilde{L} t}^{(\bf{g})} \geq 0$\;. We note that such a definition does not give a unique vector, however, we will use this to construct $N_{t_1}^{t_2} (q_0) $ in Definition \ref{definitionofNandofNtruncatedbtweentwots} and to define $ \hat{L}^{\nu} $ in Definition \ref{definitionofwidetildeLfordivergencetheoremuse}.
 \end{definition}

\begin{definition}\label{definitionofNandofNtruncatedbtweentwots}
We define $N_{t_1}^{t_2} (q_0) $ as a boundary made by the following:

For a $t_N \geq 0 $\;, we take a curve in the plane $(t, r)$\;, at an $\Om \in \SSS^{n-1}$ fixed, that starts at $(t=t_N\,,\, r=0) $ and extends to the future of $\Sigma_{t_0}$\;, and where $ \widetilde{L}^{\nu} $ (defined in Definition \ref{definitionofLtilde}) is orthogonal (with respect to the Minkowski metric $m$) to that curve at each point. For a $t_N \geq 0 $ large enough, depending on $q_0$\;, we have under the bootstrap assumption, and therefore under the a priori estimates in Lemmas \ref{apriordecayestimatesfrombootstrapassumption} and \ref{apriorestimateontheEinsteinYangMillsfieldswithoutgradiant}, that there exits such a curve that is contained in the region $\{(t, x) \;|\; q:= r-t \leq q_0 \}$\;.

We define $N (q_0)$ as the product of $\SSS^{n-1}$ and of such future inextendable curve, and we define $N_{t_1}^{t_2} (q_0) $ as $N (q_0) $ truncated between $t_1$ and $t_2$\;. Thus,  $N (q_0)$ depends on the choice of the starting point $(t_N, 0)$ that is the tip of $N (q_0)$\;, however we write $N(q_0)$ to refer to a choice of such $N(q_0)$\;. To lighten the notation, we write $N$\;, instead of $N(q_0)$\;, where a $q_0$ has been chosen fixed. 

\end{definition} 

\begin{definition}\label{definitionofoverlinebigC} Let $\overline{C}$ be the exterior region defined as the following:

The boundary $N$ separates space-time into two regions: one interior region that we shall call $C$ where space-like curves are contained in a compact region, and the other region is the complement of $C$\;, where space-like curves can go to spatial infinity. We define $\overline{C}$ as the complement of $C$\;. 

\end{definition} 

\begin{definition}\label{definitionofSigmaexteriorattequalconstantusingthestressenergymomentumtensorTgforwaveeq}
 We define for a given fixed $t$\;,
   \bea
\Sigma^{ext}_{t}  &:=&  \Sigma_t  \cap \overline{C} \, ,
\eea
where $C$ is defined in Definition \ref{definitionofoverlinebigC}.
\end{definition} 

  \begin{remark}\label{RemarkaboutthedifferenceofdefinitionofSigmaexterior}
However, in Definition \ref{definitionofSigmatexterasintersectedwithcomplementofinteriorofcausaldomainofK}, we gave a different definition for $\Sigma^{ext}_{t}$\;. We note that under the bootstrap assumptions, and therefore under the a priori estimates in Lemmas \ref{apriordecayestimatesfrombootstrapassumption} and \ref{apriorestimateontheEinsteinYangMillsfieldswithoutgradiant}, that these two definitions are equivalent for our purposes. In other words, in this context of our bootstrap assumptions,  the estimates that we get on $\Sigma^{ext}_{t}$ from one definition (of the two Definitions \ref{definitionofSigmatexterasintersectedwithcomplementofinteriorofcausaldomainofK}
 and \ref{definitionofSigmaexteriorattequalconstantusingthestressenergymomentumtensorTgforwaveeq}) could  be translated into the same estimates with the other definition.
  \end{remark}

\begin{definition} 

We define $n^{(\bf{m}), \nu}_{\Sigma}$ as the unit orthogonal vector (for the metric $m$) to the hypersurface $\Sigma^{ext}_{t} $, and $dv^{(\bf{m})}_\Sigma$ as the induced volume form on $\Sigma_t$.

We denote by $n^{(\bf{m}), \nu}_{N}$ an orthogonal vector (for the metric $m$) to the hypersurface $N_{t_1}^{t_2}$,  and by $dv^{(\bf{m})}_N$ a volume form on $N_{t_1}^{t_2} $ such that the divergence theorem applies.

\end{definition} 

\begin{definition}\label{definitionofwidetildeLfordivergencetheoremuse}
We define $ \hat{L}^{\nu} $ as a vector proportional to $ \widetilde{L}^{\nu} $ (defined in Definition \ref{definitionofLtilde}), and therefore orthogonal (with respect to the Minkowski metric $m$) to $N$ (that is defined in Definition \ref{definitionofNandofNtruncatedbtweentwots}), and such that $\hat{L}$ is oriented,  and with Euclidian length, in such a way that the stated divergence theorem in what follows hold true: see, for instance, Corollary \ref{WeightedenergyestimateusingHandnosmallnessyet}.
 \end{definition}

\begin{remark}
Based on the construction of $N(q_0)$ in Definition \ref{definitionofNandofNtruncatedbtweentwots}, we have that the exterior region $\overline{C} $ includes the region $\{(t, x) \;|\; q:= r- t \geq q_0 \}$\;.
\end{remark}

\begin{lemma}\label{comutingthetermthatcarriesthederivativeoftheweight}
We have
       \beaa
 && T_{ t t}^{(\bf{g})}  +  T_{  r t}^{(\bf{g})} \\
 &=&  \frac{1}{2} \Big(  |\derm_t  \Phi_{V} + \derm_r \Phi_{V} |^2     +  \de^{ij}  < ( \derm_i - \frac{x_i}{r} \derm_{r}  )\Phi_{V} , (\derm_j - \frac{x_j}{r} \derm_{r} ) \Phi_{V} > \Big)   \\
  && -2 H^{\underline{L}\a}< \derm_\a \Phi_{V} ,   \derm_t \Phi_{V} >  + \frac{1}{2}   \cdot H^{\a\b}  < \derm_\a \Phi_{V} ,   \derm_\b \Phi_{V} >  \\
 \eeaa

  \end{lemma}

\begin{proof}

We recall that
\beaa
L &=&  \pa_{t} + \pa_{r} \, =   \pa_{t} +  \frac{x^{i}}{r} \pa_{i} \, . 
\eeaa
Thus,
\bea
\notag
 T_{ t t}^{(\bf{g})}  +  T_{  r t}^{(\bf{g})} &=&T_{  L t}^{(\bf{g})}  \; . 
  \eea
Consequently, we compute
  \beaa
\notag
 T_{ t t}^{(\bf{g})}  +  T_{  r t}^{(\bf{g})} &=&  m_{L \underline{L}} \cdot T^{(\bf{g}) \; \underline{L}}_{\;\;\;\;\;\;\;\;\;\; t}   \; .
\eeaa
  Since $m_{L \underline{L}}  = -2$, we get
 
    \bea
 T_{ t t}^{(\bf{g})}  +  T_{  r t}^{(\bf{g})} &=&  - 2 \cdot T^{(\bf{g}) \; \underline{L}}_{\;\;\;\;\;\;\;\;\;\; t}   \; .
  \eea

Given the expression of the energy momentum-tensor, we obtain

  \bea
\notag
T^{(\bf{g}) \; \underline{L}}_{\;\;\;\;\;\;\;\;\;\; t}   =  g^{\underline{L}\a}< \derm_\a \Phi_{V} ,   \derm_t \Phi_{V} > - \frac{1}{2} m^{\underline{L}}_{\;\;\; t}  \cdot g^{\a\b}  < \derm_\a \Phi_{V} ,   \derm_\b \Phi_{V} > \;, \\
 \eea
  
We compute
  \beaa
  m^{\underline{L}}_{\;\;\; t}  = m^{\underline{L} \mu} \cdot m_{\mu t } =  m^{\underline{L} L } \cdot m_{L t } = m^{\underline{L} L } \cdot m_{t t } = \frac{1}{2} \; .
  \eeaa

Hence,
  \beaa
\notag
T^{(\bf{g}) \; \underline{L}}_{\;\;\;\;\;\;\;\;\;\; t}   &=&  g^{\underline{L}\a}< \derm_\a \Phi_{V} ,   \derm_t \Phi_{V} > - \frac{1}{4}   \cdot g^{\a\b}  < \derm_\a \Phi_{V} ,   \derm_\b \Phi_{V} > \\
 &=&  m^{\underline{L}\a}< \derm_\a \Phi_{V} ,   \derm_t \Phi_{V} > - \frac{1}{4}   \cdot m^{\a\b}  < \derm_\a \Phi_{V} ,   \derm_\b \Phi_{V} >  \\
  && +  H^{\underline{L}\a}< \derm_\a \Phi_{V} ,   \derm_t \Phi_{V} > - \frac{1}{4}   \cdot H^{\a\b}  < \derm_\a \Phi_{V} ,   \derm_\b \Phi_{V} >  \; ,
 \eeaa
 and therefore,
     \beaa
 T_{ t t}^{(\bf{g})}  +  T_{  r t}^{(\bf{g})} &=&  - 2 \cdot T^{(\bf{g}) \; \underline{L}}_{\;\;\;\;\;\;\;\;\;\; t}   \\
 &=& -2 m^{\underline{L}\a}< \derm_\a \Phi_{V} ,   \derm_t \Phi_{V} > + \frac{1}{2}   \cdot m^{\a\b}  < \derm_\a \Phi_{V} ,   \derm_\b \Phi_{V} >  \\
  && -2 H^{\underline{L}\a}< \derm_\a \Phi_{V} ,   \derm_t \Phi_{V} >  + \frac{1}{2}   \cdot H^{\a\b}  < \derm_\a \Phi_{V} ,   \derm_\b \Phi_{V} >  \; .
 \eeaa

However, we already computed in \cite{G5}, the part with the Minkowski metric $m$ and therefore, based on \cite{G5}, we have
 
 \beaa
 && -2 m^{\underline{L}\a}< \derm_\a \Phi_{V} ,   \derm_t \Phi_{V} > + \frac{1}{2}   \cdot m^{\a\b}  < \derm_\a \Phi_{V} ,   \derm_\b \Phi_{V} >  \\
  &=&  \frac{1}{2} \Big(  < \derm_t  \Phi_{V} + \derm_r \Phi_{V} , \derm_t  \Phi_{V} + \derm_r \Phi_{V} >    \\
     \notag
  && +  \de^{ij}  < ( \derm_i - \frac{x_i}{r} \derm_{r}  )\Phi_{V} , (\derm_j - \frac{x_j}{r} \derm_{r} ) \Phi_{V} > \Big)   \; .
  \eeaa
Hence, we get the announced result.
 \end{proof}
 
 We recapitulate the following corollary from \cite{G5}.
    
     \begin{corollary}\label{WeightedenergyestimateusingHandnosmallnessyet}
 For $n$ space-dimensions, we have
    \beaa  \notag
  && \int_{N_{t_1}^{t_2} }  T_{\hat{L} t}^{(\bf{g})}  \cdot  \widetilde{w} (q) \cdot dv^{(\bf{m})}_N  \\
  &&+   \int_{\Sigma^{ext}_{t_2} }  \big(    - \frac{1}{2} (m^{t t} + H^{tt} ) < \derm_t \Phi_{V} ,   \derm_t \Phi_{V}>   + \frac{1}{2} (m^{j i} +H^{j i} )   < \derm_j \Phi_{V} ,   \derm_i \Phi_{V} > \big)  \cdot   \widetilde{w} (q) \cdot d^{n}x \\
   \notag
     &=&   \int_{\Sigma^{ext}_{t_1} } \big(    - \frac{1}{2} (m^{t t} + H^{tt} ) < \derm_t \Phi_{V} ,   \derm_t \Phi_{V} >   + \frac{1}{2} (m^{j i} +H^{j i} )   < \derm_j \Phi_{V} ,   \derm_i \Phi_{V} > \big)   \cdot   \widetilde{w} (q)  \cdot d^{n}x  \\
     && -  \int_{t_1}^{t_2}  \int_{\Sigma^{ext}_{\tau} }  \Big( T_{t  t}^{(\bf{g})} +   T_{r  t}^{(\bf{g})} \Big) \cdot d\tau \cdot   \widetilde{w} ^\prime (q) d^{n}x \\
     \notag
   &&  -  \int_{t_1}^{t_2}  \int_{\Sigma^{ext}_{\tau} } \Big(  < g^{\mu\a} \derm_{\mu } \derm_\a \Phi_{V} ,   \derm_t \Phi_{V} >  +( \derm_{\mu } H^{\mu\a} ) \cdot < \derm_\a \Phi_{V} ,   \derm_t \Phi_{V} > \\
 \notag
&&- \frac{1}{2} m^{\mu}_{\;\;\; t}  \cdot  ( \derm_{\mu } H^{\a\b} ) \cdot < \derm_\a \Phi_{V} ,   \derm_\b \Phi_{V} >   \Big) \cdot d\tau \cdot   \widetilde{w} (q) d^{n}x \; .  \\
 \eeaa
 \end{corollary}

  \begin{lemma}\label{AssumptiononHforgettingthestandardnormofonecomponentsintheexpressionofenergyestimate}
  
For $ | H| \leq C < \frac{1}{n} $\,, where $n$ is the space dimension, for any vector $U\,,\, V$, we have
 \beaa
 \notag
   && - \frac{1}{2} (m^{t t} + H^{tt} ) < \derm_t \Phi_{V} ,   \derm_t \Phi_{V} >   + \frac{1}{2} (m^{j i} +H^{j i} )   < \derm_j \Phi_{V} ,   \derm_i \Phi_{V} >  \\
   &\sim & |\derm \Phi_{V} |^2 \geq 0 \; . 
  \eeaa
  
  \end{lemma}
  
  \begin{proof}
   We showed in \cite{G4}, that for $ | H| \leq C < \frac{1}{n} $\,, we have
 \beaa
 \notag
   && - \frac{1}{2} (m^{t t} + H^{tt} ) < \pa_t \Phi ,   \pa_t \Phi >   + \frac{1}{2} (m^{j i} +H^{j i} )   < \pa_j \Phi ,  \pa_i \Phi >  \\
   &\sim & |\pa \Phi |^2 \geq 0 \; , 
  \eeaa
 where for $\a\, ,\, \b \in \{ r, t, x^1, \ldots, x^n \}$\,,
\bea\label{definitionofscalarproductofpartialderivativeindirectionofwavecoordinatesoftensor}
\notag
 < \pa_\a  \Phi ,\pa_\b \Phi > &:=& \sum_{\mu,\, \nu \in \{ t, x^1, \ldots, x^n \}  } < \pa_\a  \Phi_{\mu\nu} ,\pa_\b \Phi_{\mu\nu} > \; . \\
 \eea
 However, the proof in \cite{G4} works also if we replace the partial derivative with the covariant derivative, and if we replace the scalar product in \eqref{definitionofscalarproductofpartialderivativeindirectionofwavecoordinatesoftensor} by the standard scalar product of the components $\derm_\a  \Phi_{V}$ and $\derm_\b  \Phi_{V}$\;.

  \end{proof}

Now, fix fix the space-dimension as $n=3$\;.

 \begin{lemma}\label{energyestimatewithoutestimatingthetermsthatinvolveBIGHbutbydecomposingthemcorrectlysothatonecouldgettherightestimatewithtildew}
For
\bea
| H| < \frac{1}{3} \; ,
\eea
and for $\Phi$ decaying sufficiently fast at spatial infinity, we have the following energy estimate
    \bea
   \notag
 &&     \int_{\Sigma^{ext}_{t_2} }  |\derm \Phi_{V} |^2     \cdot   \widetilde{w} (q)  \cdot d^{3}x    + \int_{N_{t_1}^{t_2} } T_{\hat{L} t}^{(\bf{g})}   \cdot    \widetilde{w} (q) \cdot dv^{(\bf{m})}_N \\
 \notag
 &&+ \int_{t_1}^{t_2}  \int_{\Sigma^{ext}_{\tau} }  \Big(    \frac{1}{2} \Big(  | \derm_t  \Phi_{V} + \derm_r \Phi_{V} |^2  +  \sum_{i=1}^3 | ( \derm_i - \frac{x_i}{r} \derm_{r}  )\Phi_{V} |^2 \Big)  \cdot d\tau \cdot   \widetilde{w} ^\prime (q) d^{3}x \\
  \notag
  &\les &       \int_{\Sigma^{ext}_{t_1} }  |\derm \Phi |^2     \cdot   \widetilde{w} (q)  \cdot d^{3}x \\
    \notag
     && +  \int_{t_1}^{t_2}  \int_{\Sigma^{ext}_{\tau} } \Big(  | H_{LL } | \cdot | \derm \Phi_{V} |^2 + |  H | \cdot | \rderm \Phi_{V} | \cdot | \derm \Phi_{V}  |  \Big)  \cdot d\tau \cdot   \widetilde{w} ^\prime (q) d^{3}x \\
     \notag
   &&  +  \int_{t_1}^{t_2}  \int_{\Sigma^{ext}_{\tau} }  |  g^{\mu\a} \derm_{\mu } \derm_\a \Phi_{V} | \cdot |  \derm_t \Phi_{V} |  \cdot d\tau \cdot   \widetilde{w} (q) d^{3}x  \\
 \notag
&& + \Big( ( | \derm H_{LL} |  + |\rderm H| ) \cdot | \derm \Phi_{V} |^2 +  | \derm H | \cdot  | \rderm \Phi_{V} | \cdot  | \derm \Phi_{V} | \Big)  \cdot d\tau \cdot   \widetilde{w} (q) d^{3}x \; .  \\
 \eea
 
 \end{lemma}
 
 \begin{proof}
 
 Injecting the result of Lemma \ref{comutingthetermthatcarriesthederivativeoftheweight} in Corollary \ref{WeightedenergyestimateusingHandnosmallnessyet}, we obtain

     \beaa  \notag
  && \int_{N_{t_1}^{t_2} }  T_{\hat{L} t}^{(\bf{g})}  \cdot  \widetilde{w}(q) \cdot dv^{(\bf{m})}_N  \\
  &&+   \int_{\Sigma^{ext}_{t_2} }  \big(    - \frac{1}{2} (m^{t t} + H^{tt} ) < \derm_t \Phi_{V} ,   \derm_t \Phi_{V} >   + \frac{1}{2} (m^{j i} +H^{j i} )   < \derm_j \Phi_{V} ,   \derm_i \Phi_{V} > \big)  \cdot   \widetilde{w} (q) \cdot d^{3}x \\
  \notag
&&    \int_{t_1}^{t_2}  \int_{\Sigma^{ext}_{\tau} }  \Big(   \frac{1}{2} \Big(  |\derm_t  \Phi_{V} + \derm_r \Phi_{V} |^2     +  \de^{ij}  < ( \derm_i - \frac{x_i}{r} \derm_{r}  )\Phi_{V} , (\derm_j - \frac{x_j}{r} \derm_{r} ) \Phi_{V} >   \Big)\\
&&  \quad \quad \quad \quad \quad \times d\tau \cdot   \widetilde{w} ^\prime (q) d^{3}x \\ 
   \notag
     &=&   \int_{\Sigma^{ext}_{t_1} } \big(    - \frac{1}{2} (m^{t t} + H^{tt} ) < \derm_t \Phi_{V} ,   \derm_t \Phi_{V} >   + \frac{1}{2} (m^{j i} +H^{j i} )   < \derm_j \Phi_{V} ,   \derm_i \Phi_{V} > \big)   \cdot   \widetilde{w} (q)  \cdot d^{3}x  \\
     && +  \int_{t_1}^{t_2}  \int_{\Sigma^{ext}_{\tau} }  \Big(  2 H^{\underline{L}\a}< \derm_\a \Phi_{V} ,   \derm_t \Phi_{V} >  - \frac{1}{2}   \cdot H^{\a\b}  < \derm_\a \Phi_{V} ,   \derm_\b \Phi_{V} > \Big) \cdot d\tau \cdot   \widetilde{w}^\prime (q) d^{3}x \\
     \notag
   &&  -  \int_{t_1}^{t_2}  \int_{\Sigma^{ext}_{\tau} } \Big(  < g^{\mu\a} \derm_{\mu } \derm_\a \Phi_{V} ,   \derm_t \Phi_{V} >  +( \derm_{\mu } H^{\mu\a} ) \cdot < \derm_\a \Phi_{V} ,   \derm_t \Phi_{V} > \\
 \notag
&&- \frac{1}{2} m^{\mu}_{\;\;\; t}  \cdot  ( \derm_{\mu } H^{\a\b} ) \cdot < \derm_\a \Phi_{V} ,   \derm_\b \Phi_{V} >   \Big) \cdot d\tau \cdot   \widetilde{w} (q) d^{3}x \; .  \\
 \eeaa
 
 Using Lemma \ref{AssumptiononHforgettingthestandardnormofonecomponentsintheexpressionofenergyestimate}, we get for $ | H| \leq C < \frac{1}{3} $\,, the following energy estimate

   \bea
   \notag
 &&     \int_{\Sigma^{ext}_{t_2} }  |\derm \Phi_{V} |^2     \cdot w(q)  \cdot d^{3}x    + \int_{N_{t_1}^{t_2} }  T_{\hat{L} t}^{(\bf{g})}   \cdot    \widetilde{w} (q) \cdot dv^{(\bf{m})}_N \\
 \notag
 &&+ \int_{t_1}^{t_2}  \int_{\Sigma^{ext}_{\tau} }  \Big(    \frac{1}{2} \Big(  | \derm_t  \Phi_{V} + \derm_r \Phi_{V} |^2  +   \sum_{i=1}^3  | ( \derm_i - \frac{x_i}{r} \derm_{r}  )\Phi_{V} |^2 \Big)  \cdot d\tau \cdot   \widetilde{w}^\prime (q) d^{3}x \\
  \notag
  &\leq &       \int_{\Sigma^{ext}_{t_1} }  |\derm \Phi_{V} |^2     \cdot \widetilde{w}(q)  \cdot d^{3}x \\
    \notag
     && +  \int_{t_1}^{t_2}  \int_{\Sigma^{ext}_{\tau} }  | 2 H^{\underline{L}\a}< \derm_\a \Phi_{V} ,   \derm_t \Phi_{V} >  - \frac{1}{2}   \cdot H^{\a\b}  < \derm_\a \Phi_{V} ,   \derm_\b \Phi_{V} >  | \cdot d\tau \cdot \widetilde{w}^\prime (q) d^{3}x \\
     \notag
   &&  +  \int_{t_1}^{t_2}  \int_{\Sigma^{ext}_{\tau} }  | < g^{\mu\a} \derm_{\mu } \derm_\a \Phi_{V} ,   \derm_t \Phi_{V} >  +( \derm_{\mu } H^{\mu\a} ) \cdot < \derm_\a \Phi_{V} ,   \derm_t \Phi_{V} > \\
 \notag
&&- \frac{1}{2} m^{\mu}_{\;\;\; t}  \cdot  ( \derm_{\mu } H^{\a\b} ) \cdot < \derm_\a \Phi_{V} ,   \derm_\b \Phi_{V} >   |  \cdot d\tau \cdot \widetilde{w}(q) d^{3}x \; .  \\
 \eea

 Since $ m^{\mu}_{\;\;\; t} = m^{\mu\a} \cdot m_{\a t}  =m^{\mu t} \cdot m_{t t} = - m^{\mu t} $, we estimate
\beaa\label{estimateonthedivergenceoftheenergymomentumtensordecomposedinnullframe}
 && | < g^{\mu\a} \derm_{\mu } \derm_\a \Phi_{V} ,   \derm_t \Phi_{V} >  +( \derm_{\mu } H^{\mu\a} ) \cdot < \derm_\a \Phi_{V} ,   \derm_t \Phi_{V} > \\
 \notag
&&- \frac{1}{2} m^{\mu}_{\;\;\; t}  \cdot  ( \derm_{\mu } H^{\a\b} ) \cdot < \derm_\a \Phi_{V} ,   \derm_\b \Phi_{V} > | \\
 &\les& |  g^{\mu\a} \derm_{\mu } \derm_\a \Phi | \cdot |  \derm_t \Phi_{V}  | \\
\notag
&& +| ( \derm_{\mu } H^{\mu\a} ) \cdot < \derm_\a \Phi_{V} ,   \derm_t \Phi_{V} > | \\
 \notag
&&+ | \frac{1}{2} m^{t t} \cdot  ( \derm_{t} H^{\a\b} ) \cdot < \derm_\a \Phi_{V} ,   \derm_\b \Phi_{V} >  |  
\notag
\eeaa

We have by decomposing simply in a null frame, or by equivalently by using Lemma 4.2 from \cite{LR10}, that 
 \bea
\notag
&& | ( \derm_{\mu } H^{\mu\a} ) \cdot < \derm_\a \Phi_{V} ,   \derm_t \Phi_{V} > | \\
 \notag
&&+ | \frac{1}{2} m^{t t} \cdot  ( \derm_{t} H^{\a\b} ) \cdot < \derm_\a \Phi_{V} ,   \derm_\b \Phi_{V} >  | \\
 \notag
&\les&  | ( \derm_{\mu } H^{\mu\a} ) \cdot < \derm_\a \Phi_{V} ,   \derm_t \Phi_{V} > | \\
 \notag
&&+ |  ( \derm_{t} H^{\a\b} ) \cdot < \derm_\a \Phi_{V} ,   \derm_\b \Phi_{V} >  | \\
 \notag
&\les& ( | \derm H_{LL} |  + |\rderm H| ) \cdot | \derm \Phi_{V} |^2 +  | \derm H | \cdot  | \rderm \Phi_{V} | \cdot  | \derm \Phi_{V} | \; .\\
\eea

Similarly, by decomposing in null frame, as in Lemma 4.2 in \cite{LR10}, we get
   \bea
\notag
&& | 2 H^{\underline{L}\a}< \derm_\a \Phi_{V} ,   \derm_t \Phi_{V} >  - \frac{1}{2}   \cdot H^{\a\b}  < \derm_\a \Phi_{V} ,   \derm_\b \Phi_{V} >  | \\
\notag
  &\leq& | H^{\underline{L}\a}< \derm_\a \Phi_{V} ,   \derm_t \Phi_{V} > |  +| H^{\a\b}  < \derm_\a \Phi_{V} ,   \derm_\b \Phi_{V} > |  \\
  \notag
 &\les& | H_{LL } | \cdot | \derm \Phi_{V} |^2 + |  H | \cdot | \rderm \Phi_{V} | \cdot | \derm \Phi_{V}  | \; . \\
 \eea

Injecting, we obtain the desired result.
\end{proof}

\begin{corollary}\label{energyestimatewithoutestimatingthetermsthatinvolveBIGHbutbydecomposingthemcorrectlysothatonecouldgettherightestimate}
For
\bea
| H| < \frac{1}{3} \; ,
\eea
and for $\ga > 0$ and $\mu < 0$\,, and for $\Phi$ decaying sufficiently fast at spatial infinity, we have the following energy estimate
    \bea
   \notag
 &&     \int_{\Sigma^{ext}_{t_2} }  |\derm \Phi_{V} |^2     \cdot   w (q)  \cdot d^{3}x    + \int_{N_{t_1}^{t_2} }   T_{\hat{L} t}^{(\bf{g})}   \cdot   w (q) \cdot dv^{(\bf{m})}_N \\
 \notag
 &&+ \int_{t_1}^{t_2}  \int_{\Sigma^{ext}_{\tau} }  \Big(    \frac{1}{2} \Big(  | \derm_t  \Phi_{V} + \derm_r \Phi_{V} |^2  +  \sum_{i=1}^3 | ( \derm_i - \frac{x_i}{r} \derm_{r}  )\Phi_{V} |^2 \Big)  \cdot d\tau \cdot   \widehat{w}^\prime (q) d^{3}x \\
  \notag
  &\les &       \int_{\Sigma^{ext}_{t_1} }  |\derm \Phi_{V} |^2     \cdot   w(q)  \cdot d^{3}x \\
    \notag
     && +  \int_{t_1}^{t_2}  \int_{\Sigma^{ext}_{\tau} } \Big(  | H_{LL } | \cdot | \derm \Phi_{V} |^2 + |  H | \cdot | \rderm \Phi_{V} | \cdot | \derm \Phi_{V}  |  \Big)  \cdot d\tau \cdot   \widehat{w} ^\prime (q) d^{3}x \\
     \notag
   &&  +  \int_{t_1}^{t_2}  \int_{\Sigma^{ext}_{\tau} }  |  g^{\mu\a} \derm_{\mu } \derm_\a \Phi_{V} | \cdot |  \derm_t \Phi_{V} |  \cdot d\tau \cdot  w(q) d^{3}x  \\
 \notag
&& + \Big( ( | \derm H_{LL} |  + |\rderm H| ) \cdot | \derm \Phi_{V} |^2 +  | \derm H | \cdot  | \rderm \Phi_{V} | \cdot  | \derm \Phi_{V} | \Big)  \cdot d\tau \cdot   w (q) d^{3}x \; .  \\
 \eea
 
 \end{corollary}

\begin{proof}
Using Lemmas \ref{equaivalenceoftildewandtildeandofderivativeoftildwandderivativeofhatw} and \ref{derivativeoftildwandrelationtotildew} and injecting in Lemma \ref{energyestimatewithoutestimatingthetermsthatinvolveBIGHbutbydecomposingthemcorrectlysothatonecouldgettherightestimatewithtildew}, we obtain the desired result.
\end{proof}

We want to estimate $ | H_{LL } | \cdot | \derm \Phi_{V} |^2 + |  H | \cdot | \rderm \Phi_{V} | \cdot | \derm \Phi_{V}  |  + ( | \derm H_{LL} |  + |\rderm H| ) \cdot | \derm \Phi_{V} |^2 +  | \derm H | \cdot  | \rderm \Phi_{V} | \cdot  | \derm \Phi_{V} |$\,. We start first by stating the following lemma.

 \begin{lemma} \label{upgradedestimateonBIGH}
 
 We have for $\ga \geq 3 \de $\,, and $0 < \de \leq \frac{1}{4}$\,, that in the exterior region, the following estimates hold for $M \leq \eps$\,, and for all $I$\,, 
                 \beaa
 \notag
 |    \Lie_{Z^I}  H   | &\leq&   C(q_0)   \cdot  c (\delta) \cdot c (\gamma) \cdot C(|I|)  \cdot E (   |I|  + 4)  \cdot \frac{\eps }{(1+t+|q|)^{1-   c (\gamma)  \cdot c (\delta)  \cdot c(|I|) \cdot E (  |I|+ 4)\cdot  \eps } }     \; , \\
      \eeaa
      and 
                    \beaa
 \notag
 && |\derm  ( \Lie_{Z^I} H )   | \\
  &\leq&   C(q_0)   \cdot  c (\delta) \cdot c (\gamma) \cdot C(|I|) \cdot E (  |I|  +4)  \cdot \frac{\eps }{(1+t+|q|)^{1-   c (\gamma)  \cdot c (\delta)  \cdot c(|I|) \cdot E ( |I|+ 4)\cdot  \eps } \cdot (1+|q|)}     \; . 
      \eeaa
 \end{lemma}
 
 \begin{proof}
Following exactly the same argument as in Lemmas \ref{aprioriestimatesonZLiederivativesofBIGH} and \ref{aprioriestimatesonZLiederivativesonbigH}, we get based on Lemma \ref{upgradedestimatesonh}, the desired result.
\end{proof}
 
 \begin{lemma}\label{estimatesonthetermsthatcontainbigHandderivativesofBigHintheenergyestimateinawaythatwecouldconcludeforneq3}
We have for $\ga \geq 3 \de $\,, and $0 < \de \leq \frac{1}{4}$\,, that in the exterior region, the following estimates hold for $M \leq \eps$\,, and for all $I$\,, 
 \beaa
   && | H_{LL } | \cdot | \derm \Phi_{V} |^2 + |  H | \cdot | \rderm \Phi_{V} | \cdot | \derm \Phi_{V}  |  \\
&\les&  C(q_0) \cdot   c (\delta) \cdot c (\gamma) \cdot E (  4)  \cdot \frac{ \eps }{ (1+t+|q|) } \cdot | \derm \Phi_{V} |^2   \\
    && +  C(q_0)   \cdot  c (\delta) \cdot c (\gamma)  \cdot E ( 4 ) \cdot \frac{\eps  }{(1+t+|q|)^{1-   c (\gamma)  \cdot c (\delta)  \cdot E (  4)\cdot  \eps } }   \cdot | \rderm \Phi_{V} |^2 \; , 
 \eeaa
and
        \beaa
   && ( | \derm H_{LL} |  + |\rderm H| ) \cdot | \derm \Phi_{V} |^2 +  | \derm H | \cdot  | \rderm \Phi_{V} | \cdot  | \derm \Phi_{V} | \\
&\les&  C(q_0) \cdot   c (\delta) \cdot c (\gamma) \cdot E (  4)  \cdot \frac{ \eps }{ (1+t+|q|) \cdot (1+|q|)  } \cdot | \derm \Phi_{V}|^2   \\
    && + C(q_0)   \cdot  c (\delta) \cdot c (\gamma)  \cdot E ( 4 ) \cdot \frac{\eps  }{(1+t+|q|)^{1-   c (\gamma)  \cdot c (\delta)  \cdot E (  4)\cdot  \eps }  \cdot (1+|q|)  }   \cdot | \rderm \Phi_{V} |^2  \; .  
 \eeaa
 
 \end{lemma}
 
 \begin{proof}
 
 \textbf{The term $ | H_{LL } | \cdot | \derm \Phi_{V} |^2 + |  H | \cdot | \rderm \Phi_{V} | \cdot | \derm \Phi_{V}  | $):}\\
 
By taking $|I| = 0$ in Lemma \ref{estimategoodcomponentspotentialandmetric}, we have in the exterior region, 
                \bea
        \notag
 |  H_{L L} |                    &\les&   C(q_0) \cdot   c (\delta) \cdot c (\gamma) \cdot E (  3)  \cdot \frac{ \eps   }{ (1+t+|q|) } \; . \\
       \eea
       
       Thus,
    \beaa
\notag
&& | H_{LL } | \cdot | \derm \Phi_{V} |^2 + |  H | \cdot | \rderm \Phi_{V} | \cdot | \derm \Phi_{V}  |  \\
&\les&  C(q_0) \cdot   c (\delta) \cdot c (\gamma) \cdot E (  3)  \cdot \frac{ \eps \cdot | \derm \Phi_{V} |^2  }{ (1+t+|q|) }  + |  H | \cdot | \rderm \Phi_{V} | \cdot | \derm \Phi_{V}  |  
 \eeaa

 However, from Lemma \ref{upgradedestimateonBIGH}, we have
                  \beaa
 \notag
 |   H   | &\leq&   C(q_0)   \cdot  c (\delta) \cdot c (\gamma)  \cdot E ( 4)  \cdot \frac{\eps }{(1+t+|q|)^{1-   c (\gamma)  \cdot c (\delta)  \cdot E (  4)\cdot  \eps } }     \; . \\
      \eeaa

Consequently,
 
     \beaa
\notag
&& | H_{LL } | \cdot | \derm \Phi_{V} |^2 + |  H | \cdot | \rderm \Phi_{V} | \cdot | \derm \Phi_{V}  |  \\
&\les&  C(q_0) \cdot   c (\delta) \cdot c (\gamma) \cdot E (  3)  \cdot \frac{ \eps \cdot | \derm \Phi_{V} |^2  }{ (1+t+|q|) } \\
&& +  C(q_0)   \cdot  c (\delta) \cdot c (\gamma)  \cdot E ( 4)  \cdot \frac{\eps }{(1+t+|q|)^{1-   c (\gamma)  \cdot c (\delta)  \cdot E (  4)\cdot  \eps } }     \cdot | \rderm \Phi_{V} | \cdot | \derm \Phi_{V}  |  
 \eeaa
 
  We estimate the second term in the right hand side of the inequality, by using $a \cdot b \les a^2 + b^2 $, and we obtain
  \beaa
  &&   C(q_0)   \cdot  c (\delta) \cdot c (\gamma)  \cdot E ( 4)  \cdot \frac{\eps }{(1+t+|q|)^{1-   c (\gamma)  \cdot c (\delta)  \cdot E (  4)\cdot  \eps } }     \cdot | \rderm \Phi_{V} | \cdot | \derm \Phi_{V}  |  \\
  &=& \big( C(q_0)   \cdot  c (\delta) \cdot c (\gamma)  \cdot E ( 4 ) \big)^{\frac{1}{2} }\cdot \frac{\sqrt{\eps} \cdot \sqrt{(1+t+|q|)} }{(1+t+|q|)^{1-   c (\gamma)  \cdot c (\delta)  \cdot E (  4)\cdot  \eps } }   \cdot | \rderm \Phi_{V} |   \\
  && \times \big( C(q_0)   \cdot  c (\delta) \cdot c (\gamma)  \cdot E ( 4) )^{\frac{1}{2} }\cdot \frac{\sqrt{\eps} }{\sqrt{(1+t+|q|)} } \big) \cdot | \derm \Phi_{V}  |   \\
    &\les&  C(q_0)   \cdot  c (\delta) \cdot c (\gamma)  \cdot E ( 4) ) \cdot \frac{\eps  }{(1+t+|q|)^{1-   c (\gamma)  \cdot c (\delta)  \cdot E (  4)\cdot  \eps } }   \cdot | \rderm \Phi_{V} |^2 \\
    && + C(q_0)   \cdot  c (\delta) \cdot c (\gamma)  \cdot E ( 4 ) \cdot \frac{\eps  }{(1+t +|q|) }  \cdot | \derm \Phi_{V} |^2
 \eeaa
 
 Finally, putting all together, we obtain
    \beaa
\notag
&&  | H_{LL } | \cdot | \derm \Phi_{V} |^2 + |  H | \cdot | \rderm \Phi_{V} | \cdot | \derm \Phi_{V}  |  \\
&\les&  C(q_0) \cdot   c (\delta) \cdot c (\gamma) \cdot E (  4)  \cdot \frac{ \eps }{ (1+t+|q|) } \cdot | \derm \Phi_{V} |^2   \\
    && +  C(q_0)   \cdot  c (\delta) \cdot c (\gamma)  \cdot E ( 4 ) \cdot \frac{\eps  }{(1+t+|q|)^{1-   c (\gamma)  \cdot c (\delta)  \cdot E (  4)\cdot  \eps } }   \cdot | \rderm \Phi_{V} |^2 \\
 \eeaa

\textbf{The term $( | \derm H_{LL} |  + |\rderm H| ) \cdot | \derm \Phi_{V} |^2 +  | \derm H | \cdot  | \rderm \Phi_{V} | \cdot  | \derm \Phi_{V} |$):}\\

Based on Lemma \ref{waveconditionestimateonzeroLiederivativeofmetric}, we have
\beaa
| \derm  H_{L L} | &\les& | \rderm  H |  + O (|H| \cdot |\derm H| ) \; .
 \eeaa
 Injecting the estimates from Lemma \ref{aprioriestimatesonZLiederivativesonbigH}, we obtain that for $\gamma > \delta  $\,, we have for all $|I|$\,, $\delta \leq \frac{1}{2}$\,,\,$\eps \leq 1$\,, in the exterior region
 \beaa
| \derm  H_{L L} | &\les&c (\delta) \cdot c (\gamma) \cdot E (  3)  \cdot \frac{\eps }{(1+t+|q|)^{2-\delta} (1+|q|)^{\de}} \\
&\les&c (\delta) \cdot c (\gamma) \cdot E (  3)  \cdot \frac{\eps }{(1+t+|q|) \cdot (1+|q|) } \\
 \eeaa

Thus, based on Lemma \ref{upgradedestimateonBIGH}, we get
      \beaa
   && ( | \derm H_{LL} |  + |\rderm H| ) \cdot | \derm \Phi_{V} |^2 +  | \derm H | \cdot  | \rderm \Phi_{V} | \cdot  | \derm \Phi_{V} | \\
&\les&C(q_0)   \cdot  c (\delta) \cdot c (\gamma)  \cdot E (    3)  \cdot \frac{\eps }{(1+t+|q|) \cdot (1+|q|) }  \cdot | \derm \Phi_{V} |^2 \\
&& +   C(q_0)   \cdot  c (\delta) \cdot c (\gamma) \cdot E (  4)  \cdot \frac{\eps }{(1+t+|q|)^{1-   c (\gamma)  \cdot c (\delta)  \cdot E (  4)\cdot  \eps } \cdot (1+|q|)}  \cdot  | \rderm \Phi_{V} | \cdot  | \derm \Phi | \\
&\les& C(q_0)   \cdot  c (\delta) \cdot c (\gamma)  \cdot E (    4)  \cdot \frac{\eps }{(1+t+|q|) \cdot (1+|q|)  }  \cdot | \derm \Phi_{V} |^2 \\
&& +   \frac{1}{(1+|q|)}  \Big(  C(q_0)   \cdot  c (\delta) \cdot c (\gamma) \cdot E (  4)  \cdot \frac{\eps }{(1+t+|q|)^{1-   c (\gamma)  \cdot c (\delta)  \cdot E (  4)\cdot  \eps } }  \cdot  | \rderm \Phi_{V} | \cdot  | \derm \Phi_{V} | \Big) 
\eeaa

Following exactly the same argument that we used to estimate the term $ | H_{LL } | \cdot | \derm \Phi_{V} |^2 + |  H | \cdot | \rderm \Phi_{V} | \cdot | \derm \Phi_{V}  | $, to estimate the second term on the right hand side of the inequality, we obtain

\beaa
&&    \frac{1}{(1+|q|)}  \Big(  C(q_0)   \cdot  c (\delta) \cdot c (\gamma) \cdot E (  4)  \cdot \frac{\eps }{(1+t+|q|)^{1-   c (\gamma)  \cdot c (\delta)  \cdot E (  4)\cdot  \eps } }  \cdot  | \rderm \Phi_{V} | \cdot  | \derm \Phi_{V} | \Big) \\
  &\les&  C(q_0)   \cdot  c (\delta) \cdot c (\gamma)  \cdot E ( 4 ) \cdot \frac{\eps  }{(1+t+|q|)^{1-   c (\gamma)  \cdot c (\delta)  \cdot E (  4)\cdot  \eps }  \cdot (1+|q|)  }   \cdot | \rderm \Phi_{V} |^2 \\
    && + C(q_0)   \cdot  c (\delta) \cdot c (\gamma)  \cdot E ( 4 ) \cdot \frac{\eps  }{(1+t +|q|) \cdot (1+|q|) }  \cdot | \derm \Phi_{V} |^2
\eeaa
 Thus, we get the stated result.
 
        \beaa
   && ( | \derm H_{LL} |  + |\rderm H| ) \cdot | \derm \Phi_{V} |^2 +  | \derm H | \cdot  | \rderm \Phi_{V} | \cdot  | \derm \Phi_{V} | \\
&\les&  C(q_0) \cdot   c (\delta) \cdot c (\gamma) \cdot E (  4)  \cdot \frac{ \eps }{ (1+t+|q|) \cdot (1+|q|)  } \cdot | \derm \Phi_{V} |^2   \\
    && + C(q_0)   \cdot  c (\delta) \cdot c (\gamma)  \cdot E ( 4 ) \cdot \frac{\eps  }{(1+t+|q|)^{1-   c (\gamma)  \cdot c (\delta)  \cdot E (  4)\cdot  \eps }  \cdot (1+|q|)  }   \cdot | \rderm \Phi_{V} |^2  \; .  
 \eeaa
\end{proof}
      
 Now, we would like to use the decay estimates on the metric $H$ to prove an suitable energy estimate that would allow us to close the bootstrap argument.

 \begin{corollary}\label{theenergyestimatewithcontrolontangentialderivativesusinghatwoveroneplusq}
For
\bea
| H| < \frac{1}{3} \; ,
\eea
and for $\ga > 0$ and $\mu < 0$\,, and for $\Phi$ a tensor of arbitrary order, say for simplicity of notation that it is a one-tensor, decaying sufficiently fast at spatial infinity, we have the following energy estimate for any vector $V$\;,
  \beaa
   \notag
 &&     \int_{\Sigma^{ext}_{t_2} }  |\derm \Phi_{V} |^2     \cdot w(q)  \cdot d^{3}x    + \int_{N_{t_1}^{t_2} }   T_{\hat{L} t}^{(\bf{g})}  (\Phi_{V}) \cdot  w(q) \cdot dv^{(\bf{m})}_N \\
 \notag
 &&+ \int_{t_1}^{t_2}  \int_{\Sigma^{ext}_{\tau} }  \Big(    \frac{1}{2} \Big(  | \derm_t  \Phi_{V} + \derm_r \Phi_{V} |^2  +  \de^{ij}  | ( \derm_i - \frac{x_i}{r} \derm_{r}  )\Phi_{V} |^2 \Big)  \cdot d\tau \cdot  \frac{\widehat{w} (q)}{(1+|q|)} d^{3}x \\
  \notag
  &\les &       \int_{\Sigma^{ext}_{t_1} }  |\derm \Phi_{V} |^2     \cdot w(q)  \cdot d^{3}x \\
    \notag
   &&  +  \int_{t_1}^{t_2}  \int_{\Sigma^{ext}_{\tau} }  |  g^{\mu\a} \derm_{\mu } \derm_\a \Phi_{V} | \cdot |  \derm_t \Phi_{V} |  \cdot d\tau \cdot w(q) d^{3}x  \\
 \notag
&& + \int_{t_1}^{t_2}  \int_{\Sigma^{ext}_{\tau} }  \Big(   C(q_0) \cdot   c (\delta) \cdot c (\gamma) \cdot E (  4)  \cdot \frac{ \eps }{ (1+\tau+|q|)\cdot (1+|q|)  } \cdot | \derm \Phi_{V} |^2   \\
    && +  C(q_0)   \cdot  c (\delta) \cdot c (\gamma)  \cdot E ( 4)  \cdot \frac{\eps  }{(1+\tau+|q|)^{1-   c (\gamma)  \cdot c (\delta)  \cdot E (  4)\cdot  \eps } \cdot (1+|q|) }   \cdot | \rderm \Phi_{V} |^2 \Big)  \cdot d\tau \cdot  w(q) d^{3}x \; . 
 \eeaa 
 \end{corollary}
 
 \begin{proof}
By injecting the result of Lemmas \ref{estimatesonthetermsthatcontainbigHandderivativesofBigHintheenergyestimateinawaythatwecouldconcludeforneq3} in Corollary \ref{energyestimatewithoutestimatingthetermsthatinvolveBIGHbutbydecomposingthemcorrectlysothatonecouldgettherightestimate}, we obtain
   \beaa
   \notag
 &&     \int_{\Sigma^{ext}_{t_2} }  |\derm \Phi_{V} |^2     \cdot w(q)  \cdot d^{3}x    + \int_{N_{t_1}^{t_2} }  T_{\hat{L} t}^{(\bf{g})}   \cdot  w(q) \cdot dv^{(\bf{m})}_N \\
 \notag
 &&+ \int_{t_1}^{t_2}  \int_{\Sigma^{ext}_{\tau} }  \Big(    \frac{1}{2} \Big(  | \derm_t  \Phi_{V} + \derm_r \Phi_{V} |^2  +  \de^{ij}  | ( \derm_i - \frac{x_i}{r} \derm_{r}  )\Phi_{V} |^2 \Big)  \cdot d\tau \cdot  \widehat{w}^\prime (q) d^{3}x \\
  \notag
  &\les &       \int_{\Sigma^{ext}_{t_1} }  |\derm \Phi_{V} |^2     \cdot w(q)  \cdot d^{3}x \\
    \notag
   &&  +  \int_{t_1}^{t_2}  \int_{\Sigma^{ext}_{\tau} }  |  g^{\mu\a} \derm_{\mu } \derm_\a \Phi_{V} | \cdot |  \derm_t \Phi_{V} |  \cdot d\tau \cdot w(q) d^{3}x  \\
 \notag
&& + \int_{t_1}^{t_2}  \int_{\Sigma^{ext}_{\tau} }  \Big(   C(q_0) \cdot   c (\delta) \cdot c (\gamma) \cdot E (  4)  \cdot \frac{ \eps }{ (1+\tau+|q|) } \cdot | \derm \Phi_{V} |^2   \\
    && +  C(q_0)   \cdot  c (\delta) \cdot c (\gamma)  \cdot E ( 4)  \cdot \frac{\eps  }{(1+\tau+|q|)^{1-   c (\gamma)  \cdot c (\delta)  \cdot E (  4)\cdot  \eps } }   \cdot | \rderm \Phi_{V} |^2 \Big)  \cdot d\tau \cdot  \widehat{w}^\prime (q)  d^{3}x  \\
     \notag
&& + \int_{t_1}^{t_2}  \int_{\Sigma^{ext}_{\tau} }  \Big(   C(q_0) \cdot   c (\delta) \cdot c (\gamma) \cdot E (  4)  \cdot \frac{ \eps }{ (1+\tau+|q|)\cdot (1+|q|)  } \cdot | \derm \Phi_{V} |^2   \\
    && +  C(q_0)   \cdot  c (\delta) \cdot c (\gamma)  \cdot E ( 4)  \cdot \frac{\eps  }{(1+\tau+|q|)^{1-   c (\gamma)  \cdot c (\delta)  \cdot E (  4)\cdot  \eps } \cdot (1+|q|) }   \cdot | \rderm \Phi_{V} |^2 \Big)  \cdot d\tau \cdot w(q) d^{3}x \; . 
 \eeaa 
Based on Lemma \ref{derivativeoftildwandrelationtotildew}, we have,
\bea
\widehat{w}^{\prime}(q) \sim \frac{\widehat{w}(q)}{(1+|q|)} \; .
\eea
and by injecting in the previous estimate, we get 
  \beaa
   \notag
 &&     \int_{\Sigma^{ext}_{t_2} }  |\derm \Phi_{V} |^2     \cdot w(q)  \cdot d^{3}x    + \int_{N_{t_1}^{t_2} }   T_{\hat{L} t}^{(\bf{g})}  \cdot  w(q) \cdot dv^{(\bf{m})}_N \\
 \notag
 &&+ \int_{t_1}^{t_2}  \int_{\Sigma^{ext}_{\tau} }  \Big(    \frac{1}{2} \Big(  | \derm_t  \Phi_{V} + \derm_r \Phi_{V} |^2  +  \de^{ij}  | ( \derm_i - \frac{x_i}{r} \derm_{r}  )\Phi_{V} |^2 \Big)  \cdot d\tau \cdot  \frac{\widehat{w} (q)}{(1+|q|)} d^{3}x \\
  \notag
  &\les &       \int_{\Sigma^{ext}_{t_1} }  |\derm \Phi_{V} |^2     \cdot w(q)  \cdot d^{3}x \\
    \notag
   &&  +  \int_{t_1}^{t_2}  \int_{\Sigma^{ext}_{\tau} }  |  g^{\mu\a} \derm_{\mu } \derm_\a \Phi_{V} | \cdot |  \derm_t \Phi_{V} |  \cdot d\tau \cdot w(q) d^{3}x  \\
 \notag
&& + \int_{t_1}^{t_2}  \int_{\Sigma^{ext}_{\tau} }  \Big(   C(q_0) \cdot   c (\delta) \cdot c (\gamma) \cdot E (  4)  \cdot \frac{ \eps }{ (1+\tau+|q|)\cdot (1+|q|)  } \cdot | \derm \Phi_{V} |^2   \\
    && +  C(q_0)   \cdot  c (\delta) \cdot c (\gamma)  \cdot E ( 4)  \cdot \frac{\eps  }{(1+\tau+|q|)^{1-   c (\gamma)  \cdot c (\delta)  \cdot E (  4)\cdot  \eps } \cdot (1+|q|) }   \cdot | \rderm \Phi_{V} |^2 \Big)  \cdot d\tau \cdot ( w(q) + \widehat{w}(q)  ) d^{3}x \; . 
 \eeaa 
 Since for $\mu < 0$\,, we have $\widehat{w}(q)  \leq w(q)$\,, we conclude the proof.
 
\end{proof}

\begin{lemma}\label{Theveryfinal }
For  $\eps$ small enough, depending on $q_0$\,, on $\de$\,, on $\ga$\;, and on $\mu < 0$\;, the following energy estimate holds for $\ga > 0$ and for $\Phi$ decaying sufficiently fast at spatial infinity,
  \beaa
   \notag
 &&     \int_{\Sigma^{ext}_{t_2} }  |\derm \Phi_{V} |^2     \cdot w(q)  \cdot d^{3}x    + \int_{N_{t_1}^{t_2} }  T_{\hat{L} t}^{(\bf{g})} (\Phi_{V})  \cdot  w(q) \cdot dv^{(\bf{m})}_N \\
 \notag
 &&+ \int_{t_1}^{t_2}  \int_{\Sigma^{ext}_{\tau} }     | \rderm \Phi_{V} |^2   \cdot  \frac{\widehat{w} (q)}{(1+|q|)} \cdot  d^{3}x  \cdot d\tau \\
  \notag
  &\les &       \int_{\Sigma^{ext}_{t_1} }  |\derm \Phi_{V} |^2     \cdot w(q)  \cdot d^{3}x \\
    \notag
   &&  +  \int_{t_1}^{t_2}  \int_{\Sigma^{ext}_{\tau} }  \frac{(1+\tau )}{\eps} \cdot |  g^{\mu\a} \derm_{\mu } \derm_\a \Phi_{V}|^2    \cdot w(q) \cdot  d^{3}x  \cdot d\tau \\
 \notag
&& + \int_{t_1}^{t_2}  \int_{\Sigma^{ext}_{\tau} }   C(q_0) \cdot   c (\delta) \cdot c (\gamma) \cdot E (  4)  \cdot \frac{ \eps }{ (1+\tau)  } \cdot | \derm \Phi_{V} |^2     \cdot  w(q)\cdot  d^{3}x  \cdot d\tau \; . 
 \eeaa 
 
 \end{lemma}

\begin{proof}
We examine the terms in Corollary \ref{theenergyestimatewithcontrolontangentialderivativesusinghatwoveroneplusq}.

\textbf{The term $| \derm_t  \Phi_{V } + \derm_r \Phi_{V} |^2  +  \de^{ij}  | ( \derm_i - \frac{x_i}{r} \derm_{r}  )\Phi_{V} |^2$:}\\

We have
\beaa
 | \derm_t  \Phi_{V} + \derm_r \Phi_{V} |^2  &=& |\derm_{L} \Phi_{V}  |^2 \; .
\eeaa

Using Definition \ref{defrestrictedderivativesintermsofZ}, we compute
\beaa
 && \de^{ij}  < ( \derm_i - \frac{x_i}{r} \derm_{r}  )\Phi_{V} , (\derm_j - \frac{x_j}{r} \derm_{r} ) \Phi_{V}> \\
  &=&  \de^{ij}  < \derm_{\rpa_{i}} \Phi_{V} , \derm_{\rpa_{i}}  \Phi_{V} >  := \sum_{i=1}^3  | \derm_{\rpa_{i}} \Phi_{V} |^2 \; .
\eeaa

We showed, in Lemma \ref{restrictedderivativesintermsofZ}, that
 \beaa
\rpa_i  &=&   \frac{x^j}{r^2}   Z_{ij} \, ,
 \eeaa
 and that we also have
 \bea
e_A = \frac{1}{r}C^{ij}_A Z_{ij} \, ,
\eea
where $C^{ij}_A$ are bounded function on $\SSS^2$.
Thus,
\beaa
&& | \derm_t  \Phi_{V } + \derm_r \Phi_{V } |^2  +  \de^{ij}  | ( \derm_i - \frac{x_i}{r} \derm_{r}  )\Phi_{V } |^2 \\
 &=&   |\derm_{L}  \Phi_{V } |^2 +  \sum_{i=1}^{3}   | \frac{x^j}{r^2}  \derm_{ Z_{ij} } \Phi_{V } |^2 \\
  &=&   |\derm_{L}  \Phi_{V}  |^2 +  \sum_{i=1}^{3}   | \frac{1}{r} \cdot \frac{x^j}{r}   \derm_{ Z_{ij} } \Phi_{V}  |^2 \\
  &=&  | \derm_{L}  \Phi_{V} |^2 +   \sum_{A=1}^{2}  | \frac{1}{r} \cdot C^{ij}_A \derm_{Z_{ij}}  \Phi_{V}  |^2 \\
    &=&    |\derm_{L}  \Phi_{V} |^2 +  \sum_{A=1}^{2}   | \derm_{e_A} \Phi_{V}  |^2 \\
    &\sim&    |\rderm  \Phi_{V} |^2 \; .
  \eeaa

\textbf{Control on $ \int_{t_1}^{t_2}  \int_{\Sigma^{ext}_{\tau} }     |\rderm \Phi_{V} |^2  \cdot  \frac{\widehat{w} (q)}{(1+|q|)} \cdot d^{n}x \cdot d\tau $:} \\

 Now, we want to control the space-time integral on the right hand side, that is $$\int_{t_1}^{t_2}  \int_{\Sigma^{ext}_{\tau} }    C(q_0)   \cdot  c (\delta) \cdot c (\gamma)  \cdot E ( 4) ) \cdot \frac{\eps  }{(1+t+|q|)^{1-   c (\gamma)  \cdot c (\delta)  \cdot E (  4)\cdot  \eps } \cdot (1+|q|) }   \cdot | \rderm \Phi_{V} |^2   \cdot  w(q) \cdot d^{n}x\ \cdot d\tau \; , $$ by absorbing it into the space-time integral on the left hand side, that is $$ \int_{t_1}^{t_2}  \int_{\Sigma^{ext}_{\tau} }     |\rderm \Phi_{V} |^2   \cdot  \frac{\widehat{w} (q)}{(1+|q|)} \cdot d^{n}x \cdot d\tau \;.$$ This would end up giving us an energy estimate with control on the latter.

 The term in the integrand on the left hand side is
 $$
\frac{\widehat{w} (q)}{(1+|q|)}  \cdot    |\rderm \Phi_{V} |^2 \; ,
$$
 
and the term in the integrand on the right hand side is
 $$
  \frac{\eps  }{(1+t+|q|)^{1-   c (\gamma)  \cdot c (\delta)  \cdot E (  4)\cdot  \eps } \cdot (1+|q|) }   \cdot | \rderm \Phi_{V} |^2   \cdot  w(q) 
  $$

 For $q > 0$, we know that $w (q) = \widehat{w} (q) $ and therefore,
  \beaa
  && C(q_0)   \cdot  c (\delta) \cdot c (\gamma)  \cdot E ( 4)  \cdot  \frac{\eps  }{(1+t+|q|)^{1-   c (\gamma)  \cdot c (\delta)  \cdot E (  4)\cdot  \eps } \cdot (1+|q|) }   \cdot | \rderm \Phi_{V} |^2   \cdot  w(q) \\
  &=& C(q_0)   \cdot  c (\delta) \cdot c (\gamma)  \cdot E ( 4) )\cdot  \frac{\eps  }{(1+t+|q|)^{1-   c (\gamma)  \cdot c (\delta)  \cdot E (  4)\cdot  \eps }  }   \cdot   \frac{\widehat{w} (q)}{(1+|q|)}  \cdot    |\rderm \Phi_{V} |^2 \\
   &\leq& C(q_0)   \cdot  c (\delta) \cdot c (\gamma)  \cdot E ( 4)  \cdot  \eps   \cdot   \frac{\widehat{w} (q)}{(1+|q|)}  \cdot    |\rderm \Phi_{V} |^2 \\
    &\leq& C(q_0)   \cdot  c (\delta) \cdot c (\gamma)    \cdot  \eps   \cdot   \frac{\widehat{w} (q)}{(1+|q|)}  \cdot    |\rderm \Phi_{V} |^2 \\
    && \text{(since $E(4) \leq 1$).}
  \eeaa
Thus, for $\eps$ small enough depending on $q_0$\,, $\de$\,, and $\ga$\,, the term on the right hand side could be absorbed into the space-time integral on the left hand side.
  
  For $q < 0$, we know that $\frac{\widehat{w} (q)}{(1+|q|)} = \frac{1}{(1+|q|)^{1-2\mu}}$ and  $\frac{w (q)}{(1+|q|)} = \frac{1}{(1+|q|)}$. Thus,
  \beaa
  && | C(q_0)   \cdot  c (\delta) \cdot c (\gamma)  \cdot E ( 4)  \cdot  \frac{\eps  }{(1+t+|q|)^{1-   c (\gamma)  \cdot c (\delta)  \cdot E (  4)\cdot  \eps } \cdot (1+|q|) }   \cdot | \rderm \Phi_{V} |^2   \cdot  w(q) | \\
  &\leq& C(q_0)   \cdot  c (\delta) \cdot c (\gamma)  \cdot E ( 4 )\cdot  \frac{\eps  }{(1+|q|)^{2-   c (\gamma)  \cdot c (\delta)  \cdot E (  4)\cdot  \eps }  }   \cdot    |\rderm \Phi_{V} |^2 \; .
  \eeaa  

Whereas 
\beaa
\frac{\widehat{w} (q)}{(1+|q|)}  \cdot    |\rderm \Phi_{V} |^2 &=& \frac{1 }{(1+|q|)^{1-  2\mu }  } \cdot    |\rderm \Phi_{V} |^2 \; .
\eeaa
We have for $q_0 \leq q \leq 0$, that for $\eps$ small enough depending on $q_0$, on $\de$\,, on $\ga$\,, on $E ( 4) $ and on $\mu$, that
\beaa
C(q_0)   \cdot  c (\delta) \cdot c (\gamma)  \cdot E ( 4 )\cdot  \frac{\eps  }{(1+|q|)^{2-   c (\gamma)  \cdot c (\delta)  \cdot E (  4)\cdot  \eps }  }  \leq  \frac{1 }{(1+|q|)^{1-  2\mu }  } \; .
\eeaa

Consequently, for all $q$, if $\eps$ small enough, depending on $q_0$, $\de$\,, $\ga$, $E ( 4) $ and on $\mu$\;, we would then be able to absorb the space-time integral in study into the left hand side of the inequality in Corollary \ref{theenergyestimatewithcontrolontangentialderivativesusinghatwoveroneplusq}, and therefore, we obtain the following energy estimate
  \bea\label{energyestimatewithoutdecomposingtheproductucttimederivativeoffieldtimessourceterms}
   \notag
 &&     \int_{\Sigma^{ext}_{t_2} }  |\derm \Phi_{V} |^2     \cdot w(q)  \cdot d^{3}x    + \int_{N_{t_1}^{t_2} }  \big( T_{\hat{L} t}^{(\bf{g})} \big)  \cdot  w(q) \cdot dv^{(\bf{m})}_N \\
 \notag
 &&+ \int_{t_1}^{t_2}  \int_{\Sigma^{ext}_{\tau} }     | \rderm \Phi_{V} |^2   \cdot  \frac{\widehat{w} (q)}{(1+|q|)} d^{3}x  \cdot d\tau \\
  \notag
  &\les &       \int_{\Sigma^{ext}_{t_1} }  |\derm \Phi |^2     \cdot w(q)  \cdot d^{3}x \\
    \notag
   &&  +  \int_{t_1}^{t_2}  \int_{\Sigma^{ext}_{\tau} }  |  g^{\mu\a} \derm_{\mu } \derm_\a \Phi_{V} | \cdot |  \derm_t \Phi_{V} |   \cdot w(q) d^{3}x  \cdot d\tau \\
 \notag
&& + \int_{t_1}^{t_2}  \int_{\Sigma^{ext}_{\tau} }  \Big(   C(q_0) \cdot   c (\delta) \cdot c (\gamma) \cdot E (  4)  \cdot \frac{ \eps }{ (1+\tau+|q|)\cdot (1+|q|)  } \cdot | \derm \Phi_{V} |^2   \Big)   \cdot  w(q) d^{3}x  \cdot d\tau \; .  \\
 \eea 
 
\begin{remark}\label{removaloftheassumptiononthesmallnessonHbyinclduingitinthesmalnnessofepsilon}
Note that for $\eps$ small enough, the assumption $| H| < \frac{1}{3}$ of Corollary \ref{theenergyestimatewithcontrolontangentialderivativesusinghatwoveroneplusq} is granted and therefore, it is no more an assumption for this lemma.
\end{remark}

\textbf{ The term $\int_{t_1}^{t_2}  \int_{\Sigma^{ext}_{\tau} }  |  g^{\mu\a} \derm_{\mu } \derm_\a \Phi_{V} | \cdot |  \derm_t \Phi_{V} |   \cdot w(q) d^{n}x  \cdot d\tau$:}\\

Using the inequality $a\cdot b \les a^2 + b^2$, we get
\beaa
&& \int_{t_1}^{t_2}  \int_{\Sigma^{ext}_{\tau} }  |  g^{\mu\a} \derm_{\mu } \derm_\a \Phi_{V} | \cdot |  \derm_t \Phi_{V} |   \cdot w(q) d^{3}x  \cdot d\tau \\
&=&\int_{t_1}^{t_2}  \int_{\Sigma^{ext}_{\tau} }  \frac{\sqrt{(1+\tau )}}{\sqrt{\eps}} \sqrt{w} \cdot | g^{\la\a} \derm_{\la}   \derm_{\a} \Phi_{V}  |\cdot  \frac{\sqrt{\eps}}{\sqrt{(1+\tau )}} \cdot |\derm \Phi_{V} | \, \sqrt{w}  \cdot  d^{3}x \cdot d\tau \\
&\les&\int_{t_1}^{t_2}  \int_{\Sigma^{ext}_{\tau} }  \frac{ \epsilon  \cdot |\derm  \Phi_{V} |^{2}}{1+\tau} \cdot  d^{3}x\cdot  d\tau \\
&&+ \int_{t_1}^{t_2}  \int_{\Sigma^{ext}_{\tau} }   \frac{(1+\tau )}{\eps} \cdot | g^{\la\a} \derm_{\la}   \derm_{\a} \Phi_{V} |^2 \, w  \cdot  d^{3}x \cdot d\tau \; .
\eeaa

Injecting in the energy estimate \eqref{energyestimatewithoutdecomposingtheproductucttimederivativeoffieldtimessourceterms}, we conclude the proof of the announced energy estimate.

\end{proof}

\section{The commutator term revisited}

\begin{lemma}\label{estimateonthecommutatortermusingtheproductsandusingbootstrap}

Let $\Phi = A $ or $\Phi = h^1$, and $\ga \geq 3 \de $\,, $0 < \de \leq \frac{1}{4}$\,, and $M \leq \eps$\,. We have the following estimate on the commutator term,

 \beaa
\notag
&&  \frac{(1+ t )}{\eps} \cdot  | \Lie_{Z^I}  ( g^{\la\mu} \derm_{\la}   \derm_{\mu}     \Phi_{V} ) - g^{\la\mu}    \derm_{\la}   \derm_{\mu}  (  \Lie_{Z^I} \Phi_{V}  ) |^2  \\
  \notag
   &\les& \frac{(1+ t )}{\eps} \cdot    \sum_{|K| \leq |I|  - 1}  | g^{\la\mu} \cdot \derm_{\la}   \derm_{\mu} (  \Lie_{Z^{K}}  \Phi_{V} ) |^2 \\
&& +  C(q_0)   \cdot  c (\delta) \cdot c (\gamma) \cdot C(  \lfloor \frac{|I|}{2} \rfloor ) \cdot E (    \lfloor \frac{|I|}{2} \rfloor   +4)  \\
&& \times \Big(  \sum_{ |K| \leq |I| }   \frac{\eps   }{(1+t+|q|)^{3-   c (\gamma)  \cdot c (\delta)  \cdot c(  \lfloor \frac{|I|}{2} \rfloor ) \cdot E (  \lfloor \frac{|I|}{2} \rfloor + 4)\cdot  \eps } \cdot (1+|q|)^2}  \, \cdot  |    \Lie_{Z^K}  H   |^2  \\
&& + \sum_{ |K| \leq |I| }  \frac{\eps }{(1+t+|q|) \cdot (1+|q|)^{2-   c (\gamma)  \cdot c (\delta)  \cdot c( \lfloor \frac{|I|}{2} \rfloor ) \cdot E (  \lfloor \frac{|I|}{2} \rfloor + 4)\cdot  \eps } }  \, \cdot | \derm ( \Lie_{Z^K}  \Phi )  |^2 \Big)   \\
&&+  C(q_0)   \cdot  c (\delta) \cdot c (\gamma) \cdot C(|I|) \cdot E (    \lfloor \frac{|I|}{2} \rfloor   +3 )  \cdot \frac{\eps }{(1+t+|q|)^{1-   c (\gamma)  \cdot c (\delta)  \cdot c(|I|) \cdot E (\lfloor \frac{|I|}{2} \rfloor  + 2)\cdot  \eps } \cdot (1+|q|)^2} \\
&& \times    \sum_{  |K| \leq |I| -1 }  | \derm ( \Lie_{Z^K}  \Phi )  |^2  \\
  && + \sum_{  |K| \leq |I| }      C(q_0)   \cdot C(|I|) \cdot E (  \lfloor \frac{|I|}{2} \rfloor  +3)  \cdot \frac{\eps^2 \cdot  | \Lie_{Z^{K}} H_{L  L} |^2 }{(1+t+|q|)^{2- 2 \de } \cdot (1+|q|)^{4+2\gamma}} \;.
\eeaa

\end{lemma}

\begin{proof}

From Lemma \ref{Thecommutationformulawithpossibilityofsperationoftangentialcomponentsaswell}, we have the following estimate on the commutator term: for any $V \in \cal T$,
 \beaa
\notag
&&| \Lie_{Z^I}  ( g^{\la\mu} \derm_{\la}   \derm_{\mu}     \Phi_{V} ) - g^{\la\mu}    \derm_{\la}   \derm_{\mu}  (  \Lie_{Z^I} \Phi_{V}  ) |  \\
  \notag
   &\les&  \sum_{|K| < |I| }  | g^{\la\mu} \cdot \derm_{\la}   \derm_{\mu} (  \Lie_{Z^{K}}  \Phi_{V} ) | \\
   \notag
&&+  \frac{1}{(1+t+|q|)}  \cdot \sum_{|K|\leq |I|,}\,\, \sum_{|J|+(|K|-1)_+\le |I|} \,\,\, | \Lie_{Z^{J}} H |\, \cdot | \derm ( \Lie_{Z^K}  \Phi )  | \\
   \notag
&& +   \frac{1}{(1+|q|)}  \cdot \sum_{|K|\leq |I|,}\,\, \sum_{|J|+(|K|-1)_+\le |I|} \,\,\, | \Lie_{Z^{J}} H_{L  L} |\, \cdot  | \derm ( \Lie_{Z^K}  \Phi  )  |  \; ,
\eeaa
  where $(|K|-1)_+=|K|-1$ if $|K|\geq 1$ and $(|K|-1)_+=0$ if $|K|=0$. We now study each term on the right hand side of the estimate.\\
  
\textbf{The term $\sum_{|K|\leq |I|,}\,\, \sum_{|J|+(|K|-1)_+\le |I|} \,\,\, | \Lie_{Z^{J}} H |\, \cdot | \derm ( \Lie_{Z^K}  \Phi )  |$:}\\

We estimate the term with the good factor $ \frac{1}{(1+t+|q|)}$, as follows,
\beaa
&&  \sum_{|K|\leq |I|,}\,\, \sum_{|J|+(|K|-1)_+\le |I|} \,\,\, | \Lie_{Z^{J}} H |\, \cdot | \derm ( \Lie_{Z^K}  \Phi )  | \\
&\leq& \sum_{ |K| \leq   \lfloor \frac{|I|}{2} \rfloor   ,\; |J| \leq |I| }   | \Lie_{Z^{J}} H |\, \cdot | \derm ( \Lie_{Z^K}  \Phi )  | + \sum_{ |J| \leq   \lfloor \frac{|I|}{2} \rfloor   ,\; |K| \leq |I| }  | \Lie_{Z^{J}} H |\, \cdot | \derm ( \Lie_{Z^K}  \Phi )  | \;.
\eeaa 
Let either $\Phi = A $ or $\Phi = h^1$, and $\ga \geq 3 \de $\,, $0 < \de \leq \frac{1}{4}$\,, and $M \leq \eps$\,. For all $|K| \leq   \lfloor \frac{|I|}{2} \rfloor $, based on Lemma \ref{upgradedestimateonLiederivativesoffields}, we have
            \bea
 \notag
&& |\derm  ( \Lie_{Z^K} \Phi)   |    \\
\notag
 &\leq&   C(q_0)   \cdot  c (\delta) \cdot c (\gamma) \cdot C(  \lfloor \frac{|I|}{2} \rfloor ) \cdot E (    \lfloor \frac{|I|}{2} \rfloor   +4)  \cdot \frac{\eps   }{(1+t+|q|)^{1-   c (\gamma)  \cdot c (\delta)  \cdot c(  \lfloor \frac{|I|}{2} \rfloor ) \cdot E (  \lfloor \frac{|I|}{2} \rfloor + 4)\cdot  \eps } \cdot (1+|q|)}     \; , \\
      \eea
      and based on Lemma  \ref{upgradedestimateonBIGH}, we have
                 \beaa
 \notag
 |    \Lie_{Z^K}  H   | &\leq&   C(q_0)   \cdot  c (\delta) \cdot c (\gamma) \cdot C( \lfloor \frac{|I|}{2} \rfloor )  \cdot E (   \lfloor \frac{|I|}{2} \rfloor + 4)  \cdot \frac{\eps }{(1+t+|q|)^{1-   c (\gamma)  \cdot c (\delta)  \cdot c( \lfloor \frac{|I|}{2} \rfloor ) \cdot E (  \lfloor \frac{|I|}{2} \rfloor + 4)\cdot  \eps } }     \; . \\
      \eeaa
Injecting, we obtain
\beaa
&&  \sum_{|K|\leq |I|,}\,\, \sum_{|J|+(|K|-1)_+\le |I|} \,\,\, | \Lie_{Z^{J}} H |\, \cdot | \derm ( \Lie_{Z^K}  \Phi )  | \\
&\leq& C(q_0)   \cdot  c (\delta) \cdot c (\gamma) \cdot C(  \lfloor \frac{|I|}{2} \rfloor ) \cdot E (    \lfloor \frac{|I|}{2} \rfloor   +4)  \\
&& \times \Big(  \sum_{ |K| \leq |I|}   \frac{\eps   }{(1+t+|q|)^{1-   c (\gamma)  \cdot c (\delta)  \cdot c(  \lfloor \frac{|I|}{2} \rfloor ) \cdot E (  \lfloor \frac{|I|}{2} \rfloor + 4)\cdot  \eps } \cdot (1+|q|)}  \, \cdot  |    \Lie_{Z^K}  H   |  \\
&& + \sum_{  |K| \leq |I| }  \frac{\eps }{(1+t+|q|)^{1-   c (\gamma)  \cdot c (\delta)  \cdot c( \lfloor \frac{|I|}{2} \rfloor ) \cdot E (  \lfloor \frac{|I|}{2} \rfloor + 4)\cdot  \eps } }  \, \cdot | \derm ( \Lie_{Z^K}  \Phi )  |  \Big)  \; .
\eeaa 

Thus,
\beaa
&&   \frac{1}{(1+t+|q|)} \cdot \sum_{|K|\leq |I|,}\,\, \sum_{|J|+(|K|-1)_+\le |I|} \,\,\, | \Lie_{Z^{J}} H |\, \cdot | \derm ( \Lie_{Z^K}  \Phi )  | \\
&\leq& C(q_0)   \cdot  c (\delta) \cdot c (\gamma) \cdot C(  \lfloor \frac{|I|}{2} \rfloor ) \cdot E (    \lfloor \frac{|I|}{2} \rfloor   +4)  \\
&& \times \Big(  \sum_{ |K| \leq |I| }   \frac{\eps   }{(1+t+|q|)^{2-   c (\gamma)  \cdot c (\delta)  \cdot c(  \lfloor \frac{|I|}{2} \rfloor ) \cdot E (  \lfloor \frac{|I|}{2} \rfloor + 4)\cdot  \eps } \cdot (1+|q|)}  \, \cdot  |    \Lie_{Z^K}  H   |  \\
&& + \sum_{ |K| \leq |I| }  \frac{\eps }{(1+t+|q|)^{2-   c (\gamma)  \cdot c (\delta)  \cdot c( \lfloor \frac{|I|}{2} \rfloor ) \cdot E (  \lfloor \frac{|I|}{2} \rfloor + 4)\cdot  \eps } }  \, \cdot | \derm ( \Lie_{Z^K}  \Phi )  | \Big)   \; .
\eeaa

\textbf{The term $\sum_{|K|\leq |I|,}\,\, \sum_{|J|+(|K|-1)_+\le |I|} \,\,\, | \Lie_{Z^{J}} H_{L  L} |\, \cdot  | \derm ( \Lie_{Z^K}  \Phi  )  |$:}\\

We now examine the term with the bad factor (non-decaying in time $t$), that is $\frac{1}{(1+|q|)}$.

We have based on Lemma \ref{upgradedestimateonLiederivativesoffields}, that 
\beaa
&& \sum_{|K|\leq |I|,}\,\, \sum_{|J|+(|K|-1)_+\le |I|} \,\,\, | \Lie_{Z^{J}} H_{L  L} |\, \cdot | \derm ( \Lie_{Z^K}  \Phi  )  | \\
 &=& \sum_{ |J| \leq  1   ,\; |K| = |I| } | \Lie_{Z^{J}} H_{L  L} |\, \cdot  | \derm ( \Lie_{Z^K}  \Phi  )  |\\
  && + \sum_{|K|\leq |I|-1,}\,\, \sum_{|J|+(|K|-1)_+\le |I|}   | \Lie_{Z^{J}} H_{L  L} |\, \cdot | \derm ( \Lie_{Z^K}  \Phi  )  | \\
     &=& \sum_{ |J| \leq  1   ,\; |K| = |I| } | \Lie_{Z^{J}} H_{L  L} |\, \cdot | \derm ( \Lie_{Z^K}  \Phi  )  | \\
     && + \sum_{ |K| \leq   \lfloor \frac{|I|}{2} \rfloor  +1  ,\; |J| \leq |I| }      | \Lie_{Z^{J}} H_{L  L} |\, \cdot  | \derm ( \Lie_{Z^K}  \Phi  )  |\\
     && +\sum_{ |K| \leq  |I|-1  ,\; |J| \leq \lfloor \frac{|I|}{2} \rfloor }  | \Lie_{Z^{J}} H_{L  L} |\, \cdot | \derm ( \Lie_{Z^K}  \Phi  )  |\\
    && \text{(where the last sum vanishes if $|I|  < 1$).}
    \eeaa

  To deal with the term $ \sum_{ |K| \leq   \lfloor \frac{|I|}{2} \rfloor  +1  ,\; |J| \leq |I| }      | \Lie_{Z^{J}} H_{L  L} |\, \cdot | \derm ( \Lie_{Z^K}  \Phi  )  | $, for either $\Phi =A$ or $\Phi= h^1$, we point out that we have based on the a priori estimates of Lemma \ref{apriordecayestimatesfrombootstrapassumption}, that for all $|K| \leq   \lfloor \frac{|I|}{2} \rfloor  +1$\;,

                  \beaa
 \notag
 && |\derm  ( \Lie_{Z^K} \Phi )   | \\
  &\leq&    C(q_0) \cdot C ( |I| ) \cdot E (  \lfloor \frac{|I|}{2} \rfloor  + 3)  \cdot \frac{\eps }{(1+t+|q|)^{1-\delta} \cdot  (1+|q|)^{1+\gamma}}    \; .
      \eeaa
      Thus,
       \beaa
     &&  \sum_{ |K| \leq   \lfloor \frac{|I|}{2} \rfloor  +1  ,\; |J| \leq |I| }      | \Lie_{Z^{J}} H_{L  L} |\, \cdot  | \derm ( \Lie_{Z^K}  \Phi  )  |  \\
       &\les&    \sum_{  |K| \leq |I| }      C(q_0)   \cdot C(|I|) \cdot E (  \lfloor \frac{|I|}{2} \rfloor  +3)  \cdot \frac{\eps \cdot  | \Lie_{Z^{K}} H_{L  L} | }{(1+t+|q|)^{1-   \de } \cdot (1+|q|)^{1+\gamma}}  \; .
       \eeaa

      To deal with the term $\sum_{ |K| \leq  |I|-1  ,\; |J| \leq \lfloor \frac{|I|}{2} \rfloor }  | \Lie_{Z^{J}} H_{L  L} |\, \cdot | \derm ( \Lie_{Z^K}  \Phi  )  |$\;, we point out that based on Lemma \ref{estimategoodcomponentspotentialandmetric}, we have for all $ |J| \leq   \lfloor \frac{|I|}{2} \rfloor $\;, 

         \beaa
        \notag
|  \Lie_{Z^J} H_{L L} |                    &\les&   \int\limits_{s,\,\Om=const} \sum_{ |M| \leq   \lfloor \frac{|I|}{2} \rfloor   -2}  |\derm  ( \Lie_{Z^M} H )  | \\
\notag
&& + C(q_0) \cdot   c (\delta) \cdot c (\gamma) \cdot C (   |I|  ) \cdot E (   \lfloor \frac{|I|}{2} \rfloor  + 3)  \cdot \frac{ \eps   }{ (1+t+|q|) } \; \\
&& \text{(where the sum vanishes if $ \lfloor \frac{|I|}{2} \rfloor   < 2$),}
       \eeaa
       and for  $ |J| \leq  1 $, 
                \beaa
        \notag
|  \Lie_{Z^J} H_{L L} |                    &\les&   C(q_0) \cdot   c (\delta) \cdot c (\gamma)  \cdot E (  4)  \cdot \frac{ \eps   }{ (1+t+|q|) } \; . \\
       \eeaa
       
          Hence,
          
          \beaa
&& \sum_{|K|\leq |I|,}\,\, \sum_{|J|+(|K|-1)_+\le |I|} \,\,\, | \Lie_{Z^{J}} H_{L  L} |\, \cdot | \derm ( \Lie_{Z^K}  \Phi  )  | \\
     &\les& \sum_{  |K| = |I| }  C(q_0) \cdot   c (\delta) \cdot c (\gamma)  \cdot E (  4)  \cdot \frac{ \eps   }{ (1+t+|q|) }  \, \cdot | \derm ( \Lie_{Z^K}  \Phi  )  | \\
&&+ C(q_0) \cdot   c (\delta) \cdot c (\gamma) \cdot C (    |I|   ) \cdot E (   \lfloor \frac{|I|}{2} \rfloor   + 3)  \cdot   \frac{ \eps   }{ (1+t+|q|)  } \, \cdot \sum_{  |K| \leq |I| -1}  | \derm ( \Lie_{Z^K}  \Phi  )  | \\
&&+  \Big(  \int\limits_{s,\,\Om=const} \sum_{ |J| \leq   \lfloor \frac{|I|}{2} \rfloor   -2}  |\derm  ( \Lie_{Z^J} H )  |  \Big) \cdot    \sum_{  |K| \leq |I| -1 } | \derm ( \Lie_{Z^K}  \Phi  )  | \\
 && + \sum_{  |K| \leq |I| }      C(q_0)   \cdot C(|I|) \cdot E (  \lfloor \frac{|I|}{2} \rfloor  +3)  \cdot \frac{\eps \cdot  | \Lie_{Z^{K}} H_{L  L} | }{(1+t+|q|)^{1-  \de } \cdot (1+|q|)^{1+\gamma}} \;.
\eeaa

Therefore,

          \beaa
&& \frac{1}{(1+|q|)} \cdot   \sum_{|K|\leq |I|,}\,\, \sum_{|J|+(|K|-1)_+\le |I|} \,\,\, | \Lie_{Z^{J}} H_{L  L} |\, \cdot  | \derm ( \Lie_{Z^K}  \Phi  )  |\\
     &\les& \sum_{  |K| \leq |I| }  C(q_0) \cdot   c (\delta) \cdot c (\gamma)  \cdot E (  4)  \cdot \frac{ \eps   }{ (1+t+|q|)  \cdot (1+|q|) }  \, \cdot  | \derm ( \Lie_{Z^K}  \Phi  )  | \\
&&+ C(q_0) \cdot   c (\delta) \cdot c (\gamma) \cdot C (    |I|   ) \cdot E (    \lfloor \frac{|I|}{2} \rfloor    + 3)  \cdot   \frac{ \eps   }{ (1+t+|q|) \cdot (1+|q|)   } \, \cdot \sum_{  |K| \leq |I| -1}  | \derm ( \Lie_{Z^K}  \Phi  )  | \\
&&+  \frac{1}{(1+|q|)} \cdot  \Big(  \int\limits_{s,\,\Om=const} \sum_{ |J| \leq    \lfloor \frac{|I|}{2} \rfloor  -2}  |\derm  ( \Lie_{Z^J} H )  |  \Big) \cdot    \sum_{  |K| \leq |I| -1 } | \derm ( \Lie_{Z^K}  \Phi  )  |\\
 && + \sum_{  |K| \leq |I| }      C(q_0)   \cdot C(|I|) \cdot E (  \lfloor \frac{|I|}{2} \rfloor  +3)  \cdot \frac{\eps \cdot  | \Lie_{Z^{K}} H_{L  L} | }{(1+t+|q|)^{1-  \de } \cdot (1+|q|)^{2+\gamma}} \;.
\eeaa

Using Lemma \ref{upgradedestimateonBIGH}, we have for $\ga \geq 3 \de $\,, and $0 < \de \leq \frac{1}{4}$\,, that in the exterior region, the following estimates hold for $M \leq \eps$\,, and for all $|J| \leq    \lfloor \frac{|I|}{2} \rfloor   -2$\,, 
                    \beaa
 \notag
 && |\derm  ( \Lie_{Z^J} H )   | \\
  &\leq&   C(q_0)   \cdot  c (\delta) \cdot c (\gamma) \cdot C(|I|) \cdot E (  \lfloor \frac{|I|}{2} \rfloor   +2)  \cdot \frac{\eps }{(1+t+|q|)^{1-   c (\gamma)  \cdot c (\delta)  \cdot c(|I|) \cdot E (  \lfloor \frac{|I|}{2} \rfloor + 2)\cdot  \eps } \cdot (1+|q|)}     \; . 
      \eeaa
Therefore, using the asymptotic behaviour of the initial data that we showed in \cite{G4} (see Remark \ref{remarkabouttheasymptoticbehaviouroftheboundarytermonhyperplaneprescribedtqual0}), we obtain
                    \beaa
 \notag
 && \int\limits_{s,\,\Om=const} \sum_{ |J| \leq    |I|  -2} |\derm  ( \Lie_{Z^J} H )   | \\
  &\leq&   C(q_0)   \cdot  c (\delta) \cdot c (\gamma) \cdot C(|I|) \cdot E (  \lfloor \frac{|I|}{2} \rfloor  +2)  \cdot \frac{\eps }{(1+t+|q|)^{1-   c (\gamma)  \cdot c (\delta)  \cdot c(|I|) \cdot E (  \lfloor \frac{|I|}{2} \rfloor + 2)\cdot  \eps } }     \; . 
      \eeaa
      
Thus, we get
          \beaa
&& \frac{1}{(1+|q|)} \cdot   \sum_{|K|\leq |I|,}\,\, \sum_{|J|+(|K|-1)_+\le |I|} \,\,\, | \Lie_{Z^{J}} H_{L  L} |\, \cdot  | \derm ( \Lie_{Z^K}  \Phi  )  | \\
     &\les& \sum_{  |K| \leq |I| }  C(q_0) \cdot   c (\delta) \cdot c (\gamma)  \cdot E (  4)  \cdot \frac{ \eps   }{ (1+t+|q|)  \cdot (1+|q|) }  \, \cdot  | \derm ( \Lie_{Z^K}  \Phi  )  |\\
&&+ C(q_0) \cdot   c (\delta) \cdot c (\gamma) \cdot C (    |I|   ) \cdot E (  \lfloor \frac{|I|}{2} \rfloor    + 3)  \cdot   \frac{ \eps   }{ (1+t+|q|) \cdot (1+|q|)   } \, \cdot \sum_{  |K| \leq |I| -1}  | \derm ( \Lie_{Z^K}  \Phi  )  | \\
&&+  C(q_0)   \cdot  c (\delta) \cdot c (\gamma) \cdot C(|I|) \cdot E (  \lfloor \frac{|I|}{2} \rfloor   +2)  \cdot \frac{\eps }{(1+t+|q|)^{1-   c (\gamma)  \cdot c (\delta)  \cdot c(|I|) \cdot E ( \lfloor \frac{|I|}{2} \rfloor  + 2)\cdot  \eps } \cdot (1+|q|)} \\
&& \times    \sum_{  |K| \leq |I| -1 } | \derm ( \Lie_{Z^K}  \Phi  )  | \\
 && + \sum_{  |K| \leq |I| }      C(q_0)   \cdot C(|I|) \cdot E (  \lfloor \frac{|I|}{2} \rfloor  +3)  \cdot \frac{\eps \cdot  | \Lie_{Z^{K}} H_{L  L} | }{(1+t+|q|)^{1-  \de } \cdot (1+|q|)^{2+\gamma}} \;.
\eeaa

Finally, we obtain
          \beaa
&& \frac{1}{(1+|q|)} \cdot   \sum_{|K|\leq |I|,}\,\, \sum_{|J|+(|K|-1)_+\le |I|} \,\,\, | \Lie_{Z^{J}} H_{L  L} |\, \cdot  | \derm ( \Lie_{Z^K}  \Phi  )  | \\
     &\les& \sum_{  |K| \leq |I| }  C(q_0) \cdot   c (\delta) \cdot c (\gamma)  \cdot E (  4)  \cdot \frac{ \eps   }{ (1+t+|q|)  \cdot (1+|q|) }  \, \cdot   | \derm ( \Lie_{Z^K}  \Phi  )  | \\
&&+  C(q_0)   \cdot  c (\delta) \cdot c (\gamma) \cdot C(|I|) \cdot E (  \lfloor \frac{|I|}{2} \rfloor   +3)  \cdot \frac{\eps }{(1+t+|q|)^{1-   c (\gamma)  \cdot c (\delta)  \cdot c(|I|) \cdot E (\lfloor \frac{|I|}{2} \rfloor  + 2)\cdot  \eps } \cdot (1+|q|)} \\
&& \times    \sum_{  |K| \leq |I| -1 } | \derm ( \Lie_{Z^K}  \Phi  )  | \\
 && + \sum_{  |K| \leq |I| }      C(q_0)   \cdot C(|I|) \cdot E (  \lfloor \frac{|I|}{2} \rfloor  +3)  \cdot \frac{\eps \cdot  | \Lie_{Z^{K}} H_{L  L} | }{(1+t+|q|)^{1-  \de } \cdot (1+|q|)^{2+\gamma}} \;.
\eeaa

\textbf{The result:}\\

Putting all together, we obtain the following estimate for the commutator term

 \beaa
\notag
&&| \Lie_{Z^I}  ( g^{\la\mu} \derm_{\la}   \derm_{\mu}     \Phi_{V} ) - g^{\la\mu}    \derm_{\la}   \derm_{\mu}  (  \Lie_{Z^I} \Phi_{V}  ) |  \\
  \notag
   &\les&  \sum_{|K| < |I| }  | g^{\la\mu} \cdot \derm_{\la}   \derm_{\mu} (  \Lie_{Z^{K}}  \Phi_{V} ) | \\
&& +  C(q_0)   \cdot  c (\delta) \cdot c (\gamma) \cdot C(  \lfloor \frac{|I|}{2} \rfloor ) \cdot E (    \lfloor \frac{|I|}{2} \rfloor   +4)  \\
&& \times \Big(  \sum_{ |K| \leq |I|}   \frac{\eps   }{(1+t+|q|)^{2-   c (\gamma)  \cdot c (\delta)  \cdot c(  \lfloor \frac{|I|}{2} \rfloor ) \cdot E (  \lfloor \frac{|I|}{2} \rfloor + 4)\cdot  \eps } \cdot (1+|q|)}  \, \cdot  |    \Lie_{Z^K}  H   |  \\
&& + \sum_{ |K| \leq |I| }  \frac{\eps }{(1+t+|q|)^{2-   c (\gamma)  \cdot c (\delta)  \cdot c( \lfloor \frac{|I|}{2} \rfloor ) \cdot E (  \lfloor \frac{|I|}{2} \rfloor + 4)\cdot  \eps } }  \, \cdot | \derm ( \Lie_{Z^K}  \Phi )  | \Big)   \\
   && + \sum_{  |K| \leq |I| }  C(q_0) \cdot   c (\delta) \cdot c (\gamma)  \cdot E (  4)  \cdot \frac{ \eps   }{ (1+t+|q|)  \cdot (1+|q|) }  \, \cdot  | \derm ( \Lie_{Z^K}  \Phi  )  | \\
&&+  C(q_0)   \cdot  c (\delta) \cdot c (\gamma) \cdot C(|I|) \cdot E (   \lfloor \frac{|I|}{2} \rfloor   +3)  \cdot \frac{\eps }{(1+t+|q|)^{1-   c (\gamma)  \cdot c (\delta)  \cdot c(|I|) \cdot E ( |I|+ 2)\cdot  \eps } \cdot (1+|q|)} \\
&& \times    \sum_{  |K| \leq |I| -1 } | \derm ( \Lie_{Z^K}  \Phi  )  |\\
 && + \sum_{  |K| \leq |I| }      C(q_0)   \cdot C(|I|) \cdot E (  \lfloor \frac{|I|}{2} \rfloor  +3)  \cdot \frac{\eps \cdot  | \Lie_{Z^{K}} H_{L  L} | }{(1+t+|q|)^{1-  \de } \cdot (1+|q|)^{2+\gamma}} \; .
\eeaa

This gives
 \beaa
\notag
&&| \Lie_{Z^I}  ( g^{\la\mu} \derm_{\la}   \derm_{\mu}     \Phi_{V} ) - g^{\la\mu}    \derm_{\la}   \derm_{\mu}  (  \Lie_{Z^I} \Phi_{V}  ) |^2  \\
  \notag
   &\les&  \sum_{|K| \leq |I|  - 1}  | g^{\la\mu} \cdot \derm_{\la}   \derm_{\mu} (  \Lie_{Z^{K}}  \Phi_{V} ) |^2 \\
&& +  C(q_0)   \cdot  c (\delta) \cdot c (\gamma) \cdot C(  \lfloor \frac{|I|}{2} \rfloor ) \cdot E (    \lfloor \frac{|I|}{2} \rfloor   +4)  \\
&& \times \Big(  \sum_{ |K| \leq |I| }   \frac{\eps^2   }{(1+t+|q|)^{4-   c (\gamma)  \cdot c (\delta)  \cdot c(  \lfloor \frac{|I|}{2} \rfloor ) \cdot E (  \lfloor \frac{|I|}{2} \rfloor + 4)\cdot  \eps } \cdot (1+|q|)^2}  \, \cdot  |    \Lie_{Z^K}  H   |^2  \\
&& + \sum_{ |K| \leq |I| }  \frac{\eps^2 }{(1+t+|q|)^{4-   c (\gamma)  \cdot c (\delta)  \cdot c( \lfloor \frac{|I|}{2} \rfloor ) \cdot E (  \lfloor \frac{|I|}{2} \rfloor + 4)\cdot  \eps } }  \, \cdot | \derm ( \Lie_{Z^K}  \Phi )  |^2 \Big)   \\
   && + \sum_{  |K| \leq |I| }  C(q_0) \cdot   c (\delta) \cdot c (\gamma)  \cdot E (  4)  \cdot \frac{ \eps^2   }{ (1+t+|q|)^2  \cdot (1+|q|)^2 }  \, \cdot  | \derm ( \Lie_{Z^K}  \Phi  )  |^2 \\
&&+  C(q_0)   \cdot  c (\delta) \cdot c (\gamma) \cdot C(|I|) \cdot E (   \lfloor \frac{|I|}{2} \rfloor   +3)  \cdot \frac{\eps^2 }{(1+t+|q|)^{2-   c (\gamma)  \cdot c (\delta)  \cdot c(|I|) \cdot E (\lfloor \frac{|I|}{2} \rfloor  + 2)\cdot  \eps } \cdot (1+|q|)^2} \\
&& \times    \sum_{  |K| \leq |I| -1 }  | \derm ( \Lie_{Z^K}  \Phi  )  |^2 \\
 && + \sum_{  |K| \leq |I| }      C(q_0)   \cdot C(|I|) \cdot E (  \lfloor \frac{|I|}{2} \rfloor  +3)  \cdot \frac{\eps^2 \cdot  | \Lie_{Z^{K}} H_{L  L} |^2 }{(1+t+|q|)^{2- 2 \de } \cdot (1+|q|)^{4+2\gamma}} \;.
\eeaa

In particular, we get 

 \beaa
\notag
&&  \frac{(1+ t )}{\eps} \cdot  | \Lie_{Z^I}  ( g^{\la\mu} \derm_{\la}   \derm_{\mu}     \Phi_{V} ) - g^{\la\mu}    \derm_{\la}   \derm_{\mu}  (  \Lie_{Z^I} \Phi_{V}  ) |^2  \\
  \notag
   &\les& \frac{(1+ t )}{\eps} \cdot    \sum_{|K| \leq |I|  - 1}  | g^{\la\mu} \cdot \derm_{\la}   \derm_{\mu} (  \Lie_{Z^{K}}  \Phi_{V} ) |^2 \\
&& +  C(q_0)   \cdot  c (\delta) \cdot c (\gamma) \cdot C(  \lfloor \frac{|I|}{2} \rfloor ) \cdot E (    \lfloor \frac{|I|}{2} \rfloor   +4)  \\
&& \times \Big(  \sum_{ |K| \leq |I| }   \frac{\eps   }{(1+t+|q|)^{3-   c (\gamma)  \cdot c (\delta)  \cdot c(  \lfloor \frac{|I|}{2} \rfloor ) \cdot E (  \lfloor \frac{|I|}{2} \rfloor + 4)\cdot  \eps } \cdot (1+|q|)^2}  \, \cdot  |    \Lie_{Z^K}  H   |^2  \\
&& + \sum_{ |K| \leq |I| }  \frac{\eps }{(1+t+|q|)^{3-   c (\gamma)  \cdot c (\delta)  \cdot c( \lfloor \frac{|I|}{2} \rfloor ) \cdot E (  \lfloor \frac{|I|}{2} \rfloor + 4)\cdot  \eps } }  \, \cdot | \derm ( \Lie_{Z^K}  \Phi )  |^2 \Big)   \\
   && + \sum_{  |K| \leq |I| }  C(q_0) \cdot   c (\delta) \cdot c (\gamma)  \cdot E (  4)  \cdot \frac{ \eps   }{ (1+t+|q|)  \cdot (1+|q|)^2 }  \, \cdot   | \derm ( \Lie_{Z^K}  \Phi  )  |^2 \\
&&+  C(q_0)   \cdot  c (\delta) \cdot c (\gamma) \cdot C(|I|) \cdot E (   \lfloor \frac{|I|}{2} \rfloor   +3)  \cdot \frac{\eps }{(1+t+|q|)^{1-   c (\gamma)  \cdot c (\delta)  \cdot c(|I|) \cdot E ( \lfloor \frac{|I|}{2} \rfloor  + 2)\cdot  \eps } \cdot (1+|q|)^2} \\
&& \times    \sum_{  |K| \leq |I| -1 } | \derm ( \Lie_{Z^K}  \Phi  )  |^2 \\
 && + \sum_{  |K| \leq |I| }      C(q_0)   \cdot C(|I|) \cdot E (  \lfloor \frac{|I|}{2} \rfloor  +3)  \cdot \frac{\eps^2 \cdot  | \Lie_{Z^{K}} H_{L  L} |^2 }{(1+t+|q|)^{2- 2 \de } \cdot (1+|q|)^{4+2\gamma}} \;.
\eeaa

\end{proof}

\subsection{Hardy type ineuqlaity on the commutator}\

We recall from \cite{G4}, that we showed the following Hardy type inequality, that we sum up in the following Corollary.

\begin{corollary}\label{HardytypeinequalityforintegralstartingatROm}
Let $w$ defined as in Definition \ref{defoftheweightw}, where $\ga > 0$.
Let  $\Phi$ a tensor that decays fast enough at spatial infinity for all time $t$\,, such that
\bea
 \int_{\SSS^{2}} \lim_{r \to \infty} \Big( \frac{r^{2}}{(1+t+r)^{a} \cdot (1+|q|) } \cdot w(q) \cdot | \Phi_{V} |^2   \Big)  d\si^{2} (t ) &=& 0 \; .
\eea
Let $R(\Om)  \geq 0 $\,, be a function of $\Om \in \SSS^{n-1}$\,. Then, since $\ga \neq 0$\,, we have for $0 \leq a \leq 2$\,, that 
\bea
\notag
 &&   \int_{\SSS^{2}} \int_{r=R(\Om)}^{r=\infty} \frac{r^{2}}{(1+t+r)^{a}} \cdot   \frac{w (q)}{(1+|q|)^2} \cdot  |\Phi_{V} |^2   \cdot dr  \cdot d\si^{2}  \\
 \notag
 &\leq& c(\ga) \cdot  \int_{\SSS^{n-1}} \int_{r=R(\Om)}^{r=\infty}  \frac{ r^{2}}{(1+t+r)^{a}} \cdot w(q) \cdot  | \pa_r\Phi_{V} |^2  \cdot  dr  \cdot d\si^{2}  \; , \\
 \eea
 where the constant $c(\ga)$ does not depend on $R(\Om)$\,. In particular, we have the following estimate in the exterior,
 \bea
\notag
  \int_{\Sigma^{ext}_{\tau} }  \frac{1}{(1+t+r)^{a}} \cdot   \frac{w (q)}{(1+|q|)^2} \cdot  |\Phi_{V} |^2 &\leq& c(\ga) \cdot  \int_{\Sigma^{ext}_{\tau} }   \frac{ 1 }{(1+t+r)^{a}} \cdot w(q) \cdot  | \pa_r\Phi_{V} |^2     \; , \\
 \eea
  where $\int_{\Sigma^{ext}_{\tau} } $ is taken with respect to the measure $r^2 \cdot dr  \cdot d\si^{2} $\,. And, in particular, if 
 \beaa
 \int_{\SSS^{2}} \lim_{r \to \infty} \Big( \frac{r^{2}}{(1+t+r)^{a} \cdot (1+|q|) } \cdot w(q) \cdot | \Phi |^2   \Big)  d\si^{2} (t ) &=& 0 \; .
\eeaa
then,
  \bea
\notag
 \int_{\Sigma^{ext}_{\tau} }  \frac{1}{(1+t+r)^{a}} \cdot   \frac{w (q)}{(1+|q|)^2} \cdot  |\Phi |^2  &\leq& c(\ga) \cdot  \int_{\Sigma^{ext}_{\tau} }   \frac{ 1 }{(1+t+r)^{a}} \cdot w(q) \cdot  | \derm \Phi |^2    \; . \\
 \eea

 \end{corollary}

 \begin{lemma}\label{Hardytypeineeuqlityappliedtothecommutatortermplusatimeintegralformula}
For $\eps$ small enough, depending on $\ga$\,, $\de$ and $|I|$\,, we have
                                           \beaa
   \notag
&&  \int_{\Sigma^{ext}_{\tau} }   \frac{(1+ t )}{\eps} \cdot  | \Lie_{Z^I}  ( g^{\la\mu} \derm_{\la}   \derm_{\mu}     \Phi_{V} ) - g^{\la\mu}    \derm_{\la}   \derm_{\mu}  (  \Lie_{Z^I} \Phi_{V}  ) |^2 \Big)  \cdot w(q)   \\
  &\leq&      \int_{\Sigma^{ext}_{\tau} }   \Big[       \frac{(1+ t )}{\eps} \cdot    \sum_{|K| \leq |I|  - 1}  | g^{\la\mu} \cdot \derm_{\la}   \derm_{\mu} (  \Lie_{Z^{K}}  \Phi_{V} ) |^2 \\
&& +  C(q_0)   \cdot  c (\delta) \cdot c (\gamma) \cdot C(  \lfloor \frac{|I|}{2} \rfloor ) \cdot E (    \lfloor \frac{|I|}{2} \rfloor   +4)  \\
&& \times  \sum_{ |K| \leq |I| }  \Big(   \frac{\eps   }{(1+t+|q|)}  \, \cdot  |  \derm (   \Lie_{Z^K}  h )   |^2  +  \frac{\eps }{(1+t+|q|) }  \, \cdot | \derm ( \Lie_{Z^K}  \Phi)  |^2 \Big)   \\
&&+  C(q_0)   \cdot  c (\delta) \cdot c (\gamma) \cdot C(|I|) \cdot E (    \lfloor \frac{|I|}{2} \rfloor   +3)  \cdot \frac{\eps }{(1+t+|q|)^{1-   c (\gamma)  \cdot c (\delta)  \cdot c(|I|) \cdot E (\lfloor \frac{|I|}{2} \rfloor  + 2)\cdot  \eps } \cdot (1+|q|)^2} \\
&& \times    \sum_{  |K| \leq |I| -1 } | \derm ( \Lie_{Z^K}  \Phi  )  |^2  \Big] \cdot w(q)\\
&& +   \int_{\Sigma^{ext}_{\tau} }   \sum_{  |K| \leq |I| }       \Big[   C(q_0)   \cdot C(|I|) \cdot E (  \lfloor \frac{|I|}{2} \rfloor  +3)  \cdot \frac{\eps \cdot  | \Lie_{Z^{K}} H_{L  L} |^2 }{(1+t+|q|)^{1- 2 \de } \cdot (1+|q|)^{4+2\gamma}}  \Big] \cdot w(q) \; .
\eeaa
\end{lemma}

\begin{proof}

Using Lemma \ref{estimateonthecommutatortermusingtheproductsandusingbootstrap}, we get for $\eps$ small enough, depending on $\ga$\,, $\de$ and $|I|$\,, the following energy estimate on the commutator term

  \beaa
\notag
&&  \frac{(1+ t )}{\eps} \cdot  | \Lie_{Z^I}  ( g^{\la\mu} \derm_{\la}   \derm_{\mu}     \Phi_{V} ) - g^{\la\mu}    \derm_{\la}   \derm_{\mu}  (  \Lie_{Z^I}  \Phi_{V}  ) |^2  \\
  \notag
   &\les& \frac{(1+ t )}{\eps} \cdot    \sum_{|K| \leq |I|  - 1}  | g^{\la\mu} \cdot \derm_{\la}   \derm_{\mu} (  \Lie_{Z^{K}}   \Phi_{V}  ) |^2 \\
&& +  C(q_0)   \cdot  c (\delta) \cdot c (\gamma) \cdot C(  \lfloor \frac{|I|}{2} \rfloor ) \cdot E (    \lfloor \frac{|I|}{2} \rfloor   +4)  \\
&& \times \Big(  \sum_{ |K| \leq |I| }   \frac{\eps   }{(1+t+|q|)\cdot (1+|q|)^2}  \, \cdot  |    \Lie_{Z^K}  H   |^2  + \sum_{ |K| \leq |I| }  \frac{\eps }{(1+t+|q|) }  \, \cdot | \derm ( \Lie_{Z^K}  \Phi  )  |^2 \Big)   \\
&&+  C(q_0)   \cdot  c (\delta) \cdot c (\gamma) \cdot C(|I|) \cdot E (   \lfloor \frac{|I|}{2} \rfloor   +3)  \cdot \frac{\eps }{(1+t+|q|)^{1-   c (\gamma)  \cdot c (\delta)  \cdot c(|I|) \cdot E ( \lfloor \frac{|I|}{2} \rfloor  + 2)\cdot  \eps } \cdot (1+|q|)^2} \\
&& \times    \sum_{  |K| \leq |I| -1 } | \derm ( \Lie_{Z^K}  \Phi  )  |^2 \\
 && + \sum_{  |K| \leq |I| }      C(q_0)   \cdot C(|I|) \cdot E (  \lfloor \frac{|I|}{2} \rfloor  +3)  \cdot \frac{\eps \cdot  | \Lie_{Z^{K}} H_{L  L} |^2 }{(1+t+|q|)^{1- 2 \de } \cdot (1+|q|)^{4+2\gamma}} \;.
\eeaa

\begin{remark}\label{remarkonwhyepsilonsmalldesonotdependontheconstantsEofkwhichareboundedby1}
Note that one can choose $\eps$ small not depending on the constants $E(k)\,, k \in \N$\,, since these are bounded by $E(k) \leq 1$\,.
\end{remark}

                                               \beaa
   \notag
&&  \int_{\Sigma^{ext}_{\tau} }   \frac{(1+ t )}{\eps} \cdot  | \Lie_{Z^I}  ( g^{\la\mu} \derm_{\la}   \derm_{\mu}     \ \Phi_{V}) - g^{\la\mu}    \derm_{\la}   \derm_{\mu}  (  \Lie_{Z^I}  \Phi_{V} ) |^2 \Big)  \cdot w(q)   \\
 &\leq&      \int_{\Sigma^{ext}_{\tau} }   \Big[       \frac{(1+ t )}{\eps} \cdot    \sum_{|K| \leq |I|  - 1}  | g^{\la\mu} \cdot \derm_{\la}   \derm_{\mu} (  \Lie_{Z^{K}}   \Phi_{V} ) |^2 \\
&& +  C(q_0)   \cdot  c (\delta) \cdot c (\gamma) \cdot C(  \lfloor \frac{|I|}{2} \rfloor ) \cdot E (    \lfloor \frac{|I|}{2} \rfloor   +4)  \\
&& \times  \sum_{ |K| \leq |I| }  \Big(   \frac{\eps   }{(1+t+|q|)}  \, \cdot  |  \derm (   \Lie_{Z^K}  H )   |^2  +  \frac{\eps }{(1+t+|q|) }  \, \cdot | \derm ( \Lie_{Z^K}  \Phi)  |^2 \Big)   \\
&&+  C(q_0)   \cdot  c (\delta) \cdot c (\gamma) \cdot C(|I|) \cdot E (    \lfloor \frac{|I|}{2} \rfloor   +3)  \cdot \frac{\eps }{(1+t+|q|)^{1-   c (\gamma)  \cdot c (\delta)  \cdot c(|I|) \cdot E ( \lfloor \frac{|I|}{2} \rfloor  + 2)\cdot  \eps } \cdot (1+|q|)^2} \\
&& \times    \sum_{  |K| \leq |I| -1 } | \derm ( \Lie_{Z^K}  \Phi  )  |^2  \Big] \cdot w(q)  \\
&& +   \int_{\Sigma^{ext}_{\tau} }   \sum_{  |K| \leq |I| }       \Big[   C(q_0)   \cdot C(|I|) \cdot E (  \lfloor \frac{|I|}{2} \rfloor  +3)  \cdot \frac{\eps \cdot  | \Lie_{Z^{K}} H_{L  L} |^2 }{(1+t+|q|)^{1- 2 \de } \cdot (1+|q|)^{4+2\gamma}}  \Big] \cdot w(q) 
\eeaa
\beaa
 &\leq&      \int_{\Sigma^{ext}_{\tau} }   \Big[      \frac{(1+ t )}{\eps} \cdot    \sum_{|K| \leq |I|  - 1}  | g^{\la\mu} \cdot \derm_{\la}   \derm_{\mu} (  \Lie_{Z^{K}}   \Phi_{V} ) |^2 \\
&& +  C(q_0)   \cdot  c (\delta) \cdot c (\gamma) \cdot C(  \lfloor \frac{|I|}{2} \rfloor ) \cdot E (    \lfloor \frac{|I|}{2} \rfloor   +4)  \\
&& \times  \sum_{ |K| \leq |I| }  \Big(   \frac{\eps   }{(1+t+|q|)}  \, \cdot  |  \derm (   \Lie_{Z^K}  H )   |^2  +  \frac{\eps }{(1+t+|q|) }  \, \cdot | \derm ( \Lie_{Z^K}  \Phi )  |^2 \Big)   \\
&&+  C(q_0)   \cdot  c (\delta) \cdot c (\gamma) \cdot C(|I|) \cdot E (    \lfloor \frac{|I|}{2} \rfloor   +3)  \cdot \frac{\eps }{(1+t+|q|)^{1-   c (\gamma)  \cdot c (\delta)  \cdot c(|I|) \cdot E ( \lfloor \frac{|I|}{2} \rfloor  + 2)\cdot  \eps } \cdot (1+|q|)^2} \\
&& \times    \sum_{  |K| \leq |I| -1 } | \derm ( \Lie_{Z^K}  \Phi  )  |^2  \Big] \cdot w(q) \\
&& +   \int_{\Sigma^{ext}_{\tau} }   \sum_{  |K| \leq |I| }       \Big[   C(q_0)   \cdot C(|I|) \cdot E (  \lfloor \frac{|I|}{2} \rfloor  +3)  \cdot \frac{\eps \cdot  | \Lie_{Z^{K}} H_{L  L} |^2 }{(1+t+|q|)^{1- 2 \de } \cdot (1+|q|)^{4+2\gamma}}  \Big] \cdot w(q) \; .
\eeaa

Thus, we get the stated result.

\end{proof}

\section{The “bad” term $\Lie_{Z^I} \big( A_{e_a}  \cdot     \derm A_{e_a} \big)$ for the wave equation of the Yang-Mills potential}
 
We want to establish an energy estimate on the “good” term $\Lie_{Z^I} A_{e_{a}}$ that we could use to establish a pointwise estimate on the “bad” term $\Lie_{Z^I} \big( A_{e_a}  \cdot     \derm A_{e_a} \big)$.

\begin{lemma}\label{estimateonttimesthesquareofthesourcesofwaveoperatorongoodcomponenentsofEinsteeinYangMillspoentialA}
   
For $\ga \geq 3 \de $\;, for $0 < \de \leq \frac{1}{4}$\;, and for $\eps$ small enough depending on $\ga$\;, $\de$\;, and $|I|$\;, we have

                                                 \beaa
   \notag
&& \frac{(1+ t )}{\eps} \cdot  |  \Lie_{Z^I}  g^{\la\mu} \derm_{\la}   \derm_{\mu}   A_{\cal T}  |^2  \\
   &\les&          \sum_{|K| \leq |I |}  \Big[   \;          O \big(   C(q_0)   \cdot  c (\delta) \cdot c (\gamma) \cdot C(|I|) \cdot E ( |I| + 5)  \cdot \frac{\eps   \cdot  | \derm ( \Lie_{Z^K} h^1 ) |^2  }{(1+t+|q|)}   \big) \\
   && + O \big(    C(q_0)   \cdot  c (\delta) \cdot c (\gamma) \cdot C(|I|) \cdot E ( \lfloor \frac{|I|}{2} \rfloor +5)  \cdot \frac{\eps  \cdot  |\derm ( \Lie_{Z^K} A ) |^2  }{(1+t+|q|) }  \big) \\
   &&  +  C(q_0)   \cdot  c (\delta) \cdot c (\gamma) \cdot C(|I|) \cdot E (\lfloor \frac{|I|}{2} \rfloor + 4)  \cdot \frac{\eps   \cdot    | \rderm ( \Lie_{Z^K} A )  |^2  }{(1+t+|q|)^{1-      c (\gamma)  \cdot c (\delta)  \cdot c(|I|) \cdot E ( \lfloor \frac{|I|}{2} \rfloor + 4) \cdot \eps } \cdot   (1+|q|) }    \\
      \notag
        &&    +     C(q_0)   \cdot  c (\delta) \cdot c (\gamma) \cdot C(|I|) \cdot E ( \lfloor \frac{|I|}{2} \rfloor + 4)  \cdot \frac{\eps  \cdot   | \rderm  ( \Lie_{Z^K} h^1 ) |^2  }{(1+t+|q|)^{1-      c (\gamma)  \cdot c (\delta)  \cdot c(|I|) \cdot E (\lfloor \frac{|I|}{2} \rfloor+ 4) \cdot \eps } \cdot (1+|q|)^{2 }}    \\
      &&    +       O \big(     C(q_0)   \cdot  c (\delta) \cdot c (\gamma) \cdot C(|I|) \cdot E ( \lfloor \frac{|I|}{2} \rfloor+ 5)  \cdot \frac{\eps \cdot |  \Lie_{Z^K} A  |^2  }{(1+t+|q|) \cdot (1+|q|)^{2  }}   \\
         &&    +         O \big( C(q_0)   \cdot  c (\delta) \cdot c (\gamma) \cdot C(|I|) \cdot E (\lfloor \frac{|I|}{2} \rfloor + 4)  \cdot \frac{\eps \cdot | \Lie_{Z^K} h^1 |^2  }{(1+t+|q|) \cdot (1+|q|)^{2 }}     \big)   \;  \Big] \\
            &&    +     C(q_0)   \cdot  c (\delta) \cdot c (\gamma) \cdot C(|I|) \cdot E ( \lfloor \frac{|I|}{2} \rfloor + 5)  \cdot \frac{\eps^3 }{(1+t+|q|)^{5-      c (\gamma)  \cdot c (\delta)  \cdot c(|I|) \cdot E ( \lfloor \frac{|I|}{2} \rfloor  + 5) \cdot \eps } \cdot (1+|q|)^{2+2\gamma - 4\de }}    \; .
\eeaa
\end{lemma}

\begin{proof}
Based on Lemma \ref{TheactualusefulstrzuctureofthesourcetermsforthewaveequationonpoentialAusingbootstrap}, we have
                                                       \beaa
   \notag
&&  |  \Lie_{Z^I}  g^{\la\mu} \derm_{\la}   \derm_{\mu}   A_{\cal T}  |  \\
   &\les&      \sum_{|K| \leq |I |}  \Big[   \;            O \big(   C(q_0)   \cdot  c (\delta) \cdot c (\gamma) \cdot C(|I|) \cdot E ( \lfloor \frac{|I|}{2} \rfloor + 5)  \cdot \frac{\eps   \cdot  | \derm ( \Lie_{Z^K} h^1 ) | }{(1+t+|q|)^{2-      c (\gamma)  \cdot c (\delta)  \cdot c(|I|) \cdot E (  \lfloor \frac{|I|}{2} \rfloor + 5) \cdot \eps } \cdot (1+|q|)^{\gamma - 2\de }}   \big) \\
   && + O \big(    C(q_0)   \cdot  c (\delta) \cdot c (\gamma) \cdot C(|I|) \cdot E (   \lfloor \frac{|I|}{2} \rfloor   +5)  \cdot \frac{\eps  \cdot  |\derm ( \Lie_{Z^K} A ) | }{(1+t+|q|)^{2-   c (\gamma)  \cdot c (\delta)  \cdot c(|I|) \cdot E ( \lfloor \frac{|I|}{2} \rfloor + 5)\cdot  \eps } }  \big) \\
   &&  +  C(q_0)   \cdot  c (\delta) \cdot c (\gamma) \cdot C(|I|) \cdot E ( \lfloor \frac{|I|}{2} \rfloor  + 4)  \cdot \frac{\eps   \cdot    | \rderm ( \Lie_{Z^K} A )  | }{(1+t+|q|)^{1-      c (\gamma)  \cdot c (\delta)  \cdot c(|I|) \cdot E (  \lfloor \frac{|I|}{2} \rfloor + 4) \cdot \eps } \cdot  (1+|q|)^{ \frac{1}{2}}  }    \\
      \notag
        &&    +     C(q_0)   \cdot  c (\delta) \cdot c (\gamma) \cdot C(|I|) \cdot E ( \lfloor \frac{|I|}{2} \rfloor + 4)  \cdot \frac{\eps  \cdot   | \rderm  ( \Lie_{Z^K} h^1 ) | }{(1+t+|q|)^{1-      c (\gamma)  \cdot c (\delta)  \cdot c(|I|) \cdot E ( \lfloor \frac{|I|}{2} \rfloor + 4) \cdot \eps } \cdot (1+|q|)^{1+\gamma - 2\de }}    \\
      &&    +       O \big(     C(q_0)   \cdot  c (\delta) \cdot c (\gamma) \cdot C(|I|) \cdot E ( \lfloor \frac{|I|}{2} \rfloor  + 5)  \cdot \frac{\eps \cdot |  \Lie_{Z^K} A  | }{(1+t+|q|)^{2-      c (\gamma)  \cdot c (\delta)  \cdot c(|I|) \cdot E ( \lfloor \frac{|I|}{2} \rfloor + 5) \cdot \eps } \cdot (1+|q|)^{\gamma - 2\de }}   \\
         &&    +         O \big( C(q_0)   \cdot  c (\delta) \cdot c (\gamma) \cdot C(|I|) \cdot E ( \lfloor \frac{|I|}{2} \rfloor + 4)  \cdot \frac{\eps \cdot | \Lie_{Z^K} h^1 | }{(1+t+|q|)^{2-      c (\gamma)  \cdot c (\delta)  \cdot c(|I|) \cdot E (  \lfloor \frac{|I|}{2} \rfloor + 4) \cdot \eps } \cdot (1+|q|)^{1+2\gamma - 4\de }}     \big)  \;  \Big]  \\
            &&    +     C(q_0)   \cdot  c (\delta) \cdot c (\gamma) \cdot C(|I|) \cdot E ( \lfloor \frac{|I|}{2} \rfloor + 5)  \cdot \frac{\eps^2 }{(1+t+|q|)^{3-      c (\gamma)  \cdot c (\delta)  \cdot c(|I|) \cdot E ( \lfloor \frac{|I|}{2} \rfloor  + 5) \cdot \eps } \cdot (1+|q|)^{1+\gamma - 2\de }}    \;.
\eeaa

                                                       \beaa
   \notag
&&  |  \Lie_{Z^I}  g^{\la\mu} \derm_{\la}   \derm_{\mu}   A_{\cal T} |^2  \\
   &\les&       \sum_{|K| \leq |I |}  \Big[   \;           O \big(   C(q_0)   \cdot  c (\delta) \cdot c (\gamma) \cdot C(|I|) \cdot E (  \lfloor \frac{|I|}{2} \rfloor + 5)  \cdot \frac{\eps   \cdot  | \derm ( \Lie_{Z^K} h^1 ) |^2  }{(1+t+|q|)^{4-      c (\gamma)  \cdot c (\delta)  \cdot c(|I|) \cdot E (  \lfloor \frac{|I|}{2} \rfloor + 5) \cdot \eps } \cdot (1+|q|)^{2\gamma - 4\de }}   \big) \\
   && + O \big(    C(q_0)   \cdot  c (\delta) \cdot c (\gamma) \cdot C(|I|) \cdot E (  \lfloor \frac{|I|}{2} \rfloor  +5)  \cdot \frac{\eps  \cdot  |\derm ( \Lie_{Z^K} A ) |^2  }{(1+t+|q|)^{4-   c (\gamma)  \cdot c (\delta)  \cdot c(|I|) \cdot E (  \lfloor \frac{|I|}{2} \rfloor + 5)\cdot  \eps } }  \big) \\
   &&  +  C(q_0)   \cdot  c (\delta) \cdot c (\gamma) \cdot C(|I|) \cdot E ( \lfloor \frac{|I|}{2} \rfloor  + 4)  \cdot \frac{\eps   \cdot    | \rderm ( \Lie_{Z^K} A )  |^2  }{(1+t+|q|)^{2-      c (\gamma)  \cdot c (\delta)  \cdot c(|I|) \cdot E ( \lfloor \frac{|I|}{2} \rfloor + 4) \cdot \eps } \cdot  (1+|q|)  }    \\
      \notag
        &&    +     C(q_0)   \cdot  c (\delta) \cdot c (\gamma) \cdot C(|I|) \cdot E (  \lfloor \frac{|I|}{2} \rfloor  + 4)  \cdot \frac{\eps  \cdot   | \rderm  ( \Lie_{Z^K} h^1 ) |^2  }{(1+t+|q|)^{2-      c (\gamma)  \cdot c (\delta)  \cdot c(|I|) \cdot E (  \lfloor \frac{|I|}{2} \rfloor + 4) \cdot \eps } \cdot (1+|q|)^{2+2\gamma - 4\de }}    \\
      &&    +       O \big(     C(q_0)   \cdot  c (\delta) \cdot c (\gamma) \cdot C(|I|) \cdot E (  \lfloor \frac{|I|}{2} \rfloor + 5)  \cdot \frac{\eps \cdot |  \Lie_{Z^K} A  |^2  }{(1+t+|q|)^{4-      c (\gamma)  \cdot c (\delta)  \cdot c(|I|) \cdot E (  \lfloor \frac{|I|}{2} \rfloor + 5) \cdot \eps } \cdot (1+|q|)^{2\gamma - 4\de }}   \\
         &&    +         O \big( C(q_0)   \cdot  c (\delta) \cdot c (\gamma) \cdot C(|I|) \cdot E (  \lfloor \frac{|I|}{2} \rfloor + 4)  \cdot \frac{\eps \cdot | \Lie_{Z^K} h^1 |^2  }{(1+t+|q|)^{4-      c (\gamma)  \cdot c (\delta)  \cdot c(|I|) \cdot E (  \lfloor \frac{|I|}{2} \rfloor + 4) \cdot \eps } \cdot (1+|q|)^{2+4\gamma - 8\de }}     \big)  \;  \Big]  \\
            &&    +     C(q_0)   \cdot  c (\delta) \cdot c (\gamma) \cdot C(|I|) \cdot E ( \lfloor \frac{|I|}{2} \rfloor + 5)  \cdot \frac{\eps^4 }{(1+t+|q|)^{6-      c (\gamma)  \cdot c (\delta)  \cdot c(|I|) \cdot E ( \lfloor \frac{|I|}{2} \rfloor  + 5) \cdot \eps } \cdot (1+|q|)^{2+ 2 \gamma - 4\de }}    \;.
\eeaa

Hence,
                                                       \beaa
   \notag
&& \frac{(1+ t )}{\eps} \cdot  |  \Lie_{Z^I}  g^{\la\mu} \derm_{\la}   \derm_{\mu}   A_{\cal T}  |^2  \\
   &\les&         \sum_{|K| \leq |I |}  \Big[   \;         O \big(   C(q_0)   \cdot  c (\delta) \cdot c (\gamma) \cdot C(|I|) \cdot E ( \lfloor \frac{|I|}{2} \rfloor + 5)  \cdot \frac{\eps   \cdot  | \derm ( \Lie_{Z^K} h^1 ) |^2  }{(1+t+|q|)^{3-      c (\gamma)  \cdot c (\delta)  \cdot c(|I|) \cdot E (  \lfloor \frac{|I|}{2} \rfloor + 5) \cdot \eps } \cdot (1+|q|)^{2\gamma - 4\de }}   \big) \\
   && + O \big(    C(q_0)   \cdot  c (\delta) \cdot c (\gamma) \cdot C(|I|) \cdot E (  \lfloor \frac{|I|}{2} \rfloor  +5)  \cdot \frac{\eps  \cdot  |\derm ( \Lie_{Z^K} A ) |^2  }{(1+t+|q|)^{3-   c (\gamma)  \cdot c (\delta)  \cdot c(|I|) \cdot E ( \lfloor \frac{|I|}{2} \rfloor + 5)\cdot  \eps } }  \big) \\
   &&  +  C(q_0)   \cdot  c (\delta) \cdot c (\gamma) \cdot C(|I|) \cdot E (  \lfloor \frac{|I|}{2} \rfloor + 4)  \cdot \frac{\eps   \cdot    | \rderm ( \Lie_{Z^K} A )  |^2  }{(1+t+|q|)^{1-      c (\gamma)  \cdot c (\delta)  \cdot c(|I|) \cdot E ( \lfloor \frac{|I|}{2} \rfloor + 4) \cdot \eps } \cdot   (1+|q|) }    \\
      \notag
        &&    +     C(q_0)   \cdot  c (\delta) \cdot c (\gamma) \cdot C(|I|) \cdot E ( \lfloor \frac{|I|}{2} \rfloor  + 4)  \cdot \frac{\eps  \cdot   | \rderm  ( \Lie_{Z^K} h^1 ) |^2  }{(1+t+|q|)^{1-      c (\gamma)  \cdot c (\delta)  \cdot c(|I|) \cdot E (  \lfloor \frac{|I|}{2} \rfloor + 4) \cdot \eps } \cdot (1+|q|)^{2+2\gamma - 4\de }}    \\
      &&    +       O \big(     C(q_0)   \cdot  c (\delta) \cdot c (\gamma) \cdot C(|I|) \cdot E (  \lfloor \frac{|I|}{2} \rfloor + 5)  \cdot \frac{\eps \cdot |  \Lie_{Z^K} A  |^2  }{(1+t+|q|)^{3-      c (\gamma)  \cdot c (\delta)  \cdot c(|I|) \cdot E ( \lfloor \frac{|I|}{2} \rfloor + 5) \cdot \eps } \cdot (1+|q|)^{2\gamma - 4\de }}   \\
         &&    +         O \big( C(q_0)   \cdot  c (\delta) \cdot c (\gamma) \cdot C(|I|) \cdot E (  \lfloor \frac{|I|}{2} \rfloor + 4)  \cdot \frac{\eps \cdot | \Lie_{Z^K} h^1 |^2  }{(1+t+|q|)^{3-      c (\gamma)  \cdot c (\delta)  \cdot c(|I|) \cdot E ( \lfloor \frac{|I|}{2} \rfloor + 4) \cdot \eps } \cdot (1+|q|)^{2+4\gamma - 8\de }}     \big)  \;  \Big]   \\
                 &&    +     C(q_0)   \cdot  c (\delta) \cdot c (\gamma) \cdot C(|I|) \cdot E ( \lfloor \frac{|I|}{2} \rfloor + 5)  \cdot \frac{\eps^3 }{(1+t+|q|)^{5-      c (\gamma)  \cdot c (\delta)  \cdot c(|I|) \cdot E ( \lfloor \frac{|I|}{2} \rfloor  + 5) \cdot \eps } \cdot (1+|q|)^{2+ 2 \gamma - 4\de }}    \;.   
\eeaa

In fact, we have

                                                       \beaa
   \notag
&& \frac{(1+ t )}{\eps} \cdot  |  \Lie_{Z^I}  g^{\la\mu} \derm_{\la}   \derm_{\mu}   A_{\cal T}  |^2  \\
   &\les&            \sum_{|K| \leq |I |}  \Big[   \;        O \big(   C(q_0)   \cdot  c (\delta) \cdot c (\gamma) \cdot C(|I|) \cdot E (\lfloor \frac{|I|}{2} \rfloor+ 5)  \cdot \frac{\eps   \cdot  | \derm ( \Lie_{Z^K} h^1 ) |^2  }{(1+t+|q|)}   \big) \\
   && + O \big(    C(q_0)   \cdot  c (\delta) \cdot c (\gamma) \cdot C(|I|) \cdot E ( \lfloor \frac{|I|}{2} \rfloor +5)  \cdot \frac{\eps  \cdot  |\derm ( \Lie_{Z^K} A ) |^2  }{(1+t+|q|) }  \big) \\
   &&  +  C(q_0)   \cdot  c (\delta) \cdot c (\gamma) \cdot C(|I|) \cdot E (\lfloor \frac{|I|}{2} \rfloor + 4)  \cdot \frac{\eps   \cdot    | \rderm ( \Lie_{Z^K} A )  |^2  }{(1+t+|q|)^{1-      c (\gamma)  \cdot c (\delta)  \cdot c(|I|) \cdot E (\lfloor \frac{|I|}{2} \rfloor+ 4) \cdot \eps } \cdot   (1+|q|) }    \\
      \notag
        &&    +     C(q_0)   \cdot  c (\delta) \cdot c (\gamma) \cdot C(|I|) \cdot E (\lfloor \frac{|I|}{2} \rfloor + 4)  \cdot \frac{\eps  \cdot   | \rderm  ( \Lie_{Z^K} h^1 ) |^2  }{(1+t+|q|)^{1-      c (\gamma)  \cdot c (\delta)  \cdot c(|I|) \cdot E ( \lfloor \frac{|I|}{2} \rfloor+ 4) \cdot \eps } \cdot (1+|q|)^{2+2\gamma - 4\de }}    \\
      &&    +       O \big(     C(q_0)   \cdot  c (\delta) \cdot c (\gamma) \cdot C(|I|) \cdot E (\lfloor \frac{|I|}{2} \rfloor + 5)  \cdot \frac{\eps \cdot |  \Lie_{Z^K} A  |^2  }{(1+t+|q|) \cdot (1+|q|)^{2-      c (\gamma)  \cdot c (\delta)  \cdot c(|I|) \cdot E ( \lfloor \frac{|I|}{2} \rfloor+ 5) \cdot \eps + 2\gamma - 4\de }}   \\
         &&    +         O \big( C(q_0)   \cdot  c (\delta) \cdot c (\gamma) \cdot C(|I|) \cdot E ( \lfloor \frac{|I|}{2} \rfloor + 4)  \cdot \frac{\eps \cdot | \Lie_{Z^K} h^1 |^2  }{(1+t+|q|) \cdot (1+|q|)^{2+4\gamma - 8\de }}     \big) \;  \Big] \\
               &&    +     C(q_0)   \cdot  c (\delta) \cdot c (\gamma) \cdot C(|I|) \cdot E ( \lfloor \frac{|I|}{2} \rfloor + 5)  \cdot \frac{\eps^3 }{(1+t+|q|)^{5-      c (\gamma)  \cdot c (\delta)  \cdot c(|I|) \cdot E ( \lfloor \frac{|I|}{2} \rfloor  + 5) \cdot \eps } \cdot (1+|q|)^{2+ 2 \gamma - 4\de }}     \;. 
\eeaa
Now, for $\ga \geq 2 \de + \de$\;, we have $2\ga - 4  \de \geq 2 \de $ and therefore, for $\eps$ small enough, depending on $\ga$\,, $\de$ and on $|I|$ (and not on $E ( |I|+ 5)$ -- see Remark \ref{remarkonwhyepsilonsmalldesonotdependontheconstantsEofkwhichareboundedby1}), we have
\beaa
2\ga - 4  \de -   c (\gamma)  \cdot c (\delta)  \cdot c(|I|) \cdot E ( |I|+ 5) \cdot \eps \geq 2\de - c (\gamma)  \cdot c (\delta)  \cdot c(|I|)  \cdot \eps \geq 0 \; ,
\eeaa
and therefore,

                                                       \beaa
   \notag
&& \frac{(1+ t )}{\eps} \cdot  |  \Lie_{Z^I}  g^{\la\mu} \derm_{\la}   \derm_{\mu}   A_{\cal T}  |^2  \\
   &\les&          \sum_{|K| \leq |I |}  \Big[   \;          O \big(   C(q_0)   \cdot  c (\delta) \cdot c (\gamma) \cdot C(|I|) \cdot E ( |I| + 5)  \cdot \frac{\eps   \cdot  | \derm ( \Lie_{Z^K} h^1 ) |^2  }{(1+t+|q|)}   \big) \\
   && + O \big(    C(q_0)   \cdot  c (\delta) \cdot c (\gamma) \cdot C(|I|) \cdot E ( \lfloor \frac{|I|}{2} \rfloor +5)  \cdot \frac{\eps  \cdot  |\derm ( \Lie_{Z^K} A ) |^2  }{(1+t+|q|) }  \big) \\
   &&  +  C(q_0)   \cdot  c (\delta) \cdot c (\gamma) \cdot C(|I|) \cdot E (\lfloor \frac{|I|}{2} \rfloor + 4)  \cdot \frac{\eps   \cdot    | \rderm ( \Lie_{Z^K} A )  |^2  }{(1+t+|q|)^{1-      c (\gamma)  \cdot c (\delta)  \cdot c(|I|) \cdot E ( \lfloor \frac{|I|}{2} \rfloor + 4) \cdot \eps } \cdot   (1+|q|) }    \\
      \notag
        &&    +     C(q_0)   \cdot  c (\delta) \cdot c (\gamma) \cdot C(|I|) \cdot E ( \lfloor \frac{|I|}{2} \rfloor + 4)  \cdot \frac{\eps  \cdot   | \rderm  ( \Lie_{Z^K} h^1 ) |^2  }{(1+t+|q|)^{1-      c (\gamma)  \cdot c (\delta)  \cdot c(|I|) \cdot E (\lfloor \frac{|I|}{2} \rfloor+ 4) \cdot \eps } \cdot (1+|q|)^{2 }}    \\
      &&    +       O \big(     C(q_0)   \cdot  c (\delta) \cdot c (\gamma) \cdot C(|I|) \cdot E ( \lfloor \frac{|I|}{2} \rfloor+ 5)  \cdot \frac{\eps \cdot |  \Lie_{Z^K} A  |^2  }{(1+t+|q|) \cdot (1+|q|)^{2  }}   \\
         &&    +         O \big( C(q_0)   \cdot  c (\delta) \cdot c (\gamma) \cdot C(|I|) \cdot E (\lfloor \frac{|I|}{2} \rfloor + 4)  \cdot \frac{\eps \cdot | \Lie_{Z^K} h^1 |^2  }{(1+t+|q|) \cdot (1+|q|)^{2 }}     \big)   \;  \Big] \\
               &&    +     C(q_0)   \cdot  c (\delta) \cdot c (\gamma) \cdot C(|I|) \cdot E ( \lfloor \frac{|I|}{2} \rfloor + 5)  \cdot \frac{\eps^3 }{(1+t+|q|)^{5-      c (\gamma)  \cdot c (\delta)  \cdot c(|I|) \cdot E ( \lfloor \frac{|I|}{2} \rfloor  + 5) \cdot \eps } \cdot (1+|q|)^{2+ 2 \gamma - 4\de }}    \;. 
\eeaa

\end{proof}

\begin{lemma}\label{estimateonthecontributionofhzerointhesourcetermsforthewaveequationontheYangMillspoential}
We have
\bea
\notag
 && \int_{\Sigma^{ext}_{\tau} }    C(q_0)   \cdot  c (\delta) \cdot c (\gamma) \cdot C(|I|) \cdot E ( \lfloor \frac{|I|}{2} \rfloor + 5)  \\
 \notag
           && \times \frac{\eps^3 }{(1+\tau+|q|)^{5-      c (\gamma)  \cdot c (\delta)  \cdot c(|I|) \cdot E ( \lfloor \frac{|I|}{2} \rfloor  + 5) \cdot \eps } \cdot (1+|q|)^{2+ 2 \gamma - 4\de }}    \cdot w(q) \\
           \notag
 &\leq &     C(q_0)   \cdot  c (\delta) \cdot c (\gamma) \cdot C(|I|) \cdot E ( \lfloor \frac{|I|}{2} \rfloor + 5)  \cdot  \frac{\eps^3 }{(1+t+|q|)^{3 -4\de -   c (\gamma)  \cdot c (\delta)  \cdot c(|I|) \cdot E ( \lfloor \frac{|I|}{2} \rfloor  + 5) \cdot \eps } }   \; . \\
 \eea

For $0 < \de \leq \frac{1}{4}$\;, we have

      \bea
                 \notag
           &&         \int_{t_1 }^{t}   \int_{\Sigma^{ext}_{\tau} }    C(q_0)   \cdot  c (\delta) \cdot c (\gamma) \cdot C(|I|) \cdot E ( \lfloor \frac{|I|}{2} \rfloor + 5)  \\
                      \notag
           && \times \frac{\eps^3 }{(1+\tau+|q|)^{5-      c (\gamma)  \cdot c (\delta)  \cdot c(|I|) \cdot E ( \lfloor \frac{|I|}{2} \rfloor  + 5) \cdot \eps } \cdot (1+|q|)^{2+ 2 \gamma - 4\de }}    \cdot w(q)   \cdot d\tau  \\
                      \notag
     &\leq&                \int_{t_1 }^{t}      C(q_0)   \cdot  c (\delta) \cdot c (\gamma) \cdot C(|I|) \cdot E ( \lfloor \frac{|I|}{2} \rfloor + 5)  \cdot  \frac{\eps^3 }{(1+t+|q|)^{3 -4\de -   c (\gamma)  \cdot c (\delta)  \cdot c(|I|) \cdot E ( \lfloor \frac{|I|}{2} \rfloor  + 5) \cdot \eps } }  \cdot d\tau \\
     \notag
         &\leq &     C(q_0)   \cdot  c (\delta) \cdot c (\gamma) \cdot C(|I|) \cdot E ( \lfloor \frac{|I|}{2} \rfloor + 5)  \cdot  \frac{\eps^3 }{(1+t+|q|)^{1 -   c (\gamma)  \cdot c (\delta)  \cdot c(|I|) \cdot E ( \lfloor \frac{|I|}{2} \rfloor  + 5) \cdot \eps } }   \; . \\
           \eea

\end{lemma}

\begin{proof}

We estimate
\beaa
           &&       \int_{\Sigma^{ext}_{\tau} }    C(q_0)   \cdot  c (\delta) \cdot c (\gamma) \cdot C(|I|) \cdot E ( \lfloor \frac{|I|}{2} \rfloor + 5)  \\
           && \times \frac{\eps^3 }{(1+t+|q|)^{5-      c (\gamma)  \cdot c (\delta)  \cdot c(|I|) \cdot E ( \lfloor \frac{|I|}{2} \rfloor  + 5) \cdot \eps } \cdot (1+|q|)^{2+ 2 \gamma - 4\de }}    \cdot w(q)     \\
           &\leq &   \int_{r=0}^{\infty}   \int_{\th=0}^{\pi}  \int_{\phi=0}^{2\pi}    C(q_0)   \cdot  c (\delta) \cdot c (\gamma) \cdot C(|I|) \cdot E ( \lfloor \frac{|I|}{2} \rfloor + 5)  \\
           && \times \frac{\eps^3 }{(1+t+|q|)^{5-      c (\gamma)  \cdot c (\delta)  \cdot c(|I|) \cdot E ( \lfloor \frac{|I|}{2} \rfloor  + 5) \cdot \eps } \cdot (1+|q|)^{2+ 2 \gamma - 4\de }}    \cdot w(q)   \cdot r^2 \sin( \th )dr d\th d\phi \; .
\eeaa

We have
\beaa
 w(q)   \leq (1+|q|)^{1+2\ga} \; ,
\eeaa
and
\beaa
r^2 \leq (1+t+|q|)^2 \;.
\eeaa
Thus,
\beaa
           &&       \int_{\Sigma^{ext}_{\tau} }    C(q_0)   \cdot  c (\delta) \cdot c (\gamma) \cdot C(|I|) \cdot E ( \lfloor \frac{|I|}{2} \rfloor + 5)  \\
           && \times \frac{\eps^3 }{(1+t+|q|)^{5-      c (\gamma)  \cdot c (\delta)  \cdot c(|I|) \cdot E ( \lfloor \frac{|I|}{2} \rfloor  + 5) \cdot \eps } \cdot (1+|q|)^{2+ 2 \gamma - 4\de }}    \cdot w(q)     \\
           &\leq &   \int_{r=0}^{\infty}   \int_{\th=0}^{\pi}  \int_{\phi=0}^{2\pi}    C(q_0)   \cdot  c (\delta) \cdot c (\gamma) \cdot C(|I|) \cdot E ( \lfloor \frac{|I|}{2} \rfloor + 5)  \\
           && \times \frac{\eps^3 }{(1+t+|q|)^{3-      c (\gamma)  \cdot c (\delta)  \cdot c(|I|) \cdot E ( \lfloor \frac{|I|}{2} \rfloor  + 5) \cdot \eps } \cdot (1+|q|)^{1- 4\de }}    \cdot dr d\th d\phi \\
              &\leq &   \int_{r=0}^{\infty}   \int_{\th=0}^{\pi}  \int_{\phi=0}^{2\pi}   C(q_0)   \cdot  c (\delta) \cdot c (\gamma) \cdot C(|I|) \cdot E ( \lfloor \frac{|I|}{2} \rfloor + 5) \\
           && \times \frac{\eps^3 }{(1+t+|q|)^{3 -4\de -   c (\gamma)  \cdot c (\delta)  \cdot c(|I|) \cdot E ( \lfloor \frac{|I|}{2} \rfloor  + 5) \cdot \eps } \cdot (1+|q|)^{1+  c (\gamma)  \cdot c (\delta)  \cdot c(|I|) \cdot E ( \lfloor \frac{|I|}{2} \rfloor  + 5) \cdot \eps }}    \cdot dr d\th d\phi \\
  &\leq &     C(q_0)   \cdot  c (\delta) \cdot c (\gamma) \cdot C(|I|) \cdot E ( \lfloor \frac{|I|}{2} \rfloor + 5)  \cdot  \frac{\eps^3 }{(1+t+|q|)^{3 -4\de -   c (\gamma)  \cdot c (\delta)  \cdot c(|I|) \cdot E ( \lfloor \frac{|I|}{2} \rfloor  + 5) \cdot \eps } }   
           \eeaa
Therefore,
\beaa
           &&         \int_{t_1 }^{t}   \int_{\Sigma^{ext}_{\tau} }    C(q_0)   \cdot  c (\delta) \cdot c (\gamma) \cdot C(|I|) \cdot E ( \lfloor \frac{|I|}{2} \rfloor + 5)  \\
           && \times \frac{\eps^3 }{(1+\tau+|q|)^{5-      c (\gamma)  \cdot c (\delta)  \cdot c(|I|) \cdot E ( \lfloor \frac{|I|}{2} \rfloor  + 5) \cdot \eps } \cdot (1+|q|)^{2+ 2 \gamma - 4\de }}    \cdot w(q)   \cdot d\tau  \\
         &\leq &     C(q_0)   \cdot  c (\delta) \cdot c (\gamma) \cdot C(|I|) \cdot E ( \lfloor \frac{|I|}{2} \rfloor + 5)  \cdot  \frac{\eps^3 }{(1+t+|q|)^{2 -4\de -   c (\gamma)  \cdot c (\delta)  \cdot c(|I|) \cdot E ( \lfloor \frac{|I|}{2} \rfloor  + 5) \cdot \eps } }   \; .
           \eeaa
           Consequently, for $0 < \de \leq \frac{1}{4}$\, we get
      \beaa
           &&         \int_{t_1 }^{t}   \int_{\Sigma^{ext}_{\tau} }    C(q_0)   \cdot  c (\delta) \cdot c (\gamma) \cdot C(|I|) \cdot E ( \lfloor \frac{|I|}{2} \rfloor + 5)  \\
           && \times \frac{\eps^3 }{(1+\tau+|q|)^{3-      c (\gamma)  \cdot c (\delta)  \cdot c(|I|) \cdot E ( \lfloor \frac{|I|}{2} \rfloor  + 5) \cdot \eps } \cdot (1+|q|)^{2+ 2 \gamma - 4\de }}    \cdot w(q)   \cdot d\tau  \\
         &\leq &     C(q_0)   \cdot  c (\delta) \cdot c (\gamma) \cdot C(|I|) \cdot E ( \lfloor \frac{|I|}{2} \rfloor + 5)  \cdot  \frac{\eps^3 }{(1+t+|q|)^{1 -   c (\gamma)  \cdot c (\delta)  \cdot c(|I|) \cdot E ( \lfloor \frac{|I|}{2} \rfloor  + 5) \cdot \eps } }   \; .\\
                    \eeaa
This concludes the proof of the lemma.\\

\end{proof}

\begin{corollary}\label{HardytypeinequalityaplliedtothesoircetermsoftheEinsteinYangMillspoentialofgoodcomponents}
For $\ga \geq 3 \de $\,, for $0 < \de \leq \frac{1}{4}$\,, and for $\eps$ small enough, depending on $\ga$\,, $\de$ and on $|I|$, we have

                                                  \beaa
   \notag
&&  \int_{\Sigma^{ext}_{\tau} }   \frac{(1+ t )}{\eps} \cdot  |  \Lie_{Z^I}  g^{\la\mu} \derm_{\la}   \derm_{\mu}   A_{\cal T}  |^2    \cdot w(q)   \\
&\les&     \sum_{|K| \leq |I |}  \Big[   \;     \int_{\Sigma^{ext}_{\tau} }   \Big(      O \big(   C(q_0)   \cdot  c (\delta) \cdot c (\gamma) \cdot C(|I|) \cdot E (\lfloor \frac{|I|}{2} \rfloor + 5)  \cdot \frac{\eps   \cdot  | \derm ( \Lie_{Z^K} h^1 ) |^2  }{(1+t+|q|)}   \big) \\
   && + O \big(    C(q_0)   \cdot  c (\delta) \cdot c (\gamma) \cdot C(|I|) \cdot E ( \lfloor \frac{|I|}{2} \rfloor+5)  \cdot \frac{\eps  \cdot  |\derm ( \Lie_{Z^K} A ) |^2  }{(1+t+|q|) }  \big) \\
   &&  +  C(q_0)   \cdot  c (\delta) \cdot c (\gamma) \cdot C(|I|) \cdot E (\lfloor \frac{|I|}{2} \rfloor+ 4)  \cdot \frac{\eps   \cdot    | \rderm ( \Lie_{Z^K} A )  |^2  }{(1+t+|q|)^{1-      c (\gamma)  \cdot c (\delta)  \cdot c(|I|) \cdot E ( \lfloor \frac{|I|}{2} \rfloor+ 4) \cdot \eps } \cdot   (1+|q|) }    \\
      \notag
        &&    +     C(q_0)   \cdot  c (\delta) \cdot c (\gamma) \cdot C(|I|) \cdot E (\lfloor \frac{|I|}{2} \rfloor + 4)  \cdot \frac{\eps  \cdot   | \rderm  ( \Lie_{Z^K} h^1 ) |^2  }{(1+t+|q|)^{1-      c (\gamma)  \cdot c (\delta)  \cdot c(|I|) \cdot E (\lfloor \frac{|I|}{2} \rfloor+ 4) \cdot \eps } \cdot (1+|q|)^{2 }} \Big)    \cdot w(q)   \;  \Big]     \\
       &&    +      C(q_0)   \cdot  c (\delta) \cdot c (\gamma) \cdot C(|I|) \cdot E ( \lfloor \frac{|I|}{2} \rfloor + 5)  \cdot  \frac{\eps^3 }{(1+t+|q|)^{3 -4\de -   c (\gamma)  \cdot c (\delta)  \cdot c(|I|) \cdot E ( \lfloor \frac{|I|}{2} \rfloor  + 5) \cdot \eps } }   \; .
  \eeaa      
        
\end{corollary}

\begin{proof}
Integrating and applying the Hardy type inequality of Corollary \ref{HardytypeinequalityforintegralstartingatROm}, with $a = 1$\,, to Lemma \ref{estimateonttimesthesquareofthesourcesofwaveoperatorongoodcomponenentsofEinsteeinYangMillspoentialA}, we get the result for the terms that contain $A$ and $h^1$. Regarding the term generated from $h^0$\;, we use Lemma \ref{estimateonthecontributionofhzerointhesourcetermsforthewaveequationontheYangMillspoential} and we get the desired result. \\            
\end{proof}

\begin{lemma}\label{estimateonboththesourcetermsforgoodcomponentofAandthecommutatorterminclduingthetfactorstimesquareandthetimeintegral}
For $\ga \geq 3 \de $\,, for $0 < \de \leq \frac{1}{4}$\,, and for $\eps$ small enough, depending on $\ga$\,, $\de$ and on $|I|$, we have

                                                 \beaa
   \notag
&&  \int_{\Sigma^{ext}_{\tau} } \Big(  \frac{(1+ t )}{\eps} \cdot  |  \Lie_{Z^I}  g^{\la\mu} \derm_{\la}   \derm_{\mu}   A_{\cal T}  |^2  \\
&& +  \frac{(1+ t )}{\eps} \cdot  | \Lie_{Z^I}  ( g^{\la\mu} \derm_{\la}   \derm_{\mu}     A_{\cal T} ) - g^{\la\mu}    \derm_{\la}   \derm_{\mu}  (  \Lie_{Z^I} A_{\cal T} ) |^2 \Big)  \cdot w(q)   \\
&\les&      \sum_{|K| \leq |I |}  \Big[   \;     \int_{\Sigma^{ext}_{\tau} }     \Big[   \;       O \big(   C(q_0)   \cdot  c (\delta) \cdot c (\gamma) \cdot C(|I|) \cdot E (\lfloor \frac{|I|}{2} \rfloor + 5)  \cdot \frac{\eps   \cdot  | \derm ( \Lie_{Z^K} h^1 ) |^2  }{(1+t+|q|)}   \big) \\
   && + O \big(    C(q_0)   \cdot  c (\delta) \cdot c (\gamma) \cdot C(|I|) \cdot E (  \lfloor \frac{|I|}{2} \rfloor +5)  \cdot \frac{\eps  \cdot  |\derm ( \Lie_{Z^K} A ) |^2  }{(1+t+|q|) }  \big) \\
   &&  +  C(q_0)   \cdot  c (\delta) \cdot c (\gamma) \cdot C(|I|) \cdot E ( \lfloor \frac{|I|}{2} \rfloor+ 4)  \cdot \frac{\eps   \cdot    | \rderm ( \Lie_{Z^K} A )  |^2  }{(1+t+|q|)^{1-      c (\gamma)  \cdot c (\delta)  \cdot c(|I|) \cdot E ( \lfloor \frac{|I|}{2} \rfloor+ 4) \cdot \eps } \cdot   (1+|q|) }    \\
      \notag
        &&    +     C(q_0)   \cdot  c (\delta) \cdot c (\gamma) \cdot C(|I|) \cdot E (\lfloor \frac{|I|}{2} \rfloor + 4)  \cdot \frac{\eps  \cdot   | \rderm  ( \Lie_{Z^K} h^1 ) |^2  }{(1+t+|q|)^{1-      c (\gamma)  \cdot c (\delta)  \cdot c(|I|) \cdot E (\lfloor \frac{|I|}{2} \rfloor+ 4) \cdot \eps } \cdot (1+|q|)^{2 }}     \Big] \cdot w(q)  \; \Big]    \\
       &&  +      \int_{\Sigma^{ext}_{\tau} }   \Big[     C(q_0)   \cdot  c (\delta) \cdot c (\gamma) \cdot C(|I|) \cdot E (   \lfloor \frac{|I|}{2} \rfloor   +3)  \cdot \frac{\eps }{(1+t+|q|)^{1-   c (\gamma)  \cdot c (\delta)  \cdot c(|I|) \cdot E (\lfloor \frac{|I|}{2} \rfloor  + 2)\cdot  \eps } \cdot (1+|q|)^2} \\
&& \times    \sum_{  |K| \leq |I| -1 }| \derm ( \Lie_{Z^K}  A )  |^2    \Big] \cdot w(q) \\
&& +   \int_{\Sigma^{ext}_{\tau} }   \sum_{  |K| \leq |I| }       \Big[   C(q_0)   \cdot C(|I|) \cdot E (  \lfloor \frac{|I|}{2} \rfloor  +3)  \cdot \frac{\eps \cdot  | \Lie_{Z^{K}} H_{L  L} |^2 }{(1+t+|q|)^{1- 2 \de } \cdot (1+|q|)^{4+2\gamma}}  \Big] \cdot w(q) \\
      &&    +      C(q_0)   \cdot  c (\delta) \cdot c (\gamma) \cdot C(|I|) \cdot E ( \lfloor \frac{|I|}{2} \rfloor + 5)  \cdot  \frac{\eps^3 }{(1+t+|q|)^{3 -4\de -   c (\gamma)  \cdot c (\delta)  \cdot c(|I|) \cdot E ( \lfloor \frac{|I|}{2} \rfloor  + 5) \cdot \eps } }    \; .
\eeaa

\end{lemma}

\begin{proof}

Finally, putting the result Corollary \ref{HardytypeinequalityaplliedtothesoircetermsoftheEinsteinYangMillspoentialofgoodcomponents} and of Lemma \ref{Hardytypeineeuqlityappliedtothecommutatortermplusatimeintegralformula}  together, we obtain
                                                    \beaa
   \notag
&&  \int_{\Sigma^{ext}_{\tau} } \Big(  \frac{(1+ t )}{\eps} \cdot  |  \Lie_{Z^I}  g^{\la\mu} \derm_{\la}   \derm_{\mu}   A_{\cal T}  |^2  \\
&& +  \frac{(1+ t )}{\eps} \cdot  | \Lie_{Z^I}  ( g^{\la\mu} \derm_{\la}   \derm_{\mu}     A_{\cal T} ) - g^{\la\mu}    \derm_{\la}   \derm_{\mu}  (  \Lie_{Z^K} A_{\cal T} ) |^2 \Big)  \cdot w(q)   \\
&\les&    \sum_{|K| \leq |I |}  \Big[   \;     \int_{\Sigma^{ext}_{\tau} }   \Big(      O \big(   C(q_0)   \cdot  c (\delta) \cdot c (\gamma) \cdot C(|I|) \cdot E (  \lfloor \frac{|I|}{2} \rfloor + 5)  \cdot \frac{\eps   \cdot  | \derm ( \Lie_{Z^K} h^1 ) |^2  }{(1+t+|q|)}   \big) \\
   && + O \big(    C(q_0)   \cdot  c (\delta) \cdot c (\gamma) \cdot C(|I|) \cdot E (   \lfloor \frac{|I|}{2} \rfloor  +5)  \cdot \frac{\eps  \cdot  |\derm ( \Lie_{Z^K} A ) |^2  }{(1+t+|q|) }  \big) \\
   &&  +  C(q_0)   \cdot  c (\delta) \cdot c (\gamma) \cdot C(|I|) \cdot E (  \lfloor \frac{|I|}{2} \rfloor + 4)  \cdot \frac{\eps   \cdot    | \rderm ( \Lie_{Z^K} A )  |^2  }{(1+t+|q|)^{1-      c (\gamma)  \cdot c (\delta)  \cdot c(|I|) \cdot E (  \lfloor \frac{|I|}{2} \rfloor+ 4) \cdot \eps } \cdot   (1+|q|) }    \\
      \notag
        &&    +     C(q_0)   \cdot  c (\delta) \cdot c (\gamma) \cdot C(|I|) \cdot E (  \lfloor \frac{|I|}{2} \rfloor + 4)  \cdot \frac{\eps  \cdot   | \rderm  ( \Lie_{Z^K} h^1 ) |^2  }{(1+t+|q|)^{1-      c (\gamma)  \cdot c (\delta)  \cdot c(|I|) \cdot E (  \lfloor \frac{|I|}{2} \rfloor+ 4) \cdot \eps } \cdot (1+|q|)^{2 }} \Big)     \;  \Big]    \cdot w(q)    \\
       &&  +      \int_{\Sigma^{ext}_{\tau} }   \Big(       \frac{(1+ t )}{\eps} \cdot    \sum_{|K| \leq |I|  - 1}  | m^{\la\mu} \cdot \derm_{\la}   \derm_{\mu} (  \Lie_{Z^{K}}  A_{e_{a}} ) |^2  \cdot w(q)  \\
&& +  C(q_0)   \cdot  c (\delta) \cdot c (\gamma) \cdot C( |I|) \cdot E (  |I| +5 )  \\
&& \times  \sum_{ |K| \leq |I| }  \Big(   \frac{\eps   }{(1+t+|q|)}  \, \cdot  |  \derm (   \Lie_{Z^K}  h )   |^2  +  \frac{\eps }{(1+t+|q|) }  \, \cdot | \derm ( \Lie_{Z^K}  A)  |^2 \Big)   \cdot w(q)  \\
&&+  C(q_0)   \cdot  c (\delta) \cdot c (\gamma) \cdot C(|I|) \cdot E (    \lfloor \frac{|I|}{2} \rfloor   +3)  \cdot \frac{\eps }{(1+t+|q|)^{1-   c (\gamma)  \cdot c (\delta)  \cdot c(|I|) \cdot E ( \lfloor \frac{|I|}{2} \rfloor  + 2)\cdot  \eps } \cdot (1+|q|)^2} \\
&& \times    \sum_{  |K| \leq |I| -1 }  | \derm ( \Lie_{Z^K}  A  )  |^2    \cdot w(q)  \Big)\\
&& +   \int_{\Sigma^{ext}_{\tau} }   \sum_{  |K| \leq |I| }       \Big[   C(q_0)   \cdot C(|I|) \cdot E (  \lfloor \frac{|I|}{2} \rfloor  +3)  \cdot \frac{\eps \cdot  | \Lie_{Z^{K}} H_{L  L} |^2 }{(1+t+|q|)^{1- 2 \de } \cdot (1+|q|)^{4+2\gamma}}  \Big] \cdot w(q) \\
 &&    +      C(q_0)   \cdot  c (\delta) \cdot c (\gamma) \cdot C(|I|) \cdot E ( \lfloor \frac{|I|}{2} \rfloor + 5)  \cdot  \frac{\eps^3 }{(1+t+|q|)^{3 -4\de -   c (\gamma)  \cdot c (\delta)  \cdot c(|I|) \cdot E ( \lfloor \frac{|I|}{2} \rfloor  + 5) \cdot \eps } }   \; . 
\eeaa

We get

                                                  \beaa
   \notag
&&  \int_{\Sigma^{ext}_{\tau} } \Big(  \frac{(1+ t )}{\eps} \cdot  |  \Lie_{Z^I}  g^{\la\mu} \derm_{\la}   \derm_{\mu}   A_{\cal T}  |^2  \\
&& +  \frac{(1+ t )}{\eps} \cdot  | \Lie_{Z^I}  ( g^{\la\mu} \derm_{\la}   \derm_{\mu}     A_{\cal T} ) - g^{\la\mu}    \derm_{\la}   \derm_{\mu}  (  \Lie_{Z^K} A_{\cal T} ) |^2 \Big)  \cdot w(q)   \\
&\les&     \sum_{|K| \leq |I |}  \Big[   \;        \int_{\Sigma^{ext}_{\tau} }   \Big(      O \big(   C(q_0)   \cdot  c (\delta) \cdot c (\gamma) \cdot C(|I|) \cdot E (\lfloor \frac{|I|}{2} \rfloor+ 5)  \cdot \frac{\eps   \cdot  | \derm ( \Lie_{Z^K} h^1 ) |^2  }{(1+t+|q|)}   \big) \\
   && + O \big(    C(q_0)   \cdot  c (\delta) \cdot c (\gamma) \cdot C(|I|) \cdot E (  \lfloor \frac{|I|}{2} \rfloor  +5)  \cdot \frac{\eps  \cdot  |\derm ( \Lie_{Z^K} A ) |^2  }{(1+t+|q|) }  \big) \\
   &&  +  C(q_0)   \cdot  c (\delta) \cdot c (\gamma) \cdot C(|I|) \cdot E ( \lfloor \frac{|I|}{2} \rfloor + 4)  \cdot \frac{\eps   \cdot    | \rderm ( \Lie_{Z^K} A )  |^2  }{(1+t+|q|)^{1-      c (\gamma)  \cdot c (\delta)  \cdot c(|I|) \cdot E ( \lfloor \frac{|I|}{2} \rfloor  + 4) \cdot \eps } \cdot   (1+|q|) }    \\
      \notag
        &&    +     C(q_0)   \cdot  c (\delta) \cdot c (\gamma) \cdot C(|I|) \cdot E ( \lfloor \frac{|I|}{2} \rfloor+ 4)  \cdot \frac{\eps  \cdot   | \rderm  ( \Lie_{Z^K} h^1 ) |^2  }{(1+t+|q|)^{1-      c (\gamma)  \cdot c (\delta)  \cdot c(|I|) \cdot E ( \lfloor \frac{|I|}{2} \rfloor+ 4) \cdot \eps } \cdot (1+|q|)^{2 }} \Big)   \;  \Big]   \cdot w(q)    \\
       &&  +      \int_{\Sigma^{ext}_{\tau} }   \Big(       \frac{(1+ t )}{\eps} \cdot    \sum_{|K| \leq |I|  - 1}  | g^{\la\mu} \cdot \derm_{\la}   \derm_{\mu} (  \Lie_{Z^{K}}  A_{e_{a}} ) |^2 \\
&&+  C(q_0)   \cdot  c (\delta) \cdot c (\gamma) \cdot C(|I|) \cdot E (    \lfloor \frac{|I|}{2} \rfloor   +3)  \cdot \frac{\eps }{(1+t+|q|)^{1-   c (\gamma)  \cdot c (\delta)  \cdot c(|I|) \cdot E (\lfloor \frac{|I|}{2} \rfloor  + 2)\cdot  \eps } \cdot (1+|q|)^2} \\
&& \times    \sum_{  |K| \leq |I| -1 }  | \derm ( \Lie_{Z^K}  A )  |^2   \Big) \cdot w(q) \\
&& +   \int_{\Sigma^{ext}_{\tau} }   \sum_{  |K| \leq |I| }       \Big[   C(q_0)   \cdot C(|I|) \cdot E (  \lfloor \frac{|I|}{2} \rfloor  +3)  \cdot \frac{\eps \cdot  | \Lie_{Z^{K}} H_{L  L} |^2 }{(1+t+|q|)^{1- 2 \de } \cdot (1+|q|)^{4+2\gamma}}  \Big] \cdot w(q) \\
      &&    +      C(q_0)   \cdot  c (\delta) \cdot c (\gamma) \cdot C(|I|) \cdot E ( \lfloor \frac{|I|}{2} \rfloor + 5)  \cdot  \frac{\eps^3 }{(1+t+|q|)^{3 -4\de -   c (\gamma)  \cdot c (\delta)  \cdot c(|I|) \cdot E ( \lfloor \frac{|I|}{2} \rfloor  + 5) \cdot \eps } }   \; .
\eeaa

But the term $    \frac{(1+ t )}{\eps} \cdot    \sum_{|K| \leq |I|  - 1}  | g^{\la\mu} \cdot \derm_{\la}   \derm_{\mu} (  \Lie_{Z^{K}}  A_{e_{a}} ) |^2 $ has a good structure that is similar, even the same, as for $    \frac{(1+ t )}{\eps} \cdot    \sum_{|K| \leq |I| }  | g^{\la\mu} \cdot \derm_{\la}   \derm_{\mu} (  \Lie_{Z^{K}}  A_{e_{a}} ) |^2 $ that we already treated (except that we have in it $ |I|  - 1$ Lie derivatives instead of $|I|$ Lie derivatives).

Thus,
                                                  \beaa
   \notag
&&  \int_{\Sigma^{ext}_{\tau} } \Big(  \frac{(1+ t )}{\eps} \cdot  |  \Lie_{Z^I}  g^{\la\mu} \derm_{\la}   \derm_{\mu}   A_{\cal T}  |^2  \\
&& +  \frac{(1+ t )}{\eps} \cdot  | \Lie_{Z^I}  ( g^{\la\mu} \derm_{\la}   \derm_{\mu}     A_{\cal T} ) - g^{\la\mu}    \derm_{\la}   \derm_{\mu}  (  \Lie_{Z^I} A_{\cal T}  ) |^2 \Big)  \cdot w(q)   \\
&\les&      \sum_{|K| \leq |I |}  \Big[   \;     \int_{\Sigma^{ext}_{\tau} }     \Big[   \;       O \big(   C(q_0)   \cdot  c (\delta) \cdot c (\gamma) \cdot C(|I|) \cdot E (\lfloor \frac{|I|}{2} \rfloor + 5)  \cdot \frac{\eps   \cdot  | \derm ( \Lie_{Z^K} h^1 ) |^2  }{(1+t+|q|)}   \big) \\
   && + O \big(    C(q_0)   \cdot  c (\delta) \cdot c (\gamma) \cdot C(|I|) \cdot E (  \lfloor \frac{|I|}{2} \rfloor +5)  \cdot \frac{\eps  \cdot  |\derm ( \Lie_{Z^K} A ) |^2  }{(1+t+|q|) }  \big) \\
   &&  +  C(q_0)   \cdot  c (\delta) \cdot c (\gamma) \cdot C(|I|) \cdot E ( \lfloor \frac{|I|}{2} \rfloor+ 4)  \cdot \frac{\eps   \cdot    | \rderm ( \Lie_{Z^K} A )  |^2  }{(1+t+|q|)^{1-      c (\gamma)  \cdot c (\delta)  \cdot c(|I|) \cdot E ( \lfloor \frac{|I|}{2} \rfloor+ 4) \cdot \eps } \cdot   (1+|q|) }    \\
      \notag
        &&    +     C(q_0)   \cdot  c (\delta) \cdot c (\gamma) \cdot C(|I|) \cdot E (\lfloor \frac{|I|}{2} \rfloor + 4)  \cdot \frac{\eps  \cdot   | \rderm  ( \Lie_{Z^K} h^1 ) |^2  }{(1+t+|q|)^{1-      c (\gamma)  \cdot c (\delta)  \cdot c(|I|) \cdot E (\lfloor \frac{|I|}{2} \rfloor+ 4) \cdot \eps } \cdot (1+|q|)^{2 }}     \Big] \cdot w(q)  \; \Big]    \\
       &&  +      \int_{\Sigma^{ext}_{\tau} }   \Big[     C(q_0)   \cdot  c (\delta) \cdot c (\gamma) \cdot C(|I|) \cdot E (   \lfloor \frac{|I|}{2} \rfloor   +3)  \cdot \frac{\eps }{(1+t+|q|)^{1-   c (\gamma)  \cdot c (\delta)  \cdot c(|I|) \cdot E ( \lfloor \frac{|I|}{2} \rfloor  + 2)\cdot  \eps } \cdot (1+|q|)^2} \\
&& \times    \sum_{  |K| \leq |I| -1 }   | \derm ( \Lie_{Z^K}  A )  |^2      \Big] \cdot w(q)  \\
&& +   \int_{\Sigma^{ext}_{\tau} }   \sum_{  |K| \leq |I| }       \Big[   C(q_0)   \cdot C(|I|) \cdot E (  \lfloor \frac{|I|}{2} \rfloor  +3)  \cdot \frac{\eps \cdot  | \Lie_{Z^{K}} H_{L  L} |^2 }{(1+t+|q|)^{1- 2 \de } \cdot (1+|q|)^{4+2\gamma}}  \Big] \cdot w(q) \\
      &&    +      C(q_0)   \cdot  c (\delta) \cdot c (\gamma) \cdot C(|I|) \cdot E ( \lfloor \frac{|I|}{2} \rfloor + 5)  \cdot  \frac{\eps^3 }{(1+t+|q|)^{3 -4\de -   c (\gamma)  \cdot c (\delta)  \cdot c(|I|) \cdot E ( \lfloor \frac{|I|}{2} \rfloor  + 5) \cdot \eps } }    \; .
\eeaa

\end{proof}

\subsection{Control on the weighted $L^2$ norm of $\derm ( \Lie_{Z^I}A_{e_{a}} ) $}\

\begin{lemma}\label{ControlontheweightedLtwonormofdermAea}
For $\ga \geq 3 \de $\,, for $0 < \de \leq \frac{1}{4}$\,, and for $\eps$ small, enough depending on $q_0$\;, on $\ga$\;, on $\de$\;, on $|I|$ and on $\mu$\;, we have

      \beaa
   \notag
 &&     \int_{\Sigma^{ext}_{t_2} }  |\derm ( \Lie_{Z^I}  A_{e_{a}} )  |^2     \cdot w(q)  \cdot d^{3}x \\  
 &\les&      \int_{\Sigma^{ext}_{t_1} }  |\derm  ( \Lie_{Z^I}  A_{e_{a}} )|^2     \cdot w(q)  \cdot d^{3}x \\
&& +   \sum_{|K| \leq |I |}   \int_{t_1}^{t_2}  \Big[    \;   \int_{\Sigma^{ext}_{\tau} }    \Big[       O \big(   C(q_0)   \cdot  c (\delta) \cdot c (\gamma) \cdot C(|I|) \cdot E ( \lfloor \frac{|I|}{2} \rfloor + 5)  \cdot \frac{\eps   \cdot  | \derm ( \Lie_{Z^K} h^1 ) |^2  }{(1+\tau+|q|)}   \big) \\
   && + O \big(    C(q_0)   \cdot  c (\delta) \cdot c (\gamma) \cdot C(|I|) \cdot E (   \lfloor \frac{|I|}{2} \rfloor +5)  \cdot \frac{\eps  \cdot  |\derm ( \Lie_{Z^K} A ) |^2  }{(1+\tau+|q|) }  \big)    \; \Big]    \cdot w(q)  \cdot d^{3}x \Big]   \cdot d\tau  \\
      \notag
        &&  +     \int_{t_1}^{t_2}  \int_{\Sigma^{ext}_{\tau} }   \Big[     C(q_0)   \cdot  c (\delta) \cdot c (\gamma) \cdot C(|I|) \cdot E (  \lfloor \frac{|I|}{2} \rfloor   +3)  \cdot \frac{\eps }{(1+\tau+|q|)^{1-   c (\gamma)  \cdot c (\delta)  \cdot c(|I|) \cdot E (\lfloor \frac{|I|}{2} \rfloor  + 2)\cdot  \eps } \cdot (1+|q|)^2} \\
&& \times    \sum_{  |K| \leq |I| -1 }  | \derm ( \Lie_{Z^K}  A )  |^2    \Big] \cdot w(q) \cdot d^{3}x  \cdot d\tau \\
&& +   \int_{t_1}^{t_2}   \int_{\Sigma^{ext}_{\tau} }   \sum_{  |K| \leq |I| }        \Big[   C(q_0)   \cdot C(|I|) \cdot E (  \lfloor \frac{|I|}{2} \rfloor  +3)  \cdot \frac{\eps \cdot  | \Lie_{Z^{K}} H_{L  L} |^2 }{(1+\tau+|q|)^{1- 2 \de } \cdot (1+|q|)^{4+2\gamma}}  \Big] \cdot w(q)  \cdot d^{3}x  \cdot d\tau \\
 && +   C(q_0)   \cdot  c (\delta) \cdot c (\gamma) \cdot C(|I|) \cdot E ( \lfloor \frac{|I|}{2} \rfloor + 5)  \cdot  \frac{\eps^3 }{(1+t+|q|)^{1 -   c (\gamma)  \cdot c (\delta)  \cdot c(|I|) \cdot E ( \lfloor \frac{|I|}{2} \rfloor  + 5) \cdot \eps } }  \; .
\eeaa

\end{lemma}

\begin{proof}

Based on Lemma \ref{Theveryfinal }, for $\eps$ small enough, depending on $q_0$, $\de$\,, $\ga$\;, and $\mu$\;, the following energy estimate holds for $\ga > 0$ and for $\Phi$ decaying sufficiently fast at spatial infinity,
  \bea
   \notag
 &&     \int_{\Sigma^{ext}_{t_2} }  |\derm \Phi_{e_{a}}  |^2     \cdot w(q)  \cdot d^{3}x    + \int_{N_{t_1}^{t_2} }  T_{\hat{L} t}^{(\bf{g})} (\Phi_{e_{a}} ) \cdot  w(q) \cdot dv^{(\bf{m})}_N \\
 \notag
 &&+ \int_{t_1}^{t_2}  \int_{\Sigma^{ext}_{\tau} }     | \rderm \Phi_{e_{a}}  |^2   \cdot  \frac{\widehat{w} (q)}{(1+|q|)} \cdot  d^{3}x  \cdot d\tau \\
  \notag
  &\les &       \int_{\Sigma^{ext}_{t_1} }  |\derm \Phi_{e_{a}}  |^2     \cdot w(q)  \cdot d^{3}x \\
    \notag
   &&  +  \int_{t_1}^{t_2}  \int_{\Sigma^{ext}_{\tau} }  \frac{(1+\tau )}{\eps} \cdot |  g^{\mu\a} \derm_{\mu } \derm_\a \Phi_{e_{a}}  |^2    \cdot w(q) \cdot  d^{3}x  \cdot d\tau \\
 \notag
&& + \int_{t_1}^{t_2}  \int_{\Sigma^{ext}_{\tau} }   C(q_0) \cdot   c (\delta) \cdot c (\gamma) \cdot E (  4)  \cdot \frac{ \eps }{ (1+\tau)  } \cdot | \derm \Phi_{e_{a}}  |^2     \cdot  w(q)\cdot  d^{3}x  \cdot d\tau \; . \\
 \eea

Performing the energy estimate on $  \Lie_{Z^I} A_{e_{a}}$\;, we get 

  \beaa
   \notag
 &&     \int_{\Sigma^{ext}_{t_2} }  |\derm ( \Lie_{Z^I}  A_{e_{a}} )  |^2     \cdot w(q)  \cdot d^{3}x    + \int_{N_{t_1}^{t_2} }  T_{\hat{L} t}^{(\bf{g})}  \cdot  w(q) \cdot dv^{(\bf{m})}_N \\
 \notag
 &&+ \int_{t_1}^{t_2}  \int_{\Sigma^{ext}_{\tau} }     | \rderm  ( \Lie_{Z^I}  A_{e_{a}} ) |^2   \cdot  \frac{\widehat{w} (q)}{(1+|q|)} \cdot  d^{3}x  \cdot d\tau \\
  \notag
  &\les &       \int_{\Sigma^{ext}_{t_1} }  |\derm  ( \Lie_{Z^I}  A_{e_{a}} )|^2     \cdot w(q)  \cdot d^{3}x \\
    \notag
   &&  +  \int_{t_1}^{t_2}  \int_{\Sigma^{ext}_{\tau} }  \frac{(1+\tau )}{\eps} \cdot |  g^{\mu\a} \derm_{\mu } \derm_\a  ( \Lie_{Z^I}  A_{e_{a}} )  |^2    \cdot w(q) \cdot  d^{3}x  \cdot d\tau \\
 \notag
&& + \int_{t_1}^{t_2}  \int_{\Sigma^{ext}_{\tau} }   C(q_0) \cdot   c (\delta) \cdot c (\gamma) \cdot E (  4)  \cdot \frac{ \eps }{ (1+\tau)  } \cdot | \derm  ( \Lie_{Z^I}  A_{e_{a}} ) |^2     \cdot  w(q)\cdot  d^{3}x  \cdot d\tau \; . 
 \eeaa 
 
 However, on one hand, we have
 \beaa
     \notag
   &&  |  g^{\mu\a} \derm_{\mu } \derm_\a  ( \Lie_{Z^I}  A_{e_{a}} )  |^2   \\
 \notag
 &\les& |  \Lie_{Z^I} \big( g^{\mu\a} \derm_{\mu } \derm_\a   A_{e_{a}} \big)  |^2   \\
 && +  | \Lie_{Z^I}  ( g^{\la\mu} \derm_{\la}   \derm_{\mu}     A_{e_{a}} ) - g^{\la\mu}    \derm_{\la}   \derm_{\mu}  (  \Lie_{Z^I} A_{e_{a}}  ) |^2 \; , 
 \eeaa
 
for which we have an estimate, using Lemma \ref{estimateonboththesourcetermsforgoodcomponentofAandthecommutatorterminclduingthetfactorstimesquareandthetimeintegral}, and Lemma \ref{estimateonthecontributionofhzerointhesourcetermsforthewaveequationontheYangMillspoential}, and we an inject that estimate in the energy estimate above.

Now, since for $\eps$ small, enough depending on $\ga$, on $\de$, on $\mu$ and on $|I|$ (recall that $E ( k) \leq 1$ so we do not have dependance on the energy -- see Remark \ref{remarkonwhyepsilonsmalldesonotdependontheconstantsEofkwhichareboundedby1}),  we have
\beaa
\frac{\eps     }{(1+t+|q|)^{1-      c (\gamma)  \cdot c (\delta)  \cdot c(|I|) \cdot E ( |I|+ 4) \cdot \eps } \cdot   (1+|q|) }  &\leq &  \frac{\widehat{w} (q)}{(1+|q|)} \; .
  \eeaa
  
 Thus, we can absorb the tangential derivatives into the left hand side of the inequality, and obtain
      \beaa
   \notag
 &&     \int_{\Sigma^{ext}_{t_2} }  |\derm ( \Lie_{Z^I}  A_{e_{a}} )  |^2     \cdot w(q)  \cdot d^{3}x \\
  \notag
  &\les &       \int_{\Sigma^{ext}_{t_1} }  |\derm  ( \Lie_{Z^I}  A_{e_{a}} )|^2     \cdot w(q)  \cdot d^{3}x \\
    \notag
&& +   \sum_{|K| \leq |I |}   \int_{t_1}^{t_2}  \Big[   \; \int_{\Sigma^{ext}_{\tau} }   C(q_0) \cdot   c (\delta) \cdot c (\gamma) \cdot E (  4)  \cdot \frac{ \eps }{ (1+\tau)  } \cdot | \derm  ( \Lie_{Z^K}  A_{e_{a}} ) |^2     \cdot  w(q)\cdot  d^{3}x  \\
&& +    \int_{\Sigma^{ext}_{\tau} }    \Big[       O \big(   C(q_0)   \cdot  c (\delta) \cdot c (\gamma) \cdot C(|I|) \cdot E ( \lfloor \frac{|I|}{2} \rfloor + 5)  \cdot \frac{\eps   \cdot  | \derm ( \Lie_{Z^K} h^1 ) |^2  }{(1+\tau+|q|)}   \big) \\
   && + O \big(    C(q_0)   \cdot  c (\delta) \cdot c (\gamma) \cdot C(|I|) \cdot E (  \lfloor \frac{|I|}{2} \rfloor  +5)  \cdot \frac{\eps  \cdot  |\derm ( \Lie_{Z^K} A ) |^2  }{(1+\tau+|q|) }  \big)  \; \Big]  \cdot w(q)  \cdot d^{3}x \Big]  \cdot d\tau  \;   \\
         &&  +      \int_{t_1}^{t_2}  \int_{\Sigma^{ext}_{\tau} }   \Big[     C(q_0)   \cdot  c (\delta) \cdot c (\gamma) \cdot C(|I|) \cdot E (   \lfloor \frac{|I|}{2} \rfloor   +3)  \cdot \frac{\eps }{(1+\tau+|q|)^{1-   c (\gamma)  \cdot c (\delta)  \cdot c(|I|) \cdot E (\lfloor \frac{|I|}{2} \rfloor  + 2)\cdot  \eps } \cdot (1+|q|)^2} \\
&& \times    \sum_{  |K| \leq |I| -1 }   | \derm ( \Lie_{Z^K}  A )  |^2     \Big]  \cdot w(q)  \cdot d^{3}x \cdot d\tau  \\
&& +   \int_{t_1}^{t_2}   \int_{\Sigma^{ext}_{\tau} }   \sum_{  |K| \leq |I| }        \Big[   C(q_0)   \cdot C(|I|) \cdot E (  \lfloor \frac{|I|}{2} \rfloor  +3)  \cdot \frac{\eps \cdot  | \Lie_{Z^{K}} H_{L  L} |^2 }{(1+\tau+|q|)^{1- 2 \de } \cdot (1+|q|)^{4+2\gamma}}  \Big] \cdot w(q)  \cdot d^{3}x  \cdot d\tau \\
  && +  C(q_0)   \cdot  c (\delta) \cdot c (\gamma) \cdot C(|I|) \cdot E ( \lfloor \frac{|I|}{2} \rfloor + 5)  \cdot  \frac{\eps^3 }{(1+t+|q|)^{1 -   c (\gamma)  \cdot c (\delta)  \cdot c(|I|) \cdot E ( \lfloor \frac{|I|}{2} \rfloor  + 5) \cdot \eps } }  \; .
\eeaa

We have  $ \lfloor \frac{|I|}{2} \rfloor  +5 \geq 4$ and therefore,

      \beaa
   \notag
 &&     \int_{\Sigma^{ext}_{t_2} }  |\derm ( \Lie_{Z^I}  A_{e_{a}} )  |^2     \cdot w(q)  \cdot d^{3}x \\  
 &\les&      \int_{\Sigma^{ext}_{t_1} }  |\derm  ( \Lie_{Z^I}  A_{e_{a}} )|^2     \cdot w(q)  \cdot d^{3}x \\
&& +   \sum_{|K| \leq |I |}   \int_{t_1}^{t_2}  \Big[    \;   \int_{\Sigma^{ext}_{\tau} }    \Big[       O \big(   C(q_0)   \cdot  c (\delta) \cdot c (\gamma) \cdot C(|I|) \cdot E ( \lfloor \frac{|I|}{2} \rfloor + 5)  \cdot \frac{\eps   \cdot  | \derm ( \Lie_{Z^K} h^1 ) |^2  }{(1+\tau+|q|)}   \big) \\
   && + O \big(    C(q_0)   \cdot  c (\delta) \cdot c (\gamma) \cdot C(|I|) \cdot E (   \lfloor \frac{|I|}{2} \rfloor +5)  \cdot \frac{\eps  \cdot  |\derm ( \Lie_{Z^K} A ) |^2  }{(1+\tau+|q|) }  \big)    \; \Big]    \cdot w(q)  \cdot d^{3}x \Big]   \cdot d\tau  \\
      \notag
        &&  +     \int_{t_1}^{t_2}  \int_{\Sigma^{ext}_{\tau} }   \Big[     C(q_0)   \cdot  c (\delta) \cdot c (\gamma) \cdot C(|I|) \cdot E (  \lfloor \frac{|I|}{2} \rfloor   +3)  \cdot \frac{\eps }{(1+\tau+|q|)^{1-   c (\gamma)  \cdot c (\delta)  \cdot c(|I|) \cdot E (\lfloor \frac{|I|}{2} \rfloor  + 2)\cdot  \eps } \cdot (1+|q|)^2} \\
&& \times    \sum_{  |K| \leq |I| -1 }  | \derm ( \Lie_{Z^K}  A )  |^2    \Big] \cdot w(q) \cdot d^{3}x  \cdot d\tau \\
&& +   \int_{t_1}^{t_2}   \int_{\Sigma^{ext}_{\tau} }   \sum_{  |K| \leq |I| }        \Big[   C(q_0)   \cdot C(|I|) \cdot E (  \lfloor \frac{|I|}{2} \rfloor  +3)  \cdot \frac{\eps \cdot  | \Lie_{Z^{K}} H_{L  L} |^2 }{(1+\tau+|q|)^{1- 2 \de } \cdot (1+|q|)^{4+2\gamma}}  \Big] \cdot w(q)  \cdot d^{3}x  \cdot d\tau \\
 && +   C(q_0)   \cdot  c (\delta) \cdot c (\gamma) \cdot C(|I|) \cdot E ( \lfloor \frac{|I|}{2} \rfloor + 5)  \cdot  \frac{\eps^3 }{(1+t+|q|)^{1 -   c (\gamma)  \cdot c (\delta)  \cdot c(|I|) \cdot E ( \lfloor \frac{|I|}{2} \rfloor  + 5) \cdot \eps } }  \; .
\eeaa

\end{proof}

\subsection{Improved pointwise decay on $\derm ( \Lie_{Z^I}A_{e_{a}} ) $}\

The goal so far, was to obtain a decay in $\frac{1}{t} $ while keeping a decay in $q$. LEt us point out that if we wanted to use the weighted estimate by Lindblad and Rodnianski (see \eqref{LinfinitynormestimateongradientderivedbyLondbladRodnianski}) with a real weight in $q$, this would imply a Grönwall inequality on the quantity that we want to control, with a factor of $\frac{1}{t} $, which would imply necessarily an $\eps$ loss in the decay, i.e. a decay as $\frac{1}{t^{1-\eps}} $\;, which we cannot afford in order to close our argument. Hence, we resorted here to a different method, which is to control the energy first for the good components, as we did, and then use Klainerman-Sobolev inequality again to inject all in the energy estimate for the full component. We shall explain.

Let us start first by using the Klainerman-Sobolev in the exterior to estimate this special component $|\derm ( \Lie_{Z^I}  A_{e_{a}} ) | $. 

\begin{lemma}\label{upgradedestimateontheA_acomponentsthatconservesthedecayasoneoverttimesintegralsoftheenergysoastousegeneralizedgronwall}

For $\ga \geq 3 \de $\,, for $0 < \de \leq \frac{1}{4}$\,, and for $\eps$ small, enough depending on $q_0$\,, on $\ga$, on $\de$, on $|I|$ and on $\mu < 0$\;, we have

      \beaa
   \notag
 &&   |\derm ( \Lie_{Z^I}  A_{e_{a}} )  |  \\
 &\les&   \frac{1}{(1+t+|q|) \cdot (1+|q|)^{1+\ga} } \cdot \Big[ \;  \sum_{|J|\leq |I| + 2 }   \int_{\Sigma^{ext}_{t_1} }  |\derm  ( \Lie_{Z^J}  A_{e_{a}} )|^2     \cdot w(q)  \cdot d^{3}x \\
&& +   \sum_{|K| \leq |I |}   \int_{t_1}^{t}  \Big[    \;   \int_{\Sigma^{ext}_{\tau} }    \Big[       O \big(   C(q_0)   \cdot  c (\delta) \cdot c (\gamma) \cdot C(|I|) \cdot E ( \lfloor \frac{|I|}{2} \rfloor + 6)  \cdot \frac{\eps   \cdot  | \derm ( \Lie_{Z^K} h^1 ) |^2  }{(1+\tau+|q|)}   \big) \\
   && + O \big(    C(q_0)   \cdot  c (\delta) \cdot c (\gamma) \cdot C(|I|) \cdot E (   \lfloor \frac{|I|}{2} \rfloor +6)  \cdot \frac{\eps  \cdot  |\derm ( \Lie_{Z^K} A ) |^2  }{(1+\tau+|q|) }  \big)    \; \Big]    \cdot w(q)  \cdot d^{3}x  \Big]   \cdot d\tau  \\
      \notag
        &&  +     \int_{t_1}^{t}  \int_{\Sigma^{ext}_{\tau} }   \Big[     C(q_0)   \cdot  c (\delta) \cdot c (\gamma) \cdot C(|I|) \cdot E ( \lfloor \frac{|I|}{2} \rfloor  +4)  \cdot \frac{\eps }{(1+\tau+|q|)^{1-   c (\gamma)  \cdot c (\delta)  \cdot c(|I|) \cdot E (  \lfloor \frac{|I|}{2} \rfloor  + 4)\cdot  \eps } \cdot (1+|q|)^2} \\
&& \times    \sum_{  |K| \leq |I| -1 }  | \derm ( \Lie_{Z^K}  A )  |^2     \Big]  \cdot w(q)  \cdot d^{3}x  \cdot d\tau \\
&& +   \int_{t_1}^{t}   \int_{\Sigma^{ext}_{\tau} }   \sum_{  |K| \leq |I| }       \Big[   C(q_0)   \cdot C(|I|) \cdot E (  \lfloor \frac{|I|}{2} \rfloor  +5)  \cdot \frac{\eps \cdot  | \Lie_{Z^{K}} H_{L  L} |^2 }{(1+\tau+|q|)^{1- 2 \de } \cdot (1+|q|)^{4+2\gamma}}  \Big] \cdot w(q)  \cdot d^{3}x  \cdot d\tau   \\
 && +   C(q_0)   \cdot  c (\delta) \cdot c (\gamma) \cdot C(|I|) \cdot E ( \lfloor \frac{|I|}{2} \rfloor + 5)  \cdot  \frac{\eps^3 }{(1+t+|q|)^{1 -   c (\gamma)  \cdot c (\delta)  \cdot c(|I|) \cdot E ( \lfloor \frac{|I|}{2} \rfloor  + 5) \cdot \eps } }  \; \Big]^{\frac{1}{2}}  \; ,
\eeaa
and 

      \beaa
   \notag
 &&   |\Lie_{Z^I}  A_{e_{a}}   |  \\
 &\les&   \frac{1}{(1+t+|q|) \cdot (1+|q|)^{\ga} } \cdot \Big[ \;   \sum_{|J|\leq |I| + 2 }  \int_{\Sigma^{ext}_{t_1} }  |\derm  ( \Lie_{Z^J}  A_{e_{a}} )|^2     \cdot w(q)  \cdot d^{3}x \\
&& +   \sum_{|K| \leq |I |}   \int_{t_1}^{t}  \Big[    \;   \int_{\Sigma^{ext}_{\tau} }    \Big[       O \big(   C(q_0)   \cdot  c (\delta) \cdot c (\gamma) \cdot C(|I|) \cdot E ( \lfloor \frac{|I|}{2} \rfloor + 6)  \cdot \frac{\eps   \cdot  | \derm ( \Lie_{Z^K} h^1 ) |^2  }{(1+\tau+|q|)}   \big) \\
   && + O \big(    C(q_0)   \cdot  c (\delta) \cdot c (\gamma) \cdot C(|I|) \cdot E (   \lfloor \frac{|I|}{2} \rfloor +6)  \cdot \frac{\eps  \cdot  |\derm ( \Lie_{Z^K} A ) |^2  }{(1+\tau+|q|) }  \big)    \; \Big]    \cdot w(q)  \cdot d^{3}x  \Big]   \cdot d\tau  \\
      \notag
        &&  +     \int_{t_1}^{t}  \int_{\Sigma^{ext}_{\tau} }   \Big[     C(q_0)   \cdot  c (\delta) \cdot c (\gamma) \cdot C(|I|) \cdot E (\lfloor \frac{|I|}{2} \rfloor   +4)  \cdot \frac{\eps }{(1+\tau+|q|)^{1-   c (\gamma)  \cdot c (\delta)  \cdot c(|I|) \cdot E (  \lfloor \frac{|I|}{2} \rfloor  + 4)\cdot  \eps } \cdot (1+|q|)^2} \\
&& \times    \sum_{  |K| \leq |I| -1 }  | \derm ( \Lie_{Z^K}  A )  |^2    \Big]  \cdot w(q)  \cdot d^{3}x  \cdot d\tau \\
&& +   \int_{t_1}^{t}   \int_{\Sigma^{ext}_{\tau} }   \sum_{  |K| \leq |I| }       \Big[   C(q_0)   \cdot C(|I|) \cdot E (  \lfloor \frac{|I|}{2} \rfloor  +5)  \cdot \frac{\eps \cdot  | \Lie_{Z^{K}} H_{L  L} |^2 }{(1+\tau+|q|)^{1- 2 \de } \cdot (1+|q|)^{4+2\gamma}}  \Big] \cdot w(q)  \cdot d^{3}x  \cdot d\tau \\
 && +   C(q_0)   \cdot  c (\delta) \cdot c (\gamma) \cdot C(|I|) \cdot E ( \lfloor \frac{|I|}{2} \rfloor + 5)  \cdot  \frac{\eps^3 }{(1+t+|q|)^{1 -   c (\gamma)  \cdot c (\delta)  \cdot c(|I|) \cdot E ( \lfloor \frac{|I|}{2} \rfloor  + 5) \cdot \eps } }   \; \Big]^{\frac{1}{2}}  \; .
\eeaa

\end{lemma}

\begin{proof}

 \textbf{Estimate on $ \derm ( \Lie_{Z^I}  A_{e_{a}}  ) $:}\\ 
In the same way that we worked out the a prior estimates (see \cite{G4} for calculations), we obtain the following estimate using the Klainerman-Sobolev inequality for $\derm ( \Lie_{Z^I}  A_{e_{a}} )  $.

\bea
\notag
|\derm ( \Lie_{Z^I}  A_{e_{a}} )  |  \cdot (1+t+|q|)^\cdot \big[ (1+|q|) \cdot w(q)\big]^{1/2} \leq
C\sum_{|J|\leq |I| + 2 }    \int_{\Sigma^{ext}_{t_2} }  |\derm ( \Lie_{Z^J}  A_{e_{a}} )  |^2     \cdot w(q)  \cdot d^{3}x \;  . \\
\eea
where here the $L^2(\Sigma^{ext}_{t} )$ norm is taken on $\Sigma^{ext}_{t} $ slice.

This gives that in the exterior, we have

\bea
\notag
|\derm ( \Lie_{Z^I}  A_{e_{a}} )  |   \les \frac{1}{(1+t+|q|) \cdot (1+|q|)^{1+\ga} }
\cdot \Big( \sum_{|J|\leq |I| + 2 }    \int_{\Sigma^{ext}_{t_2} }  |\derm ( \Lie_{Z^J}  A_{e_{a}} )  |^2     \cdot w(q)  \cdot d^{3}x \Big)^{\frac{1}{2}} \;  . \\
\eea
  Using the estimate that we showed in Lemma \ref{ControlontheweightedLtwonormofdermAea}, we get 
  
      \bea\label{boundontheenergyforONLYtheAeacomponentoftheYangMillspotential}
   \notag
 &&     \int_{\Sigma^{ext}_{t} }  |\derm ( \Lie_{Z^I}  A_{e_{a}} )  |^2     \cdot w(q)  \cdot d^{3}x \\  
    \notag
 &\les&      \int_{\Sigma^{ext}_{t_1} }  |\derm  ( \Lie_{Z^I}  A_{e_{a}} )|^2     \cdot w(q)  \cdot d^{3}x \\
    \notag
&& +   \sum_{|K| \leq |I |}   \int_{t_1}^{t}  \Big[    \;   \int_{\Sigma^{ext}_{\tau} }    \Big[       O \big(   C(q_0)   \cdot  c (\delta) \cdot c (\gamma) \cdot C(|I|) \cdot E ( \lfloor \frac{|I|}{2} \rfloor + 5)  \cdot \frac{\eps   \cdot  | \derm ( \Lie_{Z^K} h^1 ) |^2  }{(1+\tau+|q|)}   \big) \\
   \notag
   && + O \big(    C(q_0)   \cdot  c (\delta) \cdot c (\gamma) \cdot C(|I|) \cdot E (   \lfloor \frac{|I|}{2} \rfloor +5)  \cdot \frac{\eps  \cdot  |\derm ( \Lie_{Z^K} A ) |^2  }{(1+\tau+|q|) }  \big)    \; \Big]  \cdot w(q)  \cdot d^{3}x    \Big]   \cdot d\tau  \\
      \notag
        &&  +     \int_{t_1}^{t}  \int_{\Sigma^{ext}_{\tau} }   \Big[     C(q_0)   \cdot  c (\delta) \cdot c (\gamma) \cdot C(|I|) \cdot E (   \lfloor \frac{|I|}{2} \rfloor   +3)  \cdot \frac{\eps }{(1+\tau+|q|)^{1-   c (\gamma)  \cdot c (\delta)  \cdot c(|I|) \cdot E (  \lfloor \frac{|I|}{2} \rfloor  + 2)\cdot  \eps } \cdot (1+|q|)^2} \\
           \notag
&& \times    \sum_{  |K| \leq |I| -1 }  | \derm ( \Lie_{Z^K}  A )  |^2     \Big]  \cdot w(q)  \cdot d^{3}x   \cdot d\tau  \\
   \notag
&& +   \int_{t_1}^{t}   \int_{\Sigma^{ext}_{\tau} }   \sum_{  |K| \leq |I| }           \Big[   C(q_0)   \cdot C(|I|) \cdot E (  \lfloor \frac{|I|}{2} \rfloor  +5)  \cdot \frac{\eps \cdot  | \Lie_{Z^{K}} H_{L  L} |^2 }{(1+\tau+|q|)^{1- 2 \de } \cdot (1+|q|)^{4+2\gamma}}  \Big] \cdot w(q)  \cdot d^{3}x  \cdot d\tau \\
\notag
 && +   C(q_0)   \cdot  c (\delta) \cdot c (\gamma) \cdot C(|I|) \cdot E ( \lfloor \frac{|I|}{2} \rfloor + 5)  \cdot  \frac{\eps^3 }{(1+t+|q|)^{1 -   c (\gamma)  \cdot c (\delta)  \cdot c(|I|) \cdot E ( \lfloor \frac{|I|}{2} \rfloor  + 5) \cdot \eps } }   \; . \\
\eea

      \beaa
   \notag
 &&   |\derm ( \Lie_{Z^I}  A_{e_{a}} )  |  \\
 &\les&   \frac{1}{(1+t+|q|) \cdot (1+|q|)^{1+\ga} } \cdot \Big[ \;  \sum_{|J|\leq |I| + 2 }  \int_{\Sigma^{ext}_{t_1} }   |\derm  ( \Lie_{Z^J}  A_{e_{a}} )|^2     \cdot w(q)  \cdot d^{3}x \\
&& +   \sum_{|K| \leq |I |}   \int_{t_1}^{t}  \Big[    \;   \int_{\Sigma^{ext}_{\tau} }    \Big[       O \big(   C(q_0)   \cdot  c (\delta) \cdot c (\gamma) \cdot C(|I|) \cdot E ( \lfloor \frac{|I|+2}{2} \rfloor + 5)  \cdot \frac{\eps   \cdot  | \derm ( \Lie_{Z^K} h^1 ) |^2  }{(1+t+|q|)}   \big) \\
   && + O \big(    C(q_0)   \cdot  c (\delta) \cdot c (\gamma) \cdot C(|I|) \cdot E (   \lfloor \frac{|I| +2}{2} \rfloor +5)  \cdot \frac{\eps  \cdot  |\derm ( \Lie_{Z^K} A ) |^2  }{(1+\tau+|q|) }  \big)    \; \Big]   \cdot w(q)  \cdot d^{3}x    \Big]  \cdot d\tau  \\
      \notag
        &&  +     \int_{t_1}^{t}  \int_{\Sigma^{ext}_{\tau} }    \Big[     C(q_0)   \cdot  c (\delta) \cdot c (\gamma) \cdot C(|I|) \cdot E (  \lfloor \frac{|I|+2}{2} \rfloor   +3 )  \cdot \frac{\eps }{(1+\tau+|q|)^{1-   c (\gamma)  \cdot c (\delta)  \cdot c(|I|) \cdot E (  \lfloor \frac{|I|}{2} \rfloor + 4)\cdot  \eps } \cdot (1+|q|)^2} \\
&& \times    \sum_{  |K| \leq |I| -1 }  | \derm ( \Lie_{Z^K}  A )  |^2     \Big]  \cdot w(q)  \cdot d^{3}x  \cdot d\tau  \\
&& +   \int_{t_1}^{t}   \int_{\Sigma^{ext}_{\tau} }   \sum_{  |K| \leq |I| }        \Big[   C(q_0)   \cdot C(|I|) \cdot E (  \lfloor \frac{|I|}{2} \rfloor  +5)  \cdot \frac{\eps \cdot  | \Lie_{Z^{K}} H_{L  L} |^2 }{(1+\tau+|q|)^{1- 2 \de } \cdot (1+|q|)^{4+2\gamma}}  \Big] \cdot w(q)  \cdot d^{3}x  \cdot d\tau \\
 && +   C(q_0)   \cdot  c (\delta) \cdot c (\gamma) \cdot C(|I|) \cdot E ( \lfloor \frac{|I|}{2} \rfloor + 5)  \cdot  \frac{\eps^3 }{(1+t+|q|)^{1 -   c (\gamma)  \cdot c (\delta)  \cdot c(|I|) \cdot E ( \lfloor \frac{|I|}{2} \rfloor  + 5) \cdot \eps } }  \; \Big]^{\frac{1}{2}}  \; .
\eeaa
\beaa
 &\les&   \frac{1}{(1+t+|q|) \cdot (1+|q|)^{1+\ga} } \cdot \Big[ \;   \sum_{|J|\leq |I| + 2 }  \int_{\Sigma^{ext}_{t_1} }  |\derm  ( \Lie_{Z^J}  A_{e_{a}} )|^2     \cdot w(q)  \cdot d^{3}x \\
&& +   \sum_{|K| \leq |I |}   \int_{t_1}^{t}  \Big[    \;   \int_{\Sigma^{ext}_{\tau} }    \Big[       O \big(   C(q_0)   \cdot  c (\delta) \cdot c (\gamma) \cdot C(|I|) \cdot E ( \lfloor \frac{|I|}{2} \rfloor + 6)  \cdot \frac{\eps   \cdot  | \derm ( \Lie_{Z^K} h^1 ) |^2  }{(1+\tau+|q|)}   \big) \\
   && + O \big(    C(q_0)   \cdot  c (\delta) \cdot c (\gamma) \cdot C(|I|) \cdot E (   \lfloor \frac{|I|}{2} \rfloor +6)  \cdot \frac{\eps  \cdot  |\derm ( \Lie_{Z^K} A ) |^2  }{(1+\tau+|q|) }  \big)    \; \Big]    \cdot w(q)  \cdot d^{3}x   \Big]  \cdot d\tau  \\
      \notag
        &&  +     \int_{t_1}^{t}  \int_{\Sigma^{ext}_{\tau} }   \Big[     C(q_0)   \cdot  c (\delta) \cdot c (\gamma) \cdot C(|I|) \cdot E (\lfloor \frac{|I|}{2} \rfloor   +4)  \cdot \frac{\eps }{(1+\tau+|q|)^{1-   c (\gamma)  \cdot c (\delta)  \cdot c(|I|) \cdot E (  \lfloor \frac{|I|}{2} \rfloor  + 4)\cdot  \eps } \cdot (1+|q|)^2} \\
&& \times    \sum_{  |K| \leq |I| -1 }  | \derm ( \Lie_{Z^K}  A )  |^2     \Big]  \cdot w(q)  \cdot d^{3}x  \cdot d\tau  \\
&& +   \int_{t_1}^{t}   \int_{\Sigma^{ext}_{\tau} }   \sum_{  |K| \leq |I| }        \Big[   C(q_0)   \cdot C(|I|) \cdot E (  \lfloor \frac{|I|}{2} \rfloor  +5)  \cdot \frac{\eps \cdot  | \Lie_{Z^{K}} H_{L  L} |^2 }{(1+\tau+|q|)^{1- 2 \de } \cdot (1+|q|)^{4+2\gamma}}  \Big]  \cdot w(q)  \cdot d^{3}x  \cdot d\tau \\
 && +   C(q_0)   \cdot  c (\delta) \cdot c (\gamma) \cdot C(|I|) \cdot E ( \lfloor \frac{|I|}{2} \rfloor + 5)  \cdot  \frac{\eps^3 }{(1+t+|q|)^{1 -   c (\gamma)  \cdot c (\delta)  \cdot c(|I|) \cdot E ( \lfloor \frac{|I|}{2} \rfloor  + 5) \cdot \eps } }  \; \Big]^{\frac{1}{2}}  \; .
\eeaa

\textbf{Estimate on $ \Lie_{Z^I}  A_{e_{a}}   $:}\\
In order to estimate first $ | \pa_{r} A_{e_{a}}  |$, so that we could then integrate along $s= constant$ and $\Om = constant$, we will use the special fact that
      \bea\label{specialfactaboutAeacomponentsthatexpolitsthefactthatitisnotonlyAtangentialbutspecialone}
   &&   \derm_{r} e_a = 0 \\
      \notag
 &&   \text{  (see  Appendix of \cite{G2}),}
      \eea
 and therefore
      \bea\label{partialderivativeindirectionofrofAeacomponenentiscovariantderivativeofAeacomponent}
      \notag
    \pa_r  \Lie_{Z^I} A_{e_{a}} &=& \derm_r  ( \Lie_{Z^I} A_{e_{a}} )  +  \Lie_{Z^I} A (\derm_{r} e_a)  \\
     &=& \derm_r  ( \Lie_{Z^I} A_{e_{a}} ) \; .
     \eea

   Thus, integrating, as we did in \cite{G4}, we get
       \beaa
\notag
| \Lie_{Z^I} A_{e_{a}}  (t, | x | \cdot \Om) | &=& | \Lie_{Z^I} A_{e_{a}} \big(0, ( t + | x |) \cdot \Om \big) |  + \int_{| x | }^{t + | x |  } \pa_r | (\Lie_{Z^I} A_{e_{a}}  (t + | x | - r,  r  \cdot \Om ) ) | dr \\
 &\leq& | \Lie_{Z^I} A_{e_{a}}  \big(0, ( t + | x | ) \cdot \Om \big) |  + \int_{| x | }^{t + | x |  } | \pa_r  (\Lie_{Z^I} A_{e_{a}}  (t + | x | - r,  r  \cdot \Om ) ) | dr \\
  &\leq& | \Lie_{Z^I} A_{e_{a}}  \big(0, ( t + | x | ) \cdot \Om \big)|  + \int_{| x | }^{t + | x |  } | \derm_r  (\Lie_{Z^I} A_{e_a} (t + | x | - r,  r  \cdot \Om ) ) | dr \; .
\eeaa
     However, we have shown that for $0 \leq \de \leq \frac{1}{4}$, $\ga > \de$, we have
     \beaa
 | \derm_r  ( \Lie_{Z^I} A_{e_a} )  |&\les&   | \derm (\Lie_{Z^I} A_{e_a} )  |   \\
      &\les& \frac{1}{(1+t+|q|) \cdot (1+|q|)^{1+\ga} }
\cdot \Big( \sum_{|J|\leq |I| + 2 }    \int_{\Sigma^{ext}_{t_2} }  |\derm ( \Lie_{Z^J}  A_{e_{a}} )  |^2     \cdot w(q)  \cdot d^{3}x \Big)^{\frac{1}{2}} \; .
\eeaa     
Hence, proceeding the integration, as we detailed in \cite{G4}, by integrating at a fixed $\Om \in \SSS^2$, from $(t, | x | \cdot \Om)$ along the line $(\tau, r \cdot \Om)$ such that $r+\tau =  | x | +t $ (i.e. along a fixed null coordinate $s:=\tau+r $) till we reach the hyperplane $\tau=0$, as in \eqref{definitionequationforintegralalongthenullcoordinateplusboundaryterm}, and using the asymptotic behaviour of the initial data (see Remark \ref{remarkabouttheasymptoticbehaviouroftheboundarytermonhyperplaneprescribedtqual0}), and using the bound on the energy in \eqref{boundontheenergyforONLYtheAeacomponentoftheYangMillspotential}, we get that in the stated result.

\end{proof}

\section{Estimate on the energy for the full components of the potential $A$}

\begin{lemma}\label{estimateonthebadtermproductA_atimesdermA_ausingbootstrapsoastoclosethegronwallinequalityonenergy}

For $\ga \geq 3 \de $\,, for $0 < \de \leq \frac{1}{4}$\,, and for $\eps$ small, enough depending on $q_0$\,, on $\ga$, on $\de$,  on $|I|$ and on $\mu < 0$\;, we have

  \beaa
\notag
&&  \frac{ (1+t)}{\eps} \cdot \Big( \sum_{|K| +|J| \leq |I|}   | \Lie_{Z^K}  A_{e_{a}}   | \cdot |\derm ( \Lie_{Z^J}  A_{e_{a}} )  | \Big)^{2}  \\
\notag
  &\les&   C(q_0) \cdot c (\gamma) \cdot  C ( |I| ) \cdot E (  \lfloor \frac{|I|}{2} \rfloor  + 2) \\
  && \times  \frac{\eps}{(1+t+|q|)^{3-2\de} \cdot (1+|q|)^{2+4\ga} } \cdot \Big[ \;  \sum_{|J|\leq |I| + 2 }   \int_{\Sigma^{ext}_{t_1} }  |\derm  ( \Lie_{Z^J}  A_{e_{a}} )|^2     \cdot w(q)  \cdot d^{3}x \\
&& +   \sum_{|K| \leq |I |}   \int_{t_1}^{t}  \Big[    \;   \int_{\Sigma^{ext}_{\tau} }    \Big[       O \big(   C(q_0)   \cdot  c (\delta) \cdot c (\gamma) \cdot C(|I|) \cdot E ( \lfloor \frac{|I|}{2} \rfloor + 6)  \cdot \frac{\eps   \cdot  | \derm ( \Lie_{Z^K} h ) |^2  }{(1+\tau+|q|)}   \big) \\
   && + O \big(    C(q_0)   \cdot  c (\delta) \cdot c (\gamma) \cdot C(|I|) \cdot E (   \lfloor \frac{|I|}{2} \rfloor +6)  \cdot \frac{\eps  \cdot  |\derm ( \Lie_{Z^K} A ) |^2  }{(1+\tau+|q|) }  \big)    \; \Big]    \cdot w(q)  \cdot d^{3}x   \Big]  \cdot d\tau  \\
      \notag
        &&  +     \int_{t_1}^{t}  \int_{\Sigma^{ext}_{\tau} }   \Big[     C(q_0)   \cdot  c (\delta) \cdot c (\gamma) \cdot C(|I|) \cdot E (\lfloor \frac{|I|}{2} \rfloor   +4)  \cdot \frac{\eps }{(1+\tau+|q|)^{1-   c (\gamma)  \cdot c (\delta)  \cdot c(|I|) \cdot E (  \lfloor \frac{|I|}{2} \rfloor + 4)\cdot  \eps } \cdot (1+|q|)^2} \\
&& \times    \sum_{  |K| \leq |I| -1 }  | \derm ( \Lie_{Z^K}  A )  |^2    \Big]  \cdot w(q)  \cdot d^{3}x  \cdot d\tau  \\
&& +   \int_{t_1}^{t}   \int_{\Sigma^{ext}_{\tau} }   \sum_{  |K| \leq |I| }        \Big[   C(q_0)   \cdot C(|I|) \cdot E (  \lfloor \frac{|I|}{2} \rfloor  +5)  \cdot \frac{\eps \cdot  | \Lie_{Z^{K}} H_{L  L} |^2 }{(1+\tau+|q|)^{1- 2 \de } \cdot (1+|q|)^{4+2\gamma}}  \Big] \cdot w(q)  \cdot d^{3}x  \cdot d\tau  \\
 && +   C(q_0)   \cdot  c (\delta) \cdot c (\gamma) \cdot C(|I|) \cdot E ( \lfloor \frac{|I|}{2} \rfloor + 5)  \cdot  \frac{\eps^3 }{(1+t+|q|)^{1 -   c (\gamma)  \cdot c (\delta)  \cdot c(|I|) \cdot E ( \lfloor \frac{|I|}{2} \rfloor  + 5) \cdot \eps } }  \; \Big]  \; .
\eeaa

\end{lemma}

\begin{proof}
We start first by decomposing the sum in the following way,

\bea\label{decompositionofsumforbadtermA_adermA_atouseboudednessonthelowerorderterm}
\notag
&& \sum_{|K| +|J| \leq |I|}   | \Lie_{Z^K}  A_{e_{a}}   | \cdot |\derm ( \Lie_{Z^J}  A_{e_{a}} )  | \\
\notag
  &\leq&  \sum_{ |J| \leq   \lfloor \frac{|I|}{2} \rfloor   ,\; |K| \leq |I| } | \Lie_{Z^K}  A_{e_{a}}   | \cdot |\derm ( \Lie_{Z^J}  A_{e_{a}} )  |  + \sum_{ |J| \leq   \lfloor \frac{|I|}{2} \rfloor   ,\; |K| \leq |I| } | \Lie_{Z^J}  A_{e_{a}}   | \cdot |\derm ( \Lie_{Z^K}  A_{e_{a}} )  |  \; . \\
\eea

However, for all $|J| \leq  \lfloor \frac{|I|}{2} \rfloor $, we have from the a priori estimates established in Lemma \ref{apriorestimateontheEinsteinYangMillsfieldswithoutgradiant}, we have

\beaa
 | \Lie_{Z^J}  A_{e_{a}}   |    &\les&  C(q_0) \cdot c (\gamma) \cdot  C ( |I| ) \cdot E (  \lfloor \frac{|I|}{2} \rfloor  + 2) \cdot  \frac{\eps }{(1+t+|q|)^{1-\delta} \cdot  (1+|q|)^{\gamma}}\; ,
 \eeaa
 and from Lemma \ref{apriordecayestimatesfrombootstrapassumption}, we have
 \beaa
  |\derm ( \Lie_{Z^J}  A_{e_{a}} )  |  &\les& C(q_0) \cdot C ( |I| ) \cdot E ( \lfloor \frac{|I|}{2} \rfloor + 2 )  \cdot \frac{\eps }{(1+t+|q|)^{1-\delta} \cdot  (1+|q|)^{1+\gamma}} \; .
  \eeaa
  Thus, injecting in \eqref{decompositionofsumforbadtermA_adermA_atouseboudednessonthelowerorderterm}, we obtain
  \beaa
\notag
&& \sum_{|K| +|J| \leq |I|}   | \Lie_{Z^K}  A_{e_{a}}   | \cdot |\derm ( \Lie_{Z^J}  A_{e_{a}} )  | \\
\notag
  &\les&   C(q_0)  \cdot  C ( |I| ) \cdot E (  \lfloor \frac{|I|}{2} \rfloor  + 2)  \cdot \frac{\eps }{(1+t+|q|)^{1-\delta} \cdot  (1+|q|)^{1+\gamma}} \cdot  \sum_{ |K|  \leq |I| }   | \Lie_{Z^K}  A_{e_{a}}   |   \\
&&  + C(q_0) \cdot c (\gamma) \cdot  C ( |I| ) \cdot E (  \lfloor \frac{|I|}{2} \rfloor  + 2)  \cdot \frac{\eps }{(1+t+|q|)^{1-\delta} \cdot  (1+|q|)^{\gamma}}  \cdot  \sum_{ |K|  \leq |I| }    |\derm ( \Lie_{Z^K}  A_{e_{a}} )  |  \\
  &\les&   C(q_0)  \cdot c (\gamma) \cdot  C ( |I| ) \cdot E (  \lfloor \frac{|I|}{2} \rfloor  + 2)  \cdot \frac{\eps }{(1+t+|q|)^{1-\delta} } \\
  &&  \times \Big[  \frac{1}{ (1+|q|)^{1+\gamma}}  \cdot  \sum_{ |K|  \leq |I| }   | \Lie_{Z^K}  A_{e_{a}}   | + \frac{1 }{ (1+|q|)^{\gamma}}  \cdot  \sum_{ |K|  \leq |I| }   | \Lie_{Z^K}  A_{e_{a}}   |      \Big]   \; . \\
\eeaa

Now, using Lemma \ref{upgradedestimateontheA_acomponentsthatconservesthedecayasoneoverttimesintegralsoftheenergysoastousegeneralizedgronwall}, we get
  \beaa
\notag
&& \sum_{|K| +|J| \leq |I|}   | \Lie_{Z^K}  A_{e_{a}}   | \cdot |\derm ( \Lie_{Z^J}  A_{e_{a}} )  | \\
\notag
  &\les&   C(q_0) \cdot c (\gamma) \cdot  C ( |I| ) \cdot E (  \lfloor \frac{|I|}{2} \rfloor  + 2)  \cdot \frac{\eps }{(1+t+|q|)^{1-\delta} } \\
  && \times  \frac{\eps}{(1+t+|q|) \cdot (1+|q|)^{1+2\ga} } \cdot \Big[ \;  \sum_{|J|\leq |I| + 2 }   \int_{\Sigma^{ext}_{t_1} }  |\derm  ( \Lie_{Z^J}  A_{e_{a}} )|^2     \cdot w(q)  \cdot d^{3}x \\
&& +   \sum_{|K| \leq |I |}   \int_{t_1}^{t}  \Big[    \;   \int_{\Sigma^{ext}_{\tau} }    \Big[       O \big(   C(q_0)   \cdot  c (\delta) \cdot c (\gamma) \cdot C(|I|) \cdot E ( \lfloor \frac{|I|}{2} \rfloor + 6)  \cdot \frac{\eps   \cdot  | \derm ( \Lie_{Z^K} h ) |^2  }{(1+\tau+|q|)}   \big) \\
   && + O \big(    C(q_0)   \cdot  c (\delta) \cdot c (\gamma) \cdot C(|I|) \cdot E (   \lfloor \frac{|I|}{2} \rfloor +6)  \cdot \frac{\eps  \cdot  |\derm ( \Lie_{Z^K} A ) |^2  }{(1+\tau+|q|) }  \big)    \; \Big]    \cdot w(q)  \cdot d^{3}x    \Big]  \cdot d\tau  \\
      \notag
        &&  +     \int_{t_1}^{t}  \int_{\Sigma^{ext}_{\tau} }   \Big[     C(q_0)   \cdot  c (\delta) \cdot c (\gamma) \cdot C(|I|) \cdot E (\lfloor \frac{|I|}{2} \rfloor   +4)  \cdot \frac{\eps }{(1+\tau+|q|)^{1-   c (\gamma)  \cdot c (\delta)  \cdot c(|I|) \cdot E ( \lfloor \frac{|I|}{2} \rfloor  + 4)\cdot  \eps } \cdot (1+|q|)^2} \\
&& \times    \sum_{  |K| \leq |I| -1 }  | \derm ( \Lie_{Z^K}  A )  |^2    \Big]  \cdot w(q)  \cdot d^{3}x  \cdot d\tau \\
&& +   \int_{t_1}^{t}   \int_{\Sigma^{ext}_{\tau} }   \sum_{  |K| \leq |I| }         \Big[   C(q_0)   \cdot C(|I|) \cdot E (  \lfloor \frac{|I|}{2} \rfloor  +5)  \cdot \frac{\eps \cdot  | \Lie_{Z^{K}} H_{L  L} |^2 }{(1+\tau+|q|)^{1- 2 \de } \cdot (1+|q|)^{4+2\gamma}}  \Big] \cdot w(q)  \cdot d^{3}x  \cdot d\tau\\
 && +   C(q_0)   \cdot  c (\delta) \cdot c (\gamma) \cdot C(|I|) \cdot E ( \lfloor \frac{|I|}{2} \rfloor + 5)  \cdot  \frac{\eps^3 }{(1+t+|q|)^{1 -   c (\gamma)  \cdot c (\delta)  \cdot c(|I|) \cdot E ( \lfloor \frac{|I|}{2} \rfloor  + 5) \cdot \eps } }  \; \Big]^{\frac{1}{2}}  \; .
\eeaa

Squaring and multiplying by $(1+t)$ gives the desired result. 

\end{proof}

\begin{lemma}\label{estimateonthebadtermproductALtimesdermAusingbootstrapsoastoclosethegronwallinequalityonenergy}
We have for $\ga \geq 3 \de $\,, and $0 < \de \leq \frac{1}{4}$\,,

                                         \beaa
\notag
&&   \frac{ (1+t)}{\eps}  \cdot \Big(  \sum_{|K| \leq |I|}  | \Lie_{Z^K}    \big(  A_L   \cdot     \derm A    \big)    |   \Big)^{2} \\
&\les&  C(q_0)   \cdot  c (\delta) \cdot c (\gamma)  \cdot E ( 4)   \cdot \sum_{|K| = |I| }  \Big(  \frac{ \eps  \cdot     |  \derm ( \Lie_{Z^K} A )  |^2   }{ (1+t+|q|)  \cdot  (1+|q|)^{2\gamma - 4\de}  } \\
&&  +    \frac{\eps \cdot   |  \Lie_{Z^K} A_L |^2   }{(1+t+|q|)^{1-      c (\gamma)  \cdot c (\delta)   \cdot E ( 4) \cdot \eps } \cdot (1+|q|)^{2+2\gamma - 4\de }}       \Big) \\
&& +  C(q_0)   \cdot  c (\delta) \cdot c (\gamma) \cdot C( \lfloor \frac{|I|-1}{2} \rfloor) \cdot E ( \lfloor \frac{|I|-1}{2} \rfloor + 4) \\
&& \times \sum_{ |K| \leq |I| -1}    \Big(    \frac{\eps  \cdot     |  \derm ( \Lie_{Z^K} A )  |^2 }{(1+t+|q|)^{1-      c (\gamma)  \cdot c (\delta)  \cdot c( \lfloor \frac{|I|-1}{2} \rfloor) \cdot E ( \lfloor \frac{|I|-1}{2} \rfloor+ 4) \cdot \eps } \cdot (1+|q|)^{2\gamma - 4\de }}    \\
&& + \frac{\eps \cdot  | \Lie_{Z^K}  A |^2  }{(1+t+|q|)^{1-      c (\gamma)  \cdot c (\delta)  \cdot c(\lfloor \frac{|I|-1}{2} \rfloor ) \cdot E (\lfloor \frac{|I|-1}{2} \rfloor + 4) \cdot \eps } \cdot (1+|q|)^{2+2\gamma - 4\de }}   \Big) \;. 
\eeaa   

\end{lemma}

\begin{proof}
Using Corollary \ref{thestructureofthebadtermALdermAusingboostrapassumptionanddecompistionofthesumandlowerordertermsexhibitedtodealwithAL}, we get that for $\ga \geq 3 \de $\,, and $0 < \de \leq \frac{1}{4}$\,,
                                         \beaa
\notag
&& | \Lie_{Z^I}    \big(  A_L   \cdot     \derm A    \big)   | \\
&\les&  C(q_0)   \cdot  c (\delta) \cdot c (\gamma)  \cdot E ( 4)  \\
&& \times \sum_{|K| = |I| }  \Big(  \frac{ \eps  \cdot     |  \derm ( \Lie_{Z^K} A )  |   }{ (1+t+|q|)  \cdot  (1+|q|)^{\gamma - 2\de}  } \\
&&  +    \frac{\eps \cdot   |  \Lie_{Z^K} A_L |   }{(1+t+|q|)^{1-      c (\gamma)  \cdot c (\delta)   \cdot E ( 4) \cdot \eps } \cdot (1+|q|)^{1+\gamma - 2\de }}     \Big) \\
&& +  C(q_0)   \cdot  c (\delta) \cdot c (\gamma) \cdot C( \lfloor \frac{|I|-1}{2} \rfloor) \cdot E ( \lfloor \frac{|I|-1}{2} \rfloor + 4) \\
&& \times \sum_{ |K| \leq |I| -1}    \Big(    \frac{\eps  \cdot     |  \derm ( \Lie_{Z^K} A )  | }{(1+t+|q|)^{1-      c (\gamma)  \cdot c (\delta)  \cdot c( \lfloor \frac{|I|-1}{2} \rfloor) \cdot E ( \lfloor \frac{|I|-1}{2} \rfloor+ 4) \cdot \eps } \cdot (1+|q|)^{\gamma - 2\de }}    \\
&& + \frac{\eps \cdot  | \Lie_{Z^K}  A |  }{(1+t+|q|)^{1-      c (\gamma)  \cdot c (\delta)  \cdot c(\lfloor \frac{|I|-1}{2} \rfloor ) \cdot E (\lfloor \frac{|I|-1}{2} \rfloor + 4) \cdot \eps } \cdot (1+|q|)^{1+\gamma - 2\de }}   \Big) \;. 
\eeaa

We have

                                         \beaa
\notag
&& \sum_{|K| \leq |I|}  | \Lie_{Z^K}    \big(  A_L   \cdot     \derm A    \big)    | \\
&\les&  C(q_0)   \cdot  c (\delta) \cdot c (\gamma)  \cdot E ( 4)   \cdot \sum_{|K| = |I| }  \Big(  \frac{ \eps  \cdot     |  \derm ( \Lie_{Z^K} A )  |   }{ (1+t+|q|)  \cdot  (1+|q|)^{\gamma - 2\de}  } \\
&&  +    \frac{\eps \cdot   |  \Lie_{Z^K} A_L |   }{(1+t+|q|)^{1-      c (\gamma)  \cdot c (\delta)   \cdot E ( 4) \cdot \eps } \cdot (1+|q|)^{1+\gamma - 2\de }}       \Big) \\
&& +  C(q_0)   \cdot  c (\delta) \cdot c (\gamma) \cdot C( \lfloor \frac{|I|-1}{2} \rfloor) \cdot E ( \lfloor \frac{|I|-1}{2} \rfloor + 4) \\
&& \times \sum_{ |K| \leq |I| -1}    \Big(    \frac{\eps  \cdot     |  \derm ( \Lie_{Z^K} A )  | }{(1+t+|q|)^{1-      c (\gamma)  \cdot c (\delta)  \cdot c( \lfloor \frac{|I|-1}{2} \rfloor) \cdot E ( \lfloor \frac{|I|-1}{2} \rfloor+ 4) \cdot \eps } \cdot (1+|q|)^{\gamma - 2\de }}    \\
&& + \frac{\eps \cdot  | \Lie_{Z^K}  A |  }{(1+t+|q|)^{1-      c (\gamma)  \cdot c (\delta)  \cdot c(\lfloor \frac{|I|-1}{2} \rfloor ) \cdot E (\lfloor \frac{|I|-1}{2} \rfloor + 4) \cdot \eps } \cdot (1+|q|)^{1+\gamma - 2\de }}   \Big) \;. 
\eeaa     
Thus, we obtain the stated result.

\end{proof}

\begin{lemma}\label{HardyinequalityontheestimateonthebadtermproductA_atimesdermA_ausingbootstrapsoastoclosethegronwallinequalityonenergy}
For $\ga \geq 3 \de $\,, for $0 < \de \leq \frac{1}{4}$\,, and for $\eps$ small, enough depending on $q_0$\;, on $\ga$\;, on $\de$\;, on $|I|$ and on $\mu$\;, we have

 \beaa
\notag
&&   \int_{\Sigma^{ext}_{\tau} }  \frac{ (1+t)}{\eps} \cdot \Big( \sum_{|K| +|J| \leq |I|}   | \Lie_{Z^K}  A_{e_{a}}   | \cdot |\derm ( \Lie_{Z^J}  A_{e_{a}} )  | \Big)^{2}   \cdot w(q)  \\
\notag
  &\les&   C(q_0) \cdot c (\gamma) \cdot  C ( |I| ) \cdot E (  \lfloor \frac{|I|}{2} \rfloor  + 2) \\
  && \times  \frac{\eps}{(1+t)^{1+2\de}}  \cdot \Big[ \;  \sum_{|J|\leq |I| + 2 }   \int_{\Sigma^{ext}_{t_1} }  |\derm  ( \Lie_{Z^J}  A_{e_{a}} )|^2     \cdot w(q)  \cdot d^{3}x \\
&& +   \sum_{|K| \leq |I |}   \int_{t_1}^{t}  \Big[    \;   \int_{\Sigma^{ext}_{\tau} }    \Big[       O \big(   C(q_0)   \cdot  c (\delta) \cdot c (\gamma) \cdot C(|I|) \cdot E ( \lfloor \frac{|I|}{2} \rfloor + 6)  \cdot \frac{\eps   \cdot  | \derm ( \Lie_{Z^K} h^1 ) |^2  }{(1+\tau+|q|)}   \big) \\
   && + O \big(    C(q_0)   \cdot  c (\delta) \cdot c (\gamma) \cdot C(|I|) \cdot E (   \lfloor \frac{|I|}{2} \rfloor +6)  \cdot \frac{\eps  \cdot  |\derm ( \Lie_{Z^K} A ) |^2  }{(1+\tau+|q|) }  \big)    \; \Big]    \cdot w(q)  \cdot d^{3}x   \Big]  \cdot d\tau  \\
      \notag
        &&  +     \int_{t_1}^{t}  \int_{\Sigma^{ext}_{\tau} }   \Big[     C(q_0)   \cdot  c (\delta) \cdot c (\gamma) \cdot C(|I|) \cdot E (\lfloor \frac{|I|}{2} \rfloor   +4)  \cdot \frac{\eps }{(1+\tau+|q|)^{1-   c (\gamma)  \cdot c (\delta)  \cdot c(|I|) \cdot E (  \lfloor \frac{|I|}{2} \rfloor + 4)\cdot  \eps } \cdot (1+|q|)^2} \\
&& \times    \sum_{  |K| \leq |I| -1 }  | \derm ( \Lie_{Z^K}  A )  |^2     \Big]  \cdot w(q)  \cdot d^{3}x  \cdot d\tau  \\
&& +   \int_{t_1}^{t_2}   \int_{\Sigma^{ext}_{\tau} }   \sum_{  |K| \leq |I| }         \Big[   C(q_0)   \cdot C(|I|) \cdot E (  \lfloor \frac{|I|}{2} \rfloor  +5)  \cdot \frac{\eps \cdot  | \Lie_{Z^{K}} H_{L  L} |^2 }{(1+\tau+|q|)^{1- 2 \de } \cdot (1+|q|)^{4+2\gamma}}  \Big]  \cdot w(q)  \cdot d^{3}x  \cdot d\tau \\
 && +   C(q_0)   \cdot  c (\delta) \cdot c (\gamma) \cdot C(|I|) \cdot E ( \lfloor \frac{|I|}{2} \rfloor + 5)  \cdot  \frac{\eps^3 }{(1+t)^{1 -   c (\gamma)  \cdot c (\delta)  \cdot c(|I|) \cdot E ( \lfloor \frac{|I|}{2} \rfloor  + 5) \cdot \eps } }   \; \Big]     \; .
\eeaa

\end{lemma}

\begin{proof}

From Lemma \ref{estimateonthebadtermproductA_atimesdermA_ausingbootstrapsoastoclosethegronwallinequalityonenergy}, and from the fact that
\beaa
 \frac{ w(q)  }{(1+|q|)^{2+4\ga} } \leq  \frac{ 1  }{(1+|q|)^{1+2\ga} } \;,
\eeaa
we get

  \beaa
\notag
&&   \int_{\Sigma^{ext}_{\tau} }   \frac{ (1+t)}{\eps} \cdot \Big( \sum_{|K| +|J| \leq |I|}   | \Lie_{Z^K}  A_{e_{a}}   | \cdot |\derm ( \Lie_{Z^J}  A_{e_{a}} )  | \Big)^{2}   \cdot w(q)  \\
\notag
  &\les&   C(q_0) \cdot c (\gamma) \cdot  C ( |I| ) \cdot E (  \lfloor \frac{|I|}{2} \rfloor  + 2) \\
  && \times \int_{\Sigma^{ext}_{\tau} } \Big[   \frac{\eps  \cdot w(q)  }{(1+t+|q|)^{3-2\de} \cdot (1+|q|)^{2+4\ga} } \cdot \Big[ \;   \sum_{|J|\leq |I| + 2 }  \int_{\Sigma^{ext}_{t_1} }  |\derm  ( \Lie_{Z^J}  A_{e_{a}} )|^2     \cdot w(q)  \cdot d^{3}x \\
&& +   \sum_{|K| \leq |I |}   \int_{t_1}^{t}  \Big[    \;   \int_{\Sigma^{ext}_{\tau} }    \Big[       O \big(   C(q_0)   \cdot  c (\delta) \cdot c (\gamma) \cdot C(|I|) \cdot E ( \lfloor \frac{|I|}{2} \rfloor + 6)  \cdot \frac{\eps   \cdot  | \derm ( \Lie_{Z^K} h^1 ) |^2  }{(1+\tau+|q|)}   \big) \\
   && + O \big(    C(q_0)   \cdot  c (\delta) \cdot c (\gamma) \cdot C(|I|) \cdot E (   \lfloor \frac{|I|}{2} \rfloor +6)  \cdot \frac{\eps  \cdot  |\derm ( \Lie_{Z^K} A ) |^2  }{(1+\tau+|q|) }  \big)    \; \Big]    \cdot w(q)  \cdot d^{3}x   \Big]  \cdot d\tau  \\
      \notag
        &&  +     \int_{t_1}^{t}  \int_{\Sigma^{ext}_{\tau} }   \Big[     C(q_0)   \cdot  c (\delta) \cdot c (\gamma) \cdot C(|I|) \cdot E (\lfloor \frac{|I|}{2} \rfloor   +4 )  \cdot \frac{\eps }{(1+\tau+|q|)^{1-   c (\gamma)  \cdot c (\delta)  \cdot c(|I|) \cdot E ( \lfloor \frac{|I|}{2} \rfloor  + 4)\cdot  \eps } \cdot (1+|q|)^2} \\
&& \times    \sum_{  |K| \leq |I| -1 } | \derm ( \Lie_{Z^K}  A )  |^2    \Big]  \cdot w(q)  \cdot d^{3}x  \cdot d\tau \\
&& +   \int_{t_1}^{t_2}   \int_{\Sigma^{ext}_{\tau} }   \sum_{  |K| \leq |I| }        \Big[   C(q_0)   \cdot C(|I|) \cdot E (  \lfloor \frac{|I|}{2} \rfloor  +5)  \cdot \frac{\eps \cdot  | \Lie_{Z^{K}} H_{L  L} |^2 }{(1+\tau+|q|)^{1- 2 \de } \cdot (1+|q|)^{4+2\gamma}}  \Big] \cdot w(q)  \cdot d^{3}x  \cdot d\tau  \\
 && +   C(q_0)   \cdot  c (\delta) \cdot c (\gamma) \cdot C(|I|) \cdot E ( \lfloor \frac{|I|}{2} \rfloor + 5)  \cdot  \frac{\eps^3 }{(1+t+|q|)^{1 -   c (\gamma)  \cdot c (\delta)  \cdot c(|I|) \cdot E ( \lfloor \frac{|I|}{2} \rfloor  + 5) \cdot \eps } }  \; \Big]     \; .
\eeaa

\beaa
  &\les&   C(q_0) \cdot c (\gamma) \cdot  C ( |I| ) \cdot E (  \lfloor \frac{|I|}{2} \rfloor  + 2) \\
  && \times  \frac{\eps}{(1+t)^{1+2\de}} \cdot \int_{\Sigma^{ext}_{\tau} } \Big[   \frac{1}{(1+|q|)^{3+2\ga-4\de} } \cdot \Big[ \;   \sum_{|J|\leq |I| + 2 }  \int_{\Sigma^{ext}_{t_1} }  |\derm  ( \Lie_{Z^J}  A_{e_{a}} )|^2     \cdot w(q)  \cdot d^{3}x \\
&& +   \sum_{|K| \leq |I |}   \int_{t_1}^{t}  \Big[    \;   \int_{\Sigma^{ext}_{\tau} }    \Big[       O \big(   C(q_0)   \cdot  c (\delta) \cdot c (\gamma) \cdot C(|I|) \cdot E ( \lfloor \frac{|I|}{2} \rfloor + 6)  \cdot \frac{\eps   \cdot  | \derm ( \Lie_{Z^K} h^1 ) |^2  }{(1+\tau+|q|)}   \big) \\
   && + O \big(    C(q_0)   \cdot  c (\delta) \cdot c (\gamma) \cdot C(|I|) \cdot E (   \lfloor \frac{|I|}{2} \rfloor +6)  \cdot \frac{\eps  \cdot  |\derm ( \Lie_{Z^K} A ) |^2  }{(1+\tau+|q|) }  \big)    \; \Big]    \cdot w(q)  \cdot d^{3}x   \Big]  \cdot d\tau  \\
      \notag
        &&  +     \int_{t_1}^{t}  \int_{\Sigma^{ext}_{\tau} }   \Big[     C(q_0)   \cdot  c (\delta) \cdot c (\gamma) \cdot C(|I|) \cdot E (\lfloor \frac{|I|}{2} \rfloor   +4)  \cdot \frac{\eps }{(1+\tau+|q|)^{1-   c (\gamma)  \cdot c (\delta)  \cdot c(|I|) \cdot E (  \lfloor \frac{|I|}{2} \rfloor + 4)\cdot  \eps } \cdot (1+|q|)^2} \\
&& \times    \sum_{  |K| \leq |I| -1 }  | \derm ( \Lie_{Z^K}  A )  |^2    \Big]  \cdot w(q)  \cdot d^{3}x  \cdot d\tau  \\
&& +   \int_{t_1}^{t_2}   \int_{\Sigma^{ext}_{\tau} }   \sum_{  |K| \leq |I| }         \Big[   C(q_0)   \cdot C(|I|) \cdot E (  \lfloor \frac{|I|}{2} \rfloor  +5)  \cdot \frac{\eps \cdot  | \Lie_{Z^{K}} H_{L  L} |^2 }{(1+\tau+|q|)^{1- 2 \de } \cdot (1+|q|)^{4+2\gamma}}  \Big] \cdot w(q)  \cdot d^{3}x  \cdot d\tau \\
 && +   C(q_0)   \cdot  c (\delta) \cdot c (\gamma) \cdot C(|I|) \cdot E ( \lfloor \frac{|I|}{2} \rfloor + 5)  \cdot  \frac{\eps^3 }{(1+t)^{1 -   c (\gamma)  \cdot c (\delta)  \cdot c(|I|) \cdot E ( \lfloor \frac{|I|}{2} \rfloor  + 5) \cdot \eps } }  \; \Big]     \; .
\eeaa

For $\ga \geq  3 \de$\;, we have  $2\ga - 4\de \geq  2 \de$\;, and thus, for $\de > 0$, we get
\bea
\notag
\int_{\Sigma^{ext}_{\tau} }  \frac{1}{(1+|q|)^{3+2\ga-4\de} } & \leq& \int_{\Sigma^{ext}_{\tau} }  \frac{1}{(1+|q|)^{3+2\de} }  \leq \int_{r=0}^{\infty}  \frac{1}{(1+|q|)^{1+2\ga} } dr \\
&\leq & C \; .
\eea
Therefore, we obtain
 \beaa
\notag
&&   \int_{\Sigma^{ext}_{\tau} }   \frac{ (1+t)}{\eps}  \cdot \Big( \sum_{|K| +|J| \leq |I|}   | \Lie_{Z^K}  A_{e_{a}}   | \cdot |\derm ( \Lie_{Z^J}  A_{e_{a}} )  | \Big)^{2}   \cdot w(q)  \\
\notag
  &\les&   C(q_0) \cdot c (\gamma) \cdot  C ( |I| ) \cdot E (  \lfloor \frac{|I|}{2} \rfloor  + 2) \\
  && \times  \frac{\eps}{(1+t)^{1+2\de}}  \cdot \Big[ \;    \sum_{|J|\leq |I| + 2 }  \int_{\Sigma^{ext}_{t_1} }  |\derm  ( \Lie_{Z^J}  A_{e_{a}} )|^2     \cdot w(q)  \cdot d^{3}x \\
&& +   \sum_{|K| \leq |I |}   \int_{t_1}^{t}  \Big[    \;   \int_{\Sigma^{ext}_{\tau} }    \Big[       O \big(   C(q_0)   \cdot  c (\delta) \cdot c (\gamma) \cdot C(|I|) \cdot E ( \lfloor \frac{|I|}{2} \rfloor + 6)  \cdot \frac{\eps   \cdot  | \derm ( \Lie_{Z^K} h^1 ) |^2  }{(1+\tau+|q|)}   \big) \\
   && + O \big(    C(q_0)   \cdot  c (\delta) \cdot c (\gamma) \cdot C(|I|) \cdot E (   \lfloor \frac{|I|}{2} \rfloor +6)  \cdot \frac{\eps  \cdot  |\derm ( \Lie_{Z^K} A ) |^2  }{(1+\tau+|q|) }  \big)    \; \Big]    \cdot w(q)  \cdot d^{3}x   \Big]  \cdot d\tau  \\
      \notag
        &&  +     \int_{t_1}^{t}  \int_{\Sigma^{ext}_{\tau} }   \Big[     C(q_0)   \cdot  c (\delta) \cdot c (\gamma) \cdot C(|I|) \cdot E ( \lfloor \frac{|I|}{2} \rfloor   +4 )  \cdot \frac{\eps }{(1+\tau+|q|)^{1-   c (\gamma)  \cdot c (\delta)  \cdot c(|I|) \cdot E ( \lfloor \frac{|I|}{2} \rfloor + 4)\cdot  \eps } \cdot (1+|q|)^2} \\
&& \times    \sum_{  |K| \leq |I| -1 }  | \derm ( \Lie_{Z^K}  A )  |^2   \Big]  \cdot w(q)  \cdot d^{3}x  \cdot d\tau   \\
&& +   \int_{t_1}^{t_2}   \int_{\Sigma^{ext}_{\tau} }   \sum_{  |K| \leq |I| }         \Big[   C(q_0)   \cdot C(|I|) \cdot E (  \lfloor \frac{|I|}{2} \rfloor  +5)  \cdot \frac{\eps \cdot  | \Lie_{Z^{K}} H_{L  L} |^2 }{(1+\tau+|q|)^{1- 2 \de } \cdot (1+|q|)^{4+2\gamma}}  \Big]  \cdot w(q)  \cdot d^{3}x  \cdot d\tau   \\
 && +   C(q_0)   \cdot  c (\delta) \cdot c (\gamma) \cdot C(|I|) \cdot E ( \lfloor \frac{|I|}{2} \rfloor + 5)  \cdot  \frac{\eps^3 }{(1+t)^{1 -   c (\gamma)  \cdot c (\delta)  \cdot c(|I|) \cdot E ( \lfloor \frac{|I|}{2} \rfloor  + 5) \cdot \eps } }  \; \Big]      \; .
\eeaa

\end{proof}

\begin{lemma}\label{HardytpyeestimatesontehtermscontainingALinthebadstructureofsourcestermsforAusefultoobtainenergyestimates}

We have for $\ga \geq 3 \de $\,, and $0 < \de \leq \frac{1}{4}$\,,

                                         \beaa
\notag
&&  \int_{\Sigma^{ext}_{\tau} }  \frac{ (1+t)}{\eps}  \cdot \Big(  \sum_{|K| \leq |I|}  | \Lie_{Z^K}    \big(  A_L   \cdot     \derm A    \big)    |   \Big)^{2}  \cdot w(q)  \\
&\les&   \int_{\Sigma^{ext}_{\tau} } \Big[ C(q_0)   \cdot  c (\delta) \cdot c (\gamma)  \cdot E ( 4)   \\
&& \times \sum_{|K| = |I| }  \Big(  \frac{ \eps  \cdot     |  \derm ( \Lie_{Z^K} A )  |^2   }{ (1+t+|q|)  \cdot  (1+|q|)^{2\gamma - 4\de}  }   +    \frac{\eps \cdot   | \pa \Lie_{Z^K} A_L |^2   }{(1+t+|q|)^{1-      c (\gamma)  \cdot c (\delta)   \cdot E ( 4) \cdot \eps } \cdot (1+|q|)^{2\gamma - 4\de }}       \Big) \; \Big]  \cdot w(q)  \\
&& +   \int_{\Sigma^{ext}_{\tau} } \Big[  C(q_0)   \cdot  c (\delta) \cdot c (\gamma) \cdot C( |I|) \cdot E ( \lfloor \frac{|I|-1}{2} \rfloor + 4) \\
&& \times \sum_{ |K| \leq |I| -1}    \Big(    \frac{\eps  \cdot     |  \derm ( \Lie_{Z^K} A )  |^2 }{(1+t+|q|)^{1-      c (\gamma)  \cdot c (\delta)  \cdot c( \lfloor \frac{|I|-1}{2} \rfloor) \cdot E ( \lfloor \frac{|I|-1}{2} \rfloor+ 4) \cdot \eps } \cdot (1+|q|)^{2\gamma - 4\de }}      \Big)   \; \Big]   \cdot w(q)  \;. 
\eeaa

\end{lemma}

\begin{proof}

We have from Lemma \ref{estimateonthebadtermproductALtimesdermAusingbootstrapsoastoclosethegronwallinequalityonenergy}, 
                                         \beaa
\notag
&&  \int_{\Sigma^{ext}_{\tau} }  \frac{ (1+t)}{\eps}  \cdot \Big(  \sum_{|K| \leq |I|}  | \Lie_{Z^K}    \big(  A_L   \cdot     \derm A    \big)    |   \Big)^{2}  \cdot w(q)  \\
&\les&   \int_{\Sigma^{ext}_{\tau} } \Big[ C(q_0)   \cdot  c (\delta) \cdot c (\gamma)  \cdot E ( 4)   \cdot \sum_{|K| = |I| }  \Big(  \frac{ \eps  \cdot     |  \derm ( \Lie_{Z^K} A )  |^2   }{ (1+t+|q|)  \cdot  (1+|q|)^{2\gamma - 4\de}  } \\
&&  +    \frac{\eps \cdot   |  \Lie_{Z^K} A_L |^2   }{(1+t+|q|)^{1-      c (\gamma)  \cdot c (\delta)   \cdot E ( 4) \cdot \eps } \cdot (1+|q|)^{2+2\gamma - 4\de }}       \Big) \\
&& +  C(q_0)   \cdot  c (\delta) \cdot c (\gamma) \cdot C( \lfloor \frac{|I|-1}{2} \rfloor) \cdot E ( \lfloor \frac{|I|-1}{2} \rfloor + 4) \\
&& \times \sum_{ |K| \leq |I| -1}    \Big(    \frac{\eps  \cdot     |  \derm ( \Lie_{Z^K} A )  |^2 }{(1+t+|q|)^{1-      c (\gamma)  \cdot c (\delta)  \cdot c( \lfloor \frac{|I|-1}{2} \rfloor) \cdot E ( \lfloor \frac{|I|-1}{2} \rfloor+ 4) \cdot \eps } \cdot (1+|q|)^{2\gamma - 4\de }}    \\
&& + \frac{\eps \cdot  | \Lie_{Z^K}  A |^2  }{(1+t+|q|)^{1-      c (\gamma)  \cdot c (\delta)  \cdot c(\lfloor \frac{|I|-1}{2} \rfloor ) \cdot E (\lfloor \frac{|I|-1}{2} \rfloor + 4) \cdot \eps } \cdot (1+|q|)^{2+2\gamma - 4\de }}   \Big)  \; \Big]  \cdot w(q)   \;. 
\eeaa   

Using the Hardy type inequality in the exterior from Corollary \ref{HardytypeinequalityforintegralstartingatROm}, we get that for fields decaying fast enough at spatial infinity, we have,

                                         \beaa
\notag
&&  \int_{\Sigma^{ext}_{\tau} } ( \frac{ (1+t)}{\eps}   \cdot \Big(  \sum_{|K| \leq |I|}  | \Lie_{Z^K}    \big(  A_L   \cdot     \derm A    \big)    |   \Big)^{2}  \cdot w(q)  \\
&\les&   \int_{\Sigma^{ext}_{\tau} } \Big[ C(q_0)   \cdot  c (\delta) \cdot c (\gamma)  \cdot E ( 4)   \cdot \sum_{|K| = |I| }  \Big(  \frac{ \eps  \cdot     |  \derm ( \Lie_{Z^K} A )  |^2   }{ (1+t+|q|)  \cdot  (1+|q|)^{2\gamma - 4\de}  } \\
&&  +    \frac{\eps \cdot   | \pa \Lie_{Z^K} A_L |^2   }{(1+t+|q|)^{1-      c (\gamma)  \cdot c (\delta)   \cdot E ( 4) \cdot \eps } \cdot (1+|q|)^{2\gamma - 4\de }}       \Big) \\
&& +  C(q_0)   \cdot  c (\delta) \cdot c (\gamma) \cdot C( \lfloor \frac{|I|-1}{2} \rfloor) \cdot E ( \lfloor \frac{|I|-1}{2} \rfloor + 4) \\
&& \times \sum_{ |K| \leq |I| -1}    \Big(    \frac{\eps  \cdot     |  \derm ( \Lie_{Z^K} A )  |^2 }{(1+t+|q|)^{1-      c (\gamma)  \cdot c (\delta)  \cdot c( \lfloor \frac{|I|-1}{2} \rfloor) \cdot E ( \lfloor \frac{|I|-1}{2} \rfloor+ 4) \cdot \eps } \cdot (1+|q|)^{2\gamma - 4\de }}    \\
&& + \frac{\eps \cdot  | \derm ( \Lie_{Z^K}  A ) |^2  }{(1+t+|q|)^{1-      c (\gamma)  \cdot c (\delta)  \cdot c(\lfloor \frac{|I|-1}{2} \rfloor ) \cdot E (\lfloor \frac{|I|-1}{2} \rfloor + 4) \cdot \eps } \cdot (1+|q|)^{2\gamma - 4\de }}   \Big)  \; \Big]  \cdot w(q)   \;. 
\eeaa   

Hence, we obtain the desired result.

\end{proof}

\begin{lemma}\label{estimateonthetermthatcontainspartialderivatoveofALinthettimessquareofbadcomponenentsinsourcetermsforAthatcontainstheALcomponent}
For $\ga \geq 3\de$, and for $\eps$ small depending on $|I|$\;, on $\ga$ and on $\de$\;, we have
\beaa
 && \int_{\Sigma^{ext}_{\tau} } C(q_0)   \cdot  c (\delta) \cdot c (\gamma)  \cdot E ( 4)   \cdot \sum_{|K| = |I| }  \Big(    \frac{\eps \cdot   | \pa \Lie_{Z^K} A_L |^2   }{(1+t+|q|)^{1-      c (\gamma)  \cdot c (\delta)   \cdot E ( 4) \cdot \eps } \cdot (1+|q|)^{2\gamma - 4\de }}       \Big)  \cdot w(q)   \\
 &\les & \int_{\Sigma^{ext}_{\tau} }C(q_0)   \cdot  c (\delta) \cdot c (\gamma)    \cdot  E (\lfloor \frac{|I|}{2} \rfloor +4  ) \cdot  \\
 && \times   \Big[  \sum_{|K| \leq |I| }     \frac{\eps \cdot   | \rderm ( \Lie_{Z^K} A) |^2   }{(1+t+|q|)^{1-      c (\gamma)  \cdot c (\delta)   \cdot E ( 4) \cdot \eps } \cdot (1+|q|)^{2\gamma - 4\de }}         \\
 &&+  \sum_{|K| \leq |I|-1 }     \frac{\eps \cdot   | \derm \Lie_{Z^K} A |^2   }{(1+t+|q|)^{1-      c (\gamma)  \cdot c (\delta)   \cdot E ( 4) \cdot \eps } \cdot (1+|q|)^{2\gamma - 4\de }}       \Big]  \cdot w(q)   \\
&&+   \int_{\Sigma^{ext}_{\tau} } C(q_0)   \cdot  c (\delta) \cdot c (\gamma) \cdot C(|I|) \cdot E (  \lfloor \frac{|I|}{2} \rfloor  +4)  \\
&&  \times  \frac{\eps }{(1+t+|q|)^{2 } }  \cdot  \sum_{|K| \leq |I|}  \Big[  | \derm  (  \Lie_{Z^K}  A ) |^2 +  |   \derm ( \Lie_{Z^K} h^1 )   |^2  \Big]  \cdot w(q) \\
&&+   C(q_0)   \cdot  c (\delta) \cdot c (\gamma) \cdot C(|I|) \cdot E (  \lfloor \frac{|I|}{2} \rfloor  +4)  \cdot  \frac{\eps^3 }{(1+t)^{1-     c (\gamma)  \cdot c (\delta)  \cdot c(|I|) \cdot E ( \lfloor \frac{|I|}{2} \rfloor+ 4) \cdot \eps } }      \; .
\eeaa

\end{lemma}

\begin{proof}

Using Lemma \ref{estimateonpartialderivativeofALcomponent} for the term containing the component $\pa \Lie_{Z^K} A_L$, where we showed in \eqref{partialderivativeofLeiderivativeofALcomponent}, that for all $|J| \leq |I|$, we have
  
\beaa
\notag
&& \sum_{|J| \leq |I|}  | \pa  \Lie_{Z^J}  A_{L}  | \\
 &\les& E (\lfloor \frac{|I|}{2} \rfloor  ) \cdot \Big[ \sum_{|J|\leq |I|} |   \rderm  \Lie_{Z^J} A | + \sum_{|J|\leq |I| -1} |   \derm  \Lie_{Z^J} A |  + \sum_{|K| + |M| \leq |I|} O \big(  | ( \Lie_{Z^K}  h )  | \cdot  | \derm  (  \Lie_{Z^M}  A ) | \big) \Big]  \; ,\\
\eeaa

Now, we decompose the sum as 
\beaa
&&\sum_{|K| + |M| \leq |I|} O \big(  | ( \Lie_{Z^K}  h )  | \cdot  | \derm  (  \Lie_{Z^M}  A ) | \big) \\
&\les & \sum_{|J| \leq  \lfloor \frac{|I|}{2} \rfloor, \; |K| \leq |I|} O \big(  | ( \Lie_{Z^J}  h )  | \cdot  | \derm  (  \Lie_{Z^K}  A ) | \big)  + \sum_{|J| \leq  \lfloor \frac{|I|}{2} \rfloor, \; |K| \leq |I|}  O \big(  | ( \Lie_{Z^K}  h )  | \cdot  | \derm  (  \Lie_{Z^J}  A ) | \big)
\eeaa
Using Lemma \ref{upgradedestimatesonh}, we have for $\ga \geq 3 \de $, and for $0 < \de \leq \frac{1}{4}$, in the exterior region $\overline{C}$, and for all $|J| \leq  \lfloor \frac{|I|}{2} \rfloor$\;, 
                 \beaa
 \notag
 |   \Lie_{Z^J} h   | &\leq&   C(q_0)   \cdot  c (\delta) \cdot c (\gamma) \cdot C(|I|) \cdot E (  \lfloor \frac{|I|}{2} \rfloor  +4)  \cdot \frac{\eps }{(1+t+|q|)^{1-   c (\gamma)  \cdot c (\delta)  \cdot c(|I|) \cdot E ( \lfloor \frac{|I|}{2} \rfloor+ 4)\cdot  \eps } }     \; , \\
      \eeaa
and using Lemma \ref{upgradedestimateonLiederivativesoffields}, we have for all $|J| \leq  \lfloor \frac{|I|}{2} \rfloor$\;, 
            \beaa
 \notag
&& |\derm  ( \Lie_{Z^J} A)  |    \\
\notag
&\leq&   C(q_0)   \cdot  c (\delta) \cdot c (\gamma) \cdot C(|I|) \cdot E ( \lfloor \frac{|I|}{2} \rfloor  + 4)  \cdot \frac{\eps }{(1+t+|q|)^{1-      c (\gamma)  \cdot c (\delta)  \cdot c(|I|) \cdot E ( \lfloor \frac{|I|}{2} \rfloor+ 4) \cdot \eps } \cdot (1+|q|)^{1+\gamma - 2\de }}    \; .
      \eeaa
  Injecting in the sum, we obatin
      \beaa
&&\sum_{|K| + |M| \leq |I|} O \big(  | ( \Lie_{Z^K}  h )  | \cdot  | \derm  (  \Lie_{Z^M}  A ) | \big) \\
&\les & \sum_{|K| \leq |I|}  \Big[ C(q_0)   \cdot  c (\delta) \cdot c (\gamma) \cdot C(|I|) \cdot E (  \lfloor \frac{|I|}{2} \rfloor  +4)  \cdot \frac{\eps \cdot   | \derm  (  \Lie_{Z^K}  A ) |  }{(1+t+|q|)^{1-   c (\gamma)  \cdot c (\delta)  \cdot c(|I|) \cdot E ( \lfloor \frac{|I|}{2} \rfloor+ 4)\cdot  \eps } }  \\
&& + C(q_0)   \cdot  c (\delta) \cdot c (\gamma) \cdot C(|I|) \cdot E ( \lfloor \frac{|I|}{2} \rfloor  + 4)  \cdot \frac{\eps  \cdot  |   \Lie_{Z^K} h   | }{(1+t+|q|)^{1-      c (\gamma)  \cdot c (\delta)  \cdot c(|I|) \cdot E ( \lfloor \frac{|I|}{2} \rfloor+ 4) \cdot \eps } \cdot (1+|q|)^{1+\gamma - 2\de }}  \Big] \; .
\eeaa

We can now estimate, by then, the following term
\beaa
 && \int_{\Sigma^{ext}_{\tau} }  C(q_0)   \cdot  c (\delta) \cdot c (\gamma)  \cdot E ( 4)   \cdot \sum_{|K| = |I| }  \Big(    \frac{\eps \cdot   | \pa \Lie_{Z^K} A_L |^2   }{(1+t+|q|)^{1-      c (\gamma)  \cdot c (\delta)   \cdot E ( 4) \cdot \eps } \cdot (1+|q|)^{2\gamma - 4\de }}       \Big)  \cdot w(q)   \\
 &\les & \int_{\Sigma^{ext}_{\tau} }C(q_0)   \cdot  c (\delta) \cdot c (\gamma)  \cdot E ( 4)  \cdot  E (\lfloor \frac{|I|}{2} \rfloor  ) \cdot  \\
 && \times   \Big[  \sum_{|K| \leq |I| }     \frac{\eps \cdot   | \rderm ( \Lie_{Z^K} A) |^2   }{(1+t+|q|)^{1-      c (\gamma)  \cdot c (\delta)   \cdot E ( 4) \cdot \eps } \cdot (1+|q|)^{2\gamma - 4\de }}         \\
 &&+  \sum_{|K| \leq |I|-1 }     \frac{\eps \cdot   | \derm \Lie_{Z^K} A |^2   }{(1+t+|q|)^{1-      c (\gamma)  \cdot c (\delta)   \cdot E ( 4) \cdot \eps } \cdot (1+|q|)^{2\gamma - 4\de }}       \Big]  \cdot w(q)    \\
&&+   \int_{\Sigma^{ext}_{\tau} } C(q_0)   \cdot  c (\delta) \cdot c (\gamma) \cdot C(|I|) \cdot E (  \lfloor \frac{|I|}{2} \rfloor  +4)  \\
&&  \times \sum_{|K| \leq |I|}  \Big[  \frac{\eps \cdot   | \derm  (  \Lie_{Z^K}  A ) |^2  }{(1+t+|q|)^{3-   c (\gamma)  \cdot c (\delta)  \cdot c(|I|) \cdot E ( \lfloor \frac{|I|}{2} \rfloor+ 4)\cdot  \eps }  \cdot (1+|q|)^{2\gamma - 4\de } }  \\
&& + \frac{\eps  \cdot  |   \Lie_{Z^K} h   |^2 }{(1+t+|q|)^{3-      c (\gamma)  \cdot c (\delta)  \cdot c(|I|) \cdot E ( \lfloor \frac{|I|}{2} \rfloor+ 4) \cdot \eps } \cdot (1+|q|)^{2+4\gamma - 6\de }}  \Big]  \cdot w(q)   \; .
\eeaa
We can use the fact that $E ( 4)  \cdot  E (\lfloor \frac{|I|}{2} \rfloor  ) \leq   E (\lfloor \frac{|I|}{2} \rfloor +4  )$ and write for $\eps$ small depending on $|I|$, on $\ga$ and on $\de$, and for $\ga \geq 3\de$, that 
\beaa
 && \int_{\Sigma^{ext}_{\tau} }  C(q_0)   \cdot  c (\delta) \cdot c (\gamma)  \cdot E ( 4)   \cdot \sum_{|K| = |I| }  \Big(    \frac{\eps \cdot   | \pa \Lie_{Z^K} A_L |^2   }{(1+t+|q|)^{1-      c (\gamma)  \cdot c (\delta)   \cdot E ( 4) \cdot \eps } \cdot (1+|q|)^{2\gamma - 4\de }}       \Big)  \cdot w(q)   \\
 &\les & \int_{\Sigma^{ext}_{\tau} }C(q_0)   \cdot  c (\delta) \cdot c (\gamma)    \cdot  E (\lfloor \frac{|I|}{2} \rfloor +4  )  \\
 && \times   \Big[  \sum_{|K| \leq |I| }     \frac{\eps \cdot   | \rderm ( \Lie_{Z^K} A) |^2   }{(1+t+|q|)^{1-      c (\gamma)  \cdot c (\delta)   \cdot E ( 4) \cdot \eps } \cdot (1+|q|)^{2\gamma - 4\de }}         \\
 &&+  \sum_{|K| \leq |I|-1 }     \frac{\eps \cdot   | \derm \Lie_{Z^K} A |^2   }{(1+t+|q|)^{1-      c (\gamma)  \cdot c (\delta)   \cdot E ( 4) \cdot \eps } \cdot (1+|q|)^{2\gamma - 4\de }}       \Big]  \cdot w(q)    \\
&&+   \int_{\Sigma^{ext}_{\tau} } C(q_0)   \cdot  c (\delta) \cdot c (\gamma) \cdot C(|I|) \cdot E (  \lfloor \frac{|I|}{2} \rfloor  +4)  \\
&&  \times \sum_{|K| \leq |I|}  \Big[  \frac{\eps \cdot   | \derm  (  \Lie_{Z^K}  A ) |^2  }{(1+t+|q|)^{2 } }  + \frac{\eps  \cdot  |   \Lie_{Z^K} h   |^2 }{(1+t+|q|)^{3-      c (\gamma)  \cdot c (\delta)  \cdot c(|I|) \cdot E ( \lfloor \frac{|I|}{2} \rfloor+ 4) \cdot \eps }  \cdot (1+|q|)^{2+4\gamma - 6\de }}  \Big]   \cdot w(q)   \; .
\eeaa

Using again the Hardy type inequality of Corollary \ref{HardytypeinequalityforintegralstartingatROm}, we obtain
\beaa
 && \int_{\Sigma^{ext}_{\tau} }  C(q_0)   \cdot  c (\delta) \cdot c (\gamma)  \cdot E ( 4)   \cdot \sum_{|K| = |I| }  \Big(    \frac{\eps \cdot   | \pa \Lie_{Z^K} A_L |^2   }{(1+t+|q|)^{1-      c (\gamma)  \cdot c (\delta)   \cdot E ( 4) \cdot \eps } \cdot (1+|q|)^{2\gamma - 4\de }}       \Big)   \cdot w(q)   \\
 &\les & \int_{\Sigma^{ext}_{\tau} }C(q_0)   \cdot  c (\delta) \cdot c (\gamma)    \cdot  E (\lfloor \frac{|I|}{2} \rfloor +4  )   \\
 && \times   \Big[  \sum_{|K| \leq |I| }     \frac{\eps \cdot   | \rderm ( \Lie_{Z^K} A) |^2   }{(1+t+|q|)^{1-      c (\gamma)  \cdot c (\delta)   \cdot E ( 4) \cdot \eps } \cdot (1+|q|)^{2\gamma - 4\de }}         \\
 &&+  \sum_{|K| \leq |I|-1 }     \frac{\eps \cdot   | \derm \Lie_{Z^K} A |^2   }{(1+t+|q|)^{1-      c (\gamma)  \cdot c (\delta)   \cdot E ( 4) \cdot \eps } \cdot (1+|q|)^{2\gamma - 4\de }}       \Big]  \cdot w(q)   \\
&&+   \int_{\Sigma^{ext}_{\tau} } C(q_0)   \cdot  c (\delta) \cdot c (\gamma) \cdot C(|I|) \cdot E (  \lfloor \frac{|I|}{2} \rfloor  +4)  \\
&&  \times \sum_{|K| \leq |I|}  \Big[  \frac{\eps \cdot   | \derm  (  \Lie_{Z^K}  A ) |^2  }{(1+t+|q|)^{2 } }  + \frac{\eps  \cdot  |   \derm ( \Lie_{Z^K} h )   |^2 }{(1+t+|q|)^{3-      c (\gamma)  \cdot c (\delta)  \cdot c(|I|) \cdot E ( \lfloor \frac{|I|}{2} \rfloor+ 4) \cdot \eps } \cdot (1+|q|) ^{4\gamma - 6\de  } }  \Big]   \cdot w(q)   \; .
\eeaa

Thus, decomposing again $h= h^1 + h^0$ and using the estimate of Lemma \ref{Liederivativesofsphericalsymmetricpart},
\beaa
 && \int_{\Sigma^{ext}_{\tau} } C(q_0)   \cdot  c (\delta) \cdot c (\gamma)  \cdot E ( 4)   \cdot \sum_{|K| = |I| }  \Big(    \frac{\eps \cdot   | \pa \Lie_{Z^K} A_L |^2   }{(1+t+|q|)^{1-      c (\gamma)  \cdot c (\delta)   \cdot E ( 4) \cdot \eps } \cdot (1+|q|)^{2\gamma - 4\de }}       \Big)  \cdot w(q)   \\
 &\les & \int_{\Sigma^{ext}_{\tau} }C(q_0)   \cdot  c (\delta) \cdot c (\gamma)    \cdot  E (\lfloor \frac{|I|}{2} \rfloor +4  ) \cdot  \\
 && \times   \Big[  \sum_{|K| \leq |I| }     \frac{\eps \cdot   | \rderm ( \Lie_{Z^K} A) |^2   }{(1+t+|q|)^{1-      c (\gamma)  \cdot c (\delta)   \cdot E ( 4) \cdot \eps } \cdot (1+|q|)^{2\gamma - 4\de }}         \\
 &&+  \sum_{|K| \leq |I|-1 }     \frac{\eps \cdot   | \derm \Lie_{Z^K} A |^2   }{(1+t+|q|)^{1-      c (\gamma)  \cdot c (\delta)   \cdot E ( 4) \cdot \eps } \cdot (1+|q|)^{2\gamma - 4\de }}       \Big]  \cdot w(q)   \\
&&+   \int_{\Sigma^{ext}_{\tau} } C(q_0)   \cdot  c (\delta) \cdot c (\gamma) \cdot C(|I|) \cdot E (  \lfloor \frac{|I|}{2} \rfloor  +4)  \\
&&  \times  \frac{\eps }{(1+t+|q|)^{2 } }  \cdot  \sum_{|K| \leq |I|}  \Big[  | \derm  (  \Lie_{Z^K}  A ) |^2 +  |   \derm ( \Lie_{Z^K} h^1 )   |^2  \Big]  \cdot w(q)   \\
 && + \int_{\Sigma^{ext}_{\tau} } C(q_0)   \cdot  c (\delta) \cdot c (\gamma) \cdot C(|I|) \cdot E (  \lfloor \frac{|I|}{2} \rfloor  +4)  \\
 \notag
 && \times  \frac{\eps^3 }{(1+t+|q|)^{5-      c (\gamma)  \cdot c (\delta)  \cdot c(|I|) \cdot E ( \lfloor \frac{|I|}{2} \rfloor+ 4) \cdot \eps }  \cdot (1+|q|) ^{4\gamma - 6\de  } }   \cdot w(q) \; .
\eeaa
The last term, that is the contribution of $h^0$, can be estimated in the same way as in the proof of Lemma \ref{estimateonthecontributionofhzerointhesourcetermsforthewaveequationontheYangMillspoential}, and hence, using the fact that for $\ga \geq 3\de$\, we have $2\ga - 6 \de \geq 0$\,, we obtain the desired result.
\end{proof}

\begin{lemma}\label{HardyinequalityontheestimateonthebadtermproductALtimesdermAusingbootstrapsoastoclosethegronwallinequalityonenergy}
For $\ga \geq 3\de$, and for $\eps$ small depending on $|I|$\;, on $\ga$ and on $\de$\;, we have
                                         \beaa
\notag
&&  \int_{\Sigma^{ext}_{\tau} } \frac{ (1+t)}{\eps}  \cdot \Big(  \sum_{|K| \leq |I|}  | \Lie_{Z^K}    \big(  A_L   \cdot     \derm A    \big)    |   \Big)^{2}  \cdot w(q)  \\
&\les&   \int_{\Sigma^{ext}_{\tau} } \Big[ C(q_0)   \cdot  c (\delta) \cdot c (\gamma)  \cdot C(|I|) \cdot  E ( \lfloor \frac{|I|}{2} \rfloor + 4)   \\
&& \times \sum_{|K| \leq |I| }  \Big(  \frac{ \eps  \cdot     |  \derm ( \Lie_{Z^K} A )  |^2  +  \eps \cdot |  \derm ( \Lie_{Z^K} h^1 )  |^2   }{ (1+t+|q|)   }   +    \frac{\eps \cdot   | \rderm ( \Lie_{Z^K} A ) |^2   }{(1+t+|q|)^{1-      c (\gamma)  \cdot c (\delta)   \cdot E ( 4) \cdot \eps } \cdot (1+|q|)}       \Big) \; \Big]  \cdot w(q)  \\
&& +   \int_{\Sigma^{ext}_{\tau} } \Big[  C(q_0)   \cdot  c (\delta) \cdot c (\gamma) \cdot C( |I|) \cdot E ( \lfloor \frac{|I|}{2} \rfloor + 4) \\
&& \times \sum_{ |K| \leq |I| -1}    \Big(    \frac{\eps  \cdot     |  \derm ( \Lie_{Z^K} A )  |^2 }{(1+t+|q|)^{1-      c (\gamma)  \cdot c (\delta)  \cdot c( \lfloor \frac{|I|-1}{2} \rfloor) \cdot E ( \lfloor \frac{|I|-1}{2} \rfloor+ 4) \cdot \eps } \cdot (1+|q|)^{2\gamma - 4\de }}      \Big)   \; \Big]   \cdot w(q)  \\
&&+   C(q_0)   \cdot  c (\delta) \cdot c (\gamma) \cdot C(|I|) \cdot E (  \lfloor \frac{|I|}{2} \rfloor  +4)  \cdot  \frac{\eps^3 }{(1+t)^{1-     c (\gamma)  \cdot c (\delta)  \cdot c(|I|) \cdot E ( \lfloor \frac{|I|}{2} \rfloor+ 4) \cdot \eps } }  \; .
\eeaa

\end{lemma}

\begin{proof}

We inject Lemma \ref{estimateonthetermthatcontainspartialderivatoveofALinthettimessquareofbadcomponenentsinsourcetermsforAthatcontainstheALcomponent} in Lemma \ref{HardytpyeestimatesontehtermscontainingALinthebadstructureofsourcestermsforAusefultoobtainenergyestimates} and we obtain using the fact that $E ( k_1) \leq  E (k_2  ) $ for all $k_1 \leq k_2$\;, that

                                         \beaa
\notag
&&  \int_{\Sigma^{ext}_{\tau} } \frac{ (1+t)}{\eps}  \cdot \Big(  \sum_{|K| \leq |I|}  | \Lie_{Z^K}    \big(  A_L   \cdot     \derm A    \big)    |   \Big)^{2}  \cdot w(q)  \\
&\les&   \int_{\Sigma^{ext}_{\tau} } \Big[ C(q_0)   \cdot  c (\delta) \cdot c (\gamma)  \cdot E ( \lfloor \frac{|I|}{2} \rfloor + 4)   \\
&& \times \sum_{|K| \leq |I| }  \Big(  \frac{ \eps  \cdot     |  \derm ( \Lie_{Z^K} A )  |^2   }{ (1+t+|q|)  \cdot  (1+|q|)^{2\gamma - 4\de}  }   +    \frac{\eps \cdot   | \rderm ( \Lie_{Z^K} A ) |^2   }{(1+t+|q|)^{1-      c (\gamma)  \cdot c (\delta)   \cdot E ( 4) \cdot \eps } \cdot (1+|q|)^{2\gamma - 4\de }}       \Big) \; \Big]  \cdot w(q)  \\
&& +   \int_{\Sigma^{ext}_{\tau} } \Big[  C(q_0)   \cdot  c (\delta) \cdot c (\gamma) \cdot C( \lfloor \frac{|I|-1}{2} \rfloor) \cdot E ( \lfloor \frac{|I|}{2} \rfloor + 4) \\
&& \times \sum_{ |K| \leq |I| -1}    \Big(    \frac{\eps  \cdot     |  \derm ( \Lie_{Z^K} A )  |^2 }{(1+t+|q|)^{1-      c (\gamma)  \cdot c (\delta)  \cdot c( \lfloor \frac{|I|-1}{2} \rfloor) \cdot E ( \lfloor \frac{|I|-1}{2} \rfloor+ 4) \cdot \eps } \cdot (1+|q|)^{2\gamma - 4\de }}      \Big)   \; \Big]   \cdot w(q)  \\
&&+   \int_{\Sigma^{ext}_{\tau} } C(q_0)   \cdot  c (\delta) \cdot c (\gamma) \cdot C(|I|) \cdot E (  \lfloor \frac{|I|}{2} \rfloor  +4)  \\
&&  \times  \frac{\eps }{(1+t+|q|)^{2 } }  \cdot  \sum_{|K| \leq |I|}  \Big[  | \derm  (  \Lie_{Z^K}  A ) |^2 +  |   \derm ( \Lie_{Z^K} h^1 )   |^2  \Big] \cdot w(q)   \\
&&+   C(q_0)   \cdot  c (\delta) \cdot c (\gamma) \cdot C(|I|) \cdot E (  \lfloor \frac{|I|}{2} \rfloor  +4)  \cdot  \frac{\eps^3 }{(1+t)^{1-     c (\gamma)  \cdot c (\delta)  \cdot c(|I|) \cdot E ( \lfloor \frac{|I|}{2} \rfloor+ 4) \cdot \eps } }  \; .
\eeaa

Thus, for $\ga \geq 3 \de$\;, and for $\eps$ small depending on $|I|$\;, on $\ga$ and on $\de$\,, we obtain the desired result.
\end{proof}

Now, we would like to put all of this together to get an energy estimate on the full components of $\derm A$. 

\begin{lemma}\label{finallemmafortheestimateontheenergyfortheEinsteinYangMillspoentialusingthestructureoftheequationsforbothgoodandbadcomponents}

For $\ga \geq 3 \de $\,, for $0 < \de \leq \frac{1}{4}$\,, and for $\eps$ small, enough depending on $q_0$\,, on $\ga$\;, on $\de$\;, on $|I|$ and on $\mu$\;, we have

  \beaa
   \notag
 &&     \int_{\Sigma^{ext}_{t_2} }  |\derm ( \Lie_{Z^I}  A )  |^2     \cdot w(q)  \cdot d^{3}x    + \int_{N_{t_1}^{t_2} }  T_{\hat{L} t}^{(\bf{g})}  \cdot  w(q) \cdot dv^{(\bf{m})}_N \\
 \notag
 &&+ \int_{t_1}^{t_2}  \int_{\Sigma^{ext}_{\tau} }     | \rderm  ( \Lie_{Z^I}  A ) |^2   \cdot  \frac{\widehat{w} (q)}{(1+|q|)} \cdot  d^{3}x  \cdot d\tau \\
   \notag
&\les&     \sum_{|K| \leq |I |}      \int_{t_1}^{t_2} \int_{\Sigma^{ext}_{\tau} }     \Big[       O \big(   C(q_0)   \cdot  c (\delta) \cdot c (\gamma) \cdot C(|I|) \cdot E (\lfloor \frac{|I|}{2} \rfloor + 5)  \cdot \frac{\eps   \cdot  | \derm ( \Lie_{Z^K} h^1 ) |^2  }{(1+t+|q|)}   \big) \\
   && + O \big(    C(q_0)   \cdot  c (\delta) \cdot c (\gamma) \cdot C(|I|) \cdot E ( \lfloor \frac{|I|}{2} \rfloor+5)  \cdot \frac{\eps  \cdot  |\derm ( \Lie_{Z^K} A ) |^2  }{(1+t+|q|) }  \big)   \; \Big]    \cdot w(q) \cdot dt \\  
&& +    \int_{t_1}^{t_2} \int_{\Sigma^{ext}_{\tau} } \Big[  C(q_0)   \cdot  c (\delta) \cdot c (\gamma) \cdot C( |I|) \cdot E ( \lfloor \frac{|I|}{2} \rfloor + 4) \\
&& \times \sum_{ |K| \leq |I| -1}    \Big(    \frac{\eps  \cdot     |  \derm ( \Lie_{Z^K} A )  |^2 }{(1+t+|q|)^{1-      c (\gamma)  \cdot c (\delta)  \cdot c( \lfloor \frac{|I|-1}{2} \rfloor) \cdot E ( \lfloor \frac{|I|-1}{2} \rfloor+ 4) \cdot \eps } \cdot (1+|q|)^{2\gamma - 4\de }}      \Big)   \; \Big]   \cdot w(q) \cdot dt \; .
\eeaa
 \beaa
\notag
  &&  +  C(q_0) \cdot c (\gamma) \cdot  C ( |I| ) \cdot E (  \lfloor \frac{|I|}{2} \rfloor  + 2) \\
  && \times   \int_{t_1}^{t_2}  \frac{\eps}{(1+t)^{1+2\de}}  \cdot \Big[ \;   \sum_{|J|\leq |I| + 2 }   \int_{\Sigma^{ext}_{t_1} }  |\derm  ( \Lie_{Z^J}  A_{e_{a}} )|^2     \cdot w(q)  \cdot d^{3}x \\
&& +   \sum_{|K| \leq |I |}   \int_{t_1}^{t}  \Big[    \;   \int_{\Sigma^{ext}_{\tau} }    \Big[       O \big(   C(q_0)   \cdot  c (\delta) \cdot c (\gamma) \cdot C(|I|) \cdot E ( \lfloor \frac{|I|}{2} \rfloor + 6)  \cdot \frac{\eps   \cdot  | \derm ( \Lie_{Z^K} h^1 ) |^2  }{(1+\tau+|q|)}   \big) \\
   && + O \big(    C(q_0)   \cdot  c (\delta) \cdot c (\gamma) \cdot C(|I|) \cdot E (   \lfloor \frac{|I|}{2} \rfloor +6)  \cdot \frac{\eps  \cdot  |\derm ( \Lie_{Z^K} A ) |^2  }{(1+\tau+|q|) }  \big)    \; \Big]    \cdot w(q)  \cdot d^{3}x  \Big]   \cdot d\tau  \\
      \notag
        &&  +     \int_{t_1}^{t}  \int_{\Sigma^{ext}_{\tau} }   \Big[     C(q_0)   \cdot  c (\delta) \cdot c (\gamma) \cdot C(|I|) \cdot E (\lfloor \frac{|I|}{2} \rfloor   +4)  \cdot \frac{\eps }{(1+\tau+|q|)^{1-   c (\gamma)  \cdot c (\delta)  \cdot c(|I|) \cdot E ( \lfloor \frac{|I|}{2} \rfloor + 4)\cdot  \eps } \cdot (1+|q|)^2} \\
&& \times    \sum_{  |K| \leq |I| -1 }  | \derm ( \Lie_{Z^K}  A )  |^2     \Big]  \cdot w(q)  \cdot d^{3}x  \cdot d\tau \\
&& +   \int_{t_1}^{t}   \int_{\Sigma^{ext}_{\tau} }   \sum_{  |K| \leq |I| }         \Big[   C(q_0)   \cdot C(|I|) \cdot E (  \lfloor \frac{|I|}{2} \rfloor  +5)  \cdot \frac{\eps \cdot  | \Lie_{Z^{K}} H_{L  L} |^2 }{(1+\tau+|q|)^{1- 2 \de } \cdot (1+|q|)^{4+2\gamma}}  \Big] \cdot w(q)  \cdot d^{3}x  \cdot d\tau  \; \Big]    \cdot dt   
\eeaa
                                           \beaa
   \notag
  &&  +     \int_{t_1}^{t_2}  \int_{\Sigma^{ext}_{\tau} }   \Big[       \frac{(1+ t )}{\eps} \cdot    \sum_{|K| \leq |I|  - 1}  | g^{\la\mu} \cdot \derm_{\la}   \derm_{\mu} (  \Lie_{Z^{K}}  A ) |^2 \\
&&+  C(q_0)   \cdot  c (\delta) \cdot c (\gamma) \cdot C(|I|)\cdot E (  \lfloor \frac{|I|}{2} \rfloor   +3)  \cdot \frac{\eps }{(1+t+|q|)^{1-   c (\gamma)  \cdot c (\delta)  \cdot c(|I|) \cdot E ( \lfloor \frac{|I|}{2} \rfloor  + 2)\cdot  \eps } \cdot (1+|q|)^2} \\
&& \times    \sum_{  |K| \leq |I| -1 }  | \derm ( \Lie_{Z^K}  A )  |^2   \Big] \cdot dt\\
&&+   C(q_0)   \cdot  c (\delta) \cdot c (\gamma) \cdot C(|I|) \cdot E (  \lfloor \frac{|I|}{2} \rfloor  +5)  \cdot  \eps^3  \cdot (1+t_2 )^{ c (\gamma)  \cdot c (\delta)  \cdot c(|I|) \cdot E ( \lfloor \frac{|I|}{2} \rfloor+ 5) \cdot \eps }   \; .
\eeaa

\end{lemma}

\begin{proof}

Using Lemma \ref{TheactualusefulstrzuctureofthesourcetermsforthewaveequationonpoentialAusingbootstrap} to see the structure of the full components, and then using Corollary \ref{HardytypeinequalityaplliedtothesoircetermsoftheEinsteinYangMillspoentialofgoodcomponents} to estimate the “good” terms, and using Lemma \ref{HardyinequalityontheestimateonthebadtermproductALtimesdermAusingbootstrapsoastoclosethegronwallinequalityonenergy} and Lemma \ref{HardyinequalityontheestimateonthebadtermproductA_atimesdermA_ausingbootstrapsoastoclosethegronwallinequalityonenergy} to estimate the “bad” terms, we obtain
  
  For $\ga \geq 3 \de $\,, for $0 < \de \leq \frac{1}{4}$\,, and for $\eps$ small enough, depending on $\ga$\,, $\de$ and on $|I|$\,, we have

                                                 \beaa
   \notag
&&  \int_{\Sigma^{ext}_{\tau} }   \frac{(1+ t )}{\eps} \cdot  |  \Lie_{Z^I}  g^{\la\mu} \derm_{\la}   \derm_{\mu}   A_{\cal U }  |^2    \cdot w(q)   \\
&\les&     \sum_{|K| \leq |I |}      \int_{\Sigma^{ext}_{\tau} }  \Big[         O \big(   C(q_0)   \cdot  c (\delta) \cdot c (\gamma) \cdot C(|I|) \cdot E (\lfloor \frac{|I|}{2} \rfloor + 5)  \cdot \frac{\eps   \cdot  | \derm ( \Lie_{Z^K} h^1 ) |^2  }{(1+t+|q|)}   \big) \\
   && + O \big(    C(q_0)   \cdot  c (\delta) \cdot c (\gamma) \cdot C(|I|) \cdot E ( \lfloor \frac{|I|}{2} \rfloor+5)  \cdot \frac{\eps  \cdot  |\derm ( \Lie_{Z^K} A ) |^2  }{(1+t+|q|) }  \big) \\
   &&  +  C(q_0)   \cdot  c (\delta) \cdot c (\gamma) \cdot C(|I|) \cdot E (\lfloor \frac{|I|}{2} \rfloor+ 4)  \cdot \frac{\eps   \cdot    | \rderm ( \Lie_{Z^K} A )  |^2  }{(1+t+|q|)^{1-      c (\gamma)  \cdot c (\delta)  \cdot c(|I|) \cdot E ( \lfloor \frac{|I|}{2} \rfloor+ 4) \cdot \eps } \cdot   (1+|q|) }    \\
      \notag
        &&    +     C(q_0)   \cdot  c (\delta) \cdot c (\gamma) \cdot C(|I|) \cdot E (\lfloor \frac{|I|}{2} \rfloor + 4)  \cdot \frac{\eps  \cdot   | \rderm  ( \Lie_{Z^K} h^1 ) |^2  }{(1+t+|q|)^{1-      c (\gamma)  \cdot c (\delta)  \cdot c(|I|) \cdot E (\lfloor \frac{|I|}{2} \rfloor+ 4) \cdot \eps } \cdot (1+|q|)^{2 }}   \Big]       \cdot w(q)   \\
       &&    +      C(q_0)   \cdot  c (\delta) \cdot c (\gamma) \cdot C(|I|) \cdot E ( \lfloor \frac{|I|}{2} \rfloor + 5)  \cdot  \frac{\eps^3 }{(1+t+|q|)^{3 -4\de -   c (\gamma)  \cdot c (\delta)  \cdot c(|I|) \cdot E ( \lfloor \frac{|I|}{2} \rfloor  + 5) \cdot \eps } } 
  \eeaa    
  
  \beaa  
&&  + \int_{\Sigma^{ext}_{\tau} } \Big[ C(q_0)   \cdot  c (\delta) \cdot c (\gamma)  \cdot C(|I|) \cdot  E ( \lfloor \frac{|I|}{2} \rfloor + 4)   \\
&& \times \sum_{|K| \leq |I| }  \Big(  \frac{ \eps  \cdot     |  \derm ( \Lie_{Z^K} A )  |^2  +  |  \derm ( \Lie_{Z^K} h^1 )  |^2   }{ (1+t+|q|)   }   +    \frac{\eps \cdot   | \rderm ( \Lie_{Z^K} A ) |^2   }{(1+t+|q|)^{1-      c (\gamma)  \cdot c (\delta)   \cdot E ( 4) \cdot \eps } \cdot (1+|q|)}       \Big) \; \Big]  \cdot w(q)  \\
&& +   \int_{\Sigma^{ext}_{\tau} } \Big[  C(q_0)   \cdot  c (\delta) \cdot c (\gamma) \cdot C(|I|) \cdot E ( \lfloor \frac{|I|}{2} \rfloor + 4) \\
&& \times \sum_{ |K| \leq |I| -1}    \Big(    \frac{\eps  \cdot     |  \derm ( \Lie_{Z^K} A )  |^2 }{(1+t+|q|)^{1-      c (\gamma)  \cdot c (\delta)  \cdot c( \lfloor \frac{|I|-1}{2} \rfloor) \cdot E ( \lfloor \frac{|I|-1}{2} \rfloor+ 4) \cdot \eps } \cdot (1+|q|)^{2\gamma - 4\de }}      \Big)   \; \Big]   \cdot w(q)   \\
&&+   C(q_0)   \cdot  c (\delta) \cdot c (\gamma) \cdot C(|I|) \cdot E (  \lfloor \frac{|I|}{2} \rfloor  +4)  \cdot  \frac{\eps^3 }{(1+t)^{1-     c (\gamma)  \cdot c (\delta)  \cdot c(|I|) \cdot E ( \lfloor \frac{|I|}{2} \rfloor+ 4) \cdot \eps } } \; .
\eeaa
 \beaa
\notag
  &&  +  C(q_0) \cdot c (\gamma) \cdot  C ( |I| ) \cdot E (  \lfloor \frac{|I|}{2} \rfloor  + 2) \\
  && \times  \frac{\eps}{(1+t)^{1+2\de}}  \cdot \Big[ \;   \sum_{|J|\leq |I| + 2 }  \int_{\Sigma^{ext}_{t_1} }  |\derm  ( \Lie_{Z^J}  A_{e_{a}} )|^2     \cdot w(q)  \cdot d^{3}x \\
&& +   \sum_{|K| \leq |I |}   \int_{t_1}^{t_2}  \Big[    \;   \int_{\Sigma^{ext}_{\tau} }    \Big[       O \big(   C(q_0)   \cdot  c (\delta) \cdot c (\gamma) \cdot C(|I|) \cdot E ( \lfloor \frac{|I|}{2} \rfloor + 6)  \cdot \frac{\eps   \cdot  | \derm ( \Lie_{Z^K} h^1 ) |^2  }{(1+\tau+|q|)}   \big) \\
   && + O \big(    C(q_0)   \cdot  c (\delta) \cdot c (\gamma) \cdot C(|I|) \cdot E (   \lfloor \frac{|I|}{2} \rfloor +6)  \cdot \frac{\eps  \cdot  |\derm ( \Lie_{Z^K} A ) |^2  }{(1+\tau+|q|) }  \big)    \; \Big]    \cdot w(q)  \cdot d^{3}x   \Big]  \cdot d\tau  \\
      \notag
        &&  +     \int_{t_1}^{t}  \int_{\Sigma^{ext}_{\tau} }   \Big[     C(q_0)   \cdot  c (\delta) \cdot c (\gamma) \cdot C(|I|) \cdot E (\lfloor \frac{|I|}{2} \rfloor   +4)  \cdot \frac{\eps }{(1+\tau+|q|)^{1-   c (\gamma)  \cdot c (\delta)  \cdot c(|I|) \cdot E (  \lfloor \frac{|I|}{2} \rfloor + 4)\cdot  \eps } \cdot (1+|q|)^2} \\
&& \times    \sum_{  |K| \leq |I| -1 } | \derm ( \Lie_{Z^K}  A )  |^2    \Big]  \cdot w(q)  \cdot d^{3}x  \cdot d\tau \\
&& +   \int_{t_1}^{t}   \int_{\Sigma^{ext}_{\tau} }   \sum_{  |K| \leq |I| }         \Big[   C(q_0)   \cdot C(|I|) \cdot E (  \lfloor \frac{|I|}{2} \rfloor  +5)  \cdot \frac{\eps \cdot  | \Lie_{Z^{K}} H_{L  L} |^2 }{(1+\tau+|q|)^{1- 2 \de } \cdot (1+|q|)^{4+2\gamma}}  \Big]\cdot w(q)  \cdot d^{3}x  \cdot d\tau  \\
 && +   C(q_0)   \cdot  c (\delta) \cdot c (\gamma) \cdot C(|I|) \cdot E ( \lfloor \frac{|I|}{2} \rfloor + 5)  \cdot  \frac{\eps^3 }{(1+t)^{1 -   c (\gamma)  \cdot c (\delta)  \cdot c(|I|) \cdot E ( \lfloor \frac{|I|}{2} \rfloor  + 5) \cdot \eps } }    \; \Big]     \; .   
\eeaa

Thus,

                                                 \beaa
   \notag
&&  \int_{\Sigma^{ext}_{\tau} }   \frac{(1+ t )}{\eps} \cdot  |  \Lie_{Z^I}  g^{\la\mu} \derm_{\la}   \derm_{\mu}   A_{\cal U }  |^2    \cdot w(q)   \\
&\les&     \sum_{|K| \leq |I |}     \int_{\Sigma^{ext}_{\tau} }     \Big[       O \big(   C(q_0)   \cdot  c (\delta) \cdot c (\gamma) \cdot C(|I|) \cdot E (\lfloor \frac{|I|}{2} \rfloor + 5)  \cdot \frac{\eps   \cdot  | \derm ( \Lie_{Z^K} h^1 ) |^2  }{(1+t+|q|)}   \big) \\
   && + O \big(    C(q_0)   \cdot  c (\delta) \cdot c (\gamma) \cdot C(|I|) \cdot E ( \lfloor \frac{|I|}{2} \rfloor+5)  \cdot \frac{\eps  \cdot  |\derm ( \Lie_{Z^K} A ) |^2  }{(1+t+|q|) }  \big) \\
   &&  +  C(q_0)   \cdot  c (\delta) \cdot c (\gamma) \cdot C(|I|) \cdot E (\lfloor \frac{|I|}{2} \rfloor+ 4)  \cdot \frac{\eps   \cdot    | \rderm ( \Lie_{Z^K} A )  |^2  }{(1+t+|q|)^{1-      c (\gamma)  \cdot c (\delta)  \cdot c(|I|) \cdot E ( \lfloor \frac{|I|}{2} \rfloor+ 4) \cdot \eps } \cdot   (1+|q|) }    \\
      \notag
        &&    +     C(q_0)   \cdot  c (\delta) \cdot c (\gamma) \cdot C(|I|) \cdot E (\lfloor \frac{|I|}{2} \rfloor + 4)  \cdot \frac{\eps  \cdot   | \rderm  ( \Lie_{Z^K} h^1 ) |^2  }{(1+t+|q|)^{1-      c (\gamma)  \cdot c (\delta)  \cdot c(|I|) \cdot E (\lfloor \frac{|I|}{2} \rfloor+ 4) \cdot \eps } \cdot (1+|q|)^{2 }}     \;  \Big]    \cdot w(q)  \\  
&& +   \int_{\Sigma^{ext}_{\tau} } \Big[  C(q_0)   \cdot  c (\delta) \cdot c (\gamma) \cdot C( |I|) \cdot E ( \lfloor \frac{|I|}{2} \rfloor + 4) \\
&& \times \sum_{ |K| \leq |I| -1}    \Big(    \frac{\eps  \cdot     |  \derm ( \Lie_{Z^K} A )  |^2 }{(1+t+|q|)^{1-      c (\gamma)  \cdot c (\delta)  \cdot c( \lfloor \frac{|I|-1}{2} \rfloor) \cdot E ( \lfloor \frac{|I|-1}{2} \rfloor+ 4) \cdot \eps } \cdot (1+|q|)^{2\gamma - 4\de }}      \Big)   \; \Big]   \cdot w(q)  \; .
\eeaa
 \beaa
\notag
  &&  +  C(q_0) \cdot c (\gamma) \cdot  C ( |I| ) \cdot E (  \lfloor \frac{|I|}{2} \rfloor  + 2) \\
  && \times  \frac{\eps}{(1+t)^{1+2\de}}  \cdot \Big[ \;  \sum_{|J|\leq |I| + 2 }   \int_{\Sigma^{ext}_{t_1} }  |\derm  ( \Lie_{Z^J}  A_{e_{a}} )|^2     \cdot w(q)  \cdot d^{3}x \\
&& +   \sum_{|K| \leq |I |}   \int_{t_1}^{t}  \Big[    \;   \int_{\Sigma^{ext}_{\tau} }    \Big[       O \big(   C(q_0)   \cdot  c (\delta) \cdot c (\gamma) \cdot C(|I|) \cdot E ( \lfloor \frac{|I|}{2} \rfloor + 6)  \cdot \frac{\eps   \cdot  | \derm ( \Lie_{Z^K} h^1 ) |^2  }{(1+\tau+|q|)}   \big) \\
   && + O \big(    C(q_0)   \cdot  c (\delta) \cdot c (\gamma) \cdot C(|I|) \cdot E (   \lfloor \frac{|I|}{2} \rfloor +6)  \cdot \frac{\eps  \cdot  |\derm ( \Lie_{Z^K} A ) |^2  }{(1+\tau+|q|) }  \big)    \; \Big]    \cdot w(q)  \cdot d^{3}x \Big]    \cdot d\tau  \\
      \notag
        &&  +     \int_{t_1}^{t}  \int_{\Sigma^{ext}_{\tau} }   \Big[     C(q_0)   \cdot  c (\delta) \cdot c (\gamma) \cdot C(|I|) \cdot E (\lfloor \frac{|I|}{2} \rfloor   +4)  \cdot \frac{\eps }{(1+\tau+|q|)^{1-   c (\gamma)  \cdot c (\delta)  \cdot c(|I|) \cdot E (  \lfloor \frac{|I|}{2} \rfloor + 4)\cdot  \eps } \cdot (1+|q|)^2} \\
&& \times    \sum_{  |K| \leq |I| -1 }  | \derm ( \Lie_{Z^K}  A )  |^2     \Big]  \cdot w(q)  \cdot d^{3}x  \cdot d\tau \\
&& +   \int_{t_1}^{t}   \int_{\Sigma^{ext}_{\tau} }   \sum_{  |K| \leq |I| }          \Big[   C(q_0)   \cdot C(|I|) \cdot E (  \lfloor \frac{|I|}{2} \rfloor  +5)  \cdot \frac{\eps \cdot  | \Lie_{Z^{K}} H_{L  L} |^2 }{(1+\tau +|q|)^{1- 2 \de } \cdot (1+|q|)^{4+2\gamma}}  \Big] \cdot w(q)  \cdot d^{3}x  \cdot d\tau  \; \Big]   \\
&&+   C(q_0)   \cdot  c (\delta) \cdot c (\gamma) \cdot C(|I|) \cdot E (  \lfloor \frac{|I|}{2} \rfloor  +5)  \cdot  \frac{\eps^3 }{(1+t)^{1-     c (\gamma)  \cdot c (\delta)  \cdot c(|I|) \cdot E ( \lfloor \frac{|I|}{2} \rfloor+ 5) \cdot \eps } }  \; .
\eeaa

Now, using the energy estimate of Lemma \ref{Theveryfinal }, we have dor $\eps$ small enough, depending on $q_0$, $\de$\,, $\ga$ and $\mu$\;, the following energy estimate holds for $\ga > 0$ and for $\Phi$ decaying sufficiently fast at spatial infinity,

  \beaa
   \notag
 &&     \int_{\Sigma^{ext}_{t} }  |\derm ( \Lie_{Z^I}  A )  |^2     \cdot w(q)  \cdot d^{3}x    + \int_{N_{t_1}^{t} }  T_{\hat{L} t}^{(\bf{g})}  \cdot  w(q) \cdot dv^{(\bf{m})}_N \\
 \notag
 &&+ \int_{t_1}^{t}  \int_{\Sigma^{ext}_{\tau} }     | \rderm  ( \Lie_{Z^I}  A ) |^2   \cdot  \frac{\widehat{w} (q)}{(1+|q|)} \cdot  d^{3}x  \cdot d\tau \\
  \notag
  &\les &       \int_{\Sigma^{ext}_{t_1} }  |\derm  ( \Lie_{Z^I}  A )|^2     \cdot w(q)  \cdot d^{3}x \\
    \notag
   &&  +  \int_{t_1}^{t}  \int_{\Sigma^{ext}_{\tau} }  \frac{(1+\tau )}{\eps} \cdot |  \Lie_{Z^I}  ( g^{\mu\a} \derm_{\mu } \derm_\a   A )   |^2    \cdot w(q) \cdot  d^{3}x  \cdot d\tau \\
     \notag
   &&  +  \int_{t_1}^{t}  \int_{\Sigma^{ext}_{\tau} }  \frac{(1+\tau )}{\eps} \cdot | \Lie_{Z^I}  ( g^{\la\mu} \derm_{\la}   \derm_{\mu}     A ) - g^{\la\mu}    \derm_{\la}   \derm_{\mu}  (  \Lie_{Z^I} A  ) |^2    \cdot w(q) \cdot  d^{3}x  \cdot d\tau \\
 \notag
&& + \int_{t_1}^{t}  \int_{\Sigma^{ext}_{\tau} }   C(q_0) \cdot   c (\delta) \cdot c (\gamma) \cdot E (  4)  \cdot \frac{ \eps }{ (1+\tau)  } \cdot | \derm  ( \Lie_{Z^I}  A) |^2     \cdot  w(q)\cdot  d^{3}x  \cdot d\tau \; . 
 \eeaa 
 
Using Lemma \ref{Hardytypeineeuqlityappliedtothecommutatortermplusatimeintegralformula}, we can estimate the commutator term for $\eps$ small enough, depending on $\ga$\,, $\de$ and $|I|$, and then we can inject in the energy estimate.

Now, since for $\eps$ small, enough depending on $\ga$, on $\de$, on $\mu$ and on $|I|$, we have
\bea
\frac{\eps   \cdot  w (q)  }{(1+t+|q|)^{1-      c (\gamma)  \cdot c (\delta)  \cdot c(|I|) \cdot E ( |I|+ 4) \cdot \eps } \cdot   (1+|q|) }  &\leq &  \frac{\eps \cdot \widehat{w} (q)}{(1+|q|)}  \;, 
  \eea
  we can then absorb the tangential derivatives to the left hand side and get,

  \beaa
   \notag
 &&     \int_{\Sigma^{ext}_{t_2} }  |\derm ( \Lie_{Z^I}  A )  |^2     \cdot w(q)  \cdot d^{3}x    + \int_{N_{t_1}^{t_2} }  T_{\hat{L} t}^{(\bf{g})}  \cdot  w(q) \cdot dv^{(\bf{m})}_N \\
 \notag
 &&+ \int_{t_1}^{t_2}  \int_{\Sigma^{ext}_{\tau} }     | \rderm  ( \Lie_{Z^I}  A ) |^2   \cdot  \frac{\widehat{w} (q)}{(1+|q|)} \cdot  d^{3}x  \cdot d\tau \\
   \notag
&\les&     \sum_{|K| \leq |I |}    \int_{t_1}^{t_2}   \int_{\Sigma^{ext}_{\tau} }     \Big[       O \big(   C(q_0)   \cdot  c (\delta) \cdot c (\gamma) \cdot C(|I|) \cdot E (\lfloor \frac{|I|}{2} \rfloor + 5)  \cdot \frac{\eps   \cdot  | \derm ( \Lie_{Z^K} h^1 ) |^2  }{(1+t+|q|)}   \big) \\
   && + O \big(    C(q_0)   \cdot  c (\delta) \cdot c (\gamma) \cdot C(|I|) \cdot E ( \lfloor \frac{|I|}{2} \rfloor+5)  \cdot \frac{\eps  \cdot  |\derm ( \Lie_{Z^K} A ) |^2  }{(1+t+|q|) }  \big)   \; \Big]    \cdot w(q) \cdot dt \\  
&& +    \int_{t_1}^{t_2} \int_{\Sigma^{ext}_{\tau} } \Big[  C(q_0)   \cdot  c (\delta) \cdot c (\gamma) \cdot C( |I|) \cdot E ( \lfloor \frac{|I|}{2} \rfloor + 4) \\
&& \times \sum_{ |K| \leq |I| -1}    \Big(    \frac{\eps  \cdot     |  \derm ( \Lie_{Z^K} A )  |^2 }{(1+t+|q|)^{1-      c (\gamma)  \cdot c (\delta)  \cdot c( \lfloor \frac{|I|-1}{2} \rfloor) \cdot E ( \lfloor \frac{|I|-1}{2} \rfloor+ 4) \cdot \eps } \cdot (1+|q|)^{2\gamma - 4\de }}      \Big)   \; \Big]   \cdot w(q) \cdot dt
\eeaa
 \beaa
\notag
  &&  +  C(q_0) \cdot c (\gamma) \cdot  C ( |I| ) \cdot E (  \lfloor \frac{|I|}{2} \rfloor  + 2) \\
  && \times  \int_{t_1}^{t_2}   \frac{1}{(1+t)^{1+2\de}}  \cdot \Big[ \;   \sum_{|J|\leq |I| + 2 }  \int_{\Sigma^{ext}_{t_1} }  |\derm  ( \Lie_{Z^J}  A_{e_{a}} )|^2     \cdot w(q)  \cdot d^{3}x \\
&& +   \sum_{|K| \leq |I |}   \int_{t_1}^{t}  \Big[    \;   \int_{\Sigma^{ext}_{\tau} }    \Big[       O \big(   C(q_0)   \cdot  c (\delta) \cdot c (\gamma) \cdot C(|I|) \cdot E ( \lfloor \frac{|I|}{2} \rfloor + 6)  \cdot \frac{\eps   \cdot  | \derm ( \Lie_{Z^K} h^1 ) |^2  }{(1+\tau+|q|)}   \big) \\
   && + O \big(    C(q_0)   \cdot  c (\delta) \cdot c (\gamma) \cdot C(|I|) \cdot E (   \lfloor \frac{|I|}{2} \rfloor +6)  \cdot \frac{\eps  \cdot  |\derm ( \Lie_{Z^K} A ) |^2  }{(1+\tau+|q|) }  \big)    \; \Big]    \cdot w(q)  \cdot d^{3}x   \Big]  \cdot d\tau  \\
      \notag
        &&  +     \int_{t_1}^{t}  \int_{\Sigma^{ext}_{\tau} }   \Big[     C(q_0)   \cdot  c (\delta) \cdot c (\gamma) \cdot C(|I|) \cdot E (  \lfloor \frac{|I|}{2} \rfloor  +4)  \cdot \frac{\eps }{(1+\tau+|q|)^{1-   c (\gamma)  \cdot c (\delta)  \cdot c(|I|) \cdot E (  \lfloor \frac{|I|}{2} \rfloor  + 4)\cdot  \eps } \cdot (1+|q|)^2} \\
&& \times    \sum_{  |K| \leq |I| -1 }  | \derm ( \Lie_{Z^K}  A )  |^2     \Big]  \cdot w(q)  \cdot d^{3}x  \cdot d\tau \\
&& +   \int_{t_1}^{t}   \int_{\Sigma^{ext}_{\tau} }   \sum_{  |K| \leq |I| }         \Big[   C(q_0)   \cdot C(|I|) \cdot E (  \lfloor \frac{|I|}{2} \rfloor  +5)  \cdot \frac{\eps \cdot  | \Lie_{Z^{K}} H_{L  L} |^2 }{(1+\tau +|q|)^{1- 2 \de } \cdot (1+|q|)^{4+2\gamma}}  \Big] \cdot w(q)  \cdot d^{3}x  \cdot d\tau  \; \Big]  \cdot dt 
\eeaa
                                           \beaa
   \notag
  &&  +    \int_{t_1}^{t_2}  \int_{\Sigma^{ext}_{\tau} }   \Big(       \frac{(1+ t )}{\eps} \cdot    \sum_{|K| \leq |I|  - 1}  | g^{\la\mu} \cdot \derm_{\la}   \derm_{\mu} (  \Lie_{Z^{K}}  A ) |^2 \\
&& +  C(q_0)   \cdot  c (\delta) \cdot c (\gamma) \cdot C(  \lfloor \frac{|I|}{2} \rfloor ) \cdot E (    \lfloor \frac{|I|}{2} \rfloor   +4)  \\
&& \times  \sum_{ |K| \leq |I| }  \Big(   \frac{\eps   }{(1+t+|q|)}  \, \cdot  |  \derm (   \Lie_{Z^K}  h^1 )   |^2  +  \frac{\eps }{(1+t+|q|) }  \, \cdot | \derm ( \Lie_{Z^K}  A)  |^2 \Big)   \\
&&+  C(q_0)   \cdot  c (\delta) \cdot c (\gamma) \cdot C(|I|) \cdot E (  \lfloor \frac{|I|}{2} \rfloor   +3 )  \cdot \frac{\eps }{(1+t+|q|)^{1-   c (\gamma)  \cdot c (\delta)  \cdot c(|I|) \cdot E (  \lfloor \frac{|I|}{2} \rfloor + 2)\cdot  \eps } \cdot (1+|q|)^2} \\
&& \times    \sum_{  |K| \leq |I| -1 }  | \derm ( \Lie_{Z^K}  A )  |^2   \Big)  \cdot dt  \\
&&+   C(q_0)   \cdot  c (\delta) \cdot c (\gamma) \cdot C(|I|) \cdot E (  \lfloor \frac{|I|}{2} \rfloor  +5)  \cdot  \eps^3  \cdot (1+t_2 )^{ c (\gamma)  \cdot c (\delta)  \cdot c(|I|) \cdot E ( \lfloor \frac{|I|}{2} \rfloor+ 5) \cdot \eps }   \; .
\eeaa

Combining the terms, we obtain the result.

\end{proof}

\begin{lemma}\label{estimateonthetermthatcontainsHLLintheenergywithestimatesinvolvingonlyhonesothattoestiamteallbyenergy}

For $\ga \geq 3 \de $\,, for $0 < \de \leq \frac{1}{4}$\,, and for $\eps$ small, enough depending on $q_0$\,, on $\ga$\;, on $\de$\;, and on $|I|$\;,
  
 \bea
     \notag
&&   \int_{t_1}^{t_2}   \int_{\Sigma^{ext}_{\tau} }   \sum_{  |K| \leq |I| }     \Big[   C(q_0)   \cdot C(|I|) \cdot E (  \lfloor \frac{|I|}{2} \rfloor  +5)  \cdot \frac{\eps \cdot  | \Lie_{Z^{K}} H_{L  L} |^2 }{(1+\tau+|q|)^{1- 2 \de } \cdot (1+|q|)^{4+2\gamma}}  \Big] \cdot w(q)  \cdot d^{3}x  \cdot d\tau     \\
     \notag
&\leq&   \int_{t_1}^{t_2}    \int_{\Sigma^{ext}_{\tau} }   \sum_{  |K| \leq |I| }       \Big[  C(q_0)   \cdot  c (\delta) \cdot c (\gamma) \cdot C(|I|) \cdot E (  \lfloor \frac{|I|}{2} \rfloor  +5)  \\
     \notag
&& \times  \frac{\eps \cdot  | \rderm ( \Lie_{Z^{K}} h^1 )  |^2 }{(1+\tau+|q|)^{1-  2\de} \cdot (1+|q|)^{2+2\ga}} \Big] \cdot w(q)  \cdot d^{3}x  \cdot d\tau \\
     \notag
&& +  \int_{t_1}^{t_2}    \int_{\Sigma^{ext}_{\tau} }   \sum_{  |K| \leq |I| }       \Big[  C(q_0)   \cdot  c (\delta) \cdot c (\gamma) \cdot C(|I|) \cdot E (  \lfloor \frac{|I|}{2} \rfloor  +5)  \\
     \notag
&& \times  \frac{\eps \cdot  | \derm ( \Lie_{Z^{K}} h^1 )  |^2 }{(1+\tau+|q|)^{3- 2\de -  c (\gamma)  \cdot c (\delta)  \cdot c(|I|) \cdot E (   \lfloor \frac{|I|}{2} \rfloor + 5)\cdot  \eps } \cdot (1+|q|)^{2+2\ga}} \Big] \cdot w(q)  \cdot d^{3}x  \cdot d\tau   \\
 \notag
&& + C(q_0)   \cdot  c (\delta) \cdot c (\gamma) \cdot C(|I|) \cdot E ( \lfloor \frac{|I|}{2} \rfloor + 5)  \cdot  \frac{\eps^3 }{(1+t_2+|q|)^{1 -   c (\gamma)  \cdot c (\delta)  \cdot c(|I|) \cdot E ( \lfloor \frac{|I|}{2} \rfloor  + 4) \cdot \eps  } }   \; . 
\eea

\end{lemma}

\begin{proof}

Using the Hardy-type inequality established in Corollary \ref{HardytypeinequalityforintegralstartingatROm}, we get
\beaa
&&   \int_{t_1}^{t_2}   \int_{\Sigma^{ext}_{\tau} }   \sum_{  |K| \leq |I| }          \Big[   C(q_0)   \cdot C(|I|) \cdot E (  \lfloor \frac{|I|}{2} \rfloor  +5)  \cdot \frac{\eps \cdot  | \Lie_{Z^{K}} H_{L  L} |^2 }{(1+\tau+|q|)^{1- 2 \de } \cdot (1+|q|)^{4+2\gamma}}  \Big]  \cdot w(q)  \cdot d^{3}x  \cdot d\tau  \;   \\
&\les&    \int_{t_1}^{t_2}    \int_{\Sigma^{ext}_{\tau} }   \sum_{  |K| \leq |I| }         \Big[   C(q_0)   \cdot C(|I|) \cdot E (  \lfloor \frac{|I|}{2} \rfloor  +5)  \cdot \frac{\eps \cdot  | \derm ( \Lie_{Z^{K}} H_{L  L}) |^2 }{(1+\tau+|q|)^{1- 2 \de } \cdot (1+|q|)^{2+2\gamma}}  \Big] \cdot w(q)  \cdot d^{3}x  \cdot d\tau    \; .
\eeaa

Now, based on Lemma \ref{wavecoordinatesestimateonLiederivativesZonmetric}, we have
\beaa
| \derm ( \Lie_{Z^J}  H_{\cal T L} )  | &\les&\sum_{|K| \leq |I| }  |  \rderm ( \Lie_{Z^K} H ) |  + \sum_{|K|+ |M| \leq |I|}  O (| \Lie_{Z^K} H| \cdot |\derm ( \Lie_{Z^M} H ) | ) \; .
\eeaa

To estimate the term $\sum_{|K|+ |M| \leq |I|}  O (| \Lie_{Z^K} H| \cdot |\derm ( \Lie_{Z^M} H ) | )$\;, we decompose the sum as follows
\beaa
&& \sum_{|K|+ |M| \leq |I|}  O (| \Lie_{Z^K} H| \cdot |\derm ( \Lie_{Z^M} H ) | )  \\
&\leq& \sum_{ |K| \leq   \lfloor \frac{|I|}{2} \rfloor   ,\; |J| \leq |I| }   | \Lie_{Z^{J}} H |\, \cdot | \derm ( \Lie_{Z^K}  H )  | + \sum_{ |J| \leq   \lfloor \frac{|I|}{2} \rfloor   ,\; |K| \leq |I| }  | \Lie_{Z^{J}} H |\, \cdot | \derm ( \Lie_{Z^K} H )  | \;.
\eeaa 

Based on Lemma  \ref{upgradedestimateonBIGH}, we have for all $|K| \leq   \lfloor \frac{|I|}{2} \rfloor $\,,
                 \beaa
 \notag
 |    \Lie_{Z^K}  H   | &\leq&   C(q_0)   \cdot  c (\delta) \cdot c (\gamma) \cdot C( \lfloor \frac{|I|}{2} \rfloor )  \cdot E (   \lfloor \frac{|I|}{2} \rfloor + 4)  \cdot \frac{\eps }{(1+t+|q|)^{1-   c (\gamma)  \cdot c (\delta)  \cdot c( \lfloor \frac{|I|}{2} \rfloor ) \cdot E (  \lfloor \frac{|I|}{2} \rfloor + 4)\cdot  \eps } }     \; , \\
      \eeaa
      and
                    \beaa
 \notag
 && |\derm  ( \Lie_{Z^K} H )   | \\
  &\leq&   C(q_0)   \cdot  c (\delta) \cdot c (\gamma) \cdot C( \lfloor \frac{|I|}{2} \rfloor ) \cdot E (   \lfloor \frac{|I|}{2} \rfloor  +4)  \cdot \frac{\eps }{(1+t+|q|)^{1-   c (\gamma)  \cdot c (\delta)  \cdot c(|I|) \cdot E (  \lfloor \frac{|I|}{2} \rfloor + 4)\cdot  \eps } \cdot (1+|q|)}     \; . 
      \eeaa

Thus,
\beaa
&& \sum_{|K|+ |M| \leq |I|}  O (| \Lie_{Z^K} H| \cdot |\derm ( \Lie_{Z^M} H ) | )  \\
&\leq& \sum_{ |K| \leq |I| } C(q_0)   \cdot  c (\delta) \cdot c (\gamma) \cdot C(|I|) \cdot E (   \lfloor \frac{|I|}{2} \rfloor  +4)  \cdot \frac{\eps \cdot  |    \Lie_{Z^K}  H   | }{(1+t+|q|)^{1-   c (\gamma)  \cdot c (\delta)  \cdot c(|I|) \cdot E (  \lfloor \frac{|I|}{2} \rfloor  + 4)\cdot  \eps } \cdot (1+|q|)}  \\
&& +  \sum_{ |K| \leq  |I| } C(q_0)   \cdot  c (\delta) \cdot c (\gamma) \cdot C(|I|) \cdot E (   \lfloor \frac{|I|}{2} \rfloor  +4)  \cdot \frac{\eps \cdot  |\derm  ( \Lie_{Z^K} H )   | }{(1+t+|q|)^{1-   c (\gamma)  \cdot c (\delta)  \cdot c(|I|) \cdot E (  \lfloor \frac{|I|}{2} \rfloor  + 4)\cdot  \eps } } \; .
\eeaa 
      
      Injecting, we obtain
      \beaa
&& \int_{t_1}^{t_2}   \int_{\Sigma^{ext}_{\tau} }   \sum_{  |K| \leq |I| }        \Big[   C(q_0)   \cdot C(|I|) \cdot E (  \lfloor \frac{|I|}{2} \rfloor  +5)  \cdot \frac{\eps \cdot  | \Lie_{Z^{K}} H_{L  L} |^2 }{(1+\tau+|q|)^{1- 2 \de } \cdot (1+|q|)^{4+2\gamma}}  \Big] \cdot w(q)  \cdot d^{3}x  \cdot d\tau     \\
&\leq&   \int_{t_1}^{t_2}    \int_{\Sigma^{ext}_{\tau} }   \sum_{  |K| \leq |I| }       \Big[  C(q_0)   \cdot  c (\delta) \cdot c (\gamma) \cdot C(|I|) \cdot E (  \lfloor \frac{|I|}{2} \rfloor  +5)  \\
&& \times  \frac{\eps \cdot  | \rderm ( \Lie_{Z^{K}} H )  |^2 }{(1+\tau+|q|)^{1-  2\de } \cdot (1+|q|)^{2+2\gamma}} \Big] \cdot w(q)  \cdot d^{3}x  \cdot d\tau \\
&& +  \int_{t_1}^{t_2}    \int_{\Sigma^{ext}_{\tau} }   \sum_{  |K| \leq |I| }       \Big[  C(q_0)   \cdot  c (\delta) \cdot c (\gamma) \cdot C(|I|) \cdot E (  \lfloor \frac{|I|}{2} \rfloor  +5)  \\
&& \times  \frac{\eps \cdot  | \derm ( \Lie_{Z^{K}} H )  |^2 }{(1+\tau+|q|)^{3   -  2\de - c (\gamma)  \cdot c (\delta)  \cdot c(|I|) \cdot E (   \lfloor \frac{|I|}{2} \rfloor + 5)\cdot  \eps } \cdot (1+|q|)^{2+2\gamma}} \Big] \cdot w(q)  \cdot d^{3}x  \cdot d\tau  \\
&& +  \int_{t_1}^{t_2}    \int_{\Sigma^{ext}_{\tau} }   \sum_{  |K| \leq |I| }       \Big[  C(q_0)   \cdot  c (\delta) \cdot c (\gamma) \cdot C(|I|) \cdot E (  \lfloor \frac{|I|}{2} \rfloor  +5)  \\
&& \times  \frac{\eps \cdot  | \Lie_{Z^{K}} H   |^2 }{(1+\tau+|q|)^{3-  2\de -   c (\gamma)  \cdot c (\delta)  \cdot c(|I|) \cdot E (   \lfloor \frac{|I|}{2} \rfloor + 5)\cdot  \eps } \cdot (1+|q|)^{4+2\ga}} \Big] \cdot w(q)  \cdot d^{3}x  \cdot d\tau   \; . 
\eeaa

Applying the Hardy-type inequality again, established in Corollary \ref{HardytypeinequalityforintegralstartingatROm}, we get

     \bea\label{ThetermintheenergythatcontainsbigHthatwewanttoconvertintosmallhandlaterintohone}
    \notag
&& \int_{t_1}^{t_2}   \int_{\Sigma^{ext}_{\tau} }   \sum_{  |K| \leq |I| }       \Big[   C(q_0)   \cdot C(|I|) \cdot E (  \lfloor \frac{|I|}{2} \rfloor  +5)  \cdot \frac{\eps \cdot  | \Lie_{Z^{K}} H_{L  L} |^2 }{(1+\tau+|q|)^{1- 2 \de } \cdot (1+|q|)^{4+2\gamma}}  \Big] \cdot w(q)  \cdot d^{3}x  \cdot d\tau     \\
     \notag
&\leq&   \int_{t_1}^{t_2}    \int_{\Sigma^{ext}_{\tau} }   \sum_{  |K| \leq |I| }       \Big[  C(q_0)   \cdot  c (\delta) \cdot c (\gamma) \cdot C(|I|) \cdot E (  \lfloor \frac{|I|}{2} \rfloor  +5)  \\
     \notag
&& \times  \frac{\eps \cdot  | \rderm ( \Lie_{Z^{K}} H )  |^2 }{(1+\tau+|q|)^{1-  2\de } \cdot (1+|q|)^{2+2\ga}} \Big] \cdot w(q)  \cdot d^{3}x  \cdot d\tau \\
     \notag
&& +  \int_{t_1}^{t_2}    \int_{\Sigma^{ext}_{\tau} }   \sum_{  |K| \leq |I| }       \Big[  C(q_0)   \cdot  c (\delta) \cdot c (\gamma) \cdot C(|I|) \cdot E (  \lfloor \frac{|I|}{2} \rfloor  +5)  \\
     \notag
&& \times  \frac{\eps \cdot  | \derm ( \Lie_{Z^{K}} H )  |^2 }{(1+\tau+|q|)^{3 -2\de-   c (\gamma)  \cdot c (\delta)  \cdot c(|I|) \cdot E (   \lfloor \frac{|I|}{2} \rfloor + 5)\cdot  \eps } \cdot (1+|q|)^{2+2\ga}} \Big] \cdot w(q)  \cdot d^{3}x  \cdot d\tau  \; . \\
\eea

\textbf{Converting estimates on $H$ into estimates on $h$}:\\

From Lemma \ref{linkbetweenbigHandsamllh}, and by lowering indices with respect to the metric $m$\;, we have
  \beaa
H_{\mu\nu}=-h_{\mu\nu}+ O_{\mu\nu}(h^2) \; .
\eeaa
This gives that for all $|I|$\;,
  \beaa
 \Lie_{Z^I}  H_{\mu\nu} &=& -   \Lie_{Z^I}  h_{\mu\nu}+   \sum_{|J| + |K| \leq |I|} O_{\mu\nu}(  \Lie_{Z^J}  h  \cdot  \Lie_{Z^K} h ) \\
 \eeaa
 and therefore,
   \beaa
\derm_\a ( \Lie_{Z^I}  H_{\mu\nu} )  &=& -  \derm_a \Lie_{Z^I}  h_{\mu\nu}+   \sum_{|J| + |K| \leq |I|} O_{\mu\nu}( \derm_\a ( \Lie_{Z^J}  h ) \cdot  \Lie_{Z^K} h ) \; ,
\eeaa
and
   \beaa
\rderm_\a ( \Lie_{Z^I}  H_{\mu\nu} )  &=& -  \rderm_a \Lie_{Z^I}  h_{\mu\nu}+   \sum_{|J| + |K| \leq |I|} O_{\mu\nu}( \rderm_\a ( \Lie_{Z^J}  h ) \cdot  \Lie_{Z^K} h ) \; .
\eeaa

As a result
   \beaa
| \derm_\a ( \Lie_{Z^I}  H ) |  &\leq& | \derm_\a ( \Lie_{Z^I}  h ) | +   \sum_{|J| + |K| \leq |I|}  | \derm_\a ( \Lie_{Z^J}  h ) | \cdot  | \Lie_{Z^K} h  | \; ,
\eeaa
and
   \beaa
| \rderm_\a ( \Lie_{Z^I}  H ) |  &\leq& | \rderm_\a ( \Lie_{Z^I}  h ) | +   \sum_{|J| + |K| \leq |I|}  | \rderm_\a ( \Lie_{Z^J}  h ) | \cdot  | \Lie_{Z^K} h  | \; .
\eeaa

Finally, we obatin
\beaa
&& | \derm ( \Lie_{Z^I}  H ) |  \\
&\les&  | \derm ( \Lie_{Z^I}  h ) | +   \sum_{|J| + |K| \leq |I|}  | \derm ( \Lie_{Z^J}  h ) | \cdot  | \Lie_{Z^K} h  | \\
&\les &  | \derm ( \Lie_{Z^I}  h ) |  +    \sum_{|J| \leq  \lfloor \frac{|I|}{2} \rfloor  \; , \; |K|\leq |I|}   | \derm ( \Lie_{Z^J}  h ) | \cdot  | \Lie_{Z^K} h  |   +   \sum_{|J| \leq |I| \; , \; |K| \leq  \lfloor \frac{|I|}{2} \rfloor  }   | \derm ( \Lie_{Z^J}  h ) | \cdot  | \Lie_{Z^K} h  |  \; .
  \eeaa
  and
\beaa
&& | \rderm ( \Lie_{Z^I}  H ) |  \\
&\les &  | \rderm ( \Lie_{Z^I}  h ) |  +    \sum_{|J| \leq  \lfloor \frac{|I|}{2} \rfloor  \; , \; |K|\leq |I|}   | \rderm ( \Lie_{Z^J}  h ) | \cdot  | \Lie_{Z^K} h  |   +   \sum_{|J| \leq |I| \; , \; |K| \leq  \lfloor \frac{|I|}{2} \rfloor  }   | \rderm ( \Lie_{Z^J}  h ) | \cdot  | \Lie_{Z^K} h  |  \; .
  \eeaa  
    
However, based on Lemma \ref{upgradedestimatesonh}, we have for $\ga \geq 3 \de $, and $0 < \de \leq \frac{1}{4}$, we have in the exterior region $\overline{C} \subset \{ (t, x) \; | \;  q \geq q_0 \} $, for all $|I|$, 

              \beaa
 \notag
 |\derm  ( \Lie_{Z^I} h )  |  &\leq&   C(q_0)   \cdot  c (\delta) \cdot c (\gamma) \cdot C(|I|) \cdot E (  |I|  +4)  \cdot \frac{\eps }{(1+t+|q|)^{1-   c (\gamma)  \cdot c (\delta)  \cdot c(|I|) \cdot E ( |I|+ 4)\cdot  \eps } \cdot (1+|q|)}     \; , \\
 \notag
 |   \Lie_{Z^I} h   | &\leq&   C(q_0)   \cdot  c (\delta) \cdot c (\gamma) \cdot C(|I|) \cdot E (  |I|  +4)  \cdot \frac{\eps }{(1+t+|q|)^{1-   c (\gamma)  \cdot c (\delta)  \cdot c(|I|) \cdot E ( |I|+ 4)\cdot  \eps } }     \; . 
      \eeaa

Hence,
\beaa
&& | \derm ( \Lie_{Z^I}  H ) |  \\
\notag
&\les &  | \derm ( \Lie_{Z^I}  h ) | \\
&& +    \sum_{  |K|\leq |I|}   C(q_0)   \cdot  c (\delta) \cdot c (\gamma) \cdot   C (  \lfloor \frac{|I|}{2} \rfloor   ) \cdot E ( \lfloor \frac{|I|}{2} \rfloor   + 4)   \cdot \frac{\eps  \cdot | \Lie_{Z^K} h  |  }{(1+t+|q|)^{1-   c (\gamma)  \cdot c (\delta)  \cdot c(|I|) \cdot E (  \lfloor \frac{|I|}{2} \rfloor   + 4)\cdot  \eps } \cdot (1+|q|)} \\
&&   +   \sum_{|J| \leq |I|  }    C(q_0)   \cdot  c (\delta) \cdot c (\gamma) \cdot  C ( \lfloor \frac{|I|}{2} \rfloor  ) \cdot E ( \lfloor \frac{|I|}{2} \rfloor  +4) \cdot  \frac{\eps \cdot  | \derm ( \Lie_{Z^J}  h ) | }{ (1+ t + | q | )^{1-   c (\gamma)  \cdot c (\delta)  \cdot c(|I|) \cdot E (  \lfloor \frac{|I|}{2} \rfloor   + 4)\cdot  \eps } } \\
&\les & \sum_{|J| \leq |I|  }   \Big( 1+  C(q_0)   \cdot  c (\delta) \cdot c (\gamma) \cdot C ( \lfloor \frac{|I|}{2} \rfloor  ) \cdot E ( \lfloor \frac{|I|}{2} \rfloor  + 4) \cdot \eps \Big) \cdot | \derm ( \Lie_{Z^J}  h ) |  \\
\notag
&& +    \sum_{  |K|\leq |I|}     C(q_0)   \cdot  c (\delta) \cdot c (\gamma) \cdot  C (  \lfloor \frac{|I|}{2} \rfloor   ) \cdot E ( \lfloor \frac{|I|}{2} \rfloor   + 4)   \cdot \frac{\eps  \cdot | \Lie_{Z^K} h  |  }{(1+t+|q|)^{1-   c (\gamma)  \cdot c (\delta)  \cdot c(|I|) \cdot E (  \lfloor \frac{|I|}{2} \rfloor   + 4)\cdot  \eps } \cdot (1+|q|) }  \;.\\
  \eeaa
  Consequently,
  \beaa
&& | \derm ( \Lie_{Z^I}  H ) |^2  \\
\notag
&\les & \sum_{|J| \leq |I|  }   \Big( 1+  C(q_0)   \cdot  c (\delta) \cdot c (\gamma) \cdot C ( \lfloor \frac{|I|}{2} \rfloor  ) \cdot E ( \lfloor \frac{|I|}{2} \rfloor   +4) \cdot \eps \Big)^2 \cdot | \derm ( \Lie_{Z^J}  h ) |^2  \\
\notag
&& +    \sum_{  |K|\leq |I|}    C(q_0)   \cdot  c (\delta) \cdot c (\gamma) \cdot  C (  \lfloor \frac{|I|}{2} \rfloor   ) \cdot E ( \lfloor \frac{|I|}{2} \rfloor    +4)   \cdot \frac{\eps  \cdot | \Lie_{Z^K} h  |^2 }{(1+t+|q|)^{2-   c (\gamma)  \cdot c (\delta)  \cdot c(|I|) \cdot E (  \lfloor \frac{|I|}{2} \rfloor   + 4)\cdot  \eps } \cdot (1+|q|)^2 }  \;.\\
  \eeaa
  and similarly, since based on Lemma \ref{upgradedestimatesonh}, we have
   \beaa
 \notag
&&  |  \rderm ( \Lie_{Z^I} h ) (t,x)  |  \\
\notag
&\leq&   C(q_0)   \cdot  c (\delta) \cdot c (\gamma) \cdot C(|I|) \cdot E (  |I|  +5)  \cdot \frac{\eps }{(1+t+|q|)^{2-   c (\gamma)  \cdot c (\delta)  \cdot c(|I|) \cdot E ( |I|+ 5)\cdot  \eps } }     \; , 
\eeaa
we therefore get
    \beaa
&& | \rderm ( \Lie_{Z^I}  H ) |^2  \\
\notag
&\les & \sum_{|J| \leq |I|  }   \Big( 1+  C(q_0)   \cdot  c (\delta) \cdot c (\gamma) \cdot C ( \lfloor \frac{|I|}{2} \rfloor  ) \cdot E ( \lfloor \frac{|I|}{2} \rfloor   +4) \cdot \eps \Big)^2 \cdot | \rderm ( \Lie_{Z^J}  h ) |^2  \\
\notag
&& +    \sum_{  |K|\leq |I|}    C(q_0)   \cdot  c (\delta) \cdot c (\gamma) \cdot  C (  \lfloor \frac{|I|}{2} \rfloor   ) \cdot E ( \lfloor \frac{|I|}{2} \rfloor    +4)   \cdot \frac{\eps  \cdot | \Lie_{Z^K} h  |^2 }{(1+t+|q|)^{4-   c (\gamma)  \cdot c (\delta)  \cdot c(|I|) \cdot E (  \lfloor \frac{|I|}{2} \rfloor   + 4)\cdot  \eps } \cdot (1+|q|)^2 }  \;.\\
  \eeaa

Injecting in \eqref{ThetermintheenergythatcontainsbigHthatwewanttoconvertintosmallhandlaterintohone}, and using the Hardy-type inequality of Corollary \ref{HardytypeinequalityforintegralstartingatROm}, we obtain

      \bea\label{ThetermintheenergythatcontainsbigHtestimatedbysmallhbutstillneedstoconvertintosmallhone}
     \notag
&&  \int_{t_1}^{t_2}   \int_{\Sigma^{ext}_{\tau} }   \sum_{  |K| \leq |I| }       \Big[   C(q_0)   \cdot C(|I|) \cdot E (  \lfloor \frac{|I|}{2} \rfloor  +5)  \cdot \frac{\eps \cdot  | \Lie_{Z^{K}} H_{L  L} |^2 }{(1+\tau+|q|)^{1- 2 \de } \cdot (1+|q|)^{4+2\gamma}}  \Big]\cdot w(q)  \cdot d^{3}x  \cdot d\tau  \; \Big]    \cdot dt   \\
     \notag
&\leq&   \int_{t_1}^{t_2}    \int_{\Sigma^{ext}_{\tau} }   \sum_{  |K| \leq |I| }       \Big[  C(q_0)   \cdot  c (\delta) \cdot c (\gamma) \cdot C(|I|) \cdot E (  \lfloor \frac{|I|}{2} \rfloor  +5)  \\
     \notag
&& \times  \frac{\eps \cdot  | \rderm ( \Lie_{Z^{K}} h )  |^2 }{(1+\tau+|q|)^{1- 2\de } \cdot (1+|q|)^{2+2\ga}} \Big] \cdot w(q)  \cdot d^{3}x  \cdot d\tau \\
     \notag
&& +  \int_{t_1}^{t_2}    \int_{\Sigma^{ext}_{\tau} }   \sum_{  |K| \leq |I| }       \Big[  C(q_0)   \cdot  c (\delta) \cdot c (\gamma) \cdot C(|I|) \cdot E (  \lfloor \frac{|I|}{2} \rfloor  +5)  \\
     \notag
&& \times  \frac{\eps \cdot  | \derm ( \Lie_{Z^{K}} h )  |^2 }{(1+\tau+|q|)^{3 -2\de -   c (\gamma)  \cdot c (\delta)  \cdot c(|I|) \cdot E (   \lfloor \frac{|I|}{2} \rfloor + 5)\cdot  \eps } \cdot (1+|q|)^{2+2\ga}} \Big] \cdot w(q)  \cdot d^{3}x  \cdot d\tau  \; . \\
\eea

Decomposing $h = h^1 + h^0$\;, and using Lemma \ref{Liederivativesofsphericalsymmetricpart}, we get for $M \leq \eps$\;,

 \bea
     \notag
&&    \int_{t_1}^{t_2}   \int_{\Sigma^{ext}_{\tau} }   \sum_{  |K| \leq |I| }     \Big[   C(q_0)   \cdot C(|I|) \cdot E (  \lfloor \frac{|I|}{2} \rfloor  +5)  \cdot \frac{\eps \cdot  | \Lie_{Z^{K}} H_{L  L} |^2 }{(1+\tau+|q|)^{1- 2 \de } \cdot (1+|q|)^{4+2\gamma}}  \Big] \cdot w(q)  \cdot d^{3}x  \cdot d\tau    \\
     \notag
&\leq&    \int_{t_1}^{t_2}    \int_{\Sigma^{ext}_{\tau} }   \sum_{  |K| \leq |I| }       \Big[  C(q_0)   \cdot  c (\delta) \cdot c (\gamma) \cdot C(|I|) \cdot E (  \lfloor \frac{|I|}{2} \rfloor  +5)  \\
     \notag
&& \times  \frac{\eps \cdot  | \rderm ( \Lie_{Z^{K}} h^1 )  |^2 }{(1+\tau+|q|)^{1-  2\de} \cdot (1+|q|)^{2+2\ga}} \Big] \cdot w(q)  \cdot d^{3}x  \cdot d\tau \\
     \notag
&& +  \int_{t_1}^{t_2}    \int_{\Sigma^{ext}_{\tau} }   \sum_{  |K| \leq |I| }       \Big[  C(q_0)   \cdot  c (\delta) \cdot c (\gamma) \cdot C(|I|) \cdot E (  \lfloor \frac{|I|}{2} \rfloor  +5)  \\
     \notag
&& \times  \frac{\eps \cdot  | \derm ( \Lie_{Z^{K}} h^1 )  |^2 }{(1+\tau+|q|)^{3- 2\de -  c (\gamma)  \cdot c (\delta)  \cdot c(|I|) \cdot E (   \lfloor \frac{|I|}{2} \rfloor + 5)\cdot  \eps } \cdot (1+|q|)^{2+2\ga}} \Big] \cdot w(q)  \cdot d^{3}x  \cdot d\tau  \; \Big]  \\
 \notag
&& +  \int_{t_1}^{t_2}    \int_{\Sigma^{ext}_{\tau} }     \Big[  C(q_0)   \cdot  c (\delta) \cdot c (\gamma) \cdot C(|I|) \cdot E (  \lfloor \frac{|I|}{2} \rfloor  +5)  \\
     \notag
&& \times  \frac{\eps^3 }{(1+\tau+|q|)^{5- 2\de - c (\gamma)  \cdot c (\delta)  \cdot c(|I|) \cdot E (   \lfloor \frac{|I|}{2} \rfloor + 5)\cdot  \eps } \cdot (1+|q|)^{2+2\ga}} \Big] \cdot w(q)  \cdot d^{3}x  \cdot d\tau \; . \\
\eea

Based on the same argument as in the proof of Lemma \ref{estimateonthecontributionofhzerointhesourcetermsforthewaveequationontheYangMillspoential}\,, we get that $0 < \de \leq \frac{1}{4}$\;, we have

      \bea
                 \notag
           &&         \int_{t_1 }^{t_2}   \int_{\Sigma^{ext}_{\tau} }    C(q_0)   \cdot  c (\delta) \cdot c (\gamma) \cdot C(|I|) \cdot E ( \lfloor \frac{|I|}{2} \rfloor + 5)  \\
                      \notag
           && \times \frac{\eps^3 }{(1+\tau+|q|)^{5-      c (\gamma)  \cdot c (\delta)  \cdot c(|I|) \cdot E ( \lfloor \frac{|I|}{2} \rfloor  + 5) \cdot \eps } \cdot (1+|q|)^{2+ 2 \gamma - 2\de }}    \cdot w(q) \cdot d^{3}x    \cdot d\tau  \\
         \notag
         &\leq &     C(q_0)   \cdot  c (\delta) \cdot c (\gamma) \cdot C(|I|) \cdot E ( \lfloor \frac{|I|}{2} \rfloor + 5)  \cdot  \frac{\eps^3 }{(1+t_2+|q|)^{1 -   c (\gamma)  \cdot c (\delta)  \cdot c(|I|) \cdot E ( \lfloor \frac{|I|}{2} \rfloor  + 4) \cdot \eps  } }   \; . \\
           \eea

Consequently,

 \bea
     \notag
&& \int_{t_1}^{t_2}   \int_{\Sigma^{ext}_{\tau} }   \sum_{  |K| \leq |I| }     \Big[   C(q_0)   \cdot C(|I|) \cdot E (  \lfloor \frac{|I|}{2} \rfloor  +5)  \cdot \frac{\eps \cdot  | \Lie_{Z^{K}} H_{L  L} |^2 }{(1+\tau+|q|)^{1- 2 \de } \cdot (1+|q|)^{4+2\gamma}}  \Big] \cdot w(q)  \cdot d^{3}x  \cdot d\tau  \; \Big]    \cdot dt   \\
     \notag
&\leq&  \int_{t_1}^{t_2}    \int_{\Sigma^{ext}_{\tau} }   \sum_{  |K| \leq |I| }       \Big[  C(q_0)   \cdot  c (\delta) \cdot c (\gamma) \cdot C(|I|) \cdot E (  \lfloor \frac{|I|}{2} \rfloor  +5)  \\
     \notag
&& \times  \frac{\eps \cdot  | \rderm ( \Lie_{Z^{K}} h^1 )  |^2 }{(1+\tau+|q|)^{1-  2\de} \cdot (1+|q|)^{2+2\ga}} \Big] \cdot w(q)  \cdot d^{3}x  \cdot d\tau \\
     \notag
&& +  \int_{t_1}^{t_2}    \int_{\Sigma^{ext}_{\tau} }   \sum_{  |K| \leq |I| }       \Big[  C(q_0)   \cdot  c (\delta) \cdot c (\gamma) \cdot C(|I|) \cdot E (  \lfloor \frac{|I|}{2} \rfloor  +5)  \\
     \notag
&& \times  \frac{\eps \cdot  | \derm ( \Lie_{Z^{K}} h^1 )  |^2 }{(1+\tau+|q|)^{3- 2\de -  c (\gamma)  \cdot c (\delta)  \cdot c(|I|) \cdot E (   \lfloor \frac{|I|}{2} \rfloor + 5)\cdot  \eps } \cdot (1+|q|)^{2+2\ga}} \Big] \cdot w(q)  \cdot d^{3}x  \cdot d\tau   \\
 \notag
&& + C(q_0)   \cdot  c (\delta) \cdot c (\gamma) \cdot C(|I|) \cdot E ( \lfloor \frac{|I|}{2} \rfloor + 5)  \cdot  \frac{\eps^3 }{(1+t_2+|q|)^{1 -   c (\gamma)  \cdot c (\delta)  \cdot c(|I|) \cdot E ( \lfloor \frac{|I|}{2} \rfloor  + 4) \cdot \eps  } }   \; . 
\eea

\end{proof}

  \newpage
  \subsection{The estimate on the weighted $L^2$ norm of $\derm A$}\
  
  \begin{lemma} \label{theestimateontheLtwonormofAboundedbytheenergyandonetangentialtermcontainingh1}
We have for $\ga \geq 3 \de $\,, for $0 < \de \leq \frac{1}{4}$\,, and for $\eps$ small, enough depending on $q_0$\,, on $\ga$\;, on $\de$\;, on $|I|$ and on $\mu$\;,
  \beaa
   \notag
 &&       \int_{\Sigma^{ext}_{t_2} }  |\derm ( \Lie_{Z^I}  A )  |^2     \cdot w(q)  \cdot d^{3}x     + \int_{N_{t_1}^{t_2} }  T_{\hat{L} t}^{(\bf{g})}  \cdot  w(q) \cdot dv^{(\bf{m})}_N \\
 \notag
 &&+ \int_{t_1}^{t_2}  \int_{\Sigma^{ext}_{\tau} }    | \rderm  ( \Lie_{Z^I}  A ) |^2 \cdot  \frac{\widehat{w} (q)}{(1+|q|)} \cdot  d^{3}x  \cdot dt \\
   \notag
&\les&     \sum_{|K| \leq |I |}      \int_{t_1}^{t_2}   \Big[          C(q_0)   \cdot  c (\delta) \cdot c (\gamma) \cdot C(|I|) \cdot E (\lfloor \frac{|I|}{2} \rfloor + 6)  \cdot \frac{\eps  }{(1+t)}   \cdot \E_{|K|} (t) \cdot   dt  \Big]  \\  
&& +     \sum_{|K| \leq |I | -1}      \int_{t_1}^{t_2}   \Big[          C(q_0)   \cdot  c (\delta) \cdot c (\gamma) \cdot C(|I|) \cdot E (\lfloor \frac{|I|}{2} \rfloor + 5)  \\
&& \times \frac{\eps  }{(1+t)^{1-      c (\gamma)  \cdot c (\delta)  \cdot c(|I| ) \cdot E ( \lfloor \frac{|I|}{2} \rfloor  +4) \cdot \eps } }   \cdot \E_{|K|} (t) \cdot   dt  \Big] 
\eeaa
 \beaa
  \notag
&& +      \int_{t_1}^{t_2}    \int_{\Sigma^{ext}_{\tau} }   \sum_{  |K| \leq |I| }       \Big[  C(q_0)   \cdot  c (\delta) \cdot c (\gamma) \cdot C(|I|) \cdot E (  \lfloor \frac{|I|}{2} \rfloor  +5)  \\
     \notag
&& \times  \frac{\eps \cdot  | \rderm ( \Lie_{Z^{K}} h^1 )  |^2 }{(1+t+|q|)^{1-  2\de} \cdot (1+|q|)^{2+2\ga}} \Big] \cdot w(q)  \cdot d^{3}x  \cdot dt \\
\eeaa
 \beaa
   \notag
  &&  +     \int_{t_1}^{t_2}  \int_{\Sigma^{ext}_{\tau} }   \Big(       \frac{(1+ t )}{\eps} \cdot    \sum_{|K| \leq |I|  - 1}  | g^{\la\mu} \cdot \derm_{\la}   \derm_{\mu} (  \Lie_{Z^{K}}  A ) |^2   \Big) \cdot dt \\
&&+   C(q_0)   \cdot  c (\delta) \cdot c (\gamma) \cdot C(|I|) \cdot E (  \lfloor \frac{|I|}{2} \rfloor  +5)  \cdot  \eps^3  \cdot (1+t_2)^{ c (\gamma)  \cdot c (\delta)  \cdot c(|I|) \cdot E ( \lfloor \frac{|I|}{2} \rfloor+ 5) \cdot \eps } \\
\notag
  &&   + C(q_0) \cdot c (\gamma) \cdot  C ( |I| ) \cdot E (  \lfloor \frac{|I|}{2} \rfloor  + 2)  \cdot   \eps \cdot \sum_{|J|\leq |I| + 2 }  \E_{|J|} (t_1)  \; .
\eeaa

\end{lemma}

\begin{proof}

Based on Lemma \ref{finallemmafortheestimateontheenergyfortheEinsteinYangMillspoentialusingthestructureoftheequationsforbothgoodandbadcomponents}, we have for $\ga \geq 3 \de $\,, for $0 < \de \leq \frac{1}{4}$\,, and for $\eps$ small, enough depending on $q_0$\,, on $\ga$\;, on $\de$\;, on $|I|$ and on $\mu$\;,

  \beaa
   \notag
 &&       \int_{\Sigma^{ext}_{t_2} }  |\derm ( \Lie_{Z^I}  A )  |^2     \cdot w(q)  \cdot d^{3}x     + \int_{N_{t_1}^{t_2} }  T_{\hat{L} t}^{(\bf{g})}  \cdot  w(q) \cdot dv^{(\bf{m})}_N \\
 \notag
 &&+ \int_{t_1}^{t_2}  \int_{\Sigma^{ext}_{\tau} }    | \rderm  ( \Lie_{Z^I}  A ) |^2 \cdot  \frac{\widehat{w} (q)}{(1+|q|)} \cdot  d^{3}x  \cdot d\tau \\
   \notag
&\les&     \sum_{|K| \leq |I |}      \int_{t_1}^{t_2}   \Big[          C(q_0)   \cdot  c (\delta) \cdot c (\gamma) \cdot C(|I|) \cdot E (\lfloor \frac{|I|}{2} \rfloor + 5)  \cdot \frac{\eps  }{(1+t)}   \cdot \E_{|K|} (t) \cdot   dt  \Big]  \\  
&& +     \sum_{|K| \leq |I | -1}      \int_{t_1}^{t_2}   \Big[          C(q_0)   \cdot  c (\delta) \cdot c (\gamma) \cdot C(|I|) \cdot E (\lfloor \frac{|I|}{2} \rfloor + 5)  \\
&& \times \frac{\eps  }{(1+t)^{1-      c (\gamma)  \cdot c (\delta)  \cdot c( \lfloor \frac{|I|-1}{2} \rfloor) \cdot E ( \lfloor \frac{|I|-1}{2} \rfloor+ 4) \cdot \eps } }   \cdot \E_{|K|} (t) \cdot   dt  \Big] 
\eeaa
 \beaa
\notag
  &&  +  C(q_0) \cdot c (\gamma) \cdot  C ( |I| ) \cdot E (  \lfloor \frac{|I|}{2} \rfloor  + 2) \\
  && \times   \int_{t_1}^{t_2}  \frac{\eps}{(1+t)^{1+2\de}}  \cdot \Big[ \;   \sum_{|J|\leq |I| + 2 }  \E_{|J|} (t_1)  \\
&& +   \sum_{|K| \leq |I |}   \int_{t_1}^{t}  \Big[          C(q_0)   \cdot  c (\delta) \cdot c (\gamma) \cdot C(|I|) \cdot E ( \lfloor \frac{|I|}{2} \rfloor + 6)  \cdot \frac{\eps   \cdot  \E_{|K|} (t)   }{(1+\tau )}   \Big]   \cdot d\tau  \\
      \notag
      && +   \sum_{|K| \leq |I | -1}   \int_{t_1}^{t}  \Big[          C(q_0)   \cdot  c (\delta) \cdot c (\gamma) \cdot C(|I|) \cdot E ( \lfloor \frac{|I|}{2} \rfloor + 4)  \cdot \frac{\eps   \cdot  \E_{|K|} (t)   }{(1+\tau )^{1-   c (\gamma)  \cdot c (\delta)  \cdot c(|I|) \cdot E (  \lfloor \frac{|I|}{2} \rfloor + 4)\cdot  \eps }}   \Big]   \cdot d\tau  \\
&& +   \int_{t_1}^{t}   \int_{\Sigma^{ext}_{\tau} }   \sum_{  |K| \leq |I| }        \Big[   C(q_0)   \cdot C(|I|) \cdot E (  \lfloor \frac{|I|}{2} \rfloor  +5)  \cdot \frac{\eps \cdot  | \Lie_{Z^{K}} H_{L  L} |^2 }{(1+\tau+|q|)^{1- 2 \de } \cdot (1+|q|)^{4+2\gamma}}  \Big] \cdot w(q)  \cdot d^{3}x  \cdot d\tau  \; \Big]    \cdot dt   
\eeaa
 \beaa
   \notag
  &&  +     \int_{t_1}^{t_2}  \int_{\Sigma^{ext}_{\tau} }   \Big(       \frac{(1+ t )}{\eps} \cdot    \sum_{|K| \leq |I|  - 1}  | g^{\la\mu} \cdot \derm_{\la}   \derm_{\mu} (  \Lie_{Z^{K}}  A ) |^2   \Big) \cdot dt \\
&&+   C(q_0)   \cdot  c (\delta) \cdot c (\gamma) \cdot C(|I|) \cdot E (  \lfloor \frac{|I|}{2} \rfloor  +5)  \cdot  \eps^3  \cdot (1+t_2)^{ c (\gamma)  \cdot c (\delta)  \cdot c(|I|) \cdot E ( \lfloor \frac{|I|}{2} \rfloor+ 5) \cdot \eps }\; .
\eeaa

  Since for $ \de > 0$\,, we have
\beaa
 \int_{t_1}^{t_2}  \frac{\eps}{(1+t)^{1+2\de}} dt &\leq& C \; ,
\eeaa

we obtain

  \beaa
   \notag
 &&       \int_{\Sigma^{ext}_{t_2} }  |\derm ( \Lie_{Z^I}  A )  |^2     \cdot w(q)  \cdot d^{3}x     + \int_{N_{t_1}^{t_2} }  T_{\hat{L} t}^{(\bf{g})}  \cdot  w(q) \cdot dv^{(\bf{m})}_N \\
 \notag
 &&+ \int_{t_1}^{t_2}  \int_{\Sigma^{ext}_{\tau} }    | \rderm  ( \Lie_{Z^I}  A ) |^2 \cdot  \frac{\widehat{w} (q)}{(1+|q|)} \cdot  d^{3}x  \cdot d\tau \\
   \notag
&\les&     \sum_{|K| \leq |I |}      \int_{t_1}^{t_2}   \Big[          C(q_0)   \cdot  c (\delta) \cdot c (\gamma) \cdot C(|I|) \cdot E (\lfloor \frac{|I|}{2} \rfloor + 5)  \cdot \frac{\eps  }{(1+t)}   \cdot \E_{|K|} (t) \cdot   dt  \Big]  \\  
&& +     \sum_{|K| \leq |I | -1}      \int_{t_1}^{t_2}   \Big[          C(q_0)   \cdot  c (\delta) \cdot c (\gamma) \cdot C(|I|) \cdot E (\lfloor \frac{|I|}{2} \rfloor + 5)  \\
&& \times \frac{\eps  }{(1+t)^{1-      c (\gamma)  \cdot c (\delta)  \cdot c( \lfloor \frac{|I|-1}{2} \rfloor) \cdot E ( \lfloor \frac{|I|-1}{2} \rfloor+ 4) \cdot \eps } }   \cdot \E_{|K|} (t) \cdot   dt  \Big] 
\eeaa
 \beaa
\notag
  && +   C(q_0) \cdot c (\gamma) \cdot  C ( |I| ) \cdot E (  \lfloor \frac{|I|}{2} \rfloor  + 2)  \cdot \eps \cdot \sum_{|J|\leq |I| + 2 }  \E_{|J|} (t_1)  \\
&& +   \sum_{|K| \leq |I |}   \int_{t_1}^{t_2}  \Big[          C(q_0)   \cdot  c (\delta) \cdot c (\gamma) \cdot C(|I|) \cdot E ( \lfloor \frac{|I|}{2} \rfloor + 6)  \cdot \frac{\eps   \cdot  \E_{|K|} (t)   }{(1+\tau )}   \Big]   \cdot d\tau  \\
      \notag
      && +   \sum_{|K| \leq |I | -1}   \int_{t_1}^{t_2}  \Big[          C(q_0)   \cdot  c (\delta) \cdot c (\gamma) \cdot C(|I|) \cdot E ( \lfloor \frac{|I|}{2} \rfloor + 4)  \cdot \frac{\eps   \cdot  \E_{|K|} (t)   }{(1+\tau )^{1-   c (\gamma)  \cdot c (\delta)  \cdot c(|I|) \cdot E ( \lfloor \frac{|I|}{2} \rfloor + 4)\cdot  \eps }}   \Big]   \cdot d\tau  \\
&& +   \int_{t_1}^{t_2}   \int_{\Sigma^{ext}_{\tau} }   \sum_{  |K| \leq |I| }        \Big[   C(q_0)   \cdot C(|I|) \cdot E (  \lfloor \frac{|I|}{2} \rfloor  +5)  \cdot \frac{\eps \cdot  | \Lie_{Z^{K}} H_{L  L} |^2 }{(1+\tau+|q|)^{1- 2 \de } \cdot (1+|q|)^{4+2\gamma}}  \Big] \cdot w(q)  \cdot d^{3}x  \cdot d\tau  \;   
\eeaa
 \beaa
   \notag
  &&  +     \int_{t_1}^{t_2}  \int_{\Sigma^{ext}_{\tau} }   \Big(       \frac{(1+ t )}{\eps} \cdot    \sum_{|K| \leq |I|  - 1}  | g^{\la\mu} \cdot \derm_{\la}   \derm_{\mu} (  \Lie_{Z^{K}}  A ) |^2   \Big) \cdot dt \\
&&+   C(q_0)   \cdot  c (\delta) \cdot c (\gamma) \cdot C(|I|) \cdot E (  \lfloor \frac{|I|}{2} \rfloor  +5)  \cdot  \eps^3  \cdot (1+t_2)^{ c (\gamma)  \cdot c (\delta)  \cdot c(|I|) \cdot E ( \lfloor \frac{|I|}{2} \rfloor+ 5) \cdot \eps }\; .
\eeaa

Combining the terms, we get

  \beaa
   \notag
 &&       \int_{\Sigma^{ext}_{t_2} }  |\derm ( \Lie_{Z^I}  A )  |^2     \cdot w(q)  \cdot d^{3}x     + \int_{N_{t_1}^{t_2} }  T_{\hat{L} t}^{(\bf{g})}  \cdot  w(q) \cdot dv^{(\bf{m})}_N \\
 \notag
 &&+ \int_{t_1}^{t_2}  \int_{\Sigma^{ext}_{\tau} }    | \rderm  ( \Lie_{Z^I}  A ) |^2 \cdot  \frac{\widehat{w} (q)}{(1+|q|)} \cdot  d^{3}x  \cdot d\tau \\
   \notag
&\les&     \sum_{|K| \leq |I |}      \int_{t_1}^{t_2}   \Big[          C(q_0)   \cdot  c (\delta) \cdot c (\gamma) \cdot C(|I|) \cdot E (\lfloor \frac{|I|}{2} \rfloor + 6)  \cdot \frac{\eps  }{(1+t)}   \cdot \E_{|K|} (t) \cdot   dt  \Big]  \\  
&& +     \sum_{|K| \leq |I | -1}      \int_{t_1}^{t_2}   \Big[          C(q_0)   \cdot  c (\delta) \cdot c (\gamma) \cdot C(|I|) \cdot E (\lfloor \frac{|I|}{2} \rfloor + 5)  \\
&& \times \frac{\eps  }{(1+t)^{1-      c (\gamma)  \cdot c (\delta)  \cdot c(|I| ) \cdot E ( \lfloor \frac{|I|}{2} \rfloor  +4) \cdot \eps } }   \cdot \E_{|K|} (t) \cdot   dt  \Big] 
\eeaa
 \beaa
\notag
  &&   + C(q_0) \cdot c (\gamma) \cdot  C ( |I| ) \cdot E (  \lfloor \frac{|I|}{2} \rfloor  + 2)  \cdot   \eps \cdot \sum_{|J|\leq |I| + 2 }  \E_{|J|} (t_1)  \\
  \notag
&& +   \int_{t_1}^{t_2}   \int_{\Sigma^{ext}_{\tau} }   \sum_{  |K| \leq |I| }        \Big[   C(q_0)   \cdot C(|I|) \cdot E (  \lfloor \frac{|I|}{2} \rfloor  +5)  \cdot \frac{\eps \cdot  | \Lie_{Z^{K}} H_{L  L} |^2 }{(1+\tau+|q|)^{1- 2 \de } \cdot (1+|q|)^{4+2\gamma}}  \Big] \cdot w(q)  \cdot d^{3}x  \cdot d\tau  \;   
\eeaa
 \beaa
   \notag
  &&  +     \int_{t_1}^{t_2}  \int_{\Sigma^{ext}_{\tau} }   \Big(       \frac{(1+ t )}{\eps} \cdot    \sum_{|K| \leq |I|  - 1}  | g^{\la\mu} \cdot \derm_{\la}   \derm_{\mu} (  \Lie_{Z^{K}}  A ) |^2   \Big) \cdot dt \\
&&+   C(q_0)   \cdot  c (\delta) \cdot c (\gamma) \cdot C(|I|) \cdot E (  \lfloor \frac{|I|}{2} \rfloor  +5)  \cdot  \eps^3  \cdot (1+t_2)^{ c (\gamma)  \cdot c (\delta)  \cdot c(|I|) \cdot E ( \lfloor \frac{|I|}{2} \rfloor+ 5) \cdot \eps }\; .
\eeaa

Now, based on Lemma \ref{estimateonthetermthatcontainsHLLintheenergywithestimatesinvolvingonlyhonesothattoestiamteallbyenergy}, we have for $\ga \geq 3 \de $\,, for $0 < \de \leq \frac{1}{4}$\,, and for $\eps$ small, enough depending on $q_0$\,, on $\ga$\;, on $\de$\;, and on $|I|$\;,
  
 \bea
     \notag
&&   \int_{t_1}^{t_2}   \int_{\Sigma^{ext}_{\tau} }   \sum_{  |K| \leq |I| }     \Big[   C(q_0)   \cdot C(|I|) \cdot E (  \lfloor \frac{|I|}{2} \rfloor  +5)  \cdot \frac{\eps \cdot  | \Lie_{Z^{K}} H_{L  L} |^2 }{(1+\tau+|q|)^{1- 2 \de } \cdot (1+|q|)^{4+2\gamma}}  \Big] \cdot w(q)  \cdot d^{3}x  \cdot d\tau     \\
     \notag
&\leq&   \int_{t_1}^{t_2}    \int_{\Sigma^{ext}_{\tau} }   \sum_{  |K| \leq |I| }       \Big[  C(q_0)   \cdot  c (\delta) \cdot c (\gamma) \cdot C(|I|) \cdot E (  \lfloor \frac{|I|}{2} \rfloor  +5)  \\
     \notag
&& \times  \frac{\eps \cdot  | \rderm ( \Lie_{Z^{K}} h^1 )  |^2 }{(1+\tau+|q|)^{1-  2\de} \cdot (1+|q|)^{2+2\ga}} \Big] \cdot w(q)  \cdot d^{3}x  \cdot d\tau \\
     \notag
&& +  \int_{t_1}^{t_2}    \int_{\Sigma^{ext}_{\tau} }   \sum_{  |K| \leq |I| }       \Big[  C(q_0)   \cdot  c (\delta) \cdot c (\gamma) \cdot C(|I|) \cdot E (  \lfloor \frac{|I|}{2} \rfloor  +5)  \\
     \notag
&& \times  \frac{\eps \cdot  | \derm ( \Lie_{Z^{K}} h^1 )  |^2 }{(1+\tau+|q|)^{3- 2\de -  c (\gamma)  \cdot c (\delta)  \cdot c(|I|) \cdot E (   \lfloor \frac{|I|}{2} \rfloor + 5)\cdot  \eps } \cdot (1+|q|)^{2+2\ga}} \Big] \cdot w(q)  \cdot d^{3}x  \cdot d\tau   \\
 \notag
&& + C(q_0)   \cdot  c (\delta) \cdot c (\gamma) \cdot C(|I|) \cdot E ( \lfloor \frac{|I|}{2} \rfloor + 5)  \cdot  \frac{\eps^3 }{(1+t_2+|q|)^{1 -   c (\gamma)  \cdot c (\delta)  \cdot c(|I|) \cdot E ( \lfloor \frac{|I|}{2} \rfloor  + 4) \cdot \eps  } }   \; . 
\eea

 \bea
     \notag
&\leq&   \int_{t_1}^{t_2}    \int_{\Sigma^{ext}_{\tau} }   \sum_{  |K| \leq |I| }       \Big[  C(q_0)   \cdot  c (\delta) \cdot c (\gamma) \cdot C(|I|) \cdot E (  \lfloor \frac{|I|}{2} \rfloor  +5)  \\
     \notag
&& \times  \frac{\eps \cdot  | \rderm ( \Lie_{Z^{K}} h^1 )  |^2 }{(1+\tau+|q|)^{1-  2\de} \cdot (1+|q|)^{2+2\ga}} \Big] \cdot w(q)  \cdot d^{3}x  \cdot d\tau \\
     \notag
&& +  \int_{t_1}^{t_2}    \int_{\Sigma^{ext}_{\tau} }   \sum_{  |K| \leq |I| }       \Big[  C(q_0)   \cdot  c (\delta) \cdot c (\gamma) \cdot C(|I|) \cdot E (  \lfloor \frac{|I|}{2} \rfloor  +5)  \\
     \notag
&& \times  \frac{\eps  }{(1+\tau+|q|)^{3- 2\de -  c (\gamma)  \cdot c (\delta)  \cdot c(|I|) \cdot E (   \lfloor \frac{|I|}{2} \rfloor + 5)\cdot  \eps } \cdot (1+|q|)^{2+2\ga}} \Big] \cdot \E_{|K|} (t)  \cdot w(q)  \cdot d^{3}x  \cdot d\tau   \\
 \notag
&& + C(q_0)   \cdot  c (\delta) \cdot c (\gamma) \cdot C(|I|) \cdot E ( \lfloor \frac{|I|}{2} \rfloor + 5)  \cdot  \frac{\eps^3 }{(1+t_2+|q|)^{1 -   c (\gamma)  \cdot c (\delta)  \cdot c(|I|) \cdot E ( \lfloor \frac{|I|}{2} \rfloor  + 4) \cdot \eps  } }   \; . 
\eea

Injecting this estimate, we obtain the desired result.

\end{proof}

\section{Estimate on the full energy of the Einstein-Yang-Mills fields}

We can estimate the source terms for $h^1$ in a similar and easier manner than what we did for $A$ which had more troublesome structure, where, however, for $h^1$\;, we have an additional term (in the source terms of the wave equation on $h^1$) that is different, namely $ g^{\alpha\beta} \derm_\alpha \derm_\beta h^0 $ for which we have the following lemma.

\begin{lemma} \label{theestimateonthesphericallysymmetricpartcarryingthemass}
We have for $ M \leq  \eps^2 \leq \eps \leq 1 $\; that
\bea
\notag
&& \int_{t_1}^{t_2} \int_{\Sigma^{ext}_{\tau} } | \Lie_{Z^I} \big( g^{\alpha\beta} \derm_\alpha \derm_\beta h^0 \big) | \cdot |\derm ( \Lie_{Z^I}  h^1 ) | w\, d^{3}x d t  \\
\notag
&\leq&   E (|I|) \cdot C(|I|)  \cdot \eps^2  + E (|I|) \cdot C(|I|)\cdot  \eps^2 \cdot \!\! \sum_{|K|\leq |I|} \int_{t_1}^{t_2}   \int_{\Sigma^{ext}_{\tau} }   \frac{ \E_{|K|} (t)    }{(1+t) } dt  \; . \\
\eea

\end{lemma}

\begin{proof}
 Based on the proof of Lemma 11.4 in \cite{LR10}, we have for $M \leq \eps^2 $\;, 
\bea\label{theestimateofLindbladRodninaksionthesphericallysymmetricshcwarzchildianpart}
\notag
&& \int_{t_1}^{t_2} \int_{\Sigma^{ext}_{\tau} } | \Lie_{Z^I} \big( g^{\alpha\beta} \derm_\alpha \derm_\beta h^0 \big) | \cdot |\derm ( \Lie_{Z^I}  h^1 ) | w\, d^{3}x dt  \\
\notag
&\leq& E (|I|) \cdot C(|I|)\cdot \, \eps^2 \cdot \!\! \sum_{|K|\leq |I|}  \int_{t_1}^{t_2} \int_{\Sigma^{ext}_{\tau} } \frac{| \derm ( \Lie_{Z^K}   h^1) |^2}{(1+t)^2} w \, d^{3}x dt  \\
\notag
&& +  E (|I|)\cdot C(|I|) \cdot  \, \eps^2 \cdot \!\! \sum_{|K|\leq |I|} \int_{t_1}^{t_2} \!\!\Big(\int_{\Sigma^{ext}_{\tau} } | \derm ( \Lie_{Z^K}  h^1 )|^2 w\, d^{3}x \Big)^{1/2}\!\!\!\frac{dt }{(1+t)^{3/2}} \; . \\
\eea

We estimate
\beaa
&& \Big(\int_{\Sigma^{ext}_{\tau} } | \derm ( \Lie_{Z^K}  h^1 )|^2 w\, d^{3}x \Big)^{1/2}\!\!\!\frac{ dt}{(1+t)^{3/2}} \\
& =& \Big(\int_{\Sigma^{ext}_{\tau} }   \frac{ | \derm ( \Lie_{Z^K}  h^1 )|^2 }{(1+t) }  w\, d^{3}x \Big)^{1/2}\!\!\!\frac{ dt}{(1+t)}  \\
&\leq&  \int_{\Sigma^{ext}_{\tau} }   \frac{ | \derm ( \Lie_{Z^K}  h^1 )|^2 }{(1+t) }  w\, d^{3}x    + \!\!\!\frac{ dt}{(1+t)^2 }  \\
&& \text{(where we used $a \cdot b \les a^2 + b^2$)} .
\eeaa
Hence,

\beaa
 && E (|I|) \cdot C(|I|) \cdot \eps^2 \cdot \sum_{|K|\leq |I|} \int_{t_1}^{t_2} \!\!\Big(\int_{\Sigma^{ext}_{\tau} } | \derm ( \Lie_{Z^K}  h^1 )|^2 w\, d^{3}x \Big)^{1/2}\!\!\!\frac{ dt}{(1+t)^{3/2}}  \\
 &\leq& E (|I|)\cdot C(|I|) \cdot   \eps^2 \cdot  \sum_{|K|\leq |I|} \int_{t_1}^{t_2}   \int_{\Sigma^{ext}_{\tau} }   \frac{ | \derm ( \Lie_{Z^K}  h^1 )|^2 }{(1+t) }  w\, d^{3}x   + E (|I|) \cdot \, \eps^2 \cdot \!\! \sum_{|K|\leq |I|} \int_{t_1}^{t_2} \!\!\!\frac{ dt}{(1+t)^2 } \\
 &\leq& E (|I|) \cdot C(|I|)\cdot \eps^2 \cdot \sum_{|K|\leq |I|} \int_{t_1}^{t_2}   \int_{\Sigma^{ext}_{\tau} }   \frac{ | \derm ( \Lie_{Z^K}  h^1 )|^2 }{(1+t) }  w\, d^{3}x   +  \frac{ E (|I|) \cdot  \, \eps^2 \cdot C(|I|)  }{(1+t_1 ) } \\
 &\leq&   E (|I|) \cdot C(|I|)  \cdot \eps^2  + E (|I|) \cdot C(|I|)\cdot  \eps^2 \cdot \!\! \sum_{|K|\leq |I|} \int_{t_1}^{t_2}   \int_{\Sigma^{ext}_{\tau} }   \frac{ \E_{|K|} (t)    }{(1+t) }  dt  \; .
\eeaa
Thus, injecting in \eqref{theestimateofLindbladRodninaksionthesphericallysymmetricshcwarzchildianpart}, we get the stated result.
 \end{proof}

 \begin{lemma}\label{theenergyestimateonallcomponentsoftheEinstein-Yang-Millsfieldskeepingtermsonthelefthandsidethatareconstructedtobepositive}

We have for $\ga \geq 3 \de $\,, for $0 < \de \leq \frac{1}{4}$\,, and for $\eps$ small, enough depending on $q_0$\,, on $\ga$\;, on $\de$\;, on $|I|$ and on $\mu$\;, and for $M \leq \eps^2 \leq 1$\;, that for $|I| \geq \lfloor \frac{|I|}{2} \rfloor + 6$\;, 

  \beaa
   \notag
 &&     \E_{|I|} (t_2)     + \int_{N_{t_1}^{t_2} }  T_{\hat{L} t}^{(\bf{g})}  \cdot  w(q) \cdot dv^{(\bf{m})}_N \\
 \notag
 &&+ \int_{t_1}^{t_2}  \int_{\Sigma^{ext}_{\tau} }  \big(   | \rderm  ( \Lie_{Z^I}  A ) |^2  +  | \rderm  ( \Lie_{Z^I}  h^1 ) |^2 \big) \cdot  \frac{\widehat{w} (q)}{(1+|q|)} \cdot  d^{3}x  \cdot dt \\
   \notag
&\les&     \sum_{|K| \leq |I |}      \int_{t_1}^{t_2}   \Big[          C(q_0)   \cdot  c (\delta) \cdot c (\gamma) \cdot C(|I|) \cdot \big(    E (|I|  ) +E (  \lfloor \frac{|I|}{2} \rfloor  +6)   \big)   \cdot \frac{\eps  }{(1+t)}   \cdot \E_{|K|} (t) \cdot   dt  \Big]  \\  
&& +     \sum_{|K| \leq |I | -1}      \int_{t_1}^{t_2}   \Big[          C(q_0)   \cdot  c (\delta) \cdot c (\gamma) \cdot C(|I|) \cdot E (\lfloor \frac{|I|}{2} \rfloor + 5)  \\
&& \times \frac{\eps  }{(1+t)^{1-      c (\gamma)  \cdot c (\delta)  \cdot c(|I| ) \cdot E ( \lfloor \frac{|I|}{2} \rfloor  +4) \cdot \eps } }   \cdot \E_{|K|} (t) \cdot   dt  \Big] 
\eeaa
 \beaa
   \notag
  &&  +     \int_{t_1}^{t_2}  \int_{\Sigma^{ext}_{\tau} }   \Big(       \frac{(1+ t )}{\eps} \cdot    \sum_{|K| \leq |I|  - 1}  | g^{\la\mu} \cdot \derm_{\la}   \derm_{\mu} (  \Lie_{Z^{K}}  A ) |^2   \Big) \cdot dt \\
  \notag
  &&  +     \int_{t_1}^{t_2}  \int_{\Sigma^{ext}_{\tau} }   \Big(       \frac{(1+ t )}{\eps} \cdot    \sum_{|K| \leq |I|  - 1}  | g^{\la\mu} \cdot \derm_{\la}   \derm_{\mu} (  \Lie_{Z^{K}}  h^1 ) |^2   \Big) \cdot dt \\
&&+   C(q_0)   \cdot  c (\delta) \cdot c (\gamma) \cdot C(|I|) \cdot E ( |I|  )  \cdot  \eps^2  \cdot (1+t_2)^{ c (\gamma)  \cdot c (\delta)  \cdot c(|I|) \cdot E ( \lfloor \frac{|I|}{2} \rfloor+ 5) \cdot \eps } \\
\notag
  &&   + C(q_0) \cdot c (\gamma) \cdot  C ( |I| ) \cdot E (  \lfloor \frac{|I|}{2} \rfloor  + 2)  \cdot  \eps \cdot \sum_{|J|\leq |I| + 2 }  \E_{|J|} (t_1)  \; .
\eeaa

\end{lemma}

\begin{proof}

 Based on Lemma \ref{actualusefulstrzuctureofthesourcetermsforthewaveequationontheMetricsmallhoneusingbootstrap}, and considering the fact that we dealt already with the “bad” term for the sources of the metric in Subsection \ref{Structurebadsourcetermsmetricinssubsection} and more precisely  in Lemma \ref{Thebadtermforthemetricisestimatehereforetheclosureofenergyestimate}, we get a similar estimate $h^1$ to the one we obtained for $A$ in Lemma \ref{theestimateontheLtwonormofAboundedbytheenergyandonetangentialtermcontainingh1}, except that we have here an additional term arising from $h^0$\;, that we already estimated in Lemma \ref{theestimateonthesphericallysymmetricpartcarryingthemass}.
 
Thus, from Lemma \ref{theestimateontheLtwonormofAboundedbytheenergyandonetangentialtermcontainingh1}, and adding the similar estimate that we obtain for $h^1$\;, along with the additional term estimated in Lemma \ref{theestimateonthesphericallysymmetricpartcarryingthemass}, we get that for $\ga \geq 3 \de $\,, for $0 < \de \leq \frac{1}{4}$\,, and for $\eps$ small, enough depending on $q_0$\,, on $\ga$\;, on $\de$\;, on $|I|$ and on $\mu$\;, we obtain

  \beaa
   \notag
 &&     \E_{|I|} (t_2)     + \int_{N_{t_1}^{t_2} }  T_{\hat{L} t}^{(\bf{g})}  \cdot  w(q) \cdot dv^{(\bf{m})}_N \\
 \notag
 &&+ \int_{t_1}^{t_2}  \int_{\Sigma^{ext}_{\tau} }  \big(   | \rderm  ( \Lie_{Z^I}  A ) |^2  +  | \rderm  ( \Lie_{Z^I}  h^1 ) |^2 \big) \cdot  \frac{\widehat{w} (q)}{(1+|q|)} \cdot  d^{3}x  \cdot dt \\
   \notag
&\les&     \sum_{|K| \leq |I |}      \int_{t_1}^{t_2}   \Big[          C(q_0)   \cdot  c (\delta) \cdot c (\gamma) \cdot C(|I|) \cdot E (\lfloor \frac{|I|}{2} \rfloor + 6)  \cdot \frac{\eps  }{(1+t)}   \cdot \E_{|K|} (t) \cdot   dt  \Big]  \\  
&& +     \sum_{|K| \leq |I | -1}      \int_{t_1}^{t_2}   \Big[          C(q_0)   \cdot  c (\delta) \cdot c (\gamma) \cdot C(|I|) \cdot E (\lfloor \frac{|I|}{2} \rfloor + 5)  \\
&& \times \frac{\eps  }{(1+t)^{1-      c (\gamma)  \cdot c (\delta)  \cdot c(|I| ) \cdot E ( \lfloor \frac{|I|}{2} \rfloor  +4) \cdot \eps } }   \cdot \E_{|K|} (t) \cdot   dt  \Big] 
\eeaa
 \beaa
  \notag
&& +      \int_{t_1}^{t_2}    \int_{\Sigma^{ext}_{\tau} }   \sum_{  |K| \leq |I| }       \Big[  C(q_0)   \cdot  c (\delta) \cdot c (\gamma) \cdot C(|I|) \cdot E (  \lfloor \frac{|I|}{2} \rfloor  +5)  \\
     \notag
&& \times  \frac{\eps \cdot  | \rderm ( \Lie_{Z^{K}} h^1 )  |^2 }{(1+t+|q|)^{1-  2\de} \cdot (1+|q|)^{2+2\ga}} \Big] \cdot w(q)  \cdot d^{3}x  \cdot dt
\eeaa
 \beaa
   \notag
  &&  +     \int_{t_1}^{t_2}  \int_{\Sigma^{ext}_{\tau} }   \Big(       \frac{(1+ t )}{\eps} \cdot    \sum_{|K| \leq |I|  - 1}  | g^{\la\mu} \cdot \derm_{\la}   \derm_{\mu} (  \Lie_{Z^{K}}  A ) |^2   \Big) \cdot dt \\
&&+   C(q_0)   \cdot  c (\delta) \cdot c (\gamma) \cdot C(|I|) \cdot E (  \lfloor \frac{|I|}{2} \rfloor  +5)  \cdot  \eps^3  \cdot (1+t_2)^{ c (\gamma)  \cdot c (\delta)  \cdot c(|I|) \cdot E ( \lfloor \frac{|I|}{2} \rfloor+ 5) \cdot \eps } \\
\notag
  &&   + C(q_0) \cdot c (\gamma) \cdot  C ( |I| ) \cdot E (  \lfloor \frac{|I|}{2} \rfloor  + 2)  \cdot  \eps \cdot \sum_{|J|\leq |I| + 2 }  \E_{|J|} (t_1)  \\
&& +   E (|I|) \cdot C(|I|)  \cdot \eps^2  + E (|I|) \cdot C(|I|)\cdot  \eps^2 \cdot \!\! \sum_{|K|\leq |I|} \int_{t_1}^{t_2}   \int_{\Sigma^{ext}_{\tau} }   \frac{ \E_{|K|} (t)    }{(1+t) }  w\, dt
\eeaa

Combining the terms, and absorbing the term containing the tangential derivatives into the left hand side of the inequality for $\eps$ small depending $q_0$\;, on $\delta$\;, on $\gamma$\;, on $|I|$\;, we get

  \beaa
   \notag
 &&     \E_{|I|} (t_2)     + \int_{N_{t_1}^{t_2} }  T_{\hat{L} t}^{(\bf{g})}  \cdot  w(q) \cdot dv^{(\bf{m})}_N \\
 \notag
 &&+ \int_{t_1}^{t_2}  \int_{\Sigma^{ext}_{\tau} }  \big(   | \rderm  ( \Lie_{Z^I}  A ) |^2  +  | \rderm  ( \Lie_{Z^I}  h^1 ) |^2 \big) \cdot  \frac{\widehat{w} (q)}{(1+|q|)} \cdot  d^{3}x  \cdot dt \\
   \notag
&\les&     \sum_{|K| \leq |I |}      \int_{t_1}^{t_2}   \Big[          C(q_0)   \cdot  c (\delta) \cdot c (\gamma) \cdot C(|I|) \cdot \big(    E (|I|  ) +E (  \lfloor \frac{|I|}{2} \rfloor  +6)   \big)   \cdot \frac{\eps  }{(1+t)}   \cdot \E_{|K|} (t) \cdot   dt  \Big]  \\  
&& +     \sum_{|K| \leq |I | -1}      \int_{t_1}^{t_2}   \Big[          C(q_0)   \cdot  c (\delta) \cdot c (\gamma) \cdot C(|I|) \cdot E (\lfloor \frac{|I|}{2} \rfloor + 5)  \\
&& \times \frac{\eps  }{(1+t)^{1-      c (\gamma)  \cdot c (\delta)  \cdot c(|I| ) \cdot E ( \lfloor \frac{|I|}{2} \rfloor  +4) \cdot \eps } }   \cdot \E_{|K|} (t) \cdot   dt  \Big] 
\eeaa
 \beaa
   \notag
  &&  +     \int_{t_1}^{t_2}  \int_{\Sigma^{ext}_{\tau} }   \Big(       \frac{(1+ t )}{\eps} \cdot    \sum_{|K| \leq |I|  - 1}  | g^{\la\mu} \cdot \derm_{\la}   \derm_{\mu} (  \Lie_{Z^{K}}  A ) |^2   \Big) \cdot dt \\
  \notag
  &&  +     \int_{t_1}^{t_2}  \int_{\Sigma^{ext}_{\tau} }   \Big(       \frac{(1+ t )}{\eps} \cdot    \sum_{|K| \leq |I|  - 1}  | g^{\la\mu} \cdot \derm_{\la}   \derm_{\mu} (  \Lie_{Z^{K}}  h^1 ) |^2   \Big) \cdot dt \\
&&+   C(q_0)   \cdot  c (\delta) \cdot c (\gamma) \cdot C(|I|) \cdot E (  \lfloor \frac{|I|}{2} \rfloor  +5)  \cdot  \eps^3  \cdot (1+t_2)^{ c (\gamma)  \cdot c (\delta)  \cdot c(|I|) \cdot E ( \lfloor \frac{|I|}{2} \rfloor+ 5) \cdot \eps } \\
\notag
  &&   + C(q_0) \cdot c (\gamma) \cdot  C ( |I| ) \cdot E (  \lfloor \frac{|I|}{2} \rfloor  + 2)  \cdot  \eps \cdot \sum_{|J|\leq |I| + 2 }  \E_{|J|} (t_1)  \\
&& +   E (|I|) \cdot C(|I|)  \cdot \eps^2   \; .
\eeaa
Hence, we get the stated result.
\end{proof}

\subsection{The closure of the bootstrap argument on the energy}\

 \begin{lemma}\label{Theinequalityontheenergytobeusedtoapplygronwallrecursively}

We have for $\ga \geq 3 \de $\,, for $0 < \de \leq \frac{1}{4}$\,, and for $\eps$ small, enough depending on $q_0$\,, on $\ga$\;, on $\de$\;, on $|I|$ and on $\mu$\;, and for $M \leq \eps^2 \leq 1$\;, that 
     \beaa
   \notag
 &&     \E_{|I|} (t_2) \\
   \notag
&\les&    C(q_0)   \cdot  c (\delta) \cdot c (\gamma) \cdot C(|I|) \cdot \big(    E (|I|  ) +E (  \lfloor \frac{|I|}{2} \rfloor  +6)   \big)     \cdot  \Big[  \sum_{|K| \leq |I |}      \int_{t_1}^{t_2}        \frac{\eps  }{(1+t)}   \cdot \E_{|K|} (t) \cdot   dt   \\  
&& +     \sum_{|K| \leq |I | -1}      \int_{t_1}^{t_2}         \frac{\eps  }{(1+t)^{1-      c (\gamma)  \cdot c (\delta)  \cdot c(|I| ) \cdot E ( \lfloor \frac{|I|}{2} \rfloor  +4) \cdot \eps } }   \cdot \E_{|K|} (t) \cdot   dt  \\
&&+   \eps^2  \cdot (1+t_2)^{ c (\gamma)  \cdot c (\delta)  \cdot c(|I|) \cdot E ( \lfloor \frac{|I|}{2} \rfloor+ 5) \cdot \eps }   +   \eps \cdot \E_{|I|+2} (t_1)  \Big]   \; .
\eeaa

\end{lemma}

   Based on Lemma \ref{theenergyestimateonallcomponentsoftheEinstein-Yang-Millsfieldskeepingtermsonthelefthandsidethatareconstructedtobepositive}, and using that  $T_{\hat{L} t}^{(\bf{g})} \geq 0$ by construction, and considering that $|I| \geq \lfloor \frac{|I|}{2} \rfloor + 6$\;, we get
        \beaa
   \notag
 &&     \E_{|I|} (t_2) \\
   \notag
&\les&    C(q_0)   \cdot  c (\delta) \cdot c (\gamma) \cdot C(|I|) \cdot E (|I|  )    \cdot  \Big[  \sum_{|K| \leq |I |}      \int_{t_1}^{t_2}        \frac{\eps  }{(1+t)}   \cdot \E_{|K|} (t) \cdot   dt   \\  
&& +     \sum_{|K| \leq |I | -1}      \int_{t_1}^{t_2}         \frac{\eps  }{(1+t)^{1-      c (\gamma)  \cdot c (\delta)  \cdot c(|I| ) \cdot E ( \lfloor \frac{|I|}{2} \rfloor  +4) \cdot \eps } }   \cdot \E_{|K|} (t) \cdot   dt  \\
&&+   \eps^2  \cdot (1+t_2)^{ c (\gamma)  \cdot c (\delta)  \cdot c(|I|) \cdot E ( \lfloor \frac{|I|}{2} \rfloor+ 5) \cdot \eps }   +  \eps \cdot \E_{|I|+2} (t_1)  \Big]  \\
   \notag
  &&  +     \int_{t_1}^{t_2}  \int_{\Sigma^{ext}_{\tau} }   \Big(       \frac{(1+ t )}{\eps} \cdot    \sum_{|K| \leq |I|  - 1}  | g^{\la\mu} \cdot \derm_{\la}   \derm_{\mu} (  \Lie_{Z^{K}}  A ) |^2   \Big) \cdot dt \\
  \notag
  &&  +     \int_{t_1}^{t_2}  \int_{\Sigma^{ext}_{\tau} }   \Big(       \frac{(1+ t )}{\eps} \cdot    \sum_{|K| \leq |I|  - 1}  | g^{\la\mu} \cdot \derm_{\la}   \derm_{\mu} (  \Lie_{Z^{K}}  h^1 ) |^2   \Big) \cdot dt  \; .
\eeaa

Since the lower order source terms  $\Big(       \frac{(1+ t )}{\eps} \cdot    \sum_{|K| \leq |I|  - 1}  | g^{\la\mu} \cdot \derm_{\la}   \derm_{\mu} (  \Lie_{Z^{K}}  A ) |^2   \Big)$ and $\Big(       \frac{(1+ t )}{\eps} \cdot    \sum_{|K| \leq |I|  - 1}  | g^{\la\mu} \cdot \derm_{\la}   \derm_{\mu} (  \Lie_{Z^{K}}  h^1 ) |^2   \Big)$ have the same structure as the higher order terms that we already estimated and appear on the left hand side of the ineuqlaity, we can therefore, estimate these terms by the other terms on the right hand side of the inequality. Consequently, we get the stated result.

    \begin{lemma}\label{Theboundontheenergybytimetbyepsilonandboostrapassumptionsandinitialdata}

We have for $\ga \geq 3 \de $\,, for $0 < \de \leq \frac{1}{4}$\,, and for $\eps$ small, enough depending on $q_0$\,, on $\ga$\;, on $\de$\;, on $|I|$ and on $\mu$\;, and for $M \leq \eps^2 \leq 1$\;, that

  \bea
   \notag
     \E_{|I|} (t_2) &\les&    (1+t_2)^{ C(q_0)   \cdot  c (\delta) \cdot c (\gamma) \cdot C( |I|   )  \cdot \big(    E (|I| ) +E (  \lfloor \frac{|I| }{2} \rfloor  +6)   \big)      \cdot \eps } \\
        \notag
     && \times C(q_0)   \cdot  c (\delta) \cdot c (\gamma) \cdot C( |I|   )  \cdot   \big(    E (|I| ) +E (  \lfloor \frac{|I| }{2} \rfloor  +6)   \big)  \cdot   \Big[  \eps     +  \E_{|I|+2} (t_1)  \Big]  \; .\\
\eea

\end{lemma}

\begin{proof}
We have from Lemma \ref{Theinequalityontheenergytobeusedtoapplygronwallrecursively}, a recursive formula in $|I|$\;, of Grönwall type inequality, on the energy $\E_{|I|}$\;. We are going to use this recursive formula to close the argument. First, this gives that for $|I| =0$\;, that the sum over $|I | -1$ does not exist, and therefore,

     \beaa
   \notag
 &&     \E_{0} (t_2) \\
   \notag
&\les&    C(q_0)   \cdot  c (\delta) \cdot c (\gamma) \cdot C\cdot E ( 6 )    \cdot  \Big[       \int_{t_1}^{t_2}   \Big[       \frac{\eps  }{(1+t)}   \cdot \E_{0} (t) \cdot   dt   +   \eps^2  \cdot (1+t_2)^{ c (\gamma)  \cdot c (\delta)  \cdot c \cdot E (5) \cdot \eps }   +  \eps \cdot \E_{2} (t_1)  \Big]  \; .
\eeaa
Applying Grönwall lemma, we obtain
     \bea\label{initialestimateontheenergyintherecursiveformula}
   \notag
     &&\E_{0} (t_2)\\
     \notag
      &\les&    (1+t_2)^{ C(q_0)   \cdot  c (\delta) \cdot c (\gamma) \cdot C\cdot E ( 6  )    \cdot \eps } \cdot   C(q_0)   \cdot  c (\delta) \cdot c (\gamma) \cdot C\cdot E ( 6 )    \\
           \notag
      && \times  \Big[  \eps^2  \cdot (1+t_2)^{ c (\gamma)  \cdot c (\delta)  \cdot c \cdot E (5) \cdot \eps }   +   \eps \cdot \E_{2} (t_1)  \Big]  \\
      \notag
     &\les&   (1+t_2)^{ C(q_0)   \cdot  c (\delta) \cdot c (\gamma) \cdot C\cdot E ( 6  )    \cdot \eps } \cdot   C(q_0)   \cdot  c (\delta) \cdot c (\gamma) \cdot C\cdot E ( 6 )   \\
           \notag
     && \times  \Big[  \eps^2   +  \eps \cdot \E_{2} (t_1)  \Big]  \; . \\
\eea

Assuming now that for all $|I| \leq  k \in \N$, we have 

  \bea\label{assumptiononthewaytheenergygrwsfromtherecursiveformula}
   \notag
     \E_{|I|} (t_2) &\les&    (1+t_2)^{ C(q_0)   \cdot  c (\delta) \cdot c (\gamma) \cdot C( |I|   )  \cdot \big(    E (|I| ) +E (  \lfloor \frac{|I| }{2} \rfloor  +6)   \big)      \cdot \eps } \\
        \notag
     && \times C(q_0)   \cdot  c (\delta) \cdot c (\gamma) \cdot C( |I|   )  \cdot   \big(    E (|I| ) +E (  \lfloor \frac{|I| }{2} \rfloor  +6)   \big)  \cdot   \Big[  \eps^2     +   \eps \cdot\E_{|I|+2} (t_1)  \Big]  \; .\\
\eea

Injecting this estimate in the inequality from Lemma \ref{Theinequalityontheenergytobeusedtoapplygronwallrecursively} to estimate the lower order terms where $|I| - 1 \leq k$\;, we obtain

     \beaa
   \notag
 &&     \E_{k+1} (t_2) \\
   \notag
&\les&    C(q_0)   \cdot  c (\delta) \cdot c (\gamma) \cdot C(k+1) \cdot \big(    E (k+1 ) +E (  \lfloor \frac{k+1}{2} \rfloor  +6)   \big)    \cdot  \Big[     \int_{t_1}^{t_2}       \frac{\eps  }{(1+t)}   \cdot \E_{k+1} (t) \cdot   dt   \\  
&& +     \sum_{|K| \leq k }      \int_{t_1}^{t_2}         \frac{\eps  }{(1+t)^{1-      c (\gamma)  \cdot c (\delta)  \cdot c(k+1 ) \cdot E ( \lfloor \frac{k+1}{2} \rfloor  +4) \cdot \eps } }   \cdot \E_{|K|} (t) \cdot   dt  \\
&&+   \eps^2  \cdot (1+t_2)^{ c (\gamma)  \cdot c (\delta)  \cdot c(k+1) \cdot E ( \lfloor \frac{k+1}{2} \rfloor+ 5) \cdot \eps }   +  \eps \cdot \E_{k+1+2} (t_1)  \Big]   \; . \\
&\les&    C(q_0)   \cdot  c (\delta) \cdot c (\gamma) \cdot C(k+1) \cdot \big(    E (k +1) +E (  \lfloor \frac{k+1}{2} \rfloor  +6)   \big)      \cdot  \Big[     \int_{t_1}^{t_2}         \frac{\eps  }{(1+t)}   \cdot \E_{k+1} (t) \cdot   dt   \\  
&& +        \int_{t_1}^{t_2}         \frac{\eps  }{(1+t)^{1-      c (\gamma)  \cdot c (\delta)  \cdot c(k +1) \cdot E ( \lfloor \frac{k+1}{2} \rfloor  +4) \cdot \eps } }  \cdot C(q_0)   \cdot  c (\delta) \cdot c (\gamma) \cdot C(k) \cdot   \big(    E (k ) +E (  \lfloor \frac{k}{2} \rfloor  +6)   \big)     \\
&& \times \Big( (1+t_2)^{ C(q_0)   \cdot  c (\delta) \cdot c (\gamma) \cdot C(k) \cdot   \big(    E (k ) +E (  \lfloor \frac{k}{2} \rfloor  +6)   \big)    \cdot \eps } \cdot    \Big[  \eps^2    +   \eps \cdot \E_{k+2} (t_1)  \Big]   \cdot   dt  \Big)  \\
&&+   \eps^2  \cdot (1+t_2)^{ c (\gamma)  \cdot c (\delta)  \cdot c(k+1) \cdot E ( \lfloor \frac{k+1}{2} \rfloor+ 5) \cdot \eps }   +  \eps \cdot \E_{k+3} (t_1)  \Big]    \\
&& \text{(based on our assumption \eqref{assumptiononthewaytheenergygrwsfromtherecursiveformula} on the estimate on the energy derived from the recursive formula)}.
\eeaa

Consequently, we get

     \beaa
   \notag
 &&     \E_{k+1} (t_2) \\
   \notag
&\les&    C(q_0)   \cdot  c (\delta) \cdot c (\gamma) \cdot C(k+1) \cdot \big(    E (k+1 ) +E (  \lfloor \frac{k+1}{2} \rfloor  +6)   \big)     \cdot      \int_{t_1}^{t_2}         \frac{\eps  }{(1+t)}   \cdot \E_{k+1} (t) \cdot   dt   \\  
&& +       C(q_0)   \cdot  c (\delta) \cdot c (\gamma) \cdot C(k+1) \cdot \big(    E (k+1 ) +E (  \lfloor \frac{k+1}{2} \rfloor  +6)   \big)  \\
&&  \times     \int_{t_1}^{t_2}         \frac{\eps \cdot C(q_0)   \cdot  c (\delta) \cdot c (\gamma) \cdot C(k) \cdot   \big(    E (k ) +E (  \lfloor \frac{k}{2} \rfloor  +6)   \big)  }{(1+t)^{1-      c (\gamma)  \cdot c (\delta)  \cdot c(k+1 ) \cdot E ( \lfloor \frac{k+1}{2} \rfloor  +4) \cdot \eps } }  \cdot   dt \\
&& \times   (1+t_2)^{ C(q_0)   \cdot  c (\delta) \cdot c (\gamma) \cdot C(k) \cdot \big(    E (k ) +E (  \lfloor \frac{k}{2} \rfloor  +6)   \big)     \cdot \eps } \cdot    \Big[  \eps^2    +   \eps \cdot \E_{k+2} (t_1)  \Big]     \\
&&+   C(q_0)   \cdot  c (\delta) \cdot c (\gamma) \cdot C(k+1) \cdot \big(    E (k+1 ) +E (  \lfloor \frac{k+1}{2} \rfloor  +6)   \big)    \\
&& \times   \Big[ \eps^2  \cdot (1+t_2)^{ c (\gamma)  \cdot c (\delta)  \cdot c(k+1) \cdot E ( \lfloor \frac{k+1}{2} \rfloor+ 5) \cdot \eps }   +  \eps \cdot \E_{k+3} (t_1)  \Big]    
\eeaa

     \beaa
   \notag
&\les&    C(q_0)   \cdot  c (\delta) \cdot c (\gamma) \cdot C(k+1) \cdot \big(    E (k+1 ) +E (  \lfloor \frac{k+1}{2} \rfloor  +6)   \big)      \cdot      \int_{t_1}^{t_2}         \frac{\eps  }{(1+t)}   \cdot \E_{k+1} (t) \cdot   dt   \\  
&& +       C(q_0)   \cdot  c (\delta) \cdot c (\gamma) \cdot C(k+1) \cdot \big(    E (k+1 ) +E (  \lfloor \frac{k+1}{2} \rfloor  +6)   \big)   \\
&& \times       \eps \cdot (1+t_2 )^{   c (\gamma)  \cdot c (\delta)  \cdot c(k +1) \cdot E ( \lfloor \frac{k+1}{2} \rfloor  +4) \cdot \eps }    \cdot (1+t_2)^{ C(q_0)   \cdot  c (\delta) \cdot c (\gamma) \cdot C(k) \cdot  \big(    E (k ) +E (  \lfloor \frac{k}{2} \rfloor  +6)   \big)   \cdot \eps } \\
&& \times    \Big[  \eps^2   +  \eps \cdot \E_{k+2} (t_1)  \Big]    \\
&&+   C(q_0)   \cdot  c (\delta) \cdot c (\gamma) \cdot C(k+1) \cdot  \big(    E (k+1 ) +E (  \lfloor \frac{k+1}{2} \rfloor  +6)   \big)   \\
&& \times   \Big[ \eps^2  \cdot (1+t_2)^{ c (\gamma)  \cdot c (\delta)  \cdot c(k+1) \cdot E ( \lfloor \frac{k+1}{2} \rfloor+ 5) \cdot \eps }   + \eps \cdot  \E_{k+3} (t_1)  \Big]    \; .
\eeaa
Hence,
     \beaa
   \notag
 &&     \E_{k+1} (t_2) \\
   \notag
&\les&    C(q_0)   \cdot  c (\delta) \cdot c (\gamma) \cdot C(k+1) \cdot \big(    E (k+1 ) +E (  \lfloor \frac{k+1}{2} \rfloor  +6)   \big)      \cdot      \int_{t_1}^{t_2}         \frac{\eps  }{(1+t)}   \cdot \E_{k+1} (t) \cdot   dt   \\  
&& +       C(q_0)   \cdot  c (\delta) \cdot c (\gamma) \cdot C(k+1) \cdot \big(    E (k+1 ) +E (  \lfloor \frac{k+1}{2} \rfloor  +6)   \big)   \cdot           \eps   \\
&& \times  (1+t_2)^{ C(q_0)   \cdot  c (\delta) \cdot c (\gamma) \cdot C(k+1) \cdot  \big(    E (k+1 ) +E (  \lfloor \frac{k+1}{2} \rfloor  +6)   \big)   \cdot \eps } \cdot    \Big[  \eps^2    +  \eps \cdot \E_{k+2} (t_1)  \Big]    \\
&&+   C(q_0)   \cdot  c (\delta) \cdot c (\gamma) \cdot C(k+1) \cdot  \big(    E (k+1 ) +E (  \lfloor \frac{k+1}{2} \rfloor  +6)   \big)   \\
&& \times   \Big[ \eps^2  \cdot (1+t_2)^{ c (\gamma)  \cdot c (\delta)  \cdot c(k+1) \cdot E ( \lfloor \frac{k+1}{2} \rfloor+ 5) \cdot \eps }   +  \eps \cdot \E_{k+3} (t_1)  \Big]    \\
&\les&    C(q_0)   \cdot  c (\delta) \cdot c (\gamma) \cdot C(k+1) \cdot \big(    E (k+1 ) +E (  \lfloor \frac{k+1}{2} \rfloor  +6)   \big)      \cdot      \int_{t_1}^{t_2}         \frac{\eps  }{(1+t)}   \cdot \E_{k+1} (t) \cdot   dt   \\  
&& +       C(q_0)   \cdot  c (\delta) \cdot c (\gamma) \cdot C(k+1) \cdot \big(    E (k+1 ) +E (  \lfloor \frac{k+1}{2} \rfloor  +6)   \big)      \\
&& \times  (1+t_2)^{ C(q_0)   \cdot  c (\delta) \cdot c (\gamma) \cdot C(k+1) \cdot  \big(    E (k+1 ) +E (  \lfloor \frac{k+1}{2} \rfloor  +6)   \big)   \cdot \eps } \cdot    \Big[  \eps^2    +  \eps \cdot \E_{k+3} (t_1)  \Big]    \; .
\eeaa

Thus, by applying Grönwall lemma again, we get that

  \bea
   \notag
     \E_{k+1} (t_2) &\les&    (1+t_2)^{ C(q_0)   \cdot  c (\delta) \cdot c (\gamma) \cdot C( k+1  )  \cdot \big(    E (k+1 ) +E (  \lfloor \frac{k+1 }{2} \rfloor  +6)   \big)      \cdot \eps } \\
        \notag
     && \times C(q_0)   \cdot  c (\delta) \cdot c (\gamma) \cdot C(k+1   )  \cdot   \big(    E (k+1 ) +E (  \lfloor \frac{k+1 }{2} \rfloor  +6)   \big)  \\
          \notag
     && \times \Big[  \eps^2     +  \eps \cdot \E_{k+3} (t_1)  \Big]  \; ,\\
\eea
which upgrades our assumption \eqref{assumptiononthewaytheenergygrwsfromtherecursiveformula} on $|I| \leq  k$\;, to actually prove it for all $|I| \leq  k +1$\;. This, along with the initialisation \eqref{initialestimateontheenergyintherecursiveformula}, proves the stated result.

\end{proof}

\begin{proposition} \label{THEpropositiontoclosetheargumenttoboundtheenergy}

Let $N \geq 11$\;. We have for $\ga \geq 3 \de $\,, for $0 < \de \leq \frac{1}{4}$\,, for $\eps$ small, enough depending on $q_0$\,, on $\ga$\;, on $\de$\;, on $N$ and on $\mu < 0$\;, and for $M \leq \eps^2 \leq 1$\;, and for $
\overline{ \E}_{N+2} \leq \eps$ (defined in \eqref{definitionoftheenergynormforinitialdata}), that under the bootstrap assumption \eqref{bootstrap}, we have
  \beaa
   \notag
     \E_{N} (t) &\leq&  \frac{E(N) }{2} \cdot \eps  \cdot  (1+t)^{ \de }  \; .
\eeaa
\end{proposition}

\begin{proof}
 
From Lemma \ref{Theboundontheenergybytimetbyepsilonandboostrapassumptionsandinitialdata}, we have
  \beaa
   \notag
     \E_{|I|} (t) &\les&    (1+t)^{ C(q_0)   \cdot  c (\delta) \cdot c (\gamma) \cdot C( |I|   )  \cdot \big(    E (|I| ) +E (  \lfloor \frac{|I| }{2} \rfloor  +6)   \big)      \cdot \eps } \\
        \notag
     && \times C(q_0)   \cdot  c (\delta) \cdot c (\gamma) \cdot C( |I|   )  \cdot   \big(    E (|I| ) +E (  \lfloor \frac{|I| }{2} \rfloor  +6)   \big)  \cdot   \Big[  \eps^2     +   \eps \cdot \E_{|I|+2} (0)  \Big]  \; .\\
\eeaa
Hence, for $ |I|  \geq \lfloor \frac{|I|}{2} \rfloor + 6$\;, we get

  \beaa
   \notag
     \E_{|I|} (t) &\les&    (1+t)^{ C(q_0)   \cdot  c (\delta) \cdot c (\gamma) \cdot C( |I|   )  \cdot    E (|I| )     \cdot \eps } \\
        \notag
     && \times C(q_0)   \cdot  c (\delta) \cdot c (\gamma) \cdot C( |I|   )  \cdot  E (|I| ) \cdot   \Big[  \eps^2     +  \eps \cdot \E_{|I|+2} (0)  \Big]  \; .\\
\eeaa
By choosing $\eps$ small enough depending on $q_0$\;, on $\delta$\;, on $\gamma$\;, on $ E (|I| )$ and on $|I| $\;, we have 
\bea\label{choice1forepsilonsmall}
 C(q_0)   \cdot  c (\delta) \cdot c (\gamma) \cdot C( |I|   )  \cdot    E (|I| )     \cdot \eps &<&  \de \; , 
\eea
and 
\bea\label{choicetwoforepsilonsmall}
  C(q_0)   \cdot  c (\delta) \cdot c (\gamma) \cdot C( |I|   )  \cdot  E (|I| ) \cdot     \eps^2   &<& \frac{E(|I|)}{4} \cdot \eps \; .
\eea

Also, by choosing the initial data $\overline{ \E}_{|I|+2}$ small enough, say such that
\bea
\overline{ \E}_{|I|+2} \leq \eps \;, 
\eea
then this implies that 
\bea
 \E_{|I|+2} (0) \leq \eps \;.
\eea
Thus, 
\bea
\notag
 && C(q_0)   \cdot  c (\delta) \cdot c (\gamma) \cdot C( |I|   )  \cdot  E (|I| ) \cdot    \eps \cdot \E_{|I|+2} (0) \\
 \notag
 &\leq&  C(q_0)   \cdot  c (\delta) \cdot c (\gamma) \cdot C( |I|   )  \cdot  E (|I| ) \cdot     \eps^2  \\
 \notag
   &<& \frac{E(|I|)}{4} \cdot \eps \\
   && \text{(based on our choice for smallness of $\eps$ in \eqref{choicetwoforepsilonsmall}).} 
\eea

Consequently, for such $\eps$ small, we have for all $ |I|  \geq \lfloor \frac{|I|}{2} \rfloor + 6$\;,
  \beaa
   \notag
     \E_{|I|} (t) &\leq&  \frac{ E(|I|) }{2} \cdot \eps  \cdot  (1+t)^{ \de }  \; .
\eeaa
Since for $N \geq 11$\;, we have $N \geq \lfloor \frac{N}{2} \rfloor + 6$, we therefore get the stated result.

\end{proof}

\end{document}